\documentclass{amsart}
\usepackage{amsmath,amssymb,amsthm}
\usepackage{amsfonts,amssymb}
\usepackage{graphicx}
\usepackage[inkscapelatex=false]{svg}
\usepackage{enumerate}
\usepackage{enumitem}
\usepackage{tikz}
\usetikzlibrary{decorations.pathreplacing}
\usepackage{tikz-cd}
\tikzcdset{scale cd/.style={every label/.append style={scale=#1},
    cells={nodes={scale=#1}}}}
\usepackage{stmaryrd}
\usepackage{tabularx}
\usepackage{diagbox}
\usepackage{commath}
\usepackage{xcolor}
\usepackage{bbm}
\usepackage{wrapfig}
\usepackage[T1]{fontenc}
\usepackage{float}
\usepackage{color}
\usepackage{verbatim}
\usepackage{caption}
\usepackage{subcaption}
\usepackage{mathrsfs}
\usepackage{mathtools}
\usepackage{yhmath}
\usepackage[all]{xy}
\usepackage{soul}
\usepackage{hyperref}
\usepackage{nameref}
\usepackage{marvosym}
\usepackage{esvect}
\usepackage{geometry}
\usepackage{tcolorbox}
    \tcbuselibrary{skins}
    \tcbuselibrary{breakable}

\DeclareRobustCommand{\SkipTocEntry}[5]{}

\let\emptyset\varnothing

\usetikzlibrary{fit}

\tikzset{overcross/.style={double, line width=1.5, white, double=#1, double distance=\knotlinewidth},
    overcross/.default={black},
    knot/.style={line width=\knotlinewidth, baseline=-.5ex}}
\newcommand{\knotlinewidth}{.7pt}

\makeatletter
\def\fixtikzforbreqn#1#2{%
  \protected\edef#1{\noexpand\ifmmode\mathchar\the\mathcode`#2 \noexpand\else#2\noexpand\fi}%
}
\fixtikzforbreqn\tikz@nonactivesemicolon;
\fixtikzforbreqn\tikz@nonactivecolon:
\fixtikzforbreqn\tikz@nonactivebar|
\fixtikzforbreqn\tikz@nonactiveexlmark!
\makeatother

\usetikzlibrary{calc}

\theoremstyle{plain}

\newtheorem{theorem}{Theorem}[section]

\newtheorem*{theorem*}{Theorem}
\newtheorem{lemma}[theorem]{Lemma}
\newtheorem{proposition}[theorem]{Proposition}

\newtheorem{corollary}[theorem]{Corollary}

\theoremstyle{definition}
\newtheorem{definition}[theorem]{Definition}

\theoremstyle{remark}
\newtheorem{remark}[theorem]{Remark}

\newtheorem*{example*}{Example}

\numberwithin{equation}{section}

\newtheoremstyle{named}{}{}{\itshape}{}{\bfseries}{.}{.5em}{\thmnote{#3}}
\theoremstyle{named}
\newtheorem*{namedtheorem}{Theorem}

\tikzset{
    partial ellipse/.style args={#1:#2:#3}{
        insert path={+ (#1:#3) arc (#1:#2:#3)}
    }
}

\newcommand*\circled[1]{\tikz[baseline=(char.base)]{
            \node[shape=circle,draw,inner sep=2pt] (char) {#1};}}

\tikzset{commutative diagrams/.cd,
mysymbol1/.style={start anchor=center,end anchor=center,draw=none}
}
\newcommand\MySymbA[2][\circled{1}]{%
  \arrow[mysymbol1]{#2}[description]{#1}}
  
\tikzset{commutative diagrams/.cd,
mysymbol2/.style={start anchor=center,end anchor=center,draw=none}
}
\newcommand\MySymbB[2][\circled{2}]{%
  \arrow[mysymbol2]{#2}[description]{#1}}

\tikzset{commutative diagrams/.cd,
mysymbol3/.style={start anchor=center,end anchor=center,draw=none}
}
\newcommand\MySymbC[2][\circled{3}]{%
  \arrow[mysymbol3]{#2}[description]{#1}}

\tikzset{commutative diagrams/.cd,
mysymbol4/.style={start anchor=center,end anchor=center,draw=none}
}
\newcommand\MySymbD[2][\circled{4}]{%
  \arrow[mysymbol4]{#2}[description]{#1}}

\tikzset{commutative diagrams/.cd,
mysymbol5/.style={start anchor=center,end anchor=center,draw=none}
}
\newcommand\MySymbAA[2][\circled{1'}]{%
  \arrow[mysymbol5]{#2}[description]{#1}}
  
\tikzset{commutative diagrams/.cd,
mysymbol6/.style={start anchor=center,end anchor=center,draw=none}
}
\newcommand\MySymbBB[2][\circled{2'}]{%
  \arrow[mysymbol6]{#2}[description]{#1}}

\tikzset{commutative diagrams/.cd,
mysymbol7/.style={start anchor=center,end anchor=center,draw=none}
}
\newcommand\MySymbCC[2][\circled{3'}]{%
  \arrow[mysymbol7]{#2}[description]{#1}}

\tikzset{commutative diagrams/.cd,
mysymbol8/.style={start anchor=center,end anchor=center,draw=none}
}
\newcommand\MySymbDD[2][\circled{4'}]{%
  \arrow[mysymbol8]{#2}[description]{#1}}

\usetikzlibrary{decorations.pathmorphing}

\makeatletter
\newcommand{\math@param}[3]{%
  \fontdimen#3
  \ifx#1\displaystyle\textfont#2
  \else\ifx#1\textstyle\textfont#2
  \else\ifx#1\scriptstyle\scriptfont#2
  \else\scriptscriptfont#2 \fi\fi\fi
}

\newdimen\@my@yshift
\NewDocumentCommand{\myvec}{m}{\mathord{\mathpalette\@myvec{#1}}}

\newcommand*{\@myvec}[2]{%
  \begin{tikzpicture}[baseline=(n.base)]
  \node (n) [inner sep=0]
    {\global\@my@yshift=-0.5\math@param{#1}{2}{5}%
     $\m@th #1#2$};
    \draw[decorate, decoration={snake, amplitude=0.3pt, segment length=2pt}]
      ([yshift=\@my@yshift]n.south west) -- ([yshift=\@my@yshift]n.south east);
  \end{tikzpicture}%
}
\makeatother

\title{Jones-Wenzl projectors and odd Khovanov homology}
\author{Dean Spyropoulos}
\date{}

\begin{document}
\begin{abstract}
The Jones-Wenzl projectors are particular elements of the Temperley-Lieb algebra essential to the construction of quantum 3-manifold invariants. As a first step toward categorifying quantum 3-manifold invariants, Cooper and Krushkal categorified these projectors. In another direction, Naisse and Putyra gave a categorification of the Temperley-Lieb algebra compatible with odd Khovanov homology, introducing new machinery called grading categories. We provide a generalization of Naisse and Putyra's work in the spirit of Bar-Natan's canopolies or Jones's planar algebras, replacing grading categories with grading multicategories. We use our setup to prove the existence and uniqueness of categorified Jones-Wenzl projectors in odd Khovanov homology. This result quickly implies the existence of a new, ``odd'' categorification of the colored Jones polynomial.
\end{abstract}
\maketitle
\tableofcontents

\section{Introduction}

The Temperley-Lieb algebras, $TL_n$, are diagrammatic algebras originating from operator algebra theory which entered low-dimensional quantum topology with the construction of the Jones polynomial via representations of the braid group \cite{MR908150}. Elements of particular importance are special idempotents of the Temperley-Lieb algebra, $p_n \in TL_n$, called Jones-Wenzl projectors. These projectors have been studied extensively, and they are vital to the construction of the colored Jones polynomials and the skein theoretic construction of the Witten-Reshetikhin-Turaev 3-manifold invariants (cf. \cite{MR1472978}, Chapter 13).

In \cite{10.1215/S0012-7094-00-10131-7}, Khovanov provided a homological invariant of links whose graded Euler characteristic $\chi$ was the Jones polynomial, initiating the study of categorification. Since then, a major motivating question has been whether Khovanov's categorification can be extended to a categorification of quantum 3-manifold invariants. It would stand to reason that the first step in replicating the procedures of the decategorified setting would be to construct categorical lifts of the Jones-Wenzl projectors, living in some categorification of the Temperley-Lieb algebra.

A categorification of the Jones-Wenzl projectors was achieved by Cooper and Krushkal in \cite{https://doi.org/10.48550/arxiv.1005.5117}. First, Bar-Natan \cite{Bar_Natan_2005} provided a categorification of the Temperley-Lieb algebra in the sense that he constructed a category $\mathrm{Kom}(n)$ whose Grothendieck group $K_0$ was isomorphic to $TL_n$. Cooper and Krushkal then prove the existence of objects $P_n^{\text{CK}}$ of $\mathrm{Kom}(n)$ which satisfy $[P_n^{\text{CK}}] = p_n$, for $[P_n^{\text{CK}}]$ the equivalence class of $P_n^{\mathrm{CK}}$ in $K_0(\mathrm{Kom}(n))$. Said another way, $\chi(P_n^{\mathrm{CK}}) = p_n$. Rozansky \cite{rozansky2010infinitetorusbraidyields} has also given a construction of categorified projectors using the Khovanov complex associated to an infinite torus braid. For recent progress toward the categorification of quantum 3-manifold invariants from Khovanov homology, see \cite{2022arXiv221005640H}.

In this paper, we initiate an investigation of similar phenomena for a different categorification of the Jones polynomial, called odd Khovanov homology. Suppose $L$ is a link. In \cite{ozsvath2003heegaard}, Ozsv\'ath and Szab\'o constructed a spectral sequence converging to the Heegaard Floer homology of the double branched cover of $L$, $\widehat{\mathrm{HF}}(\Sigma(-L); \mathbb{Z}/2\mathbb{Z})$, with $E_2$ page the (reduced) Khovanov homology of $L$, $\widetilde{\mathrm{Kh}}(L; \mathbb{Z}/2\mathbb{Z})$. In an attempt to lift the spectral sequence to $\mathbb{Z}$ coefficients, Ozsv\'ath, Rasmussen, and Szab\'o realized that the $E_2$ page could no longer be ordinary reduced Khovanov homology. Instead, they produced a new candidate, another homological link invariant categorifying the Jones polynomial, closely related to Khovanov's construction (indeed, necessarily agreeing over $\mathbb{Z}/2\mathbb{Z}$ coefficients).

Ozsv\'ath, Rasmussen, and Szab\'o's new construction \cite{Ozsv_th_2013} is called \textit{odd Khovanov homology}, which we denote by $\mathrm{Kh}_o$; to avoid confusion, the original theory of \cite{10.1215/S0012-7094-00-10131-7} has been retroactively declared \textit{even Khovanov homology}, denoted $\mathrm{Kh}_e$. While agreeing in $\mathbb{Z}/2\mathbb{Z}$ coefficients, there exist pairs of links $L_1 \not= L_2$ for which $\mathrm{Kh}_e(L_1; \mathbb{Z}) \cong \mathrm{Kh}_e(L_2; \mathbb{Z})$, but $\mathrm{Kh}_o(L_1; \mathbb{Z})\not\cong \mathrm{Kh}_o(L_2; \mathbb{Z})$, and vice-versa; see \cite{shumakovitch2018patternsoddkhovanovhomology}. We remark that spectral sequences from odd Khovanov homology to flavors of Floer homology have been discovered: Daemi \cite{daemi2015abelian} showed that there is a spectral sequence from odd Khovanov homology to the plane Floer homology of the double branched cover, and Scaduto \cite{MR3394316} showed that another spectral sequence starting at odd Khovanov homology converges to the framed instanton homology of the double branched cover.

Recall that even Khovanov homology is built from a functor $\mathcal{F}_e$ with source the category whose objects are closed 1-manifolds and whose morphisms are embedded cobordisms, and a target category of $\mathbb{K}$-modules, for some ring $\mathbb{K}$. In the literature, a functor of this form is called a ($1+1$)-dimensional TQFT. Likewise, the original definition of odd Khovanov homology is built from a (perhaps misleadingly named) ``projective TQFT''---that is, a TQFT well-defined only up to sign---of embedded cobordisms. Indeed, the TQFT of \cite{Ozsv_th_2013}, which we will denote by $\mathcal{F}_o$, depends on some additional information. Using notation which will be introduced later (\S \ref{ss:chronologicalcobordisms and coc}), this is pictured as
\begin{equation}
\label{eq:framedsplits}
\mathcal{F}_o\left(
\tikz[baseline={([yshift=-.5ex]current bounding box.center)}, scale=.5]{
	\draw  (1,2) .. controls (1,3) and (0,3) .. (0,4);
	\draw  (2,2) .. controls (2,3) and (3,3) .. (3,4);
	\draw (1,4) .. controls (1,3) and (2,3) .. (2,4);
	\draw (0,4) .. controls (0,3.75) and (1,3.75) .. (1,4);
	\draw (0,4) .. controls (0,4.25) and (1,4.25) .. (1,4);
	\draw (2,4) .. controls (2,3.75) and (3,3.75) .. (3,4);
	\draw (2,4) .. controls (2,4.25) and (3,4.25) .. (3,4);
	\draw (1,2) .. controls (1,1.75) and (2,1.75) .. (2,2);
	\draw[dashed] (1,2) .. controls (1,2.25) and (2,2.25) .. (2,2);
    \node at (0.5, 4.5) {$a_i$};
    \node at (2.5, 4.5) {$a_{i+1}$};
    \node at (1.5, 1.5) {$a$};
        \draw[<-] (1.8,3.7) -- (1.2,2.8);
}\right)
= -
\mathcal{F}_o\left(
\tikz[baseline={([yshift=-.5ex]current bounding box.center)}, scale=.5]{
	\draw  (1,2) .. controls (1,3) and (0,3) .. (0,4);
	\draw  (2,2) .. controls (2,3) and (3,3) .. (3,4);
	\draw (1,4) .. controls (1,3) and (2,3) .. (2,4);
	\draw (0,4) .. controls (0,3.75) and (1,3.75) .. (1,4);
	\draw (0,4) .. controls (0,4.25) and (1,4.25) .. (1,4);
	\draw (2,4) .. controls (2,3.75) and (3,3.75) .. (3,4);
	\draw (2,4) .. controls (2,4.25) and (3,4.25) .. (3,4);
	\draw (1,2) .. controls (1,1.75) and (2,1.75) .. (2,2);
	\draw[dashed] (1,2) .. controls (1,2.25) and (2,2.25) .. (2,2);
    \node at (0.5, 4.5) {$a_i$};
    \node at (2.5, 4.5) {$a_{i+1}$};
    \node at (1.5, 1.5) {$a$};
        \draw[->] (1.8,3.7) -- (1.2,2.8);
}\right).
\end{equation}
Moreover, $\mathcal{F}_o$ is known to be sensitive to the exchange of critical points in embedded cobordisms between 1-manifolds.

Putyra, first in his Master's thesis \cite{putyra2010cobordismschronologiesgeneralisationkhovanov} and then in \cite{putyra20152categorychronologicalcobordismsodd}, introduced a refinement of the source category so that $\mathcal{F}_o$ may be improved to a genuine functor. By a \textit{chronological cobordism}, we mean a cobordism endowed with a framed Morse function, called a \textit{chronology}, separating critical points; see \S \ref{ss:chronologicalcobordisms and coc}. The chronology induces an orientation on each unstable manifold of index 1 and 2 critical points, which we draw as an arrow (as shown in (\ref{eq:framedsplits}) for the index 1 case). Consequently, $\mathcal{F}_o$ is upgraded to a genuine functor: the equality above is reinterpreted as a relation between the maps on modules associated with two distinct chronological cobordisms. Going forward, functors from a category of chronological cobordisms to the category of $\mathbb{K}$-modules will be called \textit{chronological TQFTs}. Also introduced in \cite{putyra20152categorychronologicalcobordismsodd} is the notion of a \textit{unified} Khovanov complex, which is a complex over the ground ring
\[
R = \mathbb{Z}[X, Y, Z^{\pm1}] \big/ (X^2 = Y^2 = 1).
\]
The homology of this complex is called \textit{unified} (also called \textit{covering} or \textit{generalized} in the literature) \textit{Khovanov homology}. The unified Khovanov complex has the incredibly desirable feature of specializing to the even theory if one sets $X = Y = Z = 1$, and to the odd theory by setting $X = Z = 1$ and $Y = -1$. We use $\mathcal{F}$ (see \S \ref{ss:unified arc algebras}) to denote the chronological TQFT for unified Khovanov homology.

\addtocontents{toc}{\SkipTocEntry}
\subsection{Unified projectors}

As in Cooper and Krushkal's work, our projectors will live in a categorification of the Temperley-Lieb algebra, which we denote by $\mathrm{Chom}(n)_R^\mathscr{G}$. Specifically, $\mathrm{Chom}(n)_R^\mathscr{G}$ is the category of $\mathscr{G}$-graded $H^n$-modules all of whose entries come from flat diskular tangles (see Figure \ref{fig:multigluingsvg} for an example of a non-flat diskular tangle). The algebra $H^n$ is the $n$th unified arc algebra; we review Khovanov's arc algebras in \S \ref{sss:Khovanov arc algebra} and unified arc algebras in \S \ref{ss:unified arc algebras}. The notation ``$\mathrm{Chom}$'' is meant to impress that we think of this category like the category ``$\mathrm{Kom}$'' of \cite{Bar_Natan_2005}, but with chronological cobordisms present. The notation $\mathscr{G}$ refers to the new grading essential to this paper; we defer an introduction to $\mathscr{G}$ momentarily. The $\mathscr{G}$-grading determines an integral $q$-grading (see \S \ref{sss:qcollapse}). We let $K_0^q$ denote the Grothendieck group which remembers only the $q$-grading and not the whole $\mathscr{G}$-grading information. Then, $\mathrm{Chom}(n)_R^\mathscr{G}$ categorifies $TL_n$ in the sense that
\[
K_0^q(\mathrm{Chom}(n)_R^\mathscr{G}) \cong TL_n
\]
as $\mathbb{Z}[q,q^{-1}]$-algebras; see Definition \ref{def:chom} of \S \ref{ss:operationsviamultigluing}.

Specializing the ground ring $R$ by $X, Y, Z = 1$ defines a forgetful functor from $\mathrm{Chom}(n)_R^\mathscr{G}$ to the category $\mathrm{Kom}(H^n\mathrm{PMod})$, another categorification of $TL_n$ compatible with even Khovanov homology. This is the categorification of Khovanov, provided in \cite{Khovanov_2002}, using \textit{projective} $H^n$-modules. Indeed, we will see that the $\mathscr{G}$-grading is not essential to the even case---the objects of $\mathrm{Kom}(H^n\mathrm{PMod})$ are not $\mathscr{G}$-graded. Likewise, specializing by $X, Z = 1$ and $Y=-1$ induces a forgetful functor from $\mathrm{Chom}(n)_R^\mathscr{G}$ to what we'll denote by $\mathrm{Chom}(n)_o^\mathscr{G}$, a categorification of $TL_n$ implicit in the work of Naisse and Putyra. We call these the even and odd forgetful functors, and denote them by $\pi_e$ and $\pi_o$ respectively. Notice that the $\mathbb{Z}/2\mathbb{Z}$-reductions of both $\mathrm{Kom}(H^n\mathrm{PMod})$ and $\mathrm{Chom}(n)_o^\mathscr{G}$ agree; we denote by $\mathrm{Kom}(H^n\mathrm{PMod})_{\mathbb{Z}/2\mathbb{Z}}$ the corresponding category. The $\mathscr{G}$-grading is also nonessential to the $\mathbb{Z}/2\mathbb{Z}$-reduction. We'll denote the corresponding forgetful functors by $\mathfrak{f}$. Then $\mathfrak{f} \circ \pi_e = \mathfrak{f}\circ \pi_o$; i.e., the diagram

\[
\begin{tikzcd}
    & \mathrm{Chom}(n)_R^\mathscr{G} \arrow[dr, "\pi_o"] \arrow[dl, "\pi_e"'] & \\
    \mathrm{Kom}(H^n\mathrm{PMod}) \arrow[dr] & &\mathrm{Chom}(n)_o^\mathscr{G} \arrow[dl] \\
    & \mathrm{Kom}(H^n\mathrm{PMod}_{\mathbb{Z}/2\mathbb{Z}}) & 
\end{tikzcd}
\]
commutes. The following is proven in \S \ref{Chapter:OddCKProj} as a combination of Proposition \ref{prop:unifiedprojectoruni} and Theorem \ref{thm:existenceofprojectors}.

\begin{namedtheorem}[Theorem A]
There exist categorifications of the Jones-Wenzl projectors, called \textit{unified projectors}, $P_n$ in $\mathrm{Chom}(n)_R^\mathscr{G}$, which are unique up to chain-homotopy equivalence. By a categorification, we mean that $[P_n] \in K_0^q(\mathrm{Chom}(n)_R^\mathscr{G})$ is equal to $p_n\in TL_n$ (for a complete description, see Definition \ref{def:unifiedprojector}). On one hand, $\pi_e(P_n)$ is a categorified projector in $\mathrm{Kom}(H^n\mathrm{PMod})$, and $\pi_e(P_n) = P_n^{\mathrm{CK}}$. On the other, under the odd forgetful functor, $\pi_o(P_n)$ is a new categorification of the $n$th Jones-Wenzl projector in $\mathrm{Chom}(n)_o^\mathscr{G}$. They both agree after reduction to $\mathbb{Z}/2\mathbb{Z}$-coefficients: $\mathfrak{f}(P_n^o) = \mathfrak{f}(P_n^{\mathrm{CK}})$.
\end{namedtheorem}

We will write $P_n^o$ to denote $\pi_o(P_n)$. We remark that Cooper and Krushkal's projectors actually live in Bar-Natan's category $\mathrm{Kom}(n)$, but it is known that this category is equivalent to Khovanov's categorification of $TL_n$, $\mathrm{Kom}(H^n\mathrm{PMod})$.

Following Section 6.4 of \cite{naisse2020odd}, we define a (diskular) tangle invariant $\mathrm{Kh}_q$ in \S \ref{ss:TangleInvariant} which specializes to unified Khovanov homology when the tangle is a closed link. The caveat is that $\mathrm{Kh}_q$ lives in a category $\mathrm{Chom}(n)_R^q$; in general, $\mathrm{Kh}$ is not a tangle invariant in the category $\mathrm{Chom}(n)_R^\mathscr{G}$, so we must ``collapse'' the $\mathscr{G}$-grading to an integral $q$-grading (see \S \ref{sss:qcollapse}---the term ``collapse'' is slightly misleading). Regardless, by construction $\mathrm{Kh}_q$ specializes to the even Khovanov tangle invariant, denoted $\mathrm{Kh}_q^e$, along with an odd Khovanov tangle invariant $\mathrm{Kh}_q^o$.

\begin{figure}[ht]
\centering
\includesvg[width=7cm]{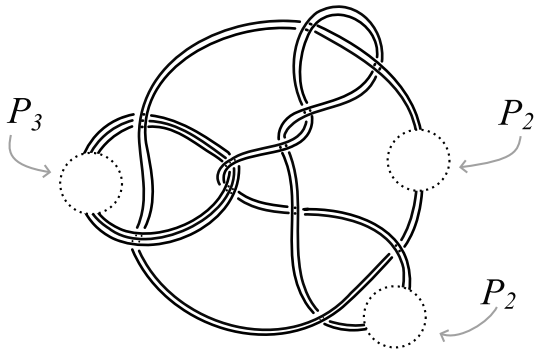}
\caption{Schematic for $\Pi^\mathbf{m}(L)$, where $L$ is the 3-component link L11n314 of the Thistelthwaite link table, and $\mathbf{m} = (3,2,2)$.}
\label{fig:coloredcomplex}
\end{figure}

In analogy with Section 5 of \cite{https://doi.org/10.48550/arxiv.1005.5117}, the existence of these tangle invariants, together with Theorem A, is immediately useful. Namely, as the Jones-Wenzl projectors are vital to the construction of the colored Jones polynomials $J(L;\mathbf{m})(q)$, the existence of categorified projectors quickly implies the existence of a categorification of the colored Jones polynomial. Using the new categorification of the Jones-Wenzl projectors (compatible with odd Khovanov homology), we construct a new, ``odd'' categorification of the colored Jones polynomial. First, if $L$ is an $n$-component link and $\mathbf{m} = (m_1,\ldots, m_n) \in \mathbb{N}^n$, denote by $T_L^\mathbf{m}$ the result of taking $m_i$ parallel copies of the $i$th component of $L$ for each $i=1,\ldots, n$ and then removing a small diskular region from each of the original components (see Figure \ref{fig:coloredcomplex}). Then, set
\[
\Pi^\mathbf{m}(L) := (P_{m_1}, \ldots, P_{m_n}) \otimes_{(H^{m_1}, \ldots, H^{m_n})} \mathrm{Kh}_q(T_L^{\mathbf{m}})
\]
where each of the $P_{m_i}$ is viewed as an object of $\mathrm{Chom}(m_i)_R^q$. This has the effect of inserting projectors into the tangle diagram $T_L^\mathbf{m}$; again, consult Figure \ref{fig:coloredcomplex} for a schematic. See \S \ref{ss:flexiblegluing} for introductory remarks regarding this tensor product.

Let $\Pi_e^\mathbf{m}(L)$ and $\Pi_o^\mathbf{m}(L)$ denote the complexes obtained by specializing $R$ by $X=Y=Z = 1$, and $X = Z = 1$, $Y=-1$ respectively. We call each of $\Pi^\mathbf{m}(L)$, $\Pi_e^\mathbf{m}(L)$, and $\Pi_o^\mathbf{m}(L)$ the \textit{unified, even, and odd $\mathbf{m}$-colored Khovanov complexes} of $L$, respectively. Finally, we define the \textit{unified, even, and odd  $\mathbf{m}$-colored Khovanov} or \textit{link homologies} of $L$ to be the homology of these complexes; we denote them by $\mathcal{H}(L; \mathbf{m})$, $\mathcal{H}_e(L; \mathbf{m})$, and $\mathcal{H}_o(L; \mathbf{m})$ respectively. We emphasize that we define the even and odd colored link homology by specializing $R$ \textit{before} taking homology. Also, notice that $\mathcal{H}(L; \mathbf{m})$ has coefficients in $R$, while $\mathcal{H}_e(L; \mathbf{m})$ and $\mathcal{H}_o(L; \mathbf{m})$ have coefficients in $\mathbb{Z}$.

Let $\chi_q$ denote the graded Euler characteristic which records only the $q$-grading associated to a particular $\mathscr{G}$-grading or $\mathscr{G}$-grading shift. Then, the following is proven in \S \ref{ss:unifiedcoloredhomology}. 

\begin{namedtheorem}[Theorem B]
For any colored link $(L;\mathbf{m})$, the chain-homotopy equivalence type 
 of the $\mathbf{m}$-colored Khovanov complex $\Pi^\mathbf{m}(L)$ is an invariant of $(L;\mathbf{m})$. Thus, the $\mathbf{m}$-colored Khovanov homologies $\mathcal{H}(L; \mathbf{m})$, $\mathcal{H}_e(L; \mathbf{m})$, and $\mathcal{H}_o(L; \mathbf{m})$ are invariants of $(L;\mathbf{m})$. Moreover, the even and odd homologies categorify the colored Jones polynomial in the sense that 
\[
\chi_q(\mathcal{H}_e(L; \mathbf{m})) = J(L;\mathbf{m})(q) = \chi_q(\mathcal{H}_o(L; \mathbf{m})).
\]
On one hand, $\mathcal{H}_e(L; \mathbf{m})$ is the colored link homology of Cooper and Krushkal. However, there are colored links $(L;\mathbf{m})$ for which $\mathcal{H}_o(L; \mathbf{m}) \not= \mathcal{H}_e(L; \mathbf{m})$, so we obtain a new categorification of the colored Jones polynomial.
\end{namedtheorem}

To see that the two categorifications are distinct, we compute the unified Khovanov homology of the full trace of $P_2$ (see \S \ref{sss:homologyoftrace} and, in particular, Equation (\ref{eq:2coloredunknot})), which coincides with the unified colored link homology of the 2-colored unknot. We obtain the even and odd colored link homologies of the 2-colored unknot by taking homology after specializing the complex of Equation (\ref{eq:2coloredunknotcomplex}) to the even and odd settings. See Table \ref{tab:evenodd} for the even (left) and odd (right) colored link homologies of the 2-colored unknot, where we have expressed the homology in terms of quantum grading $q$ and homological grading $h$. 
\begin{table}[ht]
\centering
\begin{tabular}{|c||c|c|c|c|c|c|c c}\hline 
\backslashbox{$q$}{$h$} & $~0~$ & $-1$ & $-2$ & $-3$ & $-4$ & $-5$ &  \\ \hline\hline
$2$ & $\mathbb{Z}$ &  &  &  &  &  &  \\ \hline
$0$ & $\mathbb{Z}$ &  &  &  &  &  &  \\\hline
$-2$ &  &  & $\mathbb{Z}$ &  &  &  &  \\\hline
$-4$ &  &  & $\mathbb{Z}/2$ &  &  &  &  \\\hline
$-6$ &  &  &  & $\mathbb{Z}$ & $\mathbb{Z}$ &  &  \\\hline
$-8$ &  &  &  &  & $\mathbb{Z}/2$ &  &  \\\hline
$-10$ &  &  &  &  &  & $\mathbb{Z}$ &  \\\hline
 &  &  &  &  &  &  &  $\ddots$
\end{tabular}
\qquad
\begin{tabular}{|c||c|c|c|c|c|c|c c}\hline 
\backslashbox{$q$}{$h$} & $~0~$ & $-1$ & $-2$ & $-3$ & $-4$ & $-5$ &  \\ \hline\hline
$2$ & $\mathbb{Z}$ &  &  &  &  &  &  \\ \hline
$0$ & $\mathbb{Z}$ &  &  &  &  &  &  \\\hline
$-2$ &  &  & $\mathbb{Z}$ &  &  &  &  \\\hline
$-4$ &  &  & $\mathbb{Z}$ & $\mathbb{Z}$  &  &  &  \\\hline
$-6$ &  &  &  & $\mathbb{Z}$ & $\mathbb{Z}$ &  &  \\\hline
$-8$ &  &  &  &  & $\mathbb{Z}$ & $\mathbb{Z}$ &  \\\hline
$-10$ &  &  &  &  &  & $\mathbb{Z}$ &  \\\hline
 &  &  &  &  &  &  &  $\ddots$
\end{tabular}
\caption{$\mathcal{H}_e(U;2)$ and $\mathcal{H}_o(U;2)$, respectively, where $h$ is homological grading and $q$ is quantum grading}
\label{tab:evenodd}
\end{table}

Interestingly, the 2-colored unknot has no torsion. However, using computations provided by Sch\"utz (see Theorem 8.2 and Figure 9 of \cite{Sch_tz_2022}), the odd colored homology of the 3-colored unknot, $\mathcal{H}_o(U;3)$, contains $\mathbb{Z}/3\mathbb{Z}$-torsion, whereas $\mathcal{H}_e(U;3)$ contains no $\mathbb{Z}/3\mathbb{Z}$-torsion. Finally, note that the graded Euler characteristics of both sides agree.

\addtocontents{toc}{\SkipTocEntry}
\subsection{A flexible gluing property}
\label{ss:flexiblegluing}

The majority of this paper is devoted to developing a framework for the construction and calculation of unified projectors. This will entail setting up a tangle theory that is both compatible with unified Khovanov homology and will also allow for a very flexible notion of composition for tangles. Thankfully, the work of Naisse and Putyra \cite{naisse2020odd} (to which we will keep returning) accomplishes the former goal---thus, our goal is a generalization of their work which allows for this ``more flexible gluing property.''

To be clear, recall that Khovanov's theory for knots and links has been extended to tangles via at least two methods, by both Khovanov \cite{Khovanov_2002} and Bar-Natan \cite{Bar_Natan_2002,Bar_Natan_2005} (see \S \ref{section:TLn and projectors} for a review). In the former, Khovanov extended his work to tangles with an even number of endpoints, showing that the homotopy type of the complex he associates to each tangle is an invariant of the tangle. Furthermore, for each tangle $T$, the complex $\mathrm{Kh}(T)$ has an interpretation as a graded dg-bimodule over the so-called arc algebras, $H^n$. Paramount among the properties of these bimodules is the \textit{gluing} result, which states that, for stackable tangles $T$ and $S$, 
\[
\mathrm{Kh}(T) \otimes_{H^n} \mathrm{Kh}(S) \cong \mathrm{Kh}(TS).
\]

While Khovanov and Bar-Natan were able to describe an up-to-homotopy invariant complex associated to a tangle soon after the discovery of Khovanov homology, an analogue for odd Khovanov homology remained elusive for thirteen years after its discovery. Our work will employ the first known solution, provided by Naisse and Putyra in \cite{naisse2020odd}. Before detailing their solution, we remark that, in \cite{MR4190457}, Vaz constructed a \textit{super}category and derived from it a homological invariant of tangles which \textit{super}categorified the Jones polynomial. While he proved that his invariant was distinct from even Khovanov homology, it was not evident that his theory was isomorphic to odd Khovanov homology when restricted to links until the recent work of Schelstraete and Vaz in \cite{2023arXiv231114394S}. There, Schelstraete-Vaz provided another lift of odd Khovanov homology to tangles (indeed, their work succeeded in providing the first representation theoretic construction of odd Khovanov homology) which coincided with the ``not even Khovanov homology'' of \cite{MR4190457}. Naisse and Putyra conjectured that their tangle invariant is isomorphic to Vaz's, and thus to Schelstraete-Vaz's, but this remains an open question.

Naisse and Putyra's lift of odd Khovanov homology to tangles \cite{naisse2020odd} involves the introduction of objects called \textit{grading categories} which allow one to define categories of (dg-) bimodules graded by a selected grading category. The grading category for the problem at hand is a category $\mathcal{G}$ whose morphisms are given by a pair of a flat tangle (with even inputs and even outputs) and an element of $\mathbb{Z}\times\mathbb{Z}$. Viewing the unified arc algebra as a $\mathcal{G}$-graded algebra, $H^n$ becomes graded-associative (associativity fails before this change; see \S \ref{ss:unified arc algebras} and \ref{ss:brief graded outline}, and \cite{naisse2017odd} for more detail). In the context of grading categories, it is more difficult to define what is meant by a grading shift. In order to accomplish this, Naisse and Putyra implement \textit{shifting systems} which can be assigned to a grading category; in the case of $\mathcal{G}$, a shifting system is provided by a pair of a chronological cobordism and a shift in the $\mathbb{Z}\times \mathbb{Z}$-grading. See \S  \ref{ss:brief graded outline} for a more thorough introduction to grading categories and shifting systems.

For Naisse and Putyra, all this work meant that one could mimic the constructions of Khovanov in \cite{Khovanov_2002} in a graded-associative context, yielding a tangle version of unified Khovanov homology which respects the gluing property. Continuing the analogy, the goal of the majority of this paper is to provide a generalization of the gluing result of \cite{naisse2020odd} in the spirit of Bar-Natan's canopolies \cite{Bar_Natan_2005} or of Jones's planar algebras \cite{jones1999planaralgebrasi}. While the extension is minor and well known in the even setting (see a description in Section 4 of \cite{lawson2022homotopy}), realizing the analogous result in the odd setting, in this paper, means adapting the flat tangles of Naisse and Putyra to planar arc diagrams. In particular, the grading category $\mathcal{G}$ is upgraded to what we call a \textit{grading multicategory}, denoted $\mathscr{G}$. Then, the work of Naisse and Putyra provide us with a roadmap for proving what we refer to as ``multigluing,'' Theorem \ref{thm:multigluing}. The following is a statement of multigluing in lesser generality than we prove it. Recall that $\mathcal{F}$ is the unified chronological TQFT.

\begin{namedtheorem}[Theorem C]
Suppose $T$ is a diskular tangle of type $(m_1,\ldots, m_k; n)$ (see Definition \ref{def:planararcdiagrams}) and $T_i$ is a tangle diagram in a disk with $2m_i$ points on its boundary for each $i=1, \ldots, k$. Then there is an isomorphism
\[
\left(\mathcal{F}(T_1), \ldots, \mathcal{F}(T_k)\right)\otimes_{(H^{m_1}, \ldots, H^{m_k})} \mathcal{F}(T) \cong \mathcal{F}(T(T_1,\ldots, T_k)).
\]
\end{namedtheorem}

\begin{figure}[ht]
\centering
\includesvg[width=8cm]{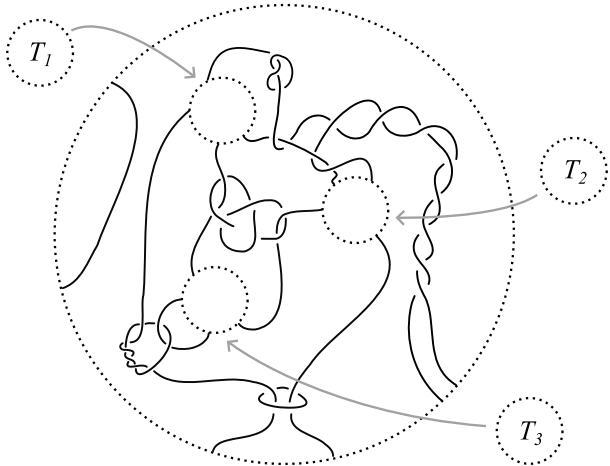}
\caption{Multigluing schematic. Here, we assume $T_1$, $T_2$, and $T_3$ are each tangles in disks with 4 points on their boundary.}
\label{fig:multigluingsvg}
\end{figure}

The notation $\otimes_{(H^{m_1}, \ldots, H^{m_k})}$, as well as the map inducing this isomorphism, is described in depth in \S \ref{s:gradingmultis and pads}. The idea of this theorem is that, given a tangle with some holes punched out, and compatible tangles $T_1, \ldots, T_k$, we can define a tensor product so that some tensor product of the dg-modules associated to $T_1, \ldots, T_k$ (denoted by $(T_1,\ldots, T_k)$) tensored with the multimodule associated to $T$ is isomorphic (as $\mathscr{G}$-graded dg-modules) to the dg-module associated to $T$ filled by the tangles $T_1,\ldots, T_k$. See Figure \ref{fig:multigluingsvg}.

\addtocontents{toc}{\SkipTocEntry}
\subsection{Other applications}

While our main motivation for this paper is a proof of existence for unified projectors and a new categorification of the colored Jones polynomial, there are other notable benefits of a more flexible gluing theorem; we will describe a few in our paper. To start, we can use Theorem \ref{thm:multigluing} to define operations on $\mathscr{G}$-graded dg-modules (e.g., a vertical stacking operation $\otimes$, juxtaposition $\sqcup$, and a partial trace $\mathrm{Tr}$) in exactly the same way as \cite{stoffregen2024joneswenzlprojectorskhovanovhomotopy}, see \S \ref{ss:operationsviamultigluing}. Defining these operations is essential as, without them, we cannot define categorified projectors. Of particular interest are our lifts of well-known adjunction statements provided by Hogancamp \cite{hogancamp2012morphismscategorifiedspinnetworks, https://doi.org/10.48550/arxiv.1405.2574}.

\begin{namedtheorem}[Theorem \ref{thm:adjunction1}]
If $A$ and $B$ are $\mathscr{G}$-graded dg-modules coming from tangle diagrams on $n-1$ and $n$ strands respectively, then
\[
\normalfont{\textsc{Hom}}_n\left(A \sqcup 1, \varphi_{\left(\tikz[baseline={([yshift=-.5ex]current bounding box.center)}, scale=.15]{
            \draw[fill=white, knot] (0,0) rectangle (4, 1.5);
            \draw[knot, red, rounded corners] (3.5, 2.125) -- (4.5, 2.125) -- (4.5, -0.625) -- (3.5, -0.625);
            \draw[knot] (0.5, 0) -- (0.5, -1.25);
            \draw[knot] (1.5, 0) -- (1.5, -1.25);
            \draw[knot] (2.5, 0) -- (2.5, -1.25);
            \draw[knot] (0.5, 1.5) -- (0.5, 2.75);
            \draw[knot] (1.5, 1.5) -- (1.5, 2.75);
            \draw[knot] (2.5, 1.5) -- (2.5, 2.75);
            \draw[knot] (3.5, 0) -- (3.5, -1.25);
            \draw[knot] (3.5, 1.5) -- (3.5, 2.75);
            \node at (2, 0.75) {\scriptsize \textit{B}};
        },\, (0,1)\right)}B\right) \cong 
    \normalfont{\textsc{Hom}}_{n-1}(A, \mathrm{Tr}(B) \{-1,0\}).
\]
\end{namedtheorem}

See the statement of Theorem \ref{thm:adjunction1} found in \S \ref{Chapter:OddCKProj} for more details. The notation $\normalfont{\textsc{Hom}}_n$ denotes the complex of maps of homogeneous bidegree; see \S \ref{ss:homomaps} and \ref{ss:dg-multimodules and hom}. Notice the $\mathscr{G}$-grading shift which is invisible to $q$-degree. We also obtain a more familiar statement, which we use in the proof of uniqueness for unified projectors:
\[
\normalfont{\textsc{Hom}}_n (A \otimes \mathcal{F}(\delta), B) \cong \normalfont{\textsc{Hom}}_n (A, B\otimes \mathcal{F}(\delta^\vee))
\]
where $\delta$ is a flat tangle. It comes as a corollary of another familiar ``duality'' statement; see Theorem \ref{thm:dualclosing}.

We remark that, in \S \ref{ss:egunifiedprojectors}, we construct $P_n$ as arising from an infinite torus braid, as in \cite{rozansky2010infinitetorusbraidyields} and Section 5 of \cite{stoffregen2024joneswenzlprojectorskhovanovhomotopy}. This description awards us with another, inductive description of $P_n$ as a filtered chain complex, inducing a spectral sequence. Explicitly, if $P_n$ and $P_{n-1}$ are projectors, we have that
\[
P_n = 
  \tikz[baseline={([yshift=-.5ex]current bounding box.center)}, scale=0.9]{
    \node at (0.85, 0.6) {$P_{n-1}$};
    \draw[knot] (0.2, 1.2) -- (0.2, 1.5);
    \draw[knot] (1.5, 1.2) -- (1.5, 1.5);
    \node at (0.85, 1.35) {$\cdots$};
    \draw[knot] (1.9, 0) to [out=-90, in=90] (0, -0.75);
    \draw[knot, overcross] (0.2, 0) -- (0.2,-0.75);
    \draw[knot, overcross] (0.5, 0) -- (0.5,-0.75);
    \draw[knot, overcross] (1.2, 0) -- (1.2,-0.75);
    \draw[knot, overcross] (1.5, 0) -- (1.5,-0.75);
    \draw[knot, overcross] (1.9, 1.5) -- (1.9, 0);
    \draw[knot] (0.2, -0.75) -- (0.2,-1.5);
    \draw[knot] (0.5, -0.75) -- (0.5,-1.5);
    \draw[knot] (1.2, -0.75) -- (1.2,-1.5);
    \draw[knot] (1.5, -0.75) -- (1.5,-1.5);
    \draw[knot, overcross] (0, -0.75) to [out=-90, in=90] (1.9, -1.5);
    \draw (0,0) rectangle (1.7, 1.2);
    \node at (0.85, -0.75) {$\cdots$};
    \node at (0.85, -1.45) {$\vdots$};
}
\]
where the wrap-around repeats indefinitely. It follows that $P_n$ is the colimit of a filtered chain complex of the following form.
\[
\begin{tikzcd}
\tikz[baseline={([yshift=-.5ex]current bounding box.center)}, scale=1]{
    \node at (0.85, 0.6) {$P_{n-1}$};
    \draw[knot] (0.2, 1.2) -- (0.2, 1.5);
    \draw[knot] (1.5, 1.2) -- (1.5, 1.5);
    \node at (0.85, 1.35) {$\cdots$};
    \draw[knot] (0.2, 0) -- (0.2,-0.75);
    \draw[knot] (0.5, 0) -- (0.5,-0.75);
    \draw[knot] (1.2, 0) -- (1.2,-0.75);
    \draw[knot] (1.5, 0) -- (1.5,-0.75);
    \draw[knot] (1.9, 1.5) -- (1.9, 0);
    \draw[knot] (1.9, 0) -- (1.9, -1.25);
    \draw[knot, overcross] (0.2, -0.75) -- (0.2,-1.25);
    \draw[knot, overcross] (0.5, -0.75) -- (0.5,-1.25);
    \draw[knot, overcross] (1.2, -0.75) -- (1.2,-1.25);
    \draw[knot, overcross] (1.5, -0.75) -- (1.5,-1.25);
    \draw (0,0) rectangle (1.7, 1.2);
    \node at (0.85, -0.75) {$\cdots$};
}
\arrow[r, hookrightarrow]
&
\tikz[baseline={([yshift=-.5ex]current bounding box.center)}, scale=1]{
    \node at (0.85, 0.6) {$P_{n-1}$};
    \draw[knot] (0.2, 1.2) -- (0.2, 1.5);
    \draw[knot] (1.5, 1.2) -- (1.5, 1.5);
    \node at (0.85, 1.35) {$\cdots$};
    \draw[knot] (1.9, 0) to [out=-90, in=90] (1.5, -0.75);
    \draw[knot] (0.2, 0) -- (0.2,-0.75);
    \draw[knot] (0.5, 0) -- (0.5,-0.75);
    \draw[knot] (1.2, 0) -- (1.2,-0.75);
    \draw[knot, overcross] (1.5,0) to[out=-90, in=90] (1.9, -0.75);
    \draw[knot] (1.9, -0.75) -- (1.9, -1.25);
    \draw[knot] (1.9, 1.5) -- (1.9, 0);
    \draw[knot, overcross] (0.2, -0.75) -- (0.2,-1.25);
    \draw[knot, overcross] (0.5, -0.75) -- (0.5,-1.25);
    \draw[knot, overcross] (1.2, -0.75) -- (1.2,-1.25);
    \draw[knot, overcross] (1.5, -0.75) -- (1.5,-1.25);
    \draw (0,0) rectangle (1.7, 1.2);
    \node at (0.85, -0.75) {$\cdots$};
}
\arrow[r, hookrightarrow]
&
\cdots
\arrow[r, hookrightarrow]
&
\tikz[baseline={([yshift=-.5ex]current bounding box.center)}, scale=1]{
    \node at (0.85, 0.6) {$P_{n-1}$};
    \draw[knot] (0.2, 1.2) -- (0.2, 1.5);
    \draw[knot] (1.5, 1.2) -- (1.5, 1.5);
    \node at (0.85, 1.35) {$\cdots$};
    \draw[knot] (1.9, 0) to [out=-90, in=90] (0.2, -0.75);
    \draw[knot, overcross] (0.2,0) to[out=-90, in=90] (0.5, -0.75);
    \draw[knot, overcross] (1.2,0) to[out=-90, in=90] (1.5, -0.75);
    \draw[knot, overcross] (1.5, 0) to[out=-90, in=90] (1.9, -0.75);
    \draw[knot, overcross] (1.9, -0.75) -- (1.9, -1.25);
    \draw[knot, overcross] (1.9, 1.5) -- (1.9, 0);
    \draw[knot, overcross] (0.2, -0.75) -- (0.2,-1.25);
    \draw[knot, overcross] (0.5, -0.75) -- (0.5,-1.25);
    \draw[knot, overcross] (1.5, -0.75) -- (1.5,-1.25);
    \draw (0,0) rectangle (1.7, 1.2);
    \node at (0.85, -0.75) {$\cdots$};
}
\arrow[r, hookrightarrow]
&
\tikz[baseline={([yshift=-.5ex]current bounding box.center)}, scale=1]{
    \node at (0.85, 0.6) {$P_{n-1}$};
    \draw[knot] (0.2, 1.2) -- (0.2, 1.5);
    \draw[knot] (1.5, 1.2) -- (1.5, 1.5);
    \node at (0.85, 1.35) {$\cdots$};
    \draw[knot] (1.9, 0) to [out=-90, in=90] (0.2, -0.75);
    \draw[knot, overcross] (0.2,0) to[out=-90, in=90] (0.5, -0.75);
    \draw[knot, overcross] (1.2,0) to[out=-90, in=90] (1.5, -0.75);
    \draw[knot, overcross] (1.5, 0) to[out=-90, in=90] (1.9, -0.75);
    \draw[knot, overcross] (1.9, -0.75) -- (1.9, -1.25);
    \draw[knot, overcross] (1.9, 1.5) -- (1.9, 0);
    \draw[knot] (0.5, -0.75) to[out=-90, in=90] (0.2, -1.25);
    \draw[knot, overcross] (0.2, -0.75) to[out=-90, in=90] (0.5, -1.25);
    \draw[knot, overcross] (1.5, -0.75) -- (1.5,-1.25);
    \draw (0,0) rectangle (1.7, 1.2);
    \node at (0.85, -0.75) {$\cdots$};
}
\arrow[r, hookrightarrow]
&
\cdots
\end{tikzcd}
\]
The filtration on $P_n$ induces one on its full trace, and using results of \S \ref{s:unifiedprops}, we conclude that $\normalfont{\textsc{Hom}}_n(P_n, P_n)$ is a filtered complex. We will investigate the associated graded of this filtration in future work.

\addtocontents{toc}{\SkipTocEntry}
\subsection{Future goals}

Further motivation for our work was provided by some questions left unanswered in this paper. We conclude the introduction by outlining a few of them.

\subsubsection*{Periodicity of projectors and a GOR conjecture}

We note (see Corollary \ref{cor:endothing}) that the existence of unified projectors (Theorem \ref{thm:existenceofprojectors}), together with an adjunction statement (Theorem \ref{thm:adjunction1}), implies that
\begin{equation}
\label{eq:homologyofendoalgebra}
H^* \normalfont{\textsc{Hom}}_n (P_n, P_n) \cong q^{-n} \mathcal{H}(U; n).
\end{equation}
In \cite{https://doi.org/10.48550/arxiv.1405.2574}, Hogancamp uses the specialization of Equation (\ref{eq:homologyofendoalgebra}) to the even setting in order to construct particular elements $U_n \in \normalfont{\textsc{Hom}}_n (P_n^{\mathrm{CK}}, P_n^{\mathrm{CK}})$ to make substantial progress toward a conjecture of Gorsky-Oblomkov-Rasmussen \cite{gorsky2013stablekhovanovhomologytorus, Gorsky_2014}. The chain maps $U_n$ take the form $t^{2-2n}q^{2n} P_n^{\mathrm{CK}} \to P_n^{\mathrm{CK}}$ and satisfy $\mathrm{Cone}(U_n) \simeq Q_n$ (for a particular complex $Q_n$), showing that $P_n^{\mathrm{CK}}$ is a periodic chain complex built from copies of $Q_n$. Interestingly, in \cite{Sch_tz_2022},  Sch\"utz computes the first few odd projectors $P_2^o$ and $P_3^o$ algorithmically and shows that, while odd $P_2^o$ is also periodic of period 2, odd $P_3^o$ is periodic of period 8, unlike even $P_3^{\mathrm{CK}}$ which has period 4 (cf. Section 4.4 of \cite{https://doi.org/10.48550/arxiv.1005.5117}). We hope to use results of this paper to prove that $P_n$ remains periodic in the unified and odd settings, and to determine the period of $P_n$ for arbitrary $n$.

\subsubsection*{Odd Khovanov spectra for tangles}

The idea for generalizing the work of Naisse-Putyra via dg-multimodules associated to diskular tangles came largely from observations of the utility of spectral multimodules in the work of Lawson-Lipshitz-Sarkar \cite{Lawson_2021, lawson2022homotopy} and Stoffregen-Willis \cite{stoffregen2024joneswenzlprojectorskhovanovhomotopy}. Now, an odd (indeed, unified) Khovanov homotopy type is known \cite{SARKAR2020107112}, but it has yet to be lifted to the setting of tangles---we hope that our work might be melded with that of \cite{Lawson_2021} and \cite{SARKAR2020107112} to produce a unified homotopy type for tangles. If this is accomplished, it is also our hope that the work here will allow for the arguments of \cite{lawson2022homotopy} to lift, proving that homotopy functoriality holds in higher generality. It is also interesting to note that the spectral projector on three strands of \cite{stoffregen2024joneswenzlprojectorskhovanovhomotopy} is periodic of period 8, like the \textit{odd} projector on three strands of \cite{Sch_tz_2022}, but unlike the three-stranded \textit{even} projector.

\subsubsection*{Functoriality}

In forthcoming work, we note that there is a natural definition of a $\mathscr{C}$-graded bar resolution and Hochschild homology for $\mathscr{C}$-graded algebras. Further, we hope that one could mimic Khovanov's proof of functoriality \cite{Khovanov_2002} (see \cite{lawson2022homotopy} for an excellent outline) in the unified setting to obtain ``up-to-unit'' functoriality for unified Khovanov homology, and thus up-to-sign functoriality for odd Khovanov homology. Perhaps more interesting (especially if aiming for a Lasagna-type invariant coming from a functorial invariant of links in $S^3$ \cite{Morrison2019InvariantsO4, morrison2024invariantssurfacessmooth4manifolds}) is the development of a functorial ``oriented model'' \cite{doi:10.1142/S0218216510007863} for odd Khovanov homology. Such a model is provided in \cite{2023arXiv231114394S}, but the question of functoriality remains open.

\addtocontents{toc}{\SkipTocEntry}
\subsection*{Outline}

Sections \ref{section:TLn and projectors} and \ref{S:odd chronologies and stuff} of this paper are preparatory and intended as introductions for the uninitiated; experts should feel free to skip them. In Section \ref{section:TLn and projectors}, we start by recalling the definition of the Temperley-Lieb algebra and its special elements, the Jones-Wenzl projectors. Then, we recall basic features of the categorifications of $TL_n$ by both of Bar-Natan \cite{Bar_Natan_2002, Bar_Natan_2005} and Khovanov \cite{Khovanov_2002}. Finally we review some facts about the first categorification of Jones-Wenzl projectors, following \cite{https://doi.org/10.48550/arxiv.1005.5117}. Section \ref{S:odd chronologies and stuff} is devoted to reviewing some of the work of Putyra regarding categories of chronological cobordisms \cite{putyra2010cobordismschronologiesgeneralisationkhovanov, putyra20152categorychronologicalcobordismsodd}. We end Section \ref{S:odd chronologies and stuff} by presenting an outline of $\mathcal{C}$-graded structures, for a grading category $\mathcal{C}$, as in \cite{naisse2020odd}---the hope is that \S \ref{ss:brief graded outline} might give the reader a bird's-eye view of the goals of Sections \ref{s:gradingmultis and pads}, \ref{S: SHIFTING SYSTEMS}, and \ref{s:tangles,multimodules,multigluing}.

Sections \ref{s:gradingmultis and pads}, \ref{S: SHIFTING SYSTEMS}, and \ref{s:tangles,multimodules,multigluing} are the technical heart of this paper, wherein we introduce grading multicategories, shifting 2-systems for those grading multicategories, and apply the general framework constructed to prove multigluing, Theorem \ref{thm:multigluing}. Again, see \S \ref{ss:brief graded outline} for a more complete outline.

In Section \ref{Chapter:tangleinvariant}, we use multigluing to obtain an invariant of (diskular) tangles, slightly generalizing a result of \cite{naisse2020odd}. As in the cited paper, the grading system is, perhaps, too sensitive for the complex associated to a (diskular) tangle diagram to be invariant under each of the Reidemeister moves (see Lemmas \ref{lem:r1}, \ref{lem:r2}, and \ref{lem:r3}). However, it is invariant up to a grading shift in which the number of saddles in the cobordism component is equal to the sum of the entries of the $\mathbb{Z} \times \mathbb{Z}$ component. Hence, we can ``collapse'' the $\mathscr{G}$-degree to an integral $q$-grading in to obtain a tangle invariant. We remark that, however slight the generalization, the added flexibility is necessary for our final result in \S \ref{ss:unifiedcoloredhomology} (additionally, we believe the differences in our proof to be notable).

Finally, in Section \ref{Chapter:OddCKProj}, we define and prove the existence and uniqueness of categorifications of the Jones-Wenzl projectors living in a category of $\mathscr{G}$-graded dg-modules, specializing to the projectors of \cite{https://doi.org/10.48550/arxiv.1005.5117}, but also to ``odd'' projectors which, prior to this paper, had only been computed up to three or so strands (cf. \cite{Sch_tz_2022}). Other highlights of this section are the proofs of the aforementioned duality and adjunction results, which we hope to be useful in future work. In conclusion, we point out that the existence of unified projectors, together with multigluing and the tangle invariant of Section \ref{Chapter:tangleinvariant}, imply the existence of a unified colored link homology, specializing to the colored link homology of, say, \cite{https://doi.org/10.48550/arxiv.1005.5117}, but also to a new, ``odd'' categorification of the colored Jones polynomial.

\addtocontents{toc}{\SkipTocEntry}
\subsection*{Acknowledgements}

I am indebted to my advisors, Efstratia Kalfagianni and Matthew Stoffregen, for patiently aiding me throughout the completion of this project, which was recommended to me by the latter. It is also a pleasure to thank Teena Gerhardt, Matthew Hedden, Adam Lowrance, Gr\'egoire Naisse, Krzysztof Putyra, David Rose, Pedro Vaz, and Michael Willis for helpful conversations, suggestions, and additional guidance. This work was supported by the NSF grants DMS-2004155 and DMS-2304033 and the NSF/RTG grant DMS-2135960.

\newpage

\section{Categorifications of \texorpdfstring{$TL_n$}{Lg} and projectors}
\label{section:TLn and projectors}

In this section, we survey attributes of the even setting which we hope to lift---in one way or another---to the odd setting. In \S \ref{ss:decategorified}, we briefly discuss the decategorified setting. In \S \ref{evenCat}, we recall the even categorifications of the Temperley-Lieb algebras due to Bar-Natan \cite{Bar_Natan_2005} and Khovanov \cite{Khovanov_2002}. We conclude by providing Cooper and Krushkal's categorification of the Jones-Wenzl projectors in \S \ref{ss:ckprojectors}, as we hope to compare their results with our work in \S \ref{Chapter:OddCKProj}.

\subsection{Temperley-Lieb algebras and Jones-Wenzl projectors}
\label{ss:decategorified}

The Temperley-Lieb algebras $TL_n$ arise naturally as the $U_q(\mathfrak{sl}_2)$-equivariant endomorphisms of $n$-fold tensor powers of the fundamental representation of $U_q(\mathfrak{sl}_2)$. As a unital $\mathbb{Z}[q,q^{-1}]$-algebra, $TL_n$ is generated by $n$ elements $1_n, e_1, \ldots, e_{n-1}$ subject to the relations
\begin{enumerate}
    \item $e_ie_j = e_je_i$ if $|i-j|\ge 2$,
    \item $e_ie_{i\pm 1}e_i = e_i$, and
    \item $e_i^2 = (q+q^{-1}) e_i$.
\end{enumerate}
The first relation is referred to as ``distant commutativity.'' We will make use of the quantum integer notation
\[
[k] = \cfrac{q^k - q^{-k}}{q-q^{-1}}
\]
so that, for example, the third relation can be rewritten $e_i^2 = [2] e_i$.

$TL_n$ can be given a diagramatic description, where the generating elements are presented by
\[
1_n = \tikz[baseline = 2ex, scale = .8]{
\draw (0,0) -- (0,1);
\draw (1,0) -- (1,1);
\node at (.5,.5) {$\cdots$};
\draw[dotted] (-.2,1) rectangle (1.2,0);
} \qquad \text{and} \qquad \,
e_i = \tikz[baseline = 2ex, scale = .8]{
\draw (0,0) -- (0,1);
\node at (.5,.5) {$\cdots$};
\draw[rounded corners = 2mm] (1,0) -- (1,.25) -- (1.5,.49) -- (2,.25) -- (2,0);
\draw[rounded corners = 2mm] (1,1) -- (1,.75) -- (1.5,.51) -- (2,.75) -- (2,1);
\node at (2.5,.5) {$\cdots$};
\draw (3,0) -- (3,1);
\draw[dotted] (-.2,1) rectangle (3.2,0);
\node[below = 1pt] at (1,0) {$i$};
\node[below = 1pt] at (2,0) {$i+1$};
}
\]
with multiplication given by top-to-bottom vertical stacking. Therefore, $TL_n$ can be viewed as the linear skein of the disk with $2n$ distinguished points on its boundary, where we regard this disk as a square with $n$ marked points on the top and $n$ marked points on the bottom. It is in this way that every $(n,n)$-tangle may be assigned an element of $TL_n$; indeed, given an oriented tangle, the relations
\[
\tikz[baseline=1.6ex, scale = .8]{
\draw[-stealth] (0,0) -- (1,1);
\draw (1,0) -- (.7,.3);
\draw[-stealth] (.3,.7) -- (0,1);
}
 = q ~ \tikz[baseline=1.6ex, scale = .8]{
\draw[rounded corners = 4mm] (0,0) -- (.45,.5) -- (0,1);
\draw[rounded corners = 4mm] (1,0) -- (.55,.5) -- (1,1);
}~ -q^2 ~ \tikz[baseline=1.6ex, scale = .8]{
\draw[rounded corners = 2mm] (0,0) -- (0,.25) -- (.5,.49) -- (1,.25) -- (1,0);
\draw[rounded corners = 2mm] (0,1) -- (0,.75) -- (.5,.51) -- (1,.75) -- (1,1);
}~ \qquad \text{and} \qquad \,
~\tikz[baseline=1.6ex, scale = .8]{
\draw[-stealth] (1,0) -- (0,1);
\draw (0,0) -- (.3, .3);
\draw[-stealth] (.7, .7) -- (1,1)
}
 = q^{-2} ~\tikz[baseline=1.6ex, scale = .8]{
\draw[rounded corners = 2mm] (0,0) -- (0,.25) -- (.5,.49) -- (1,.25) -- (1,0);
\draw[rounded corners = 2mm] (0,1) -- (0,.75) -- (.5,.51) -- (1,.75) -- (1,1);
}~ -q^{-1} ~\tikz[baseline=1.6ex, scale = .8]{
\draw[rounded corners = 4mm] (0,0) -- (.45,.5) -- (0,1);
\draw[rounded corners = 4mm] (1,0) -- (.55,.5) -- (1,1);
}
\]
yield the Jones polynomial up to normalization.

In ~\cite{Lickorish1993THESM}, it was shown that the Witten-Reshetikhin-Turaev 3-manifold invariants (\cite{cmp/1104178138, Reshetikhin1991}) may be constructed combinatorially via the Kauffman bracket. Key ingredients of this construction are the Jones-Wenzl projectors, which we recall now.

\begin{definition}
The \textit{Jones-Wenzl projectors}, denoted by $p_n$, are particular elements of $TL_n$, defined by the recursion
\[
p_1 = 1_1 ~\qquad \text{and}~ \qquad p_{n+1} = (p_n \sqcup 1) - \cfrac{[n]}{[n+1]} ~(p_n \sqcup 1) e_{n-1} (p_n \sqcup 1).
\]
It is common to depict $p_n$ by a box
\[
p_n = \tikz[baseline=.8ex, scale = .4]{
\draw (0,0) rectangle (1,1);
\draw (.5,-.5) -- (.5,0);
\draw (.5,1) -- (.5, 1.5);
\node at (.5,.5) {$n$};
}
\]
in which case the recursion appears as 
\[
\tikz[baseline={([yshift=-.5ex]current bounding box.center)}, scale=1]{
    \draw[knot] (0,0) -- (0,2);
    \draw[fill=white] (-0.3, 0.7) rectangle (0.3, 1.3);
    \node at (0,1) {$1$};
}
\quad=\quad
\tikz[baseline={([yshift=-.5ex]current bounding box.center)}, scale=1]{
    \draw[knot] (0,0) -- (0,2);
}
\qquad \text{and} \qquad
\tikz[baseline={([yshift=-.5ex]current bounding box.center)}, scale=1]{
    \draw[knot] (0,0) -- (0,2);
    \draw[fill=white] (-0.45, 0.7) rectangle (0.45, 1.3);
    \node at (0,1) {$n+1$};
}
\quad = \quad
\tikz[baseline={([yshift=-.5ex]current bounding box.center)}, scale=1]{
    \draw[knot] (0,0) -- (0,2);
    \draw[fill=white] (-0.3, 0.7) rectangle (0.3, 1.3);
    \node at (0,1) {$n$};
    \draw[knot] (0.75,0) -- (0.75,2);
}
\quad -\quad \cfrac{[n]}{[n+1]} ~
\tikz[baseline={([yshift=-.5ex]current bounding box.center)}, scale=1]{
    \draw[knot] (0.5,0) -- (0.5,2);
    \draw[knot] (1,0) -- (1, 0.65);
    \draw[knot] (1,0.65) to[out=90, in=90] (2,0.65);
    \draw[knot] (2, 0.65) -- (2,0);
    \draw[knot] (1,2) -- (1, 1.35);
    \draw[knot] (1, 1.35) to[out=-90, in=-90] (2,1.35);
    \draw[knot] (2, 1.35) -- (2,2);
    \draw[fill=white] (0.3, 0.1) rectangle (1.2, 0.7);
    \node at (0.75, 0.4) {$n$};
    \draw[fill=white] (0.3, 1.9) rectangle (1.2, 1.3);
    \node at (0.75, 1.6) {$n$};
}~.
\]
\end{definition}

The Jones-Wenzl projectors are well-studied. They may be defined equivalently as the unique elements of $TL_n$ for which
\begin{enumerate}
    \item[(JW1)] $(p_n - 1_n)$ belongs to the algebra generated by $\{e_1,\ldots, e_{n-1}\}$, and 
    \item[(JW2)] $p_n e_i = e_i p_n = 0$ for all $i=1,\ldots, n-1$.
\end{enumerate}
These properties immediately imply that the projectors are idempotents. One can also check that upon taking the Markov closure of the projectors, the Kauffman bracket evaluates them as a quantum integer:
\[
\langle \widehat{p_n} \rangle = [n+1].
\]

The purpose of listing these well-known properties of the Jones-Wenzl projectors is that their categorifications satisfy analogues in the categorified setting. We will use these properties frequently in what follows.

\subsection{Categorifications of the Temperley-Lieb algebra}
\label{evenCat}

We start by reviewing a construction of Bar-Natan ~\cite{Bar_Natan_2002, Bar_Natan_2005} which categorifies $TL_n$. Consequently, we may determine the Khovanov complex for a tangle, which turns out to be a tangle invariant up to homotopy. Afterwards, we describe another categorification of Khovanov, which has a known analogue in the odd setting. In the broader context of this paper, we wish to review Bar-Natan's categorification to motivate our grading multicategory $\mathscr{G}$, defined in Section 4.

Recall that a \textit{pre-additive category} $\mathcal{C}$ is a category such that
\begin{enumerate}
    \item for every $X,Y\in \text{ob}(\mathcal{C})$, $\text{Hom}_\mathcal{C}(X,Y)$ is an abelian group, and
    \item morphism composition distributes over the abelian group's addition rule.
\end{enumerate}
Additionally, a \textit{monoidal category} $\mathcal{C}$ is a category endowed with a functor $\otimes: \mathcal{C}\times\mathcal{C} \to \mathcal{C}$, a distinguished object $1\in\mathrm{ob}(\mathcal{C})$, and natural isomorphisms $\alpha$ (called the \textit{associator}) and left- and right-unitors $\lambda$ and $\rho$ satisfying the triangle and pentagon identities.

Given a pre-additive category $\mathcal{C}$, we may define the \textit{(split) Grothendieck group} of $\mathcal{C}$ to be the free abelian group generated by isomorphism classes in $\mathcal{C}$, with the added relation that $[A\oplus B] = [A] + [B]$:
\[
K_0(\mathcal{C}) = \mathbb{Z} \langle \mathcal{C} \rangle \left/\begin{Bmatrix}
[A]=[B]~\text{if}~A\cong B \\
[A\oplus B] = [A] + [B]
\end{Bmatrix}\right..
\]
It is common to take the Grothendieck group of pre-addivive monoidal categories---in this case, the tensor product induces an algebra structure on $K_0(\mathcal{C})$. 

For us, to categorify $TL_n$ means to define a pre-additive monoidal category $\mathcal{C}$ for which $K_0(\mathcal{C}) \cong TL_n$. Here is an outline of the construction provided by Bar-Natan.

\medskip
\noindent \textit{Step 1}: Let $\text{pre-Cob}(n)$ denote the cateory whose
\begin{itemize}
    \item objects are isotopy classes of formally $q$-graded Temperley-Lieb diagrams with $2n$ boundary points, and
    \item $\text{Hom}(q^iA, q^jB)$ is the free $\mathbb{Z}$-module spanned by isotopy classes of orientable cobordisms from $A$ to $B$.
\end{itemize}
Note that $\text{pre-Cob}(n)$ is pre-additive by definition. All of our cobordisms will be oriented upwards (from bottom to top). It is also naturally monoidal via stacking in $TL_n$. It is clear that if $C: A\to B$ and $C': A' \to B'$, then there is a cobordism $C\otimes C': A\otimes A' \to B\otimes B'$.

\begin{definition}
The \textit{degree} of a cobordism $C: q^iA \to q^j B$ is the value
\[
\deg(C) = \deg_t(C) + \deg_q(C)
\]
where
\begin{enumerate}[label=(\roman*)]
    \item $\deg_t(C) = \chi(C) - n$ is called the \textit{topological degree} of $C$, and
    \item $\deg_q(C) = j-i$ is called the \textit{quantum degree} of $C$.
\end{enumerate}
\end{definition}
\noindent It is common practice to fix $q$-gradings on the Temperley-Lieb elements so that $\deg(C)$ is always zero.

There are a few special cobordisms which we highlight here. Their frequent use necessitates additional (but commonplace) notation.
\begin{enumerate}[label = (\arabic*)]
    \item Cobordisms in this category may be decorated by dots, which correspond to hollow handle attachments up to multiplication by 2.
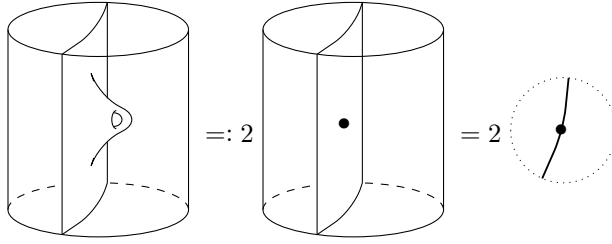
\[
    \tikz[baseline=6ex, scale = .6]{
    \draw (-2,0) arc (180:360:2 and 0.6);
    \draw[dashed] (2,0) arc (0:85:2 and 0.6);
    \draw[dashed] (-2,0) arc (180:114:2 and 0.6);
    \draw (-2,4) arc (180:360:2 and 0.6);
    \draw (2,4) arc (0:180:2 and 0.6);
    \draw (-2,0) -- (-2,4);
    \draw (2,0) -- (2,4);
    \draw[rounded corners = 3mm] (0.201251467,4+0.596954609867) -- (0,4) --  (-0.80159834,4-0.54969956);
    \draw[rounded corners = 3mm] (0.20125146,0.596954609867) -- (0,0) --  (-0.80159834,-0.54969956);
    \draw (-0.80159834,-0.54969956) -- (-0.80159834,4-0.54969956);
    \draw[] (0.20125146,0.596954609867) -- (0.20125146,1.5);
    \draw[] (0.20125146,4+0.596954609867) -- (0.20125146,2.47);
    \draw[rounded corners = 3mm] (-.1, 2.85) -- (0,2.55) -- (1,2) -- (0,1.45) -- (-.1, 1.15);
    \draw (.4,2.2) arc (90:270:.1 and .2);
    \draw (.35,1.85) arc (-90:90:.175 and .15);
} ~ =: 2~
\tikz[baseline=6ex, scale = .6]{
    \draw (-2,0) arc (180:360:2 and 0.6);
    \draw[dashed] (2,0) arc (0:85:2 and 0.6);
    \draw[dashed] (-2,0) arc (180:114:2 and 0.6);
    \draw (-2,4) arc (180:360:2 and 0.6);
    \draw (2,4) arc (0:180:2 and 0.6);
    \draw (-2,0) -- (-2,4);
    \draw (2,0) -- (2,4);
    \draw[rounded corners = 3mm] (0.201251467,4+0.596954609867) -- (0,4) --  (-0.80159834,4-0.54969956);
    \draw[rounded corners = 3mm] (0.20125146,0.596954609867) -- (0,0) --  (-0.80159834,-0.54969956);
    \draw (-0.80159834,-0.54969956) -- (-0.80159834,4-0.54969956);
    \draw[] (0.20125146,0.596954609867) -- (0.20125146,4+0.596954609867);
    \node at (-.2,1.9) {$\bullet$};
} ~ = 2~
\tikz[baseline=-1ex, scale = .7]{
    \draw[dotted] (-1,0) arc (180:360:1 and 1);
    \draw[dotted] (1,0) arc (0:180:1 and 1);
    \draw[knot, rounded corners = 3mm] (0.1006257333876,0.994924349778) -- (0,0) --  (-0.40079917208,-0.916165936749);
    \node at (-0.04,0) {$\bullet$};
}
\]
    \item Saddles will have the following shorthand.
\[
\tikz[baseline=6ex, scale = .6]{
    \draw (-2,0) arc (180:360:2 and 0.6);
    \draw[dashed] (2,0) arc (0:80:2 and 0.6);
    \draw[dashed] (-2,0) arc (180:100:2 and 0.6);
    \draw (-2,4) arc (180:360:2 and 0.6);
    \draw (2,4) arc (0:180:2 and 0.6);
    \draw (-2,0) -- (-2,4);
    \draw (2,0) -- (2,4);
    \draw[rounded corners = 3mm] (-0.395777629218,4.588134683673) -- (-.4,4) -- (-0.80159834416,4-0.54969956205);
    \draw[rounded corners = 3mm] (0.398899441995,4.587944836839) -- (.4,4) -- (0.802005173956,4-0.549646152616);
    \draw[rounded corners = 3mm] (-0.395777629218,0.588134683673) -- (0,.4) -- (0.398899441995,0.587944836839);
    \draw[rounded corners = 3mm] (-0.80159834416,-0.54969956205) -- (0,-.2) -- (0.802005173956,-0.549646152616);
    \draw (-0.395777629218,4.588134683673) -- (-0.395777629218,0.588134683673);
    \draw (-0.80159834416,4-0.54969956205) -- (-0.80159834416,-0.54969956205);
    \draw (0.398899441995,4.587944836839) -- (0.398899441995,0.587944836839); 
    \draw (0.802005173956,4-0.549646152616) -- (0.802005173956,-0.549646152616);
    \draw (-.4,4) arc (180:360:.4 and 2.5);
    \draw (0,1.5) arc (90:180:.2 and 1.8);
    \draw[dashed] (0,1.5) arc (90:0:.175 and 1);
} ~=~ 
\tikz[baseline=-1ex, scale = .7]{
    \draw[dotted] (-1,0) arc (180:360:1 and 1);
    \draw[dotted] (1,0) arc (0:180:1 and 1);
    \draw[red,thick] (0,-.43) -- (0,.43);
    \draw[knot, rounded corners=3mm] (-0.667,-0.745) -- (0,-0.3) -- (0.667,-0.745);
    \draw[knot, rounded corners=3mm] (-0.667,0.745) -- (0,0.3) -- (0.667,0.745);
}
\]
\end{enumerate}
Note that, for example
\[
\deg_t \left(\tikz[baseline=-1ex, scale = .7]{
    \draw[dotted] (-1,0) arc (180:360:1 and 1);
    \draw[dotted] (1,0) arc (0:180:1 and 1);
    \draw[knot, rounded corners = 3mm] (0.1006257333876,0.994924349778) -- (0,0) --  (-0.40079917208,-0.916165936749);
    \node at (-0.04,0) {$\bullet$};
} \right) = \chi \left( S^1\vee S^1\right) - 1 = -2
\]
since the dotted identity has the same homotopy type as the punctured torus, and
\[
\deg_t\left(\tikz[baseline=-1ex, scale = .7]{
    \draw[dotted] (-1,0) arc (180:360:1 and 1);
    \draw[dotted] (1,0) arc (0:180:1 and 1);
    \draw[red,thick] (0,-.43) -- (0,.43);
    \draw[knot, rounded corners=3mm] (-0.667,-0.745) -- (0,-0.3) -- (0.667,-0.745);
    \draw[knot, rounded corners=3mm] (-0.667,0.745) -- (0,0.3) -- (0.667,0.745);
}\right) = \chi \left( \mathbb{D}^2\right) - 2 = -1.
\]
Therefore, we will take dots to increase quantum degree 2 and saddles to increase quantum degree 1.

\medskip
\noindent\textit{Step 2}: Pass to the matrix category $\text{Mat}(\text{pre-Cob}(n))$, whose objects are vectors of objects in $\text{pre-Cob}(n)$ and whose morphisms are matrices of morphisms in $\text{pre-Cob}(n)$. Observing the defining relations in $TL_n$, to construct a category $\mathcal{C}$ for which $K_0(\mathcal{C}) \cong TL_n$, the object represented by $\bigcirc$ in $\mathcal{C}$  must be isomorphic to the sum of two empty objects in degree $\pm 1$:
\[
\bigcirc \cong q^{-1}~\emptyset ~\oplus~ q~\emptyset.
\]
We accomplish this by defining \textit{delooping} operations. Consider the morphisms
\[
\varphi: \bigcirc \xrightarrow{\begin{pmatrix} \tikz[baseline=.4ex, scale = .2]{
    \draw [domain=0:180] plot ({2*cos(\x)}, {2*sin(\x)});
    \draw (-2,0) arc (180:360:2 and 0.6);
    \draw[dashed] (2,0) arc (0:180:2 and 0.6);
} \\ \tikz[baseline=.4ex, scale = .2]{
    \draw [domain=0:180] plot ({2*cos(\x)}, {2*sin(\x)});
    \draw (-2,0) arc (180:360:2 and 0.6);
    \draw[dashed] (2,0) arc (0:180:2 and 0.6);
    \node at (0,1.25) {$\bullet$};
} \end{pmatrix} } q^{-1}~\emptyset ~\oplus~ q~\emptyset
\]
and
\[
\psi: q^{-1}~\emptyset ~\oplus~ q~\emptyset \xrightarrow{\begin{pmatrix} \tikz[baseline=.4ex, scale = .2]{
    \draw [domain=180:360] plot ({2*cos(\x)}, {2*sin(\x)});
    \draw (-2,0) arc (180:360:2 and 0.6);
    \draw (2,0) arc (0:180:2 and 0.6);
    \node at (0,-1.25) {$\bullet$};
} & \tikz[baseline=.4ex, scale = .2]{
    \draw [domain=180:360] plot ({2*cos(\x)}, {2*sin(\x)});
    \draw (-2,0) arc (180:360:2 and 0.6);
    \draw (2,0) arc (0:180:2 and 0.6);
} \end{pmatrix}} \bigcirc.
\]
We impose the isomorphism above by defining the relations implied by $\varphi \circ \psi = \text{id}_{\mathbb{Z}\otimes \mathbb{Z}}$ and $\psi \circ \varphi = \text{id}_{\bigcirc}$. On one hand,
\[
\varphi \circ \psi = \begin{pmatrix}
\tikz[baseline=.4ex, scale = .2]{
  \draw (0,0) circle (2cm);
  \draw (-2,0) arc (180:360:2 and 0.6);
  \draw[dashed] (2,0) arc (0:180:2 and 0.6);
  \node at (0,-1.25) {$\bullet$};
}& 
\tikz[baseline=.4ex, scale = .2]{
  \draw (0,0) circle (2cm);
  \draw (-2,0) arc (180:360:2 and 0.6);
  \draw[dashed] (2,0) arc (0:180:2 and 0.6);
}
\\ 
\tikz[baseline=.4ex, scale = .2]{
  \draw (0,0) circle (2cm);
  \draw (-2,0) arc (180:360:2 and 0.6);
  \draw[dashed] (2,0) arc (0:180:2 and 0.6);
  \node at (0,-1.25) {$\bullet$};
  \node at (0,1.25) {$\bullet$};
}&
\tikz[baseline=.4ex, scale = .2]{
  \draw (0,0) circle (2cm);
  \draw (-2,0) arc (180:360:2 and 0.6);
  \draw[dashed] (2,0) arc (0:180:2 and 0.6);
  \node at (0,1.25) {$\bullet$};
}
\end{pmatrix}.
\]
On the other,
\[
\psi \circ \varphi = \tikz[baseline=-4ex, scale=.2]{
    \draw [domain=180:360] plot ({2*cos(\x)}, {2*sin(\x)});
    \draw (-2,0) arc (180:360:2 and 0.6);
    \draw (2,0) arc (0:180:2 and 0.6);
    \node at (0,-1.25) {$\bullet$};
    \begin{scope}[shift={(0,-5)}]
        \draw[domain=0:180] plot ({2*cos(\x)}, {2*sin(\x)});
        \draw (-2,0) arc (180:360:2 and 0.6);
        \draw[dashed] (2,0) arc (0:180:2 and 0.6);
    \end{scope}
}
+
\tikz[baseline=-4ex, scale=.2]{
    \draw [domain=180:360] plot ({2*cos(\x)}, {2*sin(\x)});
    \draw (-2,0) arc (180:360:2 and 0.6);
    \draw (2,0) arc (0:180:2 and 0.6);
    \begin{scope}[shift={(0,-5)}]
        \draw[domain=0:180] plot ({2*cos(\x)}, {2*sin(\x)});
        \draw (-2,0) arc (180:360:2 and 0.6);
        \draw[dashed] (2,0) arc (0:180:2 and 0.6);
        \node at (0,1.25) {$\bullet$};
    \end{scope}
}.
\]
In conclusion, we define $\text{Cob}(n)$ to be the quotient of $\text{pre-Cob}(n)$ by the relations
\[
\tikz[baseline=-.5ex, scale = .2]{
  \draw (0,0) circle (2cm);
  \draw (-2,0) arc (180:360:2 and 0.6);
  \draw[dashed] (2,0) arc (0:180:2 and 0.6);
} = 0, \qquad 
\tikz[baseline=-.5ex, scale = .2]{
  \draw (0,0) circle (2cm);
  \draw (-2,0) arc (180:360:2 and 0.6);
  \draw[dashed] (2,0) arc (0:180:2 and 0.6);
  \node at (0,0) {$\bullet$};
} = 1, \qquad
\tikz[baseline=-.5ex, scale = .2]{
  \draw (0,0) circle (2cm);
  \draw (-2,0) arc (180:360:2 and 0.6);
  \draw[dashed] (2,0) arc (0:180:2 and 0.6);
  \node at (-.5,0) {$\bullet$};
  \node at (.5,0) {$\bullet$};
} = 0,~\text{and}
\]
\[\tikz[baseline=-4ex, scale=.2]{
    \draw [domain=180:360] plot ({2*cos(\x)}, {2*sin(\x)});
    \draw (-2,0) arc (180:360:2 and 0.6);
    \draw (2,0) arc (0:180:2 and 0.6);
    \node at (0,-1.25) {$\bullet$};
    \begin{scope}[shift={(0,-5)}]
        \draw[domain=0:180] plot ({2*cos(\x)}, {2*sin(\x)});
        \draw (-2,0) arc (180:360:2 and 0.6);
        \draw[dashed] (2,0) arc (0:180:2 and 0.6);
    \end{scope}
}
+
\tikz[baseline=-4ex, scale=.2]{
    \draw [domain=180:360] plot ({2*cos(\x)}, {2*sin(\x)});
    \draw (-2,0) arc (180:360:2 and 0.6);
    \draw (2,0) arc (0:180:2 and 0.6);
    \begin{scope}[shift={(0,-5)}]
        \draw[domain=0:180] plot ({2*cos(\x)}, {2*sin(\x)});
        \draw (-2,0) arc (180:360:2 and 0.6);
        \draw[dashed] (2,0) arc (0:180:2 and 0.6);
        \node at (0,1.25) {$\bullet$};
    \end{scope}
}
=
\tikz[baseline=-4ex, scale=.2]{
    \draw (-2,0) arc (180:360:2 and 0.6);
    \draw (2,0) arc (0:180:2 and 0.6);
    \draw (-2,-5) -- (-2,0);
    \draw (2,-5) -- (2,0);
    \begin{scope}[shift={(0,-5)}]
        \draw (-2,0) arc (180:360:2 and 0.6);
        \draw[dashed] (2,0) arc (0:180:2 and 0.6);
    \end{scope}
}~.
\]
The first three relations are called the \textit{sphere relations} (referred to as S0, S1, and S2 respectively), and the last relation is called the \textit{tube-cutting relation}. Interestingly, the sphere with three dots does not have an evaluation. The most general remedy is cosmetic, and it is treated as a free variable. Explicitly, in $\text{Cob}(n)$, we declare a fourth sphere relation by setting
\[
\tikz[baseline=-.5ex, scale = .2]{
  \draw (0,0) circle (2cm);
  \draw (-2,0) arc (180:360:2 and 0.6);
  \draw[dashed] (2,0) arc (0:180:2 and 0.6);
  \node at (-.8,0) {$\bullet$};
  \node at (0,0) {$\bullet$};
  \node at (.8,0) {$\bullet$};
} = \alpha.
\]
However, in what follows, we will take $\alpha$ to be zero; that is, we will replace the last sphere relation (S2) with the relation
\[
\tikz[baseline=-.35ex, scale = .3]{
    \draw (-1, -2) -- (-1, 1) -- (1,2) -- (1, -1) -- (-1, -2);
    \node at (0,0.5) {$\bullet$};
    \node at (0,-0.5) {$\bullet$};
} = 0.
\]

\begin{lemma}
There is an isomorphism of $\mathbb{Z}[q,q^{-1}]$-algebras
\[
K_0(\mathrm{Cob}(n)) \cong TL_n.
\]
\end{lemma}

\begin{proof}
Multiplication by $q$ defines an endofunctor $\text{Cob}(n) \to \text{Cob}(n)$, which in turn determines an endomorphism on $K_0(\text{Cob}(n))$, making it a $\mathbb{Z}[q,q^{-1}]$-algebra. Then the result is immediate.
\end{proof}

\medskip
\noindent \textit{Step 3}: Finally, we'd like a way to assign to a tangle in the 3-ball with $2n$ marked points some collection of objects in $\mathrm{Cob}(n)$.

\begin{definition}
Let 
\[
\mathrm{Kom}(n) = \mathrm{Kom}(\mathrm{Mat}(\mathrm{Cob}(n)))
\]
denote the category of partially bounded chain complexes of finite direct sums of objects in $\mathrm{Cob}(n)$. In this paper, we allow complexes with unbounded negative homological degree in keeping with~\cite{https://doi.org/10.48550/arxiv.1405.2574}, but opposed to, for example,~\cite{https://doi.org/10.48550/arxiv.1005.5117}.
\end{definition}

The tensor product of chain complexes extends $\otimes$ in $\mathrm{Cob}(n)$ to $\mathrm{Kom}(n)$: schematically,
\begin{align*}
C \otimes D &= \left(\begin{tikzcd}[column sep=small, scale=.915, ampersand replacement=\&]
\cdots \arrow[r] \&
\tikz[baseline=.8ex, scale = .4]{
\draw (0,0) rectangle (1,1);
\draw (.5,-.5) -- (.5,0);
\draw (.5,1) -- (.5, 1.5);
\node[scale=0.9] at (.5,.25) {$C_2$};
} \arrow[r] \& 
\tikz[baseline=.8ex, scale = .4]{
\draw (0,0) rectangle (1,1);
\draw (.5,-.5) -- (.5,0);
\draw (.5,1) -- (.5, 1.5);
\node[scale=0.9] at (.5,.25) {$C_1$};
} \arrow[r] \&
\tikz[baseline=.8ex, scale = .4]{
\draw (0,0) rectangle (1,1);
\draw (.5,-.5) -- (.5,0);
\draw (.5,1) -- (.5, 1.5);
\node[scale=0.9] at (.5,.25) {$C_0$};
}
\end{tikzcd}\right) \otimes
\left(\begin{tikzcd}[column sep=small, scale=.915, ampersand replacement=\&]
\cdots \arrow[r] \&
\tikz[baseline=.8ex, scale = .4]{
\draw (0,0) rectangle (1,1);
\draw (.5,-.5) -- (.5,0);
\draw (.5,1) -- (.5, 1.5);
\node[scale=0.9] at (.5,.25) {$D_2$};
} \arrow[r] \& 
\tikz[baseline=.8ex, scale = .4]{
\draw (0,0) rectangle (1,1);
\draw (.5,-.5) -- (.5,0);
\draw (.5,1) -- (.5, 1.5);
\node[scale=0.9] at (.5,.25) {$D_1$};
} \arrow[r] \&
\tikz[baseline=.8ex, scale = .4]{
\draw (0,0) rectangle (1,1);
\draw (.5,-.5) -- (.5,0);
\draw (.5,1) -- (.5, 1.5);
\node[scale=0.9] at (.5,.25) {$D_0$};
} 
\end{tikzcd}\right)
\\ &=
\begin{tikzcd}[scale=.95, ampersand replacement=\&]
\cdots \arrow[r] \& 
\tikz[baseline=-1.2ex, scale = .4]{
\draw (0,0) rectangle (1,1);
\draw (.5,-.5) -- (.5,0);
\draw (.5,1) -- (.5, 1.5);
\node[scale=0.9] at (.5,.25) {$C_2$};
    \begin{scope}[shift={(0,-1.5)}]
    \draw (0,0) rectangle (1,1);
    \draw (.5,-.5) -- (.5,0);
    \node[scale=0.9] at (.5,.25) {$D_0$};
\end{scope}
} \oplus 
\tikz[baseline=-1.2ex, scale = .4]{
\draw (0,0) rectangle (1,1);
\draw (.5,-.5) -- (.5,0);
\draw (.5,1) -- (.5, 1.5);
\node[scale=0.9] at (.5,.25) {$C_1$};
    \begin{scope}[shift={(0,-1.5)}]
    \draw (0,0) rectangle (1,1);
    \draw (.5,-.5) -- (.5,0);
    \node[scale=0.9] at (.5,.25) {$D_1$};
\end{scope}
} \oplus 
\tikz[baseline=-1.2ex, scale = .4]{
\draw (0,0) rectangle (1,1);
\draw (.5,-.5) -- (.5,0);
\draw (.5,1) -- (.5, 1.5);
\node[scale=0.9] at (.5,.25) {$C_0$};
    \begin{scope}[shift={(0,-1.5)}]
    \draw (0,0) rectangle (1,1);
    \draw (.5,-.5) -- (.5,0);
    \node[scale=0.9] at (.5,.25) {$D_2$};
\end{scope}
}
\arrow[r] \& 
\tikz[baseline=-1.2ex, scale = .4]{
\draw (0,0) rectangle (1,1);
\draw (.5,-.5) -- (.5,0);
\draw (.5,1) -- (.5, 1.5);
\node[scale=0.9] at (.5,.25) {$C_1$};
    \begin{scope}[shift={(0,-1.5)}]
    \draw (0,0) rectangle (1,1);
    \draw (.5,-.5) -- (.5,0);
    \node[scale=0.9] at (.5,.25) {$D_0$};
\end{scope}
} \oplus 
\tikz[baseline=-1.2ex, scale = .4]{
\draw (0,0) rectangle (1,1);
\draw (.5,-.5) -- (.5,0);
\draw (.5,1) -- (.5, 1.5);
\node[scale=0.9] at (.5,.25) {$C_0$};
    \begin{scope}[shift={(0,-1.5)}]
    \draw (0,0) rectangle (1,1);
    \draw (.5,-.5) -- (.5,0);
    \node[scale=0.9] at (.5,.25) {$D_1$};
\end{scope}
}
\arrow[r] \&
\tikz[baseline=-1.2ex, scale = .4]{
\draw (0,0) rectangle (1,1);
\draw (.5,-.5) -- (.5,0);
\draw (.5,1) -- (.5, 1.5);
\node[scale=0.9] at (.5,.25) {$C_0$};
    \begin{scope}[shift={(0,-1.5)}]
    \draw (0,0) rectangle (1,1);
    \draw (.5,-.5) -- (.5,0);
    \node[scale=0.9] at (.5,.25) {$D_0$};
\end{scope}
}
\end{tikzcd}.
\end{align*}

Indeed, passing to the homotopy category of a pre-additive category does not change the Grothendieck group up to isomorphism; see Section 2.7.1. of \cite{https://doi.org/10.48550/arxiv.1005.5117} for a full discussion.

\begin{lemma}
There is an isomorphism of $\mathbb{Z}[q]\llbracket q^{-1}\rrbracket$-algebras
\[
K_0(\mathrm{Kom}(n)) \cong TL_n.
\]
\end{lemma}

Let $\llbracket T \rrbracket$ denote the complex corresponding to a tangle $T$, obtained by the skein relations
\[
\left\llbracket \tikz[baseline=1.6ex, scale = .7]{
\draw[knot, -stealth] (0,0) -- (1,1);
\draw[knot] (1,0) -- (.7,.3);
\draw[knot, -stealth] (.3,.7) -- (0,1);
} \right\rrbracket
 = \begin{tikzcd}[]
q ~ \myvec{\tikz[baseline=1.6ex, scale = .7]{
\draw[knot, rounded corners = 4mm] (0,0) -- (.45,.5) -- (0,1);
\draw[knot, rounded corners = 4mm] (1,0) -- (.55,.5) -- (1,1);
}} \arrow{r}{\tikz[baseline=1.6ex, scale = .65]{\draw[knot, rounded corners = 4mm,-] (0,0) -- (.45,.5) -- (0,1);\draw[knot, rounded corners = 4mm,-] (1,0) -- (.55,.5) -- (1,1);\draw[red,thick,-]  (.223,.5) -- (.777,.5);}} &
q^2 ~ \tikz[baseline=1.6ex, scale = .7]{
\draw[knot, rounded corners = 2mm] (0,0) -- (0,.25) -- (.5,.49) -- (1,.25) -- (1,0);
\draw[knot, rounded corners = 2mm] (0,1) -- (0,.75) -- (.5,.51) -- (1,.75) -- (1,1);
}
\end{tikzcd}
~ \quad \text{and} \quad \,
\left\llbracket \tikz[baseline=1.6ex, scale = .7]{
\draw[knot, -stealth] (1,0) -- (0,1);
\draw[knot] (0,0) -- (.3, .3);
\draw[knot, -stealth] (.7, .7) -- (1,1)
} \right\rrbracket
 = \begin{tikzcd}[]
 q^{-2} ~\tikz[baseline=1.6ex, scale = .7]{
\draw[knot, rounded corners = 2mm] (0,0) -- (0,.25) -- (.5,.49) -- (1,.25) -- (1,0);
\draw[knot, rounded corners = 2mm] (0,1) -- (0,.75) -- (.5,.51) -- (1,.75) -- (1,1);
}~ \arrow{r}{\tikz[baseline=1.6ex, scale = .65]{\draw[knot, rounded corners = 2mm,-] (0,0) -- (0,.25) -- (.5,.49) -- (1,.25) -- (1,0);\draw[knot, rounded corners = 2mm,-] (0,1) -- (0,.75) -- (.5,.51) -- (1,.75) -- (1,1);\draw[red,thick,-] (.5,.425) -- (.5,.575);}} & 
q^{-1} ~\myvec{\tikz[baseline=1.6ex, scale = .7]{
\draw[knot, rounded corners = 4mm] (0,0) -- (.45,.5) -- (0,1);
\draw[knot, rounded corners = 4mm] (1,0) -- (.55,.5) -- (1,1);
}}
\end{tikzcd}
\]
where the $\myvec{\text{underlined}}$ term is in homological degree zero. For example, 
\[
\left\llbracket \tikz[baseline=4.5ex, scale = .8]{
\draw[knot, rounded corners = 2mm,-stealth] (0,0) -- (1,1) -- (0,2);
\draw[knot, rounded corners = 2mm] (1,0) -- (.6,.4);
\draw[knot, rounded corners = 2mm] (.4,.6) -- (0,1) -- (.4,1.4);
\draw[knot, rounded corners = 2mm,-stealth] (.6,1.6) -- (1,2);
}\right\rrbracket ~ = ~ 
\begin{tikzcd}[column sep = huge,ampersand replacement=\&]
q^{-1}~\tikz[baseline=4.5ex, scale = .8]{
\draw[knot, rounded corners = 2mm] (0,0) -- (.45, .5) -- (0,1) -- (.5,1.45) -- (1,1) -- (.55,.5) -- (1,0);
\draw[knot, rounded corners = 2mm] (0,2) -- (.5,1.55) -- (1,2);
}~ \arrow{r}{\begin{pmatrix}
\tikz[baseline=4.5ex, scale = .55]{
\draw[knot, rounded corners = 2mm,-] (0,0) -- (.45, .5) -- (0,1) -- (.5,1.45) -- (1,1) -- (.55,.5) -- (1,0);
\draw[knot, rounded corners = 2mm,-] (0,2) -- (.5,1.55) -- (1,2);
\draw[red,thick,-] (.5,1.33) -- (.5,1.67);
} &
\tikz[baseline=4.5ex, scale = .55]{
\draw[knot, rounded corners = 2mm,-] (0,0) -- (.45, .5) -- (0,1) -- (.5,1.45) -- (1,1) -- (.55,.5) -- (1,0);
\draw[knot, rounded corners = 2mm,-] (0,2) -- (.5,1.55) -- (1,2);
\draw[red,thick,-] (.33,.5) -- (.67,.5);
}
\end{pmatrix}^\top} \& 
q^0~\left(~\myvec{\tikz[baseline=4.5ex, scale = .8]{
\draw[knot, rounded corners = 2mm] (0,0) -- (.45,.5) -- (0,1) -- (.45,1.5) -- (0,2);
\draw[knot, rounded corners = 2mm] (1,0) -- (.55,.5) -- (1,1) -- (.55,1.5) -- (1,2);
} \oplus \tikz[baseline=4.5ex, scale = .8]{
\draw[knot, rounded corners = 2mm] (0,0) -- (.5,.45) -- (1,0);
\draw[knot, rounded corners = 2mm] (0,2) -- (.5,1.55) -- (1,2);
\draw[knot] (.5,1) circle (11pt);
}}~\right) \arrow{r}{\begin{pmatrix}
\tikz[baseline=4.5ex, scale = .55]{
\draw[knot, rounded corners = 2mm,-] (0,0) -- (.45,.5) -- (0,1) -- (.45,1.5) -- (0,2);
\draw[knot, rounded corners = 2mm,-] (1,0) -- (.55,.5) -- (1,1) -- (.55,1.5) -- (1,2);
\draw[red,thick,-] (.33,.5) -- (.67,.5);
} & \tikz[baseline=4.5ex, scale = .55]{
\draw[knot, rounded corners = 2mm,-] (0,0) -- (.5,.45) -- (1,0);
\draw[knot, rounded corners = 2mm,-] (0,2) -- (.5,1.55) -- (1,2);
\draw[knot] (.5,1) circle (11pt);
\draw[red,thick,-] (.5,1.4) -- (.5,1.67);
}
\end{pmatrix}} \&
q~\tikz[baseline=4.5ex, scale = .8]{
\draw[knot, rounded corners = 2mm] (0,0) -- (.5,.45) -- (1,0);
\draw[knot, rounded corners = 2mm] (0,2) -- (.45,1.5) -- (0,1) -- (.5,.55) -- (1,1) -- (.55,1.5) -- (1,2);
}
\end{tikzcd}
\]
Notice that there is a free loop in homological grading zero, hence we may apply the delooping operations to yield the complex
\begin{equation}\label{rII-diagram}
\begin{tikzcd}[column sep = large,ampersand replacement=\&]
q^{-1}~\tikz[baseline=4.5ex, scale = .8]{
\draw[knot, rounded corners = 2mm] (0,0) -- (.45, .5) -- (0,1) -- (.5,1.45) -- (1,1) -- (.55,.5) -- (1,0);
\draw[knot, rounded corners = 2mm] (0,2) -- (.5,1.55) -- (1,2);
}~ \arrow{r}{A^\top} \& 
q^{-1} ~\tikz[baseline=4.5ex, scale = .8]{
\draw[knot, rounded corners = 2mm] (0,0) -- (.5,.45) -- (1,0);
\draw[knot, rounded corners = 2mm] (0,2) -- (.5,1.55) -- (1,2);
} \oplus ~q^0~\tikz[baseline=4.5ex, scale = .8]{
\draw[knot, rounded corners = 2mm] (0,0) -- (.45,.5) -- (0,1) -- (.45,1.5) -- (0,2);
\draw[knot, rounded corners = 2mm] (1,0) -- (.55,.5) -- (1,1) -- (.55,1.5) -- (1,2);
} \oplus ~ q^1\tikz[baseline=4.5ex, scale = .8]{
\draw[knot, rounded corners = 2mm] (0,0) -- (.5,.45) -- (1,0);
\draw[knot, rounded corners = 2mm] (0,2) -- (.5,1.55) -- (1,2);
} \arrow{r}{B} \&
q~\tikz[baseline=4.5ex, scale = .8]{
\draw[knot, rounded corners = 2mm] (0,0) -- (.5,.45) -- (1,0);
\draw[knot, rounded corners = 2mm] (0,2) -- (.45,1.5) -- (0,1) -- (.5,.55) -- (1,1) -- (.55,1.5) -- (1,2);
}
\end{tikzcd}
\end{equation}
where
\[
A^\top = \begin{pmatrix}
\tikz[baseline=4.5ex, scale = .55]{
\draw[knot, rounded corners = 2mm,-] (0,0) -- (.45, .5) -- (0,1) -- (.5,1.45) -- (1,1) -- (.55,.5) -- (1,0);
\draw[knot, rounded corners = 2mm,-] (0,2) -- (.5,1.55) -- (1,2);
\draw[red,thick,-] (.5,1.33) -- (.5,1.67);
} &
\tikz[baseline=4.5ex, scale = .55]{
\draw[knot, rounded corners = 2mm,-] (0,0) -- (.45, .5) -- (0,1) -- (.5,1.45) -- (1,1) -- (.55,.5) -- (1,0);
\draw[knot, rounded corners = 2mm,-] (0,2) -- (.5,1.55) -- (1,2);
\draw[red,thick,-] (.33,.5) -- (.67,.5);
}
\end{pmatrix}^\top \circ \varphi = \begin{pmatrix}
\mathrm{id} &
\tikz[baseline=3ex, scale = .55]{
\draw[knot, rounded corners = 2mm,-] (0,0) -- (.45, .5) -- (0,1) -- (.5,1.45) -- (1,1) -- (.55,.5) -- (1,0);
\draw[knot, rounded corners = 2mm,-] (0,2) -- (.5,1.55) -- (1,2);
\draw[red,thick,-] (.5,1.33) -- (.5,1.67);
} & 
\tikz[baseline=3ex, scale = .55]{
\draw[knot, rounded corners = 2mm,-] (0,0) -- (.45, .5) -- (0,1) -- (.5,1.45) -- (1,1) -- (.55,.5) -- (1,0);
\draw[knot, rounded corners = 2mm,-] (0,2) -- (.5,1.55) -- (1,2);
\node at (.5,1.3)[circle,fill,inner sep=1pt]{};
}
\end{pmatrix}^\top
\]
and
\[
B = \psi \circ \begin{pmatrix}
\tikz[baseline=4.5ex, scale = .55]{
\draw[knot, rounded corners = 2mm,-] (0,0) -- (.45,.5) -- (0,1) -- (.45,1.5) -- (0,2);
\draw[knot, rounded corners = 2mm,-] (1,0) -- (.55,.5) -- (1,1) -- (.55,1.5) -- (1,2);
\draw[red,thick,-] (.33,.5) -- (.67,.5);
} & \tikz[baseline=4.5ex, scale = .55]{
\draw[knot, rounded corners = 2mm,-] (0,0) -- (.5,.45) -- (1,0);
\draw[knot, rounded corners = 2mm,-] (0,2) -- (.5,1.55) -- (1,2);
\draw[knot] (.5,1) circle (11pt);
\draw[red,thick,-] (.5,1.4) -- (.5,1.67);
}
\end{pmatrix} = 
\begin{pmatrix}
\tikz[baseline=3ex, scale = .55]{
\draw[knot, rounded corners = 2mm] (0,0) -- (.5,.45) -- (1,0);
\draw[knot, rounded corners = 2mm] (0,2) -- (.5,1.55) -- (1,2);
\node at (.5,1.69)[circle,fill,inner sep=1pt]{};
} & 
\tikz[baseline=3ex, scale = .55]{
\draw[knot, rounded corners = 2mm,-] (0,0) -- (.45,.5) -- (0,1) -- (.45,1.5) -- (0,2);
\draw[knot, rounded corners = 2mm,-] (1,0) -- (.55,.5) -- (1,1) -- (.55,1.5) -- (1,2);
\draw[red,thick,-] (.33,.5) -- (.67,.5);
} & \mathrm{id}
\end{pmatrix}.
\]

\subsubsection{Chain homotopy lemmas}
In~\cite{barnatan2006fast}, delooping was introduced alongside the following lemma from homotopy theory to simplify computations in Khovanov homology. 

\begin{lemma}[Simultaneous Gaussian elimination]
\label{SimGE}
Suppose $\mathcal{A}$ is a pre-additive category, 
and let $K_*$ be an object of $\mathrm{Kom}(\mathcal{A})$ of the form
\[
\begin{tikzcd} A_0 \oplus C_0 \arrow{r}{M_0} & A_1 \oplus B_1 \oplus C_1 \arrow{r}{M_1} & A_2 \oplus B_2 \oplus C_2 \arrow{r}{M_2} & \cdots \end{tikzcd}
\]
where $M_0 = \begin{pmatrix} a_0 & c_0 \\ d_0 & f_0 \\ g_0 & j_0 \end{pmatrix}$ and $M_i = \begin{pmatrix} a_i & b_i & c_i \\ d_i & e_i & f_i \\ g_i & h_i & j_i \end{pmatrix}$ for all $i>0$. If $a_{2i}:A_{2i} \to A_{2i+1}$ and $e_{2i_1}: B_{2i+1} \to B_{2i+2}$ are isomorphisms for all $i\ge 0$, then the chain complex $K_*$ is homotopy equivalent to the complex
\[
\begin{tikzcd} C_0 \arrow{r}{Q_0} & C_1 \arrow{r}{Q_1} & C_2 \arrow{r}{Q_2} & \cdots \end{tikzcd}
\]
where $\begin{cases} Q_{2i} = j_{2i} - g_{2i}a_{2i}^{-1}c_{2i} \\ Q_{2i+1} = j_{2i+1} - h_{2i+1}e_{2i+1}^{-1}f_{2i+1}\end{cases}$.
\end{lemma}

\begin{proof}
This is an application of the simpler ``Gaussian elimination,'' see~\cite{https://doi.org/10.48550/arxiv.1005.5117}. 
\end{proof}

As an application, note that we may apply simultaneous Gaussian elimination to the complex (\ref{rII-diagram}). The result is that the complex $\left\llbracket \tikz[baseline=2ex, scale = .4]{
\draw[knot, rounded corners = 2mm,-stealth] (0,0) -- (1,1) -- (0,2);
\draw[knot, rounded corners = 2mm] (1,0) -- (.6,.4);
\draw[knot, rounded corners = 2mm] (.4,.6) -- (0,1) -- (.4,1.4);
\draw[knot, rounded corners = 2mm,-stealth] (.6,1.6) -- (1,2);
}\right \rrbracket$ is homotopy equivalent (hereinafter written $\simeq$) to the chain complex $\begin{tikzcd}[column sep = tiny] 0 \arrow[r] & \myvec{\tikz[baseline=2ex, scale = .4]{
\draw[knot, rounded corners = 2mm] (0,.5) -- (.45,1) -- (0,1.5);
\draw[knot, rounded corners = 2mm] (1,.5) -- (.55,1) -- (1,1.5);
}} \arrow[r] & 0\end{tikzcd}$; i.e., the complex $\left\llbracket T \right\rrbracket$ is invariant, up to chain homotopy equivalence, under Reidemeister II moves for tangles. The following is due to Bar-Natan.

\begin{theorem}[Theorem 1 of \cite{Bar_Natan_2005}]
The homotopy class of the complex $\left \llbracket T \right\rrbracket$ regarded in $\mathrm{Kom}(n)$ is an invariant of the tangle $T$.
\end{theorem}

To conclude this subsection, we note that there is a notion of a zero object in $\mathrm{Kom}(n)$: we call a chain complex $K_*$ \textit{contractible} if $K_* \simeq 0$. The following is well known.

\begin{lemma}[Big collapse]
A chain complex $K_*$ of contractible chain complexes is, itself, contractible.
\end{lemma}

\subsubsection{Khovanov's arc algebras}
\label{sss:Khovanov arc algebra}

Another categorification, provided by Khovanov ~\cite{Khovanov_2002}, is given by the category of complexes of $H^n$-modules, where $H^n$ is the $n$th arc algebra, described below. These can be generalized to the unified setting; see \cite{naisse2017odd} for a thorough discussion. We will use arc algebras to describe odd Khovanov complexes for tangles, following \cite{putyra20152categorychronologicalcobordismsodd} and \cite{naisse2020odd}. A large portion of this paper is devoted to providing a small generalization Naisse-Putyra's construction, allowing one to perform Bar-Natanesque computations in a particular category of $H^n$-modules.

Consider the Temperley-Lieb 2-category $\mathcal{TL}$, whose
\begin{itemize}
    \item objects are natural numbers,
    \item 1-morphisms $\mathrm{Hom}_{\mathcal{TL}}(m,n)$ are isotopy classes of crossingless tangles embedded in the square with $2m$ marked points on the $[0,1]\times \{0\}$ axis and $2n$ marked points on the $[0,1]\times \{1\}$ axis, and
    \item 2-morphisms $\mathrm{Hom}_{\mathcal{TL}}(t,s)$ are cobordisms with corners from the crossingless tangle $t$ to $s$.
\end{itemize}
Write $B_m^n = \mathrm{Hom}_{\mathcal{TL}}(m,n)$. In the case that $m=0$, we write $B^n$ (respectively, $n=0$ is written $B_m$); this is the collection of \textit{crossingless matchings} of $n$ points fixed on the top axis (resp., $m$ on the bottom axis). We will write $|a| = n$ for $a\in B^n$. Composition of 1-morphisms is given by stacking: $B_n^p \times B_m^n \to B_m^p$ is given by $(s,t) \mapsto ts$. There is also a mirroring operation, $\overline{\,\cdot\,}:B_m^n\to B_n^m$, which flips tangles about the line $[0,1]\times \{1/2\}$. 

Let $a\in B^m$, $b\in B_n$, and $t\in B_m^n$. Then $atb$ is a closed 1-manifold. Let $s\in B_n^p$ and $c\in B_p$. Consider the cobordism
\[
(atb)(\overline{b}sc) \to a(ts)c
\]
given by contracting symmetric arcs of $b\overline{b}$. We denote this cobordism by $W_{abc}(t,s)$. It is minimal in the sense that its Euler characteristic is $-|b|$.

The last ingredient required for defining the arc algebra is Khovanov's Frobenius TQFT. Let $V = \mathbb{Z}\langle v_+,v_- \rangle$ denote the free abelian group generated by $v_+$ and $v_-$, and impose a grading on $V$ by $|v_+| = 1$ and $|v_-| = -1$. Consider the functor $\mathcal{F}_e: \mathrm{Pre\text{-}Cob}(0) \to \mathbb{Z}\mathrm{Mod}$ defined as follows. On objects,
\[
\mathcal{F}_e(\underbrace{\bigcirc \sqcup \cdots \sqcup \bigcirc}_{n}) = V^{\otimes n}.
\]
On morphisms, note that any cobordism may be decomposed into a sequence of merges, splits, births, and deaths, which are evaluated by $\mathcal{F}_e$ as listed below.
\begin{align*}
\mathcal{F}_e\left(\tikz[baseline=2.5ex, scale=.5]{
	\draw (0,0) .. controls (0,1) and (1,1) .. (1,2);
	\draw (1,0) .. controls (1,1) and (2,1) .. (2,0);
	\draw (3,0) .. controls (3,1) and (2,1) .. (2,2);
	\draw (0,0) .. controls (0,-.25) and (1,-.25) .. (1,0);
	\draw[dashed] (0,0) .. controls (0,.25) and (1,.25) .. (1,0);
	\draw (2,0) .. controls (2,-.25) and (3,-.25) .. (3,0);
	\draw[dashed] (2,0) .. controls (2,.25) and (3,.25) .. (3,0);
	\draw (1,2) .. controls (1,1.75) and (2,1.75) .. (2,2);
	\draw (1,2) .. controls (1,2.25) and (2,2.25) .. (2,2);
}\right) : V \otimes V \rightarrow V &= 
\begin{cases}
v_+ \otimes v_+ \mapsto v_+, & v_+ \otimes v_- \mapsto v_-, \\
v_- \otimes v_- \mapsto 0, & v_- \otimes v_+ \mapsto v_-,
\end{cases}
\\
\mathcal{F}_e\left(\tikz[baseline=9.5ex, scale=.5]{
	\draw  (1,2) .. controls (1,3) and (0,3) .. (0,4);
	\draw  (2,2) .. controls (2,3) and (3,3) .. (3,4);
	\draw (1,4) .. controls (1,3) and (2,3) .. (2,4);
	\draw (0,4) .. controls (0,3.75) and (1,3.75) .. (1,4);
	\draw (0,4) .. controls (0,4.25) and (1,4.25) .. (1,4);
	\draw (2,4) .. controls (2,3.75) and (3,3.75) .. (3,4);
	\draw (2,4) .. controls (2,4.25) and (3,4.25) .. (3,4);
	\draw (1,2) .. controls (1,1.75) and (2,1.75) .. (2,2);
	\draw[dashed] (1,2) .. controls (1,2.25) and (2,2.25) .. (2,2);
}\right) : V  \rightarrow V \otimes V &= 
\begin{cases}
v_+ \mapsto v_- \otimes v_+ + v_+ \otimes v_-,  &\\
v_- \mapsto v_- \otimes v_-, &
\end{cases}
\\
\mathcal{F}_e\left(\tikz[baseline=5.2ex, scale=.5]{
	\draw (1,2) .. controls (1,1) and (2,1) .. (2,2);
	\draw (1,2) .. controls (1,1.75) and (2,1.75) .. (2,2);
	\draw (1,2) .. controls (1,2.25) and (2,2.25) .. (2,2);
}\right) : \mathbb{Z}  \rightarrow V  &= 
\begin{cases}
1 \mapsto v_+, & 
\end{cases}
\\
\mathcal{F}_e\left(\tikz[baseline=.2ex, scale=.5]{
	\draw (1,0) .. controls (1,1) and (2,1) .. (2,0);
	\draw (1,0) .. controls (1,-.25) and (2,-.25) .. (2,0);
	\draw[dashed] (1,0) .. controls (1,.25) and (2,.25) .. (2,0);
}\right) : V  \rightarrow \mathbb{Z} &= 
\begin{cases}
v_+ \mapsto 0, & \\
v_- \mapsto 1. &
\end{cases}
\end{align*}
For example, a cylinder with a hole in it can be decomposed into a split followed by a merge. Clearly, this maps $v_+ \mapsto 2v_-$ and $v_-\mapsto 0$. So, altering $\mathrm{Cob}$ so that objects can be decorated by dots, we have that
\[
\mathcal{F}_e\left(\tikz[baseline=2.5ex, scale=.5]{
	\draw (1,0) .. controls (1,.-.25) and (2,-.25) .. (2,0);
	\draw[dashed] (1,0) .. controls (1,.25) and (2,.25) .. (2,0);
	\draw (1,0) -- (1,2);
	\draw (2,0) -- (2,2);
	\draw (1,2) .. controls (1,1.75) and (2,1.75) .. (2,2);
	\draw (1,2) .. controls (1,2.25) and (2,2.25) .. (2,2);
	\node at (1.5,1) {$\bullet$};
}\right) : V \rightarrow V = 
\begin{cases}
v_+ \mapsto v_-, & \\
v_- \mapsto 0. &
\end{cases}
\]
$\mathcal{F}_e$ extends to $\mathrm{Mat}(\mathrm{Pre\text{-}Cob}(n))$, and one can easily verify that $\mathcal{F}_e$ satisfies the each of the sphere and tube-cutting relations.

Let $t\in B_m^n$. The \textit{arc space} of $t$ is defined
\[
\mathcal{F}_e(t) = \bigoplus_{a\in B^m, b\in B_n} \mathcal{F}_e(atb).
\]
Given another tangle $s\in B_n^p$, define the composition map
\begin{align*}
&\mu[t,s]:\mathcal{F}_e(atb) \otimes \mathcal{F}_e(b'sc) \to \mathcal{F}_e(a (ts) c)\\ &\text{by}~\mu[t,s] = \begin{cases} 0 & \text{if}~\overline{b}\not=b' \\ \mathcal{F}_e(W_{abc}(t,s)) & \text{if}~\overline{b}=b' \end{cases}
\end{align*}
for $b'\in B^n$ and $c\in B_p$.

\begin{definition}
The \textit{arc algebra} $H^n$ is the arc space
\[
H^n = \mathcal{F}(1_n) = \bigoplus_{a\in B^m, b\in B_n} \mathcal{F}_e(a 1_n b)
\]
with multiplication $\mu[1_n,1_n]$.
\end{definition}

It is more work, but the category of left $H^n$-modules provides another categorification of the Temperley-Lieb algebra; see Section 5.2 of \cite{Khovanov_2002} for details.

\begin{lemma}
There is an isomorphism of $\mathbb{Z}[q,q^{-1}]$-algebras
\[
K_0(H^n\mathrm{PMod}) \cong TL_n \qquad \text{and} \qquad K_0(\mathrm{Kom}(H^n\mathrm{PMod})) \cong TL_n.
\]
for $H^n\mathrm{PMod}$ the category of \textit{projective} $H^n$-modules.
\end{lemma}

\subsection{Cooper-Krushkal projectors}
\label{ss:ckprojectors}

The first categorification of Jones-Wenzl projectors was described by Cooper and Krushkal in~\cite{https://doi.org/10.48550/arxiv.1005.5117}. Their definition mirrors that of the Jones-Wenzl projectors, and they are uniquely defined in $\mathrm{Kom}(n)$ (that is, up to homotopy equivalence). Everything presented here still holds if we replace $\mathrm{Kom}(n)$ with $\mathrm{Kom}(H^n\mathrm{PMod})$.

\begin{definition}
\label{CKprojector}
A negativiely graded chain complex $(C_*,d_*) \in \mathrm{Kom}(n)$  with degree zero differential and is called a \textit{Cooper-Krushkal projector} if it satisfies the following axioms:
\begin{enumerate}
    \item[(CK1)] $C_0 = 1_n$ and the identity does not appear in any other homological degree.
    \item[(CK2)] $C_*$ is contractible under turnbacks: for any $e_i\in TL_n$, $C_*\otimes e_i \simeq e_i \otimes C_* \simeq 0$.
\end{enumerate}
The second axiom is referred to as ``turnback killing.''
\end{definition}
Notice that, by construction, if $C\in \mathrm{Kom}(n)$ is a Cooper-Krushkal projector, then $[C]\in K_0(\mathrm{Kom}(n)) \cong TL_n$ satisfies (JW1) and (JW2), so $[C] = p_n\in TL_n$.

Like the Jones-Wenzl projectors, homotopy uniqueness of the Cooper-Krushkal projectors follows from little work. The main tool is the following generalization of idempotence (whose analogue also holds for Jones-Wenzl projectors).

\begin{proposition}
\label{GenIdempotence}
Suppose $C\in \mathrm{Kom}(m)$ and $D \in \mathrm{Kom}(n)$ are Cooper-Krushkal projectors with $0 \le m \le n$. Then
\[
C\otimes \left( D \sqcup 1_{n-m} \right) \simeq C \simeq \left(D\sqcup 1_{n-m}\right) \otimes C.
\]
\end{proposition}

Homotopy idempotence and uniqueness are then corollaries.

\begin{proof}
See Proposition 3.3 of \cite{https://doi.org/10.48550/arxiv.1005.5117}.
\end{proof}

The main theorem of~\cite{https://doi.org/10.48550/arxiv.1005.5117} is the following.

\begin{theorem}[Theorem 3.2 of \cite{https://doi.org/10.48550/arxiv.1005.5117}]
For each $n>0$, there exists a chain complex $C \in \mathrm{Kom}(n)$ that is a Cooper-Krushkal projector.
\end{theorem}

We will write $P_n^{\mathrm{CK}}$ to denote the $n$th Cooper-Krushkal projector (or a representative of it), so that $[P_n^{\mathrm{CK}}] = p_n$. We represent Cooper-Krushkal projectors via numbered boxes, as we did the Jones-Wenzl projectors. For example, here is a Jones-Wenzl projector when $n=2$:
\[
\begin{tikzcd}
\tikz[baseline=.8ex, scale = .4]{
\draw (0,0) rectangle (1,1);
\draw[knot] (.5,-.5) -- (.5,0);
\draw[knot] (.5,1) -- (.5, 1.5);
\node at (.5,.3) {$2$};
} =\quad  \cdots \arrow[r, "C_{-4}"] & 
q^{-5} \tikz[baseline=1.6ex, scale = .7]{
\draw[knot, rounded corners = 2mm] (0,0) -- (0,.25) -- (.5,.49) -- (1,.25) -- (1,0);
\draw[knot, rounded corners = 2mm] (0,1) -- (0,.75) -- (.5,.51) -- (1,.75) -- (1,1);
} \arrow[r, "C_{-3}"] & 
q^{-3} \tikz[baseline=1.6ex, scale = .7]{
\draw[knot, rounded corners = 2mm] (0,0) -- (0,.25) -- (.5,.49) -- (1,.25) -- (1,0);
\draw[knot, rounded corners = 2mm] (0,1) -- (0,.75) -- (.5,.51) -- (1,.75) -- (1,1);
} \arrow[r, "C_{-2}"] & 
q^{-1} \tikz[baseline=1.6ex, scale = .7]{
\draw[knot, rounded corners = 2mm] (0,0) -- (0,.25) -- (.5,.49) -- (1,.25) -- (1,0);
\draw[knot, rounded corners = 2mm] (0,1) -- (0,.75) -- (.5,.51) -- (1,.75) -- (1,1);
} \arrow[r, "C_{-1}"] & 
\myvec{\tikz[baseline=1.6ex, scale = .7]{
\draw[knot, rounded corners = 4mm] (0,0) -- (.45,.5) -- (0,1);
\draw[knot, rounded corners = 4mm] (1,0) -- (.55,.5) -- (1,1);
}}
\end{tikzcd}
\]
where
\[
C_i = \begin{cases} \tikz[baseline=1.6ex, scale = .65]{\draw[knot, rounded corners = 2mm,-] (0,0) -- (0,.25) -- (.5,.49) -- (1,.25) -- (1,0);\draw[knot, rounded corners = 2mm,-] (0,1) -- (0,.75) -- (.5,.51) -- (1,.75) -- (1,1);\draw[red,thick,-] (.5,.425) -- (.5,.575);} & i=-1 \\[10pt] 
\tikz[baseline=1.6ex, scale = .65]{\draw[knot, rounded corners = 2mm,-] (0,0) -- (0,.25) -- (.5,.49) -- (1,.25) -- (1,0);\draw[knot, rounded corners = 2mm,-] (0,1) -- (0,.75) -- (.5,.51) -- (1,.75) -- (1,1);\node at (.5, .58) {$\bullet$};} - \tikz[baseline=1.6ex, scale = .65]{\draw[knot, rounded corners = 2mm,-] (0,0) -- (0,.25) -- (.5,.49) -- (1,.25) -- (1,0);\draw[knot, rounded corners = 2mm,-] (0,1) -- (0,.75) -- (.5,.51) -- (1,.75) -- (1,1);\node at (.5, .38) {$\bullet$};} & i=-2k\\[10pt]
\tikz[baseline=1.6ex, scale = .65]{\draw[knot, rounded corners = 2mm,-] (0,0) -- (0,.25) -- (.5,.49) -- (1,.25) -- (1,0);\draw[knot, rounded corners = 2mm,-] (0,1) -- (0,.75) -- (.5,.51) -- (1,.75) -- (1,1);\node at (.5, .58) {$\bullet$};} + \tikz[baseline=1.6ex, scale = .65]{\draw[knot, rounded corners = 2mm,-] (0,0) -- (0,.25) -- (.5,.49) -- (1,.25) -- (1,0);\draw[knot, rounded corners = 2mm,-] (0,1) -- (0,.75) -- (.5,.51) -- (1,.75) -- (1,1);\node at (.5, .38) {$\bullet$};} & i=-2k -1\end{cases}
\]
for all positive integers $k$. It is straightforward to check that this is an element of $\mathrm{Kom}(n)$, and that it satisfies axioms (CK1) And (CK2).

This categorification succeeds in possessing many properties analogous to the original object. In particular, if $\mathrm{Tr}^n$ denotes the (complete) Markov trace applied to each entry and differential in the chain complex, we have that the graded Euler characteristic of the homology of the trace of each projector is a quantum integer; i.e.,
\[
\chi(H_*(\mathrm{Tr}^n(P_n^{\mathrm{CK}}))) = [n+1].
\]
For example, it is also straightforward to verify that, for $k$ a positive integer and $\alpha \equiv 0$,
\[
H_n(\mathrm{Tr}^2(P_2^{\mathrm{CK}})) = \begin{cases} q^2 \mathbb{Z} \oplus \mathbb{Z} & n=0  \\ 0 & n=-1 \\ q^{-4k+2} \mathbb{Z} \oplus q^{-4k} \mathbb{Z}/2\mathbb{Z} & n=-2k \\ q^{-4k-2} \mathbb{Z} & n=-2k-1\end{cases}
\]
It is interesting that the homology of $\mathrm{Tr}^n(P_n^{\mathrm{CK}})$ is not spanned only by classes which correspond to coefficients of the graded Euler characteristic. This turns out to be the case for the projectors of odd Khovanov homology as well. Moreover, the two homologies disagree (for example, there is no torsion for the odd, 2-stranded projector) but their graded Euler characteristics coincide.

\newpage

\section{The Odd setting: chronologies and \texorpdfstring{$\mathcal{G}$}{Lg}-graded structures}
\label{S:odd chronologies and stuff}

In this section, we provide a modern introduction to odd Khovanov homology. That is, rather than detailing the projective TQFT of Ozsv\'ath-Rasmussen-Szab\'o, we discuss Putyra's 2-category of chronological cobordisms and its linearlization over the ground ring $R:= \mathbb{Z}[X, Y, Z^{\pm1}] / (X^2 = Y^2 = 1)$ in \S \ref{ss:chronologicalcobordisms and coc}. In \S \ref{ss:unified arc algebras}, we attempt to mimic the constructions of \cite{Khovanov_2002}, as outlined briefly in \S \ref{sss:Khovanov arc algebra}. Here, we discover the challenges motivating the next few sections of our work: unified arc algebras are not associative in this context, and the composition maps $\mu$ are not degree-preserving. Finally, in \S \ref{ss:brief graded outline}, we give a description of the solution posed by Naisse and Putrya in \cite{naisse2020odd}. We hope that \S \ref{ss:brief graded outline} serves as a roadmap and extended outline of Sections \ref{s:gradingmultis and pads} and \ref{S: SHIFTING SYSTEMS}.

\subsection{Chronological cobordisms and changes of chronology}
\label{ss:chronologicalcobordisms and coc}

First introduced by Putyra~\cite{putyra2010cobordismschronologiesgeneralisationkhovanov, putyra20152categorychronologicalcobordismsodd}, we will proceed using the definition of chronological cobordisms provided by Sch\"utz in \cite{Sch_tz_2022}.

\begin{definition}
A \textit{chronological cobordism} between closed 1-manifolds $S_0$ and $S_1$ is a cobordism $W$ between $S_0$ and $S_1$ embedded into $\mathbb{R}^2\times [0,1]$ such that
\begin{enumerate}[label=(\roman*)]
    \item there is an $\epsilon > 0$ such that
    \[
    W \cap (\mathbb{R}^2 \times [0,1]) = S_0 \times [0,\epsilon] ~\quad~ \text{and} ~\quad~ W \cap (\mathbb{R}^2 \times [1-\epsilon, 1]) = S_1 \times [1-\epsilon, 1]
    \]
    and
    \item the height function $\tau: W \to [0,1]$ given by projection onto the third coordinate is a Morse function for which $\#\tau^{-1}(\{c\}) = 1$ whenever $c$ is a critical value of $\tau$. We call such a Morse function \textit{separative}.
\end{enumerate}
\end{definition}
Next, a \textit{framing} on a chronological cobordism is a choice of orientation of a basis for each unstable manifold $W_p \subset W$, for $p$ a critical point of $\tau$ of index 1 or 2. We will assume all chronological cobordisms to be framed. Since a framing is determined by a choice of tangent vector on each unstable manifold determined by a critical point, it is standard to visualize the framing by an arrow through critical points. We'll adapt the 2-dimensional notation to 1-dimensional diagrams appropriately; for example,
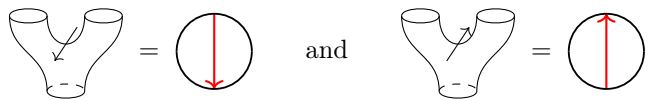
\[
\tikz[baseline=9.5ex, scale=.5]{
	\draw  (1,2) .. controls (1,3) and (0,3) .. (0,4);
	\draw  (2,2) .. controls (2,3) and (3,3) .. (3,4);
	\draw (1,4) .. controls (1,3) and (2,3) .. (2,4);
	\draw (0,4) .. controls (0,3.75) and (1,3.75) .. (1,4);
	\draw (0,4) .. controls (0,4.25) and (1,4.25) .. (1,4);
	\draw (2,4) .. controls (2,3.75) and (3,3.75) .. (3,4);
	\draw (2,4) .. controls (2,4.25) and (3,4.25) .. (3,4);
	\draw (1,2) .. controls (1,1.75) and (2,1.75) .. (2,2);
	\draw[dashed] (1,2) .. controls (1,2.25) and (2,2.25) .. (2,2);
        \draw[->] (1.8,3.7) -- (1.2,2.8);
} ~=~ \tikz[baseline=-.6ex, scale = .5]{
    \draw[red,thick,->] (0,1) -- (0,-1);
    \draw[knot] (-1,0) arc (180:360:1 and 1);
    \draw[knot] (1,0) arc (0:180:1 and 1);
    } \qquad \text{and} \qquad 
    \tikz[baseline=9.5ex, scale=.5]{
	\draw  (1,2) .. controls (1,3) and (0,3) .. (0,4);
	\draw  (2,2) .. controls (2,3) and (3,3) .. (3,4);
	\draw (1,4) .. controls (1,3) and (2,3) .. (2,4);
	\draw (0,4) .. controls (0,3.75) and (1,3.75) .. (1,4);
	\draw (0,4) .. controls (0,4.25) and (1,4.25) .. (1,4);
	\draw (2,4) .. controls (2,3.75) and (3,3.75) .. (3,4);
	\draw (2,4) .. controls (2,4.25) and (3,4.25) .. (3,4);
	\draw (1,2) .. controls (1,1.75) and (2,1.75) .. (2,2);
	\draw[dashed] (1,2) .. controls (1,2.25) and (2,2.25) .. (2,2);
        \draw[<-] (1.8,3.7) -- (1.2,2.8);
} ~=~ \tikz[baseline=-.6ex, scale = .5]{
    \draw[red,thick,<-] (0,1) -- (0,-1);
    \draw[knot] (-1,0) arc (180:360:1 and 1);
    \draw[knot] (1,0) arc (0:180:1 and 1);
    }~.
\]

Naturally, two chronological cobordisms are considered equivalent if they can be related by a diffeotopy $H_t$, $t\in [0,1]$, so that projection of $H_t(W)$ onto the third coordinate is a separative Morse function at each time $t$. This is a much more strict equivalence relation than that of the even case. To account for this, Putyra introduces the following action/relation. A \textit{change of chronology} is a diffeotopy $H_t$ such that projection of $H_t(W)$ onto the third coordinate is a generic homotopy of Morse functions, together with a smooth choice of framings on $H_t(W)$. Two changes of chronology between equivalent cobordisms are equivalent if they are homotopic in the space of oriented Igusa functions after composing with the equivalences of cobordisms; for a thorough description, consult \cite{putyra20152categorychronologicalcobordismsodd}. We write $H: W_1 \Rightarrow W_2$ for a change of chronology $H$ between chronological cobordisms $W_1$ and $W_2$.

\begin{definition}
A change of chronology $H$ on a chronological cobordism $W$ is called \textit{locally vertical} if there is a finite collection of cylinders $\{C_i\}_i$ in $\mathbb{R}^2 \times I$ such that $H$ is the identity on $W - \bigcup_i C_i$.
\end{definition}

We will use locally vertical changes of chronology frequently. Their main utility stems from the fact that they are unique up to homotopy.

\begin{proposition}[Proposition 4.4 of \cite{putyra20152categorychronologicalcobordismsodd}]
\label{PutyraHammer}
If $H$ and $H'$ are locally vertical changes of chronology (with respect to the same cylinders) with the same source and target, then they are homotopic in the space of framed diffeotopies.
\end{proposition}

There are two different ways of composing changes of chronology. First, given a sequence of cobordisms $A \xrightarrow{W} B \xrightarrow{W'} C$, and changes of chronology $H$ on $W$ and $H'$ on $W'$, there is a change of chronology $H' \circ H$ on $W' \circ W$. Second, given a sequence of changes of chronology $W \xRightarrow{H} W' \xRightarrow{H'} W''$, we will denote their composition by $H' \star H$.

On the other hand, we may completely describe the \textit{elementary chronological cobordisms} between closed 1-manifolds:
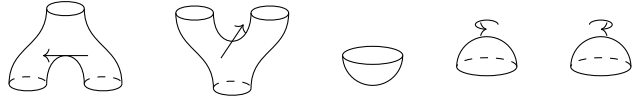
\[
\tikz[baseline=2.5ex, scale=.5]{
	\draw (0,0) .. controls (0,1) and (1,1) .. (1,2);
	\draw (1,0) .. controls (1,1) and (2,1) .. (2,0);
	\draw (3,0) .. controls (3,1) and (2,1) .. (2,2);
	\draw (0,0) .. controls (0,-.25) and (1,-.25) .. (1,0);
	\draw[dashed] (0,0) .. controls (0,.25) and (1,.25) .. (1,0);
	\draw (2,0) .. controls (2,-.25) and (3,-.25) .. (3,0);
	\draw[dashed] (2,0) .. controls (2,.25) and (3,.25) .. (3,0);
	\draw (1,2) .. controls (1,1.75) and (2,1.75) .. (2,2);
	\draw (1,2) .. controls (1,2.25) and (2,2.25) .. (2,2);
        \draw[<-] (0.9,0.75) -- (2.1,0.75);
}
\qquad
\tikz[baseline=9.5ex, scale=.5]{
	\draw  (1,2) .. controls (1,3) and (0,3) .. (0,4);
	\draw  (2,2) .. controls (2,3) and (3,3) .. (3,4);
	\draw (1,4) .. controls (1,3) and (2,3) .. (2,4);
	\draw (0,4) .. controls (0,3.75) and (1,3.75) .. (1,4);
	\draw (0,4) .. controls (0,4.25) and (1,4.25) .. (1,4);
	\draw (2,4) .. controls (2,3.75) and (3,3.75) .. (3,4);
	\draw (2,4) .. controls (2,4.25) and (3,4.25) .. (3,4);
	\draw (1,2) .. controls (1,1.75) and (2,1.75) .. (2,2);
	\draw[dashed] (1,2) .. controls (1,2.25) and (2,2.25) .. (2,2);
        \draw[<-] (1.8,3.7) -- (1.2,2.8);
}
\qquad 
\tikz[baseline=0ex, scale=.2]{
    \draw [domain=180:360] plot ({2*cos(\x)}, {2*sin(\x)});
    \draw (-2,0) arc (180:360:2 and 0.6);
    \draw (2,0) arc (0:180:2 and 0.6);
}
\qquad
\tikz[baseline=1ex, scale=.2]{
        \draw[domain=0:180] plot ({2*cos(\x)}, {2*sin(\x)});
        \draw (-2,0) arc (180:360:2 and 0.6);
        \draw[dashed] (2,0) arc (0:180:2 and 0.6);
        \draw[->] (0,2.8) [partial ellipse=0:270:5ex and 2ex];
}
\qquad
\tikz[baseline=1ex, scale=.2]{
        \draw[domain=0:180] plot ({2*cos(\x)}, {2*sin(\x)});
        \draw (-2,0) arc (180:360:2 and 0.6);
        \draw[dashed] (2,0) arc (0:180:2 and 0.6);
        \draw[<-] (0,2.8) [partial ellipse=-90:180:5ex and 2ex];
}
\]
with an additional twisting (transposing) identity cobordism. Together, these observations imply that we may decompose all changes of chronology into sequences of \textit{elementary changes of chronologies}. These are exactly those pairs of cobordisms described in the commutation chart (Figure 2) of \cite{Ozsv_th_2013}.

\begin{figure}
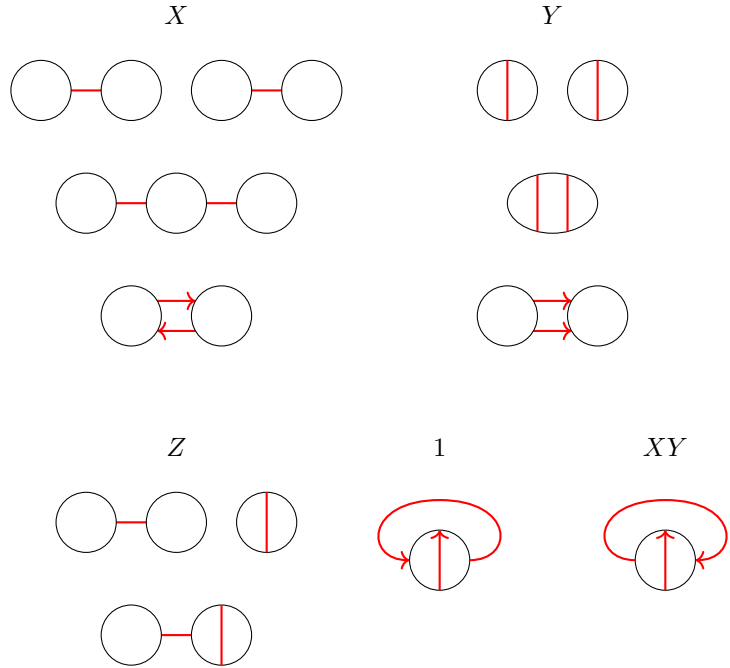

\[
\tikz{
    \node at (-3,1) {$\tikz{
  	\node at (0,0) {$X$};
  	\draw[red, thick] (-1.8,-1) -- (-.6,-1);
  	\draw[red, thick] (.6,-1) -- (1.8,-1);
  	\draw[fill=white] (-1.8,-1) circle (.4cm);
  	\draw[fill=white] (-.6,-1) circle (.4cm);
  	\draw[fill=white] (.6,-1) circle (.4cm);
  	\draw[fill=white] (1.8,-1) circle (.4cm);
  	\draw[red, thick] (-1.2,-2.5) -- (1.2,-2.5);
  	\draw[fill=white] (-1.2,-2.5) circle (.4cm);
  	\draw[fill=white] (0,-2.5) circle (.4cm);
  	\draw[fill=white] (1.2,-2.5) circle (.4cm);
  	\draw[red, thick,->] (-.6,-3.8) -- (.254,-3.8);
  	\draw[red, thick,->] (.6,-4.2) -- (-.254,-4.2);
  	\draw[fill=white] (-.6,-4) circle (.4cm);
  	\draw[fill=white] (.6,-4) circle (.4cm);
}$};
    \node at (2,1) {$\tikz{
  	\node at (0,0) {$Y$};
  	\draw[fill=white] (-.6,-1) circle (.4cm);
  	\draw[fill=white] (.6,-1) circle (.4cm);
  	\draw[red, thick] (-.6,-1.4) -- (-.6,-.6);
  	\draw[red, thick] (.6,-1.4) -- (.6,-.6);
  	\draw[fill=white] (0,-2.5) ellipse (.6cm and .4cm);
	\begin{scope}
	   	\clip (0,-2.5) ellipse (.6cm and .4cm);
  		\draw[red, thick] (-.2,-2.9) -- (-.2,-2.1);
  		\draw[red, thick] (.2,-2.9) -- (.2,-2.1);
	\end{scope}
  	\draw[red, thick,->] (-.6,-3.8) -- (.254,-3.8);
  	\draw[red, thick,->] (-.6,-4.2) -- (.254,-4.2);
  	\draw[fill=white] (-.6,-4) circle (.4cm);
  	\draw[fill=white] (.6,-4) circle (.4cm);
}$};
    \node at (-3, -4) {$\tikz{
  	\node at (0,0) {$Z$};
  	\draw[red, thick] (-1.2,-1) -- (0,-1);
  	\draw[fill=white] (-1.2,-1) circle (.4cm);
  	\draw[fill=white] (0,-1) circle (.4cm);
  	\draw[fill=white] (1.2,-1) circle (.4cm);
  	\draw[red, thick] (1.2,-1.4) -- (1.2,-.6);
  	\draw[red, thick] (-.6,-2.5) -- (.6,-2.5);
  	\draw[fill=white] (-.6,-2.5) circle (.4cm);
  	\draw[fill=white] (.6,-2.5) circle (.4cm);
  	\draw[red, thick] (.6,-2.9) -- (.6,-2.1);
  }$};
    \node at (2,-4) {$\tikz{
  	\node at (-1.5,0) {$1$};
  	\draw[fill=white] (-1.5,-1.5) circle (.4cm);
  	\draw[red, thick,->] (-1.5,-1.9) -- (-1.5,-1.1);
  	\draw[red, thick,->] (-1.1,-1.5) .. controls (-.5,-1.5) and (-.5,-.7) .. (-1.5,-.7)
  		.. controls (-2.5,-.7) and (-2.5,-1.5) .. (-1.9,-1.5);
  	\node at (1.5,0) {$XY$};
  	\draw[fill=white] (1.5,-1.5) circle (.4cm);
  	\draw[red, thick,->] (1.5,-1.9) -- (1.5,-1.1);
  	\draw[red, thick,->] (1.1,-1.5) .. controls (.5,-1.5) and (.5,-.7) .. (1.5,-.7)
  		.. controls (2.5,-.7) and (2.5,-1.5) .. (1.9,-1.5);
  	\draw[opacity=0] (.6,-2.5) circle (.4cm);
  }$};
}
\]
\caption{This is the collection of elementary changes of chronologies, together with their evaluation by $\iota$. Notice that taking $X=Z=1$ and $Y=-1$ yields the commutation chart of~\cite{Ozsv_th_2013}. Framings are omitted if evaluation by $\iota$ does not depend on them.}
\label{elementaryCoC}
\end{figure}

At this point, we have defined a 2-category whose objects are closed 1-manifolds, with chronological cobordisms as 1-morphisms and changes of chronology as 2-morphisms. This 2-category is simplified by the following procedure: for $R = \mathbb{Z}[X,Y,Z^{\pm1}]/(X^2 = Y^2 = 1)$, define the map $\iota$ which assigns to each elementary change of chronology a monomial, as pictured in Figure \ref{elementaryCoC}.\footnote{For those elementary cobordisms $H$ with $\iota(H) = Z$, it is assumed that $H$ takes a merge followed by a split to a split followed by a merge. If the opposite is true, $\iota(H) = Z^{-1}$.} Indeed,
\[
\iota(H' \circ H) = \iota(H') \iota(H) \qquad \text{and} \qquad \iota(H' \star H) = \iota(H') \iota(H)
\]
so $\iota$ assigns to every change of chronology a monomial in $R$; for more on the map $\iota$ (e.g., well-definedness and multiplicativity), see~\cite{putyra20152categorychronologicalcobordismsodd}.

Finally, as in the even case, we will eventually allow chronological cobordisms to be decorated by finitely many dots as long as each dot never shares the same level set as another dot or critical point. Precisely, let $C$ denote the critical points of $\tau$ and $D$ denote the dots on $W$. Both are taken to be finite. Then, a \textit{dotted chronological cobordism} is a chronological cobordism for which $\tau(x) \not= \tau(y)$ whenever $x,y \in C \cup D$ are distinct. In~\cite{putyra20152categorychronologicalcobordismsodd}, Putyra shows that if $H$ is a change of chronology which does nothing but move one dot past another with respect to the Morse function, then $\iota(H) = XY$.

A subtle but important distinction of the setup is the degree; define the \textit{$\mathbb{Z} \times \mathbb{Z}$-degree} of a cobordism $W$ by
\[
\abs{W} = (\#\mathrm{births} - \#\mathrm{merges} - \#\mathrm{dots}, \#\mathrm{deaths} - \#\mathrm{splits} - \#\mathrm{dots}).
\]
Note that the sum of the entries of $\abs{W}$ is the topological degree $\det_t(W)$ from \S \ref{evenCat}. Moreover, define $\lambda:(\mathbb{Z} \times \mathbb{Z})^2 \to R$ to be the bilinear map given by
\[
\lambda((x_1, y_1), (x_2, y_2)) = X^{x_1x_2} Y^{y_1y_2}Z^{x_1y_2 - y_1x_2}.
\]
Suppose $H$ is a change of chronology moving two cobordisms $W$ and $W'$ past one another; e.g., $H$ looks like
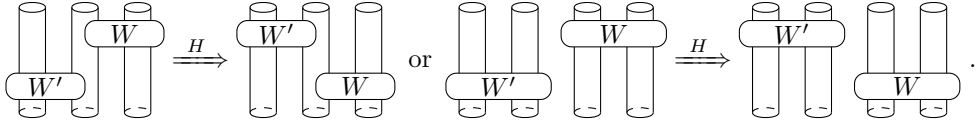
\[
\tikz[baseline=4ex, scale=.35]{
	\draw (0,0) .. controls (0,.-.25) and (1,-.25) .. (1,0);
	\draw[dashed] (0,0) .. controls (0,.25) and (1,.25) .. (1,0);
	\draw (0,0) -- (0,4);
	\draw (1,0) -- (1,4);
	\draw (0,4) .. controls (0,3.75) and (1,3.75) .. (1,4);
	\draw (0,4) .. controls (0,4.25) and (1,4.25) .. (1,4);
	\draw (2,0) .. controls (2,.-.25) and (3,-.25) .. (3,0);
	\draw[dashed] (2,0) .. controls (2,.25) and (3,.25) .. (3,0);
	\draw (2,0) -- (2,4);
	\draw (3,0) -- (3,4);
	\draw (2,4) .. controls (2,3.75) and (3,3.75) .. (3,4);
	\draw (2,4) .. controls (2,4.25) and (3,4.25) .. (3,4);
	\draw (4,0) .. controls (4,.-.25) and (5,-.25) .. (5,0);
	\draw[dashed] (4,0) .. controls (4,.25) and (5,.25) .. (5,0);
	\draw (4,0) -- (4,4);
	\draw (5,0) -- (5,4);
	\draw (4,4) .. controls (4,3.75) and (5,3.75) .. (5,4);
	\draw (4,4) .. controls (4,4.25) and (5,4.25) .. (5,4);
	\filldraw [fill=white, draw=black,rounded corners] (-.5,.5) rectangle (2.5,1.5) node[midway] { $W'$};
	\filldraw [fill=white, draw=black,rounded corners] (2.5,2.5) rectangle (5.5,3.5) node[midway] { $W$};
} 
\xRightarrow{~H~}
\tikz[baseline=4ex, scale=.35]{
	\draw (0,0) .. controls (0,.-.25) and (1,-.25) .. (1,0);
	\draw[dashed] (0,0) .. controls (0,.25) and (1,.25) .. (1,0);
	\draw (0,0) -- (0,4);
	\draw (1,0) -- (1,4);
	\draw (0,4) .. controls (0,3.75) and (1,3.75) .. (1,4);
	\draw (0,4) .. controls (0,4.25) and (1,4.25) .. (1,4);
	\draw (2,0) .. controls (2,.-.25) and (3,-.25) .. (3,0);
	\draw[dashed] (2,0) .. controls (2,.25) and (3,.25) .. (3,0);
	\draw (2,0) -- (2,4);
	\draw (3,0) -- (3,4);
	\draw (2,4) .. controls (2,3.75) and (3,3.75) .. (3,4);
	\draw (2,4) .. controls (2,4.25) and (3,4.25) .. (3,4);
	\draw (4,0) .. controls (4,.-.25) and (5,-.25) .. (5,0);
	\draw[dashed] (4,0) .. controls (4,.25) and (5,.25) .. (5,0);
	\draw (4,0) -- (4,4);
	\draw (5,0) -- (5,4);
	\draw (4,4) .. controls (4,3.75) and (5,3.75) .. (5,4);
	\draw (4,4) .. controls (4,4.25) and (5,4.25) .. (5,4);
	\filldraw [fill=white, draw=black,rounded corners] (-.5,2.5) rectangle (2.5,3.5) node[midway] { $W'$};
	\filldraw [fill=white, draw=black,rounded corners] (2.5,.5) rectangle (5.5,1.5) node[midway] { $W$};
}
\enspace \mathrm{or} \enspace
\tikz[baseline=4ex, scale=.35]{
	\draw (0,0) .. controls (0,.-.25) and (1,-.25) .. (1,0);
	\draw[dashed] (0,0) .. controls (0,.25) and (1,.25) .. (1,0);
	\draw (0,0) -- (0,4);
	\draw (1,0) -- (1,4);
	\draw (0,4) .. controls (0,3.75) and (1,3.75) .. (1,4);
	\draw (0,4) .. controls (0,4.25) and (1,4.25) .. (1,4);
	\draw (2,0) .. controls (2,.-.25) and (3,-.25) .. (3,0);
	\draw[dashed] (2,0) .. controls (2,.25) and (3,.25) .. (3,0);
	\draw (2,0) -- (2,4);
	\draw (3,0) -- (3,4);
	\draw (2,4) .. controls (2,3.75) and (3,3.75) .. (3,4);
	\draw (2,4) .. controls (2,4.25) and (3,4.25) .. (3,4);
	\filldraw [fill=white, draw=black,rounded corners] (-.5,.5) rectangle (3.5,1.5) node[midway] { $W'$};
} \ 
\tikz[baseline=4ex, scale=.35]{
	\draw (0,0) .. controls (0,.-.25) and (1,-.25) .. (1,0);
	\draw[dashed] (0,0) .. controls (0,.25) and (1,.25) .. (1,0);
	\draw (0,0) -- (0,4);
	\draw (1,0) -- (1,4);
	\draw (0,4) .. controls (0,3.75) and (1,3.75) .. (1,4);
	\draw (0,4) .. controls (0,4.25) and (1,4.25) .. (1,4);
	\draw (2,0) .. controls (2,.-.25) and (3,-.25) .. (3,0);
	\draw[dashed] (2,0) .. controls (2,.25) and (3,.25) .. (3,0);
	\draw (2,0) -- (2,4);
	\draw (3,0) -- (3,4);
	\draw (2,4) .. controls (2,3.75) and (3,3.75) .. (3,4);
	\draw (2,4) .. controls (2,4.25) and (3,4.25) .. (3,4);
	\filldraw [fill=white, draw=black,rounded corners] (-.5,2.5) rectangle (3.5,3.5) node[midway] { $W$};
}
\xRightarrow{~H~}
\tikz[baseline=4ex, scale=.35]{
	\draw (0,0) .. controls (0,.-.25) and (1,-.25) .. (1,0);
	\draw[dashed] (0,0) .. controls (0,.25) and (1,.25) .. (1,0);
	\draw (0,0) -- (0,4);
	\draw (1,0) -- (1,4);
	\draw (0,4) .. controls (0,3.75) and (1,3.75) .. (1,4);
	\draw (0,4) .. controls (0,4.25) and (1,4.25) .. (1,4);
	\draw (2,0) .. controls (2,.-.25) and (3,-.25) .. (3,0);
	\draw[dashed] (2,0) .. controls (2,.25) and (3,.25) .. (3,0);
	\draw (2,0) -- (2,4);
	\draw (3,0) -- (3,4);
	\draw (2,4) .. controls (2,3.75) and (3,3.75) .. (3,4);
	\draw (2,4) .. controls (2,4.25) and (3,4.25) .. (3,4);
	\filldraw [fill=white, draw=black,rounded corners] (-.5,2.5) rectangle (3.5,3.5) node[midway] { $W'$};
} \ \tikz[baseline=4ex, scale=.35]{
	\draw (0,0) .. controls (0,.-.25) and (1,-.25) .. (1,0);
	\draw[dashed] (0,0) .. controls (0,.25) and (1,.25) .. (1,0);
	\draw (0,0) -- (0,4);
	\draw (1,0) -- (1,4);
	\draw (0,4) .. controls (0,3.75) and (1,3.75) .. (1,4);
	\draw (0,4) .. controls (0,4.25) and (1,4.25) .. (1,4);
	\draw (2,0) .. controls (2,.-.25) and (3,-.25) .. (3,0);
	\draw[dashed] (2,0) .. controls (2,.25) and (3,.25) .. (3,0);
	\draw (2,0) -- (2,4);
	\draw (3,0) -- (3,4);
	\draw (2,4) .. controls (2,3.75) and (3,3.75) .. (3,4);
	\draw (2,4) .. controls (2,4.25) and (3,4.25) .. (3,4);
	\filldraw [fill=white, draw=black,rounded corners] (-.5,.5) rectangle (3.5,1.5) node[midway] { $W$};
}~.
\]
Then,
\[
\iota(H) = \lambda\left(\abs{W}, \abs{W'}\right).
\]
Note that this agrees with and generalizes the statement about changes of chronologies which move dots past one another. Putyra also provides the following, extremely helpful change of framing relations.
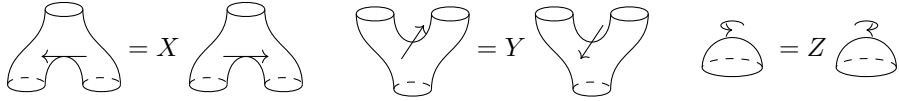
\[
\tikz[baseline=2.5ex, scale=.5]{
	\draw (0,0) .. controls (0,1) and (1,1) .. (1,2);
	\draw (1,0) .. controls (1,1) and (2,1) .. (2,0);
	\draw (3,0) .. controls (3,1) and (2,1) .. (2,2);
	\draw (0,0) .. controls (0,-.25) and (1,-.25) .. (1,0);
	\draw[dashed] (0,0) .. controls (0,.25) and (1,.25) .. (1,0);
	\draw (2,0) .. controls (2,-.25) and (3,-.25) .. (3,0);
	\draw[dashed] (2,0) .. controls (2,.25) and (3,.25) .. (3,0);
	\draw (1,2) .. controls (1,1.75) and (2,1.75) .. (2,2);
	\draw (1,2) .. controls (1,2.25) and (2,2.25) .. (2,2);
        \draw[<-] (0.9,0.75) -- (2.1,0.75);
}
 = X~
 \tikz[baseline=2.5ex, scale=.5]{
	\draw (0,0) .. controls (0,1) and (1,1) .. (1,2);
	\draw (1,0) .. controls (1,1) and (2,1) .. (2,0);
	\draw (3,0) .. controls (3,1) and (2,1) .. (2,2);
	\draw (0,0) .. controls (0,-.25) and (1,-.25) .. (1,0);
	\draw[dashed] (0,0) .. controls (0,.25) and (1,.25) .. (1,0);
	\draw (2,0) .. controls (2,-.25) and (3,-.25) .. (3,0);
	\draw[dashed] (2,0) .. controls (2,.25) and (3,.25) .. (3,0);
	\draw (1,2) .. controls (1,1.75) and (2,1.75) .. (2,2);
	\draw (1,2) .. controls (1,2.25) and (2,2.25) .. (2,2);
        \draw[->] (0.9,0.75) -- (2.1,0.75);
}
\qquad
\tikz[baseline=9.5ex, scale=.5]{
	\draw  (1,2) .. controls (1,3) and (0,3) .. (0,4);
	\draw  (2,2) .. controls (2,3) and (3,3) .. (3,4);
	\draw (1,4) .. controls (1,3) and (2,3) .. (2,4);
	\draw (0,4) .. controls (0,3.75) and (1,3.75) .. (1,4);
	\draw (0,4) .. controls (0,4.25) and (1,4.25) .. (1,4);
	\draw (2,4) .. controls (2,3.75) and (3,3.75) .. (3,4);
	\draw (2,4) .. controls (2,4.25) and (3,4.25) .. (3,4);
	\draw (1,2) .. controls (1,1.75) and (2,1.75) .. (2,2);
	\draw[dashed] (1,2) .. controls (1,2.25) and (2,2.25) .. (2,2);
        \draw[<-] (1.8,3.7) -- (1.2,2.8);
}
= Y~
\tikz[baseline=9.5ex, scale=.5]{
	\draw  (1,2) .. controls (1,3) and (0,3) .. (0,4);
	\draw  (2,2) .. controls (2,3) and (3,3) .. (3,4);
	\draw (1,4) .. controls (1,3) and (2,3) .. (2,4);
	\draw (0,4) .. controls (0,3.75) and (1,3.75) .. (1,4);
	\draw (0,4) .. controls (0,4.25) and (1,4.25) .. (1,4);
	\draw (2,4) .. controls (2,3.75) and (3,3.75) .. (3,4);
	\draw (2,4) .. controls (2,4.25) and (3,4.25) .. (3,4);
	\draw (1,2) .. controls (1,1.75) and (2,1.75) .. (2,2);
	\draw[dashed] (1,2) .. controls (1,2.25) and (2,2.25) .. (2,2);
        \draw[->] (1.8,3.7) -- (1.2,2.8);
}
\qquad
\tikz[baseline=1ex, scale=.2]{
        \draw[domain=0:180] plot ({2*cos(\x)}, {2*sin(\x)});
        \draw (-2,0) arc (180:360:2 and 0.6);
        \draw[dashed] (2,0) arc (0:180:2 and 0.6);
        \draw[->] (0,2.8) [partial ellipse=0:270:5ex and 2ex];
}
~= Z~
\tikz[baseline=1ex, scale=.2]{
        \draw[domain=0:180] plot ({2*cos(\x)}, {2*sin(\x)});
        \draw (-2,0) arc (180:360:2 and 0.6);
        \draw[dashed] (2,0) arc (0:180:2 and 0.6);
        \draw[<-] (0,2.8) [partial ellipse=-90:180:5ex and 2ex];
}
\]

In summary, we let $\mathrm{ChCob}_\bullet(0)$ (or just $\mathrm{ChCob}_\bullet$) denote the graded monoidal category whose
\begin{itemize}
    \item objects are formally $\mathbb{Z} \times \mathbb{Z}$-graded closed 1-manifolds (i.e., a pair of a closed 1-manifold and an element of $\mathbb{Z}\times\mathbb{Z}$),
    \item $\mathrm{Hom}((x_1, y_1) A, (x_2, y_2) B)$ is the free $\mathbb{Z}$-module spanned by isotopy classes of (dotted) chonological cobordisms $W$ from $A$ to $B$ with degree
    \[
    \abs{W} = (x_1 - x_2, y_1 - y_2),
    \]
    modulo the change of framing and \textit{change of chronology relations}: $W' = \iota(H) W$ for each change of chronology $H: W\Rightarrow W'$.
\end{itemize}

\subsection{Unified arc algebras}
\label{ss:unified arc algebras}

In this section, we consider the unified arc algebras $H^n$ over $R$, as provided by \cite{naisse2017odd} and \cite{naisse2020odd} (there, referred to as ``covering'' arc algebras). This is done in spirit of \cite{Khovanov_2002}, as in \S \ref{sss:Khovanov arc algebra}, using the ``chronological TQFT'' provided in \cite{putyra20152categorychronologicalcobordismsodd}. There are a number of challenges presented by this construction: for example, the unified arc algebras are non-associative, and the composition map $\mu[t,s]$ do not preserve $\mathbb{Z} \times \mathbb{Z}$-degree. The fix, as provided in \cite{naisse2020odd}, is to use the structure of a grading category, described in \S \ref{ss:brief graded outline}.

We must be a bit more careful when setting up the unified arc algebras. Still, for $a \in B^m$, $b \in B_n$, and $t\in B_m^n$, $a t b$ is a closed 1-manifold; for $s \in B_n^p$ and $c \in B_p$, we can still define a cobordism
\[
(atb)(\overline{b}sc) \to a(ts)c
\]
but we specify a chronology when we do so. The cobordism is still obtained by contracting symmetric arcs of $b\overline{b}$, and we fix the chronology by taking saddles from right-to-left and choosing the ``upwards'' framing. This is the chronological cobordism denoted by $W_{abc}(t,s)$.

Next, define the ``chronological'' TQFT $\mathcal{F}: \mathrm{ChCob}_\bullet \to R\mathrm{Mod}$. Set
\[
\mathcal{F}(\underbrace{\bigcirc \sqcup \cdots \sqcup \bigcirc}_{n}) = V^{\otimes n}.
\]
where $V = R\langle v_+, v_- \rangle$ and, now in this $\mathbb{Z}\times \mathbb{Z}$-graded landscape, we set
\[
\deg_R(v_+) = (1,0) \qquad \text{and} \qquad \deg_R(v_-) = (0,-1).
\]
so that $\deg_R(u) = (\#v_+(u), -\#v_-(u))$ where $v_\pm(u)$ denotes the collection of copies of $v_\pm$ appearing in $u$. Finally, on elementary cobordisms, set
\begin{align*}
\mathcal{F}\left(\tikz[baseline=2.5ex, scale=.5]{
	\draw (0,0) .. controls (0,1) and (1,1) .. (1,2);
	\draw (1,0) .. controls (1,1) and (2,1) .. (2,0);
	\draw (3,0) .. controls (3,1) and (2,1) .. (2,2);
	\draw (0,0) .. controls (0,-.25) and (1,-.25) .. (1,0);
	\draw[dashed] (0,0) .. controls (0,.25) and (1,.25) .. (1,0);
	\draw (2,0) .. controls (2,-.25) and (3,-.25) .. (3,0);
	\draw[dashed] (2,0) .. controls (2,.25) and (3,.25) .. (3,0);
	\draw (1,2) .. controls (1,1.75) and (2,1.75) .. (2,2);
	\draw (1,2) .. controls (1,2.25) and (2,2.25) .. (2,2);
        \draw[<-] (0.9,0.75) -- (2.1,0.75);
}\right) : V \otimes V \rightarrow V &= 
\begin{cases}
v_+ \otimes v_+ \mapsto v_+, & v_+ \otimes v_- \mapsto v_-, \\
v_- \otimes v_- \mapsto 0, & v_- \otimes v_+ \mapsto XZ v_-,~\mathrm{and}
\end{cases}
\\
\mathcal{F}\left(\tikz[baseline=9.5ex, scale=.5]{
	\draw  (1,2) .. controls (1,3) and (0,3) .. (0,4);
	\draw  (2,2) .. controls (2,3) and (3,3) .. (3,4);
	\draw (1,4) .. controls (1,3) and (2,3) .. (2,4);
	\draw (0,4) .. controls (0,3.75) and (1,3.75) .. (1,4);
	\draw (0,4) .. controls (0,4.25) and (1,4.25) .. (1,4);
	\draw (2,4) .. controls (2,3.75) and (3,3.75) .. (3,4);
	\draw (2,4) .. controls (2,4.25) and (3,4.25) .. (3,4);
	\draw (1,2) .. controls (1,1.75) and (2,1.75) .. (2,2);
	\draw[dashed] (1,2) .. controls (1,2.25) and (2,2.25) .. (2,2);
        \draw[<-] (1.8,3.7) -- (1.2,2.8);
}\right) : V  \rightarrow V \otimes V &= 
\begin{cases}
v_+ \mapsto v_- \otimes v_+ + YZ v_+ \otimes v_-,  &\\
v_- \mapsto v_- \otimes v_- &
\end{cases}
\\
\mathcal{F}\left(\tikz[baseline=5.2ex, scale=.5]{
	\draw (1,2) .. controls (1,1) and (2,1) .. (2,2);
	\draw (1,2) .. controls (1,1.75) and (2,1.75) .. (2,2);
	\draw (1,2) .. controls (1,2.25) and (2,2.25) .. (2,2);
}\right) : R  \rightarrow V  &= 
\begin{cases}
1 \mapsto v_+, & 
\end{cases}
\\
\mathcal{F}\left(\tikz[baseline=.2ex, scale=.5]{
	\draw (1,0) .. controls (1,1) and (2,1) .. (2,0);
	\draw (1,0) .. controls (1,-.25) and (2,-.25) .. (2,0);
	\draw[dashed] (1,0) .. controls (1,.25) and (2,.25) .. (2,0);
        \draw[->] (1.5,1.1) [partial ellipse=0:270:3ex and 1ex];
}\right) : V  \rightarrow R &= 
\begin{cases}
v_+ \mapsto 0, & \\
v_- \mapsto 1. &
\end{cases}
\end{align*}
To obtain a complete description, we apply the change of framing local relations. Now, notice that a cylinder with a hole evaluates to either
\[
\begin{cases}
v_+ \mapsto Z(X+Y) v_- \\ v_- \mapsto 0
\end{cases}
\qquad \text{or} \qquad
\begin{cases}
v_+ \mapsto Z (XY + 1) v_- \\ v_- \mapsto 0
\end{cases}
\]
depending on the framing. Therefore, unfortunately, we are not able to think of dots as 1/2 of a hole anymore; we define $\mathcal{F}$ on dots by setting
\[
\mathcal{F}\left(\tikz[baseline=2.5ex, scale=.5]{
	\draw (1,0) .. controls (1,.-.25) and (2,-.25) .. (2,0);
	\draw[dashed] (1,0) .. controls (1,.25) and (2,.25) .. (2,0);
	\draw (1,0) -- (1,2);
	\draw (2,0) -- (2,2);
	\draw (1,2) .. controls (1,1.75) and (2,1.75) .. (2,2);
	\draw (1,2) .. controls (1,2.25) and (2,2.25) .. (2,2);
	\node at (1.5,1) {$\bullet$};
}\right) : V \rightarrow V = 
\begin{cases}
v_+ \mapsto v_-, & \\
v_- \mapsto 0. &
\end{cases}
\]
as before. Again, it is easy to check that $\mathcal{F}$ observes the sphere and tube cutting relations.

Finally, for $t\in B_m^n$, the \textit{unified arc space} is defined
\[
\mathcal{F}(t) = \bigoplus_{a\in B^m, b\in B_n} \mathcal{F}(atb).
\]
Given another tangle $s\in B_n^p$, define the composition map
\begin{align*}
&\mu[t,s]:\mathcal{F}(atb) \otimes \mathcal{F}(b'sc) \to \mathcal{F}(a (ts) d)\\ &\text{by}~\mu[t,s] = \begin{cases} 0 & \text{if}~\overline{b}\not=b' \\ \mathcal{F}(W_{abc}(t,s)) & \text{if}~\overline{b}=b' \end{cases}
\end{align*}
where $b'\in B^n$ and $c\in B_p$. Note that, as promised, $\mu[t,s]$ \textit{does not} preserve $\mathbb{Z} \times \mathbb{Z}$-degree. 

\begin{definition}
The \textit{unified arc algebra}, which we still denote $H^n$, is the unified arc space
\[
H^n = \mathcal{F}(1_n) = \bigoplus_{a\in B^m, b\in B_m} \mathcal{F}(a1_nb)
\]
with multiplication $\mu[1_n,1_n]$.
\end{definition}

\subsection{A brief outline of \texorpdfstring{$\mathcal{C}$}{Lg}-graded structures}
\label{ss:brief graded outline}

In this subsection, we review the motivation for and construction of $\mathcal{G}$-graded $R$-modules given in \cite{naisse2020odd}. In the following sections, we provide a thorough description of a slight generalization of the procedure introduced here.

It has been shown (cf. \cite{naisse2017odd} Proposition 3.2) that the multiplication as defined above is not associative in the unified arc algebra. This presents the main difficulty---in \cite{Khovanov_2002}, Khovanov provides that
\[
\mathcal{F}(t) \otimes_{H^n} \mathcal{F}(s) \cong \mathcal{F}(ts)
\]
declaring that $(u' \cdot h) \otimes u = u' \otimes (h \cdot u)$. The assumption that multiplication in $H^n$ is associative is implicit here.

On the other hand, the failure of associativity is controlled by the cobordisms involved. Explicitly, observe the square
\[
\begin{tikzcd}[cramped, column sep=tiny]
& \begin{matrix} \mu_{abc}[t,t'](x,y) \in \mathcal{F}_c(\overline{a}tt'c) \\ z \in \mathcal{F}_c(\overline{c}t''d) \end{matrix} \ar[dr] &  \\
\begin{matrix} x \in \mathcal{F}_c(\overline{a}tb) \\ y \in \mathcal{F}_c(\overline{b}t'c) \\ z\in \mathcal{F}_c(\overline{c}t''d) \end{matrix} \ar[ur, ""] \ar[dr, ""'] & & \begin{matrix} \mu_{acd}[tt',t''](\mu_{abc}[t,t'](x,y), z) \\ \mu_{abd}[t,t't''](x, \mu_{bcd}[t',t''](y,z)) \end{matrix} \in \mathcal{F}_c(\overline{a}tt't'' d)\\
& \begin{matrix} x \in \mathcal{F}_c(\overline{a}tb) \\ \mu_{bcd}[t',t''](y,z) \in \mathcal{F}_c(\overline{b}t't''d) \end{matrix} \ar[ur] &
\end{tikzcd}
\]
In general,
\[
\mu_{acd}[tt', t''] \circ (\mu_{abc}[t,t'] \otimes 1_z) \not= \mu_{abd}[t,t't''] \circ (1_x \otimes \mu_{bcd}[t',t'']),
\]
but the failure is witnessed by the degree of the cobordisms involved: $W_{acd}(tt', t'')$ and $W_{abc}(t,t')$, and $W_{abd}(t,t't'')$ and $W_{bcd}(t',t'')$. The degree of elements also have effect.

In the literature, Majid and Albuquerque \cite{albuquerque1998quasialgebra} show that the octonions $\mathbb{O}$, while non-assoicative, admit a grading by the group $(\mathbb{Z}/2\mathbb{Z})^3$, and the gradings witness the failure of associativity. That is, they show that $\mathbb{O}$ is quasi-associative; in general, a $G$-graded $\mathbb{K}$-algebra $A$ is called \textit{quasi-associative} (or \textit{graded associative}) if there is a 3-cocycle $\alpha: G^{[3]} \to \mathbb{K}^\times$ for which
\[
a \cdot (b \cdot c) = \alpha\left(\abs{a},\abs{b},\abs{c}\right) ~ (a \cdot b) \cdot c
\]
for all homogeneous elements $a, b, c \in A$ (here, $\abs{~\cdot~}: A \to G$ is the grading).

Naisse and Putyra \cite{naisse2020odd} generalize the notion of quasi-associativity. Remarking that the 3-cocycle condition is exactly the pentagon relation for a monoidal category, their first goal is to provide similar definitions for modules and algebras \textit{graded by categories}.
\begin{definition}
By a \textit{grading category}, we will mean a category $\mathcal{C}$ endowed with a 3-cocycle $\alpha: \mathcal{C}^{[3]} \to \mathbb{K}^\times$, referred to as the \textit{associator}. Then, a \textit{$\mathcal{C}$-graded $\mathbb{K}$-module} is a $\mathbb{K}$-module $M$ which admits a decomposition
\[
M = \bigoplus_{g \in \mathcal{C}} M_g.
\]
\end{definition}
By ``$g\in \mathcal{C}$'', we just mean that $g$ is any morphism of $\mathcal{C}$. This generalizes gradings by a group by delooping: we can view any group $G$ as a category with a single object $\bullet$ with $\mathrm{End}(\bullet) = G$. A \textit{$\mathcal{C}$-graded map} $f:M\to N$ between $\mathcal{C}$-graded modules is just one which preserves grading: $f(M_g) \subset N_g$.

Define the category $\mathrm{Mod}^\mathcal{C}$ of $\mathcal{C}$-graded $\mathbb{K}$-modules with morphisms being graded maps. It is a monoidal category where the decomposision of $M' \otimes M = \bigoplus_{g \in \mathcal{C}} (M' \otimes M)_g$ is given by
\[
(M' \otimes M)_g = \bigoplus_{g = g_2 \circ g_1} M'_{g_2} \otimes_\mathbb{K} M_{g_1}
\]
for composable $g_1$ and $g_2$ (revealing a slightly different feature of the $\mathcal{C}$-graded setting). The coherence isomorphism is then given by the associator: 
\begin{align*}
(M_3 \otimes M_2) \otimes M_1 &\xrightarrow{\alpha} M_3 \otimes (M_2 \otimes M_1) \\ 
(z\otimes y) \otimes x & \mapsto \alpha(\abs{z}, \abs{y}, \abs{x})~ z \otimes (y \otimes x)
\end{align*}
for homogeneous elements $x$, $y$, and $z$. The $\mathcal{C}$-graded $\mathbb{K}$-module $\bigoplus_{X\in \mathrm{Ob}(\mathcal{C})} \mathbb{K}_{\mathrm{Id}_X}$ is the unit object, and the unitors for this tensor product may also be defined via the associator. We will describe this process explicitly in slightly more generality later on.

With this language, Naisse and Putyra are able to define $\mathcal{C}$-graded algebras and bimodules as well. First, a $\mathcal{C}$-graded $\mathbb{K}$-algebra $A$ is a $\mathcal{C}$-graded $\mathbb{K}$-module with a graded associative multiplication map $A \otimes A \to A$ such that $A_{g} \cdot A_{g'} \subset A_{g' \circ g}$, where $A_{g' \circ g} = \{0\}$ whenever $g' \circ g$ is undefined. Similarly, for two $\mathcal{C}$-graded algebras $A_1$ and $A_2$, a $\mathcal{C}$-graded $A_2$--$A_1$-bimodule $M$ is a $\mathcal{C}$-graded module $M$ with graded, $\mathbb{K}$-linear left and right actions $A_2\otimes M \to M$ and $M\otimes A_1 \to M$ satisfying the usual bimodule conditions, twisted by the associator: for example, these actions respect
\[
(y \cdot m) \cdot x = \alpha\left(\abs{y},\abs{m},\abs{x}\right) ~ y \cdot (m\cdot x)
\]
for all $y\in A_2$, $m\in M$ and $x\in A_1$. In this subsection, we will denote the category of $\mathcal{C}$-graded $A_2$--$A_1$-bimodules by $\mathrm{Bimod}^\mathcal{C}(A_2, A_1)$. The morphisms of this category will be graded maps between $A_2$--$A_1$-bimodules which preserve the left and right actions.

We employ the associator to see that, given $M'\in \mathrm{Bimod}^\mathcal{C}(A_3, A_2)$ and $M\in \mathrm{Bimod}^\mathcal{C}(A_2, A_1)$, $M'\otimes_{\mathbb{K}} M \in \mathrm{Bimod}^\mathcal{C}(A_3, A_1)$: the left and right actions are the horizontal maps making the following diagrams commute.
\[
\begin{tikzcd}[column sep = -3ex]
A_3 \otimes (M' \otimes M) \ar{rr} \ar[swap]{rd}{\alpha^{-1}} && M' \otimes M \\
&(A_3 \otimes M') \otimes M \ar[swap]{ur}{}&
\end{tikzcd}
\quad
\begin{tikzcd}[column sep = -3ex]
(M' \otimes M) \otimes A_1 \ar{rr} \ar[swap]{rd}{\alpha} && M' \otimes M \\
& M' \otimes (M \otimes A_1) \ar[swap]{ur}{}&
\end{tikzcd}
\]
Then, we can define the tensor product over the intermediary algebra $A_2$ via the coequalizer: explicitly, 
\[
M'\otimes_{A_2} M = M'\otimes_\mathbb{K} M \Big/ \left((m'\cdot x)\otimes m - \alpha\left(\abs{m'}, \abs{x}, \abs{m}\right)~ m'\otimes (x\cdot m)\right)
\]
with left $A_3$- and right $A_1$-actions induced by the ones on $M'\otimes_\mathbb{K} M$.

Now, with the goal of showing that the unified arc algebra $H^n$ is graded associative, we must build a suitable grading category $(\mathcal{G}, \alpha)$. Let $B^\bullet = \bigsqcup_{n \ge 0}B^n$ denote the collection of all crossingless matchings. Given a flat tangle $t$, we write $\widehat{t}$ or $t^\wedge$ to mean the tangle $t$ with all free loops removed; $\widehat{B}_m^n$ denotes the collection of planar tangles with no free loops. Let $\mathcal{G}$ denote the category where
\begin{itemize}
    \item $\mathrm{Ob}(\mathcal{G}) = B^\bullet$, and whose
    \item morphisms are formally $\mathbb{Z} \times \mathbb{Z}$-graded planar tangles; that is, 
    \[
    \mathrm{Hom}_\mathcal{G}(a,b) = \widehat{B}_m^n\times \mathbb{Z}^2
    \]
    for any $a\in B^m$ and $b\in B^n$.
\end{itemize}
The composition, for $(t,p) \in \mathrm{Hom}_\mathcal{G}(a,b)$ and $(t', p') \in \mathrm{Hom}_\mathcal{G}(b,c)$, is defined
\[
(t', p') \circ (t,p) = (\widehat{tt'}, p+p'+ \abs{W_{abc}(t,t')}) \in \mathrm{Hom}_\mathcal{G}(a,c).
\]
Note that, since $W_{abc}(t,t')$ consists of only saddle moves,
\[
\abs{W_{abc}(t,t')} = (-\#\mathrm{merges~in}~W_{abc}(t,t'), -\#\mathrm{splits~in}~W_{abc}(t,t')).
\]
So, it follows that the identity morphism for any crossingless matching $a \in B^m$ is $\mathrm{Id}_a = (1_m, (m,0))$. Henceforth, to make life easier, given objects $a \in B^m$ and $b\in B^n$, we'll write $atb$ when, really, we mean $at\overline{b}$.

We will omit a description of the associator until defining our own in the generalized setting---it will be apparent how to specialize ours to the current situation. Instead, we describe the way in which way elements of $H^n$, or $\mathcal{F}(t)$ in general, are $\mathcal{G}$-graded. For $u\in \mathcal{F}(atb)$, we set
\[
\deg_\mathcal{G}(u) = (\widehat{t}, \deg_R(u)) \in \mathrm{Hom}_\mathcal{G}(a,b).
\]
Hopefully this explains the choice to remove free loops from tangles: they are not involved in composition maps between arc algebras, and are extraneous information in light of the second entry of the grading. 

Secondly, this presents a solution to the first problem for unified arc algebras: $\mu_{abc}[t,s]$ preserves the $\mathcal{G}$-grading. Suppose $u \in \mathcal{F}(atb)$ and $v \in \mathcal{F}(bsc)$, so $\deg_\mathcal{G}(u) = (t,\deg_R(u)) \in \mathrm{Hom}(a,b)$ and $\deg_\mathcal{G}(v) = (s, \deg_R(v)) \in \mathrm{Hom}(b,c)$. Their composition in unified arc spaces is given by the map $\mu_{abc}[t,s]$. Recall that in the definition of the chronological TQFT $\mathcal{F}$, each merge decreases the number of copies of $v_+$ by 1, and each split increases the number of copies of $v_-$ by 1; consequently
\[
\deg_\mathcal{G}(\mu_{abc}[t,s](u,v)) = \left(\widehat{ts}, \deg_R(u) + \deg_R(v) + \abs{W_{abc}(t,s)}\right) = \deg_\mathcal{G}(v) \circ \deg_\mathcal{G}(u)
\]
as desired.

Finally, we can prove that
\[
\mu_{acd}[tt', t'']\left(\mu_{abc}[t,t'](x,y), z\right) = \alpha\left(\abs{x}, \abs{y}, \abs{z}\right) \mu_{abd}[t, t't'']\left(x, \mu_{abd}[t',t''](y,z)\right)
\]
for any $x\in \mathcal{F}(atb)$, $y \in \mathcal{F}(bt'c)$ and $z\in \mathcal{F}(c t'' d)$. In particular, Naisse and Putyra provide the following (for a discussion on unitality, see \cite{naisse2020odd} Proposition 6.2).

\begin{proposition}
$H^n$ is a unital, associative, $\mathcal{G}$-graded $R$-algebra. 
\end{proposition}

It is routine to check that, for $t \in B_m^n$, $\mathcal{F}(t)$ is an $(H^m, H^n)$-bimodule: the left $H^m$-action is given by $\mu[1_m, t]$ and the right $H^n$-action is given by $\mu[t, 1_n]$. Naisse and Putyra then provide the desired properties of these bimodules, in the sense that it mirrors results of \cite{Khovanov_2002}.

\begin{proposition}
Let $t\in B_m^n$. Then $\mathcal{F}(t)$ is an $(H^m,H^n)$-bimodule. It is also sweet as an $(H^m,H^n)$-bimodule; that is, it is projective as a left $H^m$-module and as a right $H^n$-module. Moreover, given $s\in B_n^p$, there is an isomorphism 
\[
\mathcal{F}(t) \otimes_{H^n} \mathcal{F}(s) \cong \mathcal{F}(ts)
\]
induced by $\mu[t,s]:\mathcal{F}(t) \otimes_R \mathcal{F}(s) \to \mathcal{F}(ts)$.
\end{proposition}

\subsubsection{$\mathcal{G}$-shifting system}

So far, we have successfully defined the relevant algebraic objects in the $\mathcal{G}$-graded setting. However, we have glossed over the important discussion of graded maps. In particular, given $t,s\in B_m^n$, so that $\mathcal{F}(t), \mathcal{F}(s) \in \mathrm{Ob}\left(\mathrm{Bimod}^{\mathcal{G}}(H^m, H^n)\right)$, can we describe those relevant morphisms between $\mathcal{F}(t)$ and $\mathcal{F}(s)$ in this category? Of course, any cobordism $W:t\to s$ induces a map $\mathcal{F}(W): \mathcal{F}(t) \to \mathcal{F}(s)$, but this map is clearly not graded! There must be a fix if we are to interpret cubes of resolutions with this approach; in particular, the only graded map between $\mathcal{F}\left(\tikz[baseline={([yshift=-.5ex]current bounding box.center)}, scale=.5]
{
    \begin{scope}[rotate=90]
	\draw[dotted] (3,-2) circle(0.707);
	\draw (2.5,-1.5) .. controls (2.75,-1.75) and (3.25,-1.75) .. (3.5,-1.5);
	\draw (2.5,-2.5) .. controls  (2.75,-2.25) and (3.25,-2.25) .. (3.5,-2.5);
    \end{scope}
}\right)$ and $\mathcal{F}\left(\tikz[baseline={([yshift=-.5ex]current bounding box.center)}, scale=.5]
{
	\draw[dotted] (3,-2) circle(0.707);
	\draw (2.5,-1.5) .. controls (2.75,-1.75) and (3.25,-1.75) .. (3.5,-1.5);
	\draw (2.5,-2.5) .. controls  (2.75,-2.25) and (3.25,-2.25) .. (3.5,-2.5);
}\right)$ is the zero map. The solution of Naisse and Putyra is the introduction of \textit{grading shifting functors} via a \textit{$\mathcal{G}$-shifting system}. Here is the idea of a $\mathcal{C}$-shifting system; a more precise, expanded definition is given in Section \ref{S: SHIFTING SYSTEMS}.

\begin{definition}
A \textit{$\mathcal{C}$-shifting system} is a pair $(I,\Phi)$ consisting of a monoid $(I,\bullet, e)$ and a collection $\Phi = \{\varphi_i\}_{i\in I}$ of families of maps
\[
\varphi_i = \{\varphi_i^{X,Y}: \mathsf{D}_i^{X,Y} \to \mathrm{Hom}_\mathcal{C}(X,Y)\}_{X,Y\in \mathrm{Ob}(\mathcal{C})}
\]
for $\mathsf{D}_i^{X,Y} \subset \mathrm{Hom}_\mathcal{C}(X,Y)$. These families of maps $\varphi_i$ are called \textit{$\mathcal{C}$-grading shifts}, and they are required to satisfy the property that, for each $i,j\in I$ and $X, Y \in \mathrm{Ob}(\mathcal{C})$, the following diagram commutes.
\[
\begin{tikzcd}
\mathrm{Hom}_\mathcal{C}(Y,Z) \times \mathrm{Hom}_\mathcal{C}(X,Y) \arrow[r, "\circ"] \arrow[d, "\varphi_j \times \varphi_i"'] & \mathrm{Hom}_\mathcal{C}(X,Z) \arrow[d, "\varphi_{j\bullet i}"] \\
\mathrm{Hom}_\mathcal{C}(Y,Z) \times \mathrm{Hom}_\mathcal{C}(X,Y) \arrow[r, "\circ"] & \mathrm{Hom}_\mathcal{C}(X,Z)
\end{tikzcd}
\]
\end{definition}

It is not immediate that a $\mathcal{C}$-shifting system $(I,\Phi)$ is compatible with the associator $\alpha$; a major portion of \cite{naisse2020odd}, and now our work, has to do with this observation. 

If $S = (I, \{\varphi_i\}_{i\in I})$ is a $\mathcal{C}$-shifting system compatible with $\alpha$, then for each $i\in I$, $\varphi_i:\mathrm{Mod}^\mathcal{C} \to \mathrm{Mod}^\mathcal{C}$ is a functor, called the \textit{grading shift functor}, and is defined as follows. For $M  = \bigoplus_{g\in \mathcal{C}} M_g \in \mathrm{Ob}(\mathrm{Mod}^\mathcal{C})$, put
\[
\varphi_i(M) = \bigoplus_{g\in \mathsf{D}_i} \varphi_i(M)_{\varphi_i(g)}
\]
where $\varphi_i(M)_{\varphi_i(g)} = M_g$. In other words, this grading shift functor turns elements of degree $g \in \mathsf{D}_i$ into elements of degree $\varphi_i(g)$; elements whose degree is not in $\mathsf{D}_i$ are sent to zero.

We will see that the witnesses to compatibility between a given $\mathcal{C}$-shifting system and associator imply the existence of canonical isomorphisms
\[
\varphi_j(M') \otimes \varphi_i(M) \to \varphi_{j\bullet i}(M'\otimes M).
\]
Indeed, there is a natural transformation $\varphi_j(-) \otimes \varphi_i(-) \Rightarrow \varphi_{j\bullet i}(- \otimes -)$. From here, under a certain assumption, it is easy to define shifted bimodules. In summary, this is to say that the shifting functor $\varphi_i: \mathrm{Mod}^\mathcal{C} \to \mathrm{Mod}^\mathcal{C}$ further induces a shifting functor $\varphi_i: \mathrm{Bimod}^\mathcal{C}(A_2, A_1) \to \mathrm{Bimod}^\mathcal{C}(A_2, A_1)$. The shifting functor also respects tensor products: for $M'\in \mathrm{Bimod}^\mathcal{C}(A_3, A_2)$ and $M\in \mathrm{Bimod}^\mathcal{C}(A_2, A_1)$, 
\[
\varphi_j(M') \otimes_{A_2} \varphi_i(M) \cong \varphi_{j\bullet i} (M'\otimes_{A_2} M).
\]

Returning to the situation at hand, our goal is to define a $\mathcal{G}$-shifting system (compatible with $\alpha$). The $\mathcal{G}$-shifting system we will use is given simply by weighted cobordisms $(W, v)$ where $v\in \mathbb{Z}\times \mathbb{Z}$. Explicitly, to construct the monoid in this shifting system, recall that given two cobordisms $W_1: t \to t'$ for $t,t'\in B_m^n$ and $W_2: s\to s'$ for $s,s'\in B^n_\ell$, we obtain a cobordism $W_1\bullet W_2: ts \to t's'$ by horizontal stacking.

Now, given weighted cobordisms $(W_1, v_1)$ and $(W_2, v_2)$, define $(W_1, v_1) \bullet (W_2, v_2)$ to be $(W_1 \bullet W_2, v_1 + v_2)$ whenever $W_1\bullet W_2$ is defined, and zero otherwise. The monoid of the $\mathcal{G}$-shifting system will be the collection of weighted cobordisms together with formal identity absorbing elements $\{(W,v)\} \sqcup \{e,0\}$ under the operation $\bullet$. Finally, given $t,t'\in B_m^n$ and $(W:t\to t', v)$, given any $a\in B^m$ and $b\in B^n$ we define
\[
\varphi_{(W,v)}^{a,b}(\widehat{t},p) = (\widehat{t'}, p + v + \abs{1_{a} W 1_b})
\]
where $1_{a} W 1_{b}$ is the cobordism $W$ capped off by $a\times [0,1]$ on one side and $b \times [0,1]$ (really, $\overline{b} \times [0,1]$) on the other. Since $\mathrm{Ob}(\mathcal{G}) = B^\bullet$, we can write $\varphi_{(W,v)} = \left\{\varphi_{(W,v)}^{a,b}\right\}_{a\in B^m,b\in B^n}$; we will often abuse notation and write $\varphi_{(W,v)}$ when it does not present confusion. Clearly, the domain of $\varphi_{(W,v)}^{a,b}$ is simply $\mathsf{D}_{(W,v)}^{a,b} = \{(\widehat{t},p) \in \mathrm{Hom}_\mathcal{G}(a,b) : p\in \mathbb{Z}\times\mathbb{Z}\}$. We'll write $\varphi_W$ sometimes when $v$ can be left ambiguous; however, in computations, this notation means $v = (0,0)$. Finally, for a flat tangle $t$, let $\mathbbm{1}_t$ denote the identity cobordism on $t$. Consider the collection of identity cobordisms $\mathbbm{1} = \{\mathbbm{1}_t\}_t$. Then there is an identity shift functor given by $\varphi_{\mathrm{id}} = \bigoplus_{\mathbbm{1}}\varphi_{\mathbbm{1}_t}$.

In practice, it is beneficial to view weighted cobordisms $(W,v)$ as two separate shifts; the first on a given planar tangle and the second on the $\mathbb{Z}\times\mathbb{Z}$ degree associated to that tangle. Unfortunately, to determine compatibility maps one must choose an order: we will always shift first by the chronological cobordism $W$ and second by the $\mathbb{Z}\times\mathbb{Z}$-degree. The opposite choice can also be made, and leads to small differences in the theory---for example, see Proposition \ref{FLshift}. In this way, Naisse and Putyra show that this $\mathcal{G}$-shifting system is compatible with the associator defined above; for more details, see \cite{naisse2020odd}.

Of course, there is also the possibility of vertically composing cobordisms. This is to say that the $\mathcal{G}$-shifting system may be extended to a \textit{shifting 2-system} (again, defined by Naisse-Putyra). Explicitly, in the monoid defined above, we define vertical composition in the same spirit as horizontal composition: for $W_1: t \to t'$ and $W_2: s \to s'$, 
\[
(W_2, v_2) \circ (W_1, v_1) = \begin{cases} (W_2 \circ W_1, v_2 + v_1) & \text{if}~ t' = s \\ 0 & \text{otherwise}.\end{cases}
\]
Compatibility maps are constructed via the change of chronology
\[
H: (W_2'\circ W_2) \bullet (W_1' \circ W_1) \Rightarrow (W_2' \bullet W_1') \circ (W_2 \bullet W_1).
\]
With this structure in place, we will see that any cobordism with corners $W: t \to s$ induces a \textit{graded} map $\mathcal{F}(W): \varphi_W (\mathcal{F}(t)) \to \mathcal{F}(s)$, as desired.

\newpage

\section{Grading multicategories and planar arc diagrams}
\label{s:gradingmultis and pads}

In this section, we generalize the work of Naisse and Putyra's to provide a category compatible with ``multigluing''; i.e., a framework for replacing tangles $t$ by planar arc diagrams $D$. We note that the content of this section and the next will come as little surprise to readers familiar with \cite{naisse2020odd}, outside of complications and additional structure associated to dealing with multicategories.

We start by extending the definition of $\mathcal{F}$ to planar arc diagrams, defined momentarily. In \S \ref{ss:gradingmulticategories}, we review multicategories, define grading multicategories, and construct the grading multicategory $\mathscr{G}$ utilized throughout this paper. In \S \ref{ss:Gisgradingmulti}, we verify that $\mathscr{G}$ is indeed a grading multicategory. Then, \S \ref{ss:gradingcatgeneral} is dedicated to establishing some properties of modules graded by multicategories which we use extensively. We conclude with \S \ref{ss:Ggradedmultimodules}, wherein we list consequences of observations made in \S \ref{ss:gradingcatgeneral} for $\mathscr{G}$-graded multimodules associated to planar arc diagrams by $\mathcal{F}$.

\begin{definition}
\label{def:planararcdiagrams}
An $(m_1, \ldots, m_k; n)$-\textit{planar arc diagram} $D$ is a disk $D$ with $k$ interior disks removed, together with a proper embedding of disjoint circles and closed intervals, so that there are $2m_i$ endpoints on the boundary component corresponding to the $i$th removed disk, and $2n$ endpoints on the outer boundary of $D$. Note that planar arc diagram $D$ comes with an ordering on the removed inner disks. Each boundary component carries a basepoint, disjoint from the endpoints of intervals, denoted by $\times$. We say that $D$ is oriented if the embedded circles and intervals are oriented. Both oriented and unoriented planar arc diagrams are considered up to planar isotopy. The collection of planar arc diagrams of type $(m_1, \ldots, m_k; n)$ is denoted by $\mathscr{D}_{(m_1,\ldots, m_k; n)}$. Similarly $\widehat{\mathscr{D}}_{(m_1,\ldots, m_k; n)}$ is the collection of $(m_1, \ldots, m_k; n)$-planar arc diagrams with free loops removed.
\end{definition}

\begin{wrapfigure}{r}{0pt}
\begin{tikzpicture}[scale=1.2]
    \draw[dotted] (0,0) circle (2cm);
    \node at (2,0) {$\times$};
    \draw[knot, rounded corners] (-1.685, -1.07739) -- (-1.585, -0.97739) .. controls (-1.65, 0.35) .. 
 (-1.24, 1.38472) -- (-1.34, 1.48472);
    \draw[knot, rounded corners] (-0.82, 1.82417) -- (-0.78, 1.72417) -- (-0.65, 0.38) .. controls (-0.2, 0) .. (-0.15, -1) .. controls (-0.6, -1.2) .. (-0.71, -1.86973);
    \draw[knot, rounded corners=1.5em] (-0.15, -1) -- (0.35, -0.5) -- (0.35, -1.4) -- (-0.15, -1);
    \draw[knot, rounded corners] (0.32, -1.97423) -- (0.32, -1.77423) .. controls (1.1,-1) and (0,0.4) .. (0.2, 1.3);
    \draw[knot, rounded corners=0.5em] (0.2, 1.3) -- (0.45, 1.55) .. controls (0.85, 1.2) .. (1.6, 1) -- (1.2, 0.4);
    \draw[knot, rounded corners=2em] (1.2, 0.4) -- (1.1, -0.7) -- (1.65, -1.13027);
    \draw[dotted, fill=white] (-0.65, 0.38) circle (0.3cm);
    \node at (-0.65, 0.38) {$1$};
    \node at (-0.35, 0.38) {$\times$}; 
    \draw[dotted, fill=white] (-0.15, -1) circle (0.2cm);
    \node at (-0.15, -1) {$4$};
    \node at (0.05, -1) {$\times$};
    \draw[dotted, fill=white] (0.2, 1.3) circle (0.2cm);
    \node at (0.2, 1.3) {$3$};
    \node at (0.4, 1.3) {$\times$};
    \draw[dotted, fill=white] (1.2, 0.4) circle (0.5cm);
    \node at (1.2, 0.4) {$2$};
    \node at (1.7, 0.4) {$\times$};
    \draw[knot] (-0.97, -0.73) circle (0.18cm);
\end{tikzpicture}
\end{wrapfigure}
For example, pictured to the right is an oriented $(1,1,1,2;3)$-planar arc diagram. We can compose planar arc diagrams by filling the $i$th empty region of one planar arc diagram with a $(\cdots; m_i)$ planar arc diagram. That is, given planar arc diagrams $D_i$ of type $(\ell_{i1}, \ldots, \ell_{i\alpha_i}; m_i)$ for $i=1,\ldots, k$ and $D$ of type $(m_1,\ldots, m_k; n)$, we set
\[
D \circ (D_1,\ldots, D_k) = D(D_1,\ldots, D_k; \emptyset).
\]
There is also a pairwise composition
\[
D \circ_i D_i = D(\emptyset,\ldots, D_i,\ldots, \emptyset; \emptyset).
\]
To the author's knowledge, this notation was first introduced in \cite{lawson2022homotopy} (we will adapt this definition to diskular tangles in Section \ref{ss:tangleresolution}). Note that the two notions of composition are related by
\[
D \circ (D_1,\ldots, D_k) = (\cdots ((D \circ_k D_k) \circ_{k-1} D_{k-1}) \circ_{k-2} \cdots ) \circ_1 D_1
\]
If $E$ is a planar arc diagram with an interior boundary component with $2n$ endpoints, we'll write $D(D_1,\ldots, D_k; E)$ to denote the resulting planar arc diagram. Otherwise, we frequently drop the last $\emptyset$ from the notation. This composition does not depend on order. 

On one hand, it is clear that any crossingless matching $a\in B^n$ uniquely defines a planar arc diagram of type $(; n)$. We choose the association
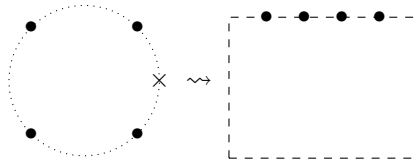
\[
\tikz[baseline={([yshift=-.55ex]current bounding box.center)}]{
    \draw[dotted] (0,0) circle (1cm);
    \node at (1,0) {$\times$};
    \node at (1.4142/2, 1.4142/2) {$\bullet$};
    \node at (-1.4142/2, 1.4142/2) {$\bullet$};
    \node at (-1.4142/2, -1.4142/2) {$\bullet$};
    \node at (1.4142/2, -1.4142/2) {$\bullet$};
}
\leadsto ~
\tikz[baseline={([yshift=-.55ex]current bounding box.center)}, scale=1.25]{
    \draw[dashed] (0,0) rectangle (2,1.5);
    \node at (.4, 1.5) {$\bullet$};
    \node at (.8, 1.5) {$\bullet$};
    \node at (1.2, 1.5) {$\bullet$};
    \node at (1.6, 1.5) {$\bullet$};
}
\]
so that, if we are being careful, the inner disks of a planar arc diagram can be filled with crossingless matchings belonging to $B^\bullet$ and can be closed on the outside by a crossingless matching belonging to $B_\bullet$.

\begin{figure}[ht]
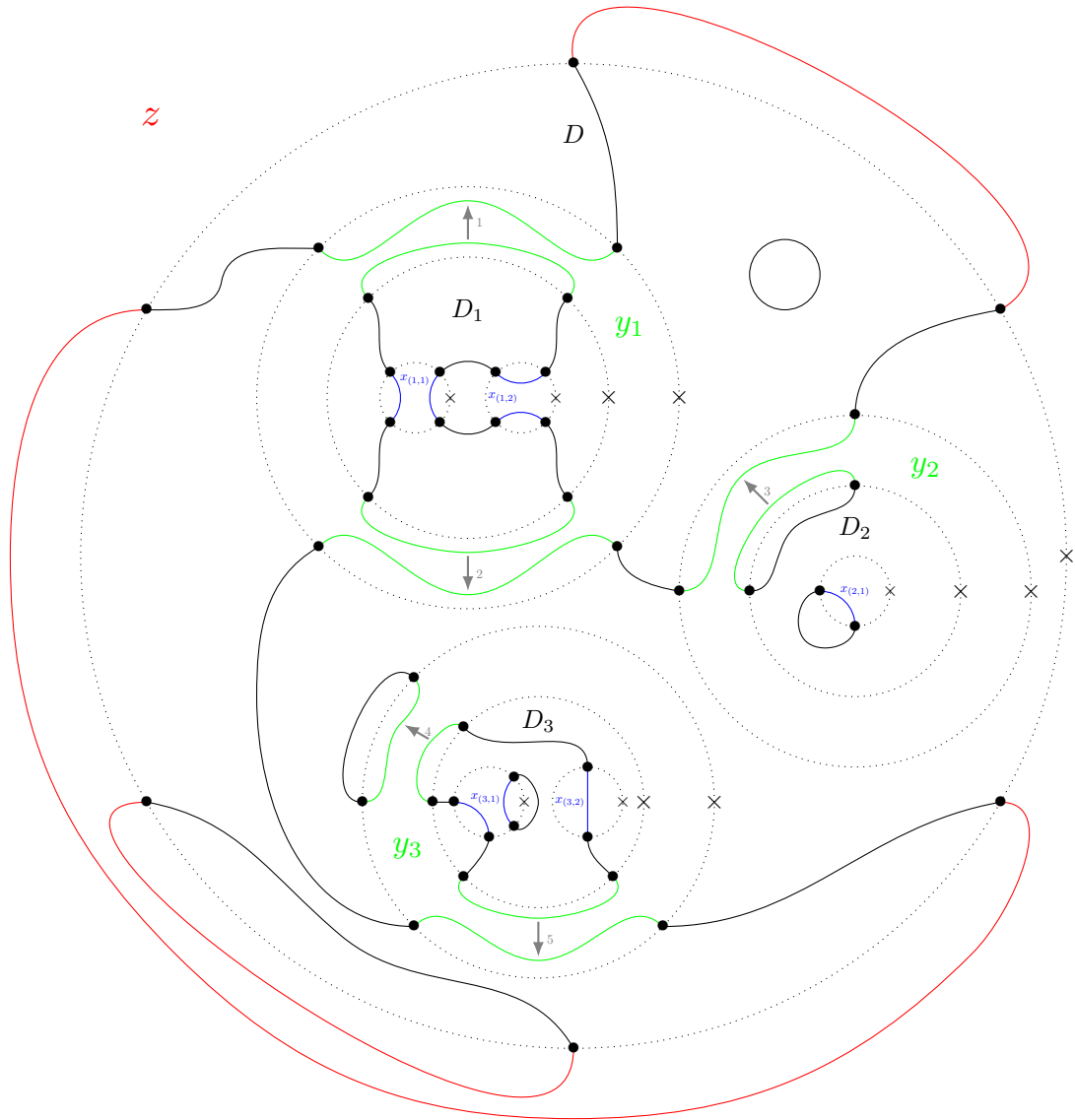

\[
\tikz[baseline={([yshift=-.5ex]current bounding box.center)}, scale=.95]{
    \draw[dotted] (0,0) circle (7cm);
        \node at (0,6) {$D$};
        \node at (7,0) {$\times$};
    \draw[dotted] (-1.5,2.25) circle (3cm);
        \node at (-1.5,3.5) {$D_1$};
        \node at (1.5,2.25) {$\times$};
    \draw[dotted] (-1.5,2.25) circle (2cm);
        \node at (0.5,2.25) {$\times$};
    \draw[dotted] (4,-.5) circle (2.5cm);
        \node at (4,.4) {$D_2$};
        \node at (5.5,-.5) {$\times$};
    \draw[dotted] (4,-.5) circle (1.5cm);
        \node at (6.5,-.5) {$\times$};
    \draw[dotted] (-.5, -3.5) circle (2.5cm);
        \node at (-.5, -2.35) {$D_3$};
        \node at (2, -3.5) {$\times$};
    \draw[dotted] (-.5, -3.5) circle (1.5cm);
        \node at (1, -3.5) {$\times$};
    \draw (-.5+1.4142/2*2.5, -3.5-1.4142/2*2.5) to[out=0,in=190] (1.732/2*7, -1/2*7);
    \draw (0,7) to[out=300, in=90] (-1.5+1.4142/2*3,2.25+1.4142/2*3);
    \draw (-1.732/2*7, 1/2*7) to[out=0,in=260] (-5,3.9) to[out=80,in=180] (-1.5-1.4142/2*3,2.25+1.4142/2*3);
    \draw (1.732/2*7, 1/2*7) to[out=190,in=90] (4,2);
    \draw (3,4) circle (.5cm);
    \draw (-1.5+1.4142/2*3,2.25-1.4142/2*3) to[out=270,in=170] (1.5,-.5);
    \draw (0,-7) to[out=120,in=320] (-3.3,-5.4) to[out=140,in=350] (-1.732/2*7, -1/2*7);
    \draw (-.5-1.4142/2*2.5, -3.5-1.4142/2*2.5) to[out=180,in=270] (-4.5, -2) to[out=90,in=210] (-1.5-1.4142/2*3,2.25-1.4142/2*3);
    \draw (-.5-1.4142/2*2.5, -3.5+1.4142/2*2.5) to[out=135,in=180] (-.5-2.5, -3.5);
    \draw[dotted] (-2.25,2.25) circle (.5cm);
    \node[scale=.75] at (-1.75,2.25) {$\times$};
    \draw[dotted] (-.75,2.25) circle (.5cm);
    \node[scale=.75] at (-.25,2.25) {$\times$};
    \draw (-2.25+1.4142/4,2.25+1.4142/4) to[out=45,in=135] (-.75-1.4142/4,2.25+1.4142/4);
    \draw (-2.25+1.4142/4,2.25-1.4142/4) to[out=315,in=225] (-.75-1.4142/4,2.25-1.4142/4);
    \draw (-.75+1.4142/4,2.25+1.4142/4) to[out=45,in=225] (-1.5+1.4142/2*2,2.25+1.4142/2*2);
    \draw (-.75+1.4142/4,2.25-1.4142/4) to[out=315,in=135] (-1.5+1.4142/2*2,2.25-1.4142/2*2);
    \draw (-2.25-1.4142/4,2.25+1.4142/4) to[out=135,in=315] (-1.5-1.4142/2*2,2.25+1.4142/2*2);
    \draw (-2.25-1.4142/4,2.25-1.4142/4) to[out=225,in=45] (-1.5-1.4142/2*2,2.25-1.4142/2*2);
    \draw[dotted] (4,-.5) circle (.5cm);
    \node[scale=.75] at (4.5,-.5) {$\times$};
    \draw (3.5,-.5) to[out=180,in=135] (4-1.4142/2,-.5-1.4142/2) to[out=315,in=270] (4,-1);
    \draw (4,1) to[out=270,in=45] (4 - 1.4142/2*1.25,-.5 + 1.4142/2*1.25) to[out=225,in=0] (2.5,-.5);
    \draw[dotted] (.2, -3.5) circle (.5cm);
    \node[scale=.75] at (.2+.5, -3.5) {$\times$};
    \draw (.2, -3.5-.5) to[out=270,in=135] (-.5+1.4142/2*1.5, -3.5-1.4142/2*1.5);
    \draw (.2, -3.5+.5) to[out=90,in=315] (-.5-1.4142/2*1.5, -3.5+1.4142/2*1.5);
    \draw[dotted] (-1.2, -3.5) circle (.5cm);
    \node[scale=.75] at (-1.2+.5, -3.5) {$\times$};
    \draw (-1.2+1.4142/4, -3.5+1.4142/4) to[out=45,in=90] (-.5,-3.5) to[out=270,in=315] (-1.2+1.4142/4, -3.5-1.4142/4);
    \draw (-1.2-.5, -3.5) -- (-.5-1.5, -3.5);
    \draw (-1.2, -3.5-.5) to[out=270,in=225] (-.5-1.4142/2*1.5, -3.5-1.4142/2*1.5);
    \draw[red] (1.732/2*7, 1/2*7) to[out=30, in=100] (0,7);
    \draw[red] (-1.732/2*7, 1/2*7) to[out=180, in=90] (-8,0) to[out=270,in=135] (-1.4142/2*8, -1.4142/2*8) to[out=315,in=180] (0,-8) to[out=0,in=225] (1.4142/2*8, -1.4142/2*8) to[out=45,in=0] (1.732/2*7, -1/2*7);
    \draw[red] (-1.732/2*7, -1/2*7) to[out=180,in=270] (0,-7);
    \node[red,scale=1.5] at (-6,6.25) {$z$};
    \draw[green] (-1.5+1.4142/2*2,2.25+1.4142/2*2) to[out=45,in=0] (-1.5,2.25+2.2) to[out=180,in=135] (-1.5-1.4142/2*2,2.25+1.4142/2*2);
    \draw[green] (-1.5+1.4142/2*3,2.25+1.4142/2*3) to[out=225,in=0] (-1.5,2.25+2.8) to[out=180,in=315] (-1.5-1.4142/2*3,2.25+1.4142/2*3);
    \draw[green] (-1.5-1.4142/2*2,2.25-1.4142/2*2) to[out=225,in=180] (-1.5,2.25-2.2) to[out=0,in=315] (-1.5+1.4142/2*2,2.25-1.4142/2*2);
    \draw[green] (-1.5-1.4142/2*3,2.25-1.4142/2*3) to[out=45,in=180] (-1.5,2.25-2.8) to[out=0,in=135] (-1.5+1.4142/2*3,2.25-1.4142/2*3);
    \node[green,scale=1.25] at (.8,3.25) {$y_1$};
    \draw[green] (4,2) to[out=270,in=45] (4 - 1.4142/2*2.3,-.5 + 1.4142/2*2.3) to[out=225,in=0] (1.5,-.5);
    \draw[green] (4,1) to[out=90,in=45] (4 - 1.4142/2*1.7,-.5 + 1.4142/2*1.7) to[out=225,in=180] (2.5,-.5);
    \node[green,scale=1.25] at (5,1.25) {$y_2$};
    \draw[green] (-.5-1.4142/2*2.5, -3.5+1.4142/2*2.5) to[out=315,in=45] (-.5-1.73205/2*2.25,-3.5+.5*2.25) to[out=225,in=0] (-.5-2.5, -3.5);
    \draw[green] (-.5-1.4142/2*1.5, -3.5+1.4142/2*1.5) to[out=135,in=45] (-.5-1.73205/2*1.75,-3.5+.5*1.75) to[out=225,in=180] (-.5-1.5, -3.5);
    \draw[green] (-.5-1.4142/2*1.5, -3.5-1.4142/2*1.5) to[out=225,in=180] (-.5, -3.5-1.65) to[out=0,in=315] (-.5+1.4142/2*1.5, -3.5-1.4142/2*1.5);
    \draw[green] (-.5-1.4142/2*2.5, -3.5-1.4142/2*2.5) to[out=45,in=180] (-.5, -3.5-2.25) to[out=0,in=135] (-.5+1.4142/2*2.5, -3.5-1.4142/2*2.5);
    \node[green,scale=1.25] at (-2.35,-4.15) {$y_3$};
    \draw[blue] (-2.25-1.4142/4,2.25+1.4142/4) to[out=315,in=45] (-2.25-1.4142/4,2.25-1.4142/4);
    \draw[blue] (-2.25+1.4142/4,2.25+1.4142/4) to[out=225,in=135] (-2.25+1.4142/4,2.25-1.4142/4);
    \node[blue,scale=.5] at (-2.25,2.5) {$x_{(1,1)}$};
    \draw[blue] (-.75-1.4142/4,2.25+1.4142/4) to[out=315,in=225] (-.75+1.4142/4,2.25+1.4142/4);
    \draw[blue] (-.75-1.4142/4,2.25-1.4142/4) to[out=45,in=135] (-.75+1.4142/4,2.25-1.4142/4);
    \node[blue,scale=.5] at (-1, 2.25) {$x_{(1,2)}$};
    \draw[blue] (3.5,-.5) to[out=0,in=90] (4,-1);
    \node[blue,scale=.5] at (4,-.5) {$x_{(2,1)}$};
    \draw[blue] (.2, -3.5-.5) -- (.2, -3.5+.5);
    \draw[blue] (-1.2-.5, -3.5) to[out=0,in=90] (-1.2, -3.5-.5);
    \draw[blue] (-1.2+1.4142/4, -3.5-1.4142/4) to[out=135,in=225] (-1.2+1.4142/4, -3.5+1.4142/4);
    \node[blue,scale=.5] at (-1.25, -3.45) {$x_{(3,1)}$};
    \node[blue,scale=.5] at (-.05, -3.5) {$x_{(3,2)}$};
    \node at (1.732/2*7, 1/2*7) {$\bullet$};
    \node at (0,7) {$\bullet$};
    \node at (-1.732/2*7, 1/2*7) {$\bullet$};
    \node at (-1.732/2*7, -1/2*7) {$\bullet$};
    \node at (0,-7) {$\bullet$};
    \node at (1.732/2*7, -1/2*7) {$\bullet$};
    \node at  (-1.5+1.4142/2*3,2.25+1.4142/2*3) {$\bullet$};
    \node at  (-1.5+1.4142/2*2,2.25+1.4142/2*2) {$\bullet$};
    \node at  (-1.5-1.4142/2*3,2.25+1.4142/2*3) {$\bullet$};
    \node at  (-1.5-1.4142/2*2,2.25+1.4142/2*2) {$\bullet$};
    \node at  (-1.5-1.4142/2*3,2.25-1.4142/2*3) {$\bullet$};
    \node at  (-1.5-1.4142/2*2,2.25-1.4142/2*2) {$\bullet$};
    \node at  (-1.5+1.4142/2*3,2.25-1.4142/2*3) {$\bullet$};
    \node at  (-1.5+1.4142/2*2,2.25-1.4142/2*2) {$\bullet$};
    \node at (4,2) {$\bullet$};
    \node at (4,1) {$\bullet$};
    \node at (1.5,-.5) {$\bullet$};
    \node at (2.5,-.5) {$\bullet$};
    \node at (-.5-1.4142/2*2.5, -3.5+1.4142/2*2.5) {$\bullet$};
    \node at (-.5-2.5, -3.5) {$\bullet$};
    \node at (-.5-1.4142/2*2.5, -3.5-1.4142/2*2.5) {$\bullet$};
    \node at (-.5-1.4142/2*1.5, -3.5+1.4142/2*1.5) {$\bullet$};
    \node at (-.5-1.5, -3.5) {$\bullet$};
    \node at (-.5-1.4142/2*1.5, -3.5-1.4142/2*1.5) {$\bullet$};
    \node at (-.5+1.4142/2*2.5, -3.5-1.4142/2*2.5) {$\bullet$};
    \node at (-.5+1.4142/2*1.5, -3.5-1.4142/2*1.5) {$\bullet$};
    \node at (-2.25+1.4142/4,2.25+1.4142/4) {$\bullet$};
    \node at (-2.25-1.4142/4,2.25+1.4142/4) {$\bullet$};
    \node at (-2.25-1.4142/4,2.25-1.4142/4) {$\bullet$};
    \node at (-2.25+1.4142/4,2.25-1.4142/4) {$\bullet$};
    \node at (-.75+1.4142/4,2.25+1.4142/4) {$\bullet$};
    \node at (-.75-1.4142/4,2.25+1.4142/4) {$\bullet$};
    \node at (-.75-1.4142/4,2.25-1.4142/4) {$\bullet$};
    \node at (-.75+1.4142/4,2.25-1.4142/4) {$\bullet$};
    \node at (3.5,-.5) {$\bullet$};
    \node at (4,-1) {$\bullet$};
    \node at (.2, -3.5-.5) {$\bullet$};
    \node at (.2, -3.5+.5) {$\bullet$};
    \node at (-1.2-.5, -3.5) {$\bullet$};
    \node at (-1.2, -3.5-.5) {$\bullet$};
    \node at (-1.2+1.4142/4, -3.5+1.4142/4) {$\bullet$};
    \node at (-1.2+1.4142/4, -3.5-1.4142/4) {$\bullet$};
    \draw[gray, thick] [-{Latex[length=2mm]}] (-1.5,4.5) -- node[scale=.5,right=1mm] {1} (-1.5,5);
    \draw[gray, thick] [-{Latex[length=2mm]}] (-1.5,0) -- node[scale=.5,right=1mm] {2} (-1.5,-.5);
    \draw[gray, thick] [-{Latex[length=2mm]}] (4 - 1.4142/2*1.75,-.5 + 1.4142/2*1.75) -- node[scale=.5,right=1mm] {3} (4 - 1.4142/2*2.25,-.5 + 1.4142/2*2.25);
    \draw[gray, thick] [-{Latex[length=2mm]}] (-.5-1.73205/2*1.8,-3.5+.5*1.8) -- node[scale=.5,right=1mm] {4} (-.5-1.73205/2*2.2,-3.5+.5*2.2);
    \draw[gray, thick] [-{Latex[length=2mm]}] (-.5, -3.5-1.7) -- node[scale=.5,right=1mm] {5} (-.5, -3.5-2.2);
}
\]
\caption{An example of a chronological coboridm $W_{\vec{x}\vec{y}z}((D_1, D_2, D_3), D)$.}
\label{fig:compositionmapsschematic}
\end{figure}

Thus, if $D$ is a $(m_1,\ldots, m_k; n)$ planar arc diagram, we define
\[
\mathcal{F}(D) = \bigoplus_{\substack{x_i\in B^{m_i}: i=1,\ldots, k \\ y\in B_n}} \mathcal{F}(D(x_1,\ldots, x_k; y))
\]
where $\mathcal{F}$ is the unified chronological TQFT. It is a $(H^{m_1}\otimes \cdots \otimes H^{m_k}, H^n)$-bimodule by the compositions
\[
\mu[(1_{m_1}, \ldots, 1_{m_k}); D] \qquad \text{and} \qquad \mu[D; 1_n].
\]
These composition maps are defined just as before: for compatible $D_i$, we define
\[
\mu[(D_1,\ldots, D_k); D]: \bigotimes_{i=1}^k \mathcal{F}(D_i) \otimes \mathcal{F}(D) \to \mathcal{F}(D(D_1\ldots, D_k))
\]
component-wise, as follows. For the time being, all tensor products are taken over $R$. Working with planar arc diagrams necessitates some burdensome notation. Notice that potentially far more closures are necessary: each $D_i$ requires, say, $\alpha_i$-many inner closures which we denote by $x_{(i,1)},\ldots, x_{(i,\alpha_i)}$, and one outer closure $y_i$. On the other hand, $D$ requires $k$ inner closures $y_1',\ldots, y_k'$ and one outer closure $z$. Let $\vec{x}$ denote the entire collection of crossingless matchings $\{x_{(1,1)}, \ldots, x_{(k,\alpha_k)}\}$. In the future, $\vec{~\cdot~}$ will always denote the entire collection of crossingless parings of that label. If $\vec{\cdot}$ has a subscript $i$, we mean all corssingless parings of that label whose first entry of their subscript is $i$; e.g., $\vec{x}_i = \{x_{(i,1)},\ldots, x_{(i,\alpha_i)}\}$. With this notation in place, we define $\mu[(D_1,\ldots, D_k); D]$ component-wise by
\[
\mu_{\vec{x}\vec{y} z}[(D_1,\ldots, D_k); D]: \bigotimes_{i=1}^k \mathcal{F}(D_i(\vec{x}_i; y_i)) \otimes \mathcal{F}(D(\vec{y'}; z)) \to \mathcal{F}(D(D_1(\vec{x}_1), \ldots, D_k(\vec{x}_k);z)
\]
where we can interpret $D_i(\vec{x}_i) := D_i(\vec{x}_i;\emptyset)$ as a crossingless matching, and
\[
\mu_{\vec{x}\vec{y} z}[(D_1,\ldots, D_k); D] = \begin{cases} 0 & \text{if}~y_i\not= y_i'~\text{for some}~i; \\ \mathcal{F}(W_{\vec{x}\vec{y}z}((D_1,\ldots, D_k); D)) & \text{if}~y_i = y_i'~\text{for all}~i.\end{cases}
\]
Elements of $\left( \bigotimes_{i=1}^k \mathcal{F}(D_i)\right) \otimes \mathcal{F}(D)$ are written $(u_1,\ldots, u_k) \otimes u$ or, frequently, $\vec{u} \otimes u$. The last thing we must do is describe the chronological cobordism $W_{\vec{x}\vec{y}z}((D_1,\ldots, D_k), D)$. This cobordism is (as one would expect, comparing to Sections \ref{sss:Khovanov arc algebra} and \ref{ss:unified arc algebras}) defined by contracting the symmetric arcs of $y_i$. The chronology is chosen by moving counter-clockwise from the basepoint of the $i$th removed disk of $D$ and contracting symmetric arcs outwardly, starting at $i=1$ and progressing to $i=k$. Use Figure \ref{fig:compositionmapsschematic} for reference. 
In this example, $W_{\vec{x} \vec{y} z} ((D_1, D_2, D_3), D)$ is the chronological cobordism obtained by contracting the symmetric arcs of $\vec{y}$ as specified by the gray arrows in the numbered order. So, it is a merge, followed by a split, and then three more merges. Notice that $W_{\vec{x} \vec{y} z} ((D_1,\ldots, D_n), D)$ has Euler characteristic $-\sum_i\abs{y_i}$ (recall that $\abs{y}=c$ whenever $y\in B^c$). As we proceed, we will use the notation  $\vec{y} D z$ to mean $D(\vec{y};z)$. This seems redundant, but it is especially helpful to write $\vec{x} (D_1,\ldots, D_k) \vec{y}$, or even $\vec{x} \vec{D}' \vec{y}$ for $\vec{D}' = (D_1, \ldots, D_k)$, rather than $(D_1(\vec{x}_1;y_1), D_2(\vec{x}_2;y_2), \ldots, D_k(\vec{x}_k; y_k))$.

Let $\vec{D}' = (D_1,\ldots, D_k)$. We will frequently refer to the chronolonological cobordism above via the (upwardly oriented) schematic 
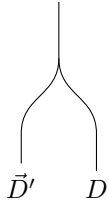
\[
\tikz{
	\draw (0,-1.9) node[below]{$\vec{D}'$}
        --
        (0,-1.5)
		.. controls (0,-1) and (.5,-1) ..
		(.5,-.5)  
		--
		(.5,0)
		--
		(.5,.25);
	\draw (1,-2) node[below]{$D$} 
		--
		(1,-1.5)
		.. controls (1,-1) and (.5,-1) ..
		(.5,-.5);
}
\]
where the trivalent vertex represents the cobordism $W_{\vec{x},\vec{y}, z}((D_1,\ldots, D_k), D)$. In following sections, we'll have to consider the compositions of such cobordisms, but it is not immediately clear how the chronology is defined. Let $\vec{D}'' = (D_{(1,1)}, \ldots, D_{(1, \alpha_1)}, \ldots, D_{(k,\alpha_k)})$, so that $\vec{D}_i''$ are the planar arc diagrams filling $D_i$. While we can interpret 
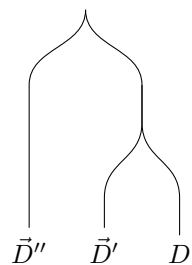
\[
\tikz[xscale=1]{
	\draw (-.5,.1) node[below]{$\vec{D}''$}
        --
        (-.5,2)
		.. controls (-.5,2.5) and (.25,2.5) ..
		(.25,3);
	\draw (.5,.1) node[below]{$\vec{D}'$}
        --
        (.5,.5)
		.. controls (.5,1) and (1,1) ..
		(1,1.5)
		--
		(1,2)
		.. controls (1,2.5) and (.25,2.5) ..
		(.25,3);
	\draw (1.5,0) node[below]{$D$}
		--
		(1.5,.5)
		.. controls (1.5,1) and (1,1) ..
		(1,1.5)
		--
		(1,2);
}
\]
as a chronological cobordism using our rules above, we'd like to consider compositions of the form
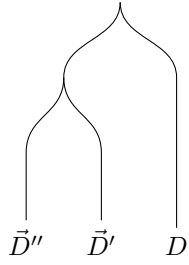
\[
\tikz[xscale=1]{
	\draw (0,.1) node[below]{$\vec{D}''$}
        --
        (0,1)
		.. controls (0,1.5) and (.5,1.5) ..
		(.5,2)
		.. controls (.5,2.5) and (1.25,2.5) ..
		(1.25,3); 
	\draw (1,.1) node[below]{$\vec{D}'$}
        --
        (1,.5)
		--
		(1,1)
		.. controls (1,1.5) and (.5,1.5) ..
		(.5,2);
	\draw (2, 0) node[below]{$D$}
		--
		(2,2)
		.. controls (2,2.5) and (1.25,2.5) ..
		(1.25,3);
}
\]
as well. In the latter, notice that the leftmost trivalent vertex is a \textit{collection} of chronological cobordisms, $W_{\vec{w}_i\vec{x}_iy_i}((D_{(i,1)},\ldots, D_{(i, \alpha_i)}))$. So, we will define the order of these chronological cobordisms to follow the index $i=1,\ldots, k$---the same idea applies to larger compositions. Denote the composition of these chronological cobordisms by $W_{\vec{w},\vec{x}, \vec{y}}(\vec{D}'', \vec{D}')$, and the corresponding map as $\mu[\vec{D}'', \vec{D}']$. To be explicit, 
\[
\mu[\vec{D}'', \vec{D}']: \left(\bigotimes_{i=1}^k \bigotimes_{j=1}^{\alpha_i} \mathcal{F}(D_{ij})\right) \otimes \left(\bigotimes_{i=1}^k \mathcal{F}(D_i)\right) \to \bigotimes_{i=1}^k \mathcal{F}(D_i(D_{i1}, \ldots, D_{i\alpha_i}))
\]
interpreting $\mu[\vec{D}'', \vec{D}'] = \bigotimes_{i=1}^k \mu[(D_{i1}, \ldots, D_{i\alpha_i}), D_i]$. We shorten the expression above to 
\[
\mu[\vec{D}'', \vec{D}']: \mathcal{F}(\vec{D}'') \otimes \mathcal{F}(\vec{D}') \to \mathcal{F}(\vec{D}'(\vec{D}'')).
\]

Finally, while we will almost always use the composition maps $\mu_{\vec{x}, \vec{y}, z}(\vec{D}, D)$ moving forward, we note that the flexibility of planar arc diagrams allows for a few more composition maps. First, note that one may fill the $i$th hole of $D$ by $D_i$, leaving the other holes unchanged, by considering the composition $\mu[(1_{n_1},\ldots, D_i, \ldots, 1_{n_k}); D]$.  On the other hand, we could also define a composition map which only fills one hole of $D$ without reference to the others. Consider the map 
\[
\mu[D_i; D]: \mathcal{F}(D_i) \otimes \mathcal{F}(D) \to \mathcal{F}(D(\emptyset,\ldots, D_i, \ldots, \emptyset))
\]
defined componentwise as
\begin{align*}
\mu_{(y_1', \ldots, \vec{x}_i, \ldots, y_k'), y_i, z}[D_i; D]: &\mathcal{F}(D_i(\vec{x}_i; y_i)) \otimes \mathcal{F}(D(\vec{y}'; z)) \to \mathcal{F}(D(y_1',\ldots, D_i(\vec{x}_i), \ldots, y_k'))\\
\mu_{(y_1', \ldots, \vec{x}_i, \ldots, y_k'), y_i, z}[D_i; D]&= \begin{cases}0 & \text{if}~ y_i' \not= y_i \\ \mathcal{F}(W_{(y_1', \ldots, \vec{x}_i, \ldots, y_k'), y_i, z}(D_i; D) ) & \text{if}~ y_i' = y_i\end{cases}
\end{align*}
where $W_{(y_1', \ldots, \vec{x}_i, \ldots, y_k'), y_i, z}(D_i; D)$ is the chronological cobordism which simply contracts symmetric arcs of $y_i\overline{y_i}$ counter-clockwise with respect to the basepoint, with closures specified by the other indices. Then, notice that
\begin{align*}
\mu[(D_1,\ldots, D_k); D] = \mu[D_k; D(D_1,\ldots, D_{k-1}, \emptyset)]  \circ\cdots & \circ \left(\mu[D_2; D(D_1, \emptyset, \ldots, \emptyset)] \otimes \mathrm{Id}_{D_3} \otimes \cdots \otimes \mathrm{Id}_{D_k}\right) \\ & \circ \left(\mu[D_1; D] \otimes \mathrm{Id}_{D_2} \otimes \cdots \otimes \mathrm{Id}_{D_k}\right)
\end{align*}
where $\mathrm{Id}_{D_i}$ means the identity on elements living in components corresponding to closures of $D_i$.

\subsection{(Grading) multicategories}
\label{ss:gradingmulticategories}

Recall that a \textit{(small) multicategory} $\mathscr{C}$ consists of
\begin{enumerate}
    \item a set of objects $\mathrm{Ob}(\mathscr{C})$,
    \item for each $k \ge 0$ and objects $x_1,\ldots, x_k, y\in \mathrm{Ob}(\mathscr{C})$, a set $\mathrm{Hom}(x_1,\ldots, x_k; y)$ of \textit{multimorphisms} from $(x_1,\ldots, x_k)$ to $y$,
    \item a composition map 
    \begin{align*}
        \mathrm{Hom}(y_1,\ldots, y_k; z) \times  \left(\mathrm{Hom}(x_{11}, \ldots, x_{1\alpha_1}; x_1) \times \cdots \times \mathrm{Hom}(x_{k1}, \ldots, x_{k \alpha_k}; x_k)\right) & \to \mathrm{Hom}(x_{11}, \ldots, x_{k\alpha_k}; z),
    \end{align*}
    and
    \item a distinguished element $\mathrm{Id}_x \in \mathrm{Hom}(x;x)$ for each $x\in \mathrm{Ob}(x)$ called the \textit{identity} of $x$
\end{enumerate}
defined so that composition is associative, in the sense that the following diagram commutes:
\[
      \xymatrix{
        {\begin{array}{l}\mathrm{Hom}(y_1,\dots,y_k;z)\\\times\prod_{i=1}^k\mathrm{Hom}(x_{i1},\dots,x_{i\alpha_i};y_i)\\\times\prod_{i=1}^k\prod_{j=1}^{\alpha_i}\mathrm{Hom}(w_{ij1},\dots,w_{ij\beta_{ij}};x_{ij})
          \end{array}}\ar[r]\ar[d]
        &
        {\begin{array}{l}\mathrm{Hom}(x_{11},\dots,x_{k\alpha_k};z)\\\times\prod_{i=1}^k\prod_{j=1}^{\alpha_i}\mathrm{Hom}(w_{ij1},\dots,w_{ij\beta_{ij}};x_{ij})\end{array}}\ar[d]\\
        {\begin{array}{l}\mathrm{Hom}(y_1,\dots,y_k;z)\\\times\prod_{i=1}^k\mathrm{Hom}(w_{i11},\dots,w_{i\alpha_i\beta_{i\alpha_i}};y_i)\end{array}}
        \ar[r]
        &
        \mathrm{Hom}(w_{111},\dots,w_{k \alpha_k \beta_{k\alpha_k}};z).
      }
    \]
In addition, we require that the identity elements are both right and left identities for composition. Proceeding, for a multimorphism $f: (x_1, \ldots, x_k) \to y$, we set $\mathrm{dom}(f) := (x_1,\ldots, x_k)$ and $\mathrm{codom}(f) := y$.

\begin{example*}
Planar arc diagrams comprise a multicategory important to the work that follows. Let $p\mathbb{T}$ denote the multicategory whose
\begin{itemize}
    \item objects are the natural numbers, including zero,
    \item $\mathrm{Hom}_{p\mathbb{T}}(m_1,\ldots, m_k; n)$ is the collection of $(m_1,\ldots, m_k; n)$ planar arc diagrams, which we will denote by $\mathscr{D}_{(m_1,\ldots, m_k; n)}$.
\end{itemize}
Composition in $p\mathbb{T}$ is composition of planar arc diagrams, as defined at the beginning of this section. It follows immediately that $p\mathbb{T}$ is a multicategory with identity elements $1_n$, which is just a circle with $n$ marked points times the interval. Note that we can view $p\mathbb{T}$ as a multicategory \textit{enriched in categories} since $\mathscr{D}_{(m_1,\ldots, m_k; n)}$ can be viewed as a category whose morphisms are (chronological) cobordisms between planar arc diagrams of type $(m_1,\ldots, m_k; n)$. 
\end{example*}

A very similar multicategory, $\mathscr{G}$, will be the main object of study for the rest of this section. The objects of $\mathscr{G}$ will be crossingless matchings rather than natural numbers, but the more striking difference between $\mathscr{G}$ and $p\mathbb{T}$ is the composition rule. We define $\mathscr{G}$ now.

\begin{definition}
Define the multicategory $\mathscr{G}$ whose
\begin{itemize}
    \item objects are crossingless matchings, $\mathrm{Ob}(\mathscr{G}) = B^\bullet$;
    \item for crossingless matchings $x_i\in B^{m_i}$, $i=1,\ldots, k$, and $y\in B^n$, set
    \[\mathrm{Hom}_{\mathscr{G}} (x_1,\ldots, x_k; y) = \widehat{\mathscr{D}}_{(m_1,\ldots, m_k; n)} \times \mathbb{Z}^2.\]
\end{itemize}
Then, composition
\[\mathrm{Hom}(y_1,\ldots, y_k; z) \times \left(\prod_{i=1}^k \mathrm{Hom}(x_{i1},\ldots, x_{i \alpha_i}; y_i)\right) \to \mathrm{Hom}(x_{11},  \ldots,  x_{k\alpha_k}; z)\]
is defined by
\[
(\widehat{D}, p) \circ \left( (\widehat{D_1}, p_1) , \cdots, (\widehat{D_k}, p_k) \right) = \left(D(D_1,\ldots, D_k; \emptyset)^{\wedge}, p + \sum_{i=1}^k p_i + \abs{W_{\vec{x}\vec{y}z}((D_1,\ldots, D_k);D)}\right)
\]
where $D(D_1,\ldots, D_k; \emptyset)^{\wedge}$ means $D(D_1,\ldots, D_k; \emptyset)$ with all closed loops removed. Finally, the distinguished identity element $\mathrm{Id}_x$ associated to each crossingless matching $x$ is given by $(1_{\abs{x}}, (\abs{x}, 0)) \in \mathrm{Hom}(x;x)$.
\end{definition}

\begin{proposition}
\label{prop: G is multicategory}
$\mathscr{G}$ is a multicategory; in particular, composition in $\mathscr{G}$ is associative.
\end{proposition}

\begin{proof}
Consider the following compositions of multimorphisms.
\begin{equation}
\label{eg.tree}
\tikz[baseline={([yshift=-.5ex]current bounding box.center)}, scale=1]{
\node at (-5.3,4) {$(w_{111},$};
\node (w1) at (-4,4) {$\cdots,$};
\node at (-2.6,4) {$w_{11\beta_{11}})$};
\node at (-1,4) {$\times\cdots\times$};
\node at (.5,4) {$(w_{1 \alpha_1 1},$};
\node (w2) at (1.7,4) {$\cdots,$};
\node at (3.3,4) {$w_{1 \alpha_1 \beta_{1 \alpha_1}})$};
\node at (-5.3,-1.5) {$(w_{k 1 1},$};
\node (w3) at (-4,-1.5) {$\cdots,$};
\node at (-2.6,-1.5) {$w_{k 1 \beta_{k 1}})$};
\node at (-1,-1.5) {$\times\cdots\times$};
\node at (.5,-1.5) {$(w_{k \alpha_k 1},$};
\node (w4) at (1.7,-1.5) {$\cdots,$};
\node at (3.3,-1.5) {$w_{k \alpha_k \beta_{k \alpha_k}})$};
\node (x3) at (-.1,0) {$(x_{k 1},$};
\node at (1,0) {$\cdots,$};
\node (x4) at (2.1,0) {$x_{k \alpha_k} )$};
\node (x1) at (-.1,2.5) {$(x_{1 1},$};
\node at (1,2.5) {$\cdots,$};
\node (x2) at (2.1,2.5) {$x_{1 \alpha_1} )$};
\node at (1,1.75) {$\times$};
\node at (1,0.58) {$\times$};
\node at (1,1.25) {$\vdots$};
\node (y1) at (5.2,1.25) {$(y_1,$};
\node at (6,1.25) {$\cdots,$};
\node (y2) at (6.8,1.25) {$y_k)$};
\node (z) at (8.5, 1.25) {$z$};
\draw[->] (w1) to[out=270,in=90] node[pos=.5,below,arrows=-]
        {$D_{11}$} (x1);
\draw[->] (w2) to[out=270,in=90] node[pos=.5,right,arrows=-]
        {$D_{1\alpha_1}$} (x2);
\draw[->] (w3) to[out=90,in=270] node[pos=.5,above,arrows=-]
        {$D_{k 1}$} (x3);
\draw[->] (w4) to[out=90,in=270] node[pos=.5,right,arrows=-]
        {$D_{k \alpha_k}$} (x4);
\draw[->] (x2) to[out=0,in=90] node[pos=.5,above,arrows=-]
        {$D_1$} (y1);
\draw[->] (x4) to[out=0,in=270] node[pos=.5,below,arrows=-]
        {$D_k$} (y2);
\draw[->] (y2) to[out=0,in=180] node[pos=.5,above,arrows=-]
        {$D$} (z);
}
\end{equation}
Our goal is to verify the associativity of these compositions in $\mathscr{G}$; i.e.,
\[
\prod_{i=1}^k\prod_{j=1}^{\alpha_i}(D_{ij}, p_{ij}) \circ \left( \prod_{i=1}^k (D_i, p_i) \circ (D,p) \right) = \left(\prod_{i=1}^k\prod_{j=1}^{\alpha_i}(D_{ij}, p_{ij}) \circ \prod_{i=1}^k (D_i, p_i)\right) \circ (D,p).
\]
In either case, the composition yields
\[
D\left(
            D_1(D_{11},\ldots, D_{1\alpha_1}),
            D_2(D_{21},\ldots, D_{2 \alpha_2}), \ldots,
            D_k(D_{k 1},\ldots, D_{k\alpha_k}) \right)^{\wedge}
\]
in the first coordinate. In the former case, the composition yields 
\begin{equation}
\label{eq:411LHS}
\begin{split}
p  + \sum_{i=1}^k p_i + \sum_{i=1}^k\sum_{j=1}^{\alpha_i} p_{ij} & + \abs{W_{\vec{x}\vec{y}z}((D_1,\ldots, D_k);D)} \\ &+ \abs{W_{\vec{w}\vec{x} z} ((D_{11},\ldots, D_{k \alpha_k}); D(D_1,\ldots, D_k))}
\end{split}
\end{equation}
in the second coordinate. In the latter case, the composition yields
\begin{equation}
\label{eq:411RHS}
\begin{split}
p + \sum_{i=1}^k p_i + \sum_{i=1}^k \sum_{j=1}^{\alpha_i} p_{ij} & + \sum_{i=1}^k \abs{W_{\vec{w}_i\vec{x}_i y_i} ((D_{i1}, \ldots, D_{i \alpha_i}); D_i)} \\ 
&+ \abs{W_{\vec{w}\vec{y}z} ((D_1(D_{11}, \ldots, D_{1 \alpha_1}), \ldots, D_k(D_{k1}, \ldots, D_{k \alpha_k})); D)}
\end{split}
\end{equation}
in the second coordinate since, for each $i=1,\ldots, k$,
\[
\prod_{j=1}^{\alpha_i}(D_{ij}, p_{ij}) \circ (D_i, p_i) = \left(D_i(D_{i1},\ldots,D_{i\alpha_i}), p_i + \sum_{j=1}^{\alpha_i} p_{ij} + \abs{W_{\vec{w}_i \vec{x}_i y_i} ((D_{i1},\ldots, D_{i\alpha_i}),D_i}\right).
\]
The values (\ref{eq:411LHS}) and (\ref{eq:411RHS}) are equivalent since the total number of merges and splits of the sequence of cobordisms is unchanged; otherwise, the minimality condition on the Euler characteristic is contradicted.
\end{proof}

By a multipath, we mean a sequence of collections of composable multimorphisms. Explicitly, a \textit{multipath of length $n$} is a sequence of sequences of multimorphisms
\[
\left(
(f_{i_1}^1)_{i_1}, (f_{i_1 i_2}^2)_{i_1i_2}, \ldots, (f_{i_1i_2\ldots i_n}^n)_{i_1i_2\ldots i_n}
\right)
\]
with ranges $i_1 = 1,\ldots, k$, $i_2 = 1,\ldots, k_{i_1}$, up to $i_n = 1,\ldots, k_{i_1i_2\ldots i_{n-1}}$ such that 
\[
\mathrm{dom}(f_{i_1\ldots i_t}^t) = \left(\mathrm{codom}(f_{i_1\ldots i_{t}1})^{t+1}, \ldots, \mathrm{codom}(f_{i_1\ldots i_t k_{i_1\ldots i_t}}^{t+1})\right)
\]
for each $t=1, \ldots, n$. Denote by $\mathscr{C}^{[n]}$ the collection of multipaths of length $n$. As we proceed, we frequently confound terminology and refer to the sequence of multimorphisms obtained by taking the composites of a multipath as a multipath. For example, suppose that
\[
\left((f_{i_1}^1), (f_{i_1i_2}^2), (f_{i_1i_2i_3}^3)\right) \in \mathscr{C}^{[3]}
\]
with $i_1 = 1,\ldots, k$, $i_2 = 1, \ldots k_{i_1}$, and $i_3 = 1, \ldots, k_{i_1i_2}$. We'll denote by
\begin{equation}
\label{eq:multipathcomp22}
(f_{i_1}^1) \circ (f_{i_1i_2}^2) \circ (f_{i_1i_2i_3}^3)
\end{equation}
the sequence of composites
\begin{align*}
f_1^1 & \circ \left( \left(f_{11}^2 \circ (f_{111}^3, \ldots, f_{11k_{11}}^3)\right), \ldots, \left(f_{1k_1}^2 \circ (f_{1k_11}^3, \ldots, f_{1k_1k_{1k_1}}^3)\right)\right),
\\
& f_2^1 \circ \left( \left(f_{21}^2 \circ (f_{211}^3, \ldots, f_{21k_{21}}^3)\right) ,\ldots, \left(f_{2k_2}^2 \circ (f_{2k_21}^3, \ldots, f_{2k_2 k_{2k_2}}^3)\right) \right), \ldots,
\\
& f_k^1 \circ \left(\left(f_{k1}^2 \circ (f_{k11}^3, \ldots, f_{k1k_{k1}}^3)\right), \ldots, \left(f_{kk_k} \circ (f_{kk_k 1}^3, \ldots, f_{kk_k k_{kk_k}}^3\right)\right).
\end{align*}
Then, this sequence is frequently referred to as a multipath of length 3, when it is really a composite of such a multipath. Finally, distilling notation further, we'll write $\vec{f}^1 := (f_1^1,\ldots, f_k^1)$, $\vec{f}_i^2 := (f_{i1}^2, \ldots, f_{ik_i}^2)$ and similarly for $\vec{f}_{ij}^3$, and write the sequence of composites of multimorphisms (\ref{eq:multipathcomp22}) as 
\begin{equation}
\label{eq:finmultipath3}
\vec{f}^1\circ 
\left(\prod_{i=1}^k \vec{f}^2_i \right)
\circ
\left(\prod_{i=1}^k \prod_{j=1}^{k_i}\vec{f}_{ij}^3\right).
\end{equation}
We will replace $k_i$ with the notation $\alpha_i$, and similarly the notation $k_{ij}$ with the notation $\beta_{ij}$. This runs the risk of presenting confusion in light of the associator and compatibility maps introduced momentarily---we hope that the meaning of notation is clear presented in context.

We remark that if $\vec{f}^1$ is a single multimorphism, then the multipath (\ref{eq:finmultipath3}) can be pictured as (\ref{eg.tree}) from the previous proof. In general, $\vec{f}$ may consist of many multimorphisms, and we can think of a multipath as a collection of such diagrams---in other words, multipaths can be viewed as trees and forests.

\begin{definition}
A \emph{grading multicategory} is pair $(\mathscr{C},\alpha)$ where $\mathscr{C}$ is a multicategory and $\alpha : \mathscr{C}^{[3]} \rightarrow \mathbb{K}^\times$ is a 3-cocycle, meaning that for all
\[
\left(\vec{f}
,
\left(\prod_{i=1}^k \vec{g}_i \right)
,
\left(\prod_{i=1}^k \prod_{j=1}^{\alpha_i}\vec{h}_{ij}\right)
,
\left(\prod_{i=1}^k \prod_{j=1}^{\alpha_i}\prod_{k=1}^{\beta_{i j}} \vec{\ell}_{ijk}\right)\right)  \in \mathscr{C}^{[4]}
\]
(shortened to $\vec{f}, \vec{g}, \vec{h}, \vec{\ell} \in \mathscr{C}^{[4]}$), $\alpha$ satisfies the expression
\[
d\alpha(\vec{\ell}, \vec{h}, \vec{g}, \vec{f}) := 
\alpha(\vec{\ell},\vec{h}, \vec{g}) \alpha(\vec{\ell}, \vec{h}, \vec{f}\vec{g})^{-1} \alpha(\vec{\ell}, \vec{g}\vec{h}, \vec{f}) \alpha(\vec{h}\vec{\ell}, \vec{g}, \vec{f})^{-1} \alpha(\vec{h}, \vec{g}, \vec{f}) = 1.
\]
We call such an $\alpha$ an \emph{associator}.
\end{definition}

\subsection{\texorpdfstring{$\mathscr{G}$}{Lg} as a grading multicategory}
\label{ss:Gisgradingmulti}

Our goal is to show that there exists a suitable associator $\alpha$ endowing $\mathscr{G}$ with the structure of a grading multicategory. We will define $\alpha$ to be the product of two values associated to changes of chronologies, one explicit and the other implicit. 

We'll use the notation $\vec{D}$, $\vec{D}'$, and so on to denote collections of planar arc diagrams which form a multipath in $p\mathbb{T}$. If $\vec{D}$ is a single planar arc diagram $D$, and $\vec{D}' = (D_1,\ldots, D_n)$, then their composition, which will in general be denoted $\vec{D}(\vec{D}')$, is denoted $D(D_1,\ldots, D_n)$. In the general setting, the constituents of a multipath $g, g', g'', g''' \in \mathscr{G}^{[4]}$ will be written
\begin{align*}
g &= (\vec{D}, \vec{p}) = \prod_i (D_i, p_i)\\
g' &= (\vec{D}', \vec{p}') = \prod_{i,j} (D_{ij}, p_{ij})\\
g'' &= (\vec{D}'', \vec{p}'') = \prod_{i,j,k} (D_{ijk}, p_{ijk})\\
g''' &= (\vec{D}''', \vec{p}''') = \prod_{i,j,k,\ell} (D_{ijk\ell}, p_{ijk\ell}).
\end{align*}
On one hand, our indexing notation allows us to write $\vec{D}_i' = \prod_j (D_{ij}, p_{ij})$. Then, $\vec{D}(\vec{D}')$ denotes the collection $(D_1(\vec{D}_1'), \ldots, D_n(\vec{D}_n')$. We could also use our indexing notation to write, for example, $g'''$ as $\prod_{i,j,k} (\vec{D}_{i,j,k}''', \vec{p}_{i,j,k}''')$. Finally, we will denote by $P$ the sum of the entries of $\vec{p}$ (that is, $P = \sum_i p_i$) and similarly for the other cases; e.g., $P''' = \vec{p}''' \cdot \langle 1,\ldots, 1\rangle = \sum_{i,j,k,\ell} p_{ijk\ell}$.

As we proceed, we will make use of the following lemma. It is implicit in the proof of Proposition \ref{prop: G is multicategory}, but we restate it here.

\begin{lemma}
\label{lemma:degree of multipaths}
For any multipath of planar arc diagrams $\vec{D}$, $\vec{D}'$, and $\vec{D}''$ as above,
\[
\abs{W_{\vec{x}\vec{y}\vec{z}}(\vec{D}', \vec{D})} + \abs{W_{\vec{w}\vec{x}\vec{z}}(\vec{D}'', \vec{D}(\vec{D}'))} = \abs{W_{\vec{w}\vec{x}\vec{y}}(\vec{D}'', \vec{D}')} + \abs{W_{\vec{w}\vec{y}\vec{z}}(\vec{D}'(\vec{D}''), \vec{D})}.
\]
\end{lemma}

\begin{proof}
The compositions $W_{\vec{w}\vec{x}\vec{z}}(\vec{D}'', \vec{D}(\vec{D}')) \circ W_{\vec{x}\vec{y}\vec{z}}(\vec{D}', \vec{D})$ and $W_{\vec{w}\vec{y}\vec{z}}(\vec{D}'(\vec{D}''), \vec{D}) \circ W_{\vec{w}\vec{x}\vec{y}}(\vec{D}'', \vec{D}')$ (we assume the first cobordisms in both composites are the identity elsewhere) have the same source and target. Thus they are isotopic cobordisms---if this were not the case, the minimality condition on the Euler characteristic would be contradicted.
\end{proof}

To construct our associator, consider the change of chronology
\[
\begin{tikzcd}[row sep=huge]
& \vec{w} \vec{D}(\vec{D}'(\vec{D}'')) \vec{z} & \\
\vec{w}\vec{D}'' \vec{x} \otimes \vec{x}D(\vec{D}') \vec{z} \arrow[ur, "W_{\vec{w}\vec{x}\vec{z}}(\vec{D}''\text{,~} \vec{D}(\vec{D}'))"] & \xRightarrow{H_\alpha} & \vec{w} \vec{D}'(\vec{D}'') \vec{y} \otimes \vec{y} \vec{D} \vec{z} \arrow[ul, "W_{\vec{w}\vec{y}\vec{z}} (\vec{D}'(\vec{D}'')\text{,~} \vec{D})"'] \\
& \vec{w} \vec{D}'' \vec{x} \otimes \vec{x} \vec{D}' \vec{y} \otimes \vec{y} \vec{D} \vec{z} \arrow[ur, "W_{\vec{w}\vec{x}\vec{y}} (\vec{D}''\text{,~}\vec{D}') \otimes \mathbbm{1}_{\vec{y}\vec{D}\vec{z}}"'] \arrow[ul, "\mathbbm{1}_{\vec{w} \vec{D}'' \vec{x}} \otimes W_{\vec{x}\vec{y}\vec{z}}(\vec{D}'\text{,~} \vec{D})"]  &
\end{tikzcd}
\]
Define $\alpha_1(g'', g', g)$ to be the evaluation of this change of chronology $\iota(H_\alpha)$---notice that this component of $\alpha$ does not see the second coordinates of its inputs. Secondly, take
\begin{align*}
\alpha_2(g'', g', g) &= \lambda \left(\abs{W_{\vec{x}\vec{y}\vec{z}}(\vec{D}', \vec{D})}, \sum_{i,j,k} p_{ijk} \right) \\ 
&= \lambda \left(\,\abs{W_{\vec{x}\vec{y}\vec{z}}(\vec{D}', \vec{D})}, P'' \right) .
\end{align*}
Then, set
\[
\alpha = \alpha_2\alpha_1.
\]

\begin{remark}
\label{remark:associator}
This definition is clearly motivated by and generalizes the associator presented in \cite{naisse2020odd}. A property we will use frequently is that the degree of cobordisms decomposes into a sum of constituents; notice, for example, that
\begin{align*} 
\alpha_2(g''', g'', g') &= \lambda\left(\,\abs{W_{\vec{w}\vec{x}\vec{y}}(\vec{D}'', \vec{D}')}, P''' \right) \\ &= \lambda\left(\sum_{i} \abs{W_{\vec{w}_i\vec{x}_i\vec{y}_i}(\vec{D}_i'', \vec{D}_i)}, P''' \right)
\end{align*}
(We could rewrite the last line as $\prod_i \lambda\left(\abs{W_{\vec{w}_i\vec{x}_i\vec{y}_i}(\vec{D}_i'', \vec{D}_i)}, P'''\right)$, invoking the bilinearity of $\lambda$, although there might be slight confusion with this rewriting since $P'''$ is a sum involving the index $i$---indeed, the second coordinates of each term in this product are equivalent.) Finally, we remark that we can view $\alpha$ as coming from the following sequence of schematics, just as in \cite{naisse2020odd} (pictured for the case we have just described).
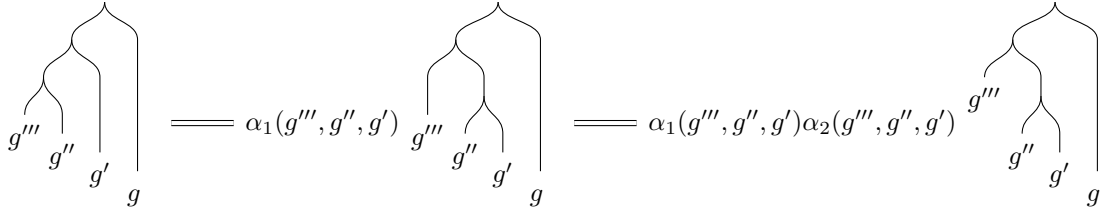
\[
\begin{tikzcd}
\tikz[scale=.5, baseline=2ex]{
	\draw (0,1) node[below]{$g'''$} .. controls (0,1.5) and (.5,1.5) .. (.5,2) .. controls (.5,2.5) and (1.25,2.5) .. (1.25,3) .. controls (1.25,3.5) and (2.125,3.5) .. (2.125,4);
	\draw (1,.5) node[below]{$g''$} -- (1,1) .. controls (1,1.5) and (.5,1.5) .. (.5,2);
	\draw (2, 0) node[below]{$g'$} -- (2,2) .. controls (2,2.5) and (1.25,2.5) .. (1.25,3);
	\draw (3,-.5) node[below]{$g\phantom{'}$} -- (3,3) .. controls (3,3.5) and (2.125,3.5) .. (2.125,4);
} \arrow[r,equals]
&
\alpha_1(g''', g'', g')\tikz[scale=.5, baseline=2ex]{
	\draw (0,1) node[below]{$g'''$} -- (0,2) .. controls (0,2.5) and (.75,2.5) .. (.75,3) .. controls (.75,3.5) and (1.875,3.5) .. (1.875,4);
	\draw (1,.5) node[below]{$g''$} .. controls (1,1) and (1.5,1) .. (1.5,1.5);
	\draw (2, 0) node[below]{$g'$} -- (2,.5) .. controls (2,1) and (1.5,1) .. (1.5,1.5) -- (1.5,2) .. controls (1.5,2.5) and (.75,2.5) .. (.75,3);
	\draw (3,-.5) node[below]{$g\phantom{'}$} -- (3,3) .. controls (3,3.5) and (1.875,3.5) .. (1.875,4);
} 
\arrow[r, equals]
&
\alpha_1(g''', g'', g')\alpha_2(g''', g'', g')\tikz[scale=.5,baseline=2ex]{
	\draw (0,2) node[below]{$g'''$} .. controls (0,2.5) and (.75,2.5) .. (.75,3) .. controls (.75,3.5) and (1.875,3.5) .. (1.875,4);
	\draw (1,.5) node[below]{$g''$} .. controls (1,1) and (1.5,1) .. (1.5,1.5);
	\draw (2, 0) node[below]{$g'$} -- (2,.5) .. controls (2,1) and (1.5,1) .. (1.5,1.5) -- (1.5,2) .. controls (1.5,2.5) and (.75,2.5) .. (.75,3);
	\draw (3,-.5) node[below]{$g\phantom{'}$} -- (3,3) .. controls (3,3.5) and (1.875,3.5) .. (1.875,4);
}
\end{tikzcd}
\]
\end{remark}

\begin{proposition}
\label{prop:associator}
The map $\alpha: \mathscr{G}^{[3]} \to R$ is a 3-cocycle.
\end{proposition}

\begin{proof}
This proof is completely analogous to the proof of Proposition 5.4 in \cite{naisse2020odd}---we represent their proof in the context of grading multicategories. As in the original case, $d\alpha(g''',g'',g',g)$ computes the difference between the paths of the diagram below.
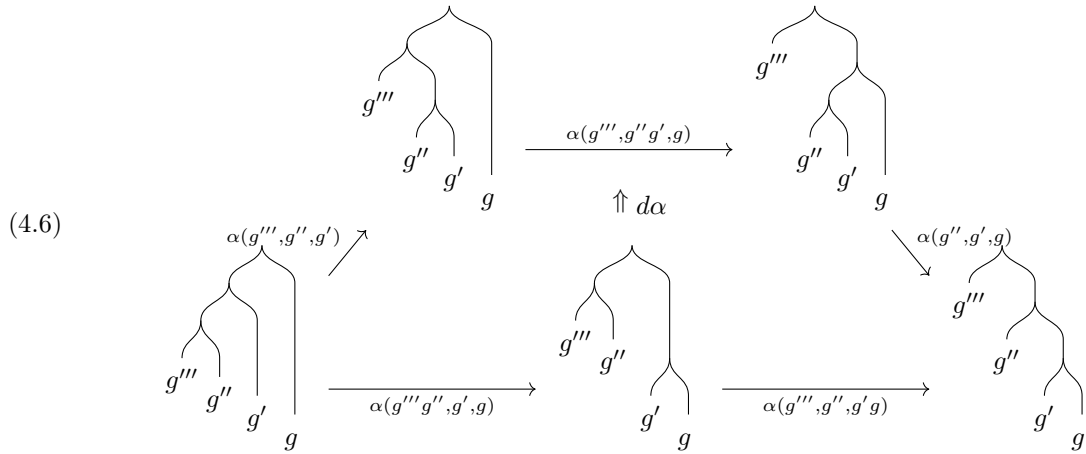
\begin{equation}
\label{multiassoc}
\begin{tikzcd}[column sep=1ex, row sep = .5ex]
&
\tikz[scale=.5]{
	\draw (0,2) node[below]{$g'''$} .. controls (0,2.5) and (.75,2.5) .. (.75,3) .. controls (.75,3.5) and (1.875,3.5) .. (1.875,4);
	\draw (1,.5) node[below]{$g''$} .. controls (1,1) and (1.5,1) .. (1.5,1.5);
	\draw (2, 0) node[below]{$g'$} -- (2,.5) .. controls (2,1) and (1.5,1) .. (1.5,1.5) -- (1.5,2) .. controls (1.5,2.5) and (.75,2.5) .. (.75,3);
	\draw (3,-.5) node[below]{$g\phantom{'}$} -- (3,3) .. controls (3,3.5) and (1.875,3.5) .. (1.875,4);
}
\ar{rr}{\alpha(g''',g''g',g)}
&{}&
\tikz[scale=.5]{
	\draw (0,3) node[below]{$g'''$} .. controls (0,3.5) and (1.125,3.5) .. (1.125,4);
	\draw (1,.5) node[below]{$g''$} .. controls (1,1) and (1.5,1) .. (1.5,1.5);
	\draw (2, 0) node[below]{$g'$} -- (2,.5) .. controls (2,1) and (1.5,1) .. (1.5,1.5) .. controls (1.5,2) and (2.25,2) ..(2.25,2.5);
	\draw (3,-.5) node[below]{$g\phantom{'}$} -- (3,1.5) .. controls (3,2) and (2.25,2) .. (2.25,2.5) -- (2.25,3) .. controls (2.25,3.5) and (1.125,3.5) .. (1.125,4);
}
\ar{dr}{\alpha(g'',g',g)}
&
\\
\tikz[scale=.5]{
	\draw (0,1) node[below]{$g'''$} .. controls (0,1.5) and (.5,1.5) .. (.5,2) .. controls (.5,2.5) and (1.25,2.5) .. (1.25,3) .. controls (1.25,3.5) and (2.125,3.5) .. (2.125,4);
	\draw (1,.5) node[below]{$g''$} -- (1,1) .. controls (1,1.5) and (.5,1.5) .. (.5,2);
	\draw (2, 0) node[below]{$g'$} -- (2,2) .. controls (2,2.5) and (1.25,2.5) .. (1.25,3);
	\draw (3,-.5) node[below]{$g\phantom{'}$} -- (3,3) .. controls (3,3.5) and (2.125,3.5) .. (2.125,4);
}
\ar{ur}{\alpha(g''',g'',g')}
\ar[swap]{rr}{\alpha(g'''g'',g',g)}
&&
\tikz[scale=.5]{
	\draw (0,2) node[below]{$g'''$} .. controls (0,2.5) and (.5,2.5) .. (.5,3) .. controls (.5,3.5) and (1.5,3.5) .. (1.5,4);
	\draw (1,1.5) node[below]{$g''$} -- (1,2) .. controls (1,2.5) and (.5,2.5) .. (.5,3);
	\draw (2, 0) node[below]{$g'$} .. controls (2,.5) and (2.5,.5) .. (2.5,1);
	\draw (3,-.5) node[below]{$g\phantom{'}$} -- (3,0) .. controls (3,.5) and (2.5,.5) .. (2.5,1) -- (2.5,3) .. controls (2.5,3.5) and (1.5,3.5) .. (1.5,4);
}
\ar[swap]{rr}{\alpha(g''',g'',g'g)}
 \ar[phantom]{u}{ \quad \rotatebox{90}{$\Rightarrow$}\ d\alpha}
&&
\tikz[scale=.5]{
	\draw (0,3) node[below]{$g'''$} .. controls (0,3.5) and (.875,3.5) .. (.875,4);
	\draw (1,1.5) node[below]{$g''$} .. controls (1,2) and (1.75,2) .. (1.75,2.5);
	\draw (2, 0) node[below]{$g'$} .. controls (2,.5) and (2.5,.5) .. (2.5,1);
	\draw (3,-.5) node[below]{$g\phantom{'}$} -- (3,0) .. controls (3,.5) and (2.5,.5) .. (2.5,1) -- (2.5,1.5) .. controls (2.5,2) and (1.75,2) .. (1.75,2.5) -- (1.75,3) .. controls (1.75,3.5) and (.875,3.5) .. (.875,4);
}
\end{tikzcd}
\end{equation}
On one hand, we are comparing two locally vertical changes of chronology with the same source and target, so the following diagram commutes by Proposition \ref{PutyraHammer}.
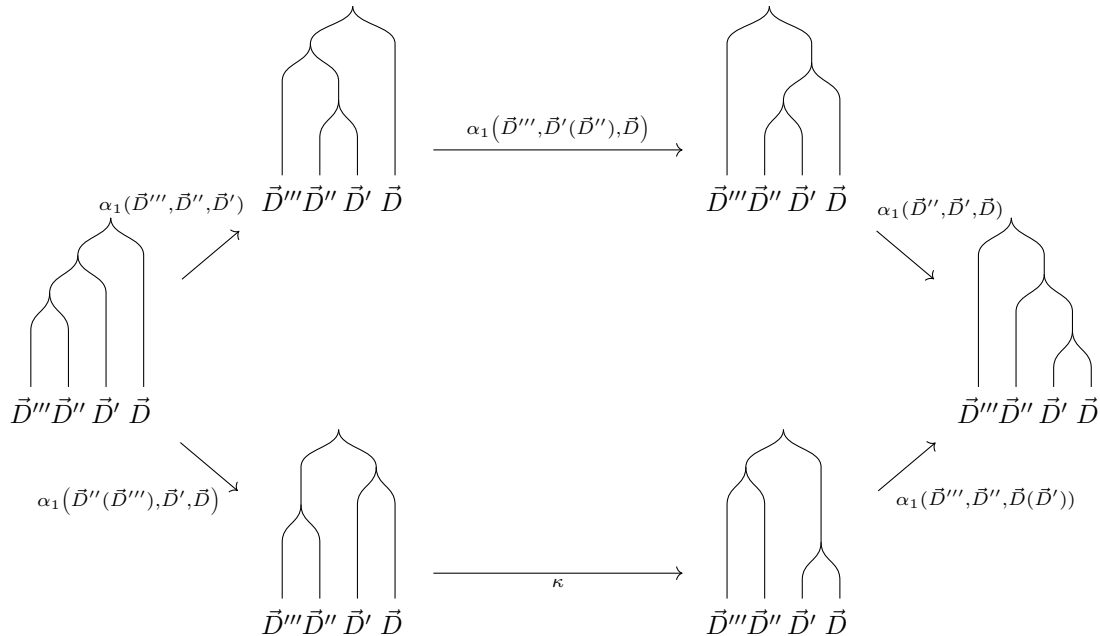
\[
\begin{tikzcd}[column sep=5ex, row sep = -2ex]
&
\tikz[scale=.5]{
	\draw (0,-.5)node[below]{$\vec{D}'''$}
		--
		(0,2) 
		.. controls (0,2.5) and (.75,2.5) ..
		(.75,3)
		.. controls (.75,3.5) and (1.875,3.5) ..
		(1.875,4);
	\draw (1,-.5) node[below]{$\vec{D}''$}
		--
		(1,.5) 
		.. controls (1,1) and (1.5,1) .. 
		(1.5,1.5);
	\draw (2, -.5) node[below]{$\vec{D}'$}
		--
		(2,.5)
		.. controls (2,1) and (1.5,1) .. 
		(1.5,1.5)
		--
		(1.5,2)
		.. controls (1.5,2.5) and (.75,2.5) ..
		(.75,3);
	\draw (3,-.5) node[below]{$\vec{D}\phantom{'}$}
		--
		(3,3)
		.. controls (3,3.5) and (1.875,3.5) ..
		(1.875,4);
}
\ar{rrrr}{\alpha_1\left(\vec{D}''', \vec{D}'(\vec{D}''), \vec{D} \right)}
&&{}&&
\tikz[scale=.5]{
	\draw (0,-.5) node[below]{$\vec{D}'''$}
		--
		(0,3)
		.. controls (0,3.5) and (1.125,3.5) ..
		(1.125,4);
	\draw (1,-.5) node[below]{$\vec{D}''$}
		--
		(1,.5)
		.. controls (1,1) and (1.5,1) .. 
		(1.5,1.5);
	\draw (2, -.5) node[below]{$\vec{D}'$}
		--
		(2,.5)
		.. controls (2,1) and (1.5,1) .. 
		(1.5,1.5)
		.. controls (1.5,2) and (2.25,2) ..
		(2.25,2.5);
	\draw (3,-.5) node[below]{$\vec{D}\phantom{'}$}
		--
		(3,1.5)
		.. controls (3,2) and (2.25,2) ..
		(2.25,2.5)
		-- 
		(2.25,3)
		.. controls (2.25,3.5) and (1.125,3.5) ..
		(1.125,4);
}
\ar[pos=-.1]{dr}{\alpha_1(\vec{D}'', \vec{D}', \vec{D})}
&
\\
\tikz[scale=.5]{
	\draw (0,-.5) node[below]{$\vec{D}'''$}
		--
		(0,1)
		.. controls (0,1.5) and (.5,1.5) ..
		(.5,2)
		.. controls (.5,2.5) and (1.25,2.5) ..
		(1.25,3)
		.. controls (1.25,3.5) and (2.125,3.5) ..
		(2.125,4);
	\draw (1,-.5) node[below]{$\vec{D}''$}
		--
		(1,1)
		.. controls (1,1.5) and (.5,1.5) ..
		(.5,2);
	\draw (2, -.5) node[below]{$\vec{D}'$}
		--
		(2,2)
		.. controls (2,2.5) and (1.25,2.5) ..
		(1.25,3);
	\draw (3,-.5) node[below]{$\vec{D}\phantom{'}$}
		--
		(3,3)
		.. controls (3,3.5) and (2.125,3.5) ..
		(2.125,4);
}
\ar[pos=1.2]{ur}{\alpha_1(\vec{D}''',\vec{D}'',\vec{D}')}
\ar[swap,pos=.8]{dr}{\alpha_1\left(\vec{D}''(\vec{D}'''), \vec{D}', \vec{D}\right)}
&&&&&&
\tikz[scale=.5]{
	\draw (0,-.5) node[below]{$\vec{D}'''$}
		--
		(0,3)
		.. controls (0,3.5) and (.875,3.5) ..
		(.875,4);
	\draw (1,-.5) node[below]{$\vec{D}''$}
		--
		(1,1.5)
	 	.. controls (1,2) and (1.75,2) ..
	 	(1.75,2.5);
	\draw (2, -.5) node[below]{$\vec{D}'$}
		--
		(2,0)
		.. controls (2,.5) and (2.5,.5) ..
		(2.5,1);
	\draw (3,-.5) node[below]{$\vec{D}\phantom{'}$}
		--
		(3,0)
		.. controls (3,.5) and (2.5,.5) ..
		(2.5,1)
		--
		(2.5,1.5)
		.. controls (2.5,2) and (1.75,2) ..
		(1.75,2.5)
		--
		(1.75,3)
		.. controls (1.75,3.5) and (.875,3.5) ..
		(.875,4);
}
\\
&
\tikz[scale=.5]{
	\draw (0,-.5) node[below]{$\vec{D}'''$}
		--
		(0,1)
		.. controls (0,1.5) and (.5,1.5) ..
		(.5,2)
		-- 
		(.5,3)
		.. controls  (.5,3.5) and (1.5,3.5) ..
		(1.5,4);
	\draw (1,-.5) node[below]{$\vec{D}''$}
		-- 
		(1,1)
	 	.. controls (1,1.5) and (.5,1.5) ..
	 	(.5,2);
	\draw (2, -.5) node[below]{$\vec{D}'$}
		--
		(2,2)
		.. controls (2,2.5) and (2.5,2.5) .. 
		(2.5,3);
	\draw (3,-.5) node[below]{$\vec{D}\phantom{'}$}
		--
		(3,2)
		.. controls (3,2.5) and (2.5,2.5) .. 
		(2.5,3)
		.. controls  (2.5,3.5) and (1.5,3.5) ..
		(1.5,4);
}
\ar[swap]{rrrr}{\kappa}
&&{} 
&&
\tikz[scale=.5]{
	\draw (0,-.5) node[below]{$\vec{D}'''$}
		--
		(0,2)
		.. controls (0,2.5) and (.5,2.5) ..
		(.5,3)
		.. controls (.5,3.5) and (1.5,3.5) ..
		(1.5,4);
	\draw (1,-.5) node[below]{$\vec{D}''$}
		-- 
		(1,2)
	 	.. controls (1,2.5) and (.5,2.5) ..
	 	(.5,3);
	\draw (2, -.5) node[below]{$\vec{D}'$}
		--
		(2,0)
		.. controls (2,.5) and (2.5,.5) ..
		(2.5,1);
	\draw (3,-.5) node[below]{$\vec{D}\phantom{'}$}
		--
		(3,0)
		.. controls (3,.5) and (2.5,.5) ..
		(2.5,1)
		--
		(2.5,3)
		.. controls (2.5,3.5) and (1.5,3.5) ..
		(1.5,4);
}
\ar[swap,pos=.2]{ur}{\alpha_1(\vec{D}''',\vec{D}'',\vec{D}(\vec{D}'))}
&
\end{tikzcd}
\] 
Since the corresponding change of chronology consists only of the sliding of two chronological cobordisms past one another, we know by work in Section \ref{ss:chronologicalcobordisms and coc} that $\kappa$ is
\[
\lambda\left(\,\abs{W_{\vec{x}\vec{y}\vec{z}}(\vec{D}', \vec{D})}, \abs{W_{\vec{v}\vec{w}\vec{x}}(\vec{D}''', \vec{D}'')}\right).
\]
Thus, the contribution of $\alpha_1$ in equation (\ref{multiassoc}) is
\[
\text{top} = \kappa \text{~bot}.
\]

On the other hand, we can compute and compare the contributions of $\alpha_2$ on the top and bottom path of (\ref{multiassoc}). The top path evaluates to
\begin{align*}
\lambda\left(\,\abs{W_{\vec{w}\vec{x}\vec{y}}(\vec{D}'', \vec{D}')}, P''' \right) \cdot \lambda\left(\,\abs{W_{\vec{w},\vec{y}, \vec{z}} (\vec{D}'(\vec{D}''), \vec{D})}, P''' \right) \cdot \lambda\left(\,\abs{W_{\vec{x}\vec{y}\vec{z}}(\vec{D}', \vec{D})}, P''\right)
\end{align*}
or, applying bilinearity of $\lambda$,
\begin{equation}
\label{***}
\lambda\left(\,\abs{W_{\vec{w}\vec{x}\vec{y}}(\vec{D}'', \vec{D}')} + \abs{W_{\vec{w}\vec{y}\vec{z}} (\vec{D}'(\vec{D}''), \vec{D})}, P'''\right)\cdot \lambda\left(\,\abs{W_{\vec{x}\vec{y}\vec{z}}(\vec{D}', \vec{D})}, P'' \right).
\end{equation}
The bottom path is slightly trickier to evaluate, since the second coordinate of $\alpha_2(g'''g'', g', g)$ requires a computation. As in the proof of Proposition \ref{prop: G is multicategory}, this comes from summing the second coordinates of $g'''$ and $g''$ and the cobordisms among their coordinates; explicitly,
\[
\alpha_2(g'''g'', g', g) = \lambda\left(\abs{W_{\vec{x}\vec{y}\vec{z}}(\vec{D}', \vec{D})}, P''' + P'' + \sum_{i,j} \abs{W_{\vec{v}_{(i,j)} \vec{w}_{ij} \vec{x}_{ij}} (\vec{D}_{ij}''', \vec{D}''_{ij})}\right).
\]
The last summation in the second coordinate can be rewritten as in Remark \ref{remark:associator}: we find that the bottom path evaluates to
\[
\lambda\left(\,\abs{W_{\vec{x}\vec{y}\vec{z}}(\vec{D}', \vec{D})}, P''' + P'' + \abs{W_{\vec{v} \vec{w} \vec{x}} (\vec{D}''', \vec{D}'')}\right) \cdot \lambda\left(\abs{W_{\vec{w}\vec{x}\vec{z}}(\vec{D}'', \vec{D}(\vec{D}'))}, P''' \right).
\]
Decomposing via bilinearity yields
\begin{align*}
& \lambda\left(\,\abs{W_{\vec{x}\vec{y}\vec{z}}(\vec{D}', \vec{D})}, P''' \right) \cdot \lambda\left(\,\abs{W_{\vec{x}\vec{y}\vec{z}}(\vec{D}', \vec{D})}, P'' \right)
\cdot \lambda\left(\,\abs{W_{\vec{x}\vec{y}\vec{z}}(\vec{D}', \vec{D})}, \abs{W_{\vec{v} \vec{w} \vec{x}} (\vec{D}''', \vec{D}'')}\right) \\ & \qquad \cdot \lambda\left(\,\abs{W_{\vec{w}\vec{x}\vec{z}}(\vec{D}'', \vec{D}(\vec{D}'))}, P''' \right).
\end{align*}
Combining the first and last term, and reordering suggestively, gives the product
\begin{align*}
&\lambda\left(\,\abs{W_{\vec{x}\vec{y}\vec{z}}(\vec{D}', \vec{D})}, \abs{W_{\vec{v} \vec{w} \vec{x}} (\vec{D}''', \vec{D}'')}\right) \cdot \lambda\left(\,\abs{W_{\vec{x}\vec{y}\vec{z}}(\vec{D}', \vec{D})} + \abs{W_{\vec{w}\vec{x}\vec{z}}(\vec{D}'', \vec{D}(\vec{D}'))}, P''' \right) \\ & \qquad \cdot \lambda\left(\,\abs{W_{\vec{x}\vec{y}\vec{z}}(\vec{D}', \vec{D})}, P'' \right).
\end{align*}
In this rewriting, the first term is $\kappa$. Moreover, by Lemma \ref{lemma:degree of multipaths}, the first coordinate of the second term is equivalent to the fist coordinate of the first term of (\ref{***}). Thus, the overall contribution of $\alpha_2$ in equation (\ref{multiassoc}) is
\[
\kappa \text{~top} = \text{bot}.
\]
Together, this provides that $d\alpha = 1$, as desired.
\end{proof}

\subsection{Generalities on modules graded by grading multicategories}
\label{ss:gradingcatgeneral}

Before proceeding with the grading multicategory at hand, we note generalities of $\mathscr{C}$-graded modules. That is, we consider the ways in which results of Section 4 of \cite{naisse2020odd} lift to the setting of grading multicategories. Throughout, $\mathscr{C}$ is a grading multicategory with associator $\alpha$ over a unital, commutative ring $\mathbb{K}$.

By a \textit{$\mathscr{C}$-graded $\mathbb{K}$-module}, we mean a $\mathbb{K}$-module $M$ with decomposition $M = \bigoplus_{g\in \mathrm{Mor}(\mathscr{C})} M_g$ where $g$ is a multimorphism of $\mathscr{C}$ (not an object). As before, we write $\abs{x} = g$ whenever $x \in M_g$. This generalizes the notion of grading by a category, introduced in \cite{naisse2020odd}, which in turn generalized the notion of grading by a group (take the category consisting of a single element $\star$ and $\mathrm{End}(\star) = G$). Of course, we are interested in the case $\mathscr{C} = \mathscr{G}$ and $\mathbb{K} = R$.

Tensor products in this setting are rather odd in the sense that their graded structure has a few different interpretations. This choice should be clear given the context. In one case, if $M$ and $M'$ are two $\mathscr{C}$-graded $\mathbb{K}$-modules, then we can define
\[
M'\otimes M = \bigoplus_{h\in\mathrm{Mor}(\mathscr{C})} (M'\otimes M)_h
\]
where
\[
(M' \otimes M)_h = \bigoplus_{h = g\circ g'} M'_{g'} \otimes_{\mathbb{K}} M_g.
\] 
Notice that this definition does not make full use of the flexibility offered by a grading multicategory. On the other hand, for $\mathscr{C}$-graded modules $M_1,\ldots, M_k$, $M$, we can view the tensor product over $\mathbb{K}$ as $\mathscr{C}$-graded by defining
\[
(M_1\otimes \cdots \otimes M_k) \otimes M = \bigoplus_{h \in \mathrm{Mor}(\mathscr{C})} [(M_1 \otimes \cdots \otimes M_k) \otimes M]_h
\]
where
\[
[(M_1 \otimes \cdots \otimes M_k) \otimes M]_h = \bigoplus_{h = g \circ (g_1, \ldots, g_k)} (M_{1, g_1} \otimes \cdots \otimes M_{k, g_k}) \otimes M_{g}.
\]
Notice that $M_1\otimes\cdots \otimes M_k$ is interpreted as a collection of $\mathscr{C}$-graded modules, but \textit{not} as a $\mathscr{C}$-graded module itself. Rather, $M_1\otimes \cdots \otimes M_k$ in the above scenario is viewed as $\mathscr{C}^k = \mathscr{C} \times \cdots \times \mathscr{C}$-graded in the sense that
\[
M_1 \otimes \cdots \otimes M_k = \bigoplus_{(g_1,\ldots, g_k)\in \mathrm{Mor}(\mathscr{C}^k)} M_{1, g_1} \otimes \cdots \otimes M_{k, g_k}.
\]
We will always abbreviate $M_1\otimes \cdots \otimes M_k$ (without interpretation as a $\mathscr{C}$-graded module itself) by $(M_1, \ldots, M_k)$ or, more succinctly, $\vec{M}$ to avoid confusion. For example, the above scenario will be written $(M_1, \ldots, M_k) \otimes M$ or, succinctly, $\vec{M} \otimes M$. Likewise, by $\vec{M}' \otimes \vec{M}$ we mean $(\vec{M}'_1\otimes M_1, \ldots, \vec{M}'_k \otimes M_k)$.

Denote by $\mathrm{Mod}_\mathbb{K}^{\mathscr{C}}$, or just $\mathrm{Mod}^\mathscr{C}$, the category of $\mathscr{C}$-graded $\mathbb{K}$-modules, whose morphisms are $\mathbb{K}$-linear maps which preserve grading. That is, for $f: M\to N$, we have $f(M_g) \subset N_g$ for each $g$. We call such maps \textit{$\mathscr{C}$-graded}, or just \textit{graded}. The associator of the grading multicategory $\mathscr{C}$ provides a coherence isomorphism
\[
\begin{tikzcd}[row sep=tiny]
(\vec{M}'' \otimes \vec{M}') \otimes M \arrow[r] & \vec{M}'' \otimes (\vec{M}' \otimes M) \\
(\vec{m}'' \otimes \vec{m}') \otimes m \arrow[r, mapsto] & \alpha\left(\abs{\vec{m}''}, \abs{\vec{m}'}, \abs{m}\right) \vec{m}'' \otimes (\vec{m}' \otimes m)
\end{tikzcd}
\]
where $\vec{m}'$ (and, similarly, $\vec{m}''$) is comprised of tensored homogeneous elements $m_i\in (M_i)_{\abs{m_i}}$, and $\abs{\vec{m}'} = (\abs{m_1}, \ldots, \abs{m_k})$ is the corresponding collection of multimorphisms (that is, $\mathscr{C}$-gradings).

Since the number of modules involved in a tensor product can vary, we have a collection of unit objects, one for each $k$, all defined as the tensor product of a single module: let $\mathbbm{1}$ denote the $\mathscr{C}$-graded $\mathbb{K}$-module $\bigoplus_{X\in \mathrm{Ob}(\mathscr{C})} (\mathbb{K})_{1_X}$. Then, $\mathbbm{1}^{\otimes k}$ is a unit object in the sense that there are (graded) isomorphisms (i.e., left- and right-unitors)
\[
\mathcal{L}: \mathbbm{1}^{\otimes k} \otimes M \cong M \qquad \mathrm{and} \qquad \mathcal{R}: M \otimes \mathbbm{1} \cong M
\]
of $\mathscr{C}$-graded modules which satisfy the triangle identity
\[
\begin{tikzcd}
\left((M_1 ,\ldots,  M_k) \otimes \mathbbm{1}^{\otimes k}\right) \otimes M  \arrow[rr, "\alpha"] \arrow[dr, "\prod_i \mathcal{R}_i \otimes \mathrm{id}_M"'] & & (M_1 ,\ldots, M_k) \otimes \left(\mathbbm{1}^{\otimes k} \otimes M\right) \arrow[dl, "\prod_i \mathrm{id}_{M_i} \otimes \mathcal{L}"] \\
& (M_1 ,\ldots,  M_k) \otimes M &
\end{tikzcd}
\]
where $\mathcal{R}_i$ means the right unitor applied to $M_i$. The left- and right-unitors we pick are determined by the associator: if each $m_i$ in $m_1\otimes \cdots\otimes m_k \in M_1\otimes \cdots\otimes M_k$ is homogeneous (with, say, $\abs{m_i}: (x_{i1},\ldots, x_{i\alpha_i}) \mapsto y_i'$), and similarly for $c_1\otimes \cdots \otimes c_k \in \mathbbm{1}^{\otimes k}$ and $m\in M$ (with, say, $\abs{m}: (y_1,\ldots, y_k) \mapsto z)$, we can choose left-unitor given by
\[
(c_1\otimes \cdots \otimes c_k) \otimes m \mapsto \alpha((1_{y_1},\ldots, 1_{y_k}), (1_{y_1},\ldots, 1_{y_k}), \abs{m})^{-1} c_1\cdots c_k m
\]
and right-unitor by
\[
(m_1\otimes \cdots \otimes m_k) \otimes (c_1 \otimes \cdots \otimes c_k) \mapsto \alpha((\abs{m_1},\ldots, \abs{m_k}), (1_{y_1},\ldots, 1_{y_k}), (1_{y_1},\ldots, 1_{y_k}))m_1 c_1 \otimes \cdots \otimes m_k c_k.
\]
To see why this satisfies the triangle identity, take $y_i = y_i'$ so that $\abs{\vec{m}'}$ and $\abs{m}$ are composable multimorphisms, and consider the path of length 4 given by 
\[
\vec{x} \xrightarrow{(\abs{m_1}, \ldots, \abs{m_k})} \vec{y} \xrightarrow{(1_{y_1},\ldots, 1_{y_k})} \vec{y} \xrightarrow{(1_{y_1},\ldots, 1_{y_k})} \vec{y} \xrightarrow{\abs{m}} z
\]
Then, the cocycle condition of $\alpha$ establishes that 
\begin{align*}
1 &= d\alpha (\abs{\vec{m}'}, 1_{\vec{y}}, 1_{\vec{y}}, \abs{m}) \\
&= \alpha(\abs{\vec{m}'}, 1_{\vec{y}}, 1_{\vec{y}}) \alpha(\abs{\vec{m}'}, 1_{\vec{y}}, \abs{m})^{-1} \alpha(1_{\vec{y}}, 1_{\vec{y}}, \abs{m}).
\end{align*}
This gives the triangle identity after re-arranging.

Since the cocyle requirement of the associator of a grading multicategory is exactly the pentagonal relation of monoidal categories, it follows from the work above that $\mathrm{Mod}_\mathbb{K}^{\mathscr{C}}$ has a structure resembling a monoidal category.

Finally, we briefly describe two important types of  $\mathscr{C}$-graded modules: algebras and multimodules. A \textit{$\mathscr{C}$-graded algebra} is a $\mathscr{C}$-graded $\mathbb{K}$-module $A = \bigoplus_{g\in \mathrm{Mor}(\mathscr{C})} A_g$, supported only in gradings $g$ which are single-input multimorphisms (i.e., morphisms) of $\mathscr{C}$, with a $\mathbb{K}$-linear multiplication map $\mu: A\otimes A \to A$ and a unit $1_X \in A_{\mathrm{Id}_X}$ for each $X \in \mathrm{Ob}(\mathscr{C})$ such that
\begin{enumerate}[label=(\roman*)]
    \item $\mu$ is graded: $\mu(A_{g'}, A_g) \subset A_{g \circ g'}$ for all $g', g\in \mathscr{C}$,
    \item $\mu$ is graded-associative: $\mu(\mu(z,y), x) = \alpha(\abs{z}, \abs{y}, \abs{x})\mu(z, \mu(y, x))$, and
    \item $\mu(1_Y, x) = \mathcal{L}(\mathrm{Id}_Y, \abs{x})~x$ and $\mu(x, 1_X) = \mathcal{R}(\abs{x}, \mathrm{Id}_X)~ x$ for all $x\in A_{\abs{x}:X \to Y}$.
\end{enumerate}
Before proceeding, we emphasize that $\mathscr{C}$-graded algebras are supported by single-input multimorphisms exclusively---really, $\mathscr{C}$-graded algebras are hardly different than the $\mathcal{C}$-graded algebras ($\mathcal{C}$ a category) of \cite{naisse2020odd}.

We'll write $\mu(x,y)$ as $x\cdot y$ when it is clear which multiplication is in use. Going on, we will only consider the tensor product $(A_1,\ldots, A_k)$---that is, $A_1\otimes \cdots \otimes A_k$ viewed as $\mathscr{C}^k$-graded---with multiplication $(a_1',\ldots, a_k')\cdot (a_1,\ldots, a_k)$, or, concisely, $\vec{a}' \cdot \vec{a}$, defined as $(\mu_{A_1}(a_1', a_1),\ldots, \mu_{A_k}(a_k', a_k))$.
 
Suppose $A_1,\ldots, A_k, B$ are $\mathscr{C}$-graded algebras. Then, a \textit{$\mathscr{C}$-graded $(A_1,\ldots, A_k; B)$-multimodule} is a $\mathscr{C}$-graded $\mathbb{K}$-module $M = \bigoplus_{g\in \mathrm{Mor}(\mathscr{C})} M_g$ with graded, $\mathbb{K}$-linear left and right actions
\[
\rho_L: (A_1, \ldots , A_k) \otimes M \to M 
\qquad \text{and} \qquad
\rho_R: M \otimes B \to M
\]
such that
\begin{enumerate}[label=(\roman*)]
    \item $\rho_L((\vec{a}' \cdot \vec{a}), m) = \alpha(\abs{\vec{a}'}, \abs{\vec{a}}, \abs{m}) \rho_L(\vec{a}', \rho_L(\vec{a}, m))$,
    \item $\rho_R(\rho_R(m, b'), b) = \alpha(\abs{m}, \abs{b}', \abs{b}) \rho_R(m, b' \cdot b)$,
    \item $\rho_R(\rho_L(\vec{a}, m), b) = \alpha(\abs{\vec{a}}, \abs{m}, \abs{b}) \rho_L(\vec{a}, \rho_R(m,b))$, and
    \item $\rho_L((1_Y,\ldots, 1_Y), m) = \mathcal{L}((1_Y,\ldots, 1_Y), \abs{m}) m$ and $\rho_R(m, 1_X) = \mathcal{R}(\abs{m}, \mathrm{Id}_X)$ for all $m\in M_{\abs{m}:X\to Y}$
\end{enumerate}
for all $\vec{a}', \vec{a} \in (A_1, \ldots, A_k)$, $b',b \in B$, and $m\in M$.

One should take caution: again, we are viewing $(A_1,\ldots, A_k)$ as a collection of $\mathscr{C}$-graded algebras, \textit{not} as a single $\mathscr{C}$-graded object. In particular, a $\mathscr{C}$-graded $(A_1,\ldots, A_k; B)$-multimodule is, perhaps surprisingly, \textit{not} equivalent to the notion of a $\mathscr{C}$-graded $(A_1\otimes_\mathbb{K} \cdots \otimes_\mathbb{K} A_k, B)$-bimodule. In particular, the left action $\rho_L$ is graded in the sense that
\[
\rho_L ((A_{1, g_1} \otimes \cdots \otimes A_{k, g_k}) \otimes M_g) \subset M_{g \circ (g_1,\ldots, g_k)}
\]
and \textit{not} in the sense that
\[
\rho_L ((A_{1, g_1} \otimes \cdots \otimes A_{k, g_k}) \otimes M_g) \subset M_{g \circ g_k \circ \cdots \circ g_1}.
\]
We define a $\mathscr{C}$-graded $(A,B)$-bimodule as a $\mathscr{C}$-graded $(A;B)$-multimodule for $\mathscr{C}$-graded algebras $A$ and $B$.

A graded map of $(A_1,\ldots, A_k; B)$-multimodules is a graded, $\mathbb{K}$-linear map satisfying 
\[
f(\rho_L(\vec{a}, m)) = \rho_L(\vec{a}, f(m)) \qquad \mathrm{and} \qquad f(\rho_R(m, b)) = \rho_R(f(m), b)
\]
for all $\vec{a}, m$, and $b$. Denote the category of $\mathscr{C}$-graded $(A_1,\ldots, A_k; B)$-multimodules, cumbersomely, by $\mathrm{MultiMod}^\mathscr{C}_R(A_1,\ldots, A_k; B)$.  As always, if it is clear what algebras we're working over, we denote this category by $\mathrm{MultiMod}^\mathscr{C}$. 

Take $M \in \mathrm{MultiMod}^\mathscr{C}(B_1,\ldots, B_k; C)$ and $M_i \in \mathrm{MultiMod}^\mathscr{C}(A_{i1},\ldots, A_{i \ell_i}; B_i)$ for each $i=1, \ldots, k$. Then $(M_1, \ldots, M_k) \otimes M$ has the structure of a $\mathscr{C}$-graded $(A_{11},\ldots, A_{k \ell_k}; C)$-multimodule by defining left- and right-actions so that the diagrams
\[
\begin{tikzcd}[column sep = 3ex]
(A_{11},\ldots, A_{k \ell_k}) \otimes ((M_1,\ldots, M_k) \otimes M) \ar{rrr} \ar[swap]{d}{\alpha^{-1}} &&& (M_1,\ldots, M_k) \otimes M \\
((A_{11},\ldots, A_{k \ell_k}) \otimes (M_1,\ldots, M_k)) \otimes M \ar[swap]{urrr}{\prod\rho_L \otimes 1}&&&
\end{tikzcd}
\]
and
\[
\begin{tikzcd}[column sep = 3ex]
((M_1, \ldots, M_k)\otimes M) \otimes C \ar{rrr} \ar[swap]{d}{\alpha} &&& (M_1,\ldots, M_k)\otimes M \\
(M_1,\ldots, M_k) \otimes (M \otimes C) \ar[swap]{urrr}{1 \otimes \rho_R}&&
\end{tikzcd}
\]
commute, interpreting $((A_{11},\ldots, A_{k \ell_k}) \otimes (M_1,\ldots, M_k))$ as
\[
((A_{11},\ldots, A_{1 \ell_1})\otimes M_1, \ldots, (A_{k1},\ldots, A_{k \ell_k})\otimes M_k).
\]
Explicitly, the left action is given by
\[
(a_{11},\ldots, a_{k \ell_k}) \cdot (\vec{m} \otimes m) := \alpha^{-1}(\abs{\vec{a}}, \abs{\vec{m}}, \abs{m}) (\rho_L^1((a_{11},\ldots, a_{1\ell_1}), m_1), \ldots, \rho_L^k((a_{k1},\ldots, a_{k\ell_k}), m_k)) \otimes m
\]
where $\rho_L^i$ is meant to denote the left action for the multimodule $M_i$. The right is just
\[
(\vec{m} \otimes m) \cdot c := \alpha(\abs{\vec{m}}, \abs{m}, \abs{c}) \vec{m} \otimes (\rho_R(m,c)).
\]

Finally, we note that the tensor product of $(M_1,\ldots, M_k)$ with $M$ over $\mathscr{C}$-graded algebras $(B_1, \ldots, B_k)$, denoted $(M_1,\ldots, M_k) \otimes_{(B_1,\ldots, B_k)} M$, is defined as
\[
(M_1,\ldots, M_k) \otimes M \big/ \left((\rho_R^1(m_1,b_1),\ldots, \rho_R^k(m_k, b_k))\otimes m - \alpha(\abs{\vec{m}},\abs{\vec{b}},\abs{m}) (m_1,\ldots,m_k) \otimes \rho_L((b_1,\ldots, b_k),m)\right)
\]
where $\rho_R^i$ is meant to denote the right action for the multimodule $M_i$. This is to say that the tensor product of $(M_1,\ldots, M_k)$ with $M$ over $(B_1,\ldots, B_k)$ is defined as the coequializer of the diagram
\[
\begin{tikzcd}
\left((M_1, \ldots, M_k) \otimes (B_1,\ldots, B_k)\right) \otimes M \arrow[dd, "\alpha"] \arrow[drr, "\prod \rho_R^i \otimes 1_M"]&& \\
&& (M_1, \ldots, M_k) \otimes M\\
(M_1, \ldots, M_k) \otimes \left( (B_1,\ldots, B_k) \otimes M\right) \arrow[urr, "\prod 1_{M_i} \otimes \rho_L"]\ && 
\end{tikzcd}
\]
in the category of $\mathscr{C}$-graded modules. Given $f:M \to N$ and $f_i: M_i \to N_i$ for all $i=1,\ldots, k$, we define the tensor product of maps
\[
(f_1,\ldots, f_k) \otimes f: (M_1,\ldots, M_k) \otimes_{(B_1,\ldots, B_k)} M \to (N_1, \ldots, N_i) \otimes_{(B_1,\ldots, B_k)} N
\]
by $\left((f_1,\ldots, f_k) \otimes f\right)\left((m_1,\ldots, m_k) \otimes m\right) = (f_1(m_1), \ldots, f_k(m_k)) \otimes f(m)$.

\subsection{\texorpdfstring{$\mathscr{G}$}{Lg}-graded arc modules}
\label{ss:Ggradedmultimodules}

If $D$ is a planar arc diagram of type $(m_1,\ldots, m_k; n)$, $\mathcal{F}(D)$ is a $\mathscr{G}$-graded $R$-multimodule where, for $u \in \mathcal{F}(D(x_1,\ldots, x_k; y)) \subset \mathcal{F}(D)$,
\[
\deg_\mathscr{G}(u) = (\widehat{D}, \deg_R(u)) \in \mathrm{Hom}_{\mathscr{G}}(x_1,\ldots, x_k; y).
\]
Then, the following lemmas are apparent.

\begin{lemma}
\label{GCompMaps}
The composition maps $\mu[(D_1,\ldots, D_k); D]$ preserve $\mathscr{G}$-grading.
\end{lemma}

\begin{proof}
This is by definitions: recall the composition maps
\[
\mu[(D_1,\ldots, D_k); D]: \left(\mathcal{F}(D_1), \ldots, \mathcal{F}(D_k)\right) \otimes \mathcal{F}(D) \to \mathcal{F}(D(D_1,\ldots, D_k))
\]
from the beginning of this section. Now, an element $(u_1,\ldots, u_k) \otimes u$ living in the source has degree
\begin{align*}
(D^\wedge, \deg_R(u)) \circ &\left((D_1^\wedge, \deg_R(u_1)), \ldots, (D_k^\wedge, \deg_R(u_k))\right) \\ &= \left(D(D_1,\ldots, D_k)^\wedge, \deg_R(u) + \sum_{i=1}^k \deg_R(u_i) + \abs{W_{\vec{x}\vec{y}z}((D_1,\ldots, D_k); D)}\right)
\end{align*}
where $\abs{u_i}: \vec{x}_i \to y_i$ and $\abs{u}: \vec{y} \to z$. On the other hand, 
$\deg_{\mathscr{G}}\left(\mu[(D_1,\ldots, D_k); D] ((u_1,\ldots, u_k) \otimes u)\right)$ is, by the definition of the degree of cobordisms, the second coordinate of the pair above.
\end{proof}

\begin{lemma}
\label{composition associativity}
For $u_{ij}\in \mathcal{F}(D_{ij})$, $u_i \in \mathcal{F}(D_i)$, and $u\in \mathcal{F}(D)$,
\[
\mu[\vec{D}'(\vec{D}''), D]\left(\mu[\vec{D}'', \vec{D}'](\vec{u}'', \vec{u}'), u\right) = \alpha\left(\abs{\vec{u}''}, \abs{\vec{u}'}, \abs{u}\right) \mu[\vec{D}'', D(\vec{D}')]\left(\vec{u}'', \mu[\vec{D}', D](\vec{u}', u)\right).
\]
\end{lemma}

\begin{proof}
This is immediate by the construction of the $\mu$ composition maps and the associator $\alpha$, recalling that $\mathrm{ChCob}_\bullet$ has the relation that $W' = \iota(H) W$ for each change of chronology $H: W \Rightarrow W'$.
\end{proof}

\begin{proposition}
\label{GHn}
The arc algebra $\mathcal{F}(1_n) = H^n$ is unital and associative as a $\mathscr{G}$-graded $R$-algebra.
\end{proposition}

\begin{proof}
Recall that the multiplication in $H^n$ is $\mu[1_n, 1_n]$, so Lemma \ref{GCompMaps} implies that the multiplication in $H^n$ is $\mathscr{G}$-graded, while Lemma \ref{composition associativity} implies that it is graded associative. Since we defined the left- and right-unitors via the associator, the third requirement of $\mathscr{G}$-graded algebras is also satisfied by Lemma \ref{composition associativity}, and we conclude that $H^n$ is a $\mathscr{G}$-graded algebra. Associativity follows from Lemma \ref{composition associativity} as well; for a proof of unitality, see the proof of Proposition 6.2 in \cite{naisse2020odd}.
\end{proof}

\begin{proposition}
Suppose $D$ is a planar arc diagram of type $(m_1,\ldots, m_k; n)$. Then $\mathcal{F}(D)$ is a $\mathscr{G}$-graded $(H^{m_1}, \ldots, H^{m_k}; H^n)$-multimodule with left action
\[
\rho_L^D = \mu[(1_{m_1},\ldots, 1_{m_k}), D]: (H^{m_1}, \ldots, H^{m_k}) \otimes \mathcal{F}(D) \to \mathcal{F}(D)\]
and right action
\[
\rho_R^D = \mu[D, 1_n]: \mathcal{F}(D) \otimes H^n \to \mathcal{F}(D).
\]
\end{proposition}

\begin{proof}
Just as the previous proposition, this follows by applying Lemmas \ref{composition associativity} and \ref{GCompMaps}, now knowing that $H^n$ is a $\mathscr{G}$-graded algebra for each $n$.
\end{proof}

Recall that if $\vec{D} = (D_1,\ldots, D_k)$ is a collection of planar arc diagrams of type $(\ell_{i1},\ldots, \ell_{i\alpha_i}; m_i)$ for each $i=1,\ldots,k$, then each of $\mathcal{F}(D_i)$ in $\mathcal{F}(\vec{D}) = \left(\mathcal{F}(D_1),\ldots ,\mathcal{F}(D_k)\right)$ is a $\mathscr{G}$-graded $\left(H^{\ell_{i1}}, \ldots, H^{\ell_{i\alpha_i}}; H^{m_i}\right)$-multimodule with left-aciton $\mu[(1_{\ell_{i1}}, \ldots, 1_{\ell{i\alpha_i}}); D_i]$ and right action $\mu[D_i; 1_{m_i}]$. Then, using results of \S \ref{ss:gradingcatgeneral}, we can view $\left(\mathcal{F}(D_1),\ldots, \mathcal{F}(D_k)\right) \otimes \mathcal{F}(D)$ as a type $\left(H^{\ell_{11}}, \ldots, H^{\ell_{k \alpha_k}}; H^m\right)$-multimodule. Similarly, comparing with the general case, we can define the tensor product $\mathcal{F}(\vec{D}) \otimes_{(H^{m_1}, \ldots, H^{m_k})} \mathcal{F}(D)$ as $\mathcal{F}(\vec{D}) \otimes \mathcal{F}(D)$ quotiented by 
\begin{equation}
\label{padcoequ}
(\mu[D_1, 1_{m_1}](u_1, x_1), \ldots, \mu[D_k, 1_{m_k}](u_k, x_k)) \otimes u - \alpha(\abs{\vec{u}}, \abs{\vec{x}}, \abs{u}) (u_1,\ldots, u_k) \otimes \mu[(1_{m_1}, \ldots, 1_{m_k}); D](\vec{x}, u)
\end{equation}
for $\vec{u}\in \mathcal{F}(\vec{D}')$, $\vec{x}\in (H^{n_1},\ldots, H^{n_k})$, and $u\in \mathcal{F}(D)$.

Mimicking \cite{Khovanov_2002}, we note each of the following. See also Section 6.1 of \cite{naisse2020odd}. The proofs of these statements are essentially identical to those found in Sections 2.6 and 2.7 of Khovanov's paper, and would take us too far afield to prove here---we leave them to the reader.

\begin{proposition}
$\mathcal{F}(D)$ is sweet: it is projective as a left $(H^{m_1},\ldots, H^{m_k})$-module and as a right $H^n$-module.
\end{proposition}

\begin{proposition}
\label{prop:compisomorphism}
If $D_i$ is a planar arc diagram of type $(\ell_{i1},\ldots, \ell_{i\alpha_i}; m_i)$ for each $i=1,\ldots, k$ and $D$ is a planar arc diagram of type $(m_1,\ldots, m_k; n)$, then there is an isomorphism of $\mathscr{G}$-graded $(H^{\ell_{11}},\ldots, H^{\ell_{k\alpha_k}}, H^n)$-multimodules
\[
\left(\bigotimes_{i=1}^k \mathcal{F}(D_i)\right) \otimes_{(H^{m_1}, \ldots, H^{m_k})} \mathcal{F}(D) \cong \mathcal{F}(D(D_1,\ldots, D_k; \emptyset))
\]
induced by $\mu[(D_1,\ldots, D_k), D]$. (The first collection of tensor products in the formula above are taken over $R$.)
\end{proposition}

We note that the sweetness proposition is important for the proof of the latter; again, see Sections 2.6 and 2.7 of \cite{Khovanov_2002}. Note also that $\mu[\vec{D}', D]: \mathcal{F}(\vec{D}') \otimes \mathcal{F}(D) \to \mathcal{F}(D(\vec{D}'))$ induces a maps $\mathcal{F}(\vec{D}') \otimes_{(H^{m_1},\ldots, H^{m_k})} \mathcal{F}(D) \to \mathcal{F}(D(\vec{D}'))$ by the universal property of the coequalizer. To see this, use Lemma \ref{composition associativity}: for $\vec{u}'\in \mathcal{F}(\vec{D}')$, $\vec{x}\in (H^{m_1}, \ldots,  H^{m_k})$, and $u\in \mathcal{F}(D)$, we have that
\[
\mu[\vec{D}', D] \left(\mu[\vec{D}', (1_{m_1}, \ldots, 1_{m_k})](\vec{u}', \vec{x}), u\right) = \alpha(\abs{\vec{u}'}, \abs{\vec{x}}, \abs{u}) \mu[\vec{D}', D] \left( \vec{u}', \mu[(1_{m_1}, \ldots,  1_{m_k}), D] (x,u) \right).
\]
Then, compare with equation (\ref{padcoequ}).

Sometimes, we will write ``$\otimes_{H}$'' as shorthand when its meaning is clear given context. For example, in the lemma below, the ``$\otimes_H$'' on the left means ``$\otimes_{(H^{\ell_{11}}, \ldots, H^{\ell_{k\alpha_k}})}$'' and the ``$\otimes_H$'' on the right means ``$\otimes_{(H^{m_1}, \ldots, H^{m_k})}$.'' We will also sometimes write ``$\mu[\vec{D}', D]$'' to mean ``the isomorphism of $\mathscr{G}$-graded bimodules induced by $\mu[\vec{D}', D]$.''

\begin{lemma}
The following diagram commutes for all $\vec{D}''$, $\vec{D}'$, and $D$. 
\[
\begin{tikzcd}[column sep = huge]
\left( \mathcal{F}(\vec{D}'') \otimes_H \mathcal{F}(\vec{D}') \right) \otimes_H \mathcal{F}(D) \arrow[r, "\mu\mathrm{[}\vec{D}''\mathrm{,} \vec{D}'\mathrm{]} \otimes 1"] \arrow[dd, "\alpha"'] & \mathcal{F}(\vec{D}'(\vec{D}'')) \otimes_H \mathcal{F}(D) \arrow[dr, "\mu\mathrm{[}\vec{D}'(\vec{D}'')\mathrm{,} D \mathrm{]}"] & \\
& & \mathcal{F}(D(\vec{D}'(\vec{D}''))) \\
\mathcal{F}(\vec{D}'') \otimes_H \left( \mathcal{F}(\vec{D}') \otimes_H \mathcal{F}(D) \right) \arrow[r, "1 \otimes \mu\mathrm{[}\vec{D}'\mathrm{,} D \mathrm{]}"] & \mathcal{F}(\vec{D}'') \otimes_H \mathcal{F}(D(\vec{D}')) \arrow[ur, "\mu\mathrm{[}\vec{D}''\mathrm{,} D(\vec{D}')\mathrm{]}"']& 
\end{tikzcd}
\]
\end{lemma}
\label{lem:multicatpentagon}

\begin{proof}
This is immediate from the definition of $\alpha$ and $\mu$, following Lemma \ref{composition associativity} in the language of Proposition \ref{prop:compisomorphism}.
\end{proof}

\newpage

\section{\texorpdfstring{$\mathscr{C}$}{Lg}-shifting systems and cobordisms}
\label{S: SHIFTING SYSTEMS}

Usually, grading shifts for graded algebraic objects are defined by way of the additive structure of $\mathbb{Z}$. This raises the question of how one should define grading shifts in a $\mathscr{C}$-graded setting. In the $\mathscr{G}$-graded case, we will see that the naive guess (i.e., a chronological cobordism in the first entry plus a $\mathbb{Z} \times \mathbb{Z}$-shift in the second) works out. We find the general definition of a $\mathscr{C}$-shifting system to be rather dense or enigmatic, so we introduce the more concrete $\mathscr{G}$-shifting system alongside the general definition, hoping it gives a helpful model for the reader. These definitions are provided in \S \ref{ss:gradingshiftsforG}, wherein we also describe the compatibility conditions required of a shifting system associated to a particular grading category. In \S \ref{ss:generalshifting}, we address generalities of shifting systems before a deep-dive into the theory of homogeneous maps for $\mathscr{C}$-graded multimodules (what does it mean for a map $f: M \to N$ of $\mathscr{C}$-graded multimodules to be homogeneous?) in \S \ref{ss:homomaps}. This includes the extension of our shifting system to a so-called ``shifting 2-system'' so that, in our context, we can interpret a composition of grading shifts as related to the grading shift associated to a composition of chronological cobordisms. Finally, $\mathscr{G}$-shifting systems are peculiar in the sense that changes of chronology induce natural transformations of grading shifts, which we detail in \S \ref{ss:gradingcocs}.

\subsection{A system of grading shifting functors for \texorpdfstring{$\mathscr{G}$}{Lg}}
\label{ss:gradingshiftsforG}

Suppose $\Delta: D \to D'$ is a chronological cobordism of planar arc diagrams $D, D'\in \mathscr{D}_{(m_1,\ldots, m_k; n)}$. If $x_i\in B^{m_i}$ for all $i = 1, \ldots, k$ and $y\in B_n$, then $\Delta$ induces a map from some subset of $\mathrm{Hom}_\mathscr{G}(x_1,\ldots, x_k; y)$ to $\mathrm{Hom}_\mathscr{G}(x_1,\ldots, x_k; y)$. Explicitly, given $v\in \mathbb{Z} \times \mathbb{Z}$, the pair $(\Delta, v)$ induces a map
\[
\varphi_{(\Delta, v)}: \{(D,p)\in \mathrm{Hom}_{\mathscr{G}}(x_1,\ldots, x_k; y)\} \to \mathrm{Hom}_{\mathscr{G}}(x_1,\ldots, x_k; y)
\]
defined
\[
\varphi_{(\Delta, v)}(D,p) = (D', p + v + \abs{\Delta(1_{x_1}, \ldots, 1_{x_k}; 1_{y})})
\]
where $\Delta(1_{x_1}, \ldots, 1_{x_k}; 1_{y})$ is the cobordism $\Delta$ corked by thickenings of the relevant crossingless matchings. We will see that any cobordism of planar arc diagrams (potentially paired with a $\mathbb{Z}\times\mathbb{Z}$-degree, in which case we call the cobordism \textit{weighted}) constitutes what we will call a \textit{$\mathscr{G}$-grading shift}.

In general, a collection $((\Delta_1, v_1) , \ldots, (\Delta_k, v_k))$ of chronological cobordisms of planar arc diagrams induces a grading shift on $((D_1, p_1), \ldots, (D_k, p_k))$. Viewing the former as a disjoint union of chronological cobordisms, there is ambiguity as to what chronology to pick. Hereafter, we fix a chronology which applies $\Delta_1$ on its component, then $\Delta_2, \ldots, $ then $\Delta_k$, followed by the identity cobordism weighed by $v_1$ on its component, then $v_2, \ldots,$ and finally $v_k$. A picture is probably more descriptive of the situation:
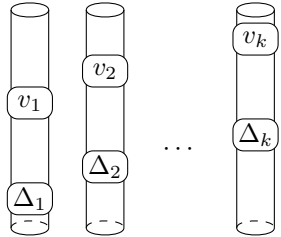
\[
\tikz[baseline=5ex, scale=.5]{
	\draw (0,-0.1) .. controls (0,.-.35) and (1,-.35) .. (1,-0.1);
	\draw[dashed] (0,-0.1) .. controls (0,.15) and (1,.15) .. (1,-0.1);
	\draw (0,-0.1) -- (0,5.7);
	\draw (1,-0.1) -- (1,5.7);
	\draw (0,5.7) .. controls (0,5.45) and (1,5.45) .. (1,5.7);
	\draw (0,5.7) .. controls (0,5.95) and (1,5.95) .. (1,5.7);
	\draw (2,-0.1) .. controls (2,.-.35) and (3,-.35) .. (3,-0.1);
	\draw[dashed] (2,-0.1) .. controls (2,.15) and (3,.15) .. (3,-0.1);
	\draw (2,-0.1) -- (2,5.7);
	\draw (3,-0.1) -- (3,5.7);
	\draw (2,5.7) .. controls (2,5.45) and (3,5.45) .. (3,5.7);
	\draw (2,5.7) .. controls (2,5.95) and (3,5.95) .. (3,5.7);
        \node at (4.5, 2) {$\cdots$};
	\draw (6,-0.1) .. controls (6,-.35) and (7,-.35) .. (7,-0.1);
	\draw[dashed] (6,-0.1) .. controls (6,.15) and (7,.15) .. (7,-0.1);
	\draw (6,-0.1) -- (6,5.7);
	\draw (7,-0.1) -- (7,5.7);
	\draw (6,5.7) .. controls (6,5.45) and (7,5.45) .. (7,5.7);
	\draw (6,5.7) .. controls (6,5.95) and (7,5.95) .. (7,5.7);
        \filldraw [fill=white, draw=black,rounded corners] (-.1,.25) rectangle (1.1, 1.1) node[midway] {$\Delta_1$};
        \filldraw [fill=white, draw=black,rounded corners] (1.9,1.1) rectangle (3.1, 1.95) node[midway] {$\Delta_2$};
        \filldraw [fill=white, draw=black,rounded corners] (5.9,1.95) rectangle (7.1, 2.8) node[midway] {$\Delta_k$};
        \filldraw [fill=white, draw=black,rounded corners] (-0.1,2.8) rectangle (1.1, 3.65) node[midway] {$v_1$};
        \filldraw [fill=white, draw=black,rounded corners] (1.9,3.65) rectangle (3.1, 4.5) node[midway] {$v_2$};
        \filldraw [fill=white, draw=black,rounded corners] (5.9,4.5) rectangle (7.1, 5.35) node[midway] {$v_k$};
}
\]
This is the chronology we mean when we write $(\vec{\Delta}, \vec{v})$. We choose this particular chronology so that our arguments remain analogous to those found in \cite{naisse2020odd}. Later on, we'll denote $\Delta(1_{x_1}, \ldots, 1_{x_k}; 1_{y})$ by $1_{\vec{x}} \Delta 1_{y}$. Again, this is especially helpful when dealing with a collection of cobordisms $(\Delta_1, \ldots, \Delta_n)$. The degree $\abs{1_{\vec{x}} (\Delta_1,\ldots, \Delta_n) 1_{\vec{y}}}$ is defined as the sum $\sum_{i=1}^n \abs{1_{\vec{x}_i} \Delta_i 1_{y_i}}$.

Now, for each $i=1,\ldots, k$, suppose $\Delta_i: D_i \to D_i'$ is a chronological chobordism for $D_i, D_i'\in \mathscr{D}_{(\ell_{i1},\ldots \ell_{i\alpha_i}; m_i)}$. We denote by $(\Delta_1,\ldots, \Delta_k) \bullet \Delta$ the chronological cobordism
\[
(\Delta_1,\ldots, \Delta_k) \bullet \Delta: D(D_1,\ldots, D_k) \to D'(D_1',\ldots, D_k')
\]
with chronology, as usual, dependant on indexing (first $\Delta$, then $\Delta_1$, and so on). If each of these cobordisms has a $\mathbb{Z}\times \mathbb{Z}$ weight, we'll set
\[
((\Delta_1, v_1),\ldots, (\Delta_k, v_k)) \bullet (\Delta, v) = \left((\Delta_1,\ldots, \Delta_k)\bullet \Delta, v + \sum_i v_i\right)
\]
for cases like the one above. Otherwise, $((\Delta_1, v_1),\ldots, (\Delta_k, v_k)) \bullet (\Delta, v) =0$. This multiplication defines what we will call a \textit{multimonoid}, whose elements are cobordisms of planar arc diagrams together with a neutral element $e$ and absorbing element $0$, with composition defined as above.

The collection of maps induced by cobordisms of planar arc diagrams $\{\varphi_{(\Delta, v)}\}$ constitutes a generalization of a shifting system, in the sense of \cite{naisse2020odd}. Explicitly, suppose $\mathscr{C}$ is a grading multicategory; a \textit{$\mathscr{C}$-grading shift} $\varphi$ is a collection of maps
\[
\varphi = \{\varphi^{\vec{X} \to Y}: \mathsf{D}^{\vec{X} \to Y} \subset \mathrm{Hom}_{\mathscr{C}}(\vec{X}; Y) \to \mathrm{Hom}_{\mathscr{C}}(\vec{X}; Y) \}_{\vec{X}, Y\in \mathrm{Ob}(\mathscr{C})}
\]
where $\vec{X} = (X_1,\ldots, X_k)$ for $X_1,\ldots, X_k \in \mathrm{Ob}(\mathscr{C})$. We write $\varphi(g)$ to mean $\varphi^{\vec{X}\to Y}(g)$ whenever $g \in \mathsf{D}^{\vec{X} \to Y}$. We write $\mathsf{D}$ to stand for ``domain'', and use the sans serif font to differentiate these from our notation for planar arc diagrams. In addition, let $\Sigma_{\text{min}}$ denote the category obtained from $\mathscr{C}$ by purging all multimorphisms besides the commuting endomorphisms: that is,
\begin{itemize}
    \item $\mathrm{Ob}(\Sigma_{\text{min}}) = \mathrm{Ob}(\mathscr{C})$ and
    \item $\mathrm{Hom}_{\Sigma_{\text{min}}}(X_1,\ldots, X_k; Y) = \begin{cases} \emptyset & k >1~\text{or}~X_1 \not=Y \\ \mathrm{Z}(\mathrm{End}_{\mathscr{C}} (Y; Y)) & \text{otherwise}\end{cases}$
\end{itemize}
Where $\mathrm{Z}$ stands for the center. Finally, by a \textit{multimonoid} $\mathscr{I}$, we mean a set equipped with an associative multiplication law
\[
\bullet: \mathscr{I}^k \times \mathscr{I} \to \mathscr{I}
\]
for each $k\ge 1$, and a neutral element $e$ so that $e^k \bullet i = i$ for each $k$ and $i \bullet e = i$ for all $i \in \mathscr{I}$. A multimonoid may also have an absorbing element $0$, so that $(j_1,\ldots, j_k) \bullet i = 0$ if any of $j_1,\ldots, j_k, i$ are $0$.

\begin{definition}
Suppose $\Sigma$ is a wide subcategory of $\mathscr{C}$ with at least all the morphisms of $\Sigma_{\text{min}}$. A \textit{$\mathscr{C}$-shifting system $S = (\mathscr{I},\Phi)$ relative $\Sigma$} for a grading multicategory $\mathscr{C}$ is a multimonoid $\mathscr{I}$ and a collection of $\mathscr{C}$-grading shifts $\Phi = \{\varphi_i\}_{i\in\mathscr{I}}$
such that
\begin{itemize}
    \item $\varphi_e$, called the \textit{neutral shift}, has
    \[
    \mathsf{D}_e^{\vec{X}\to Y} = \mathrm{Hom}_{\Sigma} (\vec{X};Y) \qquad \text{and} \qquad \varphi_e^{\vec{X} \to Y} = \mathrm{Id}_{\mathsf{D}_e^{\vec{X} \to Y}};
    \]
    \item given $\varphi_{j_1}^{(x_{11},\ldots, x_{1\alpha_1}; y_1)}, \ldots, \varphi_{j_n}^{(x_{n1},\ldots, x_{n\alpha_n};y_n)}$ and $\varphi_i^{(y_1,\ldots, y_n; z)}$, we have that 
    \[
    \mathsf{D}_i^{(y_1,\ldots, y_n;z)} \circ \prod_{k=1}^n \mathsf{D}_{j_k}^{(x_{k1},\ldots, x_{k\alpha_k}; y_k)} \subset 
    \mathsf{D}_{(j_1,\ldots, j_n) \bullet i}^{(x_{11}, \ldots, x_{n\alpha_n}; z)}
    \]
    and the diagram
    \[
    \begin{tikzcd}[row sep = large]
        \mathsf{D}_i^{(y_1,\ldots, y_n;z)} \times \prod_{k=1}^n\mathsf{D}_{j_k}^{(x_{k1},\ldots, x_{k\alpha_k}; y_k)} \arrow[r, "\circ"] \arrow[d, "(\varphi_i\text{,~}\prod_k \varphi_{j_k})"'] & \mathsf{D}_{(j_1,\ldots, j_n) \bullet i}^{(x_{11}, \ldots, x_{n\alpha_n}; z)} \arrow[d, "\varphi_{(j_1,\ldots, j_n) \bullet i}"] \\
        \mathrm{Hom}(y_1,\ldots, y_n; z) \times \prod_{k=1}^n \mathrm{Hom}(x_{k1},\ldots, x_{k\alpha_k};y_k) \arrow[r, "\circ"] & \mathrm{Hom}(x_{11},\ldots, x_{n\alpha_n}; z)
    \end{tikzcd}
    \]
    commutes;
    \item there is a subset $\mathscr{I}_{\mathrm{id}} \subset \mathscr{I}$ such that for all $k$ and all $X_1,\ldots, X_k, Y\in \mathrm{Ob}(\mathscr{C})$ there is a partition
    \[
    \mathrm{Hom}_{\mathscr{C}}(X_1,\ldots, X_k; Y) = \bigsqcup_{i\in \mathscr{I}_{\mathrm{id}}} \mathsf{D}_i^{(X_1,\ldots, X_k) \to Y}
    \]
    for which $\varphi_i = \mathrm{Id}_{\mathsf{D}_i^{\vec{X}\to Y}}$ for all $i \in \mathscr{I}_{\mathrm{id}}$;
    \item if $\mathscr{I}$ contains an absorbing element $0$, then $\varphi_0$, called the \textit{null shift}, always has $\mathsf{D}_0^{\vec{X}\to Y} = \emptyset$.
\end{itemize}
\end{definition}

\begin{remark}
We will frequently write $\mathsf{D}_{\vec{i}}^{\vec{X}\to \vec{Y}}$, or just $\mathsf{D}_{\vec{i}}$, to denote $\prod_{\ell} \mathsf{D}_{i_\ell}^{\vec{X}_\ell \to Y_\ell}$. Then, writing $g \in\mathsf{D}_{\vec{i}}$ means $g$ is an ordered tuple of morphisms as one expects. Similarly, $\varphi_{\vec{i}}(g)$ is understood component-wise. Also, we note that $\varphi_e$ is assumed only to preserve $\Sigma$. We refer the reader to Remark 4.10 of \cite{naisse2020odd} for a more detailed discussion.
\end{remark}

For example, take $\mathscr{I}$ to be the multimonoid $\{(\Delta, v)\}_{\Delta, v} \sqcup \{e, 0\}$ with multiplication $\bullet$ defined above. Taking $\mathscr{C} = \mathscr{G}$, notice that $\Sigma_{\text{min}}$ is the subcategory whose objects are crossingless matchings and whose morphisms are identity $(n;n)$-planar arc diagrams $(1_n, p): a \to a$ for $a \in B^n$, viewed only as endomorphisms. We will take $\Sigma$ to be the slightly larger category which allows for morphisms $(1_n, p):a \to b$ for potentially distinct $a, b\in B^n$. Using the notation of the above definition, to a chronological cobordism of planar arc diagrams $\Delta: D \to D'$, $D, D'\in \mathscr{D}_{(m_1,\ldots, m_k; n)}$ and $v\in \mathbb{Z}\times\mathbb{Z}$, we have a $\mathscr{G}$-grading shift $\varphi_{(\Delta, v)}$ so that for any crossingless matchings $x_1,\ldots, x_k, y$ with $\abs{x_i} = m_i$ and $\abs{y} = n$, 
\[
\mathsf{D}_{(\Delta, v)}^{(x_1,\ldots, x_k) \to y} = \{(D^\wedge, p) \in \mathrm{Hom}_{\mathscr{G}}(x_1,\ldots, x_k; y) : p\in \mathbb{Z}\times\mathbb{Z}\}.
\]

\begin{proposition}
\label{prop:G-shifting}
The multimonoid $\mathscr{I} = \{(\Delta, v)\}_{\Delta, v} \sqcup \{e,0\}$ together with the induced $\mathscr{G}$-grading shifts $\{\varphi_{i}\}_{i\in \mathscr{I}}$ form a $\mathscr{G}$-shifting system.
\end{proposition}

\begin{proof}
We define $\varphi_e$ and $\varphi_0$ so that the first and last points are satisfied. The second point is straightforward. Finally, for the third point, we take $\mathscr{I}_{\mathrm{id}} = \{(\mathbbm{1}_{D^\wedge}, (0,0)): D~\text{is a planar arc diagram}\}$, where $\mathbbm{1}_D$ is the identity cobordism on $D$.
\end{proof}

The definition of a $\mathscr{C}$-shifting system made no reference to the associator of the grading category $\mathscr{C}$. We say that a $\mathscr{C}$-shifting system $S$ is \textit{compatible} with the associator $\alpha$ of $\mathscr{C}$ if there is a family of maps
\[
\beta_{(\vec{j}_1,\ldots, \vec{j}_n) , \vec{i}}^{\vec{x}\vec{y}\vec{z}}: \prod_{k=1}^n\mathsf{D}_{\vec{j}_k}^{(\vec{x}_{k1},\ldots, \vec{x}_{k\alpha_k}; \vec{y}_k)} \times\mathsf{D}_{\vec{i}}^{(\vec{y}_1,\ldots, \vec{y}_n;\vec{z})} \to \mathbb{K}^{\times},
\]
for each $\vec{w}, \vec{x}, \vec{y}, \vec{z}$ consisting of objects of $\mathscr{C}$ and $\vec{i}, \vec{j}$ consisting of objects in $\mathscr{I}$, called \textit{compatibility maps} granted they satisfy the relations
\begin{align}
\label{eq:shiftingsystcomp}
\begin{split}
\alpha(g'',g',g)&
\beta^{\vec{w}\vec{x}\vec{z}}_{\vec{k}\bullet \vec{j,}\vec{i}}(g''g',g)
\beta^{\vec{x}\vec{y}\vec{z}}_{\vec{k},\vec{j}}(g'',g') \\
&= 
\beta^{\vec{w}\vec{y}\vec{z}}_{\vec{k}, \vec{j}\bullet \vec{i}}(g'', g'g)
\beta^{\vec{w}\vec{x}\vec{y}}_{\vec{j},\vec{i}}(g',g)
\alpha\left(\varphi^{\vec{y}\vec{z}}_{\vec{k}}(g''), \varphi^{\vec{x}\vec{y}}_{\vec{j}}(g'), \varphi^{\vec{w}\vec{x}}_{\vec{i}}(g)\right),
\end{split}
\end{align}
for all valid $g'', g, g$ and $\vec{i},\vec{j}, \vec{k}$, and $\beta_{e,e} = \beta_{(e,\ldots, e),(e,\ldots, e)} = 1$. Diagrammatically, this is to say that the following picture commutes (here, the boxed number $n$ refers to the $\mathscr{C}$-grading shift, and we suppress burdensome indices).
\begin{equation}
\label{beta_CompatCondition}
\begin{tikzcd}[row sep=huge, column sep=huge]
\tikz[yscale=.5,xscale=.75]{
	\draw (0,0) node[below,scale=.75]{$g''$}
        -- (0,1) 
		.. controls (0,1.5) and (.5,1.5) ..
		(.5,2)
		.. controls (.5,2.5) and (1.25,2.5) ..
		(1.25,3);
	\draw (1,-1.5) node[below,scale=.75]{$g'$}
		--
		(1,-1)
		-- 
		(1,-1)
		--
		(1,1)
		.. controls (1,1.5) and (.5,1.5) ..
		(.5,2);
	\draw (2,-3) node[below,scale=.75]{$g\phantom{'}$}
		--
		(2, -2.5)
		--
		(2,2)
		.. controls (2,2.5) and (1.25,2.5) ..
		(1.25,3);
        \node[fill=white,draw,rounded corners,scale=.75] at (0,0.35) {3};
        \node[fill=white,draw,rounded corners,scale=.75] at (1,-1) {2};
        \node[fill=white,draw,rounded corners,scale=.75] at (2, -2.5) {1};
}
\arrow[r, "\beta(g''\text{,~}g')"] \arrow[d, "\alpha(\varphi_3(g'')\text{,~}\varphi_2(g')\text{,~}\varphi_1(g))"']
&
\tikz[yscale=.5,xscale=.75]{
	\draw (0,0) node[below,scale=.75]{$g''$}
		--
		(0,0) 
		.. controls (0,.5) and (.5,.5) ..
		(.5,1)
		--
		(.5,2)
		.. controls (.5,2.5) and (1.25,2.5) ..
		(1.25,3);
	\draw (1,-1.5) node[below,scale=.75]{$g'$}
		-- 
		(1,0)   
		--
		(1,0)
		.. controls (1,.5) and (.5,.5) ..
		(.5,1);
	\draw (2,-3) node[below,scale=.75]{$g\phantom{'}$}
		--
		(2, -2.5)
		--
		(2,2)
		.. controls (2,2.5) and (1.25,2.5) ..
		(1.25,3);
        \node[fill=white,draw,rounded corners,scale=.75] at (.5,1.3) {$3\bullet2$};
        \node[fill=white,draw,rounded corners,scale=.75] at (2, -2.5) {$1$};
}
\arrow[r, "\beta(g''g'\text{,~}g)"]
&
\tikz[yscale=.5,xscale=.75]{
	\draw (0,-1) node[below,scale=.75]{$g''$}
		--
		(0,-1) 
		.. controls (0,-.5) and (.5,-.5) ..
		(.5,0)
		.. controls (.5,.5) and (1.25,.5) ..
		(1.25,1)
		--
		(1.25,3);
	\draw (1,-1.5) node[below,scale=.75]{$g'$}
		--
		(1,-1)
		.. controls (1,-.5) and (.5,-.5) ..
		(.5,0);
	\draw (2,-2) node[below,scale=.75]{$g\phantom{'}$}
		--
		(2, 0)
		.. controls (2,.5) and (1.25,.5) ..
		(1.25,1);
        \node[fill=white,draw,rounded corners,scale=.75] at (1.25,1.75) {$3\bullet2\bullet1$};
}
\arrow[d, "\alpha(g''\text{,~}g'\text{,~}g)"]
\\
\tikz[yscale=.5,xscale=.75]{
	\draw (0,1) node[below,scale=.75]{$g''$}
		--
		(0,1.5)
		--
		(0,2)
		.. controls (0,2.5) and (.75,2.5) ..
		(.75,3);
	\draw (1,-1.5) node[below,scale=.75]{$g'$}
		--
		(1,-1)
		-- 
		(1,-1)
		--
		(1,0)
		.. controls (1,0.5) and (1.5,0.5) ..
		(1.5,1);
	\draw (2,-3) node[below,scale=.75]{$g\phantom{'}$}
		--
		(2, -2.5)
		--
		(2,0)
		.. controls (2,0.5) and (1.5,0.5) ..
		(1.5,1)
        --
        (1.5,2)
		.. controls (1.5,2.5) and (.75,2.5) .. 
		(.75,3);
        \node[fill=white,draw,rounded corners,scale=.75] at (0,1.35) {$3$};
        \node[fill=white,draw,rounded corners,scale=.75] at (1,-1) {$2$};
        \node[fill=white,draw,rounded corners,scale=.75] at (2, -2.5) {$1$};
}
\arrow[r, "\beta(g'\text{,~}g)"]
&
\tikz[yscale=.5,xscale=.75]{
	\draw (0,1.5) node[below,scale=.75]{$g''$}
		--
		(0,2) 
		--
		(0,2.5)
		.. controls (0,3) and (.75,3) ..
		(.75,3.5);
	\draw (1,-1.5) node[below,scale=.75]{$g'$}
		--
		(1,-.5)
		.. controls (1,0) and (1.5,0) ..
		(1.5,.5);
	\draw (2,-3) node[below,scale=.75]{$g\phantom{'}$}
		--
		(2, -.5) 
		.. controls (2,0) and (1.5,0) ..
		(1.5,.5)
		--
		(1.5,1)
		--
		(1.5,.5)
		--
		(1.5,2.5)  
		.. controls (1.5,3) and (.75,3) .. 
		(.75,3.5);
        \node[fill=white,draw,rounded corners,scale=.75] at (0,2) {$3$};
        \node[fill=white,draw,rounded corners,scale=.75] at (1.5,.75) {$2\bullet1$};
}
\arrow[r, "\beta(g''\text{,~}g'g)"]
&
\tikz[yscale=.5,xscale=.75]{
	\draw (0,-.5) node[below,scale=.75]{$g''$}
		--
		(0,.5) 
		.. controls (0,1) and (.75,1) ..
		(.75,1.5);
        \draw (1.5,-1) -- (1.5,0.5);
        \draw (1,-2.5) node[below,scale=.75]{$g'$}
		--
		(1,-2)
		.. controls (1,-1.5) and (1.5,-1.5) ..
		(1.5,-1);
	\draw (2,-3) node[below,scale=.75]{$g\phantom{'}$}
		--
		(2, -2) 
		.. controls (2,-1.5) and (1.5,-1.5) ..
		(1.5,-1)
        --
        (1.5,0.5)
		.. controls (1.5,1) and (.75,1) .. 
		(.75,1.5)
		--
		(.75,3.5);
        \node[fill=white,draw,rounded corners,scale=.75] at (.75,2.15) {$3\bullet2\bullet1$};
}
\end{tikzcd}
\end{equation}

For $(\mathscr{G},\alpha)$, we will define the compatibility maps $\beta$ in a way analogous to the presentation in \cite{naisse2020odd}. Suppose 
\[
g' = (\vec{D}', \vec{p}') =  \prod_i(\vec{D}_i', \vec{p}_i') = \prod_{i,j} (D_{ij}, p_{ij})
\]
and 
\[
g = (\vec{D}, \vec{p}) = \prod_i (D_i, p_i)
\]
constitute a multipath of length two; $g\circ g' \in \mathscr{G}^{[2]}$. Again, denote by $P\in\mathbb{Z}\times \mathbb{Z}$ the sum of the entries of $\vec{p}$ and $P'$ the sum of the entries of $\vec{p}'$. Write $\vec{D} = (D_1,\ldots, D_k)$, and suppose that $(\vec{\Delta}, \vec{v}) = ((\Delta_1, v_1), \ldots, (\Delta_k, v_k))$ is a collection of cobordisms for $g$. We'll write
\[
(\vec{\Delta}, \vec{v})(\vec{D}, \vec{p}) = (\vec{\Delta}(\vec{D}), \vec{v} + \vec{p})
\]
where $\vec{\Delta}(\vec{D}) = (\Delta_1(D_1),\ldots, \Delta_k(D_k))$ and $\Delta_i(D_i)$ denotes the boundary of $\Delta_i$ other than $D_i$. Finally, $V$ and $V'\in \mathbb{Z}\times \mathbb{Z}$ will denote the sums of the entries of $\vec{v}$ and $\vec{v}'$ respectively.

The value $\beta$ will be defined as the product of four values. First, consider the change of chronology 
\[
\begin{tikzcd}[scale=1.25, row sep = huge]
& \vec{x} \vec{\Delta}(\vec{D}) \left(\vec{\Delta}'(\vec{D}')\right) \vec{z} & \\
\vec{x} \vec{D}(\vec{D}') \vec{z} \arrow[ur, "1_{\vec{x}} \left( \vec{\Delta}' \bullet \vec{\Delta} \right) 1_{\vec{z}}"] & \xRightarrow{H_\beta} & \vec{x} \vec{\Delta}'(\vec{D}') \vec{y} \otimes \vec{y} \vec{\Delta}(\vec{D}) \vec{z} \arrow[ul, "W_{\vec{x}\vec{y}\vec{z}} (\vec{\Delta'}(\vec{D}')\text{,~}\vec{\Delta}(\vec{D}))"'] \\
& \vec{x} \vec{D}' \vec{y} \otimes \vec{y} \vec{D} \vec{z} \arrow[ur, "1_{\vec{x}} \vec{\Delta'} 1_{\vec{y}} \sqcup 1_{\vec{y}} \vec{\Delta} 1_{\vec{z}}"']  \arrow[lu, "W_{\vec{x}\vec{y}\vec{z}} (\vec{D}'\text{,~}\vec{D})"]& 
\end{tikzcd}
\]
Set $\beta_1 = \iota(H_\beta)$. Then, set
\begin{align*}
    \beta_2 &= \lambda\left(~\abs{W_{\vec{x}\vec{y}\vec{z}}(\vec{\Delta}'(\vec{D}'), \vec{\Delta}(\vec{D}))}, V' + V \right), \\
    \beta_3 &= \lambda\left(~\abs{1_{\vec{x}}\vec{\Delta}' 1_{\vec{y}}}, V\right),~\text{and}\\
    \beta_4 &= \lambda\left(P',~ \abs{1_{\vec{y}} \vec{\Delta} 1_{\vec{z}}} + V\right).
\end{align*}
We define
\[
\beta = \beta_4\beta_3\beta_2\beta_1.
\]
Naisse and Putyra describe this shift diagrammatically as follows.
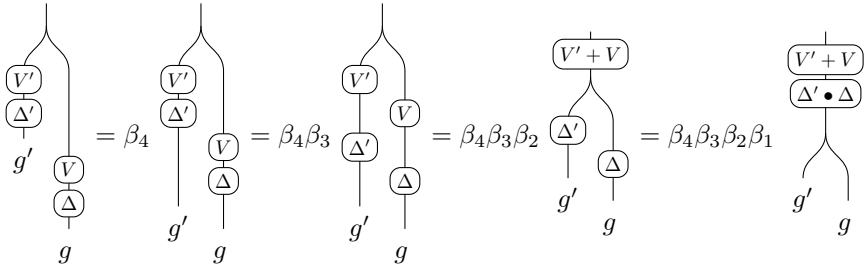
\[
\tikz[baseline={([yshift=-.5ex]current bounding box.center)}, scale=0.6, yscale=1]{
	\draw (0,1) node[below]{$g'$}
		--
		(0,1.75) node[near end,,fill=white,draw,rounded corners,scale=.75]{$\Delta'$}
		--
		(0,2.5) node[near end,,fill=white,draw,rounded corners,scale=.75]{$V'$}
		.. controls (0,3) and (.5,3) ..
		(.5,3.5) -- (.5,4);
	\draw (1,-1) node[below]{$g\phantom{'}$}
		--
		(1,-.25) node[near end,,fill=white,draw,rounded corners,scale=.75]{$\Delta$}
		-- 
		(1,.5) node[near end,,fill=white,draw,rounded corners,scale=.75]{$V$}
		-- 
		(1,2.5)
		.. controls (1,3) and (.5,3) ..
		(.5,3.5);
}
=
\beta_4 \ 
\tikz[baseline={([yshift=-.5ex]current bounding box.center)}, scale=0.6, yscale=1]{
	\draw (0,-.5) node[below]{$g'$}
		--
		(0,1)
		--
		(0,1.75) node[near end,,fill=white,draw,rounded corners,scale=.75]{$\Delta'$}
		--
		(0,2.5) node[near end,,fill=white,draw,rounded corners,scale=.75]{$V'$}
		.. controls (0,3) and (.5,3) ..
		(.5,3.5) -- (.5,4);
	\draw (1,-1) node[below]{$g\phantom{'}$}
		--
		(1,-.5)
		--
		(1,.25) node[near end,,fill=white,draw,rounded corners,scale=.75]{$\Delta$}
		-- 
		(1,1) node[near end,,fill=white,draw,rounded corners,scale=.75]{$V$}
		-- 
		(1,2.5)
		.. controls (1,3) and (.5,3) ..
		(.5,3.5);
}
=
\beta_4 \beta_3 \ 
\tikz[baseline={([yshift=-.5ex]current bounding box.center)}, scale=0.6, yscale=1]{
	\draw (0,-.5) node[below]{$g'$}
		--
		(0,.25)
		--
		(0,1)node[near end,,fill=white,draw,rounded corners,scale=.75]{$\Delta'$}
		--
		(0,1.75) 
		--
		(0,2.5) node[near end,,fill=white,draw,rounded corners,scale=.75]{$V'$}
		.. controls (0,3) and (.5,3) ..
		(.5,3.5) -- (.5,4);
	\draw (1,-1) node[below]{$g\phantom{'}$}
		--
		(1,-.5)
		--
		(1,.25) node[near end,,fill=white,draw,rounded corners,scale=.75]{$\Delta$}
		-- 
		(1,1) 
		-- 
		(1,1.75) node[near end,,fill=white,draw,rounded corners,scale=.75]{$V$}
		-- 
		(1,2.5)
		.. controls (1,3) and (.5,3) ..
		(.5,3.5);
}
=
\beta_4 \beta_3 \beta_2 \
 \tikz[baseline={([yshift=-.5ex]current bounding box.center)}, scale=0.6, yscale=1]{
	\draw (0,-.75) node[below]{$g'$}
		--
		(0,-.25)
		--
 		(0,.5)  node[near end,,fill=white,draw,rounded corners,scale=.75]{$\Delta'$}
		.. controls (0,1) and (.5,1) ..
		(.5,1.5)  
		--
		(.5,2.5)  node[midway,,fill=white,draw,rounded corners,scale=.75]{$V' + V$};
	\draw (1,-1.25) node[below]{$g\phantom{'}$} 
		--
		(1,-1)
		--
		(1,-.25)  node[near end,,fill=white,draw,rounded corners,scale=.75]{$\Delta$}
		-- 
		(1,.5)
		.. controls (1,1) and (.5,1) ..
		(.5,1.5);
}
=
\beta_4\beta_3\beta_2\beta_1 \
 \tikz[baseline={([yshift=-.5ex]current bounding box.center)}, scale=0.6, yscale=1]{
	\draw (0,-1.5) node[below]{$g'$}
		.. controls (0,-1) and (.5,-1) ..
		(.5,-.5)  
		--
		(.5,0)
		--
		(.5,.75) node[midway,,fill=white,draw,rounded corners,scale=.75]{$\Delta'\bullet \Delta$}
		--
		(.5,1.5)  node[midway,,fill=white,draw,rounded corners,scale=.75]{$V'+V$}
		--
		(.5,1.75);
	\draw (1,-2) node[below]{$g\phantom{'}$} 
		--
		(1,-1.5)
		.. controls (1,-1) and (.5,-1) ..
		(.5,-.5);
}
\]

\begin{lemma}
\label{lem:computational2}
If $D \in \mathscr{D}_{(m_1,\ldots, m_k; n)}$ and $\vec{D} = (D_1,\ldots, D_k)$ for $D_i \in \mathscr{D}_{(\ell_{i1},\ldots, \ell_{i\alpha_i}; m_i)}$. Moreover, let $\Delta$ be a cobordism on $D$ and $\vec{\Delta}'$ be a collection of cobordisms for $\vec{D}'$. Then, for any $\vec{x} \in \prod_{i,j} B_{\ell_{ij}}$, $\vec{y} \in \prod_i B_{m_i}$ and $z\in B_n$, we have that
\[
\abs{1_{\vec{y}} \Delta 1_z} + \abs{1_{\vec{x}} \vec{\Delta}' 1_{\vec{y}}} + \abs{W_{\vec{x}\vec{y}z} (\vec{\Delta}'(\vec{D}'), \Delta(D))} = \abs{W_{\vec{x}\vec{y}z} (\vec{D}', D)} + \abs{1_{\vec{x}} (\vec{\Delta}' \bullet \Delta) 1_z}
\]
\end{lemma}

Notice that this lemma immediately applies to the cases when $D$ is actually a collection $\vec{D} = \prod_i D_i$ and $\vec{D}' = \prod_{i,j} D_{ij}$.

\begin{proof}
Exactly as in \cite{naisse2020odd}, there is a diffeomorphism between the cobordisms below, so they must have the same degree.
\[
\tikz[yscale=.5,xscale=.75, baseline=1ex]{
	\draw (1,-1.5) node[below]{$\vec{D}'$}
		--
		(1,-.5)
		-- 
		(1,.5)
		--
		(1,1)
		.. controls (1,1.5) and (1.5,1.5) ..
		(1.5,2);
	\draw (2,-1.5) node[below]{$\vec{D}$}
		--
		(2, -.5)
		--
		(2,1)
		.. controls (2,1.5) and (1.5,1.5) ..
		(1.5,2)
		--
		(1.5,2.5);
        \node[fill=white,draw,rounded corners,scale=.8] at (1,.3) {$\vec{\Delta}'$};
        \node[fill=white,draw,rounded corners,scale=.8] at (2, -.7) {$\vec{\Delta}$};
}
\ \cong \ 
\tikz[yscale=.5,xscale=.75, baseline=1ex]{
	\draw (1,-1.5) node[below]{$\vec{D}'$}
		--
		(1,-1)
		.. controls (1,-.5) and (1.5,-.5) ..
		(1.5,0);
	\draw (2,-1.5) node[below]{$\vec{D}$}
		--
		(2, -1)   
		.. controls (2,-.5) and (1.5,-.5) ..
		(1.5,0)
		--
		(1.5,2.5);
        \node[fill=white,draw,rounded corners,scale=.8] at (1.5,1.3){$\vec{\Delta}' \bullet \vec{\Delta}$};
}
\]
\end{proof}

\begin{proposition}
\label{prop:comp1}
The $\mathscr{G}$-shifting system of Proposition \ref{prop:G-shifting} is compatible with the associator of Proposition \ref{prop:associator} through the compatibility map $\beta$ defined above.
\end{proposition}

\begin{proof}
The following follows closely the proof of \cite{naisse2020odd}. We will show that $\beta$ satisfies equation (\ref{eq:shiftingsystcomp}). The first step is to document the contributions of $\alpha_1$ and $\beta_1$. To do this, consider the following diagram of cobordisms. 
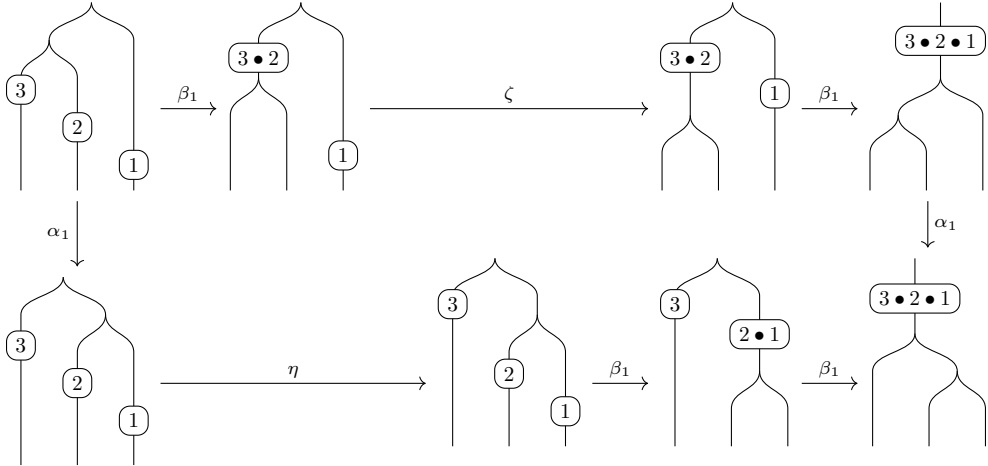
\[
\begin{tikzcd}[column sep=5ex, scale=1.1]
\tikz[yscale=.5,xscale=.75]{
	\draw (0,-2)
		--
		(0,0)
		--
		(0,1)   
		.. controls (0,1.5) and (.5,1.5) ..
		(.5,2)
		.. controls (.5,2.5) and (1.25,2.5) ..
		(1.25,3);
	\draw (1,-2)
		--
		(1,-1)
		-- 
		(1,0)  
		--
		(1,1)
		.. controls (1,1.5) and (.5,1.5) ..
		(.5,2);
	\draw (2,-2)
		--
		(2, -1)   
		--
		(2,2)
		.. controls (2,2.5) and (1.25,2.5) ..
		(1.25,3);
        \node[fill=white,draw,rounded corners,scale=.8] at (0,0.5) {$3$};
        \node[fill=white,draw,rounded corners,scale=.8] at (1,-0.5) {$2$};
        \node[fill=white,draw,rounded corners,scale=.8] at (2,-1.5) {$1$};
}
\ar{r}{\beta_1}
\ar[swap]{d}{\alpha_1}
&
\tikz[yscale=.5,xscale=.75]{
	\draw (0,-2)
		--
		(0,0) 
		.. controls (0,.5) and (.5,.5) ..
		(.5,1)
		--
		(.5,2)
		.. controls (.5,2.5) and (1.25,2.5) ..
		(1.25,3);
	\draw (1,-2)
		-- 
		(1,0)   
		--
		(1,0)
		.. controls (1,.5) and (.5,.5) ..
		(.5,1);
	\draw (2,-2)
		--
		(2, -1)
		--
		(2,2)
		.. controls (2,2.5) and (1.25,2.5) ..
		(1.25,3);
        \node[fill=white,draw,rounded corners,scale=.8] at (0.5,1.35) {$3 \bullet 2$};
        \node[fill=white,draw,rounded corners,scale=.8] at (2,-1.25) {$1$};
}
\ar{rr}{\zeta}
&
\hspace{2ex}
&
\tikz[yscale=.5,xscale=.75]{
	\draw (0,-2)
		--
		(0,-1) 
		.. controls (0,-.5) and (.5,-.5) ..
		(.5,0)
		--
		(.5,1)
		--
		(.5,2)
		.. controls (.5,2.5) and (1.25,2.5) ..
		(1.25,3);
	\draw (1,-2)
		--
		(1,-1)
		.. controls (1,-.5) and (.5,-.5) ..
		(.5,0);
	\draw (2,-2)
		--
		(2, 0) 
		--
		(2,1)
		--
		(2,2)
		.. controls (2,2.5) and (1.25,2.5) ..
		(1.25,3);
        \node[fill=white,draw,rounded corners,scale=.8] at (0.5,1.35) {$3 \bullet 2$};
        \node[fill=white,draw,rounded corners,scale=.8] at (2,0.4) {$1$};
}
\ar{r}{\beta_1}
&
\tikz[yscale=.5,xscale=.75]{
	\draw (0,-2)
		--
		(0,-1) 
		.. controls (0,-.5) and (.5,-.5) ..
		(.5,0)
		.. controls (.5,.5) and (1.25,.5) ..
		(1.25,1)
		--
		(1.25,3);
	\draw (1,-2)
		--
		(1,-1)
		.. controls (1,-.5) and (.5,-.5) ..
		(.5,0);
	\draw (2,-2)
		--
		(2, 0)
		.. controls (2,.5) and (1.25,.5) ..
		(1.25,1);
        \node[fill=white,draw,rounded corners,scale=.8] at (1.25,1.8) {$3 \bullet 2 \bullet 1$};
}
\ar{d}{\alpha_1}
\\
\tikz[yscale=.5,xscale=.75]{
	\draw (0,-2)
		--
		(0,1)
		--
		(0,1)
		--
		(0,2)
		.. controls (0,2.5) and (.75,2.5) ..
		(.75,3);
	\draw (1,-2)
		--
		(1,-1)
		-- 
		(1,0)
		--
		(1,1)
		.. controls (1,1.5) and (1.5,1.5) ..
		(1.5,2);
	\draw (2,-2)
		--
		(2, -1)
		--
		(2,1)
		.. controls (2,1.5) and (1.5,1.5) ..
		(1.5,2)
		.. controls (1.5,2.5) and (.75,2.5) .. 
		(.75,3);
        \node[fill=white,draw,rounded corners,scale=.8] at (0,1) {$3$};
        \node[fill=white,draw,rounded corners,scale=.8] at (1,0) {$2$};
        \node[fill=white,draw,rounded corners,scale=.8] at (2,-1) {$1$};
}
\ar{rr}{\eta}
&
&
\tikz[yscale=.5,xscale=.75]{
	\draw (0,-1.5)
		--
		(0,1) 
		--
		(0,2.5)
		.. controls (0,3) and (.75,3) ..
		(.75,3.5);
	\draw (1,-1.5)
		--
		(1,-.5)
		-- 
		(1,.5)
		--
		(1,1)
		.. controls (1,1.5) and (1.5,1.5) ..
		(1.5,2);
	\draw (2,-1.5)
		--
		(2, -.5)
		--
		(2,1)
		.. controls (2,1.5) and (1.5,1.5) ..
		(1.5,2)
		--
		(1.5,2.5)
		.. controls (1.5,3) and (.75,3) .. 
		(.75,3.5);
        \node[fill=white,draw,rounded corners,scale=.8] at (0,2.1) {$3$};
        \node[fill=white,draw,rounded corners,scale=.8] at (1,.25) {$2$};
        \node[fill=white,draw,rounded corners,scale=.8] at (2,-0.75) {$1$};
}
\ar{r}{\beta_1}
&
\tikz[yscale=.5,xscale=.75]{
	\draw (0,-1.5)
		--
		(0,1) 
		--
		(0,2.5)
		.. controls (0,3) and (.75,3) ..
		(.75,3.5);
	\draw (1,-1.5)
		--
		(1,-.5)
		.. controls (1,0) and (1.5,0) ..
		(1.5,.5);
	\draw (2,-1.5)
		--
		(2, -.5) 
		.. controls (2,0) and (1.5,0) ..
		(1.5,.5)
		--
		(1.5,1)
		--
		(1.5,1.5)
		--
		(1.5,2.5)  
		.. controls (1.5,3) and (.75,3) .. 
		(.75,3.5);
        \node[fill=white,draw,rounded corners,scale=.8] at (0,2.1) {$3$};
        \node[fill=white,draw,rounded corners,scale=.8] at (1.5, 1.3) {$2 \bullet 1$};
}
\ar{r}{\beta_1}
&
\tikz[yscale=.5,xscale=.75]{
	\draw (0,-1.5)
		--
		(0,.5) 
		.. controls (0,1) and (.75,1) ..
		(.75,1.5);
	\draw (1,-1.5)
		--
		(1,-.5)
		.. controls (1,0) and (1.5,0) ..
		(1.5,.5);
	\draw (2,-1.5)
		--
		(2, -.5) 
		.. controls (2,0) and (1.5,0) ..
		(1.5,.5)  
		.. controls (1.5,1) and (.75,1) .. 
		(.75,1.5)
		--
		(.75,3.5);
        \node[fill=white,draw,rounded corners,scale=.8] at (0.75,2.25) {$3 \bullet 2 \bullet 1$};
}
\end{tikzcd}
\]
where
\[
\zeta = \lambda\left(~\abs{W_{\vec{w}\vec{x}\vec{y}}(\vec{D}''\text{,~}\vec{D}')}\text{,~}\abs{1_{\vec{y}}\vec{\Delta}1_{\vec{z}}}\right) 
\qquad \text{and} \qquad
\eta = \lambda\left(~\abs{W_{\vec{x}\vec{y}\vec{z}}(\vec{\Delta}'(\vec{D}')\text{,~} \vec{\Delta}(\vec{D}))}\text{,~} \abs{1_{\vec{w}}\vec{\Delta}'' 1_{\vec{x}}} \right).
\]
Comparing with diagram (\ref{beta_CompatCondition}), we see that
\[
\zeta \times \text{(Left of (\ref{eq:shiftingsystcomp}))} = \eta \times \text{(Right of (\ref{eq:shiftingsystcomp}))}.
\]

Next, we have to compare the contributions of $\alpha_2, \beta_2, \beta_3$, and $\beta_4$. First, for the left side of (\ref{eq:shiftingsystcomp}) (or, the top-and-then-down path in (\ref{beta_CompatCondition})), we have the following.
\begin{align*}
& \lambda \left(~\abs{W_{\vec{w}\vec{x}\vec{y}} (\vec{\Delta}''(\vec{D}''), \vec{\Delta}'(\vec{D}')}, V'' + V' \right) \cdot \lambda\left(~\abs{1_{\vec{w}} \vec{\Delta}'' 1_{\vec{y}}}, V'\right) \cdot \lambda\left(P'', \abs{1_{\vec{x}}\vec{\Delta}' 1_{\vec{y}}} + V'\right) & \beta_{2,3,4} (g'',g')\\
& \cdot \lambda\left(~\abs{W_{\vec{w}\vec{y}\vec{z}}((\vec{\Delta}''\bullet \vec{\Delta}')(\vec{D}'(\vec{D}'')), \vec{\Delta}(\vec{D}))}, V'' + V' + V\right) \cdot \underbrace{\lambda\left(~\abs{1_{\vec{w}} (\vec{\Delta}''\bullet\vec{\Delta}')1_{\vec{y}}}, V\right)}_{(*)} &  \\ & \qquad \cdot \underbrace{\lambda\left(P' + P'' + \abs{W_{\vec{w}\vec{x}\vec{y}}(\vec{D}'', \vec{D}')}, \abs{1_{\vec{y}} \vec{\Delta} 1_{\vec{z}}} + V \right)}_{(**)} & \beta_{2,3,4} (g''g', g) \\ 
& \cdot \lambda \left(~\abs{W_{\vec{x}\vec{y}\vec{z}} (\vec{D}', \vec{D})}, P''\right) & \alpha_2(g'', g', g)
\end{align*}
Turn your attention to the term marked $(*)$. By application of Lemma \ref{lem:computational2}, we can write
\begin{align*}
\lambda\left(~\abs{1_{\vec{w}} (\vec{\Delta}''\bullet\vec{\Delta}')1_{\vec{y}}}, V\right) &= \lambda\left(~\abs{1_{\vec{x}} \vec{\Delta}' 1_{\vec{y}}} + \abs{1_{\vec{w}} \vec{\Delta}'' 1_{\vec{x}}}, V\right) \cdot \lambda\left(~\abs{W_{\vec{w}\vec{x}\vec{y} } (\vec{\Delta}''(\vec{D}''), \vec{\Delta}'(\vec{D}'))}, V\right) \\ & \cdot \lambda\left(~\abs{W_{\vec{w}\vec{x}\vec{y}}(\vec{D}'', \vec{D}')}, V\right)^{-1}.
\end{align*}
On the other hand, using linearity, we can rewrite the term marked $(**)$ as
\begin{align*}
\lambda\left(P' + P'' + \abs{W_{\vec{w}\vec{x}\vec{y}}(\vec{D}'', \vec{D}')}, \abs{1_{\vec{y}} \vec{\Delta} 1_{\vec{z}}} + V\right) &= \lambda\left(P' + P'' + \abs{W_{\vec{w}\vec{x}\vec{y}}(\vec{D}'', \vec{D}')}, \abs{1_{\vec{y}} \vec{\Delta} 1_{\vec{z}}}\right) \\ & \cdot \lambda\left(P' + P'', V\right) \cdot \lambda\left(~\abs{W_{\vec{w}\vec{x}\vec{y}}(\vec{D}'', \vec{D}')}, V\right)
\end{align*}
The last terms in the past two expansions cancel, and we can rewrite the contributions of $\alpha_2, \beta_2, \beta_3$, and $\beta_4$ on the left side of (\ref{eq:shiftingsystcomp}) as 
\begin{align*}
& \lambda \left(~\abs{W_{\vec{w}\vec{x}\vec{y}} (\vec{\Delta}''(\vec{D}''), \vec{\Delta}'(\vec{D}')}, V'' + V' +V \right) \cdot \lambda\left(~\abs{1_{\vec{w}} \vec{\Delta}'' 1_{\vec{y}}}, V'\right) \cdot \lambda\left(P'', \abs{1_{\vec{x}}\vec{\Delta}' 1_{\vec{y}}} + V'\right)\\
& \cdot \lambda\left(~\abs{W_{\vec{w}\vec{y}\vec{z}}((\vec{\Delta}''\bullet \vec{\Delta}')(\vec{D}'(\vec{D}'')), \vec{\Delta}(\vec{D}))}, V'' + V' + V\right) \cdot \lambda\left(\abs{1_{\vec{x}} \vec{\Delta}' 1_{\vec{y}}} + \abs{1_{\vec{w}} \vec{\Delta}'' 1_{\vec{x}}}, V~\right) \\ & \qquad \cdot \lambda\left(P' + P'', V\right) \cdot \lambda\left( P' + P'', \abs{1_{\vec{y}} \vec{\Delta} 1_{\vec{z}}}\right) \cdot \underbrace{\lambda\left(~\abs{W_{\vec{w}\vec{x}\vec{y}}(\vec{D}'', \vec{D}') }, \abs{1_{\vec{y}} \vec{\Delta} 1_{\vec{z}}}\right)}_{\zeta}  \\ 
& \cdot \lambda \left(~\abs{W_{\vec{x}\vec{y}\vec{z}} (\vec{D}', \vec{D})}, P''\right).
\end{align*}
The process above could be described as ``simplifying'' $\beta_{2,3,4}(g''g', g)$.

Likewise, on the right side of (\ref{eq:shiftingsystcomp}) (or, the down-and-then-bottom path in (\ref{beta_CompatCondition})), we have the following.
\begin{align*}
& \underbrace{\lambda\left(~\abs{W_{\vec{x}\vec{y}\vec{z}} (\vec{\Delta}'(\vec{D}'), \vec{\Delta}(\vec{D}))}, P'' + V'' + \abs{1_{\vec{w}} \vec{\Delta}'' 1_{\vec{x}}}\right)}_{(\star)}  & \alpha_2(\varphi(g''), \varphi(g'), \varphi(g)) \\
&\cdot \underbrace{\lambda\left(~\abs{W_{\vec{x}\vec{y}\vec{z}}(\vec{\Delta}'(\vec{D}'), \vec{\Delta}(\vec{D}))}, V' + V \right)}_{(\star\star)} \cdot \lambda\left(~\abs{1_{\vec{x}}\vec{\Delta}' 1_{\vec{y}}}, V\right).
& \\
& \qquad \cdot \lambda\left(P',~ \abs{1_{\vec{y}} \vec{\Delta} 1_{\vec{z}}} + V\right).
& \beta_{2,3,4}(g', g) \\
& \cdot \lambda\left(~ \abs{W_{\vec{w}\vec{x}\vec{z}}(\vec{\Delta}''(\vec{D}''), (\vec{\Delta}' \bullet \vec{\Delta}) (\vec{D}(\vec{D}')))}, V'' + V' + V \right) & \\
& \qquad \cdot \lambda\left(~\abs{1_{\vec{w}} \vec{\Delta}'' 1_{\vec{x}}}, V + V' \right) \cdot \underbrace{\lambda\left(P'', \abs{1_{\vec{x}} (\vec{\Delta}' \bullet \vec{\Delta}) 1_{\vec{z}}} + V + V'\right)}_{(\star\star\star)} & \beta_{2,3,4}(g'', g'g)
\end{align*}
Notice that the term $(\star)$ has three ``parts.'' The $V''$ part can be absorbed into the term $(\star\star)$; the rest can be written
\[
\lambda\left(~\abs{W_{\vec{x}\vec{y}\vec{z}} (\vec{\Delta}'(\vec{D}'), \vec{\Delta}(\vec{D}))}, P''\right) \cdot \lambda\left(~\abs{W_{\vec{x}\vec{y}\vec{z}} (\vec{\Delta}'(\vec{D}'), \vec{\Delta}(\vec{D}))}, \abs{1_{\vec{w}} \vec{\Delta}'' 1_{\vec{x}}}\right).
\]
The $(\star\star\star)$ term decomposes into parts
\[
\lambda\left(P'', \abs{1_{\vec{x}} (\vec{\Delta}' \bullet \vec{\Delta}) 1_{\vec{z}}}\right) \cdot \lambda\left(P'', V + V'\right).
\]
Again, applying Lemma \ref{lem:computational2}, we can write
\begin{align*}
\lambda\left(P'', \abs{1_{\vec{x}} (\vec{\Delta}' \bullet \vec{\Delta}) 1_{\vec{z}}}\right) &= \lambda\left(P'', \abs{1_{\vec{y}} \vec{\Delta} 1_{\vec{z}}} + \abs{1_{\vec{x}} \vec{\Delta}' 1_{\vec{y}}}\right) \cdot \lambda\left(P'', \abs{W_{\vec{x}\vec{y}\vec{z}}(\vec{\Delta}'(\vec{D}'), \vec{\Delta}(\vec{D}))}\right) \\ & \cdot \lambda\left(P'', -\abs{W_{\vec{x}\vec{y}\vec{z}}(\vec{D}', \vec{D})}\right).
\end{align*}
The middle term after the equality cancels with first term in the rewriting of $(\star)$. The last term can be rewritten as $\lambda\left(~\abs{W_{\vec{x}\vec{y}\vec{z}}(\vec{D}', \vec{D})}, P''\right)$. All together, this means that we can rewrite the contributions of $\alpha_2$, $\beta_2$, $\beta_3$, and $\beta_4$ on the right side of (\ref{eq:shiftingsystcomp}) as 
\begin{align*}
& \underbrace{\lambda\left(~\abs{W_{\vec{x}\vec{y}\vec{z}} (\vec{\Delta}'(\vec{D}'), \vec{\Delta}(\vec{D}))}, \abs{1_{\vec{w}} \vec{\Delta}'' 1_{\vec{x}}}\right)}_{\eta} \\
&\cdot \lambda\left(~\abs{W_{\vec{x}\vec{y}\vec{z}}(\vec{\Delta}'(\vec{D}'), \vec{\Delta}(\vec{D}))}, V'' + V' + V \right) \cdot \lambda\left(~\abs{1_{\vec{x}}\vec{\Delta}' 1_{\vec{y}}}, V\right) \cdot \lambda\left(P',~ \abs{1_{\vec{y}} \vec{\Delta} 1_{\vec{z}}} + V\right). \\
& \cdot \lambda\left(~ \abs{W_{\vec{w}\vec{x}\vec{z}}(\vec{\Delta}''(\vec{D}''), (\vec{\Delta}' \bullet \vec{\Delta}) (\vec{D}(\vec{D}')))}, V'' + V' + V \right) \cdot \lambda\left(~\abs{1_{\vec{w}} \vec{\Delta}'' 1_{\vec{x}}}, V + V' \right) \\
& \qquad \cdot \lambda\left(P'', V + V'\right) \cdot \lambda\left(P'', \abs{1_{\vec{y}} \vec{\Delta} 1_{\vec{z}}} + \abs{1_{\vec{x}} \vec{\Delta}' 1_{\vec{y}}}\right) \cdot \lambda\left(~\abs{W_{\vec{x}\vec{y}\vec{z}}(\vec{D}', \vec{D})}, P''\right)
\end{align*}

Compare the simplifications of contributions from each side. One one hand, 
\[
\lambda \left(~\abs{W_{\vec{w}\vec{x}\vec{y}} (\vec{\Delta}''(\vec{D}''), \vec{\Delta}'(\vec{D}')}, V'' + V' +V \right) \cdot \lambda\left(~\abs{W_{\vec{w}\vec{y}\vec{z}}((\vec{\Delta}''\bullet \vec{\Delta}')(\vec{D}'(\vec{D}'')), \vec{\Delta}(\vec{D}))}, V'' + V' + V\right)
\]
is equal to
\[
\lambda\left(~\abs{W_{\vec{x}\vec{y}\vec{z}}(\vec{\Delta}'(\vec{D}'), \vec{\Delta}(\vec{D}))}, V'' + V' + V \right) \cdot \lambda\left(~ \abs{W_{\vec{w}\vec{x}\vec{z}}(\vec{\Delta}''(\vec{D}''), (\vec{\Delta}' \bullet \vec{\Delta}) (\vec{D}(\vec{D}')))}, V'' + V' + V \right)
\]
by Lemma \ref{lemma:degree of multipaths}. On the other hand, careful observation reveals that, via bilinearity of $\lambda$ alone, the two collections of terms apart from these, and the terms marked $\zeta$ and $\eta$, are equivalent. The conclusion is that 
\[
\eta \times \text{(Left of (\ref{eq:shiftingsystcomp}))} = \zeta \times \text{(Right of (\ref{eq:shiftingsystcomp}))}.
\]
This completes the proof.
\end{proof}

As we proceed, for simplicity of exposition (and because it is the only situation which matters in our application) we will only consider multipaths which end in a single multimorphism; we have shown in the previous arguments how the situation is generalized without problem.

\subsection{Generalities on shifting systems for grading multicategories}
\label{ss:generalshifting}

We conclude this discussion by detailing the generalities of $\mathscr{C}$-shifting systems. These are results of \cite{naisse2020odd} which lift to the setting of grading multicategories. Throughout, let $\mathscr{C}$ be a grading multicategory with associator $\alpha$, and $S = \{\mathscr{I}, \{\varphi_i\}_{i\in\mathscr{I}}\}$ a $\mathscr{C}$-shifting system compatible with $\alpha$ through compatibility maps $\beta$.

Just as in the non-multi setting, we define for each $i \in \mathscr{I}$ a \textit{grading shift functor} $\varphi_i: \mathrm{Mod}^\mathscr{C} \to \mathrm{Mod}^\mathscr{C}$ by putting
\[
\varphi_i(M) = \bigoplus_{g\in \mathsf{D}_i} \varphi_i(M)_{\varphi_i(g)}
\]
for $\varphi_i(M)_{\varphi_i(g)}:= M_g$; that is, $\varphi_i$ sends elements in degree $g\in \mathsf{D}_i$ to elements in degree $\varphi_i(g)$, and elements whose degree does not belong to $\mathsf{D}_i$ to zero. Sometimes we call $\varphi_i$ a \textit{$\mathscr{C}$-grading shift} or just a \textit{grading shift}.

Now, if $M, M_1, \ldots, M_k$ are $\mathscr{C}$-graded modules, there is a canonical isomorphism
\begin{equation}
\label{eq:compatibilityCI}
\begin{split}
\beta_{(j_1,\ldots, j_k), i}:& \left(\varphi_{j_1}(M_1), \ldots, \varphi_{j_k}(M_k)\right) \otimes \varphi_i(M) \xrightarrow{\sim} \varphi_{(j_1,\ldots, j_k)\bullet i} \left((M_1,\ldots, M_k) \otimes M\right)\\
\text{by}~&(m_1, \ldots, m_k) \otimes m \mapsto \beta_{(j_1,\ldots, j_k), i}(\abs{\vec{m}}, \abs{m}) (m_1,\ldots, m_k) \otimes m.
\end{split}
\end{equation}
The compatibility requirement, equation (\ref{eq:shiftingsystcomp}), ensures that this isomorphism is compatible with the coherence isomorphism given by $\alpha$. Moreover, since grading shift functors do not have effect on graded maps, the compatibility maps $\beta_{\vec{j}, i}$ define natural isomorphisms (denoted by the same symbol) of multifunctors 
\begin{equation}
\label{eq:compatibilityNI}
\beta_{\vec{j}, i}: \left(\varphi_{j_1}(-), \ldots, \varphi_{j_k}(-)\right) \otimes \varphi_i(-) \xrightarrow{\sim} \varphi_{\vec{j}, i} \left((-,\ldots, -), -\right)
\end{equation}
for all $j_1,\ldots, j_k, i\in \mathscr{I}$.

We define the \textit{identity shift functor} $\varphi_{\mathrm{Id}}$ as $\bigoplus_{i\in \mathscr{I}_{\mathrm{Id}}} \varphi_i$; thus, $\varphi_{\mathrm{id}}(M) \cong M$. In general, the identity shift and the neutral shift are not the same (see, for example, \cite{naisse2020odd}, Remark 4.10). We'll consider the set $\widetilde{\mathscr{I}}$, defined to be $\mathscr{I} \sqcup \{\mathrm{Id}\}$. We do not think of  $\widetilde{\mathscr{I}}$ as a multimonoid---writing it this way just helps to simplify notation. For example, we will write $\varphi_{j\bullet \mathrm{Id}}$ to mean $\bigoplus_{i\in \mathscr{I}_{\mathrm{Id}}}  \varphi_{j\bullet i}$. Similarly, $\varphi_{\vec{\mathrm{Id}}\bullet i}$ means $\bigoplus_{\vec{j}\in J} \varphi_{\vec{j} \bullet i}$ where $J = \{(j_1,\ldots, j_k): j_\ell \in \mathscr{I}_{\mathrm{Id}}~\text{for all}~\ell=1,\ldots, k\}$. To extend the compatibility maps $\beta$ to $\widetilde{\mathscr{I}}$, define $\beta_{\vec{\mathrm{Id}}, i}(g', g) = \beta_{\vec{j}, i}(g', g)$ where $g' \in \mathsf{D}_{\vec{j}}$ and $\vec{j} \in \mathscr{I}_{\mathrm{Id}}$; similarly $\beta_{j, \mathrm{Id}}(g', g) = \beta_{j,i}(g', g)$ where $g\in \mathsf{D}_i$, $i\in \mathscr{I}_{\mathrm{Id}}$. Lastly, we fix $\beta_{\vec{\mathrm{Id}}, \mathrm{Id}} = 1$. 

\subsubsection{Shifting multimodules}

To continue in the general setting, we must make the following assumption.

\begin{tcolorbox}[breakable, enhanced,colback=yellow!5!white,boxrule=0pt,frame hidden,
borderline west={.5mm}{0mm}{black}]
\textit{Assumption}: From hereafter, all $\mathscr{C}$-graded algebras $A$ are supported only in $\Sigma$; that is, $A_g = 0$ whenever $g\not\in \mathrm{Hom}_{\Sigma}$
\end{tcolorbox}

Thus, for $\mathscr{C}$-algebras $A$ which satisfy this assumption, we have that $\varphi_e(A) \cong A$ (really, $\varphi_e(A) = A$, since $\varphi_e$ acts as the identity wherever defined). Recall that, since
\[
H^n = \mathcal{F}(1_n) = \bigoplus_{a,b\in B^n} \mathcal{F}(a 1_n \overline{b})
\]
any $m\in H^n$ has degree $\deg_\mathscr{G}(m) = (1_n, \deg_R(m)): a \to b$; that is, arc algebras are $\mathscr{G}$-graded algebras supported only in $\Sigma$. 

If $M$ is a $\mathscr{C}$-graded $(A_1,\ldots, A_k; B)$-multimodule, and $\varphi_i$ is a $\mathscr{C}$-grading shifting functor, then we can view $\varphi_i(M)$ as a $\mathscr{C}$-graded $(A_1,\ldots, A_k; B)$-multimodule by defining left- and right-acitons
\[
  \begin{split}
    &\varphi_i\rho_L: (A_1,\ldots, A_k) \otimes \varphi_i(M) \to \varphi_i(M)\\
    \text{by}~&\varphi_i\rho_L(\vec{a}, \varphi_i(m)) = \beta_{(e,\ldots, e), i}(\abs{\vec{a}}, \abs{m}) \varphi_i(\rho_L(\vec{a},m))
  \end{split}
\quad \text{and} \quad
  \begin{split}
    &\varphi_i\rho_R: \varphi_i(M) \otimes B \to \varphi_i(M)\\
    \text{by}~&\varphi_i\rho_R(\varphi_i(m),b) = \beta_{i,e}(\abs{m}, \abs{b}) \varphi_i(\rho_R(m,b))
  \end{split}
\]
In other words, $\varphi_i\rho_L$ and $\varphi_i\rho_R$ are defined as the composites
\[
\begin{tikzcd}
(A_1,\ldots, A_k) \otimes \varphi_i(M) \arrow[rr, dotted] \arrow[d, "\star"] & & \varphi_i(M) \arrow[dd, equals] \\
\left(\varphi_e(A_1), \ldots, \varphi_e(A_k)\right) \otimes \varphi_i(M) \arrow[d, "\beta_{\vec{e},i}"]& & \\
\varphi_{(e,\ldots,e)\bullet i}\left((A_1,\ldots, A_k)\otimes M\right) \arrow[rr, "\rho_L"] & & \varphi_{(e,\ldots,e)\bullet i}(M)
\end{tikzcd}
\quad \text{and} \quad
\begin{tikzcd}
\varphi_i(M) \otimes B \arrow[rr, dotted] \arrow[d, "\star"] & & \varphi_i(M) \arrow[dd, equals] \\
\varphi_i(M) \otimes \varphi_e(B) \arrow[d, "\beta_{i, e}"]& & \\
\varphi_{i \bullet e}\left(M \otimes B\right) \arrow[rr, "\rho_R"] & & \varphi_{i\bullet e}(M)
\end{tikzcd}
\]
where the maps labeled $\star$ are isomorphisms thanks to the assumption from the start of the section.

We'll breifly describe why $\varphi_i(M)$ is indeed a $\mathscr{C}$-graded multimodule. First, notice that $\varphi_i\rho_L$ and $\varphi_i\rho_R$ are both graded maps. To illustrate for the left action, if $(a_1,\ldots, a_k) \otimes m$ has grading $g \circ (g_1,\ldots, g_k)$ in $(A_1,\ldots, A_k) \otimes M$, it has grading $\varphi_i(g) \circ (g_1, \ldots, g_k)$ in $(A_1,\ldots, A_k) \otimes \varphi_i(M)$. Thanks to the assumption from the start of the section, $g_i = \varphi_e(g_i)$ since all algebras in sight are supported only in $\Sigma$, and thus $\varphi_e$ acts as the identity map. Applying the natural isomorphism (\ref{eq:compatibilityNI}) provides the desired result. To see that requirements (i)-(iv) of the definition of $\mathscr{C}$-graded multimodules holds, one must simply apply equation (\ref{eq:shiftingsystcomp}) and $\beta_{(e,\ldots,e),(e,\ldots,e)} = 1$ in each of the scenarios.

Thus, grading shift functors are also functors for categories of multimodules. In conclusion, we have the following.

\begin{proposition}
Suppose $M\in \mathrm{MultiMod}^\mathscr{C}(B_1,\ldots, B_k; C)$ and $M_i \in \mathrm{MultiMod}^\mathscr{C}(A_{i1}, \ldots, A_{i\alpha_i}; B_i)$ for each $i=1,\ldots, k$. Then, for each $i, j_1, \ldots, j_k \in \mathscr{I}$, there is an isomorphism of $\mathscr{C}$-graded $(A_{11},\ldots, A_{k\alpha_k})$-multimodules
\[
\beta_{(j_1,\ldots, j_k), i}: \left(\varphi_{j_1}(M_1),\ldots, \varphi_{j_k}(M_k)\right) \otimes_{(B_1,\ldots, B_k)} \varphi_i(M) \xrightarrow{\sim} \varphi_{(j_1,\ldots, j_k)\bullet i}\left((M_1,\ldots, M_k) \otimes_{(B_1,\ldots, B_k)} M\right)
\]
induced by the canonical isomorphism (\ref{eq:compatibilityCI}).
\end{proposition}

\begin{proof}
We direct the reader to \cite{naisse2020odd} Proposition 4.18 for a complete proof; the one here is completely analogous.
\end{proof}

\subsection{Homogeneous maps}
\label{ss:homomaps}

One of the goals of this paper is to prove an adjunction for unified Khovanov homology, generalizing Theorem 2.31 of \cite{https://doi.org/10.48550/arxiv.1405.2574}. This means we must define $\mathrm{HOM}$-complexes which, in our case, necessitates defining what is meant by maps of homogeneous $\mathscr{G}$-degree. This opens a whole can of worms, which most of the rest of this section is devoted to describing. We proceed with the same assumptions as before: $(\mathscr{C},\alpha)$ is a grading multicategory, and $S = (\mathscr{I}, \{\varphi_i\}_{i\in \mathscr{I}})$ is a shifting system compatible with $\alpha$ through maps $\beta$. Moreover, all $\mathscr{C}$-graded algebras are assumed to be supported entirely in $\Sigma$ so that previous results hold.

\begin{definition}
\label{def:purelyhomo}
Suppose $M$ and $N$ are $\mathscr{C}$-graded $(A_1,\ldots, A_k; B)$-multimodules. A $\mathbb{K}$-linear map $f: M \to N$ is called \textit{purely homogeneous of degree $i$} (for $i \in \mathscr{I} \sqcup \{\mathrm{Id}\}$) if, for all $m\in M$,
\begin{enumerate}[label = (\roman*)]
    \item $f(m) = 0$ if $\abs{m} \not\in \mathsf{D}_i$, 
    \item $\abs{f(m)} = \varphi_i(\abs{m})$ if $\abs{m} \in \mathsf{D}_i$,
    \item $\rho_L(\vec{a}, f(m)) = \beta_{(e,\ldots, e), i}(\abs{\vec{a}}, \abs{m}) f( \rho_L(\vec{a}, m)$ for all $\vec{a} \in (A_1,\ldots, A_k)$, and
    \item $\rho_R(f(m), b) = \beta_{i,e} (\abs{m}, \abs{b}) f(\rho_R(m,b))$ for all $b\in B$.
\end{enumerate}
A map $f: M \to M$ is called \textit{homogeneous} if it is a finite sum of purely homogeneous maps, written $f = \sum_{j} f^j$. We'll write $\abs{f} = i$ if $f$ is a purely homogeneous map of degree $i$.
\end{definition}
Importantly, we do not require that a purely homogeneous map preserve $\mathscr{C}$-degree; however, every purely homogeneous map of degree $i$, $f: M\to N$, induces a graded one, $\widetilde{f}: \varphi_i(M) \to N$, by setting $\widetilde{f}(\varphi_i(m)) = f(m)$.

Using the shifting system and compatibility maps, we can define the tensor product of homogeneous maps. Let $f_i: M_i \to N_i$ for $i=1, \ldots, k$ and $f: M \to N$ be (not necessarily purely) homogeneous maps of $(A_{i1},\ldots, A_{i\alpha_i}; B_i)$-multimodules and $(B_1,\ldots, B_k; C)$-multimodules respectively. Then, define
\[
(f_1, \ldots, f_k) \otimes f: (M_1, \ldots, M_k) \otimes M \to (N_1,\ldots, N_k) \otimes N
\]
by setting $(f_1,\ldots, f_k) \otimes f = \sum_j [(f_1, \ldots, f_k)\otimes f]^j$ where
\[
[(f_1,\ldots, f_k)\otimes f]^j ((m_1,\ldots, m_k) \otimes m)  = \sum_{(i_1,\ldots, i_k)\bullet i = j} \beta_{\abs{\vec{f}}, \abs{f}} (\abs{\vec{m}}, \abs{m})^{-1} \left(f_1^{i_1}(m_1), \ldots, f_k^{i_k}(m_k)\right) \otimes f^i(m)
\]
for all homogeneous elements $\vec{m} \in (M_1,\ldots, M_k)$, $m\in M$.

First, notice that homogeneous maps behave well with respect to this tensor product (or, horizontal composition).

\begin{proposition}
\label{prop:C-graded tensor}
If $f_1,\ldots, f_k, f$ are purely homogeneous maps of degrees $i_1, \ldots, i_k$ and $i$ respectively, then $(f_1, \ldots, f_k) \otimes f$ is purely homogeneous of degree $(i_1,\ldots, i_k) \bullet i$.
\end{proposition}

\begin{proof}
For requirement (i), recall that $\abs{(m_1,\ldots, m_k) \otimes m} = g \circ (g_1,\ldots, g_k)$. The assumption that $\abs{\vec{m} \otimes m} \not\in \mathsf{D}_{\vec{i} \bullet i}$ implies that either $g \not\in \mathsf{D}_i$, hence $f(m) = 0$ since $f$ is homogeneous of degree $i$, or $g_\ell \not \in \mathsf{D}_{i_\ell}$ for some $\ell$, in which case $f(m_\ell) = 0$ for the same reason. Thus $\left((f_1, \ldots, f_k) \otimes f\right)(\vec{m} \otimes m) = 0.$

For (ii), we compute
\begin{align*}
\abs{\left((f_1, \ldots, f_k) \otimes f\right)(\vec{m} \otimes m)} &= \beta_{\vec{i}, i} \left(\abs{\vec{m}}, \abs{m}\right)^{-1} \abs{\left(f_1(m_1), \ldots, f_k(m_k)\right) \otimes f(m)} \\
&= \beta_{\vec{i}, i} \left(\abs{\vec{m}}, \abs{m}\right)^{-1} \left(\varphi_{i_1}(\abs{m_1}), \ldots, \varphi_{i_k}(\abs{m_k})\right) \circ \varphi_i(\abs{m})\\
&= \varphi_{\vec{i} \bullet i} \left(\abs{\vec{m} \otimes m}\right)
\end{align*}
as desired.

For (iii), 
\begin{align*}
    &\rho_L\left(\vec{a}, \left((f_1, \ldots, f_k) \otimes f\right) (\vec{m} \otimes m)\right) = \beta_{\vec{i}, i} \left(\abs{\vec{m}}, \abs{m} \right)^{-1} \rho_L \left(\vec{a}, \left((f_1(m_1), \ldots, f_k(m_k) \right) \otimes f(m) \right)
    \\
    &= \beta_{\vec{i}, i} \left(\abs{\vec{m}}, \abs{m} \right)^{-1} \alpha(\abs{\vec{a}}, \underbrace{\abs{\vec{f}(\vec{m})}}_{=\varphi_{\vec{i}}(\abs{\vec{m}})}, \underbrace{\abs{f(m)}}_{=\varphi_i(\abs{m})})^{-1} \left(\rho_L^1\left(\vec{a}_1, f_1(m_1)\right), \ldots, \rho_L^k\left(\vec{a}_k, f_k(m_k)\right)\right) \otimes f(m)
    \\
    &= \beta_{\vec{i}, i} \left(\abs{\vec{m}}, \abs{m} \right)^{-1} \alpha(\abs{\vec{a}},\varphi_{\vec{i}}(\abs{\vec{m}}), \varphi_i(\abs{m}))^{-1} \beta_{\vec{e}, \vec{i}}(\abs{\vec{a}}, \abs{\vec{m}})\left( f_1\left(\rho_L^1(\vec{a}_1, m_1)\right), \ldots, f_k\left(\rho_L^k(\vec{a}_k, m_k)\right)\right) \otimes f(m)
    \\
    &= \beta_{\vec{e}, \vec{i}\bullet i}(\abs{\vec{a}}, \abs{\vec{m}} \circ \abs{m}) \alpha(\abs{\vec{a}}, \abs{\vec{m}}, \abs{m})^{-1} \beta_{\vec{e} \bullet \vec{i}, i}(\abs{\vec{a}}\circ \abs{\vec{m}}, \abs{m})^{-1} \left( f_1\left(\rho_L^1(\vec{a}_1, m_1)\right), \ldots, f_k\left(\rho_L^k(\vec{a}_k, m_k)\right)\right) \otimes f(m)
    \\
    &= \beta_{\vec{e}, \vec{i}\bullet i}(\abs{\vec{a}}, \abs{\vec{m}} \circ \abs{m}) \alpha(\abs{\vec{a}}, \abs{\vec{m}}, \abs{m})^{-1} \left((f_1,\ldots, f_k) \otimes f \right) \left( \left(\rho_L^1(\vec{a}_1, m_1), \ldots, \rho_L^k(\vec{a}_k, m_k)\right) \otimes m \right)
    \\
    &= \beta_{\vec{e}, \vec{i}\bullet i}(\abs{\vec{a}}, \abs{\vec{m}} \circ \abs{m}) \left((f_1,\ldots, f_k) \otimes f \right) \left(\rho_L(\vec{a}, \vec{m} \otimes m) \right)
\end{align*}
The first equality is by definition and $\mathbb{K}$-linearity of the left action. The second equality is by the definition of the $\mathscr{C}$-graded multimodule left-action on $(M_1,\ldots, M_k) \otimes M$. The third equality follows from the assumption that $\abs{f_\ell} = i_\ell$. The fourth equality follows from equation (\ref{eq:shiftingsystcomp}). Finally, the fifth and sixth equalities follow from unraveling definitions; in particular, the fifth follows since $\abs{\vec{a}} \circ \abs{\vec{m}} = \abs{\left(\rho_L^1(\vec{a}_1, m_1), \ldots, \rho_L^k(\vec{a}_k, m_k\right)}$ and
the sixth invokes the $\mathbb{K}$-linearity of $(f_1,\ldots, f_k) \otimes f$. Finally, this gives us the desired result since $\abs{\vec{m}} \circ \abs{m} = \abs{\vec{m} \otimes m}$: we have
\[
\rho_L\left(\vec{a}, \left((f_1, \ldots, f_k) \otimes f\right) (\vec{m} \otimes m)\right) = \beta_{\vec{e}, \vec{i}\bullet i}(\abs{\vec{a}}, \abs{\vec{m}\otimes m}) \left((f_1,\ldots, f_k) \otimes f \right) \left(\rho_L(\vec{a}, \vec{m} \otimes m) \right).
\]
Showing that (iv) holds is completely analogous (and easier)---we leave it to the reader.
\end{proof}

On the other hand, we do not yet have a method for composing grading shifts vertically, so that we cannot define the composition of homogeneous maps. We introduce the fix in the following section.

\subsubsection{Extension to a shifting 2-system}

As before, we will consider the $\mathscr{G}$-graded situation and then present generalities. Thankfully, the extension of a $\mathscr{C}$-shifting system to a $\mathscr{C}$-shifting 2-system is almost exactly like the categorically-graded situation.

Suppose $\Delta_1: D_1 \to D_1'$ and $\Delta_2: D_2 \to D_2'$ are cobordisms of planar arc diagrams, so that $(\Delta_1, v_1)$ and $(\Delta_2, v_2)$ induce grading shift functors for any $v_1, v_2 \in \mathbb{Z} \oplus \mathbb{Z}$; that is, they belong to the multimonoid $\mathscr{I}$ of $\mathscr{G}$. Define a binary operation, which we call vertical composition, 
\[
\circ: \mathscr{I} \times \mathscr{I} \to \mathscr{I}
\]
by stacking: set
\[
(\Delta_2, v_2) \circ (\Delta_1, v_1) = \begin{cases} (\Delta_2 \circ \Delta_1, v_2 + v_1) & \text{if}~D_1' = D_2,~\text{and}\\ 0 & \text{otherwise}.\end{cases}
\]
In our multivariable setting, vertical composition must be extended to a family of vertical compositions for each $k\ge 1$,
\[
\circ: \mathscr{I}^k \times \mathscr{I}^k \to \mathscr{I}^k.
\]
So, if $\vec{\Delta}_i = (\Delta_{i1}, \ldots, \Delta_{ik})$ for $i=1,2$ and $D_{1j}\xrightarrow{\Delta_{1j}} D_{1j}'$ and $D_{2j} \xrightarrow{\Delta_{2j}} D_{2j}'$ for $j=1,\ldots, k$, we set 
\[
(\vec{\Delta}_2, \vec{v}_2) \circ (\vec{\Delta}_1, \vec{v}_1) = \begin{cases} (\vec{\Delta}_2 \circ \vec{\Delta}_1, \vec{v}_2 + \vec{v}_1) & \text{if}~D_{1j}' = D_{2j}~\text{for all}~j,~\text{and}\\ 0 & \text{otherwise.}\end{cases}
\] 
The nonzero term on the right is given a chronology as follows.
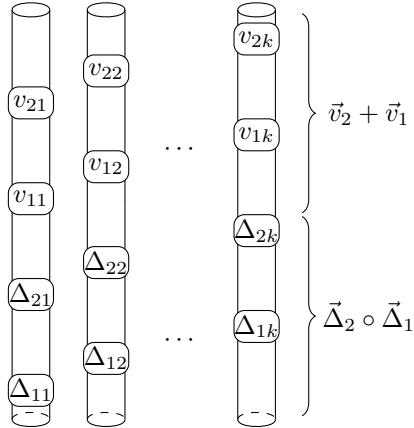
\[
\tikz[baseline=5ex, scale=.5]{
	\draw (0,-0.1) .. controls (0,.-.35) and (1,-.35) .. (1,-0.1);
	\draw[dashed] (0,-0.1) .. controls (0,.15) and (1,.15) .. (1,-0.1);
	\draw (0,-0.1) -- (0,5.7);
	\draw (1,-0.1) -- (1,5.7);
	\draw (2,-0.1) .. controls (2,.-.35) and (3,-.35) .. (3,-0.1);
	\draw[dashed] (2,-0.1) .. controls (2,.15) and (3,.15) .. (3,-0.1);
	\draw (2,-0.1) -- (2,5.7);
	\draw (3,-0.1) -- (3,5.7);
        \node at (4.5, 2) {$\cdots$};
	\draw (6,-0.1) .. controls (6,-.35) and (7,-.35) .. (7,-0.1);
	\draw[dashed] (6,-0.1) .. controls (6,.15) and (7,.15) .. (7,-0.1);
	\draw (6,-0.1) -- (6,5.7);
	\draw (7,-0.1) -- (7,5.7);
        %
        \filldraw [fill=white, draw=black,rounded corners] (-.1,.25) rectangle (1.1, 1.1) node[midway] {$\Delta_{11}$};
        \filldraw [fill=white, draw=black,rounded corners] (1.9,1.1) rectangle (3.1, 1.95) node[midway] {$\Delta_{12}$};
        \filldraw [fill=white, draw=black,rounded corners] (5.9,1.95) rectangle (7.1, 2.8) node[midway] {$\Delta_{1k}$};
        \filldraw [fill=white, draw=black,rounded corners] (-0.1,2.8) rectangle (1.1, 3.65) node[midway] {$\Delta_{21}$};
        \filldraw [fill=white, draw=black,rounded corners] (1.9,3.65) rectangle (3.1, 4.5) node[midway] {$\Delta_{22}$};
        \filldraw [fill=white, draw=black,rounded corners] (5.9,4.5) rectangle (7.1, 5.35) node[midway] {$\Delta_{2k}$};
    \begin{scope}[shift={(0,5.1)}]
	\draw (0,.59) -- (0,5.7);
	\draw (1,.59) -- (1,5.7);
	\draw (0,5.7) .. controls (0,5.45) and (1,5.45) .. (1,5.7);
	\draw (0,5.7) .. controls (0,5.95) and (1,5.95) .. (1,5.7);
	\draw (2,.59) -- (2,5.7);
	\draw (3,.59) -- (3,5.7);
	\draw (2,5.7) .. controls (2,5.45) and (3,5.45) .. (3,5.7);
	\draw (2,5.7) .. controls (2,5.95) and (3,5.95) .. (3,5.7);
        \node at (4.5, 2) {$\cdots$};
	\draw (6,.59) -- (6,5.7);
	\draw (7,.59) -- (7,5.7);
	\draw (6,5.7) .. controls (6,5.45) and (7,5.45) .. (7,5.7);
	\draw (6,5.7) .. controls (6,5.95) and (7,5.95) .. (7,5.7);
        \filldraw [fill=white, draw=black,rounded corners] (-.1,.25) rectangle (1.1, 1.1) node[midway] {$v_{11}$};
        \filldraw [fill=white, draw=black,rounded corners] (1.9,1.1) rectangle (3.1, 1.95) node[midway] {$v_{12}$};
        \filldraw [fill=white, draw=black,rounded corners] (5.9,1.95) rectangle (7.1, 2.8) node[midway] {$v_{1k}$};
        \filldraw [fill=white, draw=black,rounded corners] (-0.1,2.8) rectangle (1.1, 3.65) node[midway] {$v_{21}$};
        \filldraw [fill=white, draw=black,rounded corners] (1.9,3.65) rectangle (3.1, 4.5) node[midway] {$v_{22}$};
        \filldraw [fill=white, draw=black,rounded corners] (5.9,4.5) rectangle (7.1, 5.35) node[midway] {$v_{2k}$};
    \end{scope}
        \draw [decorate,decoration={brace,amplitude=5pt,mirror,raise=4ex}]
  (6.5,0) -- (6.5,5.3) node[midway,xshift=4.25em]{$\vec{\Delta}_2 \circ \vec{\Delta}_1$};
        \draw [decorate,decoration={brace,amplitude=5pt,mirror,raise=4ex}]
  (6.5,5.4) -- (6.5,10.7) node[midway,xshift=4.25em]{$\vec{v}_2 + \vec{v}_1$};
}
\]
Again, we choose this particular chronology so that the arguments of \cite{naisse2020odd} lift to our setting.

In general, if $S = \{\mathscr{I}, \{\varphi_i\}\}$ is already a $\mathscr{C}$-shifting system, equipping $\mathscr{I}$ with a vertical composition map $\mathscr{I}^k \times \mathscr{I}^k \to \mathscr{I}^k$ of this form constitutes what is called a \textit{$\mathscr{C}$-shifting 2-system}, granted it satisfies the following requirements:
\begin{enumerate}[label = (\roman*)]
    \item $e \circ e = e$,
    \item $\mathsf{D}_{j\circ i} = \mathsf{D}_i \cap \varphi_i^{-1} (\mathsf{D}_j)$,
    \item $\varphi_{j\circ i} = \varphi_j|_{\varphi_i(\mathsf{D}_i) \cap \mathsf{D}_j} \circ \varphi_i|_{\mathsf{D}_{j\circ i}}$, and
    \item $\varphi_{\left((j_1,\ldots, j_k) \circ (i_1,\ldots, i_k)\right) \bullet (j\circ i)} = \varphi_{\left((j_1,\ldots, j_k)\bullet j\right) \circ \left((i_1, \ldots, i_k) \bullet i\right)}$ for all $j_1, \ldots, j_k, j, i_1,\ldots, i_k, i \in \mathscr{I}$.
\end{enumerate}
The first three requirements are written in the single-input case to ignore burdensome notation and should be extended to the $k$-input cases. To elucidate the above requirements notice that (in particular, if $j\circ i$ is nonzero) $\varphi_j$ and $\varphi_i$ must be defined on (frequently distinct) subsets of the same hom-set. In the $\mathscr{G}$-graded case, this causes no confusion: $\mathsf{D}_{j\circ i} = \mathsf{D}_i$ since cobordisms which start at $D_1$ and factor through $D_1' = D_2$ still start at $D_1$. In general, we should be a little more careful: 
\[
\begin{tikzcd}
\mathsf{D}_i \arrow[r, "\varphi_i"] & \mathrm{Hom}_{\mathscr{C}}(X_1,\ldots, X_k; Y) \supset \mathsf{D}_j \arrow[r, "\varphi_j"] & \mathrm{Hom}_{\mathscr{C}}(X_1,\ldots, X_k; Y)
\end{tikzcd}
\]
so, in general, $\varphi_{j\circ i}$ is defined only on the subset $\mathsf{D}_i \cap \varphi_i^{-1}(\mathsf{D}_j)$, as in (ii) and (iii). Condition (iv) just ensures that vertical composition and horizontal composition play nicely together---(iv) obviously holds in the $\mathscr{G}$-setting for weighted cobordisms of planar arc diagrams.

For completeness, we include a description of compatibility maps. We say that a $\mathscr{C}$-shifting 2-system $S = \{\mathscr{I}, \{\varphi_i\}_{i\in\mathscr{I}}\}$ is \textit{compatible} with the associator $\alpha$ of $\mathscr{C}$ if there are $(\beta, \gamma, \Xi)$ such that the underlying $\mathscr{C}$-shifting system is compatible with $\alpha$ through $\beta$, and $\gamma$ and $\Xi$ are as follows. First, $\gamma$ stands for a collection of maps 
\[
\gamma_{i,j}^{\vec{X}\to Y}: \mathsf{D}_i^{\vec{X}\to Y} \to \mathbb{K}^\times
\]
for all $i,j\in \mathscr{I}$ and $\vec{X}, Y\in \mathscr{C}$ satisfying $\gamma_{i,j} = 1$ whenever $i\in \mathscr{I}_{\mathrm{id}}, j\in \mathscr{I}_{\mathrm{id}},$ or $i = j = e$. More generally, we construct multivariable functions
\[
\gamma_{\vec{i}, \vec{j}}^{\vec{X} \to \vec{Y}}: \prod_{\ell=1}^k \mathsf{D}_{i_\ell}^{\vec{X}_i \to Y_i} \to \mathbb{K}^\times
\]
with analogous requirements ($\gamma_{\vec{i}, \vec{j}} = 1$ whenever each entry of $\vec{i}$ belongs to $\mathscr{I}_{\mathrm{Id}}$, each entry of $\vec{j}$ belongs to $\mathscr{I}_{\mathrm{Id}}$, or $\vec{i} = \vec{j} = \vec{e}$). We do not require that $\gamma_{\vec{i}, \vec{j}}^{\vec{X} \to \vec{Y}} = \gamma_{i_1, j_1}^{\vec{X}_1 \to Y_1}\cdots \gamma_{i_k, j_k}^{\vec{X}_k \to Y_k}$. For example, this is not the case for the $\mathscr{G}$-graded setting, at least the way we've set things up. Second, $\Xi$ stands for a collection of invertible scalars
\[
\Xi_{\substack{i, \vec{i} \\ j, \vec{j}}}^{\vec{X} \to \vec{Y} \to Z} \in \mathbb{K}^\times
\]
satisfying (i) $\Xi_{\substack{i, \vec{i} \\ j, \vec{j}}} = 1$ whenever $(\vec{j} \circ \vec{i}) \bullet (j \circ i) = (\vec{j} \bullet j) \circ (\vec{i} \bullet i)$ and (ii) $\Xi_{\substack{i, \vec{i} \\ j, \vec{j}}}$ is invariant when exchanging elements of $\mathscr{I}_\mathrm{id}$ out with other elements of $\mathscr{I}_\mathrm{id}$. We often write $\Xi_{\substack{i, \vec{i} \\ j, \vec{j}}}(g'g)$ for $\Xi_{\substack{i, \vec{i} \\ j, \vec{j}}}^{\vec{X} \to \vec{Y} \to Z}$ when $\vec{X}\xrightarrow{g'} \vec{Y} \xrightarrow{g} Z$, or $\Xi_{\substack{i, \vec{i} \\ j, \vec{j}}}(g)$ for $\Xi_{\substack{i, \vec{i} \\ j, \vec{j}}}^{\vec{X} \to Y}$ when $\vec{X} \xrightarrow{g} Y$. Finally, we say that the shifting 2-system is compatible with $\alpha$ through $(\beta, \gamma, \Xi)$ if, in addition, the two following equations hold. The first reads
\begin{equation}
\label{eq:shiftingsystvcomp}
\gamma_{\vec{i} \bullet i, \vec{j} \bullet j}(g'g) 
\beta_{i, \vec{i}}(g',g) 
\beta_{j ,\vec{j}}(\varphi_{\vec{i}}(g'), \varphi_i(g))
=
\Xi_{\substack{i, \vec{i} \\ j, \vec{j}}}(g'g) \beta_{j\circ i, \vec{j} \circ \vec{i}}(g',g) 
\gamma_{\vec{i}, \vec{j}}(g') 
\gamma_{i,j}(g),
\end{equation}
for all $g' \in \mathsf{D}_{\vec{i}}^{\vec{Z} \to \vec{Y}}$ and $g \in \mathsf{D}_i^{\vec{X} \to Y}$.
Again, this looks awful, but this is just to say that $\gamma$ and $\Xi$ are chosen so that the following diagram commutes.
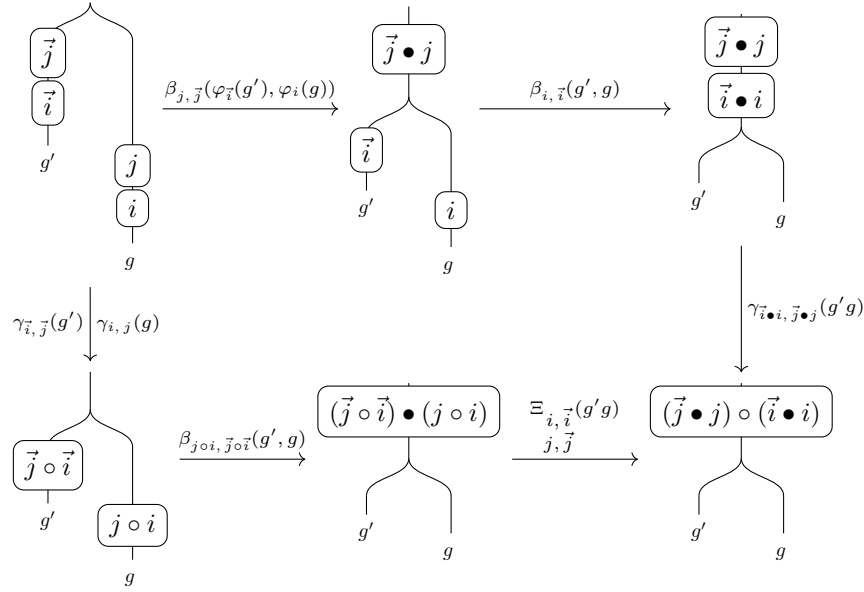
\[
\begin{tikzcd}[column sep=huge, row sep=large, scale=1.1]
\tikz[yscale=.5,xscale=.75]{
	\draw 
		(.5,2)
		.. controls (.5,2.5) and (1.25,2.5) ..
		(1.25,3);
	\draw (.5,-.8) node[below,scale=.75]{$g'$}
		-- 
		(.5,-.8)
		--
		(.5,1)
		.. controls (.5,1.5) and (.5,1.5) ..
		(.5,2);
	\draw (2,-3.4) node[below,scale=.75]{$g\phantom{'}$}
		--
		(2, -2.5)
		--
		(2,2)
		.. controls (2,2.5) and (1.25,2.5) ..
		(1.25,3);
        \node[fill=white,draw,rounded corners,scale=1] at (.5,1.4) {$\vec{j}$};
        \node[fill=white,draw,rounded corners,scale=1] at (.5,0) {$\vec{i}$};
        \node[fill=white,draw,rounded corners,scale=1] at (2, -1.5) {$j$};
        \node[fill=white,draw,rounded corners,scale=1] at (2, -2.7) {$i$};
}
\arrow[r, "\beta_{j\text{,}\,\vec{j}}(\varphi_{\vec{i}}(g')\text{,}\, \varphi_i(g))"] \arrow[d, "\gamma_{\vec{i}\text{,}\, \vec{j}}(g')"', "\gamma_{i\text{,}\,j}(g)"]
&
\tikz[yscale=.5,xscale=.75,yshift=-2.5cm]{
        \draw (1.25, 3) -- (1.25, 5.4);
	\draw 
		(.5,2)
		.. controls (.5,2.5) and (1.25,2.5) ..
		(1.25,3);
	\draw (.5,.5) node[below,scale=.75]{$g'$}
		--
		(.5,1)
		.. controls (.5,1.5) and (.5,1.5) ..
		(.5,2);
	\draw (2,-1) node[below,scale=.75]{$g\phantom{'}$}
		--
		(2,2)
		.. controls (2,2.5) and (1.25,2.5) ..
		(1.25,3);
        \node[fill=white,draw,rounded corners,scale=1] at (.5,1.25) {$\vec{i}$};
        \node[fill=white,draw,rounded corners,scale=1] at (2, -.25) {$i$};
        \node[fill=white,draw,rounded corners,scale=1] at (1.25, 4) {$\vec{j}\bullet j$};
}
\arrow[r, "\beta_{i\text{,}\, \vec{i}}(g'\text{,}\, g)"]
&
\tikz[yscale=.5,xscale=.75,yshift=-3.3cm]{
        \draw (1.25, 3) -- (1.25, 6);
	\draw 
		(.5,2)
		.. controls (.5,2.5) and (1.25,2.5) ..
		(1.25,3);
	\draw (.5,1.5) node[below,scale=.75]{$g'$}
		.. controls (.5,1.5) and (.5,1.5) ..
		(.5,2);
	\draw (2,1) node[below,scale=.75]{$g\phantom{'}$}
		--
		(2,2)
		.. controls (2,2.5) and (1.25,2.5) ..
		(1.25,3);
        \node[fill=white,draw,rounded corners,scale=1] at (1.25, 3.5) {$\vec{i}\bullet i$};
        \node[fill=white,draw,rounded corners,scale=1] at (1.25, 5) {$\vec{j}\bullet j$};
}
\arrow[d, "\gamma_{\vec{i} \bullet i\text{,}\, \vec{j} \bullet j}(g'g)"]
\\
\tikz[yscale=.5,xscale=.75,yshift=-1.5cm]{
        \draw (1.25, 3) -- (1.25, 4);
	\draw 
		(.5,2)
		.. controls (.5,2.5) and (1.25,2.5) ..
		(1.25,3);
	\draw (.5,.5) node[below,scale=.75]{$g'$}
		--
		(.5,1)
		.. controls (.5,1.5) and (.5,1.5) ..
		(.5,2);
	\draw (2,-1) node[below,scale=.75]{$g\phantom{'}$}
		--
		(2,2)
		.. controls (2,2.5) and (1.25,2.5) ..
		(1.25,3);
        \node[fill=white,draw,rounded corners,scale=1] at (.5,1.25) {$\vec{j} \circ \vec{i}$};
        \node[fill=white,draw,rounded corners,scale=1] at (2, -.25) {$j\circ i$}; 
} 
\arrow[r, "\beta_{j\circ i\text{,}\, \vec{j} \circ \vec{i}}(g'\text{,}\,g)"]
&
\tikz[yscale=.5,xscale=.75,yshift=-2.8cm]{
        \draw (1.25, 3) -- (1.25, 5);
	\draw 
		(.5,2)
		.. controls (.5,2.5) and (1.25,2.5) ..
		(1.25,3);
	\draw (.5,1.5) node[below,scale=.75]{$g'$}
		.. controls (.5,1.5) and (.5,1.5) ..
		(.5,2);
	\draw (2,1) node[below,scale=.75]{$g\phantom{'}$}
		--
		(2,2)
		.. controls (2,2.5) and (1.25,2.5) ..
		(1.25,3);
        \node[fill=white,draw,rounded corners,scale=1] at (1.25, 4) {$(\vec{j}\circ \vec{i}) \bullet (j \circ i)$};
}
\arrow[r, "\Xi_{\substack{i\text{,}\, \vec{i} \\ j\text{,}\, \vec{j}}}(g'g)"]
&
\tikz[yscale=.5,xscale=.75,yshift=-2.8cm]{
        \draw (1.25, 3) -- (1.25, 5);
	\draw 
		(.5,2)
		.. controls (.5,2.5) and (1.25,2.5) ..
		(1.25,3);
	\draw (.5,1.5) node[below,scale=.75]{$g'$}
		.. controls (.5,1.5) and (.5,1.5) ..
		(.5,2);
	\draw (2,1) node[below,scale=.75]{$g\phantom{'}$}
		--
		(2,2)
		.. controls (2,2.5) and (1.25,2.5) ..
		(1.25,3);
        \node[fill=white,draw,rounded corners,scale=1] at (1.25, 4) {$(\vec{j}\bullet j)\circ(\vec{i}\bullet i)$};
}
\end{tikzcd}
\]
The second requirement reads
\begin{equation}
\label{eq:gammarec}
\gamma_{i, k\circ j}(g) \gamma_{j,k}(\varphi_i(g)) = \gamma_{j\circ i, k}(g) \gamma_{i,j}(g)
\end{equation}
for all $g\in \mathsf{D}_i^{\vec{X}\to Y}$ which, a little more obviously, is to say the following diagram commutes.
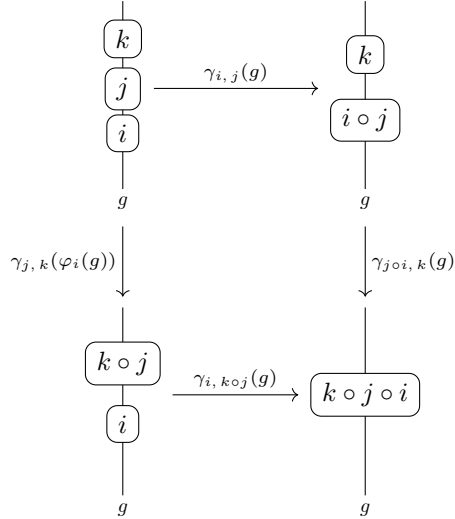
\[
\begin{tikzcd}[column sep=huge, row sep=large, scale=1.1]
\tikz[yscale=.5,xscale=.75]{
    \draw (0,-2.5) node[below,scale=.75]{$g$} -- (0,2.5);
    \node[fill=white,draw,rounded corners,scale=1] at (0,-1.25) {$i$};
    \node[fill=white,draw,rounded corners,scale=1] at (0,0) {$j$};
    \node[fill=white,draw,rounded corners,scale=1] at (0,1.25) {$k$};
}
\arrow[r, "\gamma_{i\text{,}\,j}(g)"] \arrow[d, "\gamma_{j\text{,}\,k}(\varphi_i(g))"']
&
\tikz[yscale=.5,xscale=.75]{
    \draw (0,-2.5) node[below,scale=.75]{$g$} -- (0,2.5);
    \node[fill=white,draw,rounded corners,scale=1] at (0,-0.833) {$i\circ j$};
    \node[fill=white,draw,rounded corners,scale=1] at (0,0.833) {$k$};
}
\arrow[d, "\gamma_{j\circ i\text{,}\, k}(g)"]
\\
\tikz[yscale=.5,xscale=.75]{
    \draw (0,-2.5) node[below,scale=.75]{$g$} -- (0,2.5);
    \node[fill=white,draw,rounded corners,scale=1] at (0,-0.833) {$i$};
    \node[fill=white,draw,rounded corners,scale=1] at (0,0.833) {$k\circ j$};
}
\arrow[r, "\gamma_{i\text{,}\, k\circ j}(g)"]
&
\tikz[yscale=.5,xscale=.75]{
    \draw (0,-2.5) node[below,scale=.75]{$g$} -- (0,2.5);
    \node[fill=white,draw,rounded corners,scale=1] at (0,0) {$k\circ j \circ i$};
}
\end{tikzcd}
\]

So, to extend the $\mathscr{G}$-shifting system we have to a $\mathscr{G}$-shifting 2-system, we choose compatibility maps
\[
\begin{split}
\gamma_{(\Delta_1, v_1), (\Delta_2, v_2)}^{(x_1,\ldots, x_k) \to y} (D,p) = \lambda\left(\abs{1_{\vec{x}} \Delta_2 1_y}, v_1\right)
\end{split}
\qquad \text{or, in general,} \qquad
\begin{split}
\gamma_{(\vec{\Delta}_1, \vec{v}_1),(\vec{\Delta}_2, \vec{v}_2)}^{\vec{x}\to \vec{y}}(\vec{D}, \vec{p}) = \lambda\left(\abs{1_{\vec{x}} \vec{\Delta}_2 1_{\vec{y}}}, V_1 \right)
\end{split}
\]
where $V_1$ is the sum of entries in $\vec{v}_1$, and
\[
\Xi_{\substack{(\Delta_1, v_1), (\vec{\Delta}_1, \vec{v}_1) \\ (\Delta_2, v_2), (\vec{\Delta}_2, \vec{v}_2)}}^{\vec{x} \to \vec{y} \to z} = \iota({}_{\vec{x}}H_z) \lambda(V_1, v_2)
\]
where $H:(\vec{\Delta}_2 \circ \vec{\Delta}_1)\bullet (\Delta_2 \circ \Delta_1) \Rightarrow (\vec{\Delta}_2 \bullet \Delta_2) \circ (\vec{\Delta}_1 \bullet \Delta_1)$ and $V_1$, as before, means the sum of the entries of $\vec{v}_1$. We refer to the first factor of $\Xi$ by $\Xi_1$ and the second factor by $\Xi_2$. Of course, the definitions above only hold if the cobordisms involved are vertically composable with respect to the chosen order; otherwise, these maps are zero.

To understand where these choices come from, notice that $(\Delta_2, v_2) \circ (\Delta_1, v_1)$ can be rewritten schematically as 
\[
\tikz[yscale=.5,xscale=.75,baseline=0ex]{
    \draw (0,-3) node[below,scale=.75]{$g$} -- (0,3);
    \node[fill=white,draw,rounded corners,scale=1] at (0,-1.8) {$\Delta_1$};
    \node[fill=white,draw,rounded corners,scale=1] at (0,-0.6) {$v_1$};
    \node[fill=white,draw,rounded corners,scale=1] at (0,0.6) {$\Delta_2$};
    \node[fill=white,draw,rounded corners,scale=1] at (0,1.8) {$v_2$};
}
= \lambda(\abs{1_{\vec{x}}\Delta_2 1_y}, v_1)
\tikz[yscale=.5,xscale=.75,baseline=0ex]{
    \draw (0,-3) node[below,scale=.75]{$g$} -- (0,3);
    \node[fill=white,draw,rounded corners,scale=1] at (0,-1.8) {$\Delta_1$};
    \node[fill=white,draw,rounded corners,scale=1] at (0,-0.6) {$\Delta_2$};
    \node[fill=white,draw,rounded corners,scale=1] at (0,0.6) {$v_1$};
    \node[fill=white,draw,rounded corners,scale=1] at (0,1.8) {$v_2$};
}
\]
for $g: \vec{x} \to y$. That is, $(\Delta_2, v_2) \circ (\Delta_1, v_1) = \lambda(\abs{1_{\vec{x}}\Delta_2 1_y}, v_1) (\Delta_2 \circ \Delta_1, v_2 + v_1)$, so we hope $\gamma$ has the form above. For $\Xi$, we can start by recognizing that, schematically (and thanks to our chronology conventions), $(\vec{\Delta}_2 \circ \vec{\Delta}_1) \bullet (\Delta_2 \circ \Delta_1)$ looks like
\[
\tikz[yscale=.5,xscale=.75,baseline=0ex]{
    \draw (0,-3) node[below,scale=.75]{$g$} -- (0,3);
    \node[fill=white,draw,rounded corners,scale=1] at (0,-2.3) {$\vec{1} \bullet \Delta_1$};
    \node[fill=white,draw,rounded corners,scale=1] at (0,-0.78) {$\vec{1} \bullet \Delta_2$};
    \node[fill=white,draw,rounded corners,scale=1] at (0,0.78) {$\vec{\Delta}_1 \bullet 1$};
    \node[fill=white,draw,rounded corners,scale=1] at (0,2.3) {$\vec{\Delta}_2 \bullet 1$};
}
\]
Where ``$1$'' and ``$\vec{1}$'' just stand for the identity cobordism on their respective components (in particular, an element of $\mathscr{I}_{\text{Id}}$).On the other hand, $(\vec{\Delta}_2 \bullet \Delta_2) \circ (\vec{\Delta}_1 \bullet \Delta_1)$ looks like 
\[
\tikz[yscale=.5,xscale=.75,baseline=0ex]{
    \draw (0,-3) node[below,scale=.75]{$g$} -- (0,3);
    \node[fill=white,draw,rounded corners,scale=1] at (0,-2.3) {$\vec{1} \bullet \Delta_1$};
    \node[fill=white,draw,rounded corners,scale=1] at (0,0.78) {$\vec{1} \bullet \Delta_2$};
    \node[fill=white,draw,rounded corners,scale=1] at (0,-0.78) {$\vec{\Delta}_1 \bullet 1$};
    \node[fill=white,draw,rounded corners,scale=1] at (0,2.3) {$\vec{\Delta}_2 \bullet 1$};
}
\]
So, we have that
\[
(\vec{\Delta}_2 \circ \vec{\Delta}_1) \bullet (\Delta_2 \circ \Delta_1)
= \iota({}_{\vec{x}}H_{z}) (\vec{\Delta}_2 \bullet \Delta_2) \circ (\vec{\Delta}_1 \bullet \Delta_1)
\]
where $H:(\vec{\Delta}_2 \circ \vec{\Delta}_2)\bullet (\Delta_2 \circ \Delta_1) \Rightarrow (\vec{\Delta}_2 \bullet \Delta_2) \circ (\vec{\Delta}_1 \bullet \Delta_1)$ is the locally vertical change of chronology which simply pushes $\Delta_2$ past the cobordisms involved in $\vec{\Delta}_1$. So that $\Xi$ satisfies equation (\ref{eq:shiftingsystvcomp}), we must also multiply by $\lambda(V_1, v_2)$.

\begin{proposition}
The $\mathscr{G}$-shifting 2-system $S$ defined above is compatible with $\alpha$ through $(\beta, \gamma, \Xi)$.
\end{proposition}

\begin{proof}
We know that the underlying shifting system is compatible with $\alpha$ through $\beta$ by Proposition \ref{prop:comp1}. Since 
\[
\mathscr{I}_{\mathrm{id}} = \{(\mathbbm{1}_{D^\wedge}, (0,0)): D~\text{is a planar arc diagram}\},
\]
it is clear that $\gamma$ and $\Xi$ as chosen satisfy preliminary requirements; all we need to do is verify equations (\ref{eq:shiftingsystvcomp}) and (\ref{eq:gammarec}). Verifying (\ref{eq:gammarec}) is easy: computing both sides yields
\[
\lambda(\abs{1_{\vec{x}}(\Delta_3 \circ \Delta_2)1_y}, v_1) \lambda(\abs{1_{\vec{x}} \Delta_3 1_y}, v_2) = \lambda(\abs{1_{\vec{x}} \Delta_3 1_y}, v_1 + v_2) \lambda(\abs{1_{\vec{x}} \Delta_2 1_y}, v_1)
\]
which is true since $\abs{1_{\vec{x}}(\Delta_3 \circ \Delta_2)1_y} = \abs{1_{\vec{x}} \Delta_3 1_y} + \abs{1_{\vec{x}}\Delta_2 1_y}$; applying bilinearity shows that both sides are equal to $\lambda(\abs{1_{\vec{x}} \Delta_3 1_y}, v_1)  \lambda(\abs{1_{\vec{x}} \Delta_3 1_y}, v_2)  \lambda(\abs{1_{\vec{x}} \Delta_2 1_y}, v_1)$. 

Verifying equation (\ref{eq:shiftingsystvcomp}) looks a lot like the proofs of Propositions \ref{prop:associator} and \ref{prop:comp1}. Again, start by considering the contributions of $\beta_1$ and $\Xi_1$ only. To do this, one can consider the two sequences of changes of chronology encoded by the diagram below.
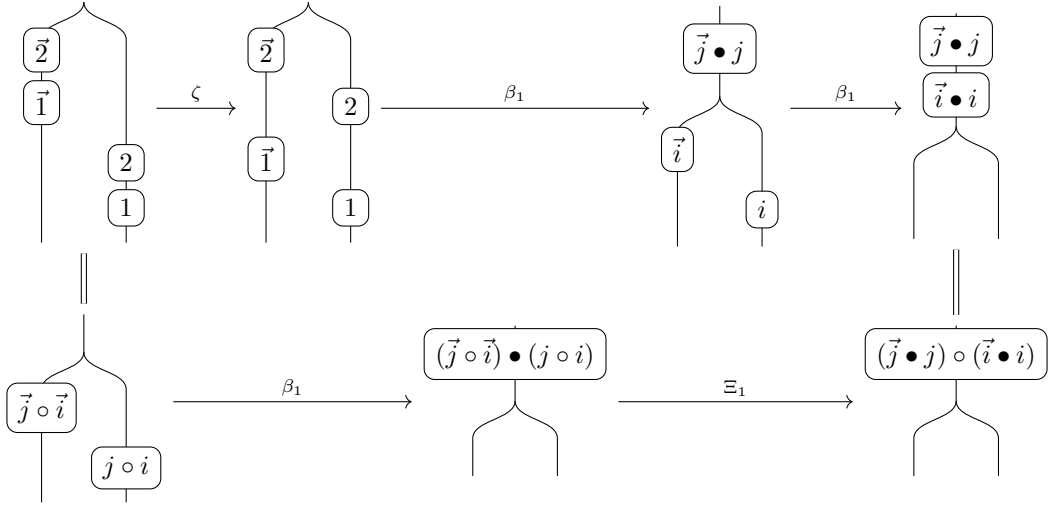
\[
\begin{tikzcd}[column sep = small]
\tikz[yscale=.5,xscale=.75]{
	\draw 
		(.5,2)
		.. controls (.5,2.5) and (1.25,2.5) ..
		(1.25,3);
	\draw (.5,-3.4)
		-- 
		(.5,-.8)
		--
		(.5,1)
		.. controls (.5,1.5) and (.5,1.5) ..
		(.5,2);
	\draw (2,-3.4)
		--
		(2, -2.5)
		--
		(2,2)
		.. controls (2,2.5) and (1.25,2.5) ..
		(1.25,3);
        \node[fill=white,draw,rounded corners,scale=1] at (.5,1.4) {$\vec{2}$};
        \node[fill=white,draw,rounded corners,scale=1] at (.5,0) {$\vec{1}$};
        \node[fill=white,draw,rounded corners,scale=1] at (2, -1.5) {$2$};
        \node[fill=white,draw,rounded corners,scale=1] at (2, -2.7) {$1$};
}
\arrow[d, equal]
\arrow[rr, "\zeta"]
&&
\tikz[yscale=.5,xscale=.75]{
	\draw 
		(.5,2)
		.. controls (.5,2.5) and (1.25,2.5) ..
		(1.25,3);
	\draw (.5,-3.4)
		-- 
		(.5,-.8)
		--
		(.5,1)
		.. controls (.5,1.5) and (.5,1.5) ..
		(.5,2);
	\draw (2,-3.4)
		--
		(2, -2.5)
		--
		(2,2)
		.. controls (2,2.5) and (1.25,2.5) ..
		(1.25,3);
        \node[fill=white,draw,rounded corners,scale=1] at (.5,1.4) {$\vec{2}$};
        \node[fill=white,draw,rounded corners,scale=1] at (2,0) {$2$};
        \node[fill=white,draw,rounded corners,scale=1] at (.5, -1.5) {$\vec{1}$};
        \node[fill=white,draw,rounded corners,scale=1] at (2, -2.7) {$1$};
}
\arrow[rr, "\beta_1"]
&&
\tikz[yscale=.5,xscale=.75,yshift=-2.5cm]{
        \draw (1.25, 3) -- (1.25, 5.4);
	\draw 
		(.5,2)
		.. controls (.5,2.5) and (1.25,2.5) ..
		(1.25,3);
	\draw (.5,-1)
		--
		(.5,1)
		.. controls (.5,1.5) and (.5,1.5) ..
		(.5,2);
	\draw (2,-1)
		--
		(2,2)
		.. controls (2,2.5) and (1.25,2.5) ..
		(1.25,3);
        \node[fill=white,draw,rounded corners,scale=1] at (.5,1.25) {$\vec{i}$};
        \node[fill=white,draw,rounded corners,scale=1] at (2, -.25) {$i$};
        \node[fill=white,draw,rounded corners,scale=1] at (1.25, 4) {$\vec{j}\bullet j$};
}
\arrow[rr, "\beta_1"]
&&
\tikz[yscale=.5,xscale=.75,yshift=-3.3cm]{
        \draw (1.25, 3) -- (1.25, 6);
	\draw 
		(.5,2)
		.. controls (.5,2.5) and (1.25,2.5) ..
		(1.25,3);
	\draw (.5,0)
		.. controls (.5,1.5) and (.5,1.5) ..
		(.5,2);
	\draw (2,0)
		--
		(2,2)
		.. controls (2,2.5) and (1.25,2.5) ..
		(1.25,3);
        \node[fill=white,draw,rounded corners,scale=1] at (1.25, 3.5) {$\vec{i}\bullet i$};
        \node[fill=white,draw,rounded corners,scale=1] at (1.25, 5) {$\vec{j}\bullet j$};
}
\arrow[d, equal]
\\
\tikz[yscale=.5,xscale=.75,yshift=-1.5cm]{
        \draw (1.25, 3) -- (1.25, 4);
	\draw 
		(.5,2)
		.. controls (.5,2.5) and (1.25,2.5) ..
		(1.25,3);
	\draw (.5,-1)
		--
		(.5,1)
		.. controls (.5,1.5) and (.5,1.5) ..
		(.5,2);
	\draw (2,-1)
		--
		(2,2)
		.. controls (2,2.5) and (1.25,2.5) ..
		(1.25,3);
        \node[fill=white,draw,rounded corners,scale=1] at (.5,1.25) {$\vec{j} \circ \vec{i}$};
        \node[fill=white,draw,rounded corners,scale=1] at (2, -.25) {$j\circ i$}; 
} 
\arrow[rrr, "\beta_1"]
&&&
\tikz[yscale=.5,xscale=.75,yshift=-2.8cm]{
        \draw (1.25, 3) -- (1.25, 5);
	\draw 
		(.5,2)
		.. controls (.5,2.5) and (1.25,2.5) ..
		(1.25,3);
	\draw (.5,1)
		.. controls (.5,1.5) and (.5,1.5) ..
		(.5,2);
	\draw (2,1)
		--
		(2,2)
		.. controls (2,2.5) and (1.25,2.5) ..
		(1.25,3);
        \node[fill=white,draw,rounded corners,scale=1] at (1.25, 4) {$(\vec{j}\circ \vec{i}) \bullet (j \circ i)$};
}
\arrow[rrr, "\Xi_1"]
&&&
\tikz[yscale=.5,xscale=.75,yshift=-2.8cm]{
        \draw (1.25, 3) -- (1.25, 5);
	\draw 
		(.5,2)
		.. controls (.5,2.5) and (1.25,2.5) ..
		(1.25,3);
	\draw (.5,1)
		.. controls (.5,1.5) and (.5,1.5) ..
		(.5,2);
	\draw (2,1)
		--
		(2,2)
		.. controls (2,2.5) and (1.25,2.5) ..
		(1.25,3);
        \node[fill=white,draw,rounded corners,scale=1] at (1.25, 4) {$(\vec{j}\bullet j)\circ(\vec{i}\bullet i)$};
}
\end{tikzcd}
\]
where
\[
\zeta = \lambda\left(\abs{1_{\vec{x}} \vec{\Delta}_1 1_{\vec{y}}}, \abs{1_{\vec{y}} \Delta_2 1_z}\right).
\]
Thus, by Proposition \ref{PutyraHammer}, we see that the contributions of $\beta_1$ and $\Xi_1$ in equation (\ref{eq:shiftingsystvcomp}) is
\[
\zeta \times \text{(Left of (\ref{eq:shiftingsystvcomp}))} = \text{(Right of (\ref{eq:shiftingsystvcomp}))}.
\]

The remainder of the proof is computing the contributions of $\beta_2,\beta_3,\beta_4$, $\gamma$, and $\Xi_2$. The contributions of these on the left-hand side of (\ref{eq:shiftingsystvcomp}) are
\begin{align*}
&\lambda\left(\abs{ W_{\vec{x}\vec{y}z} ((\vec{\Delta}_2 \circ \vec{\Delta}_1) (\vec{D}), (\Delta_2 \circ \Delta_1)(D))}, V_2 + v_2\right) \cdot \lambda\left(\abs{1_{\vec{x}} \vec{\Delta}_2 1_{\vec{y}}}, v_2\right) & (\beta_{2,3,4})_{j, \vec{j}}(\varphi_{\vec{i}}(g'), \varphi_i(g)) \\
& \cdot \underbrace{\lambda\left(\abs{1_{\vec{x}} \vec{\Delta}_1 1_{\vec{y}}} + P + V_1, \abs{1_{\vec{y}} \Delta_2 1_{z}} + v_2\right)}_{(*)} &
\\
&\lambda\left(\abs{W_{\vec{x}\vec{y}z}(\vec{\Delta}_1(\vec{D}), \Delta_1(D))}, V_1 + v_1\right) \cdot \lambda\left(\abs{1_{\vec{x}} \vec{\Delta}_1 1_{\vec{y}}}, v_1\right) & (\beta_{2,3,4})_{i, \vec{i}}(g',g)
\\
& \cdot \lambda\left(P, \abs{1_{\vec{y}}\Delta_1 1_z} + v_1\right)
\\
& \underbrace{\lambda\left(\abs{1_{\vec{x}} (\vec{\Delta}_2 \bullet \Delta_2) 1_{z}}, V_1 + v_1\right)}_{(**)} & \gamma_{\vec{i}\bullet i, \vec{j} \bullet j}(g'g)
\end{align*}
We rewrite this product by expanding $(*)$ via bilinearity, expanding $(**)$ via Lemma \ref{lem:computational2} and bilinearity, and then performing the obvious cancellations; the result is the following.
\begin{align*}
&\lambda\left(\abs{ W_{\vec{x}\vec{y}z} ((\vec{\Delta}_2 \circ \vec{\Delta}_1) (\vec{D}), (\Delta_2 \circ \Delta_1)(D))}, V_2 + v_2\right)
\cdot \lambda\left(\abs{1_{\vec{x}} \vec{\Delta}_2 1_{\vec{y}}}, v_2\right)
\\
& \cdot \lambda\left(P, \abs{1_{\vec{y}} \Delta_2 1_{z}} + v_2\right) \cdot \underbrace{\lambda\left(\abs{1_{\vec{x}} \vec{\Delta}_1 1_{\vec{y}}}, \abs{1_{\vec{y}} \Delta_2 1_{z}}\right)}_{\zeta} \cdot \lambda\left(\abs{1_{\vec{x}} \vec{\Delta}_1 1_{\vec{y}}}, v_2\right) \cdot \lambda(V_1, v_2)\\
& \lambda\left(\abs{1_{\vec{x}} \vec{\Delta}_1 1_{\vec{y}}}, v_1\right) \cdot \lambda\left(P, \abs{1_{\vec{y}}\Delta_1 1_z} + v_1\right)
\\
& \lambda\left(\abs{W_{\vec{x}\vec{y}z}((\vec{\Delta}_2\circ \vec{\Delta}_1)(\vec{D}), (\Delta_2 \circ \Delta_1)(D))}, V_1 + v_1\right) \cdot \lambda\left(\abs{1_{\vec{y}} \Delta_2 1_z}, v_1\right) \cdot \lambda\left(\abs{1_{\vec{x}} \vec{\Delta}_2 1_{\vec{y}}}, V_1\right) \cdot \lambda\left(\abs{1_{\vec{x}} \vec{\Delta}_2 1_{\vec{y}}}, v_1\right)
\end{align*}

On the other hand, the contributions of $\beta_2, \beta_3, \beta_4, \gamma$ and $\Xi_2$ on the right-hand side of (\ref{eq:shiftingsystvcomp}) are
\begin{align*}
&\lambda\left(\abs{1_{\vec{x}}\vec{\Delta}_2 1_{\vec{y}}}, V_1\right) \cdot \lambda\left(\abs{1_{\vec{y}} \Delta_2 1_z}, v_1\right) & \gamma_{\vec{i}, \vec{j}}(g') \cdot \gamma_{i,j}(g)
\\
&\lambda\left(\abs{W_{\vec{x}\vec{y}z} ((\vec{\Delta}_2 \circ \vec{\Delta}_1)(\vec{D}), (\Delta_2 \circ \Delta_1)(D))}, V_1 + V_2 + v_1 + v_2\right)  & (\beta_{2,3,4})_{j\circ i, \vec{j}\circ \vec{i}}(g', g)
\\
&\cdot \lambda\left(\abs{1_{\vec{x}}(\vec{\Delta}_2 \circ \vec{\Delta}_1) 1_{\vec{y}}}, v_2 + v_1\right) \cdot \lambda\left(P, \abs{1_{\vec{y}(\Delta_2 \circ \Delta_1) 1_{z}}} + v_2 + v_2\right) &
\\
& \lambda\left(V_1, v_2\right) & \Xi_2
\end{align*}
Comparing the updated form of the left-hand side with this, we see that everything cancels except for the $\zeta$ term present in the former. Thus, we conclude that the contributions of $\beta_2, \beta_3, \beta_4, \gamma$ and $\Xi_2$ in equation (\ref{eq:shiftingsystvcomp}) is
\[
\text{(Left of (\ref{eq:shiftingsystvcomp}))} = \zeta \times \text{(Right of (\ref{eq:shiftingsystvcomp}))},
\]
which completes the proof.
\end{proof}

\subsubsection{$\mathscr{C}$-graded vertical composition}

As before, we construct natural isomorphisms $\varphi_j \circ \varphi_i \Rightarrow \varphi_{j\circ i}$ or $\varphi_{\vec{j}} \circ \varphi_{\vec{i}} \Rightarrow \varphi_{\vec{j}\circ \vec{i}}$ given by
\[
\begin{split}
(\varphi_j \circ \varphi_i)(M) &\to \varphi_{j\circ i}(M)\\
m &\mapsto \gamma_{i,j}(\abs{m}) m
\end{split}
\qquad\text{or, in general,}\qquad
\begin{split}
(\varphi_{\vec{j}} \circ \varphi_{\vec{i}})(\vec{M}) &\to \varphi_{\vec{j}\circ \vec{i}}(\vec{M})\\
\vec{m} &\mapsto \gamma_{\vec{i},\vec{j}}(\abs{\vec{m}}) \vec{m}
\end{split}
\]
respectively, and $\varphi_{(\vec{j} \circ \vec{i}) \bullet (j \circ i)} \Rightarrow \varphi_{(\vec{j} \bullet j) \circ (\vec{i} \bullet i)}$ by
\begin{align*}
(\varphi_{(\vec{j} \circ \vec{i}) \bullet (j \circ i)})(M) &\to \varphi_{(\vec{j} \bullet j) \circ (\vec{i} \bullet i)}(M)\\
m &\mapsto \Xi_{\substack{i, \vec{i} \\ j, \vec{j}}}(\abs{m}) m
\end{align*}
for all homogeneous $m\in M$. In terms of these natural isomorphisms, equations (\ref{eq:shiftingsystvcomp}) and (\ref{eq:gammarec}) translate to mean that the diagram
\[
\begin{tikzcd}[scale cd=.92, row sep = large]
\varphi_{\vec{j}} \circ \varphi_{\vec{i}} (M_1,\ldots, M_k) \otimes \varphi_j \circ \varphi_i(M)
\arrow[r, "\beta_{j, \vec{j}}"]
\arrow[d, "\gamma_{\vec{i}, \vec{j}} \otimes \gamma_{i,j}"]
&
\varphi_{\vec{j} \bullet j} \left(\varphi_{\vec{i}} (M_1,\ldots, M_k) \otimes \varphi_i(M)\right)
\arrow[r, "\beta_{i, \vec{i}}"]
&
\varphi_{\vec{j} \bullet j} \circ \varphi_{\vec{i} \bullet i} \left((M_1,\ldots, M_k) \otimes M\right)
\arrow[d, "\gamma_{\vec{i}\bullet i, \vec{j} \bullet j}"]
\\
\varphi_{\vec{j}\circ \vec{i}} (M_1,\ldots, M_k) \otimes \varphi_{j \circ i} (M)
\arrow[r, "\beta_{j\circ i, \vec{j} \circ \vec{i}}"]
&
\varphi_{(\vec{j} \circ \vec{i}) \bullet (j\circ i)} \left((M_1,\ldots, M_k)\otimes M\right)
\arrow[r, "\Xi_{\substack{i, \vec{i} \\ j, \vec{j}}}"]
&
\varphi_{(\vec{j} \bullet j) \circ (\vec{i} \bullet i)} \left((M_1,\ldots, M_k)\otimes M\right)
\end{tikzcd}
\]
commutes for all $\mathscr{C}$-graded multimodules $M_1, \ldots, M_k, M$, and 
\[
\begin{tikzcd}
\varphi_k \circ \varphi_j \circ \varphi_i (M)
\arrow[r, "\gamma_{i,j}"]
\arrow[d, "\gamma_{j,k}"']
&
\varphi_k \circ \varphi_{j\circ i} (M)
\arrow[d, "\gamma_{j\circ i, k}"]
\\
\varphi_{k\circ j, i}(M)
\arrow[r, "\gamma_{i, k\circ j}"]
&
\varphi_{k\circ j\circ i}(M)
\end{tikzcd}
\]
commutes for each $\mathscr{C}$-graded multimodule $M$.

Before moving on, we note that a shifting 2-system may be extended to $\widetilde{\mathscr{I}}$. In particular, since $\varphi_{\mathrm{Id}}$ acts as the identity, we can extend vertical composition itself by declaring $\mathrm{Id} \circ i = i = i \circ \mathrm{Id}$. If $\mathrm{Id}$ appears in the subscript of $\Xi$, it can be replaced by an compatible element in $\mathscr{I}_{\mathrm{Id}}$.

Finally, we can properly define a vertical composition of homogeneous maps. Suppose $f: M \to N$ is homogeneous of degree $i$, and $g: N \to L$ is homogeneous of degree $j$. Define their $\mathscr{C}$-graded composition as
\[
\left(g \circ_\mathscr{C} f\right)(m) = \gamma_{i,j}(\abs{m})^{-1} \left(g\circ f\right)(m).
\]

\begin{proposition}
\label{prop:C-graded comp}
With the assumptions above, $g\circ_\mathscr{C} f$ is purely homogeneous of degree $j\circ i$.
\end{proposition}

\begin{proof}
Requirement (i) of Definition \ref{def:purelyhomo} follows easily since $\mathsf{D}_{j\circ i} = \mathsf{D}_i \cap \varphi_i^{-1}(\mathsf{D_j})$. Additionally, $\abs{f} = i$ so $\abs{f(m)} = \varphi_i(\abs{m})$, and $\abs{g} = j$ so $\abs{g(f(m))} = \varphi_j \circ \varphi_i (\abs{m})$. Thus,
\[
\abs{\left(g \circ_\mathscr{C} f\right)(m)} = \gamma_{i,j}(\abs{m})^{-1} \varphi_j \circ \varphi_i (\abs{m}) = \varphi_{j\circ i}(\abs{m}),
\]
so (ii) is satisfied.

For (iii) we claim that, for any $\vec{a} = (a_1,\ldots, a_k) \in (A_1,\ldots, A_k)$,
\[
\gamma_{i,j}(\abs{m})^{-1}\beta_{\vec{e}, i}(\abs{\vec{a}}, \abs{m}) \beta_{\vec{e}, j}(\abs{\vec{a}}, \varphi_i(\abs{m})) = \beta_{\vec{e}, j\circ i}(\abs{\vec{a}}, \abs{m}) \gamma_{i,j}\left(\abs{\rho_L(\vec{a},m)}\right)^{-1}
\]
where $\vec{e} = (e, \ldots, e)$ as usual. The desired result follows easily from here, since 
\begin{align*}
\rho_L\left(\vec{a}, (g\circ_\mathscr{C} f)(m)\right) &= \rho_L\left(\vec{a}, \gamma_{i,j}(\abs{m})^{-1} g(f(m))\right) & \text{by definition,}
\\
&= \gamma_{i,j}(\abs{m})^{-1} \rho_L\left(\vec{a}, g(f(m))\right) & \text{by}~\mathbb{K}\text{-linearity of}~\rho_L,
\\
&= \gamma_{i,j}(\abs{m})^{-1} \beta_{\vec{e}, j}(\abs{\vec{a}}, \underbrace{\abs{f(m)}}_{=\varphi_i(\abs{m})}) g\left(\rho_L(\vec{a}, f(m))\right) & \text{since}~\abs{g} = j,
\\
&= \gamma_{i,j}(\abs{m})^{-1} \beta_{\vec{e}, j}(\abs{\vec{a}}, \varphi_i(\abs{m})) \beta_{\vec{e}, i}(\abs{\vec{a}}, \abs{m}) g\left(f\left(\rho_L(\vec{a}, m)\right)\right) & \text{since}~\abs{f} = i~\&~g~\text{is}~\mathbb{K}\text{-linear},
\\
&= \beta_{\vec{e}, j\circ i} (\abs{\vec{a}}, \abs{m}) \gamma_{i,j}\left(\abs{\rho_L(\vec{a},m)}\right)^{-1} (g\circ f)\left(\rho_L(\vec{a}, m)\right)& \text{by the claim, and}
\\
&= \beta_{\vec{e}, j\circ i}(\abs{\vec{a}}, \abs{m}) \left(g \circ_\mathscr{C} f\right) (\rho_L(\vec{a}, m)) & \text{by definition}.
\end{align*}
To prove the claim, we apply equation (\ref{eq:shiftingsystvcomp}) when $\vec{j} = \vec{i} = \vec{e}$ and $g' = \abs{\vec{a}}$ and $g = \abs{m}$; it reads
\begin{equation*}
\gamma_{\vec{e} \bullet i, \vec{e} \bullet j}(\abs{\vec{a}}\circ \abs{m}) 
\beta_{\vec{e}, i}(\abs{\vec{a}},\abs{m}) 
\beta_{\vec{e}, j}(\varphi_{\vec{e}}(\abs{\vec{a}}), \varphi_i(\abs{m}))
=
\Xi_{\substack{i, \vec{e} \\ j, \vec{e}}}(\abs{\vec{a}}\circ \abs{m}) \beta_{(\vec{e} \circ \vec{e}), j\circ i}(\abs{\vec{a}}, \abs{m}) 
\gamma_{\vec{e}, \vec{e}}(\abs{\vec{a}}) 
\gamma_{i,j}(\abs{m}).
\end{equation*}
Now, $\Xi_{\substack{i, \vec{e} \\ j, \vec{e}}} = 1$ since $(\vec{e} \circ \vec{e}) \bullet (j \circ i) = \vec{e} \bullet (j\circ i) = j \circ i = (\vec{e} \bullet j) \circ (\vec{e} \bullet i)$. Moreover, by our working assumption that all $\mathscr{C}$-graded algebras are supported entirely in $\Sigma$, $\varphi_{\vec{e}}(\abs{\vec{a}}) = \abs{\vec{a}}$. Then, noting $\gamma_{\vec{e}, \vec{e}} = 1$, the equation above may be rewritten
\[
\gamma_{i,j}(\abs{\vec{a}}\circ \abs{m})\beta_{\vec{e}, i}(\abs{\vec{a}}, \abs{m}) \beta_{\vec{e}, j}(\abs{\vec{a}}, \varphi_i(\abs{m})) = \beta_{\vec{e}, j\circ i}(\abs{\vec{a}}, \abs{m}) \gamma_{i,j}(\abs{m}).
\]
Note that $\abs{\rho_L(\vec{a}, m)} = \abs{\vec{a}} \circ \abs{m}$, since the action maps of multimodules are $\mathscr{C}$-graded---thus, rearranging provides the desired result. Requirement (iv) is proven in exactly the same manner, noting that $\Xi_{\substack{e, i \\ e, j}} = 1$.
\end{proof}

In general, suppose $M_\ell \xrightarrow{f_\ell} N_{\ell} \xrightarrow{g_\ell} L_{\ell}$ is a composition of purely homogeneous maps of degree $i_\ell$ and $j_\ell$ respectively, for $\ell = 1,\ldots, k$. We say that $\vec{f}$ is purely homogeneous of degree $\vec{i}$, and $\vec{g}$ is purely homogeneous of degree $\vec{j}$. Then, for $\vec{m} \in (M_1,\ldots, M_k)$, we define 
\begin{align*}
(\vec{g} \circ_\mathscr{C} \vec{f})(\vec{m}) &= \gamma_{\vec{i}, \vec{j}}(\abs{\vec{m}})^{-1} (\vec{g} \circ \vec{f})(\vec{m})
\\
&= \gamma_{\vec{i}, \vec{j}}(\abs{\vec{m}})^{-1} (g_1(f_1(m_1)), \ldots, g_k(f_k(m_k))).
\end{align*}
The proof above extends to this situation without trouble, so $(\vec{g} \circ_\mathscr{C} \vec{f})$ is purely homogeneous of degree $\vec{j} \circ \vec{i}$.

\begin{proposition}
$\mathscr{C}$-graded vertical composition is associative.
\end{proposition}

\begin{proof}
Suppose $M\xrightarrow{f} N \xrightarrow{g} L \xrightarrow{h} K$ are purely homogeneous of degrees $\abs{f} = i $, $\abs{g} = j$, and $\abs{h} = k$. On one hand,
\[
\left(h \circ_\mathscr{C} (g\circ_\mathscr{C} f)\right) (m) = \gamma_{i,j}(\abs{m})^{-1} \left(h \circ_\mathscr{C} gf\right)(m) = \gamma_{i,j}(\abs{m})^{-1} \gamma_{j\circ i, k}(\abs{m})^{-1} hgf(m).
\]
On the other,
\[
\left(\left(h \circ_\mathscr{C} g\right) \circ_\mathscr{C} f\right) (m) = \gamma_{j,k}\left(\abs{f(m)}\right)^{-1} \left(hg \circ_\mathscr{C} f\right)(m) = \gamma_{j,k}\left(\abs{f(m)}\right)^{-1} \gamma_{i, k\circ j}(\abs{m})^{-1} hgf(m).
\]
Since $\abs{f(m)} = \varphi_i(\abs{m})$, associativity follows from equation (\ref{eq:gammarec}).
\end{proof}

Propositions \ref{prop:C-graded tensor} and \ref{prop:C-graded comp} imply that the $\mathscr{C}$-graded composition and tensor product of homogeneous maps is again a homogeneous map. The last thing we must do is check the compatibility of $\otimes$ and $\circ_\mathscr{C}$.

\begin{proposition}
\label{prop:vertical&horizontalcomp}
Suppose $f:M\to N$ and $\{f_\alpha: M_\alpha \to N_\alpha\}_{\alpha=1, \ldots, k}$ are purely homogeneous maps of degree $i$ and $i_\alpha$ respectively, and similarly $g: N \to L$ and $\{g_\beta: N_\beta \to L_\beta\}_{\beta=1, \ldots, k}$ are purely homogeneous maps of degree $j$ and $j_\beta$ respectively. Then
\[
\left(\left((g_1 \circ_\mathscr{C} f_1), \ldots, (g_k \circ_\mathscr{C} f_k) 
\right) \otimes g\circ_\mathscr{C} f \right) (\vec{m} \otimes m)
=
\Xi_{\substack{i, \vec{i} \\ j, \vec{j}}} \left((g_1,\ldots, g_k) \otimes g\right) \circ_\mathscr{C} \left((f_1,\ldots, f_k) \otimes f\right) (\vec{m} \otimes m).
\]
\end{proposition}

\begin{proof}
We will just unravel both sides of the equation above. The equality will follow from equation (\ref{eq:shiftingsystvcomp}). On one hand,
\begin{align*}
( & (( g_1 \circ_\mathscr{C} f_1), \ldots, (g_k \circ_\mathscr{C} f_k)) \otimes g\circ_\mathscr{C} f) (\vec{m} \otimes m)
\\
    &= \beta_{\vec{j}\circ \vec{i}, j\circ i}(\abs{\vec{m}}, \abs{m})^{-1} \left((g_1 \circ_\mathscr{C} f_1)(m_1), \ldots, (g_k \circ_\mathscr{C} f_k)(m_k)\right)\otimes (g \circ_\mathscr{C} f)(m)
\\
    &= \beta_{\vec{j}\circ \vec{i}, j\circ i}(\abs{\vec{m}}, \abs{m})^{-1} \left(\gamma_{i_1, j_1}(\abs{m_1})^{-1}(g_1 \circ f_1)(m_1), \ldots, \gamma_{i_k, j_k}(\abs{m_1})^{-1}(g_k \circ f_k)(m_k)\right)\otimes \gamma_{i, j}(\abs{m_1})^{-1}(g \circ f)(m)
\\
    &= \beta_{\vec{j}\circ \vec{i}, j\circ i}(\abs{\vec{m}}, \abs{m})^{-1} \gamma_{\vec{i}, \vec{j}}(\abs{\vec{m}})^{-1} \gamma_{i, j}(\abs{m})^{-1} \left((g_1 \circ f_1)(m_1), \ldots, (g_k \circ f_k)(m_k)\right)\otimes (g \circ f)(m).
\end{align*}
The first equality follows from Proposition \ref{prop:C-graded comp} since each $g_\ell \circ f_\ell$ is purely homogeneous of degree $j_\ell \circ i_\ell$, so $\left((g_1 \circ_\mathscr{C} f_1), \ldots, g_k \circ_\mathscr{C} f_k \right)$ is purely homogeneous of degree $\vec{j} \circ \vec{i}$. The second equality follows from the definition of $\circ_\mathscr{C}$, while the third is just a rewriting step. On the other hand,
\begin{align*}
(((g_1,\ldots &, g_k) \otimes g) \circ_\mathscr{C} (f_1,\ldots, f_k) \otimes f)) (\vec{m} \otimes m)
\\
    &= \gamma_{\vec{i} \bullet i, \vec{j} \bullet j}(\abs{\vec{m}\otimes m})^{-1} \left(((g_1,\ldots, g_k) \otimes g) \circ ((f_1,\ldots, f_k) \otimes f) \right) (\vec{m}\otimes m)
\\
    &= \gamma_{\vec{i} \bullet i, \vec{j} \bullet j}(\abs{\vec{m}\otimes m})^{-1} ((g_1,\ldots, g_k) \otimes g) \left( \beta_{\vec{i}, i}(\abs{\vec{m}}, \abs{m})^{-1} (f_1(m_1), \ldots, f_k(m_k)) \otimes f(m) \right)
\\
    &= \gamma_{\vec{i} \bullet i, \vec{j} \bullet j} (\abs{\vec{m}\otimes m})^{-1}\beta_{\vec{i}, i}(\abs{\vec{m}}, \abs{m})^{-1} \beta_{\vec{j}, j}(\abs{\vec{f}(\vec{m})}, \abs{f(m)})^{-1}\left(g_1(f_1(m)), \ldots, g_k(f_k(m_k))\right) \otimes g(f(m)).
\end{align*}
The first equality follows from Proposition \ref{prop:C-graded tensor}, and the second and third follow from the definition of the tensor product of homogeneous maps. As suggested, the equality follows from equation (\ref{eq:shiftingsystvcomp}), taking $g' = \abs{\vec{m}}$ and $g = \abs{m}$---we must only compensate by $\Xi_{\substack{i, \vec{i} \\ j, \vec{j}}}$.
\end{proof}

\subsection{Changes of chronology}
\label{ss:gradingcocs}

An important feature of $\mathscr{G}$-shifting systems in particular is that changes of chronology induce natural transformations of grading shift functors. Recall that we have a few different notions of composition for changes of chronolology:
\begin{itemize}
    \item to a sequence of chronological cobordisms $A \xrightarrow{W} B \xrightarrow{W'} C$ and changes of chronology $H$ on $W$ and $H'$ on $W'$, there is a change of chronology $H' \circ H$ on $W' \circ W$;
    \item A sequence of changes of chronology $W_1\xRightarrow{H_1} W_2 \xRightarrow{H_2}W_3$ is itself a change of chronology, denoted $H_2 \star H_1$.
\end{itemize}
Clearly $\circ$ and $\star$ extend to chronological cobordisms with corners $\Delta$. In this setting we obtain another way of composing changes of chronology. Suppose $\Delta, \Delta_1,\ldots, \Delta_k$ are chronological cobordisms with corners so that $(\Delta_1,\ldots, \Delta_k) \bullet \Delta$ is nonzero, and suppose $H, H_1,\ldots, H_k$ are changes of chronology on $\Delta, \Delta_1,\ldots, \Delta_k$. Then we denote by $(H_1,\ldots, H_k) \bullet H$ the change of chronology on $(\Delta_1,\ldots, \Delta_k)\bullet \Delta$ defined by applying the $H_i$ and the $H$ in order according to the chronology. Indeed, we could define the $\bullet$ operation in terms of successive applications of the $\circ$ operation after extending each change of chronology to be trivial outside of its original component.

Now, each change of chronology $H: \Delta \to \Delta'$ of chronological cobordisms with corners extends to a change of chronology without corners given appropriate crossingless matchings $x_1,\ldots, x_k, y$. The latter is denoted by
\[
{}_{\vec{x}}H_y: 1_{\vec{x}} \Delta 1_y \to 1_{\vec{x}} \Delta' 1_y.
\]
We claim that this observation induces a natural transformation of grading shift functors
\[
\varphi_H: \varphi_\Delta \Rightarrow \varphi_{\Delta'}
\]
defined on each $M \in \mathrm{Ob}(\mathrm{MultiMod}^\mathscr{G})$ by
\begin{align*}
    \varphi_H(M): \varphi_\Delta(M) &\to \varphi_{\Delta'}(M)\\
    \varphi_\Delta(m) &\mapsto \iota(H(\abs{m}))^{-1} \varphi_{\Delta'}(m)
\end{align*}
where $H(\abs{m})$ means ${}_{\vec{x}}H_y$ for $\abs{m}: \vec{x} \to y$. In general, 
\[
\varphi_{(H_1, \ldots, H_k)}: \varphi_{(\Delta_1, \ldots, \Delta_k)} \Rightarrow \varphi_{(\Delta_1', \ldots, \Delta_k')}
\]
where
\[
\varphi_{(H_1,\ldots, H_k)}: (\varphi_{\Delta_1}(M_1), \ldots, \varphi_{\Delta_k}(M_k)) \to (\varphi_{\Delta_1'}(M_1), \ldots, \varphi_{\Delta_k'}(M_k))
\]
is given by
\[
\varphi_{\vec{\Delta}}(\vec{m}) \mapsto \prod_{i=1}^k \iota(H_i(\abs{m_i}))^{-1} \varphi_{\vec{\Delta}'}(\vec{m}).
\]
We abbreviate $\prod_{i=1}^k \iota(H_i(\abs{m_i}))^{-1}$ to $\iota(\vec{H}(\abs{\vec{m}}))^{-1}$. Sometimes, we write $\varphi_H$ when we mean $\varphi_H(M)$.

\begin{proposition}
\label{prop:CoCnaturality}
The diagram
\[
\begin{tikzcd}[column sep = huge, row sep = large]
\left(\varphi_{\Delta_1}(M_1), \ldots, \varphi_{\Delta_k}(M_k)\right) \otimes \varphi_\Delta(M) \arrow[d,"(\varphi_{H_1}\text{,}\, \ldots\text{,}\, \varphi_{H_k}) \otimes \varphi_H"'] \arrow[r, "\beta_{(\Delta_1,\ldots, \Delta_k), \Delta}"] & \varphi_{(\Delta_1, \ldots, \Delta_k)\bullet \Delta} \left((M_1, \ldots, M_k) \otimes M\right) \arrow[d, "\varphi_{(H_1\text{,}\, \ldots\text{,}\, H_k) \bullet H}"]
\\
\left(\varphi_{\Delta_1'}(M_1), \ldots, \varphi_{\Delta_k'}(M_k)\right) \otimes \varphi_{\Delta'}(M) \arrow[r, "\beta_{(\Delta_1',\ldots, \Delta_k'), \Delta'}"] & \varphi_{(\Delta_1', \ldots, \Delta_k')\bullet \Delta'} \left((M_1, \ldots, M_k) \otimes M\right)
\end{tikzcd}
\]
commutes. Thus, $\varphi_H(M)$ is a map of $\mathscr{G}$-graded multimodules and, in particular, $\varphi_H$ is a natural transformation of $\mathrm{MultiMod}^{\mathscr{G}}(A_1,\ldots, A_k; B)$ functors.
\end{proposition}

\begin{proof}
Assume that $(m_1,\ldots, m_k) \in (M_1,\ldots, M_k)$ and $m\in M$ such that $\abs{m_i} \in \mathrm{Hom}(x_{i1}, \ldots, x_{i\alpha_i}; y_i)$ for each $i=1,\ldots, k$ and $\abs{m} \in \mathrm{Hom}(y_1,\ldots, y_k; z)$. Recall that $\beta$ is defined as the composite $\beta_1\beta_2\beta_3\beta_4$ and notice that $(\beta_i)_{(\Delta_1, \ldots, \Delta_k), \Delta} = (\beta_i)_{(\Delta_1', \ldots, \Delta_k'), \Delta'}$ for $i=2, 3,$ and 4. Denote by $H_{\beta}$ and $H_{\beta'}$ the changes of chronology used to define $(\beta_1)_{(\Delta_1, \ldots, \Delta_k), \Delta}$ and $(\beta_1)_{(\Delta_1', \ldots, \Delta_k'), \Delta'}$ respectively. Also, consider the changes of chronology
\[
{}_{\vec{x}}\left((H_1,\ldots, H_k) \bullet H\right)_{z}: 1_{\vec{x}} (\vec{\Delta} \bullet \Delta) 1_{z} \Rightarrow 1_{\vec{x}} (\vec{\Delta}'\bullet \Delta')1_z,
\]
which we abbreviate to $H_\bullet$, and
\[
\left({}_{\vec{x}_1}(H_1)_{y_1}, \ldots, {}_{\vec{x}_k}(H_k)_{y_k}\right) \sqcup {}_{\vec{y}} H_z : 1_{\vec{x}} \vec{\Delta} 1_{\vec{y}} \sqcup 1_{\vec{y}} \Delta 1_z \Rightarrow 1_{\vec{x}} \vec{\Delta}' 1_{\vec{y}} \sqcup 1_{\vec{y}} \Delta' 1_z,
\]
which we abbreviate to $H_{\sqcup}$. Then we have the following sequences of changes of chronology.
\[
\begin{tikzcd}
\left(1_{\vec{x}}(\vec{\Delta} \bullet \Delta)1_z\right) \circ W_{\vec{x}\vec{y}z}(\vec{D}, D) 
\arrow[r, "H_{\beta_1}"]
\arrow[d, "H_\bullet \circ \mathrm{Id}"']
& 
W_{\vec{x}\vec{y}z}(\vec{\Delta}(\vec{D}), \Delta(D)) \circ \left(1_{\vec{x}}\vec{\Delta}1_{\vec{y}} \sqcup 1_{\vec{y}} \Delta 1_z \right)
\arrow[d, "\mathrm{Id} \circ H_\sqcup"]
\\
\left(1_{\vec{x}}(\vec{\Delta}' \bullet \Delta')1_z\right) \circ W_{\vec{x}\vec{y}z}(\vec{D}, D) 
\arrow[r, "H_{\beta_1}"]
& 
W_{\vec{x}\vec{y}z}(\vec{\Delta}'(\vec{D}), \Delta'(D)) \circ \left(1_{\vec{x}}\vec{\Delta}'1_{\vec{y}} \sqcup 1_{\vec{y}} \Delta' 1_z \right)
\end{tikzcd}
\]
Then, proposition \ref{PutyraHammer} implies that
\[
\iota(H_{\sqcup}) \beta_{\vec{\Delta}, \Delta} (\abs{\vec{m}}, \abs{m})
=
\beta_{\vec{\Delta}', \Delta'}(\abs{\vec{m}}, \abs{m}) \iota(H_\bullet).
\]
On the other hand,
\[
\left((\varphi_{H_1}, \ldots, \varphi_{H_k}) \otimes \varphi_H\right) (\abs{\vec{m}}\otimes \abs{m}) = \varphi_{\vec{H}}(\abs{\vec{m}})\varphi_H(\abs{m}) = \iota(\vec{H}(\abs{\vec{m}}))^{-1} \iota(H(\abs{m}))^{-1} = \iota(H_\sqcup)^{-1}
\]
and
\[
\varphi_{(H_1,\ldots, H_k)\bullet H}(\abs{m} \circ \abs{\vec{m}}) = \iota(H_\bullet)^{-1},
\]
which concludes the proof.
\end{proof}

\begin{proposition}
We have that
\[
\varphi_{H'} \circ \varphi_{H} \cong \varphi_{H'\star H} \qquad 
\]
\end{proposition}

\begin{proof}
This is immediate, since 
$\iota({}_{\vec{x}}H'_y \circ {}_{\vec{x}}H_y) = \iota({}_{\vec{x}}H'_y) \iota({}_{\vec{x}}H_y) = \iota({}_{\vec{x}}H'_y \star {}_{\vec{x}}H_y).$
\end{proof}

\newpage

\section{Tangles, dg-multimodules, and multigluing}
\label{s:tangles,multimodules,multigluing}

In this section, we finally prove multigluing in generality upon detailing our method for associating to a diskular tangle a $\mathscr{G}$-graded dg-multimodule; see \S \ref{ss:tangleresolution}. This is preceded by \S \ref{ss:dg-multimodules and hom} wherein we define $\mathscr{C}$-graded dg-multimodules, whose differential preserves $\mathscr{C}$-degree. We also define the HOM-complex associated to $\mathscr{C}$-graded dg-multimodules, important to \S \ref{Chapter:OddCKProj}. Finally, in \S \ref{ss:dgcgraded}, we also discuss graded commutativity and analogues of Naisse-Putyra's ``dg-$\mathcal{C}$-graded'' bimodules, whose differential is $\mathscr{C}$-homogeneous. We do this mostly for completeness, and cite \S \ref{ss:dgcgraded} very sparingly in successive sections.

\subsection{\texorpdfstring{$\mathscr{C}$}{Lg}-graded dg-multimodules and related concepts}
\label{ss:dg-multimodules and hom}

We remark that we only consider the situation of $\mathscr{C}$-graded dg-multimodules over $\mathscr{C}$-graded algebras, rather than over $\mathscr{C}$-graded dg-algebras.
\begin{definition}
If $A_1,\ldots, A_k, B$ are $\mathscr{C}$-graded algebras, a \textit{$\mathscr{C}$-graded dg-$(A_1,\ldots, A_k; B)$-multimodule} $(M, d_M)$ is a $\mathbb{Z} \times \mathscr{C}$-graded $(A_1,\ldots, A_k; B)$-multimodule $M = \bigoplus_{n\in\mathbb{Z}, g\in\mathrm{Mor}(\mathscr{C})} M_g^n$ together with a $\mathbb{K}$-linear map $d_M: M\to M$ satisfying
\begin{enumerate}[label = (\roman*)]
    \item $d_M(M_g^n) \subset M_g^{n+1}$,
    \item $d_M(\rho_L(\vec{a}, m)) = \rho_L(\vec{a}, d_M(m))$,
    \item $d_M(\rho_R(m,b)) = \rho_R(d_M(m), b)$, and
    \item $d_M \circ d_M = 0$.
\end{enumerate}
for all $\vec{a} \in (A_1,\ldots, A_k), b\in B$, and $m\in M$. The $\mathbb{Z}$-grading is called the \textit{homological} grading; the homological grading of $m\in M$ is denoted $\abs{m}_h$. We assume the left and right action on a multimodule preserves homological grading; i.e., $\abs{\rho_L(\vec{a}, m)}_h = \abs{m}_h$. A map of $\mathscr{C}$-graded dg-bimodules $f: M \to N$ will always mean a $\mathbb{K}$-linear chain map (i.e., it commutes with the differentials) which preserves both homological and $\mathscr{C}$-grading.
\end{definition}

Given $\mathscr{C}$-graded dg-$(A_{i1},\ldots, A_{i \alpha_i}; B_i)$-multimodules $(M_i, d_{M_i})$ for each $i=1,\ldots, k$ and a $\mathscr{C}$-graded dg-$(B_1,\ldots, B_k; C)$-multimodules $(M, d_M)$, we define a new $\mathscr{C}$-graded dg-$(A_{11}, \ldots, A_{k\alpha_k}; C)$-multimodule
\[
\left((M_1, d_{M_1}), \ldots, (M_k, d_{M_k})\right) \otimes_{(B_1,\ldots, B_k)} (M, d_M) = \left((M_1,\ldots, M_k) \otimes_{(B_1, \ldots, B_k)} M, d_{\vec{M}\otimes M}\right)
\]
where
\[
d_{\vec{M}\otimes M} (\vec{m}\otimes m) = \sum_{i=1}^k (-1)^{\sum_{j=1}^{i-1} \abs{m_j}_h} (m_1, \ldots,  d_{M_i}(m_i) , \ldots,  m_k) \otimes m + (-1)^{\sum_{i=1}^k \abs{m_i}_h} \vec{m} \otimes d_M(m).
\]
We will frequently denote the first large summation, perhaps confusingly, by simply $d_{\vec{M}}(\vec{m})$.

\begin{proposition}
The tensor product of $\mathscr{C}$-graded dg-multimodules, as defined above, is a $\mathscr{C}$-graded dg-multimodule.
\end{proposition}

\begin{proof}
The requirement (i) is obvious. Also, it is routine (but tedious) to check requirement (iv), that $d_{\vec{M} \otimes M}^2 = 0$. To see requirements (ii) and (iii), note that $d_{M_i}$ preserves $\mathscr{C}$-grading, so for any $i$,
\[
\abs{(m_1, \ldots, d_{M_i}(m_i), \ldots, m_k)} = \abs{(m_1,\ldots, m_i, \ldots, m_k)}
\]
thus, in particular,
\[
\alpha(\abs{\vec{a}}, \abs{\vec{m}}, \abs{m}) = \alpha(\abs{\vec{a}}, \abs{\vec{m}}, \abs{d_M(m)}) = \alpha(\abs{\vec{a}}, \abs{(m_1, \ldots, d_{M_i}(m_i), \ldots, m_k)}, \abs{m}).
\]
We leave the rest of the proof to the reader.
\end{proof}

The \textit{homology} of a $\mathscr{C}$-graded dg-multimodule $(M, d_M)$ is the $\mathscr{C}\times \mathbb{Z}$-graded multimodule $H(M, d_M) = \ker(d_M)\big/ \mathrm{im}(d_M)$. We call a map of $\mathscr{C}$-graded dg-multimodules $f: (M, d_M) \to (N, d_N)$ a \textit{quasi-isomorphism} if the induced map $f_*: H(M, d_M) \to H(N, d_N)$ is an isomorphism. 

Write $\mathrm{Multimod}_{dg}^\mathscr{C}(A_1,\ldots, A_k; B)$ for the category of $\mathscr{C}$-graded $(A_1,\ldots, A_k; B)$-multimodules. The \textit{derived category of $\mathscr{C}$-graded dg-$(A_1,\ldots, A_k; B)$-multimodules}, denoted $\mathsf{Multimod}_{dg}^\mathscr{C}(A_1,\ldots, A_k; B)$ is obtained from $\mathrm{Multimod}_{dg}^\mathscr{C}(A_1,\ldots, A_k; B)$ by localizing along quasi-isomorphisms. We can form the \textit{derived tensor product} of $\mathscr{C}$-graded multimodules by replacing $M$ and/or the $M_i$ in the regular tensor product with a K-flat resolution and then taking the regular tensor product.

We define the mapping cone of a map of $\mathscr{C}$-graded dg-multimodules as follows. First, recall the homological shifting functor $[k]$ which sends the dg-multimodule $(M, d_M)$ to $(M[k], d_{M[k]})$ where $M[k]_g^n = M_g^{n-k}$, $d_{M[k]} = (-1)^k d_M$, and $M[k]$ inherits the left and right actions of $M$. Then the mapping cone of  $f: (M, d_M) \to (N, d_N)$ is the $\mathscr{C}$-graded dg-multimodule
\[
\mathrm{Cone}(f) = (M[-1] \oplus N, d_{\mathrm{Cone}(f)}) \qquad \text{where} \qquad d_{\mathrm{Cone}(f)} = \begin{pmatrix} -d_M & 0 \\ f & d_N\end{pmatrix}.
\]

We also define the HOM complex of $\mathscr{C}$-graded dg-multimodules. Suppose $M$ and $N$ are two $\mathscr{C}$-graded dg-$(A_1,\ldots, A_k; B)$-multimodules. Let $\mathrm{HOM}(M,N)$ denote the chain complex of bihomogeneous (that is, homogeneous in homological degree and purely homogeneous in $\widetilde{\mathscr{I}}$-degree) maps $f$ of arbitrary $(\mathbb{Z} \times \widetilde{\mathscr{I}})$-degree, with differential 
\[
D(f) = d_N \circ f - (-1)^{\abs{f}_h} f \circ d_M.
\]
Thus, $D$ preserves the $\widetilde{\mathscr{I}}$-degree of a bihomogeneous map, but increases the homological degree by one. For example, if $f$ has degree $(k, i) \in \mathbb{Z} \times \widetilde{\mathscr{I}}$, then the differential of $\mathrm{HOM}(M,N)$ simply takes the difference of the following paths.
\[
\begin{tikzcd}
M_g^n \arrow[r, "f"] \arrow[d, "d_M"] & N_{\varphi_i(g)}^{n+k} \arrow[d, "d_N"]
\\
M_g^{n+1} \arrow[r, "f"] & N_{\varphi_i(g)}^{n + k +1}
\end{tikzcd}
\]
Recall that each purely homogeneous map of degree $i$ induces a graded map $\widetilde{f}: \varphi_i(M) \to M$. Moreover, $\mathscr{C}$-grading preserving maps can be viewed as purely homogeneous of degree $\mathrm{Id}\in \widetilde{\mathscr{I}}$; indeed, purely homogeneous maps of degree $\mathrm{Id}$ induce maps graded maps $\varphi_{\mathrm{Id}}(M) \to N$, but, $\varphi_{\mathrm{Id}}(M) = M$. This (tautological) correspondence allows us to view the HOM complex as a bigraded abelian group
\[
\mathrm{HOM}(M, N)_i^k \cong \prod_{n\in \mathbb{Z}} \mathrm{Hom}_{\mathrm{MultiMod}^{\mathscr{C}}}(\varphi_i(M^n), N^{n+k})
\]
with differential of bidegree $(1, e)$. 

\subsection{Resolution of diskular tangles}
\label{ss:tangleresolution}

A \textit{diskular $(m_1,\ldots, m_k; n)$-tangle} is a tangle diagram $T$ in $\mathbb{D}^2 - (\mathring{D}_1 \cup \cdots \cup \mathring{D}_k)$, where each of the $D_i$ are disjoint disks lying within the interior of $\mathbb{D}^2$, each of the form $\{z\in \mathbb{D}^2: \abs{z-z_i} \le r_i\}$ for some $z_i\in \mathring{\mathbb{D}}^2$ and $r_i >0$, so that
\begin{itemize}
    \item Each $D_i$ has $2m_i$ marked points on its boundary, all disjoint from a fixed basepoint in $\partial D_i$, and
    \item $\mathbb{D}^2$ itself has $2n$ marked points on its boundary, all disjoint from a fixed basepoint on $\partial \mathbb{D}^2$.
\end{itemize}
By ``$T$ is a tangle diagram in $\mathbb{D}^2 - (\mathring{D}_1 \cup \cdots \cup \mathring{D}_k)$,'' we mean that the interval components of $T$ all have endpoints lying on the marked points of $\mathbb{D}^2 - (\mathring{D}_1 \cup \cdots \cup \mathring{D}_k)$. We view the disks $D_1, \ldots, D_k$ as ordered.

As with planar arc diagrams, if $S_i$ is a diskular $(\ell_{i1},\ldots, \ell_{i\alpha_i}; m_i)$-tangle for each $i=1,\ldots, k$, we denote by $T(S_1,\ldots, S_k)$ the diskular $(\ell_{11}, \ldots, \ell_{k\alpha_k}; n)$-tangle obtained by filling the $i$th removed disk with $S_i$, identifying distinguished points and basepoints appropriately. Again, there is also a pairwise composition, which we write as $T \circ_i S_i$, and the two are related by
\[
T(S_1, \ldots, S_k) = ( \cdots ((T\circ_k S_k) \circ_{k-1} \cdots ) \circ_1 S_1.
\]
A diskular $(; n)$-tangle is referred to as a \textit{diskular $n$-tangle}.

Let $c(T)$ denote the number of crossings in $T$ and take an ordering $\chi(T) = \{\chi_1,\ldots, \chi_{c(T)}\}$ of the crossings of $T$. Let $v = (v_1,\ldots, v_{c(T)}): \chi(T) \to \{0,1\}^{c(T)}$ be an assignment of 0 or 1 to each crossing of $T$. To each $v$, thought of as the coordinates of the vertices of the hypercube $[0,1]^{c(T)}$, we associate a planar arc diagram $T_v$ of type $(m_1, \ldots, m_k; n)$ by resolving each crossing according to the following rule.
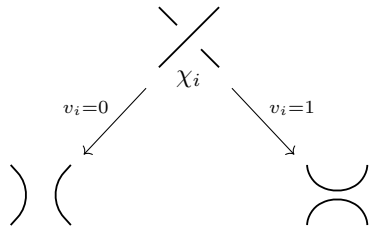
\[
\begin{tikzcd}
& \tikz[baseline=1.6ex, scale = .8]{
\draw[knot] (0,0) -- (1,1);
\draw[knot] (1,0) -- (.7,.3);
\draw[knot] (.3,.7) -- (0,1);
\node at (0.5, -0.25) {$\chi_i$};
} \arrow[dl, "v_i = 0"'] \arrow[dr, "v_i = 1"] & 
\\
\tikz[baseline=1.6ex, scale = .8]{
\draw[knot, rounded corners = 4mm] (0,0) -- (.45,.5) -- (0,1);
\draw[knot, rounded corners = 4mm] (1,0) -- (.55,.5) -- (1,1);
} & & \tikz[baseline=1.6ex, scale = .8]{
\draw[knot, rounded corners = 2mm] (0,0) -- (0,.25) -- (.5,.49) -- (1,.25) -- (1,0);
\draw[knot, rounded corners = 2mm] (0,1) -- (0,.75) -- (.5,.51) -- (1,.75) -- (1,1);
}
\end{tikzcd}
\]
We call $T_v$ a resolution of $T$. As this procedure associates planar arc diagrams to each vertex of the cube $[0,1]^{c(T)}$, we can associate to each edge a cobordism of planar arc diagrams. First, to ensure this cobordism comes with a chronology, we require that $T$ come labeled with one of
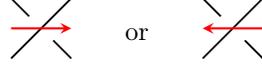
\[
\tikz[baseline=1.6ex, scale = .8]{
\draw[knot] (0,0) -- (1,1);
\draw[knot] (1,0) -- (.7,.3);
\draw[knot] (.3,.7) -- (0,1);
\draw[-stealth, red, thick] (0,0.5) -- (1, 0.5);
}
\qquad \text{or} \qquad
\tikz[baseline=1.6ex, scale = .8]{
\draw[knot] (0,0) -- (1,1);
\draw[knot] (1,0) -- (.7,.3);
\draw[knot] (.3,.7) -- (0,1);
\draw[stealth-, red, thick] (0,0.5) -- (1, 0.5);
}
\]
at each crossing. For each $v_i = 0$ in some vertex $v$, we write $v + i$ to denote the vertex which is identical to $v$ except that $(v+i)_i = 1$. Introduce a direction on the edges of the cube so that $v \rightarrow v+i$. Finally, to each of these edges, we associate the chronological cobordism
\[
W_{v,i}: T_v \to T_{v+i}
\]
obtained by putting a saddle in a small cylinder above the 0-resolution of the $i$th crossing with chronology determined by the labeling, and taking the identity everywhere outside of this cylinder.

Our goal is to assign to each diskular $(m_1,\ldots, m_k; n)$-tangle $T$ a $\mathscr{G}$-graded dg-$(H^{m_1}, \ldots, H^{m_k}; H^{n})$-multimodule $\mathcal{F}(T)$. We have already seen that to each $T_v$, $\mathcal{F}(T_v)$ is a $\mathscr{G}$-graded $(H^{m_1}, \ldots, H^{m_k}; H^{n})$-multimodule. Also, to each edge cobordism $W_{v,i}: T_v \to T_{v+i}$, we can associate a $\mathscr{G}$-graded map
\[
\mathcal{F}(W_{v,i}): \varphi_{W_{v,i}}(\mathcal{F}(T_v)) \to \mathcal{F}(T_{v+i}).
\]
We will need a slightly different graded map, achieved by constructing another family of chronological cobordisms for each $v$. Denote by $\underline{1}$ the ``all one'' vertex $(1, \ldots, 1)$. Recursively, set $W_{\underline{1}} = \mathbbm{1}_{T_{\underline{1}}}$, the identity cobordism of $T_{\underline{1}}$. For $v \not= \underline{1}$, let $\ell$ denote the lowest integer so that $v_\ell = 0$. Then, define
\[
W_v := W_{v+\ell} \circ W_{v, \ell}
\]
which has path $T_v \xrightarrow{W_{v,\ell}} T_{v+\ell} \xrightarrow{W_{v+\ell}} T_{\underline{1}}$. Additionally, notice that for each $v_j = 0$, there is a locally vertical change of chronology
\[
H_{v,j}: W_v \Rightarrow W_{v+j} \circ W_{v,j}
\]
obtained by pushing the saddle over the $j$th crossing to the beginning of the sequence of saddles.

Now, set
\[
C(T)_r = \bigoplus_{\abs{v} = r} C(T)_v[r] \qquad \text{where} \qquad C(T)_v = \varphi_{W_v}(\mathcal{F}(T_v)).
\]
Here, $r$ is the homological index of the dg-bimodule we are building. The first step in defining the differential is to associate to each edge $v \to v+j$ the $\mathscr{G}$-graded map
\[
d_{v,j} = \mathcal{F}(W_{v,j}) \circ \varphi_{H_{v,j}}(\mathcal{F}(T_v)): C(T)_v \to C(T)_{v+j}.
\]
Perhaps this doesn't seem to make sense. Indeed, there should be an intermediary $\gamma_{W_{v+j}, W_{v,j}}$ for the composition to parse:
\[
\varphi_{W_v}(\mathcal{F}(T_v)) \xrightarrow{\varphi_{H_{v,j}}(\mathcal{F}(T_v))} \varphi_{W_{v+j} \circ W_{v,j}} (\mathcal{F}(T_v)) 
\xrightarrow{\gamma_{W_{v+j}, W_{v,j}}} \varphi_{W_{v+j}}\left(\varphi_{W_{v,j}} \left(\mathcal{F}(T_v)\right) \right)
\xrightarrow{\mathcal{F}(W_{v,j})} \varphi_{W_{v+j}} (\mathcal{F}(T_{v+j})).
\]
However, by our definition of these compatibility maps, $\gamma_{W_{v+j}, W_{v,j}} = 1$ since both cobordisms involved are unweighted. Actually, the grading shifting system imposed on this grading multicategory implies that $\varphi_{W_{v+j} \circ W_{v,j}} = \varphi_{W_{v+j}} \circ \varphi_{W_{v,j}}$.

\begin{lemma}[\cite{naisse2020odd}, Lemma 6.7]
The diagram
\[
\begin{tikzcd}
& C(T)_{v+i} \arrow[dr, "d_{v+i,j}"] & \\
C(T)_v \arrow[ur, "d_{v,i}"] \arrow[dr, "d_{v,j}"'] & & C(T)_{v+i+j}\\
& C(T)_{v+j} \arrow[ur, "d_{v+j, i}"'] & 
\end{tikzcd}
\]
commutes for all $v$ and $i,j$ for which $v_i = v_j = 0$.
\end{lemma}

The proof of this Lemma is exactly as Naisse-Putyra. Indeed, the validity of this Lemma, without sign assignments, is the first meaningful benefit of working with grading (multi)categories.

Finally, define $d_r: C(T)_r \to C(T)_{r+1}$ by setting
\[
d_r|_{C(T)_v} = \sum_{\{j : v_j = 0\}} (-1)^{p(v,j)} d_{v,j}
\]
for all $v$ with $\abs{v} = r$, where $p(v,j) = \{\ell : j < \ell \le c(T)~\text{and}~v_\ell = 1\}$ counts the number of 1-resolutions occurring after the $j$th entry of $v$. In conclusion, we set
\[
\mathcal{F}(T) = \left(\bigoplus_r C(T)_r, d = \sum_r d_r\right).
\]

The following is apparent, but we write it as a proposition for future reference.

\begin{proposition}
\label{Prop:cofibration}
Suppose $T$ is a diskular tangle. Given a specified crossing of $T$, write $T_i$, for $i=0,1$, to denote the diskular tangles resulting from taking the $i$th resolution of this crossing. Write $\sigma$ to denote the saddle from $T_0$ to $T_1$. Then,
\[
\mathcal{F}(T) \cong \mathrm{Cone}
\left(\varphi_{\sigma}
\mathcal{F}(T_0)
\xrightarrow{\mathcal{F}(\sigma)}
\mathcal{F}(T_1)
\right).
\]
Equivalently, we have an exact triangle
\[
\mathcal{F}(T_1) \to \mathcal{F}(T) \to \varphi_{\sigma}\mathcal{F}(T_0)[1].
\]
\end{proposition}

Less apparent is the fact that $\mathcal{F}(T)$ is actually a $\mathscr{C}$-graded dg-multimodule.

\begin{proposition}
\label{Prop:diskularisDG}
If $T$ is a diskular $(m_1,\ldots, m_k; n)$-tangle, $\mathcal{F}(T)$ is a $\mathscr{G}$-graded dg-$(H^{m_1}, \ldots, H^{m_k}; H^n)$-multimodule.
\end{proposition}

\begin{proof}
It is clear that $d(\mathcal{F}(T)_g^\ell) \subset \mathcal{F}(T)_g^{\ell+1}$ and $d^2 = 0$ by definition. We will show that $d(\rho_L(\vec{a}, u)) = \rho_L(\vec{a}, d(u))$; the requirement for the right action follows by a similar argument. By linearity, it suffices to show that the diagram
\[
\begin{tikzcd}[column sep = large]
(A_1,\ldots, A_k) \otimes \varphi_{W_v}(\mathcal{F}(T_v)) 
\arrow[d, "1 \otimes d_{v,j}"'] \arrow[r, "\varphi_{W_v}\rho_L^v"] & \varphi_{W_v}(\mathcal{F}(T_v)) \arrow[d, "d_{v,j}"]\\
(A_1,\ldots, A_k) \otimes \varphi_{W_{v+j}}(\mathcal{F}(T_{v+j})) \arrow[r, "\varphi_{W_{v+j}}\rho_L^{v+j}"] & \varphi_{W_{v+j}}(\mathcal{F}(T_{v+j}))
\end{tikzcd}
\]
commutes, where $\rho_L^v$ denotes $\mu[(1_{m_1}, \ldots, 1_{m_k}); T_v]$. By definition of $d_{v,j}$ and left actions on shifted multimodules, this diagram factors as follows (we've refrained from labeling arrows to avoid clutter.
\[
\begin{tikzcd}[row sep = huge]
\varphi_{\vec{e}}(A_1,\ldots, A_k) \otimes \varphi_{W_v}(\mathcal{F}(T_v)) \MySymbA{dr} \arrow[d] \arrow[r] 
&
\varphi_{W_v}((A_1,\ldots, A_k) \otimes \mathcal{F}(T_v)) \MySymbB{dr} \arrow[d] \arrow[r]
& 
\varphi_{W_v}(\mathcal{F}(T_v)) \arrow[d]
\\
\varphi_{\vec{e}}(A_1,\ldots, A_k) \otimes \varphi_{W_{v+j} \circ W_{v,j}}(\mathcal{F}(T_v)) \MySymbC{dr} \arrow[d] \arrow[r] 
&
\varphi_{W_{v+j} \circ W_{v,j}} ((A_1,\ldots, A_k) \otimes \mathcal{F}(T_v)) \MySymbD{dr} \arrow[d] \arrow[r]
& 
\varphi_{W_{v+j} \circ W_{v,j}}(\mathcal{F}(T_v)) \arrow[d]
\\
\varphi_{\vec{e}}(A_1,\ldots, A_k) \otimes \varphi_{W_{v+j}}(\mathcal{F}(T_{v+j})) \arrow[r]
&
\varphi_{W_{v+j}}((A_1,\ldots, A_k) \otimes \mathcal{F}(T_{v+j})) \arrow[r]
&
\varphi_{W_{v+j}}(\mathcal{F}(T_{v+j}))
\end{tikzcd}
\]
Here, we are using the fact that $(A_1,\ldots, A_k) = \varphi_{\vec{e}}(A_1,\ldots A_k)$. We will show that the original diagram commutes by showing that squares $\circled{1}$--$\circled{4}$ commute up to constants which cancel with one another.

Square $\circled{1}$,
\[
\begin{tikzcd}[column sep = huge, row sep = huge]
\varphi_{\vec{e}}(A_1,\ldots, A_k) \otimes \varphi_{W_v}(\mathcal{F}(T_v)) \arrow[d, "1 \otimes \varphi_{H_{v,j}}(\mathcal{F}(T_v))"'] \arrow[r, "\beta_{\vec{e}, W_v}"] 
&
\varphi_{W_v}((A_1,\ldots, A_k) \otimes \mathcal{F}(T_v)) \arrow[d, "\varphi_{H_{v,j}}(\mathcal{F}(T_v))"]
\\
\varphi_{\vec{e}}(A_1,\ldots, A_k) \otimes \varphi_{W_{v+j} \circ W_{v,j}}(\mathcal{F}(T_v)) \arrow[r, "\beta_{\vec{e}, W_{v+j} \circ W_{v,j} }"] 
&
\varphi_{W_{v+j} \circ W_{v,j}} ((A_1,\ldots, A_k) \otimes \mathcal{F}(T_v)) 
\end{tikzcd}
\]
commutes on the nose by Proposition \ref{prop:CoCnaturality}, taking $\vec{\Delta} = \vec{e}$ and $\Delta = W_v$. Technically, if $h$ is the ``do nothing'' change of chronology, we are acting by $\varphi_{\vec{h}}(\vec{A})$ on the left terms, but this is clearly equal to 1. Similarly, the vertical arrow on the right should be $\varphi_{\vec{h} \bullet H_{v,j}}$.

Square $\circled{2}$,
\[
\begin{tikzcd}[column sep = huge, row sep = huge]
\varphi_{W_v}((A_1,\ldots, A_k) \otimes \mathcal{F}(T_v)) \arrow[d, "\varphi_{H_{v,j}}(\mathcal{F}(T_v))"'] \arrow[r, "\rho_L^v"]
& 
\varphi_{W_v}(\mathcal{F}(T_v)) \arrow[d, "\varphi_{H_{v,j}}(\mathcal{F}(T_v))"]
\\
\varphi_{W_{v+j} \circ W_{v,j}} ((A_1,\ldots, A_k) \otimes \mathcal{F}(T_v)) \arrow[r, "\rho_L^v"]
& 
\varphi_{W_{v+j} \circ W_{v,j}}(\mathcal{F}(T_v))
\end{tikzcd}
\]
commutes by the naturality of $\varphi_{H_{v,j}}$; again, see Proposition \ref{prop:CoCnaturality}.

Square $\circled{3}$,
\[
\begin{tikzcd}[column sep = huge, row sep = huge]
\varphi_{\vec{e}}(A_1,\ldots, A_k) \otimes \varphi_{W_{v+j} \circ W_{v,j}}(\mathcal{F}(T_v)) \arrow[d, "1 \otimes \mathcal{F}(W_{v,j})"'] \arrow[r, "\beta_{\vec{e}, W_{v+j} \circ w_{v,j}}"] 
&
\varphi_{W_{v+j} \circ W_{v,j}} ((A_1,\ldots, A_k) \otimes \mathcal{F}(T_v)) \arrow[d, "1 \otimes \mathcal{F}(W_{v,j})"]
\\
\varphi_{\vec{e}}(A_1,\ldots, A_k) \otimes \varphi_{W_{v+j}}(\mathcal{F}(T_{v+j})) \arrow[r, "\beta_{\vec{e}, W_{v+j}}"]
&
\varphi_{W_{v+j}}((A_1,\ldots, A_k) \otimes \mathcal{F}(T_{v+j}))
\end{tikzcd}
\]
commutes up to a factor of $\beta_{\vec{e}, W_{v,j}}(\abs{\vec{a}}, \abs{u})$, where we've fixed $\vec{a} \in (A_1,\ldots, A_k)$ and $u \in \mathcal{F}(T_v)$. To see this, recall that $\beta$ decomposes into 4 terms, $\beta_1$--$\beta_4$, and that here $\beta_2 = \beta_3 = 1$ for both compatibility maps since all cobordisms involved are unweighted. Otherwise, suppose $\abs{a_i}: x_i \to y_i$ and $\abs{u}: (y_1,\ldots, y_k) \to z$. Note that if $\abs{u}: \vec{y} \to x$ then $\abs{\varphi_W(u)}: \vec{y} \to x$ and $\abs{\mathcal{F}(\Delta)(u)} \vec{y} \to x$, as long as the values are nonzero. Then
\[
\left(\beta_{\vec{e}, W_{v+j} \circ W_{v,j}}(\abs{\vec{a}}, \abs{u})\right)_4 = \lambda(P', \abs{1_{\vec{y}} (W_{v+j} \circ W_{v,j}) 1_z})
\]
and
\[
\left(\beta_{\vec{e}, W_{v+j}}(\abs{\vec{a}}, \abs{\mathcal{F}(W_{v,j})(u)})\right)_4 = \lambda(P', \abs{1_{\vec{y}}(W_{v+j}) 1_z})
\]
where $P'$ is the sum of the second coordinaters of $\vec{a_i}$ for $i=1,\ldots, k$. Bilinearity of $\lambda$ implies that the contribution from the $\beta_4$ terms is
\[
\lambda(P', \abs{1_{\vec{y}} (W_{v,j}) 1_z}) \times \text{(down, then right)} = \text{(right, then down)}.
\]
On the other hand, the $\beta_1$ terms are computed via changes of chronology. Similar to before, we have that 
\[
\left(\beta_{\vec{e}, W_{v+j} \circ W_{v,j}}(\abs{\vec{a}}, \abs{u})\right)_1 = 
\left(\beta_{\vec{e}, W_{v,j}}(\abs{\vec{a}}, \abs{u})\right)_1
\times
\left(\beta_{\vec{e}, W_{v+j}}(\abs{\vec{a}}, \abs{\mathcal{F}(W_{v,j})(u)})\right)_1.
\]
The easiest way to see this is by noticing that the change of chronology on the left factors into changes of chronologies corresponding to the right terms.
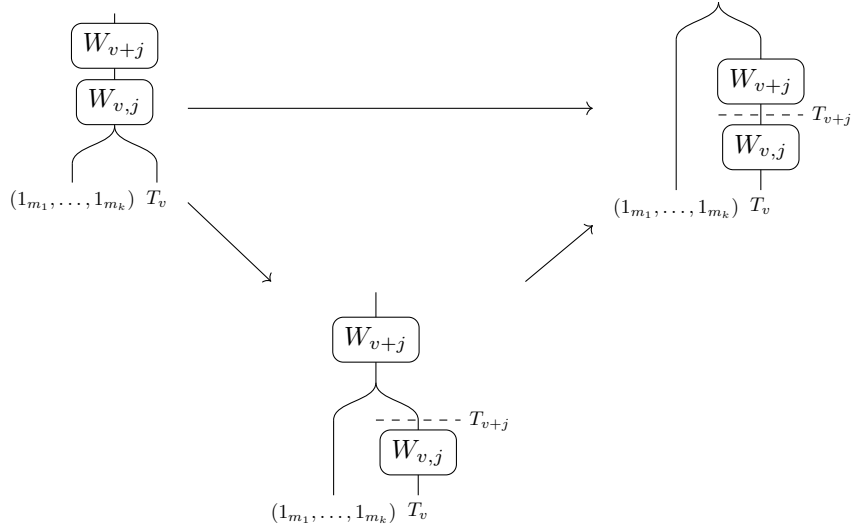
\[
\begin{tikzcd}
\tikz[yscale=.5,xscale=.75,yshift=-3.3cm]{
        \draw (1.25, 3) -- (1.25, 6);
	\draw 
		(.5,2)
		.. controls (.5,2.5) and (1.25,2.5) ..
		(1.25,3);
	\draw (.5,1.5) node[below,scale=.75]{$(1_{m_1}, \ldots, 1_{m_k})$}
		.. controls (.5,1.5) and (.5,1.5) ..
		(.5,2);
	\draw (2,1.5) node[below,scale=.75]{$T_v$}
		--
		(2,2)
		.. controls (2,2.5) and (1.25,2.5) ..
		(1.25,3);
        \node[fill=white,draw,rounded corners,scale=1] at (1.25, 3.5) {$W_{v,j}$};
        \node[fill=white,draw,rounded corners,scale=1] at (1.25, 5) {$W_{v+j}$};
} \arrow[dr] \arrow[rr]
&&
\tikz[yscale=.5,xscale=.75]{
	\draw 
		(.5,2)
		.. controls (.5,2.5) and (1.25,2.5) ..
		(1.25,3);
	\draw (.5,-2) node[below,scale=.75]{$(1_{m_1}, \ldots, 1_{m_k})$}
		-- 
		(.5,-.8)
		--
		(.5,1)
		.. controls (.5,1.5) and (.5,1.5) ..
		(.5,2);
	\draw (2,-2) node[below,scale=.75]{$T_v$}
		--
		(2,2)
		.. controls (2,2.5) and (1.25,2.5) ..
		(1.25,3);
        \draw[dashed] (1.25,0) -- (2.75,0) node[right,scale=.75]{$T_{v+j}$};
        \node[fill=white,draw,rounded corners,scale=1] at (2, .75) {$W_{v+j}$};
        \node[fill=white,draw,rounded corners,scale=1] at (2, -1) {$W_{v,j}$};
}
\\
&
\tikz[yscale=.5,xscale=.75,yshift=-2.5cm]{
        \draw (1.25, 3) -- (1.25, 5.4);
	\draw 
		(.5,2)
		.. controls (.5,2.5) and (1.25,2.5) ..
		(1.25,3);
	\draw (.5,0) node[below,scale=.75]{$(1_{m_1},\ldots, 1_{m_k})$}
		--
		(.5,1)
		.. controls (.5,1.5) and (.5,1.5) ..
		(.5,2);
	\draw (2,0) node[below,scale=.75]{$T_v$}
		--
		(2,2)
		.. controls (2,2.5) and (1.25,2.5) ..
		(1.25,3);
        \draw[dashed] (1.25,2) -- (2.75,2) node[right,scale=.75]{$T_{v+j}$};
        \node[fill=white,draw,rounded corners,scale=1] at (2, 1) {$W_{v,j}$};
        \node[fill=white,draw,rounded corners,scale=1] at (1.25, 4) {$W_{v+j}$};
} \arrow[ur]
&
\end{tikzcd}
\]
Together, this means that 
\[
\left(\beta_{\vec{e}, W_{v,j}}(\abs{\vec{a}}, \abs{u})\right)
\times \text{(down, then right)} = \text{(right, then down)}.
\]

Finally, square $\circled{4}$,
\[
\begin{tikzcd}[column sep = huge, row sep = huge]
\varphi_{W_{v+j} \circ W_{v,j}} ((A_1,\ldots, A_k) \otimes \mathcal{F}(T_v)) \arrow[d, "1 \otimes \mathcal{F}(W_{v,j})"'] \arrow[r, "\rho_L^v"]
& 
\varphi_{W_{v+j} \circ W_{v,j}}(\mathcal{F}(T_v)) \arrow[d, "\mathcal{F}(W_{v,j})"]
\\
\varphi_{W_{v+j}}((A_1,\ldots, A_k) \otimes \mathcal{F}(T_{v+j})) \arrow[r, "\rho_L^{v+j}"]
&
\varphi_{W_{v+j}}(\mathcal{F}(T_{v+j}))
\end{tikzcd}
\]
commutes up to a factor of $\beta_{\vec{e}, W_{v,j}}(\abs{\vec{a}}, \abs{u})$. To see this, recall that $\mathcal{F}(W_{i,j})$ is homogeneous of degree $W_{i,j}$, hence
\[
\rho_L^{v+j}(\vec{a}, \mathcal{F}(W_{v,j})(u)) = \beta_{\vec{e}, W_{v,j}}(\abs{\vec{a}}, \abs{u}) \mathcal{F}(W_{v,j})(\rho_L^v(\vec{a}, u)).
\]
These two contributions of $\beta_{\vec{e}, W_{v,j}}$ cancel each other out, which concludes the proof.
\end{proof}

\subsubsection{Multigluing}

Finally, we prove that $\mathcal{F}$ behaves as we hope with respect to composition of tangles; this isomorphism is referred to as multigluing.

\begin{theorem}
\label{thm:multigluing}
Suppose $T$ is a diskular $(m_1,\ldots, m_k; n)$-tangle and $T_i$ is a diskular $(\ell_{i1}, \ldots, \ell_{i\alpha_i}; m_i)$ tangle for each $i = 1,\ldots, k$. Then there is an isomorphism
\[
\left(\mathcal{F}(T_1), \ldots, \mathcal{F}(T_k)\right)\otimes_{(H^{m_1}, \ldots, H^{m_k})} \mathcal{F}(T) \cong \mathcal{F}(T(T_1,\ldots, T_k))
\]
induced by $\mu[((T_1)_{v_1}, \ldots, (T_k)_{v_k}); T_v]$.
\end{theorem}

\begin{proof}
Recall that $C(T)_v = \varphi_{W_v}\mathcal{F}(T_v)$ and  $C(T_i)_{v_i} = \varphi_{W_{v_i}} \mathcal{F}(T_i)_{v_i}$. We'll write $\otimes_{(H^{m_1},\ldots, H^{m_k})}$ as $\otimes_{\vec{m}}$. First, notice that
{\footnotesize
\begin{align*}
\left(\varphi_{W_{v_1}}\mathcal{F}(T_1)_{v_1}, \ldots, \varphi_{W_{v_k}} \mathcal{F}(T_k)_{v_k}\right) \otimes_{\vec{m}} \varphi_{W_v}\mathcal{F}(T_v)
&\xrightarrow{\beta_{(W_{v_1},\ldots, W_{v_k}), W_v}}
\varphi_{(W_{v_1}, \ldots, W_{v_k})\bullet W_v} \left((\mathcal{F}(T_1)_{v_1}, \ldots, \mathcal{F}(T_k)_{v_k})\otimes_{\vec{m}} \mathcal{F}(T_v)\right)
\\
&\xrightarrow{\mu[((T_1)_{v_1},\ldots, (T_k)_{v_k}); T_v]}
\varphi_{(W_{v_1}, \ldots, W_{v_k})\bullet W_v} \mathcal{F}(T_v((T_1)_{v_1}, \ldots, (T_k)_{v_k}))
\end{align*}
}%
is an isomorphism thanks to Proposition \ref{prop:compisomorphism}. This composition is what we mean by ``the map induced by $\mu[((T_1)_{v_1}, \ldots, (T_k)_{v_k}); T_v]$''; we will denote it by $\mu^\star$ when there is no confusion. Notice that the target of this composition can be rewritten
\[
\varphi_{W_{(v, v_1, \ldots, v_k)}} \mathcal{F}(T(T_1,\ldots, T_k)_{(v, v_1, \ldots, v_k)}) = C(T(T_1,\ldots, T_k))_{(v, v_1, \ldots, v_k)}
\]
where we've ordered the crossings of $T(T_1,\ldots, T_k)$ by the crossings of $T$ first, and then the crossings of $T_1, T_2$, and so on. So, to conclude the proof, we need only show the diagrams
\[
\begin{tikzcd}
    \left(C(T_1)_{v_1}, \ldots, C(T_i)_{v_i}, \ldots, C(T_k)_{v_k}\right) \otimes_{\vec{m}} C(T)_v 
    \arrow[r, "\mu^\star"]
    \arrow[d, "(1\text{,} \ldots\text{,} d_{v_i, j}\text{,} \ldots\text{,} 1) \otimes 1"']
    & 
    C(T(T_1,\ldots, T_k))_{(v, v_1,\ldots, v_i, \ldots, v_k)}
    \arrow[d, "d_{(v, v_1,\ldots, v_k), c + c_1 + \ldots + c_{i-1}+j}"]
    \\
    \left(C(T_1)_{v_1}, \ldots, C(T_i)_{v_i+j}, \ldots, C(T_k)_{v_k}\right)\otimes_{\vec{m}} C(T)_v
    \arrow[r, "\mu^\star"]
    &
    C(T(T_1,\ldots, T_k))_{(v, v_1,\ldots, v_i + j, \ldots, v_k)}
\end{tikzcd}
\]
(where $c$ is the number of crossings of $T$ and $c_i$ is the number of crossings of $T_i$) and
\[
\begin{tikzcd}
    \left(C(T_1)_{v_1}, \ldots, C(T_k)_{v_k}\right) \otimes_{\vec{m}} C(T)_v 
    \arrow[r, "\mu^\star"]
    \arrow[d, "(1\text{,} \ldots\text{,} 1) \otimes d_{v,j}"']
    & 
    C(T(T_1,\ldots, T_k))_{(v, v_1,\ldots, v_k)}
    \arrow[d, "d_{(v, v_1,\ldots, v_k),j}"]
    \\
    \left(C(T_1)_{v_1}, \ldots, C(T_k)_{v_k}\right)\otimes_{\vec{m}} C(T)_{v+j}
    \arrow[r, "\mu^\star"]
    &
    C(T(T_1,\ldots, T_k))_{(v+j, v_1, \ldots, v_k)}
\end{tikzcd}
\]
commute. As in the proof of Proposition \ref{Prop:diskularisDG}, we will show that each square factors into squares which commute up to values which cancel.

We introduce the following notation: we'll write
\begin{itemize}
    \item $\varphi_{W_{v_i}} = \varphi_i$ and $\mathcal{F}(T_i)_{v_i} = C_i$, so that $C(T_i)_{v_i} = \varphi_{W_{v_i}}\mathcal{F}(T_i)_{v_i}$ can be written $\varphi_iC_i$;
    \item $\varphi_{W_v} = \varphi_0$ and $\mathcal{F}(T)_{v} = C_0$, so that $C(T)_{v} = \varphi_{W_{v}}\mathcal{F}(T)_{v}$ can be written $\varphi_0C_0$;
    \item $\varphi_{W_{v_i + j}\circ W_{v_i,j}} = \varphi_{i'}$ and $\varphi_{W_{v_i + j}} = \varphi_{i''}$. Similarly, $\varphi_{W_{v + j}\circ W_{v,j}} = \varphi_{0'}$ and $\varphi_{W_{v + j}} = \varphi_{0''}$.
    \item $C = \mathcal{F}(T_v((T_1)_{v_1}, \ldots, (T_i)_{v_i}, \ldots, (T_k)_{v_k}))$, $C' = \mathcal{F}(T_v((T_1)_{v_1}, \ldots, (T_i)_{v_i+j}, \ldots, (T_k)_{v_k}))$, and $C'' = \mathcal{F}(T_{v+j}((T_1)_{v_1}, \ldots, (T_k)_{v_k}))$.
\end{itemize}
Other notation is defined accordingly; for example, $\varphi_{(W_{v_1}, \ldots, W_{v_i}, \ldots, W_{v_k})\bullet W_v}$ is rewritten $\varphi_{(1, \ldots, i, \ldots, k) \bullet 0}$, and so on. The maps involved also adapt, including the writing of $\varphi_{H_i}$ for $\varphi_{H_{v_i, j}}$, and $\mathcal{F}_i$ and $\mathcal{F}_i'$ for $\mathcal{F}(W_{v_i, j})$ and $\mathcal{F}(W_{(v, v_1,\ldots, v_k), c + c_1 + \cdots + c_{i-1} + j})$.

With this new notation, the first diagram factorizes as follows.
\[
\begin{tikzcd}[column sep = huge, row sep = huge, scale cd=.8]
    \left(\varphi_1C_1, \ldots, \varphi_iC_i, \ldots, \varphi_k C_k\right) \otimes_{\vec{m}} \varphi_0 C_0
    \MySymbA{dr}
    \arrow[r, "\beta_{(1,\ldots, i, \ldots, k), 0}"]
    \arrow[d, "(\mathbbm{1}\text{,} \ldots\text{,} \varphi_{H_i}\text{,} \ldots\text{,} \mathbbm{1})\otimes \mathbbm{1}"]
    &
    \varphi_{(1,\ldots, i,\ldots, k)\bullet 0} \left(C_0, \ldots, C_i, \ldots, C_k\right) \otimes_{\vec{m}} C_0
    \MySymbB{dr}
    \arrow[r, "\mu"]
    \arrow[d, "\varphi_{(\mathbbm{1}\text{,}\ldots\text{,} H_i\text{,}\ldots\text{,}\mathbbm{1})\bullet \mathbbm{1}}"]
    &
    \varphi_{(1,\ldots, i,\ldots, k)\bullet 0} C
    \arrow[d, "\varphi_{(\mathbbm{1}\text{,}\ldots\text{,} H_i\text{,}\ldots\text{,}\mathbbm{1})\bullet \mathbbm{1}}"']
    \\
    \left(\varphi_1C_1, \ldots, \varphi_{i'}C_i, \ldots, \varphi_k C_k\right) \otimes_{\vec{m}} \varphi_0 C_0
    \MySymbC{dr}
    \arrow[r, "\beta_{(1,\ldots, i', \ldots, k), 0}"]
    \arrow[d, "(\mathbbm{1}\text{,} \ldots\text{,} \mathcal{F}_i\text{,} \ldots\text{,} \mathbbm{1})\otimes \mathbbm{1}"]
    &
    \varphi_{(1,\ldots, i',\ldots, k)\bullet 0} \left(C_0, \ldots, C_i, \ldots, C_k\right) \otimes_{\vec{m}} C_0
    \MySymbD{dr}
    \arrow[r, "\mu"]
    \arrow[d, "(\mathbbm{1}\text{,} \ldots\text{,} \mathcal{F}_i\text{,} \ldots\text{,} \mathbbm{1})\otimes \mathbbm{1}"]
    &
    \varphi_{(1,\ldots, i',\ldots, k)\bullet 0} C
    \arrow[d, "\mathcal{F}_i'"']
    \\
    \left(\varphi_1C_1, \ldots, \varphi_{i''}C_i', \ldots, \varphi_k C_k\right) \otimes_{\vec{m}} \varphi_0 C_0
    \arrow[r, "\beta_{(1,\ldots, i'', \ldots, k), 0}"]
    &
    \varphi_{(1,\ldots, i'',\ldots, k)\bullet 0} \left(C_0, \ldots, C_i', \ldots, C_k\right) \otimes_{\vec{m}} C_0
    \arrow[r, "\mu"]
    &
    \varphi_{(1,\ldots, i'',\ldots, k)\bullet 0} C'
\end{tikzcd}
\]
Squares $\circled{1}$ and $\circled{2}$ both commute by Proposition \ref{prop:CoCnaturality}. Comparing the horizontal arrows, square $\circled{3}$ commutes up to a factor of $\beta_{(\mathbbm{1},\ldots, W_{v_i, j} ,\ldots,\mathbbm{1}), \mathbbm{1}}$, in the sense that 
\[
\beta_{(\mathbbm{1},\ldots, W_{v_i, j} ,\ldots,\mathbbm{1}), \mathbbm{1}} \times \text{(down, then right)} = \text{(right, then down)}.
\]
Notice that the $\beta_2$ and $\beta_3$ terms are both equal to 1, since all cobordisms involved are unweighted. Moreover, the $\beta_4$ term is equal to 1 since $\abs{1_{\vec{x}}\mathbbm{1}1_{y}} = (0,0)$ given any closures $\vec{x}$, $y$. Thus, the two sides differ by a value given by a single change of chronology
\[
\beta_{(\mathbbm{1},\ldots, W_{v_i, j} ,\ldots,\mathbbm{1}), \mathbbm{1}} = \iota\left(
\tikz[baseline={([yshift=-.5ex]current bounding box.center)}, scale=.75, yscale=.5,xscale=.75,yshift=-2.5cm]{
        \draw (1.25, 3) -- (1.25, 5.4);
	\draw 
		(.5,2)
		.. controls (.5,2.5) and (1.25,2.5) ..
		(1.25,3);
	\draw (.5,0) 
		--
		(.5,1)
		.. controls (.5,1.5) and (.5,1.5) ..
		(.5,2);
	\draw (2,0) 
		--
		(2,2)
		.. controls (2,2.5) and (1.25,2.5) ..
		(1.25,3);
        \node[fill=white,draw,rounded corners,scale=1] at (1.25, 4) {$W_{v_i,j}$};
}
\Rightarrow
\tikz[baseline={([yshift=-.5ex]current bounding box.center)}, scale=.75, yscale=.5,xscale=.75,yshift=-2.5cm]{
        \draw (1.25, 3) -- (1.25, 5.4);
	\draw 
		(.5,2)
		.. controls (.5,2.5) and (1.25,2.5) ..
		(1.25,3);
	\draw (.5,0) 
		--
		(.5,1)
		.. controls (.5,1.5) and (.5,1.5) ..
		(.5,2);
	\draw (2,0) 
		--
		(2,2)
		.. controls (2,2.5) and (1.25,2.5) ..
		(1.25,3);
        \node[fill=white,draw,rounded corners,scale=1] at (0.5, 1) {$W_{v_i,j}$};
}
\right)
\]
as in the definition of the compatibility maps $\beta$. Of course, we should view the $W_{v_i,j}$ on the left as $(\mathbbm{1}, \ldots, W_{v_i,j}, \ldots, \mathbbm{1}) \bullet \mathbbm{1}$, and the $W_{v_i,j}$ on the right as $(\mathbbm{1}, \ldots, W_{v_i,j}, \ldots, \mathbbm{1})$. On the other hand, square $\circled{4}$ commutes up to the value
\[
\iota\left(\tikz[baseline={([yshift=-.5ex]current bounding box.center)}, scale=.75, yscale=.5,xscale=.75,yshift=-2.5cm]{
        \draw (1.25, 3) -- (1.25, 5.4);
	\draw 
		(.5,2)
		.. controls (.5,2.5) and (1.25,2.5) ..
		(1.25,3);
	\draw (.5,0) 
		--
		(.5,1)
		.. controls (.5,1.5) and (.5,1.5) ..
		(.5,2);
	\draw (2,0) 
		--
		(2,2)
		.. controls (2,2.5) and (1.25,2.5) ..
		(1.25,3);
        \node[fill=white,draw,rounded corners,scale=1] at (0.5, 1) {$W_{v_i,j}$};
}
\Rightarrow
\tikz[baseline={([yshift=-.5ex]current bounding box.center)}, scale=.75, yscale=.5,xscale=.75,yshift=-2.5cm]{
        \draw (1.25, 3) -- (1.25, 5.4);
	\draw 
		(.5,2)
		.. controls (.5,2.5) and (1.25,2.5) ..
		(1.25,3);
	\draw (.5,0) 
		--
		(.5,1)
		.. controls (.5,1.5) and (.5,1.5) ..
		(.5,2);
	\draw (2,0) 
		--
		(2,2)
		.. controls (2,2.5) and (1.25,2.5) ..
		(1.25,3);
        \node[fill=white,draw,rounded corners,scale=1] at (1.25, 4) {$W_{v_i,j}$};
}
\right)
= 
(\beta_{(\mathbbm{1},\ldots, W_{v_i, j} ,\ldots,\mathbbm{1}), \mathbbm{1}})^{-1}
\]
in the sense that
\[
(\beta_{(\mathbbm{1},\ldots, W_{v_i, j} ,\ldots,\mathbbm{1}), \mathbbm{1}})^{-1}
\times
\text{(down, then right)}
=
\text{(right, then down)}.
\]
Thus, the former diagram commutes.

There are subtle differences in validating the commutativity of the latter diagram---in particular, the $\beta_4$ term in the analogue to square $\circled{3}$ is nontrivial. Anyway, the diagram in question factorizes as follows.
\[
\begin{tikzcd}[column sep = huge, row sep = huge, scale cd=1]
    \left(\varphi_1C_1, \ldots, \varphi_k C_k\right) \otimes_{\vec{m}} \varphi_0 C_0
    \MySymbAA{dr}
    \arrow[r, "\beta_{(1, \ldots, k), 0}"]
    \arrow[d, "(\mathbbm{1}\text{,} \ldots\text{,} \mathbbm{1})\otimes \varphi_H"]
    &
    \varphi_{(1,\ldots, k)\bullet 0} \left(C_0, \ldots, C_k\right) \otimes_{\vec{m}} C_0
    \MySymbBB{dr}
    \arrow[r, "\mu"]
    \arrow[d, "\varphi_{(\mathbbm{1}\text{,}\ldots\text{,}\mathbbm{1})\bullet H}"]
    &
    \varphi_{(1,\ldots, k)\bullet 0} C
    \arrow[d, "\varphi_{(\mathbbm{1}\text{,}\ldots\text{,}\mathbbm{1})\bullet H}"]
    \\
    \left(\varphi_1C_1, \ldots, \varphi_k C_k\right) \otimes_{\vec{m}} \varphi_{0'} C_0
    \MySymbCC{dr}
    \arrow[r, "\beta_{(1, \ldots, k), 0'}"]
    \arrow[d, "(\mathbbm{1} \text{,} \ldots\text{,} \mathbbm{1})\otimes \mathcal{F}_0"]
    &
    \varphi_{(1, \ldots, k)\bullet 0'} \left(C_0, \ldots, C_k\right) \otimes_{\vec{m}} C_0
    \MySymbDD{dr}
    \arrow[r, "\mu"]
    \arrow[d, "(\mathbbm{1} \text{,} \ldots\text{,} \mathbbm{1})\otimes \mathcal{F}_0"]
    &
    \varphi_{(1, \ldots, k)\bullet 0'} C
    \arrow[d, "\mathcal{F}_0'"]
    \\
    \left(\varphi_1C_1, \ldots, \varphi_k C_k\right) \otimes_{\vec{m}} \varphi_{0''} C_0'
    \arrow[r, "\beta_{(1, \ldots, k), 0''}"]
    &
    \varphi_{(1,\ldots, k)\bullet 0''} \left(C_0, \ldots, C_k\right) \otimes_{\vec{m}} C_0'
    \arrow[r, "\mu"]
    &
    \varphi_{(1, \ldots, k)\bullet 0''} C''
\end{tikzcd}
\]
Again, squares \circled{$1'$} and \circled{$2'$} commute thanks to Proposition \ref{prop:CoCnaturality}, and square \circled{$3'$} commutes up to a factor of $\beta_{\vec{\mathbbm{1}}, W_{v,j}}$ in the sense that
\[
\beta_{\vec{\mathbbm{1}}, W_{v,j}} \times \text{(down, then right)} = \text{(right, then down)}.
\]
Indeed, the $\beta_2$ and $\beta_3$ terms are trivial, but otherwise we have
\[
\beta_{\vec{\mathbbm{1}}, W_{v,j}} = \iota\left(
\tikz[baseline={([yshift=-.5ex]current bounding box.center)}, scale=.75, yscale=.5,xscale=.75,yshift=-2.5cm]{
        \draw (1.25, 3) -- (1.25, 5.4);
	\draw 
		(.5,2)
		.. controls (.5,2.5) and (1.25,2.5) ..
		(1.25,3);
	\draw (.5,0) 
		--
		(.5,1)
		.. controls (.5,1.5) and (.5,1.5) ..
		(.5,2);
	\draw (2,0) 
		--
		(2,2)
		.. controls (2,2.5) and (1.25,2.5) ..
		(1.25,3);
        \node[fill=white,draw,rounded corners,scale=1] at (1.25, 4) {$W_{v,j}$};
}
\Rightarrow
\tikz[baseline={([yshift=-.5ex]current bounding box.center)}, scale=.75, yscale=.5,xscale=.75,yshift=-2.5cm]{
        \draw (1.25, 3) -- (1.25, 5.4);
	\draw 
		(.5,2)
		.. controls (.5,2.5) and (1.25,2.5) ..
		(1.25,3);
	\draw (.5,0) 
		--
		(.5,1)
		.. controls (.5,1.5) and (.5,1.5) ..
		(.5,2);
	\draw (2,0) 
		--
		(2,2)
		.. controls (2,2.5) and (1.25,2.5) ..
		(1.25,3);
        \node[fill=white,draw,rounded corners,scale=1] at (2, 1) {$W_{v,j}$};
}
\right) \times \lambda\left(P', \abs{1_{\vec{y}} W_{v,j} 1_{z}}\right).
\]
Then again, the coherence isomorphisms specify that
\[
\tikz[baseline={([yshift=-.5ex]current bounding box.center)}, scale=.75, yscale=.5,xscale=.75,yshift=-2.5cm]{
        \draw (1.25, 3) -- (1.25, 5.4);
	\draw 
		(.5,2)
		.. controls (.5,2.5) and (1.25,2.5) ..
		(1.25,3);
	\draw (.5,1) 
		--
		(.5,1)
		.. controls (.5,1.5) and (.5,1.5) ..
		(.5,2);
	\draw (2,-1) 
		--
		(2,2)
		.. controls (2,2.5) and (1.25,2.5) ..
		(1.25,3);
        \node[fill=white,draw,rounded corners,scale=1] at (2, 0) {$W_{v,j}$};
}
= \beta_{\vec{\mathbbm{1}}, W_{v,j}} \times 
\tikz[baseline={([yshift=-.5ex]current bounding box.center)}, scale=.75, yscale=.5,xscale=.75,yshift=-2.5cm]{
        \draw (1.25, 3) -- (1.25, 5.4);
	\draw 
		(.5,2)
		.. controls (.5,2.5) and (1.25,2.5) ..
		(1.25,3);
	\draw (.5,1) 
		--
		(.5,1)
		.. controls (.5,1.5) and (.5,1.5) ..
		(.5,2);
	\draw (2,-1) 
		--
		(2,2)
		.. controls (2,2.5) and (1.25,2.5) ..
		(1.25,3);
        \node[fill=white,draw,rounded corners,scale=1] at (1.25, 4) {$W_{v,j}$};
}
\]
which is precisely
\[
\text{(down, then right)} = \beta_{\vec{\mathbbm{1}}, W_{v,j}} \times \text{(right, then down)}
\]
for square \circled{$4'$}. Thus the latter diagram commutes, concluding the proof.
\end{proof}

\subsection{dg-\texorpdfstring{$\mathscr{C}$}{Lg}-graded multimodules}
\label{ss:dgcgraded}
In \cite{naisse2020odd}, Naisse and Putrya provide a second notion of $\mathscr{C}$-graded dg-multimodules with differential which is $\mathscr{C}$-homogeneous rather than $\mathscr{C}$-grading preserving: they are distinguished from the former notion by calling them \textit{dg-$\mathscr{C}$-graded} multimodules. The only difference lies in the differential. This subsection is devoted to showing that the analogous objects exists in the multicategorically graded setting. However, along the way, we develop the notion of $\mathscr{C}$-commutative diagrams (see \S \ref{ss:C-commutativity}). While almost all succeeding work in this paper does not rely on anything proven in this section, we will use $\mathscr{C}$-commutative diagrams (especially Proposition \ref{prop:C-commutativity}) briefly in the discussion of duality (\S \ref{ss:duality}) and very minimally in the proof of properties of unified projectors (\S \ref{s:unifiedprops}). The author suggests proceeding to Section \ref{Chapter:tangleinvariant} and referring back to this section if necessary.

\begin{definition}
If $A_1,\ldots, A_k, B$ are $\mathscr{C}$-graded algebras, a \textit{dg-$\mathscr{C}$-graded $(A_1,\ldots, A_k; B)$-multimodule} $(M, d_M)$ is a $\mathbb{Z} \times \mathscr{C}$-graded $(A_1,\ldots, A_k; B)$-multimodule $M = \bigoplus_{n\in\mathbb{Z}, g\in\mathscr{C}} M_g^n$ together with a homogeneous differential $d_M$, written $\sum_{j}d_{M}^j$, satisfying $d_M \circ d_M = 0$.\footnote{Really, we could have written $d_M \circ_\mathscr{C} d_M = 0$, but clearly this is the case if and only if $d_M \circ d_M = 0$, since $\gamma$ does not take 0 for a value.} Again, we assume that the left and right action preserves homological grading. A map of dg-$\mathscr{C}$-graded multimodules $f: M \to N$ means a $\mathbb{K}$-linear map which preserves homological grading \textit{and $\mathscr{C}$-graded commutes with the differentials.}
\end{definition}

To be explicit, since $d_M = \sum_j d_M^j$ is homogeneous, we know that it is $\mathbb{K}$-linear, the sum is finite, and for all $m\in M$ it satisfies
\begin{enumerate}[label=(\roman*)]
    \item $d_M^j(M_g^n) \subset M_{\varphi_j(g)}^{n+1}$ if $g\in \mathsf{D}_j$ and $d_M^j(M_g^n)=0$ otherwise,
    \item $d_M^j(\rho_L(\vec{a}, m) = \beta_{\vec{e}, j}(\abs{\vec{a}}, \abs{m})^{-1} \rho_L(\vec{a}, d_M^j(m))$ for all $\vec{a}\in (A_1,\ldots, A_k)$, and
    \item $d_M^j(\rho_R(m, b)) = \beta_{j,e}(\abs{m}, \abs{b}) \rho_R(d_M^j(m), b)$ for all $b\in B$.
\end{enumerate}

Note that we do not require maps of dg-$\mathscr{C}$-graded multimodules to preserve $\mathscr{C}$-grading. The rest of this section is devoted to understanding what we mean by $\mathscr{C}$-graded commutativity; see \cite{naisse2020odd} for more details.

A \textit{commutativity system} on $\{\mathscr{I}, \Phi\}$ is a collection
\[
\mathtt{T} = \left\{\left((i,j),(i', j')\right) \in \left(\mathscr{I}^m\right)^2 \times \left(\mathscr{I}^m\right)^2 \right\}
\]
(that is, each of $i,j,i',j'$ may be $m$-vectors) such that
\begin{itemize}
    \item if $\left((i,j),(i', j')\right) \in \mathtt{T}$, then 
    \begin{enumerate}[label=(\alph*)]
        \item $\varphi_{j\circ i} = \varphi_{j' \circ i'}$, and
        \item $\left((i',j'),(i, j)\right) \in \mathtt{T}$
    \end{enumerate}
    and
    \item for any $k\ge 1$, if $\left((i_1,j_1),(i_1', j_1')\right), \ldots, \left((i_k,j_k),(i_k', j_k')\right), \left((i,j),(i', j')\right) \in \mathtt{T}$, then 
    \[
    \left(((i_1,\ldots, i_k)\bullet i,(j_1,\ldots, j_k)\bullet j),((i_1',\ldots, i_k')\bullet i',(j_1',\ldots, j_k')\bullet j')\right) \in \mathtt{T}.
    \]
\end{itemize}
We abbreviate the last requirement to $\left((\vec{i}, \vec{j}), (\vec{i}', \vec{j}')\right), \left((i,j), (i', j')\right)\in \mathtt{T} \implies (\vec{i}\bullet i, \vec{j} \bullet j)$.

For simplicity of exposition, assume $i,j,i',j'$ are single-entry. To witness the commutativity system, we introduce a collection of scalars
\[
\tau_{\substack{i, i' \\ j,j'}}^{\vec{X}\to Y}\in \mathbb{K}^\times
\]
for each $X_1,\ldots, X_k, Y\in \mathrm{Ob}(\mathscr{C})$ and $\left((i,j),(i',j')\right)\in \mathtt{T}$, satisfying
\begin{enumerate}[label = (\roman*)]
    \item $\tau_{\substack{i, i' \\ j,j'}}^{\vec{X}\to Y} = 1$ whenever $j\circ i = j' \circ i'$, and 
    \item $\left(\tau_{\substack{i, i' \\ j,j'}}^{\vec{X}\to Y} \right)^{-1} = \tau_{\substack{i', i \\ j',j}}^{\vec{X}\to Y}$ for each $\left((i,j),(i',j')\right)\in \mathtt{T}$.
\end{enumerate}
If $\left((i,j), (i', j')\right)\not\in \mathtt{T}$, then we declare $\tau_{\substack{i, i' \\ j, j'}}$ to be zero. We will write $\tau_{\substack{\vec{i}, \vec{i}' \\ \vec{j}, \vec{j}'}}^{\vec{X} \to \vec{Y}}$ for the scalar witness when $m\not=1$---the above definition of $\tau$ extends to the case where $i,j, i', j'$ are vectors, requiring $\tau_{\substack{\vec{i}, \vec{i}' \\ \vec{j}, \vec{j}'}}^{\vec{X} \to \vec{Y}} = 1$ whenever $\vec{j} \circ \vec{i} = \vec{j}' \circ \vec{i}'$, interpreted correctly. As earlier, we write $\tau_{\substack{i, i' \\ j, j'}}(g)$ to mean $\tau_{\substack{i, i' \\ j, j'}}^{\vec{X} \to Y}$ whenever $g: \vec{X} \to Y$. Finally, we say that a commutativity system $\mathtt{T}$ is \textit{compatible} with a shifting 2-system through $\tau$ if two equations are satisfied. The first is
\begin{equation}
\label{eq:CommSystComp}
\tau_{\substack{\vec{i}_1 \bullet i_1, \vec{i}_2 \bullet i_2 \\ \vec{j}_1 \bullet j_1, \vec{j}_2 \bullet j_2}} (g'g) \Xi_{\substack{i_1, \vec{i}_1 \\ j_1, \vec{j}_1}}(g'g) \beta_{\vec{j}_1 \circ \vec{i}_1, j_1 \circ i_1}(g',g) 
=
\Xi_{\substack{i_2, \vec{i}_2 \\ j_2, \vec{j}_2}}(g'g) \beta_{\vec{j}_2 \circ \vec{i}_2, j_2 \circ i_2}(g', g) \tau_{\substack{\vec{i}_1, \vec{i}_2 \\ \vec{j}_1, \vec{j}_2}}(g') \tau_{\substack{i_1, i_2 \\ j_1, j_2}}(g)
\end{equation}
which translates to the following diagram.
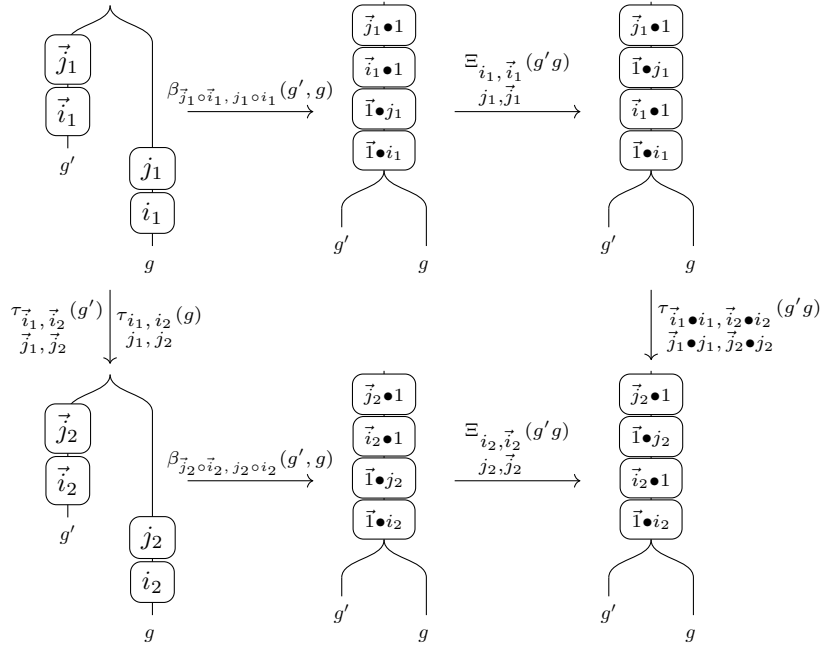
\[
\begin{tikzcd}[column sep=huge, row sep=large, scale=1.5]
\tikz[yscale=.5,xscale=.75]{
	\draw 
		(.5,2)
		.. controls (.5,2.5) and (1.25,2.5) ..
		(1.25,3);
	\draw (.5,-.8) node[below,scale=.75]{$g'$}
		-- 
		(.5,-.8)
		--
		(.5,1)
		.. controls (.5,1.5) and (.5,1.5) ..
		(.5,2);
	\draw (2,-3.4) node[below,scale=.75]{$g\phantom{'}$}
		--
		(2, -2.5)
		--
		(2,2)
		.. controls (2,2.5) and (1.25,2.5) ..
		(1.25,3);
        \node[fill=white,draw,rounded corners,scale=1] at (.5,1.3) {$\vec{j}_1$};
        \node[fill=white,draw,rounded corners,scale=1] at (.5,-0.1) {$\vec{i}_1$};
        \node[fill=white,draw,rounded corners,scale=1] at (2, -1.5) {$j_1$};
        \node[fill=white,draw,rounded corners,scale=1] at (2, -2.7) {$i_1$};
}
\arrow[r, "\beta_{\vec{j}_1 \circ \vec{i}_1\text{,}\,j_1 \circ i_1}(g'\text{,}\,g)"] \arrow[d, "\tau_{\substack{\vec{i}_1\text{,}\, \vec{i}_2 \\ \vec{j}_1\text{,}\, \vec{j}_2}}(g')"', "\tau_{\substack{i_1\text{,}\, i_2 \\ j_1\text{,}\, j_2}}(g)"]
&
\tikz[yscale=.5,xscale=.75,yshift=-4.4cm]{
        \draw (1.25, 3) -- (1.25, 7.5);
	\draw 
		(.5,2)
		.. controls (.5,2.5) and (1.25,2.5) ..
		(1.25,3);
	\draw (.5,1.5) node[below,scale=.75]{$g'$}
		.. controls (.5,1.5) and (.5,1.5) ..
		(.5,2);
	\draw (2,1) node[below,scale=.75]{$g\phantom{'}$}
		--
		(2,2)
		.. controls (2,2.5) and (1.25,2.5) ..
		(1.25,3);
         \node[fill=white,draw,rounded corners,scale=1] at (1.25, 3.35) {${\scriptstyle \vec{1}\bullet i_1}$};
        \node[fill=white,draw,rounded corners,scale=1] at (1.25, 4.4667) {${\scriptstyle \vec{1}\bullet j_1}$};
        \node[fill=white,draw,rounded corners,scale=1] at (1.25, 5.56) {${\scriptstyle \vec{i}_1 \bullet 1}$};
        \node[fill=white,draw,rounded corners,scale=1] at (1.25, 6.69) {${\scriptstyle \vec{j}_1 \bullet 1}$};
}
\arrow[r, "\Xi_{\substack{i_1\text{,}\, \vec{i}_1 \\ j_1, \vec{j}_1}}(g'g)"]
&
\tikz[yscale=.5,xscale=.75,yshift=-4.4cm]{
        \draw (1.25, 3) -- (1.25, 7.5);
	\draw 
		(.5,2)
		.. controls (.5,2.5) and (1.25,2.5) ..
		(1.25,3);
	\draw (.5,1.5) node[below,scale=.75]{$g'$}
		.. controls (.5,1.5) and (.5,1.5) ..
		(.5,2);
	\draw (2,1) node[below,scale=.75]{$g\phantom{'}$}
		--
		(2,2)
		.. controls (2,2.5) and (1.25,2.5) ..
		(1.25,3);
         \node[fill=white,draw,rounded corners,scale=1] at (1.25, 3.35) {${\scriptstyle \vec{1}\bullet i_1}$};
        \node[fill=white,draw,rounded corners,scale=1] at (1.25, 5.57) {${\scriptstyle \vec{1}\bullet j_1}$};
        \node[fill=white,draw,rounded corners,scale=1] at (1.25, 4.447) {${\scriptstyle \vec{i}_1 \bullet 1}$};
        \node[fill=white,draw,rounded corners,scale=1] at (1.25, 6.69) {${\scriptstyle \vec{j}_1 \bullet 1}$};
}
\arrow[d, "\tau_{\substack{\vec{i}_1 \bullet i_1\text{,}\, \vec{i}_2 \bullet i_2 \\ \vec{j}_1 \bullet j_1\text{,}\, \vec{j}_2 \bullet j_2}} (g'g) "]
\\
\tikz[yscale=.5,xscale=.75]{
	\draw 
		(.5,2)
		.. controls (.5,2.5) and (1.25,2.5) ..
		(1.25,3);
	\draw (.5,-.8) node[below,scale=.75]{$g'$}
		-- 
		(.5,-.8)
		--
		(.5,1)
		.. controls (.5,1.5) and (.5,1.5) ..
		(.5,2);
	\draw (2,-3.4) node[below,scale=.75]{$g\phantom{'}$}
		--
		(2, -2.5)
		--
		(2,2)
		.. controls (2,2.5) and (1.25,2.5) ..
		(1.25,3);
        \node[fill=white,draw,rounded corners,scale=1] at (.5,1.3) {$\vec{j}_2$};
        \node[fill=white,draw,rounded corners,scale=1] at (.5,-.1) {$\vec{i}_2$};
        \node[fill=white,draw,rounded corners,scale=1] at (2, -1.5) {$j_2$};
        \node[fill=white,draw,rounded corners,scale=1] at (2, -2.7) {$i_2$};
}
\arrow[r, "\beta_{\vec{j}_2 \circ \vec{i}_2\text{,}\, j_2 \circ i_2}(g'\text{,}\, g)"]
&
\tikz[yscale=.5,xscale=.75,yshift=-4.4cm]{
        \draw (1.25, 3) -- (1.25, 7.5);
	\draw 
		(.5,2)
		.. controls (.5,2.5) and (1.25,2.5) ..
		(1.25,3);
	\draw (.5,1.5) node[below,scale=.75]{$g'$}
		.. controls (.5,1.5) and (.5,1.5) ..
		(.5,2);
	\draw (2,1) node[below,scale=.75]{$g\phantom{'}$}
		--
		(2,2)
		.. controls (2,2.5) and (1.25,2.5) ..
		(1.25,3);
         \node[fill=white,draw,rounded corners,scale=1] at (1.25, 3.35) {${\scriptstyle \vec{1}\bullet i_2}$};
        \node[fill=white,draw,rounded corners,scale=1] at (1.25, 4.4667) {${\scriptstyle \vec{1}\bullet j_2}$};
        \node[fill=white,draw,rounded corners,scale=1] at (1.25, 5.56) {${\scriptstyle \vec{i}_2 \bullet 1}$};
        \node[fill=white,draw,rounded corners,scale=1] at (1.25, 6.69) {${\scriptstyle \vec{j}_2 \bullet 1}$};
}
\arrow[r, "\Xi_{\substack{i_2, \vec{i}_2 \\ j_2, \vec{j}_2}}(g'g)"]
&
\tikz[yscale=.5,xscale=.75,yshift=-4.4cm]{
        \draw (1.25, 3) -- (1.25, 7.5);
	\draw 
		(.5,2)
		.. controls (.5,2.5) and (1.25,2.5) ..
		(1.25,3);
	\draw (.5,1.5) node[below,scale=.75]{$g'$}
		.. controls (.5,1.5) and (.5,1.5) ..
		(.5,2);
	\draw (2,1) node[below,scale=.75]{$g\phantom{'}$}
		--
		(2,2)
		.. controls (2,2.5) and (1.25,2.5) ..
		(1.25,3);
         \node[fill=white,draw,rounded corners,scale=1] at (1.25, 3.35) {${\scriptstyle \vec{1}\bullet i_2}$};
        \node[fill=white,draw,rounded corners,scale=1] at (1.25, 5.57) {${\scriptstyle \vec{1}\bullet j_2}$};
        \node[fill=white,draw,rounded corners,scale=1] at (1.25, 4.447) {${\scriptstyle \vec{i}_2 \bullet 1}$};
        \node[fill=white,draw,rounded corners,scale=1] at (1.25, 6.69) {${\scriptstyle \vec{j}_2 \bullet 1}$};
}
\end{tikzcd}
\]
The second equation establishes consistency between $\tau$ and $\Xi$: we require that
\begin{equation}
\label{eq:TauXiComp}
\tau_{\substack{\vec{\mathrm{Id}}\bullet i, \vec{j} \bullet \mathrm{Id} \\ \vec{j} \bullet \mathrm{Id}, \vec{\mathrm{Id}} \bullet i}} = \Xi_{\substack{\mathrm{Id}, \vec{j} \\ i, \vec{\mathrm{Id}}}} \Xi_{\substack{i, \vec{\mathrm{Id}} \\ \mathrm{Id}, \vec{j}}}^{-1}.
\end{equation}
In particular, notice that in order to conclude that $\left((\vec{\mathrm{Id}}\bullet i, \vec{j} \bullet \mathrm{Id}), (\vec{j} \bullet \mathrm{Id}, \vec{\mathrm{Id}} \bullet i)\right) \in \mathtt{T}$, where any $\mathrm{Id}$ may be replaced by any element of $\mathscr{I}_{\mathrm{Id}}$, it is sufficient if $((i', \mathrm{Id}), (\mathrm{Id}, i')) \in \mathtt{T}$ for any $i'\in \mathscr{I}$---this will clearly be the case in the $\mathscr{G}$-graded setting. This equation translates to the following diagram.
\[
\begin{tikzcd}
\left(\vec{j} \bullet \mathrm{Id} \right) \circ \left(\vec{\mathrm{Id}} \bullet i \right) 
\arrow[rrr, "\tau_{\substack{\vec{\mathrm{Id}}\bullet i\text{,}\, \vec{j} \bullet \mathrm{Id} \\ \vec{j} \bullet \mathrm{Id}\text{,}\, \vec{\mathrm{Id}} \bullet i}}"]
\arrow[dr, "\Xi_{\substack{i\text{,}\, \vec{\mathrm{Id}} \\ \mathrm{Id}\text{,}\, \vec{j}}}^{-1}"']
&&&
\left(\vec{\mathrm{Id}} \bullet i\right) \circ \left(\vec{j} \bullet \mathrm{Id} \right)
\\
&
\left(\vec{j} \circ \vec{\mathrm{Id}}\right) \bullet \left(\mathrm{Id} \circ i \right)
\arrow[r, equals]
&
\left(\vec{\mathrm{Id}} \circ \vec{j}\right) \bullet \left(i \circ \mathrm{Id} \right)
\arrow[ur, "\Xi_{\substack{\mathrm{Id}\text{,}\, \vec{j} \\ i\text{,}\, \vec{\mathrm{Id}}}}"']
&
\end{tikzcd}
\]

\subsubsection{$\mathscr{G}$-graded commutativity} 
As before, one last time, we will describe the $\mathscr{G}$-graded setting before passing on to generalities. We will not consider dg-$\mathscr{G}$-graded multimodules explicitly, but we can construct them using the information of this section. See \cite{naisse2020odd} for more generalities of these objects. We'll write $\Delta^v$ for $(\Delta, v) \in \mathscr{I}$ to reduce the number of nested ordered pairs. We'll describe the non-vectorized setting first. Let $\mathtt{T}$ denote the collection of all pairs $\left\{\left((\Delta_1^{v_1}, \Delta_2^{v_2}), (\Delta_1'^{v_1'}, \Delta_2'^{v_2'})\right)\right\}$ for which
\begin{itemize}
    \item there exists a locally vertical change of chronology $H: \Delta_2 \circ \Delta_1 \Rightarrow \Delta_2' \circ \Delta_1'$, and
    \item $v_1 = v_2'$ and $v_2 = v_1'$.
\end{itemize}
Similarly, in the vecotrized setting, $\left((\vec{\Delta}_1^{\vec{v}_1}, \vec{\Delta}_2^{\vec{v}_2}), (\vec{\Delta}_1'^{\vec{v}_1'}, \vec{\Delta}_2'^{\vec{v}_2'})\right)$ is in $\mathtt{T}$ if there are locally vertical changes of chronology $H_\ell: \Delta_{2,\ell}\circ \Delta_{1,\ell} \Rightarrow \Delta_{2,\ell}' \circ \Delta_{1, \ell}'$ for all $\ell$ and $\vec{v}_1 = \vec{v}_2'$ and $\vec{v}_2 = \vec{v}_1'$. Notice that $\mathtt{T}$ satisfies the criteria of a commutativity system since cobordisms which differ only with respect to a locally vertical change of chronology induce the same $\mathscr{G}$-grading shift, locally vertical changes of chronology are invertible, and locally vertical changes of chronology are well behaved with respect to horizontal composition of cobordisms.

Next, we set
\[
\tau_{\substack{(\Delta_1, v_1), (\Delta_1', v_1') \\ (\Delta_2, v_2), (\Delta_2', v_2')}}^{\vec{x} \to y} = \iota({}_{\vec{x}}H_y) \lambda(v_2, v_1)
\]
again, where $H$ is the locally vertical change of chronology $H: \Delta_2 \circ \Delta_1 \Rightarrow \Delta_2' \circ \Delta_1'$. 
In the vectorized setting, we set
\[
\tau_{\substack{(\vec{\Delta}_1, \vec{v}_1), (\vec{\Delta}_1', \vec{v}_1') \\ (\vec{\Delta}_2, \vec{v}_2), (\vec{\Delta}_2', \vec{v}_2')}}^{\vec{x} \to \vec{y}} = \prod_\ell\iota({}_{\vec{x}_\ell}(H_\ell)_{y_\ell}) \lambda(V_2, V_1)
\]
where $V_2$, $V_1$ denote the sums of the entries of $\vec{v}_2$ and $\vec{v}_1$ respectively. We will write $\tau = \tau_1 \tau_2$, for $\tau_1$ the part coming from the change of chronology and $\tau_2$ the other. 

Notice that if $(\Delta_2, v_2) \circ (\Delta_1, v_1) = (\Delta_2', v_2') \circ (\Delta_1,' v_1')$, then $H$ is the identical change of chronology, and $v_1 = v_2' = v_2 = v_1'$ so $\lambda(v_2, v_1) = 1$, hence $\tau_{\substack{(\Delta_1, v_1), (\Delta_1', v_1') \\ (\Delta_2, v_2), (\Delta_2', v_2')}}^{\vec{x} \to y} = 1$. Also, if $\overline{H}: \Delta_2' \circ \Delta_1' \Rightarrow \Delta_2 \circ \Delta_1$ is also a locally vertical change of chronology (guaranteed to exist by the existence of $H$) then
\[
\tau_{\substack{(\Delta_1', v_1'), (\Delta_1, v_1) \\ (\Delta_2', v_2'), (\Delta_2, v_2)}}^{\vec{x} \to y} = \iota({}_{\vec{x}}\overline{H}_y) \lambda(v_2', v_1') = \iota({}_{\vec{x}}H_y)^{-1} \lambda(v_1, v_2) = \left(\tau_{\substack{(\Delta_1, v_1), (\Delta_1', v_1') \\ (\Delta_2, v_2), (\Delta_2', v_2')}}^{\vec{x} \to y}\right)^{-1}
\]
as desired.

\begin{proposition}
This commutativity system $\mathtt{T}$ is compatible with the $\mathscr{G}$-grading shifting 2-system defined previously, through the scalars $\tau$.
\end{proposition}

\begin{proof}
The validity of (\ref{eq:TauXiComp}) is simple: recall that $\mathscr{I}_{\mathrm{Id}}$ consists of elements $(\mathbbm{1}_{D^\wedge}, (0,0))$ for any planar arc diagram $D$. Thus
\[
\tau_{\substack{\mathrm{Id} \bullet (\Delta, v), (\vec{\Delta}, \vec{v} \bullet \mathrm{Id}) \\ (\vec{\Delta}, \vec{v}) \bullet \mathrm{Id}, \vec{\mathrm{Id}} \bullet (\Delta, v)}}^{\vec{x} \to y} = \iota({}_{\vec{x}}H_y) \lambda(V, v) 
\]
where $V$ is the sum of entries of $\vec{v}$ and $H: ((\vec{\Delta}, \vec{v})\bullet \mathrm{Id}) \circ (\vec{\mathrm{Id}} \bullet (\Delta, v)) \Rightarrow (\vec{\mathrm{Id}} \bullet (\Delta, v)) \circ ((\vec{\Delta}, \vec{v}) \bullet \mathrm{Id})$. On the other hand, 
\[
\Xi_{\substack{\mathrm{Id}, (\vec{\Delta}, \vec{v}) \\ (\Delta, v), \vec{\mathrm{Id}}}} \Xi_{\substack{(\Delta, v), \vec{\mathrm{Id}} \\ \mathrm{Id}, (\vec{\Delta}, \vec{v})}}^{-1} = \left(\iota({}_{\vec{x}}H''_y) \lambda(V,v)\right) \cdot \left(\iota({}_{\vec{x}} H'_y) \lambda((0,0)(0,0))^{-1} \right)
\]
where
\[
((\vec{\Delta}, \vec{v})\bullet \mathrm{Id}) \circ (\vec{\mathrm{Id}} \bullet (\Delta, v))
\xRightarrow{H'}
((\vec{\Delta}, \vec{v})\circ \vec{\mathrm{Id}}) \bullet (\mathrm{Id} \circ (\Delta, v))
=
(\vec{\mathrm{Id}} \circ (\vec{\Delta}, \vec{v})) \bullet ((\Delta, v) \circ \mathrm{Id})
\xRightarrow{H''}
(\vec{\mathrm{Id}} \bullet (\Delta, v)) \circ ((\vec{\Delta}, \vec{v}) \bullet \mathrm{Id}).
\]
Since $H$ and $H'' \circ H'$ are locally vertical changes of chronology with the source and target, Proposition \ref{PutyraHammer} implies that $\iota({}_{\vec{x}}H_y) = \iota({}_{\vec{x}}(H'' \circ H')_y) = \iota({}_{\vec{x}}H''_y) \iota({}_{\vec{x}}H'_y)$, so equation (\ref{eq:TauXiComp}) is satisfied.

To check equation (\ref{eq:CommSystComp}), we apply familiar arguments. Actually, the computation is fairly simple compared to the previous proofs of this type. On one hand, ignoring $\mathbb{Z} \times \mathbb{Z}$-degree to start, consider the diagram
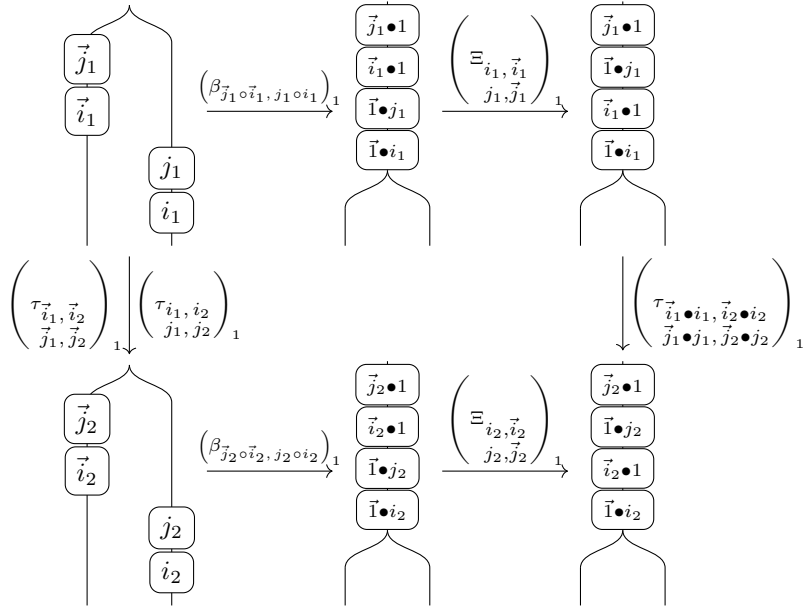
\[
\begin{tikzcd}[column sep=huge, row sep=huge, scale=1.5]
\tikz[yscale=.5,xscale=.75]{
	\draw 
		(.5,2)
		.. controls (.5,2.5) and (1.25,2.5) ..
		(1.25,3);
	\draw (.5,-3.4)
		-- 
		(.5,-.8)
		--
		(.5,1)
		.. controls (.5,1.5) and (.5,1.5) ..
		(.5,2);
	\draw (2,-3.4)
		--
		(2, -2.5)
		--
		(2,2)
		.. controls (2,2.5) and (1.25,2.5) ..
		(1.25,3);
        \node[fill=white,draw,rounded corners,scale=1] at (.5,1.3) {$\vec{j}_1$};
        \node[fill=white,draw,rounded corners,scale=1] at (.5,-0.1) {$\vec{i}_1$};
        \node[fill=white,draw,rounded corners,scale=1] at (2, -1.5) {$j_1$};
        \node[fill=white,draw,rounded corners,scale=1] at (2, -2.7) {$i_1$};
}
\arrow[r, "\left(\beta_{\vec{j}_1 \circ \vec{i}_1\text{,}\,j_1 \circ i_1}\right)_1"] \arrow[d, "\left(\tau_{\substack{\vec{i}_1\text{,}\, \vec{i}_2 \\ \vec{j}_1\text{,}\, \vec{j}_2}}\right)_1"', "\left(\tau_{\substack{i_1\text{,}\, i_2 \\ j_1\text{,}\, j_2}}\right)_1"]
&
\tikz[yscale=.5,xscale=.75,yshift=-4.4cm]{
        \draw (1.25, 3) -- (1.25, 7.5);
	\draw 
		(.5,2)
		.. controls (.5,2.5) and (1.25,2.5) ..
		(1.25,3);
	\draw (.5,1)
		.. controls (.5,1.5) and (.5,1.5) ..
		(.5,2);
	\draw (2,1)
		--
		(2,2)
		.. controls (2,2.5) and (1.25,2.5) ..
		(1.25,3);
         \node[fill=white,draw,rounded corners,scale=1] at (1.25, 3.35) {${\scriptstyle \vec{1}\bullet i_1}$};
        \node[fill=white,draw,rounded corners,scale=1] at (1.25, 4.4667) {${\scriptstyle \vec{1}\bullet j_1}$};
        \node[fill=white,draw,rounded corners,scale=1] at (1.25, 5.56) {${\scriptstyle \vec{i}_1 \bullet 1}$};
        \node[fill=white,draw,rounded corners,scale=1] at (1.25, 6.69) {${\scriptstyle \vec{j}_1 \bullet 1}$};
}
\arrow[r, "\left(\Xi_{\substack{i_1\text{,}\, \vec{i}_1 \\ j_1, \vec{j}_1}}\right)_1"]
&
\tikz[yscale=.5,xscale=.75,yshift=-4.4cm]{
        \draw (1.25, 3) -- (1.25, 7.5);
	\draw 
		(.5,2)
		.. controls (.5,2.5) and (1.25,2.5) ..
		(1.25,3);
	\draw (.5,1) 
		.. controls (.5,1.5) and (.5,1.5) ..
		(.5,2);
	\draw (2,1) 
		--
		(2,2)
		.. controls (2,2.5) and (1.25,2.5) ..
		(1.25,3);
         \node[fill=white,draw,rounded corners,scale=1] at (1.25, 3.35) {${\scriptstyle \vec{1}\bullet i_1}$};
        \node[fill=white,draw,rounded corners,scale=1] at (1.25, 5.57) {${\scriptstyle \vec{1}\bullet j_1}$};
        \node[fill=white,draw,rounded corners,scale=1] at (1.25, 4.447) {${\scriptstyle \vec{i}_1 \bullet 1}$};
        \node[fill=white,draw,rounded corners,scale=1] at (1.25, 6.69) {${\scriptstyle \vec{j}_1 \bullet 1}$};
}
\arrow[d, "\left(\tau_{\substack{\vec{i}_1 \bullet i_1\text{,}\, \vec{i}_2 \bullet i_2 \\ \vec{j}_1 \bullet j_1\text{,}\, \vec{j}_2 \bullet j_2}}\right)_1"]
\\
\tikz[yscale=.5,xscale=.75]{
	\draw 
		(.5,2)
		.. controls (.5,2.5) and (1.25,2.5) ..
		(1.25,3);
	\draw (.5,-3.4)
		-- 
		(.5,-.8)
		--
		(.5,1)
		.. controls (.5,1.5) and (.5,1.5) ..
		(.5,2);
	\draw (2,-3.4)
		--
		(2, -2.5)
		--
		(2,2)
		.. controls (2,2.5) and (1.25,2.5) ..
		(1.25,3);
        \node[fill=white,draw,rounded corners,scale=1] at (.5,1.3) {$\vec{j}_2$};
        \node[fill=white,draw,rounded corners,scale=1] at (.5,-.1) {$\vec{i}_2$};
        \node[fill=white,draw,rounded corners,scale=1] at (2, -1.5) {$j_2$};
        \node[fill=white,draw,rounded corners,scale=1] at (2, -2.7) {$i_2$};
}
\arrow[r, "\left(\beta_{\vec{j}_2 \circ \vec{i}_2\text{,}\, j_2 \circ i_2}\right)_1"]
&
\tikz[yscale=.5,xscale=.75,yshift=-4.4cm]{
        \draw (1.25, 3) -- (1.25, 7.5);
	\draw 
		(.5,2)
		.. controls (.5,2.5) and (1.25,2.5) ..
		(1.25,3);
	\draw (.5,1) 
		.. controls (.5,1.5) and (.5,1.5) ..
		(.5,2);
	\draw (2,1)
		--
		(2,2)
		.. controls (2,2.5) and (1.25,2.5) ..
		(1.25,3);
         \node[fill=white,draw,rounded corners,scale=1] at (1.25, 3.35) {${\scriptstyle \vec{1}\bullet i_2}$};
        \node[fill=white,draw,rounded corners,scale=1] at (1.25, 4.4667) {${\scriptstyle \vec{1}\bullet j_2}$};
        \node[fill=white,draw,rounded corners,scale=1] at (1.25, 5.56) {${\scriptstyle \vec{i}_2 \bullet 1}$};
        \node[fill=white,draw,rounded corners,scale=1] at (1.25, 6.69) {${\scriptstyle \vec{j}_2 \bullet 1}$};
}
\arrow[r, "\left(\Xi_{\substack{i_2, \vec{i}_2 \\ j_2, \vec{j}_2}}\right)_1"]
&
\tikz[yscale=.5,xscale=.75,yshift=-4.4cm]{
        \draw (1.25, 3) -- (1.25, 7.5);
	\draw 
		(.5,2)
		.. controls (.5,2.5) and (1.25,2.5) ..
		(1.25,3);
	\draw (.5,1)
		.. controls (.5,1.5) and (.5,1.5) ..
		(.5,2);
	\draw (2,1)
		--
		(2,2)
		.. controls (2,2.5) and (1.25,2.5) ..
		(1.25,3);
         \node[fill=white,draw,rounded corners,scale=1] at (1.25, 3.35) {${\scriptstyle \vec{1}\bullet i_2}$};
        \node[fill=white,draw,rounded corners,scale=1] at (1.25, 5.57) {${\scriptstyle \vec{1}\bullet j_2}$};
        \node[fill=white,draw,rounded corners,scale=1] at (1.25, 4.447) {${\scriptstyle \vec{i}_2 \bullet 1}$};
        \node[fill=white,draw,rounded corners,scale=1] at (1.25, 6.69) {${\scriptstyle \vec{j}_2 \bullet 1}$};
}
\end{tikzcd}
\]
where $i_1 = \Delta_1$,  $j_1 = \Delta_2$, $i_2 = \Delta_1'$, $j_2 = \Delta_2'$, and so on. The two paths trace out changes of chronology with the same source and target, so we conclude that the contributions of $\tau_1$, $\Xi_1$, and $\beta_1$ from equation (\ref{eq:CommSystComp}) agree on the nose.

On the other hand, since $\vec{\Delta}_1 \circ \vec{\Delta}_2$ and $\vec{\Delta}_1' \circ \vec{\Delta}_2'$, as well as $\Delta_2\circ \Delta_1$ and $(\Delta_2' \circ \Delta_1')$, differ only by a locally vertical changes of chronology, plus $v_1 + v_2 = v_2' + v_1'$ and $V_1 + V_2 = V_2' + V_1'$, it is easy to find that
\[
\left(\beta_{(\vec{\Delta}_2 \circ \vec{\Delta}_1, \vec{v}_1 + \vec{v}_2), (\Delta_2 \circ \Delta_1, v_1 + v_2)}\right)_{2,3,4}
=
\left(\beta_{(\vec{\Delta}_2' \circ \vec{\Delta}_1', \vec{v}_1' + \vec{v}_2'), (\Delta_2' \circ \Delta_1', v_1' + v_2')}\right)_{2,3,4}.
\]
If these conditions were not true, then the $\tau$ maps involved would be zero, and equation (\ref{eq:CommSystComp}) would hold trivially. Moreover, we compute
\[
\left(\tau_{\substack{(\vec{\Delta}_1 \bullet \Delta_1, V_1 + v_1), (\vec{\Delta}_1' \bullet \Delta_1', V_1' + v_1') \\ (\vec{\Delta}_2 \bullet \Delta_2, V_2 + v_2), (\vec{\Delta}_2' \bullet \Delta_2', V_2' + v_2')}}\right)_2 = \lambda(V_2 + v_2, V_1 + v_1) = \lambda(V_2, V_1) \cdot \lambda(V_2, v_1) \cdot \lambda(v_2, V_1) \cdot \lambda(v_2, v_1)
\]
and
\[
\left(\Xi_{\substack{(\Delta_1, v_1), (\vec{\Delta}_1, \vec{v}_1) \\ (\Delta_2, v_2), (\vec{\Delta}_1, \vec{v}_2)}}\right)_2 = \lambda(V_1, v_2)
\]
on one side, and 
\[
\left(\tau_{\substack{(\vec{\Delta}_1, \vec{v}_1), (\vec{\Delta}_1', \vec{v}_1') \\ (\vec{\Delta}_2, \vec{v}_2), (\vec{\Delta}_2', \vec{v}_2')}}\right)_2 = \lambda(V_2, V_1),
\]
\[
\left(\tau_{\substack{(\Delta_1, v_1), (\Delta_1', v_1') \\ (\Delta_2, v_2), (\Delta_2', v_2')}}\right)_2 = \lambda(v_2, v_1),
\]
and
\[
\left(\Xi_{\substack{(\Delta_1', v_1'), (\vec{\Delta}_1', \vec{v}_1') \\ (\Delta_2', v_2'), (\vec{\Delta}_1', \vec{v}_2')}}\right)_2 = \lambda(V_1', v_2') = \lambda(V_2, v_1)
\]
on the other. Since $\lambda(V_1, v_2) = \lambda(v_2, V_1)^{-1}$, these computations tell us that the contributions of $\tau_2$, $\Xi_2$, and $\beta_{2,3,4}$ from equation (\ref{eq:CommSystComp}) also agree on the nose, concluding the proof.
\end{proof}

There may be other choices of commutativity systems compatible with the $\mathscr{G}$-grading shifting 2-system. However, this doesn't matter so much: the existence of a commutativity system is more important than the commutativity system itself.

\subsubsection{Generalities of commutativity systems}
\label{ss:C-commutativity}

As before, we obtain natural transformations $\varphi_{j\circ i} \Rightarrow \varphi_{j' \circ i'}$
\begin{align*}
    \varphi_{j\circ i} (M) & \to \varphi_{j' \circ i'}(M) \\
    m & \mapsto \tau_{\substack{i, i' \\ j, j'}}(\abs{m}) m
\end{align*}
or, more generally, $\varphi_{\vec{j} \circ \vec{i}} \Rightarrow \varphi_{\vec{j}' \circ \vec{i}'}$ given by
\begin{align*}
\varphi_{\vec{j} \circ \vec{i}}(M_1,\ldots, M_k) &\to \varphi_{\vec{j}' \circ \vec{i}'}(M_1,\ldots, M_k) \\
\vec{m} &\mapsto \tau_{\substack{\vec{i}, \vec{i}' \\ \vec{j}, \vec{j}'}}(\abs{\vec{m}}) \vec{m}
\end{align*}
Then, the compatibility equations (\ref{eq:CommSystComp}) and (\ref{eq:TauXiComp}) imply the following commutative diagrams in categories of $\mathscr{G}$-graded multimodules.
\[
\begin{tikzcd}[column sep = huge, row sep = huge, scale cd=.795]
\varphi_{\vec{j}_1\circ \vec{i}_1} (M_1,\ldots, M_k) \otimes \varphi_{j_1 \circ i_1} (M) \arrow[d, "\tau_{\substack{\vec{i}_1, \vec{i}_2 \\ \vec{j}_1, \vec{j}_2}} \otimes \tau_{\substack{i_1, i_2 \\ j_1, j_2}}"'] \arrow[r, "\beta_{\vec{j}_1 \circ \vec{i}_1, j_1\circ i_1}"]
&
\varphi_{(\vec{j}_1\circ \vec{i}_1) \bullet (j_1 \circ i_1)}\left((M_1, \ldots, M_k) \otimes M\right) \arrow[r, "\Xi_{\substack{i_1, \vec{i}_1 \\ j_1, \vec{j}_1}}"]
&
\varphi_{(\vec{j}_1 \bullet j_1) \circ (\vec{i}_1 \bullet i_1)} \left((M_1,\ldots, M_k) \otimes M \right) \arrow[d, "\tau_{\substack{\vec{i}_1 \bullet i_1, \vec{i}_2 \bullet i_2 \\ \vec{j}_1\bullet j_1, \vec{j}_2 \bullet j_2}}"]
\\
\varphi_{\vec{j}_2 \circ \vec{i}_2}(M_1,\ldots, M_k) \otimes \varphi_{j_2 \circ i_2}(M) \arrow[r, "\beta_{\vec{j}_2\circ \vec{i}_2, j_2\circ i_2}"]
&
\varphi_{(\vec{j}_2 \circ \vec{i}_2) \bullet (j_2 \circ i_2)}\left((M_1,\ldots, M_k) \otimes M\right) \arrow[r, "\Xi_{\substack{i_2, \vec{i}_2 \\ j_2, \vec{j}_2}}"]
&
\varphi_{(\vec{j}_2 \bullet j_2) \circ (\vec{i}_2 \bullet i_2)}\left((M_1,\ldots, M_k) \otimes M\right)
\end{tikzcd}
\]
\[
\begin{tikzcd}[column sep = huge, row sep = huge, scale cd=1]
\varphi_{(\vec{j} \bullet \mathrm{Id}) \circ (\vec{\mathrm{Id}} \bullet i)}(M) \arrow[d, "\Xi_{\substack{i, \vec{\mathrm{Id}} \\ \mathrm{Id}, \vec{j}}}"] \arrow[rr, "\tau_{\substack{\vec{\mathrm{Id}} \bullet i, \vec{j} \bullet \mathrm{Id} \\ \vec{j} \bullet \mathrm{Id}, \vec{\mathrm{Id}} \bullet i}}"]
&&
\varphi_{(\vec{\mathrm{Id}} \bullet i) \circ (\vec{j} \bullet \mathrm{Id})} (M)
\\
\varphi_{(\vec{j} \circ \vec{\mathrm{Id}}) \bullet (\mathrm{Id} \circ i)} (M) \arrow[r, equal]
&
\varphi_{\vec{j} \bullet i}(M) \arrow[r, equal]
&
\varphi_{(\vec{\mathrm{Id}} \circ \vec{j})\bullet (i \circ \mathrm{Id})} \arrow[u, "\Xi_{\substack{\mathrm{Id}, \vec{j} \\ i, \mathrm{Id}}}"]
\end{tikzcd}
\]

Consider a diagram of purely homogeneous maps, with degrees pictured.
\[
\begin{tikzcd}
& M_{12} \arrow[dr, "f_{*2}", "j'"'] & \\
M_{11} \arrow[ur, "f_{1*}", "i'"'] \arrow[dr, "f_{*1}"', "i"] & & M_{22} \\
& M_{21} \arrow[ur, "f_{2*}"', "j"] & 
\end{tikzcd}
\]
We say that the diagram is \textit{$\mathscr{C}$-graded commutative} if $\left((i,j),(i', j')\right)\in \mathtt{T}$, and 
\[
\left(f_{2*} \circ_{\mathscr{C}} f_{*1}\right) = \tau_{\substack{i, i' \\ j, j'}} \left(f_{*2} \circ_{\mathscr{C}} f_{1*}\right).
\]
Note that $\left(f_{*2} \circ_{\mathscr{C}} f_{1*}\right)$ has degree $j' \circ i'$ and $\left(f_{2*} \circ_{\mathscr{C}} f_{*1}\right)$ has degree $j\circ i$, so $\tau_{\substack{i, i' \\ j, j'}}$ ensures their $\mathscr{C}$-degrees agree. This situation is abbreviated by including an arrow $\mathbin{\rotatebox[origin=c]{90}{$\Rightarrow$}}$ as in the following proposition.

\begin{proposition}
\label{prop:C-commutativity}
Given $\mathscr{C}$-graded commutative diagrams
\[
\begin{tikzcd}
& M_{12} \arrow[dr, "f_{*2}"] & \\
M_{11} \arrow[ur, "f_{1*}"] \arrow[dr, "f_{*1}"'] & & M_{22} \\
& M_{21} \arrow[ur, "f_{2*}"'] \arrow[uu, Rightarrow, shorten=2.5ex] & 
\end{tikzcd}
\]
and
\[
\begin{tikzcd}
& (M_{12})_1 \arrow[dr, "(f_{*2})_1"] & \\
(M_{11})_1 \arrow[ur, "(f_{1*})_1"] \arrow[dr, "(f_{*1})_1"'] & & (M_{22})_1 \\
& (M_{21})_1 \arrow[ur, "(f_{2*})_1"'] \arrow[uu, Rightarrow, shorten=2.5ex] & 
\end{tikzcd}
\quad,~\ldots,\quad
\begin{tikzcd}
& (M_{12})_k \arrow[dr, "(f_{*2})_k"] & \\
(M_{11})_k \arrow[ur, "(f_{1*})_k"] \arrow[dr, "(f_{*1})_k"'] & & (M_{22})_k \\
& (M_{21})_k \arrow[ur, "(f_{2*})_k"'] \arrow[uu, Rightarrow, shorten=2.5ex] & 
\end{tikzcd}
\]
the diagram
\[
\begin{tikzcd}[row sep=huge]
& \left((M_{12})_1, \ldots, (M_{12})_k\right)\otimes M_{12} \arrow[dr, "\vec{f}_{*2} \otimes f_{*2}"] & \\
\left((M_{11})_1, \ldots, (M_{11})_k\right)\otimes M_{11} \arrow[ur, "\vec{f}_{1*} \otimes f_{1*}"] \arrow[dr, "\vec{f}_{*1} \otimes f_{*1}"']& & \left((M_{22})_1, \ldots, (M_{22})_k\right)\otimes M_{22}\\
& \left((M_{21})_1, \ldots, (M_{21})_k\right)\otimes M_{21} \arrow[ur, "\vec{f}_{2*} \otimes f_{2*}"'] \arrow[uu, Rightarrow, shorten=2.5ex] & 
\end{tikzcd}
\]
is $\mathscr{C}$-graded commutative. Here, $\vec{f}_{1*} \otimes f_{1*}$ is shorthand for $\left((f_{1*})_1, \ldots, (f_{1*})_k\right) \otimes f_{1*}$, and so on.
\end{proposition}

\begin{proof}
This is simple, given Proposition \ref{prop:vertical&horizontalcomp} and equation (\ref{eq:CommSystComp}). We drop some notation in what follows; hopefully it is clear:
\begin{align*}
    ((\vec{f}_{*2} \otimes f_{*2}) \circ_{\mathscr{C}} (\vec{f}_{1*} \otimes f_{1*})) (\vec{m} \otimes m) 
    &= \Xi_{\substack{i', \vec{i}' \\ j', \vec{j}'}}^{-1} \left( (\vec{f}_{*2} \circ_{\mathscr{C}} \vec{f}_{1*}) \otimes (f_{*2} \circ_\mathscr{C} f_{1*}) \right) (\vec{m} \otimes m)
    \\
    &= \Xi_{\substack{i', \vec{i}' \\ j', \vec{j}'}}^{-1} \beta_{\vec{j}' \circ \vec{i}', j'\circ i'}^{-1} (\vec{f}_{*2} \circ_{\mathscr{C}} \vec{f}_{1*})(\vec{m}) \otimes (f_{*2} \circ_\mathscr{C} f_{1*}) (m)
    \\
    &= \Xi_{\substack{i', \vec{i}' \\ j', \vec{j}'}}^{-1} \beta_{\vec{j}' \circ \vec{i}', j'\circ i'}^{-1} \tau_{\substack{\vec{i}, \vec{i}' \\ \vec{j}, \vec{j}'}}^{-1} \tau_{\substack{i, i'\\ j, j'}}^{-1} (\vec{f}_{2*} \circ_{\mathscr{C}} \vec{f}_{*1})(\vec{m}) \otimes (f_{2*} \circ_\mathscr{C} f_{*1}) (m)
    \\
    &= \tau_{\substack{\vec{i} \bullet i, \vec{i}' \bullet i' \\ \vec{j} \bullet j, \vec{j}' \bullet j'}}^{-1} \Xi_{\substack{i, \vec{i} \\ j, \vec{j}}}^{-1} \beta_{\vec{j} \circ \vec{i}, j\circ i}^{-1} (\vec{f}_{2*} \circ_{\mathscr{C}} \vec{f}_{*1})(\vec{m}) \otimes (f_{2*} \circ_\mathscr{C} f_{*1}) (m)
    \\
    &= \tau_{\substack{\vec{i} \bullet i, \vec{i}' \bullet i' \\ \vec{j} \bullet j, \vec{j}' \bullet j'}}^{-1} \Xi_{\substack{i, \vec{i} \\ j, \vec{j}}}^{-1} \left((\vec{f}_{2*} \circ_{\mathscr{C}} \vec{f}_{*1}) \otimes (f_{2*} \circ_\mathscr{C} f_{*1})\right) (\vec{m} \otimes m) 
    \\
    &= \tau_{\substack{\vec{i} \bullet i, \vec{i}' \bullet i' \\ \vec{j} \bullet j, \vec{j}' \bullet j'}}^{-1} ((\vec{f}_{2*} \otimes f_{2*}) \circ_{\mathscr{C}} (\vec{f}_{*1} \otimes f_{*1})) (\vec{m} \otimes m).
\end{align*}
\end{proof}

\newpage

\section{An invariant of diskular tangles}
\label{Chapter:tangleinvariant}

In this section, we describe an invariant of diskular tangles. In \S \ref{ss:quickcomps}, we describe useful computational tools necessary to successive work, inspired by \cite{barnatan2006fast} but paying particular attention to the ``simplification'' of $\mathscr{G}$-grading shifts. We remark that, as in \cite{naisse2020odd}, our $\mathscr{G}$-grading system is a little too sensitive the $\mathscr{G}$-graded dg-multimodule we associate to a diskular tangle to be invariant under each Reidemeister move. However, we also describe a procedure (important to the results of section \S \ref{Chapter:OddCKProj}) which collapses $\mathscr{G}$-grading to a $q$-grading, in which case we obtain an honest tangle invariant. The work here is motivated by and serves as a generalization of \cite{naisse2020odd}. We write $\mathrm{Kom}(\cdot)$ to indicate the category of complexes which are bounded below in homological degree and of finite rank in each quantum or $\mathscr{G}$ degree.

\subsection{Quick computations in unified Khovanov homology}
\label{ss:quickcomps}

To begin, we will describe a few tools which will allow for quick computations in $\mathrm{Kom}\left(H^n\mathrm{Mod}_R^\mathscr{G}\right)$. In particular, we hope to use the methods introduced in \cite{barnatan2006fast}, but must develop others to deal with problems posed by $\mathscr{G}$-shifts.

\subsubsection{Delooping}

As an internal check, we can derive a formula for delooping in the current setting. A birth $\tikz[baseline={([yshift=-.5ex]current bounding box.center)}, scale=.5]{
	\draw (1,2) .. controls (1,1) and (2,1) .. (2,2);
	\draw (1,2) .. controls (1,1.75) and (2,1.75) .. (2,2);
	\draw (1,2) .. controls (1,2.25) and (2,2.25) .. (2,2);
}
: \emptyset \to \bigcirc$ induces a graded map $\mathcal{F}\left(\tikz[baseline={([yshift=-.5ex]current bounding box.center)}, scale=.5]{
	\draw (1,2) .. controls (1,1) and (2,1) .. (2,2);
	\draw (1,2) .. controls (1,1.75) and (2,1.75) .. (2,2);
	\draw (1,2) .. controls (1,2.25) and (2,2.25) .. (2,2);
}\right): \varphi_{\tikz[baseline={([yshift=-.5ex]current bounding box.center)}, scale=.3]{
	\draw (1,2) .. controls (1,1) and (2,1) .. (2,2);
	\draw (1,2) .. controls (1,1.75) and (2,1.75) .. (2,2);
	\draw (1,2) .. controls (1,2.25) and (2,2.25) .. (2,2);
}}(R) \to V$, since $\mathcal{F}(\emptyset) = R$ and $\mathcal{F}(\bigcirc) = V$. Notice that this $\mathcal{G}$-grading shift functor has only the effect of adding $(1,0)$ in the second coordinate (free loops are ignored in the first coorinate): $\varphi_{\tikz[baseline={([yshift=-.5ex]current bounding box.center)}, scale=0.3]{
	\draw (1,2) .. controls (1,1) and (2,1) .. (2,2);
	\draw (1,2) .. controls (1,1.75) and (2,1.75) .. (2,2);
	\draw (1,2) .. controls (1,2.25) and (2,2.25) .. (2,2);
}} \cong \{1,0\}$. So we have a graded map
\[\mathcal{F}\left(\tikz[baseline={([yshift=-.5ex]current bounding box.center)}, scale=0.5]{
	\draw (1,2) .. controls (1,1) and (2,1) .. (2,2);
	\draw (1,2) .. controls (1,1.75) and (2,1.75) .. (2,2);
	\draw (1,2) .. controls (1,2.25) and (2,2.25) .. (2,2);
}\right): R\{1,0\} \to V.\]
Similarly,
\[\mathcal{F}\left(\tikz[baseline={([yshift=-.5ex]current bounding box.center)}, scale=0.5]{
	\draw (1,2) .. controls (1,1) and (2,1) .. (2,2);
	\draw (1,2) .. controls (1,1.75) and (2,1.75) .. (2,2);
	\draw (1,2) .. controls (1,2.25) and (2,2.25) .. (2,2);
        \node at (1.5,1.55) {$\bullet$};
}\right): R\{0,-1\} \to V\]
is a graded map. The grading shift functors $\{u,v\}$ have clear inverses given by $\{-u,-v\}$. This fact, together with similar analysis on graded maps induced by deaths, yields the following array of graded maps:
\begin{equation}
\label{eq:deloopingdiagram}
\tikz[baseline={([yshift=-.5ex]current bounding box.center)}, scale=1]{
    \node(00) at (-2,0) {$\mathcal{F}(\bigcirc)$};
    \node(10) at (1, 1.5) {$\mathcal{F}(\emptyset)\{1,0\}$};
    \node(01) at (1,-1.5) {$\mathcal{F}(\emptyset) \{0,-1\}$};
    \node(11) at (4,0) {$\mathcal{F}(\bigcirc)$};
    \node at (1,0) {$\oplus$};
\draw[->] (00) -- node[pos=.25,above,arrows=-]
        {$
        \tikz[baseline=-4ex, scale=.2]{
            \draw[domain=0:180] plot ({2*cos(\x)}, {2*sin(\x)});
            \draw (-2,0) arc (180:360:2 and 0.6);
            \draw[dashed] (2,0) arc (0:180:2 and 0.6);
            \node at (0,0.25) {$\bullet$};
            \draw[->] (0,2.8) [partial ellipse=0:270:5ex and 2ex];
            }$}  
        (10);
\draw[->] (00) -- node[pos=.25,below,arrows=-]
        {$
        \tikz[baseline=-4ex, scale=.2]{
            \draw[domain=0:180] plot ({2*cos(\x)}, {2*sin(\x)});
            \draw (-2,0) arc (180:360:2 and 0.6);
            \draw[dashed] (2,0) arc (0:180:2 and 0.6);
            \draw[->] (0,2.8) [partial ellipse=0:270:5ex and 2ex];
            }$}  
        (01);
\draw[->] (10) -- node[pos=.75,above,arrows=-]
        {$
        \tikz[baseline=-4ex, scale=.2]{
            \draw [domain=180:360] plot ({2*cos(\x)}, {2*sin(\x)});
            \draw (-2,0) arc (180:360:2 and 0.6);
            \draw (2,0) arc (0:180:2 and 0.6);
            }$}
        (11);
\draw[->] (01) -- node[pos=.75,below=0.5em,arrows=-]
        {$
        \tikz[baseline=-4ex, scale=.2]{
            \draw [domain=180:360] plot ({2*cos(\x)}, {2*sin(\x)});
            \draw (-2,0) arc (180:360:2 and 0.6);
            \draw (2,0) arc (0:180:2 and 0.6);
            \node at (0,-0.5) {$\bullet$};
            }$}
        (11);
}
\end{equation}
It might seem pedantic, but we note that the arrows on the left-hand side of (\ref{eq:deloopingdiagram}) should be precomposed with the isomorphisms coming from natural transformations of grading shift functors
\[
\mathrm{Id} \Rightarrow \{1,0\} \circ \{-1,0\}
\]
and
\[
\mathrm{Id} \Rightarrow
 \{0, -1\} \circ \{0,1\}
\]
respectively, so that the maps on the left are graded with respect to our conventions. We will neglect writing these isomorphisms outside of special situations (e.g., the proof of Theorems \ref{thm:adjunction1} and \ref{thm:dualclosing}).

\begin{proposition}[Delooping]
$\mathcal{F}(\bigcirc) \cong \mathcal{F}(\emptyset)\{0,-1\} \oplus \mathcal{F}(\emptyset)\{1,0\}.$
\end{proposition}

\begin{proof}
This follows directly from the definition of $\mathcal{F}$. For example, the composition shown in diagram \ref{eq:deloopingdiagram} reads
\[
\mathcal{F}\left(
\tikz[baseline={([yshift=-.5ex]current bounding box.center)}, scale=0.2]{
    \draw [domain=180:360] plot ({2*cos(\x)}, {2*sin(\x)});
    \draw (-2,0) arc (180:360:2 and 0.6);
    \draw (2,0) arc (0:180:2 and 0.6);
    \begin{scope}[shift={(0,-5.5)}]
        \draw[domain=0:180] plot ({2*cos(\x)}, {2*sin(\x)});
        \draw (-2,0) arc (180:360:2 and 0.6);
        \draw[dashed] (2,0) arc (0:180:2 and 0.6);
        \node at (0,1.25) {$\bullet$};
        \draw[->] (0,2.8) [partial ellipse=0:270:5ex and 2ex];
    \end{scope}
}\right)
+
\mathcal{F}\left(
\tikz[baseline={([yshift=-.5ex]current bounding box.center)}, scale=0.2]{
    \draw [domain=180:360] plot ({2*cos(\x)}, {2*sin(\x)});
    \draw (-2,0) arc (180:360:2 and 0.6);
    \draw (2,0) arc (0:180:2 and 0.6);
    \node at (0,-1.25) {$\bullet$};
    \begin{scope}[shift={(0,-5.5)}]
        \draw[domain=0:180] plot ({2*cos(\x)}, {2*sin(\x)});
        \draw (-2,0) arc (180:360:2 and 0.6);
        \draw[dashed] (2,0) arc (0:180:2 and 0.6);
        \draw[->] (0,2.8) [partial ellipse=0:270:5ex and 2ex];
    \end{scope}
}\right)
=
\mathcal{F}\left(
\tikz[baseline={([yshift=-.5ex]current bounding box.center)}, scale=0.2]{
    \draw (-2,0) arc (180:360:2 and 0.6);
    \draw (2,0) arc (0:180:2 and 0.6);
    \draw (-2,-5.5) -- (-2,0);
    \draw (2,-5.5) -- (2,0);
    \begin{scope}[shift={(0,-5.5)}]
        \draw (-2,0) arc (180:360:2 and 0.6);
        \draw[dashed] (2,0) arc (0:180:2 and 0.6);
    \end{scope}
}\right)\,.
\]
One can verify this by checking that a dotted cylinder, followed by a positive death, and then a birth maps $v_+$ to $v_+$ and $v_-$ to zero, while a positive death, followed by a birth, and then a dotted cylinder maps $v_+$ to zero and $v_-$ to $v_-$. The other composition is also the identity: this amounts to showing that 
\[
\mathcal{F}\left(\tikz[baseline={([yshift=-.5ex]current bounding box.center)}, scale=0.2]{
  \draw (0,0) circle (2cm);
  \draw (-2,0) arc (180:360:2 and 0.6);
  \draw[dashed] (2,0) arc (0:180:2 and 0.6);
  \draw[->] (0,2.8) [partial ellipse=0:270:5ex and 2ex];
}\right) = \mathcal{F}\left( \tikz[baseline={([yshift=-.5ex]current bounding box.center)}, scale=0.2]{
  \draw (0,0) circle (2cm);
  \draw (-2,0) arc (180:360:2 and 0.6);
  \draw[dashed] (2,0) arc (0:180:2 and 0.6);
  \node at (0,1.2) {$\bullet$};
  \node at (0,-1.2) {$\bullet$};
  \draw[->] (0,2.8) [partial ellipse=0:270:5ex and 2ex];
} \right) = 0, \quad \text{and} \quad 
\mathcal{F}\left(\tikz[baseline={([yshift=-.5ex]current bounding box.center)}, scale=0.2]{
  \draw (0,0) circle (2cm);
  \draw (-2,0) arc (180:360:2 and 0.6);
  \draw[dashed] (2,0) arc (0:180:2 and 0.6);
  \node at (0,0) {$\bullet$};
  \draw[->] (0,2.8) [partial ellipse=0:270:5ex and 2ex];
}\right) = 1.
\]
That is, the tube-cutting and sphere relations hold in the category $H^n\mathrm{Mod}_R^\mathscr{G}$.

We should expect the gradings as they are since $\deg_R(v_+) = (1,0)$ and $\deg_R(v_-) = (0,-1)$, with $V = R\langle v_+, v_-\rangle$.
\end{proof}

\subsubsection{Simplifying grading shift functors}

In the even setting, delooping and Gaussian elimination allowed us to perform quick computations. To perform similar computations in the unified setting, we need to develop a system for simplifying $\mathscr{G}$-shifts. In the best cases, this means that $W$ consists of no ambiguous saddles, and is equivalent to a grading shift supported entirely in the $\mathbb{Z} \times \mathbb{Z}$ component; for example, we previously used that $\varphi_{\tikz[baseline={([yshift=-.5ex]current bounding box.center)}, scale=0.3]{
	\draw (1,2) .. controls (1,1) and (2,1) .. (2,2);
	\draw (1,2) .. controls (1,1.75) and (2,1.75) .. (2,2);
	\draw (1,2) .. controls (1,2.25) and (2,2.25) .. (2,2);
}} \cong \{1,0\}$. Usually this is not the case. Instead, given a cobordism $W:t \to t'$, we'd like to equate $\varphi_{(W,v)}$ with  $\varphi_{(\check{W}, u)}$ for some $u\in \mathbb{Z}\times \mathbb{Z}$ where $\check{W}$ is minimal. Recall that if $W$ is not minimal, it fails to be so up to some addition of tubes. Therefore, to approach the problem of simplify grading shift functors, it makes sense to ask how $\varphi_{(W,v)}$ behaves under tube-cutting.

\begin{proposition}
\label{ShiftingTubes}
Let $W:t \to t'$ be a cobordism. There is a minimal cobordism $\check{W}: t\to t'$ which is isotopic to $W$ outside of finitely many tubes. Denote the number of tubes in $W$ by $\tau_W$. Then
\[
\varphi_{(W,v)} \cong \varphi_{(\check{W}, v + \tau_W(-1,-1))}.
\]
\end{proposition}

\begin{proof}
Any tube in $W$ is either unambiguous (it is a split followed by a merge or vice versa) or it is ambiguous (it is impossible to determine the order of elementary cobordisms which constitute the tube without a given closure). Consider the (locally vertical) change of chronology $H: W \Rightarrow W'$ which changes all ambiguous tubes into unambiguous tubes, e.g.,
\[
H: 
\tikz[baseline={([yshift=-.5ex]current bounding box.center)}, scale=0.35]{
	\draw (0,0) .. controls (1.5,0) and (2,1) .. (.5,1); 
	\draw (4,0) .. controls (2.5,0) and (3,1) .. (4.5,1); 
	\filldraw[fill=white, draw=white] (.5,.55)  rectangle  (4,2);
	\draw[dashed] (0,0) .. controls (1.5,0) and (2,1) .. (.5,1); 
	\draw[dashed] (4,0) .. controls (2.5,0) and (3,1) .. (4.5,1); 
	\draw (1.375,.5) .. controls (1.375,1.5) and (3.125,1.5) .. (3.125,.5);
	\draw[yshift=1cm]  (0,2) .. controls (1.5,2) and (2,3) .. (.5,3); 
	\draw[yshift=1cm]  (4,2) .. controls (2.5,2) and (3,3) .. (4.5,3); 
	\draw[yshift=1cm]  (1.375,2.5) .. controls (1.375,1.5) and (3.125,1.5) .. (3.125,2.5);
	\draw (0,0) -- (0,3);
	\draw (.5,3.25) -- (.5,4);
	\draw[dashed] (.5,1) -- (.5,4);
	\draw (4,0) -- (4,3);
	\draw (4.5,1) -- (4.5,4);
}
\Rightarrow
\tikz[baseline={([yshift=-.5ex]current bounding box.center)}, scale=0.28]{
	\draw (0,0) .. controls (1.5,0) and (2,1) .. (.5,1); 
	\draw (4,0) .. controls (2.5,0) and (3,1) .. (4.5,1); 
	\filldraw[fill=white, draw=white] (.5,.55)  rectangle  (4.5,2);
	\draw[dashed] (0,0) .. controls (1.5,0) and (2,1) .. (.5,1); 
	\draw[dashed] (4,0) .. controls (2.5,0) and (3,1) .. (4.5,1); 
	\draw (0,0) -- (0,4);
	\draw[dashed] (.5,1) -- (.5,5);
	\draw (1.375,.5) -- (1.375,2.5);
	\draw (3.125,.5) -- (3.125,2.5);
	\draw (.5,4) -- (.5,5) -- (4.5,5) -- (4.5,4);
	\draw (1.375,2.5) .. controls (1.375,3.5) and (3.125,3.5) .. (3.125,2.5);
	\draw (0,2) -- (0,4) -- (4,4);
	\draw (8,1) -- (8.5,1) -- (8.5,5);
	\draw (5.375,2.5) .. controls (5.375,1.5) and (7.125,1.5) .. (7.125,2.5);
	\draw[dashed] (4.5,5) -- (4.5,1) -- (8.5,1) -- (8.5,3);
	\draw (4,4) -- (4,0) -- (8,0) -- (8,2);
	\draw (4,4) .. controls (5.5,4) and (6,5) .. (4.5,5); 
	\draw (8,4) .. controls (6.5,4) and (7,5) .. (8.5,5); 
	\draw (8,2) -- (8,4);
	\draw (8.5,3) -- (8.5,5);
	\draw (5.375,2.5) -- (5.375,4.5);
	\draw (7.125,2.5) -- (7.125,4.5);
}
\]
wherever ambiguous tubes are present. From our analysis earlier, there is an induced natural transformation $\varphi_H: \varphi_W\Rightarrow \varphi_{W'}$. Note that $\deg(1_{\overline{a}} W 1_b) = \deg(1_{\overline{a}} W' 1_b)$ since any tube in $W$ corresponds to the addition of $(-1,-1)$ in degree on any closure, ambiguous or not. This implies that $\varphi_{(W,v)} \cong \varphi_{(W',v)}$. Since each tube in $W'$ is unambiguous, we know that each tube in $W'$ acts as a degree $(-1,-1)$ shift, so the result follows.
\end{proof}

A consequence of this proposition is that \textit{all} grading shift functors have inverses, not just $\{u,v\}$.

\begin{corollary}
For any pair $(W:t\to t', v)$, $\varphi_{(W,v)}$ has a left inverse
\[\varphi_{(W,v)}^{-1} = \varphi_{(\overline{W}, -v + \tau_{\overline{W}\circ W}(1,1))}\]
(where $\overline{W}: t' \to t$ is the mirror image of $W$) in the sense that
\[\varphi_{(W,v)}^{-1} \circ \varphi_{(W,v)} \cong \varphi_{\mathbbm{1}_t}.\]
\end{corollary}

\begin{proof}
If the composition $\overline{W}\circ W$ produces any tubes, the contribution by these tubes on the second coordinate are killed by the addition of the term $\tau_{\overline{W}\circ W}(1,1)$.
\end{proof}

\begin{example*}
An elementary saddle cobordism $\tikz[baseline={([yshift=-.5ex]current bounding box.center)}, scale=0.45]
{
    \begin{scope}[rotate=90]
	\draw[dotted] (3,-2) circle(0.707);
	\draw (2.5,-1.5) .. controls (2.75,-1.75) and (3.25,-1.75) .. (3.5,-1.5);
	\draw (2.5,-2.5) .. controls  (2.75,-2.25) and (3.25,-2.25) .. (3.5,-2.5);
        \draw[red,thick] (3,-1.7) -- (3,-2.3);
    \end{scope}
} : \tikz[baseline={([yshift=-.5ex]current bounding box.center)}, scale=0.45]
{
    \begin{scope}[rotate=90]
	\draw[dotted] (3,-2) circle(0.707);
	\draw (2.5,-1.5) .. controls (2.75,-1.75) and (3.25,-1.75) .. (3.5,-1.5);
	\draw (2.5,-2.5) .. controls  (2.75,-2.25) and (3.25,-2.25) .. (3.5,-2.5);
    \end{scope}
} \to \tikz[baseline={([yshift=-.5ex]current bounding box.center)}, scale=0.45]
{
	\draw[dotted] (3,-2) circle(0.707);
	\draw (2.5,-1.5) .. controls (2.75,-1.75) and (3.25,-1.75) .. (3.5,-1.5);
	\draw (2.5,-2.5) .. controls  (2.75,-2.25) and (3.25,-2.25) .. (3.5,-2.5);
}$ induces the graded map 
\[
\mathcal{F}\left(\tikz[baseline={([yshift=-.5ex]current bounding box.center)}, scale=0.45]
{
    \begin{scope}[rotate=90]
	\draw[dotted] (3,-2) circle(0.707);
	\draw (2.5,-1.5) .. controls (2.75,-1.75) and (3.25,-1.75) .. (3.5,-1.5);
	\draw (2.5,-2.5) .. controls  (2.75,-2.25) and (3.25,-2.25) .. (3.5,-2.5);
        \draw[red,thick] (3,-1.7) -- (3,-2.3);
    \end{scope}
}\right) : \varphi_{\tikz[baseline={([yshift=-.5ex]current bounding box.center)}, scale=0.45]
{
    \begin{scope}[rotate=90]
	\draw[dotted] (3,-2) circle(0.707);
	\draw (2.5,-1.5) .. controls (2.75,-1.75) and (3.25,-1.75) .. (3.5,-1.5);
	\draw (2.5,-2.5) .. controls  (2.75,-2.25) and (3.25,-2.25) .. (3.5,-2.5);
        \draw[red,thick] (3,-1.7) -- (3,-2.3);
    \end{scope}
}} \mathcal{F}\left(\tikz[baseline={([yshift=-.5ex]current bounding box.center)}, scale=0.45]
{
    \begin{scope}[rotate=90]
	\draw[dotted] (3,-2) circle(0.707);
	\draw (2.5,-1.5) .. controls (2.75,-1.75) and (3.25,-1.75) .. (3.5,-1.5);
	\draw (2.5,-2.5) .. controls  (2.75,-2.25) and (3.25,-2.25) .. (3.5,-2.5);
    \end{scope}
}\right) \to \mathcal{F}\left(\tikz[baseline={([yshift=-.5ex]current bounding box.center)}, scale=0.45]
{
	\draw[dotted] (3,-2) circle(0.707);
	\draw (2.5,-1.5) .. controls (2.75,-1.75) and (3.25,-1.75) .. (3.5,-1.5);
	\draw (2.5,-2.5) .. controls  (2.75,-2.25) and (3.25,-2.25) .. (3.5,-2.5);
}\right).
\]
Consider the isomorphism induced by change of chronology:
\[
\varphi_H: \mathrm{Id} \Rightarrow \varphi_{\tikz[baseline={([yshift=-.5ex]current bounding box.center)}, scale=0.45]
{
    \begin{scope}[rotate=90]
	\draw[dotted] (3,-2) circle(0.707);
	\draw (2.5,-1.5) .. controls (2.75,-1.75) and (3.25,-1.75) .. (3.5,-1.5);
	\draw (2.5,-2.5) .. controls  (2.75,-2.25) and (3.25,-2.25) .. (3.5,-2.5);
        \draw[red,thick] (3,-1.7) -- (3,-2.3);
    \end{scope}
}}^{-1} \circ \varphi_{\tikz[baseline={([yshift=-.5ex]current bounding box.center)}, scale=0.45]
{
    \begin{scope}[rotate=90]
	\draw[dotted] (3,-2) circle(0.707);
	\draw (2.5,-1.5) .. controls (2.75,-1.75) and (3.25,-1.75) .. (3.5,-1.5);
	\draw (2.5,-2.5) .. controls  (2.75,-2.25) and (3.25,-2.25) .. (3.5,-2.5);
        \draw[red,thick] (3,-1.7) -- (3,-2.3);
    \end{scope}
}}.
\]
Then, precomposing with $\varphi_H$, the saddle can be reinterpreted as the following graded map.
\[
\mathcal{F}\left(\tikz[baseline={([yshift=-.5ex]current bounding box.center)}, scale=0.45]
{
    \begin{scope}[rotate=90]
	\draw[dotted] (3,-2) circle(0.707);
	\draw (2.5,-1.5) .. controls (2.75,-1.75) and (3.25,-1.75) .. (3.5,-1.5);
	\draw (2.5,-2.5) .. controls  (2.75,-2.25) and (3.25,-2.25) .. (3.5,-2.5);
        \draw[red,thick] (3,-1.7) -- (3,-2.3);
    \end{scope}
}\right) \circ \varphi_H : \mathcal{F}\left(\tikz[baseline={([yshift=-.5ex]current bounding box.center)}, scale=0.45]
{
    \begin{scope}[rotate=90]
	\draw[dotted] (3,-2) circle(0.707);
	\draw (2.5,-1.5) .. controls (2.75,-1.75) and (3.25,-1.75) .. (3.5,-1.5);
	\draw (2.5,-2.5) .. controls  (2.75,-2.25) and (3.25,-2.25) .. (3.5,-2.5);
    \end{scope}
}\right) \to \varphi_{\tikz[baseline={([yshift=-.5ex]current bounding box.center)}, scale=0.45]
{
    \begin{scope}[rotate=90]
	\draw[dotted] (3,-2) circle(0.707);
	\draw (2.5,-1.5) .. controls (2.75,-1.75) and (3.25,-1.75) .. (3.5,-1.5);
	\draw (2.5,-2.5) .. controls  (2.75,-2.25) and (3.25,-2.25) .. (3.5,-2.5);
        \draw[red,thick] (3,-1.7) -- (3,-2.3);
    \end{scope}
}}^{-1} \mathcal{F}\left(\tikz[baseline={([yshift=-.5ex]current bounding box.center)}, scale=0.45]
{
	\draw[dotted] (3,-2) circle(0.707);
	\draw (2.5,-1.5) .. controls (2.75,-1.75) and (3.25,-1.75) .. (3.5,-1.5);
	\draw (2.5,-2.5) .. controls  (2.75,-2.25) and (3.25,-2.25) .. (3.5,-2.5);
}\right).
\]
We compute that $\varphi_{\tikz[baseline={([yshift=-.5ex]current bounding box.center)}, scale=0.45]
{
    \begin{scope}[rotate=90]
	\draw[dotted] (3,-2) circle(0.707);
	\draw (2.5,-1.5) .. controls (2.75,-1.75) and (3.25,-1.75) .. (3.5,-1.5);
	\draw (2.5,-2.5) .. controls  (2.75,-2.25) and (3.25,-2.25) .. (3.5,-2.5);
        \draw[red,thick] (3,-1.7) -- (3,-2.3);
    \end{scope}
}}^{-1} = \varphi_{\left(\tikz[baseline={([yshift=-.5ex]current bounding box.center)}, scale=0.45]
{
	\draw[dotted] (3,-2) circle(0.707);
	\draw (2.5,-1.5) .. controls (2.75,-1.75) and (3.25,-1.75) .. (3.5,-1.5);
	\draw (2.5,-2.5) .. controls  (2.75,-2.25) and (3.25,-2.25) .. (3.5,-2.5);
        \draw[red,thick] (3,-1.7) -- (3,-2.3);
}\,,~ (1,1)\right)}$ since $\tikz[baseline={([yshift=-.5ex]current bounding box.center)}, scale=0.45]
{
	\draw[dotted] (3,-2) circle(0.707);
	\draw (2.5,-1.5) .. controls (2.75,-1.75) and (3.25,-1.75) .. (3.5,-1.5);
	\draw (2.5,-2.5) .. controls  (2.75,-2.25) and (3.25,-2.25) .. (3.5,-2.5);
        \draw[red,thick] (3,-1.7) -- (3,-2.3);
} \circ \tikz[baseline={([yshift=-.5ex]current bounding box.center)}, scale=0.45]
{
    \begin{scope}[rotate=90]
	\draw[dotted] (3,-2) circle(0.707);
	\draw (2.5,-1.5) .. controls (2.75,-1.75) and (3.25,-1.75) .. (3.5,-1.5);
	\draw (2.5,-2.5) .. controls  (2.75,-2.25) and (3.25,-2.25) .. (3.5,-2.5);
        \draw[red,thick] (3,-1.7) -- (3,-2.3);
    \end{scope}
}$ produces a tube. In general,
\[
\varphi_{\left(\tikz[baseline={([yshift=-.5ex]current bounding box.center)}, scale=0.45]
{
    \begin{scope}[rotate=90]
	\draw[dotted] (3,-2) circle(0.707);
	\draw (2.5,-1.5) .. controls (2.75,-1.75) and (3.25,-1.75) .. (3.5,-1.5);
	\draw (2.5,-2.5) .. controls  (2.75,-2.25) and (3.25,-2.25) .. (3.5,-2.5);
        \draw[red,thick] (3,-1.7) -- (3,-2.3);
    \end{scope}
}\,,~ (u,v)\right)}^{-1} = \varphi_{\left(\tikz[baseline={([yshift=-.5ex]current bounding box.center)}, scale=0.45]
{
	\draw[dotted] (3,-2) circle(0.707);
	\draw (2.5,-1.5) .. controls (2.75,-1.75) and (3.25,-1.75) .. (3.5,-1.5);
	\draw (2.5,-2.5) .. controls  (2.75,-2.25) and (3.25,-2.25) .. (3.5,-2.5);
        \draw[red,thick] (3,-1.7) -- (3,-2.3);
}\,,~ (1-u,1-v)\right)}.
\]
\end{example*}

\begin{remark}
Returning to diagram (\ref{eq:deloopingdiagram}), we see that the natural transformations of grading shifting functors actually take the forms
\[
\varphi_H: \mathrm{Id} \Rightarrow \varphi_{\tikz[baseline={([yshift=-.5ex]current bounding box.center)}, scale=0.2]{
            \draw[domain=0:180] plot ({2*cos(\x)}, {2*sin(\x)});
            \draw (-2,0) arc (180:360:2 and 0.6);
            \draw[dashed] (2,0) arc (0:180:2 and 0.6);
            \node at (0,1.25) {$\bullet$};
            }}^{-1} \circ \varphi_{\tikz[baseline={([yshift=-.5ex]current bounding box.center)}, scale=0.2]{
            \draw[domain=0:180] plot ({2*cos(\x)}, {2*sin(\x)});
            \draw (-2,0) arc (180:360:2 and 0.6);
            \draw[dashed] (2,0) arc (0:180:2 and 0.6);
            \node at (0,1.25) {$\bullet$};
            }} \cong \{1,0\} \circ \{-1,0\}
\]
and
\[
\varphi_H: \mathrm{Id} \Rightarrow
\varphi_{\tikz[baseline=-4ex, scale=.2]{
            \draw[domain=0:180] plot ({2*cos(\x)}, {2*sin(\x)});
            \draw (-2,0) arc (180:360:2 and 0.6);
            \draw[dashed] (2,0) arc (0:180:2 and 0.6);
            }}^{-1} \circ \varphi_{\tikz[baseline=-4ex, scale=.2]{
            \draw[domain=0:180] plot ({2*cos(\x)}, {2*sin(\x)});
            \draw (-2,0) arc (180:360:2 and 0.6);
            \draw[dashed] (2,0) arc (0:180:2 and 0.6);
            }} \cong \{0, -1\} \circ \{0,1\}.
\]
\end{remark}

The dismissal of free loops by the $\mathcal{G}$-shifting system leads to another possibility for simplification of grading shift functors. We will frequently use the following simplification while cooking up projectors; see \cite{naisse2020odd} for a proof.

\begin{proposition}\label{FLshift}
Suppose $W:t \to t'$ is a cobordism and $t$ contains a free loop $\ell$. Then there is a natural isomorphism 
\[\varphi_{(W, (u,v))} \cong \varphi_{(W', (u-1, v))} \quad \text{by} \quad m \mapsto \lambda((1,0), \deg(1_{\overline{a}} W 1_b)) m\]
where $W':t - \ell \to t'$ is given by gluing a birth under the free loop in $t$, and $m \in M_{g:a \to b}$. If, on the other hand, $t'$ contains the free loop, the natural isomorphism is on the nose: 
\[\varphi_{(W, (u,v))} \cong \varphi_{(W'', (u, v-1))}\]
and $W'': t \to t' - \ell$ is given by gluing a death above $W$.
\end{proposition}

\begin{example*}
Here is a way we may use the preceding proposition. Consider the $\mathscr{G}$-grading shifting map 
$\varphi_{\tikz[baseline={([yshift=-.5ex]current bounding box.center)}, scale=.35]{
        \draw (0,1) to[out=-70,in=250] (1,1);
        \draw(0,-1) to[out=70,in=110] node[pos=.5,inner sep=0](bot){} (1,-1);
        \draw (.5,0) circle (.35);
        \draw[red] (.5,-.35)--(bot);
        \draw[red] (.5, .35) -- (.5, .7);
        }}$ (the choice of chronology is unimportant). Then Proposition \ref{FLshift} says that this grading shift is isomorphic to the grading shift $\varphi_{(W', (-1,0))}$, where $W'$ is pictured below.
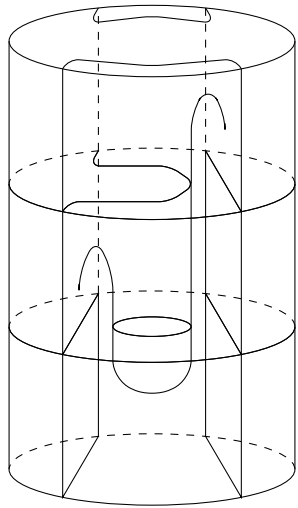
\[
\tikz[baseline=6ex, scale = .475]{
        \draw (-4,0) arc (180:360:4 and 1);
        \draw[dashed] (4,0) arc (0:180:4 and 1);
        \draw (-4,4) arc (180:360:4 and 1);
        \draw (-4,4) -- (-4,0);
        \draw (4,4) -- (4,0);
        \draw (-1.5,0.927) -- (-2.5, -0.781) -- (-2.5, -0.781+4) -- (-1.5,4.927) -- (-1.5,0.927);
        \draw (1.5,0.927) -- (2.5, -0.781) -- (2.5, -0.781+4) -- (1.5,4.927) -- (1.5,0.927);
        \draw (-1.095,4) arc (180:360:1.095 and 0.274);
        \draw (1.095,4) arc (0:180:1.095 and 0.274);
        \draw[rounded corners=7.5pt] (-1.095,4) -- (-1.095,2.5) -- (0,2) -- (1.095,2.5) -- (1.095,4);
\begin{scope}[shift={(0,4)}]
        \draw (-4,0) arc (180:360:4 and 1);
        \draw[dashed] (4,0) arc (0:180:4 and 1);
        \draw (-4,4) arc (180:360:4 and 1);
        \draw (-4,4) -- (-4,0);
        \draw (4,4) -- (4,0);
        \draw (-1.5,0.927) -- (-2.5, -0.781) -- (-2.5, -0.781+4);
        \draw[dashed] (-1.5,0.927) -- (-1.5,4.927);
        \draw[rounded corners] (-1.5,4.927) -- (-1.75, 4.5) -- (0.5,4.5) -- (1.2175,4) -- (0.5, 3.5) -- (-2.335, 3.5) -- (-2.5, -0.781+4);
        \draw (1.5,0.927) -- (2.5, -0.781) -- (2.5, -0.781+4) -- (1.5,4.927) -- (1.5,0.927);
        \draw (-1.095,0) arc (180:360:1.095 and 0.274);
        \draw (1.095,0) arc (0:180:1.095 and 0.274);
        \draw (1.095,0) -- (1.095,4);
        \draw[rounded corners=6.5pt] (-1.095,0) -- (-1.095,1.5) -- (-1.5695, 2.5) -- (-2.043, 1.5) -- (-2.043, 1.2);
\end{scope}
\begin{scope}[shift={(0,8)}]
        \draw (-4,0) arc (180:360:4 and 1);
        \draw[dashed] (4,0) arc (0:180:4 and 1);
        \draw (-4,4) arc (180:360:4 and 1);
        \draw (4,4) arc (0:180:4 and 1);
        \draw (-4,4) -- (-4,0);
        \draw (4,4) -- (4,0);
        \draw[dashed] (1.5,4.927) -- (1.5,0.927);
        \draw (1.5,0.927) -- (2.5, -0.781) -- (2.5, -0.781+4);
        \draw[rounded corners] (-1.5,0.927) -- (-1.75, 0.5) -- (0.5,0.5) -- (1.2175,0) -- (0.5, -0.5) -- (-2.335, -0.5) -- (-2.5, -0.781);
        \draw[dashed] (-1.5,0.927) -- (-1.5,4.927);
        \draw (-2.5, -0.781) -- (-2.5, -0.781+4);
        \draw[rounded corners] (-1.5, 4.927) -- (-1.75, 4.5) -- (0,4.7) -- (1.75, 4.5) -- (1.5, 4.927);
        \draw[rounded corners] (-2.5,-0.781+4) -- (-2.335, -0.5+4) -- (0,-0.7+4) -- (2.335, -0.5+4) -- (2.5, -0.781+4);
        \draw[rounded corners=6.5pt] (1.1,0) -- (1.1,2) -- (1.5715,2.75) -- (2.043,2) -- (2.043,1.7);
\end{scope}
}
\]
Of course, $W'$ is isotopic to an elementary saddle $\tikz[baseline={([yshift=-.5ex]current bounding box.center)}, scale=.35]{
    \draw (0,.5) to[out=-70,in=250] (1,.5);
    \draw (0,-.5) to[out=70,in=110] (1,-.5);
    \draw[red] (.5, -.25) -- (.5, .25);
    }$, so Proposition \ref{ShiftingTubes} allows us to conclude that
\[
\varphi_{\tikz[baseline={([yshift=-.5ex]current bounding box.center)}, scale=.35]{
        \draw (0,1) to[out=-70,in=250] (1,1);
        \draw(0,-1) to[out=70,in=110] node[pos=.5,inner sep=0](bot){} (1,-1);
        \draw (.5,0) circle (.35);
        \draw[red] (.5,-.35)--(bot);
        \draw[red] (.5, .35) -- (.5, .7);
        }}
\cong
\varphi_{\left(\tikz[baseline={([yshift=-.5ex]current bounding box.center)}, scale=.35]{
    \draw (0,.5) to[out=-70,in=250] (1,.5);
    \draw (0,-.5) to[out=70,in=110] (1,-.5);
    \draw[red] (.5, -.25) -- (.5, .25);
    }, (-1, 0) \right)}.
\]
\end{example*}

The examples of this subsection illustrates a peculiarity of computations in $\mathrm{Kom}(H^n\mathrm{Mod}_R^\mathscr{G})$---that a $\mathscr{G}$-grading shift has a few different representatives. A difficulty in coming work (cf. the proof of Lemma \ref{lem:r3}) is choosing the correct representative.

\subsection{Tangle invariant}
\label{ss:TangleInvariant}

In this subsection, we finally construct an invariant of diskular tangles. The work presented here is motivated by and serves as a generalization of Section 6.4 in \cite{naisse2020odd}.

Suppose $T$ is a diskular $(n_1,\ldots, n_k; m)$-tangle with $c$-many crossings. We will continue under the assumption that $T$ carries an orientation. Then $T$ defines an oriented $(\underbrace{2, \ldots, 2}_{c\text{ times}}, n_1, \ldots, n_k; m)$-planar arc diagram $D_T$ by replacing each crossing of $T$ with a new diskular region with four endpoints; consult the schematic below.
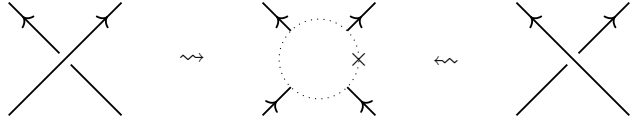
\[
\tikz[baseline={([yshift=-.5ex]current bounding box.center)}, scale=1.5]
{
        \draw[knot] (0,0) -- (1,1);
        \draw[knot, ->] (0.874, 0.874) -- (0.875, 0.875);
        \draw[knot] (1,0) -- (0.55,0.45);
        \draw[knot] (0.45, 0.55) -- (0, 1);
        \draw[knot, ->] (0.126, 0.874) -- (0.125, 0.875);
    \begin{scope}[xshift=1.125cm]
        \node at (0.5, 0.5) {$\rightsquigarrow$};
    \end{scope}
    \begin{scope}[xshift=2.25cm]
        \draw[dotted] (0.5, 0.5)circle (0.707/2);
        \node at (0.707/2 + 0.5, 0.5) {$\times$};
        \draw[knot] (0,0) -- (0.25, 0.25);
        \draw[knot, ->] (0.124, 0.124) -- (0.125, 0.125);
        \draw[knot] (1,0) -- (0.75, 0.25);
        \draw[knot, ->] (0.876, 0.124) -- (0.875, 0.125);
        \draw[knot] (0.25, 0.75) -- (0,1);
        \draw[knot, ->] (0.126, 0.874) -- (0.125, 0.875);
        \draw[knot] (0.75, 0.75) -- (1,1);
        \draw[knot, ->] (0.874, 0.874) -- (0.875, 0.875);
    \end{scope}
    \begin{scope}[xshift=3.375cm]
        \node[rotate=180] at (0.5, 0.5) {$\rightsquigarrow$};
    \end{scope}
    \begin{scope}[xshift=4.5cm]
        \draw[knot] (0,0) -- (0.45, 0.45);
        \draw[knot] (0.55, 0.55) -- (1,1);
        \draw[knot, ->] (0.874, 0.874) -- (0.875, 0.875);
        \draw[knot] (1,0) -- (0,1);
        \draw[knot, ->] (0.126, 0.874) -- (0.125, 0.875);
    \end{scope}
}
\]
Denote the crossings of $T$ by $x_1, \ldots, x_c$ and define the complex
\[
\mathrm{Kh}(T) := \left(\mathrm{Kh}(x_1), \ldots, \mathrm{Kh} (x_c)\right) \otimes_{(H^2, \ldots, H^2)} \mathcal{F}(D_T)
\]
where 
\begin{align*}
\mathrm{Kh}\left( 
\tikz[baseline={([yshift=-.5ex]current bounding box.center)}, scale=.75]
{
	\draw[dotted] (.5,.5) circle(0.707);
	\draw[knot, ->](0,0) -- (1,1);
	\fill[fill=white] (.5,.5) circle (.15);
	\draw[knot, ->](1,0) -- (0,1);
}
 \right) 
\ &:= \ 
\mathrm{Cone}\left( 
 \varphi_{\tikz[baseline=-6.5ex, scale=.45]
{
	\draw[dotted] (3,-2) circle(0.707);
        \draw[knot, red,thick] (3,-1.7) -- (3,-2.3);
	\draw[knot] (2.5,-1.5) .. controls (2.75,-1.75) and (3.25,-1.75) .. (3.5,-1.5);
	\draw[knot] (2.5,-2.5) .. controls  (2.75,-2.25) and (3.25,-2.25) .. (3.5,-2.5);
}}
\mathcal{F}\left(
\tikz[baseline={([yshift=-.5ex]current bounding box.center)}, scale=.75]
{
	\draw[dotted] (3,-2) circle(0.707);
	\draw[knot] (2.5,-1.5) .. controls (2.75,-1.75) and (3.25,-1.75) .. (3.5,-1.5);
	\draw[knot] (2.5,-2.5) .. controls  (2.75,-2.25) and (3.25,-2.25) .. (3.5,-2.5);
}
  \right) 
\xrightarrow{\mathcal{F}\left( \tikz[baseline=-6.5ex, scale=.45]
{
	\draw[dotted] (3,-2) circle(0.707);
        \draw[knot, red,thick] (3,-1.7) -- (3,-2.3);
	\draw[knot] (2.5,-1.5) .. controls (2.75,-1.75) and (3.25,-1.75) .. (3.5,-1.5);
	\draw[knot] (2.5,-2.5) .. controls  (2.75,-2.25) and (3.25,-2.25) .. (3.5,-2.5);
} \right)}
\myvec{
\mathcal{F}\left(
\tikz[baseline={([yshift=-.5ex]current bounding box.center)}, scale=.75]
{
	\draw[dotted] (-2,-2) circle(0.707);
	\draw[knot] (-1.5,-1.5) .. controls (-1.75,-1.75) and (-1.75,-2.25) .. (-1.5,-2.5);
	\draw[knot] (-2.5,-1.5) .. controls (-2.25,-1.75) and (-2.25,-2.25) ..  (-2.5,-2.5);
}  
 \right)
 }
\right) \{-1,0\},~ \text{and}
\\
\mathrm{Kh}\left(
\tikz[baseline={([yshift=-.5ex]current bounding box.center)}, scale=.75]
{
	\draw[dotted] (.5,.5) circle(0.707);
	\draw[knot, ->](1,0) -- (0,1);
	\fill[fill=white] (.5,.5) circle (.15);
	\draw[knot, ->](0,0) -- (1,1);
}
 \right) 
\ &:= \ 
\mathrm{Cone}\left( 
\myvec{
\mathcal{F}\left(
\tikz[baseline={([yshift=-.5ex]current bounding box.center)}, scale=.75]
{
	\draw[dotted] (-2,-2) circle(0.707);
	\draw[knot] (-1.5,-1.5) .. controls (-1.75,-1.75) and (-1.75,-2.25) .. (-1.5,-2.5);
	\draw[knot] (-2.5,-1.5) .. controls (-2.25,-1.75) and (-2.25,-2.25) ..  (-2.5,-2.5);
} 
  \right) 
  }
\xrightarrow{
\mathcal{F}\left( \tikz[baseline=8.25ex, scale=.45]
{
    \begin{scope}[rotate=90]
	\draw[dotted] (3,-2) circle(0.707);
        \draw[knot, red,thick] (3,-1.7) -- (3,-2.3);
	\draw[knot] (2.5,-1.5) .. controls (2.75,-1.75) and (3.25,-1.75) .. (3.5,-1.5);
	\draw[knot] (2.5,-2.5) .. controls  (2.75,-2.25) and (3.25,-2.25) .. (3.5,-2.5);
    \end{scope}
} \right)
\circ \varphi_H}
\varphi_{\tikz[baseline=8.25ex, scale=.45]
{
    \begin{scope}[rotate=90]
	\draw[dotted] (3,-2) circle(0.707);
        \draw[knot, red,thick] (3,-1.7) -- (3,-2.3);
	\draw[knot] (2.5,-1.5) .. controls (2.75,-1.75) and (3.25,-1.75) .. (3.5,-1.5);
	\draw[knot] (2.5,-2.5) .. controls  (2.75,-2.25) and (3.25,-2.25) .. (3.5,-2.5);
    \end{scope}
}} ^{-1}
\mathcal{F}\left(
\tikz[baseline={([yshift=-.5ex]current bounding box.center)}, scale=.75]
{
	\draw[dotted] (3,-2) circle(0.707);
	\draw[knot] (2.5,-1.5) .. controls (2.75,-1.75) and (3.25,-1.75) .. (3.5,-1.5);
	\draw[knot] (2.5,-2.5) .. controls  (2.75,-2.25) and (3.25,-2.25) .. (3.5,-2.5);
}
\right)
\right)\{0,1\},
\end{align*}
for $\varphi_H: \mathrm{Id} \Rightarrow \varphi_{\tikz[baseline=8.25ex, scale=.45]
{
    \begin{scope}[rotate=90]
	\draw[dotted] (3,-2) circle(0.707);
        \draw[knot, red,thick] (3,-1.7) -- (3,-2.3);
	\draw[knot] (2.5,-1.5) .. controls (2.75,-1.75) and (3.25,-1.75) .. (3.5,-1.5);
	\draw[knot] (2.5,-2.5) .. controls  (2.75,-2.25) and (3.25,-2.25) .. (3.5,-2.5);
    \end{scope}
}} ^{-1} \circ \varphi_{\tikz[baseline=8.25ex, scale=.45]
{
    \begin{scope}[rotate=90]
	\draw[dotted] (3,-2) circle(0.707);
        \draw[knot, red,thick] (3,-1.7) -- (3,-2.3);
	\draw[knot] (2.5,-1.5) .. controls (2.75,-1.75) and (3.25,-1.75) .. (3.5,-1.5);
	\draw[knot] (2.5,-2.5) .. controls  (2.75,-2.25) and (3.25,-2.25) .. (3.5,-2.5);
    \end{scope}
}}$. Recall that the $\myvec{\text{underlined}}$ entry is in homological degree zero.

The reader should compare this with the unoriented case, where we have
\[
\mathcal{F}(T) \cong \left(\mathcal{F}(x_1),\ldots, \mathcal{F}(x_c)\right) \otimes_{(H^2, \ldots, H^2)} \mathcal{F}(D_T)
\]
by Theorem \ref{thm:multigluing}. So, we would expect the following lemma.

\begin{lemma}
\label{lem:khcoherence}
For any diskular tangle $T$, there exists a shifting functor $\varphi$ and integer $\ell$ such that
\[
\mathrm{Kh}(T) \cong \varphi(\mathcal{F}(T))[\ell]
\]
\end{lemma}

\begin{proof}
Recall that the dg-multimodule associated to a single crossing is given by
\[
\mathcal{F}\left(
\tikz[baseline={([yshift=-.5ex]current bounding box.center)}, scale=.75]
{
	\draw[dotted] (.5,.5) circle(0.707);
	\draw[knot,-](1,0) -- (0,1);
	\fill[fill=white] (.5,.5) circle (.15);
	\draw[knot,-](0,0) -- (1,1);
}
 \right)  = 
 \mathrm{Cone} \left(
 \myvec{
 \varphi_{\tikz[baseline=8.25ex, scale=.45]
{
    \begin{scope}[rotate=90]
	\draw[dotted] (3,-2) circle(0.707);
        \draw[knot, red,thick] (3,-1.7) -- (3,-2.3);
	\draw[knot] (2.5,-1.5) .. controls (2.75,-1.75) and (3.25,-1.75) .. (3.5,-1.5);
	\draw[knot] (2.5,-2.5) .. controls  (2.75,-2.25) and (3.25,-2.25) .. (3.5,-2.5);
    \end{scope}
}}
\mathcal{F}\left(
\tikz[baseline={([yshift=-.5ex]current bounding box.center)}, scale=.75]
{
	\draw[dotted] (-2,-2) circle(0.707);
	\draw[knot] (-1.5,-1.5) .. controls (-1.75,-1.75) and (-1.75,-2.25) .. (-1.5,-2.5);
	\draw[knot] (-2.5,-1.5) .. controls (-2.25,-1.75) and (-2.25,-2.25) ..  (-2.5,-2.5);
} 
  \right) 
  }
\xrightarrow{
\mathcal{F}\left( \tikz[baseline=8.25ex, scale=.45]
{
    \begin{scope}[rotate=90]
	\draw[dotted] (3,-2) circle(0.707);
        \draw[knot, red,thick] (3,-1.7) -- (3,-2.3);
	\draw[knot] (2.5,-1.5) .. controls (2.75,-1.75) and (3.25,-1.75) .. (3.5,-1.5);
	\draw[knot] (2.5,-2.5) .. controls  (2.75,-2.25) and (3.25,-2.25) .. (3.5,-2.5);
    \end{scope}
} \right)}
\mathcal{F}\left(
\tikz[baseline={([yshift=-.5ex]current bounding box.center)}, scale=.75]
{
	\draw[dotted] (3,-2) circle(0.707);
	\draw[knot] (2.5,-1.5) .. controls (2.75,-1.75) and (3.25,-1.75) .. (3.5,-1.5);
	\draw[knot] (2.5,-2.5) .. controls  (2.75,-2.25) and (3.25,-2.25) .. (3.5,-2.5);
}
\right)
\right).
\]
On one hand, it is obvious that
\[
\mathrm{Kh}\left( 
\tikz[baseline={([yshift=-.5ex]current bounding box.center)}, scale=.75]
{
	\draw[dotted] (.5,.5) circle(0.707);
	\draw[knot, ->](0,0) -- (1,1);
	\fill[fill=white] (.5,.5) circle (.15);
	\draw[knot, ->](1,0) -- (0,1);
}
 \right) 
\cong 
\mathcal{F}\left( 
\tikz[baseline={([yshift=-.5ex]current bounding box.center)}, scale=.75]
{
	\draw[dotted] (.5,.5) circle(0.707);
	\draw[knot, ->](0,0) -- (1,1);
	\fill[fill=white] (.5,.5) circle (.15);
	\draw[knot, ->](1,0) -- (0,1);
}
 \right) \{-1,0\}[-1].
\]
On the other hand, the diagram
\[
\tikz[xscale=1]{
    \node(11) at (0,0) {$\mathcal{F}\left(
\tikz[baseline={([yshift=-.5ex]current bounding box.center)}, scale=.75]
{
	\draw[dotted] (-2,-2) circle(0.707);
	\draw[knot] (-1.5,-1.5) .. controls (-1.75,-1.75) and (-1.75,-2.25) .. (-1.5,-2.5);
	\draw[knot] (-2.5,-1.5) .. controls (-2.25,-1.75) and (-2.25,-2.25) ..  (-2.5,-2.5);
} 
  \right)$};
    \node(12) at (8,0) {$\varphi_{\tikz[baseline=8.25ex, scale=.45]
{
    \begin{scope}[rotate=90]
	\draw[dotted] (3,-2) circle(0.707);
        \draw[knot, red,thick] (3,-1.7) -- (3,-2.3);
	\draw[knot] (2.5,-1.5) .. controls (2.75,-1.75) and (3.25,-1.75) .. (3.5,-1.5);
	\draw[knot] (2.5,-2.5) .. controls  (2.75,-2.25) and (3.25,-2.25) .. (3.5,-2.5);
    \end{scope}
}} ^{-1}
\mathcal{F}\left(
\tikz[baseline={([yshift=-.5ex]current bounding box.center)}, scale=.75]
{
	\draw[dotted] (3,-2) circle(0.707);
	\draw[knot] (2.5,-1.5) .. controls (2.75,-1.75) and (3.25,-1.75) .. (3.5,-1.5);
	\draw[knot] (2.5,-2.5) .. controls  (2.75,-2.25) and (3.25,-2.25) .. (3.5,-2.5);
}
\right)$};
    \node(21) at (0,-3) {$\varphi_{\tikz[baseline=8.25ex, scale=.45]
{
    \begin{scope}[rotate=90]
	\draw[dotted] (3,-2) circle(0.707);
        \draw[knot, red,thick] (3,-1.7) -- (3,-2.3);
	\draw[knot] (2.5,-1.5) .. controls (2.75,-1.75) and (3.25,-1.75) .. (3.5,-1.5);
	\draw[knot] (2.5,-2.5) .. controls  (2.75,-2.25) and (3.25,-2.25) .. (3.5,-2.5);
    \end{scope}
}} ^{-1} \circ \varphi_{\tikz[baseline=8.25ex, scale=.45]
{
    \begin{scope}[rotate=90]
	\draw[dotted] (3,-2) circle(0.707);
        \draw[knot, red,thick] (3,-1.7) -- (3,-2.3);
	\draw[knot] (2.5,-1.5) .. controls (2.75,-1.75) and (3.25,-1.75) .. (3.5,-1.5);
	\draw[knot] (2.5,-2.5) .. controls  (2.75,-2.25) and (3.25,-2.25) .. (3.5,-2.5);
    \end{scope}
}} \mathcal{F}\left( \tikz[baseline={([yshift=-.5ex]current bounding box.center)}, scale=.75]
{
	\draw[dotted] (-2,-2) circle(0.707);
	\draw[knot] (-1.5,-1.5) .. controls (-1.75,-1.75) and (-1.75,-2.25) .. (-1.5,-2.5);
	\draw[knot] (-2.5,-1.5) .. controls (-2.25,-1.75) and (-2.25,-2.25) ..  (-2.5,-2.5);
} \right)$};
    \node(22) at (8, -3) {$\varphi_{\tikz[baseline=8.25ex, scale=.45]
{
    \begin{scope}[rotate=90]
	\draw[dotted] (3,-2) circle(0.707);
        \draw[knot, red,thick] (3,-1.7) -- (3,-2.3);
	\draw[knot] (2.5,-1.5) .. controls (2.75,-1.75) and (3.25,-1.75) .. (3.5,-1.5);
	\draw[knot] (2.5,-2.5) .. controls  (2.75,-2.25) and (3.25,-2.25) .. (3.5,-2.5);
    \end{scope}
}} ^{-1}
\mathcal{F}\left(
\tikz[baseline={([yshift=-.5ex]current bounding box.center)}, scale=.75]
{
	\draw[dotted] (3,-2) circle(0.707);
	\draw[knot] (2.5,-1.5) .. controls (2.75,-1.75) and (3.25,-1.75) .. (3.5,-1.5);
	\draw[knot] (2.5,-2.5) .. controls  (2.75,-2.25) and (3.25,-2.25) .. (3.5,-2.5);
}
\right)$};
    \draw[->] (11) to[out=0,in=180] node[pos=.5,above,arrows=-]
        {$\mathcal{F}\left( \tikz[baseline=8.25ex, scale=.45]
{
    \begin{scope}[rotate=90]
	\draw[dotted] (3,-2) circle(0.707);
        \draw[knot, red,thick] (3,-1.7) -- (3,-2.3);
	\draw[knot] (2.5,-1.5) .. controls (2.75,-1.75) and (3.25,-1.75) .. (3.5,-1.5);
	\draw[knot] (2.5,-2.5) .. controls  (2.75,-2.25) and (3.25,-2.25) .. (3.5,-2.5);
    \end{scope}
}\right) \circ \varphi_H$} (12);
    \draw[->] (11) to[out=270,in=90] node[pos=.5,right,arrows=-]
        {$\varphi_H$} (21);
    \draw[->] (21) to[out=0,in=180] node[pos=.5,above,arrows=-]
        {$\mathcal{F}\left( \tikz[baseline=8.25ex, scale=.45]
{
    \begin{scope}[rotate=90]
	\draw[dotted] (3,-2) circle(0.707);
        \draw[knot, red,thick] (3,-1.7) -- (3,-2.3);
	\draw[knot] (2.5,-1.5) .. controls (2.75,-1.75) and (3.25,-1.75) .. (3.5,-1.5);
	\draw[knot] (2.5,-2.5) .. controls  (2.75,-2.25) and (3.25,-2.25) .. (3.5,-2.5);
    \end{scope}
}\right)$} (22);
    \draw[double equal sign distance] (12) to[out=270,in=90] (22);
}
\]
commutes tautologically. The bottom line is exactly $\varphi_{\tikz[baseline=8.25ex, scale=.45]
{
    \begin{scope}[rotate=90]
	\draw[dotted] (3,-2) circle(0.707);
        \draw[knot, red,thick] (3,-1.7) -- (3,-2.3);
	\draw[knot] (2.5,-1.5) .. controls (2.75,-1.75) and (3.25,-1.75) .. (3.5,-1.5);
	\draw[knot] (2.5,-2.5) .. controls  (2.75,-2.25) and (3.25,-2.25) .. (3.5,-2.5);
    \end{scope}
}}^{-1}
\mathcal{F} \left(
\tikz[baseline={([yshift=-.5ex]current bounding box.center)}, scale=.75]
{
	\draw[dotted] (.5,.5) circle(0.707);
	\draw[knot, ->](1,0) -- (0,1);
	\fill[fill=white] (.5,.5) circle (.15);
	\draw[knot, ->](0,0) -- (1,1);
}
 \right)$, so we conclude that 
 \[
\mathrm{Kh} \left(
\tikz[baseline={([yshift=-.5ex]current bounding box.center)}, scale=.75]
{
	\draw[dotted] (.5,.5) circle(0.707);
	\draw[knot, ->](1,0) -- (0,1);
	\fill[fill=white] (.5,.5) circle (.15);
	\draw[knot, ->](0,0) -- (1,1);
}
 \right) \cong
 \varphi_{\tikz[baseline=8.25ex, scale=.45]
{
    \begin{scope}[rotate=90]
	\draw[dotted] (3,-2) circle(0.707);
        \draw[knot, red,thick] (3,-1.7) -- (3,-2.3);
	\draw[knot] (2.5,-1.5) .. controls (2.75,-1.75) and (3.25,-1.75) .. (3.5,-1.5);
	\draw[knot] (2.5,-2.5) .. controls  (2.75,-2.25) and (3.25,-2.25) .. (3.5,-2.5);
    \end{scope}
}}^{-1}
\mathcal{F} \left(
\tikz[baseline={([yshift=-.5ex]current bounding box.center)}, scale=.75]
{
	\draw[dotted] (.5,.5) circle(0.707);
	\draw[knot, ->](1,0) -- (0,1);
	\fill[fill=white] (.5,.5) circle (.15);
	\draw[knot, ->](0,0) -- (1,1);
}
 \right) \{0,1\}.
 \]
Then the desired result follows from the definition of $\mathrm{Kh}$ and Theorem \ref{thm:multigluing}.
\end{proof}

Unfortunately, $\mathrm{Kh}$ is not an invariant of oriented tangles in the $\mathscr{G}$-graded sense; rather, $\mathrm{Kh}$ will be an invariant of diskular tangles up to $\mathscr{G}$-grading shift (Theorem \ref{thm:almosttangleinvt}). We break the computation up into three lemmas of increasing difficulty.

\begin{remark}
Notice that invariance under planar isotopy is immediately apparent in the $\mathscr{G}$-graded setting, in contrast to \cite{naisse2020odd}, since $\mathcal{F}(D_T) \cong \mathcal{F}(D_{T'})$ if $T'$ is obtained from $T$ via planar isotopy. Moreover, in our setup, we no longer have to assume $T$ is presented in a generic position.
\end{remark}

\begin{lemma}
\label{lem:r1}
There are isomorphisms
\[
\mathrm{Kh}
\left(
\tikz[baseline={([yshift=-.5ex]current bounding box.center)}, scale=.75]
{
        \draw[dotted] (.5,.5) circle(0.707);
        \draw[knot] (1, 0.5) to[out=-90, in=-60] (0,1);
        \draw[knot, overcross] (0,0) to[out=60, in =90] (1, 0.5);

}
\right)
\cong
\mathrm{Kh}
\left(
\tikz[baseline={([yshift=-.5ex]current bounding box.center)}, scale=.75]
{
        \draw[dotted] (.5,.5) circle(0.707);
        \draw[knot] (0,0) to[out=60, in=-90] (.7, .5) to[out=90, in=-60] (0,1);
}
\right)
\cong
\mathrm{Kh}
\left(
\tikz[baseline={([yshift=-.5ex]current bounding box.center)}, scale=.75]
{
        \draw[dotted] (.5,.5) circle(0.707);
        \draw[knot] (0,0) to[out=60, in =90] (1, 0.5);
        \draw[knot, overcross] (1, 0.5) to[out=-90, in=-60] (0,1);
}
\right)
\]
in $\mathrm{Kom}(H^1\mathrm{Mod}_R^\mathscr{G})$ (here, the choice of orientation does not matter).
\end{lemma}

\begin{proof}
Picking an orientation for the right handed twist, we compute
\begin{align*}
\mathrm{Kh}
\left(
\tikz[baseline={([yshift=-.5ex]current bounding box.center)}, scale=.75]
{
        \draw[dotted] (.5,.5) circle(0.707);
        \draw[knot, ->] (1, 0.5) to[out=-90, in=-60] (0,1);
        \draw[knot, overcross] (0,0) to[out=60, in =90] (1, 0.5);
}
\right)
&
\cong
\left(
\tikz[baseline={([yshift=-.5ex]current bounding box.center)}]
{
    \node(res0) at (0,0) {$
        \mathcal{F}
        \left(
        \tikz[baseline={([yshift=-.5ex]current bounding box.center)}, scale=.85]
        {
            \draw[dotted] (0.5,0.5) circle(0.707);
            \draw[knot] (0,0) to[out=60, in=-90] (0.2, 0.5) to[out=90, in=-60] (0,1);
            \draw[knot] (0.6, 0.5) circle(0.2);
        }
        \right)$
    };
    \node(res1) at (5,0) {$
        \varphi_{
        \tikz[baseline={([yshift=-.5ex]current bounding box.center)}, scale=.45]
        {
            \draw[dotted] (0.5,0.5) circle(0.707);
            \draw[knot, red] (0.4, 0.5) -- (0.2, 0.5);
            \draw[knot] (0,0) to[out=60, in=-90] (0.2, 0.5) to[out=90, in=-60] (0,1);
            \draw[knot] (0.6, 0.5) circle(0.2);
        }
        }^{-1}
        \mathcal{F}
        \left(
        \tikz[baseline={([yshift=-.5ex]current bounding box.center)}, scale=.85]
        {
            \draw[dotted] (.5,.5) circle(0.707);
            \draw[knot] (0,0) to[out=60, in=-90] (.8, .5) to[out=90, in=-60] (0,1);
        }
        \right)$
    };
    \draw[->] (res0) to[out=0, in=180]
        node[pos=.5, above, arrows=-] {$
        \tikz[baseline={([yshift=-.5ex]current bounding box.center)}, scale=.45]
        {
            \draw[dotted] (0.5,0.5) circle(0.707);
            \draw[knot, red] (0.4, 0.5) -- (0.2, 0.5);
            \draw[knot] (0,0) to[out=60, in=-90] (0.2, 0.5) to[out=90, in=-60] (0,1);
            \draw[knot] (0.6, 0.5) circle(0.2);
        }
        \circ \varphi_H
        $
        }
        (res1);
}
\right) \{0,1\}
\\
&
\cong
\left(
\tikz[baseline={([yshift=-.5ex]current bounding box.center)}]
{
    \node at (0,0) {$\oplus$};
    \node(res0A) at (0,1.25) {$
        \mathcal{F}
        \left(
        \tikz[baseline={([yshift=-.5ex]current bounding box.center)}, scale=.85]
        {
            \draw[dotted] (0.5,0.5) circle(0.707);
            \draw[knot] (0,0) to[out=60, in=-90] (0.2, 0.5) to[out=90, in=-60] (0,1);
        }
        \right)
        \{0,-1\}$
    };
    \node(res0B) at (0,-1.25) {$
        \mathcal{F}
        \left(
        \tikz[baseline={([yshift=-.5ex]current bounding box.center)}, scale=.85]
        {
            \draw[dotted] (0.5,0.5) circle(0.707);
            \draw[knot] (0,0) to[out=60, in=-90] (0.2, 0.5) to[out=90, in=-60] (0,1);
        }
        \right)
        \{1,0\}$
    };
    \node(res1) at (5,0) {$
        \mathcal{F}
        \left(
        \tikz[baseline={([yshift=-.5ex]current bounding box.center)}, scale=.85]
        {
            \draw[dotted] (.5,.5) circle(0.707);
            \draw[knot] (0,0) to[out=60, in=-90] (.8, .5) to[out=90, in=-60] (0,1);
        }
        \right) \{1,0\}$
    };
    \draw[->] (res0A) to
    node[pos=.6, above=1mm, arrows=-] {$
        XZ\varphi_H
        $
        }
    (res1);
    \draw[->] (res0B) to
    node[pos=.6, below=1mm, arrows=-] {$
        \varphi_H
        $
        }
    (res1);
}
\right) \{0,1\}
\\
&
\cong
\mathrm{Kh}
\left(
\tikz[baseline={([yshift=-.5ex]current bounding box.center)}, scale=.75]
{
        \draw[dotted] (.5,.5) circle(0.707);
        \draw[knot] (0,0) to[out=60, in=-90] (.7, .5) to[out=90, in=-60] (0,1);
}
\right) \,.
\end{align*}
The second isomorphism is by delooping, noticing that $\varphi_{
        \tikz[baseline={([yshift=-.5ex]current bounding box.center)}, scale=.45]
        {
            \draw[dotted] (0.5,0.5) circle(0.707);
            \draw[knot, red] (0.4, 0.5) -- (0.2, 0.5);
            \draw[knot] (0,0) to[out=60, in=-90] (0.2, 0.5) to[out=90, in=-60] (0,1);
            \draw[knot] (0.6, 0.5) circle(0.2);
        }
        }^{-1}
$ is isomorphic to $\{1,0\}$ as shifting functors. Additionally, the maps are obtained from the former by precomposing with a birth or a dotted birth. The third isomorphism is by Gaussian elimination. The reader may verify that the computation for $\mathrm{Kh}
\left(
\tikz[baseline={([yshift=-.5ex]current bounding box.center)}, scale=.75]
{
        \draw[dotted] (.5,.5) circle(0.707);
        \draw[knot] (1, 0.5) to[out=-90, in=-60] (0,1);
        \draw[knot, overcross] (0,0) to[out=60, in =90] (1, 0.5);
        \draw[knot, ->] (0.01,0.01) -- (0,0);
}
\right)$ is duplicate.

Doing the same for the left handed twist,
\begin{align*}
\mathrm{Kh}
\left(
\tikz[baseline={([yshift=-.5ex]current bounding box.center)}, scale=.75]
{
        \draw[dotted] (.5,.5) circle(0.707);
        \draw[knot] (0,0) to[out=60, in =90] (1, 0.5);
        \draw[knot, overcross] (1, 0.5) to[out=-90, in=-60] (0,1);
        \draw[knot, ->] (0.01, 0.99) -- (0,1);
}
\right)
&
\cong
\left(
\tikz[baseline={([yshift=-.5ex]current bounding box.center)}]
{
    \node(res0) at (0,0) {$
        \varphi_{
        \tikz[baseline={([yshift=-.5ex]current bounding box.center)}, scale=.45]
        {
            \draw[dotted] (.5,.5) circle(0.707);
            \draw[knot, red] (0.4, .2) -- (0.4, .8);
            \draw[knot] (0,0) to[out=60, in=-90] (.8, .5) to[out=90, in=-60] (0,1);
        }
        }
        \mathcal{F}
        \left(
        \tikz[baseline={([yshift=-.5ex]current bounding box.center)}, scale=.85]
        {
            \draw[dotted] (.5,.5) circle(0.707);
            \draw[knot] (0,0) to[out=60, in=-90] (.8, .5) to[out=90, in=-60] (0,1);
        }
        \right)$
    };
    \node(res1) at (5,0) {$
        \mathcal{F}
        \left(
        \tikz[baseline={([yshift=-.5ex]current bounding box.center)}, scale=.85]
        {
            \draw[dotted] (0.5,0.5) circle(0.707);
            \draw[knot] (0,0) to[out=60, in=-90] (0.2, 0.5) to[out=90, in=-60] (0,1);
            \draw[knot] (0.6, 0.5) circle(0.2);
        }
        \right)$
    };
    \draw[->] (res0) to
    node[pos=0.5, above, arrows=-] {$
    \tikz[baseline={([yshift=-.5ex]current bounding box.center)}, scale=.45]
        {
            \draw[dotted] (.5,.5) circle(0.707);
            \draw[knot, red] (0.4, .2) -- (0.4, .8);
            \draw[knot] (0,0) to[out=60, in=-90] (.8, .5) to[out=90, in=-60] (0,1);
        }$
        }
    (res1);
}
\right) \{-1,0\}
\\
&
\cong
\left(
\tikz[baseline={([yshift=-.5ex]current bounding box.center)}]
{
    \node at (5,0) {$\oplus$};
    \node(res0) at (0,0) {$
        \mathcal{F}
        \left(
        \tikz[baseline={([yshift=-.5ex]current bounding box.center)}, scale=.85]
        {
            \draw[dotted] (.5,.5) circle(0.707);
            \draw[knot] (0,0) to[out=60, in=-90] (.8, .5) to[out=90, in=-60] (0,1);
        }
        \right)
        \{0, -1\}$
    };
    \node(res1A) at (5,1.25) {$
        \mathcal{F}
        \left(
        \tikz[baseline={([yshift=-.5ex]current bounding box.center)}, scale=.85]
        {
            \draw[dotted] (0.5,0.5) circle(0.707);
            \draw[knot] (0,0) to[out=60, in=-90] (0.2, 0.5) to[out=90, in=-60] (0,1);
        }
        \right)
        \{0,-1\}$
    };
    \node(res1B) at (5,-1.25) {$
        \mathcal{F}
        \left(
        \tikz[baseline={([yshift=-.5ex]current bounding box.center)}, scale=.85]
        {
            \draw[dotted] (0.5,0.5) circle(0.707);
            \draw[knot] (0,0) to[out=60, in=-90] (0.2, 0.5) to[out=90, in=-60] (0,1);
        }
        \right)
        \{1,0\}$
    };
    \draw[->] (res0) to
    node[pos=0.5, above=1mm, arrows=-] {$1$}
        (res1A);
    \draw[->] (res0) to
    node[pos=0.5, below=1mm, arrows=-] {$YZ$}
    (res1B);
}
\right) \{-1,0\}
\\
&
\cong
\mathrm{Kh}
\left(
\tikz[baseline={([yshift=-.5ex]current bounding box.center)}, scale=.75]
{
        \draw[dotted] (.5,.5) circle(0.707);
        \draw[knot] (0,0) to[out=60, in=-90] (.7, .5) to[out=90, in=-60] (0,1);
}
\right)
\end{align*}
follows by the same reasoning, and the computation for $\mathrm{Kh}
\left(
\tikz[baseline={([yshift=-.5ex]current bounding box.center)}, scale=.75]
{
        \draw[dotted] (.5,.5) circle(0.707);
        \draw[knot, <-] (0,0) to[out=60, in =90] (1, 0.5);
        \draw[knot, overcross] (1, 0.5) to[out=-90, in=-60] (0,1);
}
\right)$ is its doppelg\"anger.
\end{proof}

\begin{lemma}
\label{lem:r2}
There are isomoprhisms
\[
\mathrm{Kh}
\left(
\tikz[baseline={([yshift=-.5ex]current bounding box.center)}, scale=.85]
{
        \draw[dotted] (.5,.5) circle(0.707);
        \draw[knot, ->] (0,0) to[out=30, in=-90] (0.8, 0.5);
        \draw[knot] (0.8, 0.5) to[out=90, in=-30] (0,1);
        \draw[knot, overcross] (1,0) to[out=150, in=-90] (0.3, 0.5);
        \draw[knot, overcross] (0.3, 0.5) to[out=90, in=210] (1,1);
        \draw[knot, ->] (0.3, 0.499) -- (0.3, 0.501);
}
\right)
\cong
\mathrm{Kh}
\left(
\tikz[baseline={([yshift=-.5ex]current bounding box.center)}, scale=.85]
{
        \draw[dotted] (.5,.5) circle(0.707);
        \draw[knot] (0, 0) to[out=30, in=-90] (0.35, 0.5);
        \draw[knot] (0.35, 0.5) to[out=90, in=-30] (0, 1);
        \draw[knot] (1, 0) to[out=150, in=-90] (0.65, 0.5);
        \draw[knot] (0.65, 0.5) to[out=90, in=210] (1, 1);
}
\right)
\{-1,1\}
\cong
\mathrm{Kh}
\left(
\tikz[baseline={([yshift=-.5ex]current bounding box.center)}, scale=.85]
{
        \draw[dotted] (.5,.5) circle(0.707);
        \draw[knot, ->] (1,0) to[out=150, in=-90] (0.3, 0.5);
        \draw[knot] (0.3, 0.5) to[out=90, in=210] (1,1);
        \draw[knot, overcross] (0,0) to[out=30, in=-90] (0.8, 0.5);
        \draw[knot, overcross] (0.8, 0.5) to[out=90, in=-30] (0,1);
        \draw[knot, ->] (0.8, 0.499) -- (0.8, 0.501);
}
\right)
\]
and
\[
\mathrm{Kh}
\left(
\tikz[baseline={([yshift=-.5ex]current bounding box.center)}, scale=.85]
{
        \draw[dotted] (.5,.5) circle(0.707);
        \draw[knot, ->] (0,0) to[out=30, in=-90] (0.8, 0.5);
        \draw[knot] (0.8, 0.5) to[out=90, in=-30] (0,1);
        \draw[knot, overcross] (1,0) to[out=150, in=-90] (0.3, 0.5);
        \draw[knot, overcross] (0.3, 0.5) to[out=90, in=210] (1,1);
        \draw[knot, <-] (0.3, 0.499) -- (0.3, 0.501);
}
\right)
\cong
\varphi_{\left(\tikz[baseline=8.25ex, scale=.45]
{
    \begin{scope}[rotate=90]
	\draw[dotted] (3,-2) circle(0.707);
	\draw[knot] (2.5,-1.5) .. controls (2.75,-1.75) and (3.25,-1.75) .. (3.5,-1.5);
	\draw[knot] (2.5,-2.5) .. controls  (2.75,-2.25) and (3.25,-2.25) .. (3.5,-2.5);
        \draw[red, knot] (3,-1.7) -- (3,-2.3);
    \end{scope}
},\, (0,1)\right)}
\mathrm{Kh}
\left(
\tikz[baseline={([yshift=-.5ex]current bounding box.center)}, scale=.85]
{
        \draw[dotted] (.5,.5) circle(0.707);
        \draw[knot] (0, 0) to[out=30, in=-90] (0.35, 0.5);
        \draw[knot] (0.35, 0.5) to[out=90, in=-30] (0, 1);
        \draw[knot] (1, 0) to[out=150, in=-90] (0.65, 0.5);
        \draw[knot] (0.65, 0.5) to[out=90, in=210] (1, 1);
}
\right)
\cong
\mathrm{Kh}
\left(
\tikz[baseline={([yshift=-.5ex]current bounding box.center)}, scale=.85]
{
        \draw[dotted] (.5,.5) circle(0.707);
        \draw[knot] (1,0) to[out=150, in=-90] (0.3, 0.5);
        \draw[knot, <-] (0.3, 0.5) to[out=90, in=210] (1,1);
        \draw[knot, overcross] (0,0) to[out=30, in=-90] (0.8, 0.5);
        \draw[knot, overcross] (0.8, 0.5) to[out=90, in=-30] (0,1);
        \draw[knot, ->] (0.8, 0.499) -- (0.8, 0.501);
}
\right)
\]
in $\mathrm{Kom}(H^2\mathrm{Mod}_R^\mathscr{G})$. We call the first pair of isomorphisms $\mathrm{RII}_+$ and the second pair $\mathrm{RII}_-$.
\end{lemma}

\begin{proof}
By definition, $\mathrm{Kh}
\left(
\tikz[baseline={([yshift=-.5ex]current bounding box.center)}, scale=.85]
{
        \draw[dotted] (.5,.5) circle(0.707);
        \draw[knot, ->] (0,0) to[out=30, in=-90] (0.8, 0.5);
        \draw[knot] (0.8, 0.5) to[out=90, in=-30] (0,1);
        \draw[knot, overcross] (1,0) to[out=150, in=-90] (0.3, 0.5);
        \draw[knot, overcross] (0.3, 0.5) to[out=90, in=210] (1,1);
        \draw[knot, ->] (0.3, 0.499) -- (0.3, 0.501);
}
\right)$ is the complex
\[
\tikz[baseline={([yshift=-.5ex]current bounding box.center)}]
{
    \node at (5,0) {$\oplus$};
    \node(res00) at (0,1.25) {$
    \varphi_{
    \tikz[baseline={([yshift=-.5ex]current bounding box.center)}, scale=.45]
        {
        \draw[dotted] (0.5, 0.5) circle (0.707);
    %
        \draw[knot, red] (0.5, 0.4) -- (0.5, 0.2);
        \draw[knot] (0,0) to[out=45, in=135] (1,0);
        \draw[knot] (0,1) to[out=-45, in=90] (0.35, 0.75);
        \draw[knot] (0.35, 0.75) to[out=-90, in=90] (0.2, 0.5);
        \draw[knot] (0.2, 0.5) to[out=-90, in=180] (0.5, 0.4);
        \draw[knot] (0.5, 0.4) to[out=0, in=-90] (0.8, 0.5);
        \draw[knot] (0.8, 0.5) to[out=90, in=-90] (0.65, 0.75);
        \draw[knot] (0.65, 0.75) to[out=90, in=225] (1,1);
        }
        }\,
        \mathcal{F}
        \left(
        \tikz[baseline={([yshift=-.5ex]current bounding box.center)}, scale=.85]
        {
        \draw[dotted] (0.5, 0.5) circle (0.707);
    %
    %
        \draw[knot] (0,0) to[out=45, in=135] (1,0);
        \draw[knot] (0,1) to[out=-45, in=90] (0.35, 0.75);
        \draw[knot] (0.35, 0.75) to[out=-90, in=90] (0.2, 0.5);
        \draw[knot] (0.2, 0.5) to[out=-90, in=180] (0.5, 0.4);
        \draw[knot] (0.5, 0.4) to[out=0, in=-90] (0.8, 0.5);
        \draw[knot] (0.8, 0.5) to[out=90, in=-90] (0.65, 0.75);
        \draw[knot] (0.65, 0.75) to[out=90, in=225] (1,1);
        }
        \right)
    $};
    \node(res10) at (5,1.25) {$
 \mathcal{F}
\left(
\tikz[baseline={([yshift=-.5ex]current bounding box.center)}, scale=.85]
{
        \draw[dotted] (0.5, 0.5) circle (0.707);

        \draw[knot] (0,1) to[out=-45, in=90] (0.35, 0.75);
        \draw[knot] (0.35, 0.75) to[out=-90, in=90] (0.2, 0.5);
        \draw[knot] (0.2, 0.5) to[out=-90, in=90] (0.35, 0.25);
        \draw[knot] (0.35, 0.25) to[out=-90, in=45] (0,0);
        \draw[knot] (1,1) to[out=225, in=90] (0.65, 0.75);
        \draw[knot] (0.65, 0.75) to[out=-90, in=90] (0.8, 0.5);
        \draw[knot] (0.8, 0.5) to[out=-90, in=90] (0.65, 0.25);
        \draw[knot] (0.65, 0.25) to[out=-90, in=135] (1,0);
}
\right)
    $};
    \node(res01) at (5,-1.25) {$
    \varphi_{
        \tikz[baseline={([yshift=-.5ex]current bounding box.center)}, scale=.45]
        {
        \draw[dotted] (0.5, 0.5) circle (0.707);
        \draw[knot, red] (0.5, 0.4) -- (0.5, 0.2);
        \draw[knot] (0,0) to[out=45, in=135] (1,0);
        \draw[knot] (0,1) to[out=-45, in=90] (0.35, 0.75);
        \draw[knot] (0.35, 0.75) to[out=-90, in=90] (0.2, 0.5);
        \draw[knot] (0.2, 0.5) to[out=-90, in=180] (0.5, 0.4);
        \draw[knot] (0.5, 0.4) to[out=0, in=-90] (0.8, 0.5);
        \draw[knot] (0.8, 0.5) to[out=90, in=-90] (0.65, 0.75);
        \draw[knot] (0.65, 0.75) to[out=90, in=225] (1,1);
        }
    }
    \circ
    \varphi_{
        \tikz[baseline={([yshift=-.5ex]current bounding box.center)}, scale=.45]
        {
        \draw[dotted] (0.5, 0.5) circle (0.707);
        \draw[knot, red] (0.35, 0.75) -- (0.65, 0.75);
        \draw[knot] (0,0) to[out=45, in=135] (1,0);
        \draw[knot] (0,1) to[out=-45, in=90] (0.35, 0.75);
        \draw[knot] (0.35, 0.75) to[out=-90, in=90] (0.2, 0.5);
        \draw[knot] (0.2, 0.5) to[out=-90, in=180] (0.5, 0.4);
        \draw[knot] (0.5, 0.4) to[out=0, in=-90] (0.8, 0.5);
        \draw[knot] (0.8, 0.5) to[out=90, in=-90] (0.65, 0.75);
        \draw[knot] (0.65, 0.75) to[out=90, in=225] (1,1);
        }
    }^{-1}
    \,
    \mathcal{F}
    \left(
    \tikz[baseline={([yshift=-.5ex]current bounding box.center)}, scale=.85]
        {
        \draw[dotted] (0.5, 0.5) circle (0.707);

        \draw[knot] (0,0) to[out=45, in=135] (1,0);
        \draw[knot] (0,1) to[out=-45, in=-135] (1,1);
        \draw[knot] (0.5, 0.65) to[out=0, in=90] (0.8, 0.5) to[out=-90, in=0] (0.5, 0.35) to[out=180, in=-90] (0.2, 0.5) to[out=90, in=180] (0.5, 0.65);
        }
    \right)
    $};
    \node(res11) at (11, -1.25){$
    \varphi_{
    \tikz[baseline={([yshift=-.5ex]current bounding box.center)}, scale=.45]
        {
        \draw[dotted] (0.5, 0.5) circle (0.707);
        \draw[knot, red] (0.35, 0.75) -- (0.65, 0.75);
        \draw[knot] (0,1) to[out=-45, in=90] (0.35, 0.75);
        \draw[knot] (0.35, 0.75) to[out=-90, in=90] (0.2, 0.5);
        \draw[knot] (0.2, 0.5) to[out=-90, in=90] (0.35, 0.25);
        \draw[knot] (0.35, 0.25) to[out=-90, in=45] (0,0);
        \draw[knot] (1,1) to[out=225, in=90] (0.65, 0.75);
        \draw[knot] (0.65, 0.75) to[out=-90, in=90] (0.8, 0.5);
        \draw[knot] (0.8, 0.5) to[out=-90, in=90] (0.65, 0.25);
        \draw[knot] (0.65, 0.25) to[out=-90, in=135] (1,0);
        }
        }^{-1}\,
        \mathcal{F}
        \left(
        \tikz[baseline={([yshift=-.5ex]current bounding box.center)}, scale=.85]
        {
        \draw[dotted] (0.5, 0.5) circle (0.707);
        \draw[knot] (0,1) to[out=-45, in=-135] (1,1);
        \draw[knot] (0,0) to[out=45, in=-90] (0.35, 0.25);
        \draw[knot] (0.35, 0.25) to[out=90, in=-90] (0.2, 0.5);
        \draw[knot] (0.2, 0.5) to[out=90, in=180] (0.5, 0.6);
        \draw[knot] (0.5, 0.6) to[out=0, in=90] (0.8, 0.5);
        \draw[knot] (0.8, 0.5) to[out=-90, in=90] (0.65, 0.25);
        \draw[knot] (0.65, 0.25) to[out=-90, in=135] (1,0);
        }
        \right)
    $};
    \draw[->] (res00) to
        node[pos=.5,above,arrows=-] {$\tikz[baseline={([yshift=-.5ex]current bounding box.center)}, scale=.45]
        {
        \draw[dotted] (0.5, 0.5) circle (0.707);
    %
        \draw[knot, red] (0.5, 0.4) -- (0.5, 0.2);
        \draw[knot] (0,0) to[out=45, in=135] (1,0);
        \draw[knot] (0,1) to[out=-45, in=90] (0.35, 0.75);
        \draw[knot] (0.35, 0.75) to[out=-90, in=90] (0.2, 0.5);
        \draw[knot] (0.2, 0.5) to[out=-90, in=180] (0.5, 0.4);
        \draw[knot] (0.5, 0.4) to[out=0, in=-90] (0.8, 0.5);
        \draw[knot] (0.8, 0.5) to[out=90, in=-90] (0.65, 0.75);
        \draw[knot] (0.65, 0.75) to[out=90, in=225] (1,1);
        }
        $}
        (res10);
    \draw[->] (res00) to 
        node[pos=.2, below=1.5mm, arrows=-] {$
        \tikz[baseline={([yshift=-.5ex]current bounding box.center)}, scale=.45]
        {
        \draw[dotted] (0.5, 0.5) circle (0.707);
        \draw[knot, red] (0.35, 0.75) -- (0.65, 0.75);
    %
        \draw[knot] (0,0) to[out=45, in=135] (1,0);
        \draw[knot] (0,1) to[out=-45, in=90] (0.35, 0.75);
        \draw[knot] (0.35, 0.75) to[out=-90, in=90] (0.2, 0.5);
        \draw[knot] (0.2, 0.5) to[out=-90, in=180] (0.5, 0.4);
        \draw[knot] (0.5, 0.4) to[out=0, in=-90] (0.8, 0.5);
        \draw[knot] (0.8, 0.5) to[out=90, in=-90] (0.65, 0.75);
        \draw[knot] (0.65, 0.75) to[out=90, in=225] (1,1);
        }
        \circ \varphi_{H_1}
        $}
        (res01);
    \draw[->] (res10) to
        node[pos=.6, above=1.5mm, arrows=-] {$
        \tikz[baseline={([yshift=-.5ex]current bounding box.center)}, scale=.45]
        {
        \draw[dotted] (0.5, 0.5) circle (0.707);
        \draw[knot, red] (0.35, 0.75) -- (0.65, 0.75);
    %
        \draw[knot] (0,1) to[out=-45, in=90] (0.35, 0.75);
        \draw[knot] (0.35, 0.75) to[out=-90, in=90] (0.2, 0.5);
        \draw[knot] (0.2, 0.5) to[out=-90, in=90] (0.35, 0.25);
        \draw[knot] (0.35, 0.25) to[out=-90, in=45] (0,0);
        \draw[knot] (1,1) to[out=225, in=90] (0.65, 0.75);
        \draw[knot] (0.65, 0.75) to[out=-90, in=90] (0.8, 0.5);
        \draw[knot] (0.8, 0.5) to[out=-90, in=90] (0.65, 0.25);
        \draw[knot] (0.65, 0.25) to[out=-90, in=135] (1,0);
        }
        \circ \varphi_{H_2}
        $}
        (res11);
    \draw[->] (res01) to
        node[pos=.5, above, arrows=-] {$
        \tikz[baseline={([yshift=-.5ex]current bounding box.center)}, scale=.45]
        {
        \draw[dotted] (0.5, 0.5) circle (0.707);
    %
        \draw[knot, red] (0.5, 0.35) -- (0.5, 0.2);
        \draw[knot] (0,0) to[out=45, in=135] (1,0);
        \draw[knot] (0,1) to[out=-45, in=-135] (1,1);
        \draw[knot] (0.5, 0.65) to[out=0, in=90] (0.8, 0.5) to[out=-90, in=0] (0.5, 0.35) to[out=180, in=-90] (0.2, 0.5) to[out=90, in=180] (0.5, 0.65);
        }
        $}
        (res11);
}
\]
with a global shift by $\{-1,1\}$. However, up to isomorphism, we can rewrite the grading shifts on the 01 and 11 resolutions suggestively, so that the complex takes the form
\[
\tikz[baseline={([yshift=-.5ex]current bounding box.center)}]
{
    \node at (5,0) {$\oplus$};
    \node(res00) at (0,1.25) {$
    \varphi_{
    \tikz[baseline={([yshift=-.5ex]current bounding box.center)}, scale=.45]
        {
        \draw[dotted] (0.5, 0.5) circle (0.707);
    %
        \draw[knot, red] (0.5, 0.4) -- (0.5, 0.2);
        \draw[knot] (0,0) to[out=45, in=135] (1,0);
        \draw[knot] (0,1) to[out=-45, in=90] (0.35, 0.75);
        \draw[knot] (0.35, 0.75) to[out=-90, in=90] (0.2, 0.5);
        \draw[knot] (0.2, 0.5) to[out=-90, in=180] (0.5, 0.4);
        \draw[knot] (0.5, 0.4) to[out=0, in=-90] (0.8, 0.5);
        \draw[knot] (0.8, 0.5) to[out=90, in=-90] (0.65, 0.75);
        \draw[knot] (0.65, 0.75) to[out=90, in=225] (1,1);
        }
        }\,
        \mathcal{F}
        \left(
        \tikz[baseline={([yshift=-.5ex]current bounding box.center)}, scale=.85]
        {
        \draw[dotted] (0.5, 0.5) circle (0.707);
    %
    %
        \draw[knot] (0,0) to[out=45, in=135] (1,0);
        \draw[knot] (0,1) to[out=-45, in=90] (0.35, 0.75);
        \draw[knot] (0.35, 0.75) to[out=-90, in=90] (0.2, 0.5);
        \draw[knot] (0.2, 0.5) to[out=-90, in=180] (0.5, 0.4);
        \draw[knot] (0.5, 0.4) to[out=0, in=-90] (0.8, 0.5);
        \draw[knot] (0.8, 0.5) to[out=90, in=-90] (0.65, 0.75);
        \draw[knot] (0.65, 0.75) to[out=90, in=225] (1,1);
        }
        \right)
    $};
    \node(res10) at (5,1.25) {$
 \mathcal{F}
\left(
\tikz[baseline={([yshift=-.5ex]current bounding box.center)}, scale=.85]
{
        \draw[dotted] (0.5, 0.5) circle (0.707);

        \draw[knot] (0,1) to[out=-45, in=90] (0.35, 0.75);
        \draw[knot] (0.35, 0.75) to[out=-90, in=90] (0.2, 0.5);
        \draw[knot] (0.2, 0.5) to[out=-90, in=90] (0.35, 0.25);
        \draw[knot] (0.35, 0.25) to[out=-90, in=45] (0,0);
        \draw[knot] (1,1) to[out=225, in=90] (0.65, 0.75);
        \draw[knot] (0.65, 0.75) to[out=-90, in=90] (0.8, 0.5);
        \draw[knot] (0.8, 0.5) to[out=-90, in=90] (0.65, 0.25);
        \draw[knot] (0.65, 0.25) to[out=-90, in=135] (1,0);
}
\right)
    $};
    \node(res01) at (5,-1.25) {$
    \varphi_{
        \left(
        \tikz[baseline={([yshift=-.5ex]current bounding box.center)}, scale=.45]
        {
        \draw[dotted] (0.5, 0.5) circle (0.707);
        \draw[knot, red] (0.5, 0.4) -- (0.5, 0.2);
        \draw[knot] (0,0) to[out=45, in=135] (1,0);
        \draw[knot] (0,1) to[out=-45, in=90] (0.35, 0.75);
        \draw[knot] (0.35, 0.75) to[out=-90, in=90] (0.2, 0.5);
        \draw[knot] (0.2, 0.5) to[out=-90, in=180] (0.5, 0.4);
        \draw[knot] (0.5, 0.4) to[out=0, in=-90] (0.8, 0.5);
        \draw[knot] (0.8, 0.5) to[out=90, in=-90] (0.65, 0.75);
        \draw[knot] (0.65, 0.75) to[out=90, in=225] (1,1);
        }
        , \,
        (0,1)
        \right)
    }
    \,
    \mathcal{F}
    \left(
    \tikz[baseline={([yshift=-.5ex]current bounding box.center)}, scale=.85]
        {
        \draw[dotted] (0.5, 0.5) circle (0.707);

        \draw[knot] (0,0) to[out=45, in=135] (1,0);
        \draw[knot] (0,1) to[out=-45, in=-135] (1,1);
        \draw[knot] (0.5, 0.65) to[out=0, in=90] (0.8, 0.5) to[out=-90, in=0] (0.5, 0.35) to[out=180, in=-90] (0.2, 0.5) to[out=90, in=180] (0.5, 0.65);
        }
    \right)
    $};
    \node(res11) at (11, -1.25){$
    \varphi_{
    \left(
    \tikz[baseline={([yshift=-.5ex]current bounding box.center)}, scale=.45]
        {
        \draw[dotted] (0.5, 0.5) circle (0.707);
    %
        \draw[knot, red] (0.5, 0.6) -- (0.5, 0.8);
        \draw[knot] (0,1) to[out=-45, in=-135] (1,1);
        \draw[knot] (0,0) to[out=45, in=-90] (0.35, 0.25);
        \draw[knot] (0.35, 0.25) to[out=90, in=-90] (0.2, 0.5);
        \draw[knot] (0.2, 0.5) to[out=90, in=180] (0.5, 0.6);
        \draw[knot] (0.5, 0.6) to[out=0, in=90] (0.8, 0.5);
        \draw[knot] (0.8, 0.5) to[out=-90, in=90] (0.65, 0.25);
        \draw[knot] (0.65, 0.25) to[out=-90, in=135] (1,0);        
        }
        , \,
        (1,1)
        \right)
        }\,
        \mathcal{F}
        \left(
        \tikz[baseline={([yshift=-.5ex]current bounding box.center)}, scale=.85]
        {
        \draw[dotted] (0.5, 0.5) circle (0.707);
        \draw[knot] (0,1) to[out=-45, in=-135] (1,1);
        \draw[knot] (0,0) to[out=45, in=-90] (0.35, 0.25);
        \draw[knot] (0.35, 0.25) to[out=90, in=-90] (0.2, 0.5);
        \draw[knot] (0.2, 0.5) to[out=90, in=180] (0.5, 0.6);
        \draw[knot] (0.5, 0.6) to[out=0, in=90] (0.8, 0.5);
        \draw[knot] (0.8, 0.5) to[out=-90, in=90] (0.65, 0.25);
        \draw[knot] (0.65, 0.25) to[out=-90, in=135] (1,0);
        }
        \right)
    $};
    \draw[->] (res00) to
        node[pos=.5,above,arrows=-] {$\tikz[baseline={([yshift=-.5ex]current bounding box.center)}, scale=.45]
        {
        \draw[dotted] (0.5, 0.5) circle (0.707);
    %
        \draw[knot, red] (0.5, 0.4) -- (0.5, 0.2);
        \draw[knot] (0,0) to[out=45, in=135] (1,0);
        \draw[knot] (0,1) to[out=-45, in=90] (0.35, 0.75);
        \draw[knot] (0.35, 0.75) to[out=-90, in=90] (0.2, 0.5);
        \draw[knot] (0.2, 0.5) to[out=-90, in=180] (0.5, 0.4);
        \draw[knot] (0.5, 0.4) to[out=0, in=-90] (0.8, 0.5);
        \draw[knot] (0.8, 0.5) to[out=90, in=-90] (0.65, 0.75);
        \draw[knot] (0.65, 0.75) to[out=90, in=225] (1,1);
        }
        $}
        (res10);
    \draw[->] (res00) to 
        node[pos=.2, below=1.5mm, arrows=-] {$
        \tikz[baseline={([yshift=-.5ex]current bounding box.center)}, scale=.45]
        {
        \draw[dotted] (0.5, 0.5) circle (0.707);
        \draw[knot, red] (0.35, 0.75) -- (0.65, 0.75);
    %
        \draw[knot] (0,0) to[out=45, in=135] (1,0);
        \draw[knot] (0,1) to[out=-45, in=90] (0.35, 0.75);
        \draw[knot] (0.35, 0.75) to[out=-90, in=90] (0.2, 0.5);
        \draw[knot] (0.2, 0.5) to[out=-90, in=180] (0.5, 0.4);
        \draw[knot] (0.5, 0.4) to[out=0, in=-90] (0.8, 0.5);
        \draw[knot] (0.8, 0.5) to[out=90, in=-90] (0.65, 0.75);
        \draw[knot] (0.65, 0.75) to[out=90, in=225] (1,1);
        }
        \circ \varphi_{H_1}
        $}
        (res01);
    \draw[->] (res10) to
        node[pos=.6, above=1.5mm, arrows=-] {$
        \tikz[baseline={([yshift=-.5ex]current bounding box.center)}, scale=.45]
        {
        \draw[dotted] (0.5, 0.5) circle (0.707);
        \draw[knot, red] (0.35, 0.75) -- (0.65, 0.75);
    %
        \draw[knot] (0,1) to[out=-45, in=90] (0.35, 0.75);
        \draw[knot] (0.35, 0.75) to[out=-90, in=90] (0.2, 0.5);
        \draw[knot] (0.2, 0.5) to[out=-90, in=90] (0.35, 0.25);
        \draw[knot] (0.35, 0.25) to[out=-90, in=45] (0,0);
        \draw[knot] (1,1) to[out=225, in=90] (0.65, 0.75);
        \draw[knot] (0.65, 0.75) to[out=-90, in=90] (0.8, 0.5);
        \draw[knot] (0.8, 0.5) to[out=-90, in=90] (0.65, 0.25);
        \draw[knot] (0.65, 0.25) to[out=-90, in=135] (1,0);
        }
        \circ \varphi_{H_2}
        $}
        (res11);
    \draw[->] (res01) to
        node[pos=.5, above, arrows=-] {$
        \tikz[baseline={([yshift=-.5ex]current bounding box.center)}, scale=.45]
        {
        \draw[dotted] (0.5, 0.5) circle (0.707);
    %
        \draw[knot, red] (0.5, 0.35) -- (0.5, 0.2);
        \draw[knot] (0,0) to[out=45, in=135] (1,0);
        \draw[knot] (0,1) to[out=-45, in=-135] (1,1);
        \draw[knot] (0.5, 0.65) to[out=0, in=90] (0.8, 0.5) to[out=-90, in=0] (0.5, 0.35) to[out=180, in=-90] (0.2, 0.5) to[out=90, in=180] (0.5, 0.65);
        }
        $}
        (res11);
}
\]
again, with a global shift by $\{-1,1\}$. Now, by delooping, 
\[
\varphi_{
        \left(
        \tikz[baseline={([yshift=-.5ex]current bounding box.center)}, scale=.45]
        {
        \draw[dotted] (0.5, 0.5) circle (0.707);
        \draw[knot, red] (0.5, 0.4) -- (0.5, 0.2);
        \draw[knot] (0,0) to[out=45, in=135] (1,0);
        \draw[knot] (0,1) to[out=-45, in=90] (0.35, 0.75);
        \draw[knot] (0.35, 0.75) to[out=-90, in=90] (0.2, 0.5);
        \draw[knot] (0.2, 0.5) to[out=-90, in=180] (0.5, 0.4);
        \draw[knot] (0.5, 0.4) to[out=0, in=-90] (0.8, 0.5);
        \draw[knot] (0.8, 0.5) to[out=90, in=-90] (0.65, 0.75);
        \draw[knot] (0.65, 0.75) to[out=90, in=225] (1,1);
        }
        , \,
        (0,1)
        \right)
    }
    \,
    \mathcal{F}
    \left(
    \tikz[baseline={([yshift=-.5ex]current bounding box.center)}, scale=.85]
        {
        \draw[dotted] (0.5, 0.5) circle (0.707);

        \draw[knot] (0,0) to[out=45, in=135] (1,0);
        \draw[knot] (0,1) to[out=-45, in=-135] (1,1);
        \draw[knot] (0.5, 0.65) to[out=0, in=90] (0.8, 0.5) to[out=-90, in=0] (0.5, 0.35) to[out=180, in=-90] (0.2, 0.5) to[out=90, in=180] (0.5, 0.65);
        }
    \right)
    \cong
\varphi_{
        \tikz[baseline={([yshift=-.5ex]current bounding box.center)}, scale=.45]
        {
        \draw[dotted] (0.5, 0.5) circle (0.707);
        \draw[knot,red] (0.5, 0.2) -- (0.5, 0.8);
        \draw[knot] (0,0) to[out=45, in=135] (1,0);
        \draw[knot] (0,1) to[out=-45, in=-135] (1,1);
        }
    }
    \,
    \mathcal{F}
    \left(
    \tikz[baseline={([yshift=-.5ex]current bounding box.center)}, scale=.85]
        {
        \draw[dotted] (0.5, 0.5) circle (0.707);
        \draw[knot] (0,0) to[out=45, in=135] (1,0);
        \draw[knot] (0,1) to[out=-45, in=-135] (1,1);
        }
    \right)
    \oplus
    \varphi_{
        \left(
        \tikz[baseline={([yshift=-.5ex]current bounding box.center)}, scale=.45]
        {
        \draw[dotted] (0.5, 0.5) circle (0.707);
        \draw[knot,red] (0.5, 0.2) -- (0.5, 0.8);
        \draw[knot] (0,0) to[out=45, in=135] (1,0);
        \draw[knot] (0,1) to[out=-45, in=-135] (1,1);
        }
        , \,
        (1,1)
        \right)
    }
    \,
    \mathcal{F}
    \left(
    \tikz[baseline={([yshift=-.5ex]current bounding box.center)}, scale=.85]
        {
        \draw[dotted] (0.5, 0.5) circle (0.707);
        \draw[knot] (0,0) to[out=45, in=135] (1,0);
        \draw[knot] (0,1) to[out=-45, in=-135] (1,1);
        }
    \right).
\]
Moreover, the maps 
$\mathcal{F}\left(
        \tikz[baseline={([yshift=-.5ex]current bounding box.center)}, scale=.45]
        {
        \draw[dotted] (0.5, 0.5) circle (0.707);
        \draw[knot, red] (0.35, 0.75) -- (0.65, 0.75);
    %
        \draw[knot] (0,0) to[out=45, in=135] (1,0);
        \draw[knot] (0,1) to[out=-45, in=90] (0.35, 0.75);
        \draw[knot] (0.35, 0.75) to[out=-90, in=90] (0.2, 0.5);
        \draw[knot] (0.2, 0.5) to[out=-90, in=180] (0.5, 0.4);
        \draw[knot] (0.5, 0.4) to[out=0, in=-90] (0.8, 0.5);
        \draw[knot] (0.8, 0.5) to[out=90, in=-90] (0.65, 0.75);
        \draw[knot] (0.65, 0.75) to[out=90, in=225] (1,1);
        }
\right) \circ \varphi_{H_1}$
and 
$\mathcal{F}\left(
         \tikz[baseline={([yshift=-.5ex]current bounding box.center)}, scale=.45]
        {
        \draw[dotted] (0.5, 0.5) circle (0.707);
    %
        \draw[knot, red] (0.5, 0.35) -- (0.5, 0.2);
        \draw[knot] (0,0) to[out=45, in=135] (1,0);
        \draw[knot] (0,1) to[out=-45, in=-135] (1,1);
        \draw[knot] (0.5, 0.65) to[out=0, in=90] (0.8, 0.5) to[out=-90, in=0] (0.5, 0.35) to[out=180, in=-90] (0.2, 0.5) to[out=90, in=180] (0.5, 0.65);
        }
\right)$
compose with the delooping isomorphism to yield invertible maps where desired, so that Gaussian elimination tells us that the entire complex is homotopy equivalent to $ \mathrm{Kh}
\left(
\tikz[baseline={([yshift=-.5ex]current bounding box.center)}, scale=.85]
{
        \draw[dotted] (0.5, 0.5) circle (0.707);

        \draw[knot] (0,1) to[out=-45, in=90] (0.35, 0.75);
        \draw[knot] (0.35, 0.75) to[out=-90, in=90] (0.2, 0.5);
        \draw[knot] (0.2, 0.5) to[out=-90, in=90] (0.35, 0.25);
        \draw[knot] (0.35, 0.25) to[out=-90, in=45] (0,0);
        \draw[knot] (1,1) to[out=225, in=90] (0.65, 0.75);
        \draw[knot] (0.65, 0.75) to[out=-90, in=90] (0.8, 0.5);
        \draw[knot] (0.8, 0.5) to[out=-90, in=90] (0.65, 0.25);
        \draw[knot] (0.65, 0.25) to[out=-90, in=135] (1,0);
}\right)\{-1,1\}
$, 
as desired. Duplicate this work for the other side of $\mathrm{RII}_+$.

We play the exact same game for $\mathrm{RII}_-$: $\mathrm{Kh}
\left(
\tikz[baseline={([yshift=-.5ex]current bounding box.center)}, scale=.85]
{
        \draw[dotted] (.5,.5) circle(0.707);
        \draw[knot, ->] (0,0) to[out=30, in=-90] (0.8, 0.5);
        \draw[knot] (0.8, 0.5) to[out=90, in=-30] (0,1);
        \draw[knot, overcross] (1,0) to[out=150, in=-90] (0.3, 0.5);
        \draw[knot, overcross] (0.3, 0.5) to[out=90, in=210] (1,1);
        \draw[knot, <-] (0.3, 0.499) -- (0.3, 0.501);
}
\right)$ is
\[
\tikz[baseline={([yshift=-.5ex]current bounding box.center)}]
{
    \node at (5,0) {$\oplus$};
    \node(res00) at (-1,1.25) {$
    \varphi_{
    \tikz[baseline={([yshift=-.5ex]current bounding box.center)}, scale=.45]
        {
        \draw[dotted] (0.5, 0.5) circle (0.707);
        \draw[knot, red] (0.35, 0.75) -- (0.65, 0.75);
    %
        \draw[knot] (0,0) to[out=45, in=135] (1,0);
        \draw[knot] (0,1) to[out=-45, in=90] (0.35, 0.75);
        \draw[knot] (0.35, 0.75) to[out=-90, in=90] (0.2, 0.5);
        \draw[knot] (0.2, 0.5) to[out=-90, in=180] (0.5, 0.4);
        \draw[knot] (0.5, 0.4) to[out=0, in=-90] (0.8, 0.5);
        \draw[knot] (0.8, 0.5) to[out=90, in=-90] (0.65, 0.75);
        \draw[knot] (0.65, 0.75) to[out=90, in=225] (1,1);
        }
        }\,
        \mathcal{F}
        \left(
        \tikz[baseline={([yshift=-.5ex]current bounding box.center)}, scale=.85]
        {
        \draw[dotted] (0.5, 0.5) circle (0.707);
    %
    %
        \draw[knot] (0,0) to[out=45, in=135] (1,0);
        \draw[knot] (0,1) to[out=-45, in=90] (0.35, 0.75);
        \draw[knot] (0.35, 0.75) to[out=-90, in=90] (0.2, 0.5);
        \draw[knot] (0.2, 0.5) to[out=-90, in=180] (0.5, 0.4);
        \draw[knot] (0.5, 0.4) to[out=0, in=-90] (0.8, 0.5);
        \draw[knot] (0.8, 0.5) to[out=90, in=-90] (0.65, 0.75);
        \draw[knot] (0.65, 0.75) to[out=90, in=225] (1,1);
        }
        \right)
    $};
    \node(res10) at (5,1.25) {$
    \varphi_{
        \tikz[baseline={([yshift=-.5ex]current bounding box.center)}, scale=.45]
        {
        \draw[dotted] (0.5, 0.5) circle (0.707);
    %
        \draw[knot, red] (0.35, 0.75) -- (0.65, 0.75);
        \draw[knot] (0,0) to[out=45, in=135] (1,0);
        \draw[knot] (0,1) to[out=-45, in=90] (0.35, 0.75);
        \draw[knot] (0.35, 0.75) to[out=-90, in=90] (0.2, 0.5);
        \draw[knot] (0.2, 0.5) to[out=-90, in=180] (0.5, 0.4);
        \draw[knot] (0.5, 0.4) to[out=0, in=-90] (0.8, 0.5);
        \draw[knot] (0.8, 0.5) to[out=90, in=-90] (0.65, 0.75);
        \draw[knot] (0.65, 0.75) to[out=90, in=225] (1,1);
        }
    }
    \circ
    \varphi_{
        \tikz[baseline={([yshift=-.5ex]current bounding box.center)}, scale=.45]
        {
        \draw[dotted] (0.5, 0.5) circle (0.707);
        \draw[knot, red] (0.5, 0.4) -- (0.5, 0.2);
    %
        \draw[knot] (0,0) to[out=45, in=135] (1,0);
        \draw[knot] (0,1) to[out=-45, in=90] (0.35, 0.75);
        \draw[knot] (0.35, 0.75) to[out=-90, in=90] (0.2, 0.5);
        \draw[knot] (0.2, 0.5) to[out=-90, in=180] (0.5, 0.4);
        \draw[knot] (0.5, 0.4) to[out=0, in=-90] (0.8, 0.5);
        \draw[knot] (0.8, 0.5) to[out=90, in=-90] (0.65, 0.75);
        \draw[knot] (0.65, 0.75) to[out=90, in=225] (1,1);
        }
    }^{-1}
 \mathcal{F}
\left(
\tikz[baseline={([yshift=-.5ex]current bounding box.center)}, scale=.85]
{
        \draw[dotted] (0.5, 0.5) circle (0.707);

        \draw[knot] (0,1) to[out=-45, in=90] (0.35, 0.75);
        \draw[knot] (0.35, 0.75) to[out=-90, in=90] (0.2, 0.5);
        \draw[knot] (0.2, 0.5) to[out=-90, in=90] (0.35, 0.25);
        \draw[knot] (0.35, 0.25) to[out=-90, in=45] (0,0);
        \draw[knot] (1,1) to[out=225, in=90] (0.65, 0.75);
        \draw[knot] (0.65, 0.75) to[out=-90, in=90] (0.8, 0.5);
        \draw[knot] (0.8, 0.5) to[out=-90, in=90] (0.65, 0.25);
        \draw[knot] (0.65, 0.25) to[out=-90, in=135] (1,0);
}
\right)
    $};
    \node(res01) at (5,-1.25) {$
    \mathcal{F}
    \left(
    \tikz[baseline={([yshift=-.5ex]current bounding box.center)}, scale=.85]
        {
        \draw[dotted] (0.5, 0.5) circle (0.707);
        \draw[knot] (0,0) to[out=45, in=135] (1,0);
        \draw[knot] (0,1) to[out=-45, in=-135] (1,1);
        \draw[knot] (0.5, 0.65) to[out=0, in=90] (0.8, 0.5) to[out=-90, in=0] (0.5, 0.35) to[out=180, in=-90] (0.2, 0.5) to[out=90, in=180] (0.5, 0.65);
        }
    \right)
    $};
    \node(res11) at (11, -1.25){$
    \varphi_{
    \tikz[baseline={([yshift=-.5ex]current bounding box.center)}, scale=.45]
        {
        \draw[dotted] (0.5, 0.5) circle (0.707);
    %
        \draw[knot, red] (0.5, 0.35) -- (0.5, 0.2);    
        \draw[knot] (0,0) to[out=45, in=135] (1,0);
        \draw[knot] (0,1) to[out=-45, in=-135] (1,1);
        \draw[knot] (0.5, 0.65) to[out=0, in=90] (0.8, 0.5) to[out=-90, in=0] (0.5, 0.35) to[out=180, in=-90] (0.2, 0.5) to[out=90, in=180] (0.5, 0.65);
        }
        }^{-1}\,
        \mathcal{F}
        \left(
        \tikz[baseline={([yshift=-.5ex]current bounding box.center)}, scale=.85]
        {
        \draw[dotted] (0.5, 0.5) circle (0.707);
        \draw[knot] (0,1) to[out=-45, in=-135] (1,1);
        \draw[knot] (0,0) to[out=45, in=-90] (0.35, 0.25);
        \draw[knot] (0.35, 0.25) to[out=90, in=-90] (0.2, 0.5);
        \draw[knot] (0.2, 0.5) to[out=90, in=180] (0.5, 0.6);
        \draw[knot] (0.5, 0.6) to[out=0, in=90] (0.8, 0.5);
        \draw[knot] (0.8, 0.5) to[out=-90, in=90] (0.65, 0.25);
        \draw[knot] (0.65, 0.25) to[out=-90, in=135] (1,0);
        }
        \right)
    $};
    \draw[->] (res00) to
        node[pos=.5,above,arrows=-] {$\tikz[baseline={([yshift=-.5ex]current bounding box.center)}, scale=.45]
        {
        \draw[dotted] (0.5, 0.5) circle (0.707);
    %
        \draw[knot, red] (0.5, 0.4) -- (0.5, 0.2);
        \draw[knot] (0,0) to[out=45, in=135] (1,0);
        \draw[knot] (0,1) to[out=-45, in=90] (0.35, 0.75);
        \draw[knot] (0.35, 0.75) to[out=-90, in=90] (0.2, 0.5);
        \draw[knot] (0.2, 0.5) to[out=-90, in=180] (0.5, 0.4);
        \draw[knot] (0.5, 0.4) to[out=0, in=-90] (0.8, 0.5);
        \draw[knot] (0.8, 0.5) to[out=90, in=-90] (0.65, 0.75);
        \draw[knot] (0.65, 0.75) to[out=90, in=225] (1,1);
        }
        \circ \varphi_{H_1'}
        $}
        (res10);
    \draw[->] (res00) to 
        node[pos=.2, below=1.5mm, arrows=-] {$
        \tikz[baseline={([yshift=-.5ex]current bounding box.center)}, scale=.45]
        {
        \draw[dotted] (0.5, 0.5) circle (0.707);
        \draw[knot, red] (0.35, 0.75) -- (0.65, 0.75);
    %
        \draw[knot] (0,0) to[out=45, in=135] (1,0);
        \draw[knot] (0,1) to[out=-45, in=90] (0.35, 0.75);
        \draw[knot] (0.35, 0.75) to[out=-90, in=90] (0.2, 0.5);
        \draw[knot] (0.2, 0.5) to[out=-90, in=180] (0.5, 0.4);
        \draw[knot] (0.5, 0.4) to[out=0, in=-90] (0.8, 0.5);
        \draw[knot] (0.8, 0.5) to[out=90, in=-90] (0.65, 0.75);
        \draw[knot] (0.65, 0.75) to[out=90, in=225] (1,1);
        }
        $}
        (res01);
    \draw[->] (res10) to
        node[pos=.6, above=1.5mm, arrows=-] {$
        \tikz[baseline={([yshift=-.5ex]current bounding box.center)}, scale=.45]
        {
        \draw[dotted] (0.5, 0.5) circle (0.707);
        \draw[knot, red] (0.35, 0.75) -- (0.65, 0.75);
    %
        \draw[knot] (0,1) to[out=-45, in=90] (0.35, 0.75);
        \draw[knot] (0.35, 0.75) to[out=-90, in=90] (0.2, 0.5);
        \draw[knot] (0.2, 0.5) to[out=-90, in=90] (0.35, 0.25);
        \draw[knot] (0.35, 0.25) to[out=-90, in=45] (0,0);
        \draw[knot] (1,1) to[out=225, in=90] (0.65, 0.75);
        \draw[knot] (0.65, 0.75) to[out=-90, in=90] (0.8, 0.5);
        \draw[knot] (0.8, 0.5) to[out=-90, in=90] (0.65, 0.25);
        \draw[knot] (0.65, 0.25) to[out=-90, in=135] (1,0);
        }
        $}
        (res11);
    \draw[->] (res01) to
        node[pos=.5, above, arrows=-] {$
        \tikz[baseline={([yshift=-.5ex]current bounding box.center)}, scale=.45]
        {
        \draw[dotted] (0.5, 0.5) circle (0.707);
    %
        \draw[knot, red] (0.5, 0.35) -- (0.5, 0.2);
        \draw[knot] (0,0) to[out=45, in=135] (1,0);
        \draw[knot] (0,1) to[out=-45, in=-135] (1,1);
        \draw[knot] (0.5, 0.65) to[out=0, in=90] (0.8, 0.5) to[out=-90, in=0] (0.5, 0.35) to[out=180, in=-90] (0.2, 0.5) to[out=90, in=180] (0.5, 0.65);
        }
        \circ \varphi_{H_2'}
        $}
        (res11);
}
\]
with a global shift of $\{-1,1\}$. By grading shift arethmetic, we know
\[
\varphi_{
    \tikz[baseline={([yshift=-.5ex]current bounding box.center)}, scale=.45]
        {
        \draw[dotted] (0.5, 0.5) circle (0.707);
        \draw[knot, red] (0.35, 0.75) -- (0.65, 0.75);
    %
        \draw[knot] (0,0) to[out=45, in=135] (1,0);
        \draw[knot] (0,1) to[out=-45, in=90] (0.35, 0.75);
        \draw[knot] (0.35, 0.75) to[out=-90, in=90] (0.2, 0.5);
        \draw[knot] (0.2, 0.5) to[out=-90, in=180] (0.5, 0.4);
        \draw[knot] (0.5, 0.4) to[out=0, in=-90] (0.8, 0.5);
        \draw[knot] (0.8, 0.5) to[out=90, in=-90] (0.65, 0.75);
        \draw[knot] (0.65, 0.75) to[out=90, in=225] (1,1);
        }
        }
\cong \{0,-1\},
\qquad
\varphi_{
        \tikz[baseline={([yshift=-.5ex]current bounding box.center)}, scale=.45]
        {
        \draw[dotted] (0.5, 0.5) circle (0.707);
        \draw[knot, red] (0.5, 0.4) -- (0.5, 0.2);
    %
        \draw[knot] (0,0) to[out=45, in=135] (1,0);
        \draw[knot] (0,1) to[out=-45, in=90] (0.35, 0.75);
        \draw[knot] (0.35, 0.75) to[out=-90, in=90] (0.2, 0.5);
        \draw[knot] (0.2, 0.5) to[out=-90, in=180] (0.5, 0.4);
        \draw[knot] (0.5, 0.4) to[out=0, in=-90] (0.8, 0.5);
        \draw[knot] (0.8, 0.5) to[out=90, in=-90] (0.65, 0.75);
        \draw[knot] (0.65, 0.75) to[out=90, in=225] (1,1);
        }
    }^{-1}
\cong 
\varphi_{\left(
        \tikz[baseline={([yshift=-.5ex]current bounding box.center)}, scale=.45]
        {
        \draw[dotted] (0.5, 0.5) circle (0.707);
    %
        \draw[knot, red] (0.35, 0.25) -- (0.65, 0.25);
        \draw[knot] (0,1) to[out=-45, in=90] (0.35, 0.75);
        \draw[knot] (0.35, 0.75) to[out=-90, in=90] (0.2, 0.5);
        \draw[knot] (0.2, 0.5) to[out=-90, in=90] (0.35, 0.25);
        \draw[knot] (0.35, 0.25) to[out=-90, in=45] (0,0);
        \draw[knot] (1,1) to[out=225, in=90] (0.65, 0.75);
        \draw[knot] (0.65, 0.75) to[out=-90, in=90] (0.8, 0.5);
        \draw[knot] (0.8, 0.5) to[out=-90, in=90] (0.65, 0.25);
        \draw[knot] (0.65, 0.25) to[out=-90, in=135] (1,0);
        }
        \, (1,1)\right)
        }
,~\text{and}
\qquad
\varphi_{
    \tikz[baseline={([yshift=-.5ex]current bounding box.center)}, scale=.45]
        {
        \draw[dotted] (0.5, 0.5) circle (0.707);
    %
        \draw[knot, red] (0.5, 0.35) -- (0.5, 0.2);    
        \draw[knot] (0,0) to[out=45, in=135] (1,0);
        \draw[knot] (0,1) to[out=-45, in=-135] (1,1);
        \draw[knot] (0.5, 0.65) to[out=0, in=90] (0.8, 0.5) to[out=-90, in=0] (0.5, 0.35) to[out=180, in=-90] (0.2, 0.5) to[out=90, in=180] (0.5, 0.65);
        }
        }^{-1}
\cong \{1,0\},
\]
so that the complex may be rewritten
\[
\tikz[baseline={([yshift=-.5ex]current bounding box.center)}]
{
    \node at (5,0) {$\oplus$};
    \node(res00) at (-1,1.25) {$
        \mathcal{F}
        \left(
        \tikz[baseline={([yshift=-.5ex]current bounding box.center)}, scale=.85]
        {
        \draw[dotted] (0.5, 0.5) circle (0.707);
    %
    %
        \draw[knot] (0,0) to[out=45, in=135] (1,0);
        \draw[knot] (0,1) to[out=-45, in=90] (0.35, 0.75);
        \draw[knot] (0.35, 0.75) to[out=-90, in=90] (0.2, 0.5);
        \draw[knot] (0.2, 0.5) to[out=-90, in=180] (0.5, 0.4);
        \draw[knot] (0.5, 0.4) to[out=0, in=-90] (0.8, 0.5);
        \draw[knot] (0.8, 0.5) to[out=90, in=-90] (0.65, 0.75);
        \draw[knot] (0.65, 0.75) to[out=90, in=225] (1,1);
        }
        \right)
        \{0,-1\}
    $};
    \node(res10) at (5,1.25) {$
    \varphi_{\left(
        \tikz[baseline={([yshift=-.5ex]current bounding box.center)}, scale=.45]
        {
        \draw[dotted] (0.5, 0.5) circle (0.707);
    %
        \draw[knot, red] (0.35, 0.25) -- (0.65, 0.25);
        \draw[knot] (0,1) to[out=-45, in=90] (0.35, 0.75);
        \draw[knot] (0.35, 0.75) to[out=-90, in=90] (0.2, 0.5);
        \draw[knot] (0.2, 0.5) to[out=-90, in=90] (0.35, 0.25);
        \draw[knot] (0.35, 0.25) to[out=-90, in=45] (0,0);
        \draw[knot] (1,1) to[out=225, in=90] (0.65, 0.75);
        \draw[knot] (0.65, 0.75) to[out=-90, in=90] (0.8, 0.5);
        \draw[knot] (0.8, 0.5) to[out=-90, in=90] (0.65, 0.25);
        \draw[knot] (0.65, 0.25) to[out=-90, in=135] (1,0);
        }
        \, (1,0)\right)
        }
 \mathcal{F}
\left(
\tikz[baseline={([yshift=-.5ex]current bounding box.center)}, scale=.85]
{
        \draw[dotted] (0.5, 0.5) circle (0.707);

        \draw[knot] (0,1) to[out=-45, in=90] (0.35, 0.75);
        \draw[knot] (0.35, 0.75) to[out=-90, in=90] (0.2, 0.5);
        \draw[knot] (0.2, 0.5) to[out=-90, in=90] (0.35, 0.25);
        \draw[knot] (0.35, 0.25) to[out=-90, in=45] (0,0);
        \draw[knot] (1,1) to[out=225, in=90] (0.65, 0.75);
        \draw[knot] (0.65, 0.75) to[out=-90, in=90] (0.8, 0.5);
        \draw[knot] (0.8, 0.5) to[out=-90, in=90] (0.65, 0.25);
        \draw[knot] (0.65, 0.25) to[out=-90, in=135] (1,0);
}
\right)
    $};
    \node(res01) at (5,-1.25) {$
    \mathcal{F}
    \left(
    \tikz[baseline={([yshift=-.5ex]current bounding box.center)}, scale=.85]
        {
        \draw[dotted] (0.5, 0.5) circle (0.707);
        \draw[knot] (0,0) to[out=45, in=135] (1,0);
        \draw[knot] (0,1) to[out=-45, in=-135] (1,1);
        \draw[knot] (0.5, 0.65) to[out=0, in=90] (0.8, 0.5) to[out=-90, in=0] (0.5, 0.35) to[out=180, in=-90] (0.2, 0.5) to[out=90, in=180] (0.5, 0.65);
        }
    \right)
    $};
    \node(res11) at (11, -1.25){$
        \mathcal{F}
        \left(
        \tikz[baseline={([yshift=-.5ex]current bounding box.center)}, scale=.85]
        {
        \draw[dotted] (0.5, 0.5) circle (0.707);
        \draw[knot] (0,1) to[out=-45, in=-135] (1,1);
        \draw[knot] (0,0) to[out=45, in=-90] (0.35, 0.25);
        \draw[knot] (0.35, 0.25) to[out=90, in=-90] (0.2, 0.5);
        \draw[knot] (0.2, 0.5) to[out=90, in=180] (0.5, 0.6);
        \draw[knot] (0.5, 0.6) to[out=0, in=90] (0.8, 0.5);
        \draw[knot] (0.8, 0.5) to[out=-90, in=90] (0.65, 0.25);
        \draw[knot] (0.65, 0.25) to[out=-90, in=135] (1,0);
        }
        \right)
        \{1,0\}
    $};
    \draw[->] (res00) to
        node[pos=.5,above,arrows=-] {$\tikz[baseline={([yshift=-.5ex]current bounding box.center)}, scale=.45]
        {
        \draw[dotted] (0.5, 0.5) circle (0.707);
    %
        \draw[knot, red] (0.5, 0.4) -- (0.5, 0.2);
        \draw[knot] (0,0) to[out=45, in=135] (1,0);
        \draw[knot] (0,1) to[out=-45, in=90] (0.35, 0.75);
        \draw[knot] (0.35, 0.75) to[out=-90, in=90] (0.2, 0.5);
        \draw[knot] (0.2, 0.5) to[out=-90, in=180] (0.5, 0.4);
        \draw[knot] (0.5, 0.4) to[out=0, in=-90] (0.8, 0.5);
        \draw[knot] (0.8, 0.5) to[out=90, in=-90] (0.65, 0.75);
        \draw[knot] (0.65, 0.75) to[out=90, in=225] (1,1);
        }
        \circ \varphi_{H_1'}
        $}
        (res10);
    \draw[->] (res00) to 
        node[pos=.2, below=1.5mm, arrows=-] {$
        \tikz[baseline={([yshift=-.5ex]current bounding box.center)}, scale=.45]
        {
        \draw[dotted] (0.5, 0.5) circle (0.707);
        \draw[knot, red] (0.35, 0.75) -- (0.65, 0.75);
    %
        \draw[knot] (0,0) to[out=45, in=135] (1,0);
        \draw[knot] (0,1) to[out=-45, in=90] (0.35, 0.75);
        \draw[knot] (0.35, 0.75) to[out=-90, in=90] (0.2, 0.5);
        \draw[knot] (0.2, 0.5) to[out=-90, in=180] (0.5, 0.4);
        \draw[knot] (0.5, 0.4) to[out=0, in=-90] (0.8, 0.5);
        \draw[knot] (0.8, 0.5) to[out=90, in=-90] (0.65, 0.75);
        \draw[knot] (0.65, 0.75) to[out=90, in=225] (1,1);
        }
        $}
        (res01);
    \draw[->] (res10) to
        node[pos=.6, above=1.5mm, arrows=-] {$
        \tikz[baseline={([yshift=-.5ex]current bounding box.center)}, scale=.45]
        {
        \draw[dotted] (0.5, 0.5) circle (0.707);
        \draw[knot, red] (0.35, 0.75) -- (0.65, 0.75);
    %
        \draw[knot] (0,1) to[out=-45, in=90] (0.35, 0.75);
        \draw[knot] (0.35, 0.75) to[out=-90, in=90] (0.2, 0.5);
        \draw[knot] (0.2, 0.5) to[out=-90, in=90] (0.35, 0.25);
        \draw[knot] (0.35, 0.25) to[out=-90, in=45] (0,0);
        \draw[knot] (1,1) to[out=225, in=90] (0.65, 0.75);
        \draw[knot] (0.65, 0.75) to[out=-90, in=90] (0.8, 0.5);
        \draw[knot] (0.8, 0.5) to[out=-90, in=90] (0.65, 0.25);
        \draw[knot] (0.65, 0.25) to[out=-90, in=135] (1,0);
        }
        $}
        (res11);
    \draw[->] (res01) to
        node[pos=.5, above, arrows=-] {$
        \tikz[baseline={([yshift=-.5ex]current bounding box.center)}, scale=.45]
        {
        \draw[dotted] (0.5, 0.5) circle (0.707);
    %
        \draw[knot, red] (0.5, 0.35) -- (0.5, 0.2);
        \draw[knot] (0,0) to[out=45, in=135] (1,0);
        \draw[knot] (0,1) to[out=-45, in=-135] (1,1);
        \draw[knot] (0.5, 0.65) to[out=0, in=90] (0.8, 0.5) to[out=-90, in=0] (0.5, 0.35) to[out=180, in=-90] (0.2, 0.5) to[out=90, in=180] (0.5, 0.65);
        }
        \circ \varphi_{H_2'}
        $}
        (res11);
}
\]
Delooping the 01 entry and applying Gaussian elimination, we conclude that the entire complex is homotopy equivalent to $\varphi_{\left(\tikz[baseline=8.25ex, scale=.45]
{
    \begin{scope}[rotate=90]
	\draw[dotted] (3,-2) circle(0.707);
	\draw[knot] (2.5,-1.5) .. controls (2.75,-1.75) and (3.25,-1.75) .. (3.5,-1.5);
	\draw[knot] (2.5,-2.5) .. controls  (2.75,-2.25) and (3.25,-2.25) .. (3.5,-2.5);
        \draw[red, knot] (3,-1.7) -- (3,-2.3);
    \end{scope}
},\, (1,0)\right)}
\mathrm{Kh}
\left(
\tikz[baseline={([yshift=-.5ex]current bounding box.center)}, scale=.85]
{
        \draw[dotted] (.5,.5) circle(0.707);
        \draw[knot] (0, 0) to[out=30, in=-90] (0.35, 0.5);
        \draw[knot] (0.35, 0.5) to[out=90, in=-30] (0, 1);
        \draw[knot] (1, 0) to[out=150, in=-90] (0.65, 0.5);
        \draw[knot] (0.65, 0.5) to[out=90, in=210] (1, 1);
}
\right) \{-1,1\}$
; i.e., 
$\varphi_{\left(\tikz[baseline=8.25ex, scale=.45]
{
    \begin{scope}[rotate=90]
	\draw[dotted] (3,-2) circle(0.707);
	\draw[knot] (2.5,-1.5) .. controls (2.75,-1.75) and (3.25,-1.75) .. (3.5,-1.5);
	\draw[knot] (2.5,-2.5) .. controls  (2.75,-2.25) and (3.25,-2.25) .. (3.5,-2.5);
        \draw[red, knot] (3,-1.7) -- (3,-2.3);
    \end{scope}
},\, (0,1)\right)}
\mathrm{Kh}
\left(
\tikz[baseline={([yshift=-.5ex]current bounding box.center)}, scale=.85]
{
        \draw[dotted] (.5,.5) circle(0.707);
        \draw[knot] (0, 0) to[out=30, in=-90] (0.35, 0.5);
        \draw[knot] (0.35, 0.5) to[out=90, in=-30] (0, 1);
        \draw[knot] (1, 0) to[out=150, in=-90] (0.65, 0.5);
        \draw[knot] (0.65, 0.5) to[out=90, in=210] (1, 1);
}
\right)$, as desired. Again, the other side of $\mathrm{RII}_-$ is similar.
\end{proof}

\begin{remark}
\label{rem:R3noticing}
Lemma \ref{lem:r2} establishes that the grading shift coming from Reidemeister II moves is dependent on orientation. This, together with Lemma \ref{lem:r1}, implies that Reidemeister III moves must---at least, sometimes---come at the cost of a nontrivial grading shift. For example, if this was not the case, the sequence of isomorphisms
\[
\begin{tikzcd}[column sep = large, row sep = large]
\varphi_{\left(\tikz[baseline=8.25ex, scale=.45]
{
    \begin{scope}[rotate=90]
	\draw[dotted] (3,-2) circle(0.707);
	\draw[knot] (2.5,-1.5) .. controls (2.75,-1.75) and (3.25,-1.75) .. (3.5,-1.5);
	\draw[knot] (2.5,-2.5) .. controls  (2.75,-2.25) and (3.25,-2.25) .. (3.5,-2.5);
        \draw[red, knot] (3,-1.7) -- (3,-2.3);
    \end{scope}
},\, (0,1)\right)}
\mathrm{Kh}
\left(
\tikz[baseline={([yshift=-.5ex]current bounding box.center)}, scale=.85]
{
        \draw[dotted] (.5,.5) circle(0.707);
        \draw[knot] (0, 0) to[out=30, in=-90] (0.35, 0.5);
        \draw[knot] (0.35, 0.5) to[out=90, in=-30] (0, 1);
        \draw[knot] (1, 0) to[out=150, in=-90] (0.65, 0.5);
        \draw[knot] (0.65, 0.5) to[out=90, in=210] (1, 1);
}
\right)
\arrow[r, "\mathrm{RII}_-^{-1}"]
&
\mathrm{Kh}
\left(
\tikz[baseline={([yshift=-.5ex]current bounding box.center)}, scale=.85]
{
        \draw[dotted] (.5,.5) circle(0.707);
        \draw[knot, ->] (0,0) to[out=30, in=-90] (0.8, 0.5);
        \draw[knot] (0.8, 0.5) to[out=90, in=-30] (0,1);
        \draw[knot, overcross] (1,0) to[out=150, in=-90] (0.3, 0.5);
        \draw[knot, overcross] (0.3, 0.5) to[out=90, in=210] (1,1);
        \draw[knot, <-] (0.3, 0.499) -- (0.3, 0.501);
}
\right)
\arrow[r, "\mathrm{RI}^{-1}"]
&
\mathrm{Kh}
\left(
\tikz[baseline={([yshift=-.5ex]current bounding box.center)}, scale=.85]
{
        \draw[dotted] (.5,.5) circle(0.707);
        \draw[knot, ->] (0,0) to[out=30, in=-90] (0.8, 0.5);
        \draw[knot] (0.8, 0.5) to[out=90, in=-30] (0,1);
        \draw[knot, overcross] (1,0) to[out=180, in=0] (0.15, 0.7);
        \draw[knot] (0.15, 0.7) to[out=180, in=90] (0, 0.5);
        \draw[knot, ->] (0.15, 0.3) to[out=180, in=-90] (0, 0.5);
        \draw[knot, overcross] (0.15, 0.3) to[out=0, in=180] (1,1);
}
\right)
\arrow[d, "\star"]
\\
\mathrm{Kh}
\left(
\tikz[baseline={([yshift=-.5ex]current bounding box.center)}, scale=.85]
{
        \draw[dotted] (.5,.5) circle(0.707);
        \draw[knot] (0, 0) to[out=30, in=-90] (0.35, 0.5);
        \draw[knot] (0.35, 0.5) to[out=90, in=-30] (0, 1);
        \draw[knot] (1, 0) to[out=150, in=-90] (0.65, 0.5);
        \draw[knot] (0.65, 0.5) to[out=90, in=210] (1, 1);
}
\right)
\{-1,1\}
&
\mathrm{Kh}
\left(
\tikz[baseline={([yshift=-.5ex]current bounding box.center)}, scale=.85]
{
        \draw[dotted] (.5,.5) circle(0.707);
        \draw[knot, ->] (0, 0) to[out=30, in=-90] (0.2, 0.5);
        \draw[knot] (0.2, 0.5) to[out=90, in=-30] (0, 1);
        \draw[knot] (1, 0) to[out=135, in=0] (0.65, 0.8);
        \draw[knot] (0.65, 0.8) to[out=180, in=90] (0.5, 0.5);
        \draw[knot, <-] (0.5, 0.5) to[out=-90, in=180] (0.65, 0.2);
        \draw[knot, overcross] (0.65, 0.2) to[out=0, in=225] (1, 1);
    %
}
\right)
\{-1,1\}
\arrow[l, "\mathrm{RI}"']
&
\mathrm{Kh}
\left(
\tikz[baseline={([yshift=-.5ex]current bounding box.center)}, scale=.85]
{
        \draw[dotted] (.5,.5) circle(0.707);
        \draw[knot] (0, 0) to[out=30, in=-90] (0.55, 0.5);
        \draw[knot] (0.55, 0.5) to[out=90, in=-30] (0, 1);
        \draw[knot] (1, 0) to[out=150, in=0] (0.5, 0.8);
        \draw[knot, overcross] (0.5, 0.8) to[out=180, in=90] (0.25, 0.5);
        \draw[knot, overcross] (0.25, 0.5) to[out=-90, in=180] (0.5, 0.2);
        \draw[knot, overcross] (0.5, 0.2) to[out=0, in=210] (1, 1);
        \draw[knot, ->] (0.25, 0.499) -- (0.25, 0.501);
        \draw[knot, ->] (0.55, 0.499) -- (0.55, 0.501);
}
\right)
\arrow[l, "\mathrm{RII}_+"']
\end{tikzcd}
\]
would yield a contradiction. Notice that the vertical arrow is an Reidemeister III move of type
\[
\tikz[baseline={([yshift=-.5ex]current bounding box.center)}, scale=.85]
{
        \draw[dotted] (.5,.5) circle(0.707);
        \draw[knot] (0.5, -0.207) to[out=90, in=-90] (0.085, 0.5);
        \draw[knot] (0.085, 0.5) to[out=90, in=-90] (0.5, 1.207);
        \draw[knot, <-] (0.085, 0.499) -- (0.085, 0.501);
        \draw[knot, overcross] (0,1) -- (0.45, 0.55);
        \draw[knot] (0.45, 0.55) -- (0.55, 0.45);
        \draw[knot] (0.55, 0.45) -- (1,0);
        \draw[knot, <-] (0,1) -- (0.01, 0.99);
        \draw[knot, overcross] (0,0) -- (0.45, 0.45);
        \draw[knot, overcross] (0.45, 0.45) -- (0.55, 0.55);
        \draw[knot, ->] (0.55, 0.55) -- (1,1);
}
\rightsquigarrow
\tikz[baseline={([yshift=-.5ex]current bounding box.center)}, scale=.85]
{
        \draw[dotted] (.5,.5) circle(0.707);
        \draw[knot] (0.5, -0.207) to[out=90, in=-90] (0.915, 0.5);
        \draw[knot] (0.915, 0.5) to[out=90, in=-90] (0.5, 1.207);
        \draw[knot, <-] (0.915, 0.499) -- (0.915, 0.501);
        \draw[knot, <-] (0,1) -- (0.45, 0.55);
        \draw[knot] (0.45, 0.55) -- (0.55, 0.45);
        \draw[knot, overcross] (0.55, 0.45) -- (1,0);
        \draw[knot] (0,0) -- (0.45, 0.45);
        \draw[knot, overcross] (0.45, 0.45) -- (0.55, 0.55);
        \draw[knot, overcross] (0.55, 0.55) -- (1,1);
        \draw[knot, ->] (0.99, 0.99) -- (1,1);

}\,.
\]
\end{remark}

\begin{lemma}
\label{lem:r3}
We have the following isomorphisms in $\mathrm{Kom}(H^3\mathrm{Mod}_R^\mathscr{G})$:
\[
\mathrm{Kh}
\left(
\tikz[baseline={([yshift=-.5ex]current bounding box.center)}, scale=.85]
{
        \draw[dotted] (.5,.5) circle(0.707);
        \draw[knot, <-] (0,1) -- (0.45, 0.55);
        \draw[knot] (0.45, 0.55) -- (0.55, 0.45);
        \draw[knot] (0.55, 0.45) -- (1,0);
        \draw[knot] (0,0) -- (0.45, 0.45);
        \draw[knot, overcross] (0.45, 0.45) -- (0.55, 0.55);
        \draw[knot, ->] (0.55, 0.55) -- (1,1);
        \draw[knot, overcross] (0.5, -0.207) to[out=90, in=-90] (0.085, 0.5);
        \draw[knot, overcross] (0.085, 0.5) to[out=90, in=-90] (0.5, 1.207);
        \draw[knot, ->] (0.085, 0.499) -- (0.085, 0.501);
}
\right)
\cong
\mathrm{Kh}
\left(
\tikz[baseline={([yshift=-.5ex]current bounding box.center)}, scale=.85]
{
        \draw[dotted] (.5,.5) circle(0.707);
        \draw[knot, <-] (0,1) -- (0.45, 0.55);
        \draw[knot] (0.45, 0.55) -- (0.55, 0.45);
        \draw[knot] (0.55, 0.45) -- (1,0);
        \draw[knot] (0,0) -- (0.45, 0.45);
        \draw[knot, overcross] (0.45, 0.45) -- (0.55, 0.55);
        \draw[knot, ->] (0.55, 0.55) -- (1,1);
        \draw[knot, overcross] (0.5, -0.207) to[out=90, in=-90] (0.915, 0.5);
        \draw[knot, overcross] (0.915, 0.5) to[out=90, in=-90] (0.5, 1.207);
        \draw[knot, ->] (0.915, 0.499) -- (0.915, 0.501);
}
\right)
,
\qquad
\mathrm{Kh}
\left(
\tikz[baseline={([yshift=-.5ex]current bounding box.center)}, scale=.85]
{
        \draw[dotted] (.5,.5) circle(0.707);
        \draw[knot, <-] (0,1) -- (0.45, 0.55);
        \draw[knot] (0.45, 0.55) -- (0.55, 0.45);
        \draw[knot] (0.55, 0.45) -- (1,0);
        \draw[knot] (0,0) -- (0.45, 0.45);
        \draw[knot, overcross] (0.45, 0.45) -- (0.55, 0.55);
        \draw[knot, ->] (0.55, 0.55) -- (1,1);
        \draw[knot, overcross] (0.5, -0.207) to[out=90, in=-90] (0.085, 0.5);
        \draw[knot, overcross] (0.085, 0.5) to[out=90, in=-90] (0.5, 1.207);
        \draw[knot, <-] (0.085, 0.499) -- (0.085, 0.501);
}
\right)
\cong
\varphi_{
\tikz[baseline={([yshift=-.5ex]current bounding box.center)}, scale=.45]
{
        \draw[dotted] (.5,.5) circle(0.707); 
        \draw[knot, red] (0.2, 0.5) -- (0.5, 0.5);
        \draw[knot] (0,0) to[out=45, in=-45] (0,1);
        \draw[knot] (1,0) to[out=135, in=-135] (1,1);
        \draw[knot] (0.5, -0.207) -- (0.5, 1.207);
}
}
\circ
\varphi_{
\tikz[baseline={([yshift=-.5ex]current bounding box.center)}, scale=.45]
{
        \draw[dotted] (.5,.5) circle(0.707); 
        \draw[knot, red] (0.8, 0.5) -- (0.5, 0.5);
        \draw[knot] (0,0) to[out=45, in=-45] (0,1);
        \draw[knot] (1,0) to[out=135, in=-135] (1,1);
        \draw[knot] (0.5, -0.207) -- (0.5, 1.207);
}
}^{-1}
\mathrm{Kh}
\left(
\tikz[baseline={([yshift=-.5ex]current bounding box.center)}, scale=.85]
{
        \draw[dotted] (.5,.5) circle(0.707);
        \draw[knot, <-] (0,1) -- (0.45, 0.55);
        \draw[knot] (0.45, 0.55) -- (0.55, 0.45);
        \draw[knot] (0.55, 0.45) -- (1,0);
        \draw[knot] (0,0) -- (0.45, 0.45);
        \draw[knot, overcross] (0.45, 0.45) -- (0.55, 0.55);
        \draw[knot, ->] (0.55, 0.55) -- (1,1);
        \draw[knot, overcross] (0.5, -0.207) to[out=90, in=-90] (0.915, 0.5);
        \draw[knot, overcross] (0.915, 0.5) to[out=90, in=-90] (0.5, 1.207);
        \draw[knot, <-] (0.915, 0.499) -- (0.915, 0.501);

}
\right),
\]
\[
\mathrm{Kh}
\left(
\tikz[baseline={([yshift=-.5ex]current bounding box.center)}, scale=.85]
{
        \draw[dotted] (.5,.5) circle(0.707);
        \draw[knot] (0,0) -- (0.45, 0.45);
        \draw[knot] (0.45, 0.45) -- (0.55, 0.55);
        \draw[knot, ->] (0.55, 0.55) -- (1,1);
        \draw[knot, <-] (0,1) -- (0.45, 0.55);
        \draw[knot, overcross] (0.45, 0.55) -- (0.55, 0.45);
        \draw[knot] (0.55, 0.45) -- (1,0);
        \draw[knot, overcross] (0.5, -0.207) to[out=90, in=-90] (0.085, 0.5);
        \draw[knot, overcross] (0.085, 0.5) to[out=90, in=-90] (0.5, 1.207);
        \draw[knot, ->] (0.085, 0.499) -- (0.085, 0.501);
}
\right)
\cong
\mathrm{Kh}
\left(
\tikz[baseline={([yshift=-.5ex]current bounding box.center)}, scale=.85]
{
        \draw[dotted] (.5,.5) circle(0.707);
        \draw[knot] (0,0) -- (0.45, 0.45);
        \draw[knot] (0.45, 0.45) -- (0.55, 0.55);
        \draw[knot, ->] (0.55, 0.55) -- (1,1);
        \draw[knot, <-] (0,1) -- (0.45, 0.55);
        \draw[knot, overcross] (0.45, 0.55) -- (0.55, 0.45);
        \draw[knot] (0.55, 0.45) -- (1,0);
        \draw[knot, overcross] (0.5, -0.207) to[out=90, in=-90] (0.915, 0.5);
        \draw[knot, overcross] (0.915, 0.5) to[out=90, in=-90] (0.5, 1.207);
        \draw[knot, ->] (0.915, 0.499) -- (0.915, 0.501);
}
\right)
,
\qquad
\mathrm{Kh}
\left(
\tikz[baseline={([yshift=-.5ex]current bounding box.center)}, scale=.85]
{
        \draw[dotted] (.5,.5) circle(0.707);
        \draw[knot] (0,0) -- (0.45, 0.45);
        \draw[knot] (0.45, 0.45) -- (0.55, 0.55);
        \draw[knot, ->] (0.55, 0.55) -- (1,1);
        \draw[knot, <-] (0,1) -- (0.45, 0.55);
        \draw[knot, overcross] (0.45, 0.55) -- (0.55, 0.45);
        \draw[knot] (0.55, 0.45) -- (1,0);
        \draw[knot, overcross] (0.5, -0.207) to[out=90, in=-90] (0.085, 0.5);
        \draw[knot, overcross] (0.085, 0.5) to[out=90, in=-90] (0.5, 1.207);
        \draw[knot, <-] (0.085, 0.499) -- (0.085, 0.501);
}
\right)
\cong
\varphi_{
\tikz[baseline={([yshift=-.5ex]current bounding box.center)}, scale=.45]
{
        \draw[dotted] (.5,.5) circle(0.707); 
        \draw[knot, red] (0.2, 0.5) -- (0.5, 0.5);
        \draw[knot] (0,0) to[out=45, in=-45] (0,1);
        \draw[knot] (1,0) to[out=135, in=-135] (1,1);
        \draw[knot] (0.5, -0.207) -- (0.5, 1.207);
}
}
\circ
\varphi_{
\tikz[baseline={([yshift=-.5ex]current bounding box.center)}, scale=.45]
{
        \draw[dotted] (.5,.5) circle(0.707); 
        \draw[knot, red] (0.8, 0.5) -- (0.5, 0.5);
        \draw[knot] (0,0) to[out=45, in=-45] (0,1);
        \draw[knot] (1,0) to[out=135, in=-135] (1,1);
        \draw[knot] (0.5, -0.207) -- (0.5, 1.207);
}
}^{-1}
\mathrm{Kh}
\left(
\tikz[baseline={([yshift=-.5ex]current bounding box.center)}, scale=.85]
{
        \draw[dotted] (.5,.5) circle(0.707);
        \draw[knot] (0,0) -- (0.45, 0.45);
        \draw[knot] (0.45, 0.45) -- (0.55, 0.55);
        \draw[knot, ->] (0.55, 0.55) -- (1,1);
        \draw[knot, <-] (0,1) -- (0.45, 0.55);
        \draw[knot, overcross] (0.45, 0.55) -- (0.55, 0.45);
        \draw[knot] (0.55, 0.45) -- (1,0);
        \draw[knot, overcross] (0.5, -0.207) to[out=90, in=-90] (0.915, 0.5);
        \draw[knot, overcross] (0.915, 0.5) to[out=90, in=-90] (0.5, 1.207);
        \draw[knot, <-] (0.915, 0.499) -- (0.915, 0.501);

}
\right),
\]
\[
\mathrm{Kh}
\left(
\tikz[baseline={([yshift=-.5ex]current bounding box.center)}, scale=.85]
{
        \draw[dotted] (.5,.5) circle(0.707);
        \draw[knot] (0.5, -0.207) to[out=90, in=-90] (0.085, 0.5);
        \draw[knot] (0.085, 0.5) to[out=90, in=-90] (0.5, 1.207);
        \draw[knot, ->] (0.085, 0.499) -- (0.085, 0.501);
        \draw[knot, overcross] (0,1) -- (0.45, 0.55);
        \draw[knot, <-] (0,1) -- (0.01, 0.99);
        \draw[knot] (0.45, 0.55) -- (0.55, 0.45);
        \draw[knot] (0.55, 0.45) -- (1,0);
        \draw[knot, overcross] (0,0) -- (0.45, 0.45);
        \draw[knot, overcross] (0.45, 0.45) -- (0.55, 0.55);
        \draw[knot, ->] (0.55, 0.55) -- (1,1);
}
\right)
\cong
\mathrm{Kh}
\left(
\tikz[baseline={([yshift=-.5ex]current bounding box.center)}, scale=.85]
{
        \draw[dotted] (.5,.5) circle(0.707);
        \draw[knot] (0.5, -0.207) to[out=90, in=-90] (0.915, 0.5);
        \draw[knot] (0.915, 0.5) to[out=90, in=-90] (0.5, 1.207);
        \draw[knot, ->] (0.915, 0.499) -- (0.915, 0.501);
        \draw[knot, <-] (0,1) -- (0.45, 0.55);
        \draw[knot] (0.45, 0.55) -- (0.55, 0.45);
        \draw[knot, overcross] (0.55, 0.45) -- (1,0);
        \draw[knot] (0,0) -- (0.45, 0.45);
        \draw[knot, overcross] (0.45, 0.45) -- (0.55, 0.55);
        \draw[knot, overcross] (0.55, 0.55) -- (1,1);
        \draw[knot, ->] (0.99, 0.99) -- (1,1);
}
\right)
,
\qquad
\mathrm{Kh}
\left(
\tikz[baseline={([yshift=-.5ex]current bounding box.center)}, scale=.85]
{
        \draw[dotted] (.5,.5) circle(0.707);
        \draw[knot] (0.5, -0.207) to[out=90, in=-90] (0.085, 0.5);
        \draw[knot] (0.085, 0.5) to[out=90, in=-90] (0.5, 1.207);
        \draw[knot, <-] (0.085, 0.499) -- (0.085, 0.501);
        \draw[knot, overcross] (0,1) -- (0.45, 0.55);
        \draw[knot] (0.45, 0.55) -- (0.55, 0.45);
        \draw[knot] (0.55, 0.45) -- (1,0);
        \draw[knot, <-] (0,1) -- (0.01, 0.99);
        \draw[knot, overcross] (0,0) -- (0.45, 0.45);
        \draw[knot, overcross] (0.45, 0.45) -- (0.55, 0.55);
        \draw[knot, ->] (0.55, 0.55) -- (1,1);
}
\right)
\cong
\varphi_{
\tikz[baseline={([yshift=-.5ex]current bounding box.center)}, scale=.45]
{
        \draw[dotted] (.5,.5) circle(0.707); 
        \draw[knot, red] (0.2, 0.5) -- (0.5, 0.5);
        \draw[knot] (0,0) to[out=45, in=-45] (0,1);
        \draw[knot] (1,0) to[out=135, in=-135] (1,1);
        \draw[knot] (0.5, -0.207) -- (0.5, 1.207);
}
}
\circ
\varphi_{
\tikz[baseline={([yshift=-.5ex]current bounding box.center)}, scale=.45]
{
        \draw[dotted] (.5,.5) circle(0.707); 
        \draw[knot, red] (0.8, 0.5) -- (0.5, 0.5);
        \draw[knot] (0,0) to[out=45, in=-45] (0,1);
        \draw[knot] (1,0) to[out=135, in=-135] (1,1);
        \draw[knot] (0.5, -0.207) -- (0.5, 1.207);
}
}^{-1}
\mathrm{Kh}
\left(
\tikz[baseline={([yshift=-.5ex]current bounding box.center)}, scale=.85]
{
        \draw[dotted] (.5,.5) circle(0.707);
        \draw[knot] (0.5, -0.207) to[out=90, in=-90] (0.915, 0.5);
        \draw[knot] (0.915, 0.5) to[out=90, in=-90] (0.5, 1.207);
        \draw[knot, <-] (0.915, 0.499) -- (0.915, 0.501);
        \draw[knot, <-] (0,1) -- (0.45, 0.55);
        \draw[knot] (0.45, 0.55) -- (0.55, 0.45);
        \draw[knot, overcross] (0.55, 0.45) -- (1,0);
        \draw[knot] (0,0) -- (0.45, 0.45);
        \draw[knot, overcross] (0.45, 0.45) -- (0.55, 0.55);
        \draw[knot, overcross] (0.55, 0.55) -- (1,1);
        \draw[knot, ->] (0.99, 0.99) -- (1,1);

}
\right),
\]
\[
\mathrm{Kh}
\left(
\tikz[baseline={([yshift=-.5ex]current bounding box.center)}, scale=.85]
{
        \draw[dotted] (.5,.5) circle(0.707);
        \draw[knot] (0.5, -0.207) to[out=90, in=-90] (0.085, 0.5);
        \draw[knot] (0.085, 0.5) to[out=90, in=-90] (0.5, 1.207);
        \draw[knot, ->] (0.085, 0.499) -- (0.085, 0.501);
        \draw[knot, overcross] (0,0) -- (0.45, 0.45);
        \draw[knot] (0.45, 0.45) -- (0.55, 0.55);
        \draw[knot, ->] (0.55, 0.55) -- (1,1);
        \draw[knot, overcross] (0,1) -- (0.45, 0.55);
        \draw[knot, <-] (0,1) -- (0.01, 0.99);
        \draw[knot, overcross] (0.45, 0.55) -- (0.55, 0.45);
        \draw[knot] (0.55, 0.45) -- (1,0);
}
\right)
\cong
\mathrm{Kh}
\left(
\tikz[baseline={([yshift=-.5ex]current bounding box.center)}, scale=.85]
{
        \draw[dotted] (.5,.5) circle(0.707);
        \draw[knot] (0.5, -0.207) to[out=90, in=-90] (0.915, 0.5);
        \draw[knot] (0.915, 0.5) to[out=90, in=-90] (0.5, 1.207);
        \draw[knot, ->] (0.915, 0.499) -- (0.915, 0.501);
        \draw[knot] (0,0) -- (0.45, 0.45);
        \draw[knot] (0.45, 0.45) -- (0.55, 0.55);
        \draw[knot, overcross] (0.55, 0.55) -- (1,1);
        \draw[knot, ->] (0.99, 0.99) -- (1,1);
        \draw[knot, <-] (0,1) -- (0.45, 0.55);
        \draw[knot, overcross] (0.45, 0.55) -- (0.55, 0.45);
        \draw[knot, overcross] (0.55, 0.45) -- (1,0);
}
\right)
,\qquad
\mathrm{Kh}
\left(
\tikz[baseline={([yshift=-.5ex]current bounding box.center)}, scale=.85]
{
        \draw[dotted] (.5,.5) circle(0.707);
        \draw[knot] (0.5, -0.207) to[out=90, in=-90] (0.085, 0.5);
        \draw[knot] (0.085, 0.5) to[out=90, in=-90] (0.5, 1.207);
        \draw[knot, <-] (0.085, 0.499) -- (0.085, 0.501);
        \draw[knot, overcross] (0,0) -- (0.45, 0.45);
        \draw[knot] (0.45, 0.45) -- (0.55, 0.55);
        \draw[knot, ->] (0.55, 0.55) -- (1,1);
        \draw[knot, overcross] (0,1) -- (0.45, 0.55);
        \draw[knot, overcross] (0.45, 0.55) -- (0.55, 0.45);
        \draw[knot] (0.55, 0.45) -- (1,0);
        \draw[knot, <-] (0,1) -- (0.01, 0.99);
}
\right)
\cong
\varphi_{
\tikz[baseline={([yshift=-.5ex]current bounding box.center)}, scale=.45]
{
        \draw[dotted] (.5,.5) circle(0.707); 
        \draw[knot, red] (0.2, 0.5) -- (0.5, 0.5);
        \draw[knot] (0,0) to[out=45, in=-45] (0,1);
        \draw[knot] (1,0) to[out=135, in=-135] (1,1);
        \draw[knot] (0.5, -0.207) -- (0.5, 1.207);
}
}
\circ
\varphi_{
\tikz[baseline={([yshift=-.5ex]current bounding box.center)}, scale=.45]
{
        \draw[dotted] (.5,.5) circle(0.707); 
        \draw[knot, red] (0.8, 0.5) -- (0.5, 0.5);
        \draw[knot] (0,0) to[out=45, in=-45] (0,1);
        \draw[knot] (1,0) to[out=135, in=-135] (1,1);
        \draw[knot] (0.5, -0.207) -- (0.5, 1.207);
}
}^{-1}
\mathrm{Kh}
\left(
\tikz[baseline={([yshift=-.5ex]current bounding box.center)}, scale=.85]
{
        \draw[dotted] (.5,.5) circle(0.707);
        \draw[knot] (0.5, -0.207) to[out=90, in=-90] (0.915, 0.5);
        \draw[knot] (0.915, 0.5) to[out=90, in=-90] (0.5, 1.207);
        \draw[knot, <-] (0.915, 0.499) -- (0.915, 0.501);
        \draw[knot] (0,0) -- (0.45, 0.45);
        \draw[knot] (0.45, 0.45) -- (0.55, 0.55);
        \draw[knot, overcross] (0.55, 0.55) -- (1,1);
        \draw[knot, ->] (0.99, 0.99) -- (1,1);
        \draw[knot, <-] (0,1) -- (0.45, 0.55);
        \draw[knot, overcross] (0.45, 0.55) -- (0.55, 0.45);
        \draw[knot, overcross] (0.55, 0.45) -- (1,0);
}
\right).
\]
\end{lemma}

\begin{proof}
We will describe the proof by illustrating one of the isomorphisms on the left-hand side and its counterpart on the right-hand side. Each computation is slightly different, but we hope that this discussion sates the reader, or illuminates the procedure enough so that they might check the others on their own.

The idea for any isomorphism on the left-hand side is to expand each complex and apply Gaussian elimination carefully. If Gaussian elimination is done properly, the two complexes are isotopic. If we do the same procedure for complexes appearing on the right-hand side, we will find that the entries of the complex are isotopic, but the grading shifts disagree. In this case, we will argue that one is taken to the other by applying the grading shifts provided in the statement of the Lemma.

Observe the complex associated to $\mathrm{Kh}
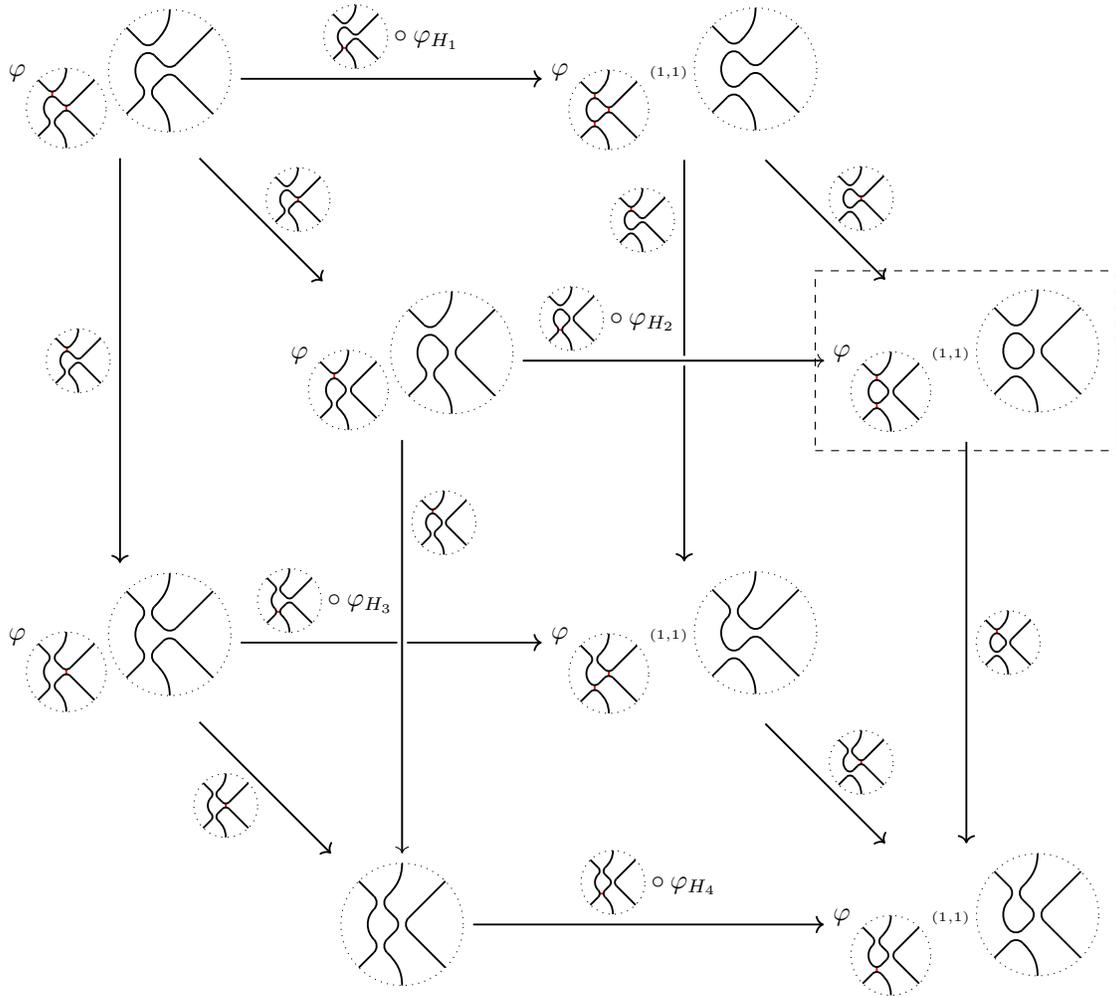
\left(
\tikz[baseline={([yshift=-.5ex]current bounding box.center)}, scale=.85]
{
        \draw[dotted] (.5,.5) circle(0.707);
        \draw[knot] (0.5, -0.207) to[out=90, in=-90] (0.085, 0.5);
        \draw[knot] (0.085, 0.5) to[out=90, in=-90] (0.5, 1.207);
        \draw[knot, ->] (0.085, 0.499) -- (0.085, 0.501);
        \draw[knot, overcross] (0,0) -- (0.45, 0.45);
        \draw[knot] (0.45, 0.45) -- (0.55, 0.55);
        \draw[knot, ->] (0.55, 0.55) -- (1,1);
        \draw[knot, overcross] (0,1) -- (0.45, 0.55);
        \draw[knot, <-] (0,1) -- (0.01, 0.99);
        \draw[knot, overcross] (0.45, 0.55) -- (0.55, 0.45);
        \draw[knot] (0.55, 0.45) -- (1,0);
}
\right)$.

\[
\tikz[scale=1.5]
{
    \node(000) at (0,0) {$
    \varphi_{
    \tikz[baseline={([yshift=-.5ex]current bounding box.center)}, scale=.75]
{
    \draw[dotted] (.5,.5) circle(0.707);
        \draw[knot, red] (0.5, 0.45) -- (0.5, 0.55);
        \draw[knot, red] (0.25, 0.8) -- (0.25, 0.7);
        \draw[knot] (0.5, 1.207) to[out=-90, in=45] (0.325, 0.825);
        \draw[knot] (0.175, 0.675) to[out=225, in=135] (0.175, 0.325);
        \draw[knot] (0.325, 0.175) to[out=-45, in=90] (0.5, -0.207);
        \draw[knot] (0,0) to (0.175, 0.175);
        \draw[knot] (0.325, 0.325) to (0.425, 0.425);
        \draw[knot] (0.575, 0.575) to (1,1);
        \draw[knot] (0,1) to (0.175, 0.825);
        \draw[knot] (0.325, 0.675) to (0.425, 0.575);
        \draw[knot] (0.575, 0.425) to (1,0);
        \draw[knot] (0.175, 0.175) to[out=45, in=-45] (0.175, 0.325);
        \draw[knot] (0.325, 0.175) to[out=135, in=-135] (0.325, 0.325);
        \draw[knot] (0.425, 0.425) to[out=45, in=135] (0.575, 0.425);
        \draw[knot] (0.425, 0.575) to[out=-45, in=-135] (0.575, 0.575);
        \draw[knot] (0.175, 0.825) to[out=-45, in=-135] (0.325, 0.825);
        \draw[knot] (0.175, 0.675) to[out=45, in=135] (0.325, 0.675);
}
    }
    \tikz[baseline={([yshift=-.5ex]current bounding box.center)}, scale=1.15]
{
    \draw[dotted] (.5,.5) circle(0.707);
        \draw[knot] (0.5, 1.207) to[out=-90, in=45] (0.325, 0.825);
        \draw[knot] (0.175, 0.675) to[out=225, in=135] (0.175, 0.325);
        \draw[knot] (0.325, 0.175) to[out=-45, in=90] (0.5, -0.207);
        \draw[knot] (0,0) to (0.175, 0.175);
        \draw[knot] (0.325, 0.325) to (0.425, 0.425);
        \draw[knot] (0.575, 0.575) to (1,1);
        \draw[knot] (0,1) to (0.175, 0.825);
        \draw[knot] (0.325, 0.675) to (0.425, 0.575);
        \draw[knot] (0.575, 0.425) to (1,0);
        \draw[knot] (0.175, 0.175) to[out=45, in=-45] (0.175, 0.325);
        \draw[knot] (0.325, 0.175) to[out=135, in=-135] (0.325, 0.325);
        \draw[knot] (0.425, 0.425) to[out=45, in=135] (0.575, 0.425);
        \draw[knot] (0.425, 0.575) to[out=-45, in=-135] (0.575, 0.575);
        \draw[knot] (0.175, 0.825) to[out=-45, in=-135] (0.325, 0.825);
        \draw[knot] (0.175, 0.675) to[out=45, in=135] (0.325, 0.675);
}
    $};
    \node(100) at (5,0) {$
        \varphi_{
    \tikz[baseline={([yshift=-.5ex]current bounding box.center)}, scale=.75]
{
    \draw[dotted] (.5,.5) circle(0.707);
        \draw[knot, red] (0.25, 0.2) -- (0.25, 0.3);
        \draw[knot, red] (0.5, 0.45) -- (0.5, 0.55);
        \draw[knot, red] (0.25, 0.8) -- (0.25, 0.7);
        \draw[knot] (0.5, 1.207) to[out=-90, in=45] (0.325, 0.825);
        \draw[knot] (0.175, 0.675) to[out=225, in=135] (0.175, 0.325);
        \draw[knot] (0.325, 0.175) to[out=-45, in=90] (0.5, -0.207);
        \draw[knot] (0,0) to (0.175, 0.175);
        \draw[knot] (0.325, 0.325) to (0.425, 0.425);
        \draw[knot] (0.575, 0.575) to (1,1);
        \draw[knot] (0,1) to (0.175, 0.825);
        \draw[knot] (0.325, 0.675) to (0.425, 0.575);
        \draw[knot] (0.575, 0.425) to (1,0);
        \draw[knot] (0.175, 0.175) to[out=45, in=135] (0.325, 0.175);
        \draw[knot] (0.175, 0.325) to[out=-45, in=-135] (0.325, 0.325);
        \draw[knot] (0.425, 0.425) to[out=45, in=135] (0.575, 0.425);
        \draw[knot] (0.425, 0.575) to[out=-45, in=-135] (0.575, 0.575);
        \draw[knot] (0.175, 0.825) to[out=-45, in=-135] (0.325, 0.825);
        \draw[knot] (0.175, 0.675) to[out=45, in=135] (0.325, 0.675);
}^{(1,1)}
\,
    }
    \tikz[baseline={([yshift=-.5ex]current bounding box.center)}, scale=1.15]
{
    \draw[dotted] (.5,.5) circle(0.707);
        \draw[knot] (0.5, 1.207) to[out=-90, in=45] (0.325, 0.825);
        \draw[knot] (0.175, 0.675) to[out=225, in=135] (0.175, 0.325);
        \draw[knot] (0.325, 0.175) to[out=-45, in=90] (0.5, -0.207);
        \draw[knot] (0,0) to (0.175, 0.175);
        \draw[knot] (0.325, 0.325) to (0.425, 0.425);
        \draw[knot] (0.575, 0.575) to (1,1);
        \draw[knot] (0,1) to (0.175, 0.825);
        \draw[knot] (0.325, 0.675) to (0.425, 0.575);
        \draw[knot] (0.575, 0.425) to (1,0);
        \draw[knot] (0.175, 0.175) to[out=45, in=135] (0.325, 0.175);
        \draw[knot] (0.175, 0.325) to[out=-45, in=-135] (0.325, 0.325);
        \draw[knot] (0.425, 0.425) to[out=45, in=135] (0.575, 0.425);
        \draw[knot] (0.425, 0.575) to[out=-45, in=-135] (0.575, 0.575);
        \draw[knot] (0.175, 0.825) to[out=-45, in=-135] (0.325, 0.825);
        \draw[knot] (0.175, 0.675) to[out=45, in=135] (0.325, 0.675);
}
    $};
    \node(010) at (2.5, -2.5) {$
        \varphi_{
    \tikz[baseline={([yshift=-.5ex]current bounding box.center)}, scale=0.75]
{
    \draw[dotted] (.5,.5) circle(0.707);
        \draw[knot, red] (0.25, 0.8) -- (0.25, 0.7);
        \draw[knot] (0.5, 1.207) to[out=-90, in=45] (0.325, 0.825);
        \draw[knot] (0.175, 0.675) to[out=225, in=135] (0.175, 0.325);
        \draw[knot] (0.325, 0.175) to[out=-45, in=90] (0.5, -0.207);
        \draw[knot] (0,0) to (0.175, 0.175);
        \draw[knot] (0.325, 0.325) to (0.425, 0.425);
        \draw[knot] (0.575, 0.575) to (1,1);
        \draw[knot] (0,1) to (0.175, 0.825);
        \draw[knot] (0.325, 0.675) to (0.425, 0.575);
        \draw[knot] (0.575, 0.425) to (1,0);
        \draw[knot] (0.175, 0.175) to[out=45, in=-45] (0.175, 0.325);
        \draw[knot] (0.325, 0.175) to[out=135, in=-135] (0.325, 0.325);
        \draw[knot] (0.425, 0.425) to[out=45, in=-45] (0.425, 0.575); 
        \draw[knot] (0.575, 0.425) to[out=135, in=-135] (0.575, 0.575);
        \draw[knot] (0.175, 0.825) to[out=-45, in=-135] (0.325, 0.825);
        \draw[knot] (0.175, 0.675) to[out=45, in=135] (0.325, 0.675);
}
    }
    \tikz[baseline={([yshift=-.5ex]current bounding box.center)}, scale=1.15]
{
    \draw[dotted] (.5,.5) circle(0.707);
        \draw[knot] (0.5, 1.207) to[out=-90, in=45] (0.325, 0.825);
        \draw[knot] (0.175, 0.675) to[out=225, in=135] (0.175, 0.325);
        \draw[knot] (0.325, 0.175) to[out=-45, in=90] (0.5, -0.207);
        \draw[knot] (0,0) to (0.175, 0.175);
        \draw[knot] (0.325, 0.325) to (0.425, 0.425);
        \draw[knot] (0.575, 0.575) to (1,1);
        \draw[knot] (0,1) to (0.175, 0.825);
        \draw[knot] (0.325, 0.675) to (0.425, 0.575);
        \draw[knot] (0.575, 0.425) to (1,0);
        \draw[knot] (0.175, 0.175) to[out=45, in=-45] (0.175, 0.325);
        \draw[knot] (0.325, 0.175) to[out=135, in=-135] (0.325, 0.325);
        \draw[knot] (0.425, 0.425) to[out=45, in=-45] (0.425, 0.575); 
        \draw[knot] (0.575, 0.425) to[out=135, in=-135] (0.575, 0.575);
        \draw[knot] (0.175, 0.825) to[out=-45, in=-135] (0.325, 0.825);
        \draw[knot] (0.175, 0.675) to[out=45, in=135] (0.325, 0.675);
}
    $};
    \node(110) at (7.5, -2.5) {$
        \varphi_{
    \tikz[baseline={([yshift=-.5ex]current bounding box.center)}, scale=.75]
{
    \draw[dotted] (.5,.5) circle(0.707);
        \draw[knot, red] (0.25, 0.2) -- (0.25, 0.3);
        \draw[knot, red] (0.25, 0.8) -- (0.25, 0.7);
        \draw[knot] (0.5, 1.207) to[out=-90, in=45] (0.325, 0.825);
        \draw[knot] (0.175, 0.675) to[out=225, in=135] (0.175, 0.325);
        \draw[knot] (0.325, 0.175) to[out=-45, in=90] (0.5, -0.207);
        \draw[knot] (0,0) to (0.175, 0.175);
        \draw[knot] (0.325, 0.325) to (0.425, 0.425);
        \draw[knot] (0.575, 0.575) to (1,1);
        \draw[knot] (0,1) to (0.175, 0.825);
        \draw[knot] (0.325, 0.675) to (0.425, 0.575);
        \draw[knot] (0.575, 0.425) to (1,0);
        \draw[knot] (0.175, 0.175) to[out=45, in=135] (0.325, 0.175);
        \draw[knot] (0.175, 0.325) to[out=-45, in=-135] (0.325, 0.325);
        \draw[knot] (0.425, 0.425) to[out=45, in=-45] (0.425, 0.575); 
        \draw[knot] (0.575, 0.425) to[out=135, in=-135] (0.575, 0.575);
        \draw[knot] (0.175, 0.825) to[out=-45, in=-135] (0.325, 0.825);
        \draw[knot] (0.175, 0.675) to[out=45, in=135] (0.325, 0.675);
}^{(1,1)}
\,
    }
    \tikz[baseline={([yshift=-.5ex]current bounding box.center)}, scale=1.15]
{
    \draw[dotted] (.5,.5) circle(0.707);
        \draw[knot] (0.5, 1.207) to[out=-90, in=45] (0.325, 0.825);
        \draw[knot] (0.175, 0.675) to[out=225, in=135] (0.175, 0.325);
        \draw[knot] (0.325, 0.175) to[out=-45, in=90] (0.5, -0.207);
        \draw[knot] (0,0) to (0.175, 0.175);
        \draw[knot] (0.325, 0.325) to (0.425, 0.425);
        \draw[knot] (0.575, 0.575) to (1,1);
        \draw[knot] (0,1) to (0.175, 0.825);
        \draw[knot] (0.325, 0.675) to (0.425, 0.575);
        \draw[knot] (0.575, 0.425) to (1,0);
        \draw[knot] (0.175, 0.175) to[out=45, in=135] (0.325, 0.175);
        \draw[knot] (0.175, 0.325) to[out=-45, in=-135] (0.325, 0.325);
        \draw[knot] (0.425, 0.425) to[out=45, in=-45] (0.425, 0.575); 
        \draw[knot] (0.575, 0.425) to[out=135, in=-135] (0.575, 0.575);
        \draw[knot] (0.175, 0.825) to[out=-45, in=-135] (0.325, 0.825);
        \draw[knot] (0.175, 0.675) to[out=45, in=135] (0.325, 0.675);
}
    $};
    \node(001) at (0,-5) {$
        \varphi_{
    \tikz[baseline={([yshift=-.5ex]current bounding box.center)}, scale=.75]
{
    \draw[dotted] (.5,.5) circle(0.707);
        \draw[knot, red] (0.5, 0.45) -- (0.5, 0.55);
        \draw[knot] (0.5, 1.207) to[out=-90, in=45] (0.325, 0.825);
        \draw[knot] (0.175, 0.675) to[out=225, in=135] (0.175, 0.325);
        \draw[knot] (0.325, 0.175) to[out=-45, in=90] (0.5, -0.207);
        \draw[knot] (0,0) to (0.175, 0.175);
        \draw[knot] (0.325, 0.325) to (0.425, 0.425);
        \draw[knot] (0.575, 0.575) to (1,1);
        \draw[knot] (0,1) to (0.175, 0.825);
        \draw[knot] (0.325, 0.675) to (0.425, 0.575);
        \draw[knot] (0.575, 0.425) to (1,0);
        \draw[knot] (0.175, 0.175) to[out=45, in=-45] (0.175, 0.325);
        \draw[knot] (0.325, 0.175) to[out=135, in=-135] (0.325, 0.325);
        \draw[knot] (0.425, 0.425) to[out=45, in=135] (0.575, 0.425);
        \draw[knot] (0.425, 0.575) to[out=-45, in=-135] (0.575, 0.575);
        \draw[knot] (0.175, 0.675) to[out=45, in=-45] (0.175, 0.825);
        \draw[knot] (0.325, 0.675) to[out=135, in=-135] (0.325, 0.825);
}
    }
    \tikz[baseline={([yshift=-.5ex]current bounding box.center)}, scale=1.15]
{
    \draw[dotted] (.5,.5) circle(0.707);
        \draw[knot] (0.5, 1.207) to[out=-90, in=45] (0.325, 0.825);
        \draw[knot] (0.175, 0.675) to[out=225, in=135] (0.175, 0.325);
        \draw[knot] (0.325, 0.175) to[out=-45, in=90] (0.5, -0.207);
        \draw[knot] (0,0) to (0.175, 0.175);
        \draw[knot] (0.325, 0.325) to (0.425, 0.425);
        \draw[knot] (0.575, 0.575) to (1,1);
        \draw[knot] (0,1) to (0.175, 0.825);
        \draw[knot] (0.325, 0.675) to (0.425, 0.575);
        \draw[knot] (0.575, 0.425) to (1,0);
        \draw[knot] (0.175, 0.175) to[out=45, in=-45] (0.175, 0.325);
        \draw[knot] (0.325, 0.175) to[out=135, in=-135] (0.325, 0.325);
        \draw[knot] (0.425, 0.425) to[out=45, in=135] (0.575, 0.425);
        \draw[knot] (0.425, 0.575) to[out=-45, in=-135] (0.575, 0.575);
        \draw[knot] (0.175, 0.675) to[out=45, in=-45] (0.175, 0.825);
        \draw[knot] (0.325, 0.675) to[out=135, in=-135] (0.325, 0.825);
}
    $};
    \node(101) at (5,-5) {$
        \varphi_{
    \tikz[baseline={([yshift=-.5ex]current bounding box.center)}, scale=.75]
{
    \draw[dotted] (.5,.5) circle(0.707);
        \draw[knot, red] (0.25, 0.2) -- (0.25, 0.3);
        \draw[knot, red] (0.5, 0.45) -- (0.5, 0.55);
        \draw[knot] (0.5, 1.207) to[out=-90, in=45] (0.325, 0.825);
        \draw[knot] (0.175, 0.675) to[out=225, in=135] (0.175, 0.325);
        \draw[knot] (0.325, 0.175) to[out=-45, in=90] (0.5, -0.207);
        \draw[knot] (0,0) to (0.175, 0.175);
        \draw[knot] (0.325, 0.325) to (0.425, 0.425);
        \draw[knot] (0.575, 0.575) to (1,1);
        \draw[knot] (0,1) to (0.175, 0.825);
        \draw[knot] (0.325, 0.675) to (0.425, 0.575);
        \draw[knot] (0.575, 0.425) to (1,0);
        \draw[knot] (0.175, 0.175) to[out=45, in=135] (0.325, 0.175);
        \draw[knot] (0.175, 0.325) to[out=-45, in=-135] (0.325, 0.325);
        \draw[knot] (0.425, 0.425) to[out=45, in=135] (0.575, 0.425);
        \draw[knot] (0.425, 0.575) to[out=-45, in=-135] (0.575, 0.575);
        \draw[knot] (0.175, 0.675) to[out=45, in=-45] (0.175, 0.825);
        \draw[knot] (0.325, 0.675) to[out=135, in=-135] (0.325, 0.825);
}^{(1,1)}
\,
    }
    \tikz[baseline={([yshift=-.5ex]current bounding box.center)}, scale=1.15]
{
    \draw[dotted] (.5,.5) circle(0.707);
        \draw[knot] (0.5, 1.207) to[out=-90, in=45] (0.325, 0.825);
        \draw[knot] (0.175, 0.675) to[out=225, in=135] (0.175, 0.325);
        \draw[knot] (0.325, 0.175) to[out=-45, in=90] (0.5, -0.207);
        \draw[knot] (0,0) to (0.175, 0.175);
        \draw[knot] (0.325, 0.325) to (0.425, 0.425);
        \draw[knot] (0.575, 0.575) to (1,1);
        \draw[knot] (0,1) to (0.175, 0.825);
        \draw[knot] (0.325, 0.675) to (0.425, 0.575);
        \draw[knot] (0.575, 0.425) to (1,0);
        \draw[knot] (0.175, 0.175) to[out=45, in=135] (0.325, 0.175);
        \draw[knot] (0.175, 0.325) to[out=-45, in=-135] (0.325, 0.325);
        \draw[knot] (0.425, 0.425) to[out=45, in=135] (0.575, 0.425);
        \draw[knot] (0.425, 0.575) to[out=-45, in=-135] (0.575, 0.575);
        \draw[knot] (0.175, 0.675) to[out=45, in=-45] (0.175, 0.825);
        \draw[knot] (0.325, 0.675) to[out=135, in=-135] (0.325, 0.825);
}
    $};
    \node(011) at (2.5, -7.5) {$
    \tikz[baseline={([yshift=-.5ex]current bounding box.center)}, scale=1.15]
{
    \draw[dotted] (.5,.5) circle(0.707);
        \draw[knot] (0.5, 1.207) to[out=-90, in=45] (0.325, 0.825);
        \draw[knot] (0.175, 0.675) to[out=225, in=135] (0.175, 0.325);
        \draw[knot] (0.325, 0.175) to[out=-45, in=90] (0.5, -0.207);
        \draw[knot] (0,0) to (0.175, 0.175);
        \draw[knot] (0.325, 0.325) to (0.425, 0.425);
        \draw[knot] (0.575, 0.575) to (1,1);
        \draw[knot] (0,1) to (0.175, 0.825);
        \draw[knot] (0.325, 0.675) to (0.425, 0.575);
        \draw[knot] (0.575, 0.425) to (1,0);
        \draw[knot] (0.175, 0.175) to[out=45, in=-45] (0.175, 0.325);
        \draw[knot] (0.325, 0.175) to[out=135, in=-135] (0.325, 0.325);
        \draw[knot] (0.425, 0.425) to[out=45, in=-45] (0.425, 0.575); 
        \draw[knot] (0.575, 0.425) to[out=135, in=-135] (0.575, 0.575);
        \draw[knot] (0.175, 0.675) to[out=45, in=-45] (0.175, 0.825);
        \draw[knot] (0.325, 0.675) to[out=135, in=-135] (0.325, 0.825);
}
    $};
    \node(111) at (7.5, -7.5) {$
        \varphi_{
    \tikz[baseline={([yshift=-.5ex]current bounding box.center)}, scale=.75]
{
    \draw[dotted] (.5,.5) circle(0.707);
        \draw[knot, red] (0.25, 0.2) -- (0.25, 0.3);
        \draw[knot] (0.5, 1.207) to[out=-90, in=45] (0.325, 0.825);
        \draw[knot] (0.175, 0.675) to[out=225, in=135] (0.175, 0.325);
        \draw[knot] (0.325, 0.175) to[out=-45, in=90] (0.5, -0.207);
        \draw[knot] (0,0) to (0.175, 0.175);
        \draw[knot] (0.325, 0.325) to (0.425, 0.425);
        \draw[knot] (0.575, 0.575) to (1,1);
        \draw[knot] (0,1) to (0.175, 0.825);
        \draw[knot] (0.325, 0.675) to (0.425, 0.575);
        \draw[knot] (0.575, 0.425) to (1,0);
        \draw[knot] (0.175, 0.175) to[out=45, in=135] (0.325, 0.175);
        \draw[knot] (0.175, 0.325) to[out=-45, in=-135] (0.325, 0.325);
        \draw[knot] (0.425, 0.425) to[out=45, in=-45] (0.425, 0.575); 
        \draw[knot] (0.575, 0.425) to[out=135, in=-135] (0.575, 0.575);
        \draw[knot] (0.175, 0.675) to[out=45, in=-45] (0.175, 0.825);
        \draw[knot] (0.325, 0.675) to[out=135, in=-135] (0.325, 0.825);
}^{(1,1)}
\,
    }
    \tikz[baseline={([yshift=-.5ex]current bounding box.center)}, scale=1.15]
{
    \draw[dotted] (.5,.5) circle(0.707);
        \draw[knot] (0.5, 1.207) to[out=-90, in=45] (0.325, 0.825);
        \draw[knot] (0.175, 0.675) to[out=225, in=135] (0.175, 0.325);
        \draw[knot] (0.325, 0.175) to[out=-45, in=90] (0.5, -0.207);
        \draw[knot] (0,0) to (0.175, 0.175);
        \draw[knot] (0.325, 0.325) to (0.425, 0.425);
        \draw[knot] (0.575, 0.575) to (1,1);
        \draw[knot] (0,1) to (0.175, 0.825);
        \draw[knot] (0.325, 0.675) to (0.425, 0.575);
        \draw[knot] (0.575, 0.425) to (1,0);
        \draw[knot] (0.175, 0.175) to[out=45, in=135] (0.325, 0.175);
        \draw[knot] (0.175, 0.325) to[out=-45, in=-135] (0.325, 0.325);
        \draw[knot] (0.425, 0.425) to[out=45, in=-45] (0.425, 0.575); 
        \draw[knot] (0.575, 0.425) to[out=135, in=-135] (0.575, 0.575);
        \draw[knot] (0.175, 0.675) to[out=45, in=-45] (0.175, 0.825);
        \draw[knot] (0.325, 0.675) to[out=135, in=-135] (0.325, 0.825);
}
    $};
    \draw[knot, ->] (000) to
    node[pos=0.5, above, arrows=-] {$
\tikz[baseline={([yshift=-.5ex]current bounding box.center)}, scale=0.6]
{
    \draw[dotted] (.5,.5) circle(0.707);
        \draw[knot, red] (0.2, 0.25) -- (0.3, 0.25);
        \draw[knot] (0.5, 1.207) to[out=-90, in=45] (0.325, 0.825);
        \draw[knot] (0.175, 0.675) to[out=225, in=135] (0.175, 0.325);
        \draw[knot] (0.325, 0.175) to[out=-45, in=90] (0.5, -0.207);
        \draw[knot] (0,0) to (0.175, 0.175);
        \draw[knot] (0.325, 0.325) to (0.425, 0.425);
        \draw[knot] (0.575, 0.575) to (1,1);
        \draw[knot] (0,1) to (0.175, 0.825);
        \draw[knot] (0.325, 0.675) to (0.425, 0.575);
        \draw[knot] (0.575, 0.425) to (1,0);
        \draw[knot] (0.175, 0.175) to[out=45, in=-45] (0.175, 0.325);
        \draw[knot] (0.325, 0.175) to[out=135, in=-135] (0.325, 0.325);
        \draw[knot] (0.425, 0.425) to[out=45, in=135] (0.575, 0.425);
        \draw[knot] (0.425, 0.575) to[out=-45, in=-135] (0.575, 0.575);
        \draw[knot] (0.175, 0.825) to[out=-45, in=-135] (0.325, 0.825);
        \draw[knot] (0.175, 0.675) to[out=45, in=135] (0.325, 0.675);
} \circ \varphi_{H_1}
    $}
    (100);
    \draw[knot, ->] (000) to
    node[pos=0.8, above=2mm, arrows=-] {$
\tikz[baseline={([yshift=-.5ex]current bounding box.center)}, scale=0.6]
{
    \draw[dotted] (.5,.5) circle(0.707);
        \draw[knot, red] (0.5, 0.45) -- (0.5, 0.55);
        \draw[knot] (0.5, 1.207) to[out=-90, in=45] (0.325, 0.825);
        \draw[knot] (0.175, 0.675) to[out=225, in=135] (0.175, 0.325);
        \draw[knot] (0.325, 0.175) to[out=-45, in=90] (0.5, -0.207);
        \draw[knot] (0,0) to (0.175, 0.175);
        \draw[knot] (0.325, 0.325) to (0.425, 0.425);
        \draw[knot] (0.575, 0.575) to (1,1);
        \draw[knot] (0,1) to (0.175, 0.825);
        \draw[knot] (0.325, 0.675) to (0.425, 0.575);
        \draw[knot] (0.575, 0.425) to (1,0);
        \draw[knot] (0.175, 0.175) to[out=45, in=-45] (0.175, 0.325);
        \draw[knot] (0.325, 0.175) to[out=135, in=-135] (0.325, 0.325);
        \draw[knot] (0.425, 0.425) to[out=45, in=135] (0.575, 0.425);
        \draw[knot] (0.425, 0.575) to[out=-45, in=-135] (0.575, 0.575);
        \draw[knot] (0.175, 0.825) to[out=-45, in=-135] (0.325, 0.825);
        \draw[knot] (0.175, 0.675) to[out=45, in=135] (0.325, 0.675);
}
    $}
    (010);
    \draw[knot, ->] (000) to
    node[pos=0.5, left, arrows=-] {$
\tikz[baseline={([yshift=-.5ex]current bounding box.center)}, scale=0.6]
{
    \draw[dotted] (.5,.5) circle(0.707);
        \draw[knot, red] (0.25, 0.8) -- (0.25, 0.7);
        \draw[knot] (0.5, 1.207) to[out=-90, in=45] (0.325, 0.825);
        \draw[knot] (0.175, 0.675) to[out=225, in=135] (0.175, 0.325);
        \draw[knot] (0.325, 0.175) to[out=-45, in=90] (0.5, -0.207);
        \draw[knot] (0,0) to (0.175, 0.175);
        \draw[knot] (0.325, 0.325) to (0.425, 0.425);
        \draw[knot] (0.575, 0.575) to (1,1);
        \draw[knot] (0,1) to (0.175, 0.825);
        \draw[knot] (0.325, 0.675) to (0.425, 0.575);
        \draw[knot] (0.575, 0.425) to (1,0);
        \draw[knot] (0.175, 0.175) to[out=45, in=-45] (0.175, 0.325);
        \draw[knot] (0.325, 0.175) to[out=135, in=-135] (0.325, 0.325);
        \draw[knot] (0.425, 0.425) to[out=45, in=135] (0.575, 0.425);
        \draw[knot] (0.425, 0.575) to[out=-45, in=-135] (0.575, 0.575);
        \draw[knot] (0.175, 0.825) to[out=-45, in=-135] (0.325, 0.825);
        \draw[knot] (0.175, 0.675) to[out=45, in=135] (0.325, 0.675);
}
    $}
    (001);
    \draw[knot, ->] (100) to
    node[pos=0.8, above=2mm, arrows=-] {$
\tikz[baseline={([yshift=-.5ex]current bounding box.center)}, scale=0.6]
{
    \draw[dotted] (.5,.5) circle(0.707);
        \draw[knot, red] (0.5, 0.45) -- (0.5, 0.55);
        \draw[knot] (0.5, 1.207) to[out=-90, in=45] (0.325, 0.825);
        \draw[knot] (0.175, 0.675) to[out=225, in=135] (0.175, 0.325);
        \draw[knot] (0.325, 0.175) to[out=-45, in=90] (0.5, -0.207);
        \draw[knot] (0,0) to (0.175, 0.175);
        \draw[knot] (0.325, 0.325) to (0.425, 0.425);
        \draw[knot] (0.575, 0.575) to (1,1);
        \draw[knot] (0,1) to (0.175, 0.825);
        \draw[knot] (0.325, 0.675) to (0.425, 0.575);
        \draw[knot] (0.575, 0.425) to (1,0);
        \draw[knot] (0.175, 0.175) to[out=45, in=135] (0.325, 0.175);
        \draw[knot] (0.175, 0.325) to[out=-45, in=-135] (0.325, 0.325);
        \draw[knot] (0.425, 0.425) to[out=45, in=135] (0.575, 0.425);
        \draw[knot] (0.425, 0.575) to[out=-45, in=-135] (0.575, 0.575);
        \draw[knot] (0.175, 0.825) to[out=-45, in=-135] (0.325, 0.825);
        \draw[knot] (0.175, 0.675) to[out=45, in=135] (0.325, 0.675);
}
    $}
    (110);
    \draw[knot, ->] (100) to
    node[pos=0.15, left, arrows=-] {$
\tikz[baseline={([yshift=-.5ex]current bounding box.center)}, scale=0.6]
{
    \draw[dotted] (.5,.5) circle(0.707);
        \draw[knot, red] (0.25, 0.8) -- (0.25, 0.7);
        \draw[knot] (0.5, 1.207) to[out=-90, in=45] (0.325, 0.825);
        \draw[knot] (0.175, 0.675) to[out=225, in=135] (0.175, 0.325);
        \draw[knot] (0.325, 0.175) to[out=-45, in=90] (0.5, -0.207);
        \draw[knot] (0,0) to (0.175, 0.175);
        \draw[knot] (0.325, 0.325) to (0.425, 0.425);
        \draw[knot] (0.575, 0.575) to (1,1);
        \draw[knot] (0,1) to (0.175, 0.825);
        \draw[knot] (0.325, 0.675) to (0.425, 0.575);
        \draw[knot] (0.575, 0.425) to (1,0);
        \draw[knot] (0.175, 0.175) to[out=45, in=135] (0.325, 0.175);
        \draw[knot] (0.175, 0.325) to[out=-45, in=-135] (0.325, 0.325);
        \draw[knot] (0.425, 0.425) to[out=45, in=135] (0.575, 0.425);
        \draw[knot] (0.425, 0.575) to[out=-45, in=-135] (0.575, 0.575);
        \draw[knot] (0.175, 0.825) to[out=-45, in=-135] (0.325, 0.825);
        \draw[knot] (0.175, 0.675) to[out=45, in=135] (0.325, 0.675);
}
    $}
    (101);
    \draw[knot, ->] (001) to
    node[pos=0.28, above, arrows=-] {$
\tikz[baseline={([yshift=-.5ex]current bounding box.center)}, scale=0.6]
{
    \draw[dotted] (.5,.5) circle(0.707);
        \draw[knot, red] (0.2, 0.25) -- (0.3, 0.25);
        \draw[knot] (0.5, 1.207) to[out=-90, in=45] (0.325, 0.825);
        \draw[knot] (0.175, 0.675) to[out=225, in=135] (0.175, 0.325);
        \draw[knot] (0.325, 0.175) to[out=-45, in=90] (0.5, -0.207);
        \draw[knot] (0,0) to (0.175, 0.175);
        \draw[knot] (0.325, 0.325) to (0.425, 0.425);
        \draw[knot] (0.575, 0.575) to (1,1);
        \draw[knot] (0,1) to (0.175, 0.825);
        \draw[knot] (0.325, 0.675) to (0.425, 0.575);
        \draw[knot] (0.575, 0.425) to (1,0);
        \draw[knot] (0.175, 0.175) to[out=45, in=-45] (0.175, 0.325);
        \draw[knot] (0.325, 0.175) to[out=135, in=-135] (0.325, 0.325);
        \draw[knot] (0.425, 0.425) to[out=45, in=135] (0.575, 0.425);
        \draw[knot] (0.425, 0.575) to[out=-45, in=-135] (0.575, 0.575);
        \draw[knot] (0.175, 0.675) to[out=45, in=-45] (0.175, 0.825);
        \draw[knot] (0.325, 0.675) to[out=135, in=-135] (0.325, 0.825);
} \circ \varphi_{H_3}
    $}
    (101);
    \draw[knot, ->] (001) to
    node[pos=0.2, below=2mm, arrows=-] {$
\tikz[baseline={([yshift=-.5ex]current bounding box.center)}, scale=0.6]
{
    \draw[dotted] (.5,.5) circle(0.707);
        \draw[knot, red] (0.5, 0.45) -- (0.5, 0.55);
        \draw[knot] (0.5, 1.207) to[out=-90, in=45] (0.325, 0.825);
        \draw[knot] (0.175, 0.675) to[out=225, in=135] (0.175, 0.325);
        \draw[knot] (0.325, 0.175) to[out=-45, in=90] (0.5, -0.207);
        \draw[knot] (0,0) to (0.175, 0.175);
        \draw[knot] (0.325, 0.325) to (0.425, 0.425);
        \draw[knot] (0.575, 0.575) to (1,1);
        \draw[knot] (0,1) to (0.175, 0.825);
        \draw[knot] (0.325, 0.675) to (0.425, 0.575);
        \draw[knot] (0.575, 0.425) to (1,0);
        \draw[knot] (0.175, 0.175) to[out=45, in=-45] (0.175, 0.325);
        \draw[knot] (0.325, 0.175) to[out=135, in=-135] (0.325, 0.325);
        \draw[knot] (0.425, 0.425) to[out=45, in=135] (0.575, 0.425);
        \draw[knot] (0.425, 0.575) to[out=-45, in=-135] (0.575, 0.575);
        \draw[knot] (0.175, 0.675) to[out=45, in=-45] (0.175, 0.825);
        \draw[knot] (0.325, 0.675) to[out=135, in=-135] (0.325, 0.825);
}
    $}
    (011);
    \draw[knot, ->] (110) to
    node[pos=0.5, right, arrows=-] {$
\tikz[baseline={([yshift=-.5ex]current bounding box.center)}, scale=0.6]
{
    \draw[dotted] (.5,.5) circle(0.707);
        \draw[knot, red] (0.25, 0.8) -- (0.25, 0.7);
        \draw[knot] (0.5, 1.207) to[out=-90, in=45] (0.325, 0.825);
        \draw[knot] (0.175, 0.675) to[out=225, in=135] (0.175, 0.325);
        \draw[knot] (0.325, 0.175) to[out=-45, in=90] (0.5, -0.207);
        \draw[knot] (0,0) to (0.175, 0.175);
        \draw[knot] (0.325, 0.325) to (0.425, 0.425);
        \draw[knot] (0.575, 0.575) to (1,1);
        \draw[knot] (0,1) to (0.175, 0.825);
        \draw[knot] (0.325, 0.675) to (0.425, 0.575);
        \draw[knot] (0.575, 0.425) to (1,0);
        \draw[knot] (0.175, 0.175) to[out=45, in=135] (0.325, 0.175);
        \draw[knot] (0.175, 0.325) to[out=-45, in=-135] (0.325, 0.325);
        \draw[knot] (0.425, 0.425) to[out=45, in=-45] (0.425, 0.575); 
        \draw[knot] (0.575, 0.425) to[out=135, in=-135] (0.575, 0.575);
        \draw[knot] (0.175, 0.825) to[out=-45, in=-135] (0.325, 0.825);
        \draw[knot] (0.175, 0.675) to[out=45, in=135] (0.325, 0.675);
}
    $}
    (111);
    \draw[knot, ->] (101) to
    node[pos=0.8, above=2mm, arrows=-] {$
\tikz[baseline={([yshift=-.5ex]current bounding box.center)}, scale=0.6]
{
    \draw[dotted] (.5,.5) circle(0.707);
        \draw[knot, red] (0.5, 0.45) -- (0.5, 0.55);
        \draw[knot] (0.5, 1.207) to[out=-90, in=45] (0.325, 0.825);
        \draw[knot] (0.175, 0.675) to[out=225, in=135] (0.175, 0.325);
        \draw[knot] (0.325, 0.175) to[out=-45, in=90] (0.5, -0.207);
        \draw[knot] (0,0) to (0.175, 0.175);
        \draw[knot] (0.325, 0.325) to (0.425, 0.425);
        \draw[knot] (0.575, 0.575) to (1,1);
        \draw[knot] (0,1) to (0.175, 0.825);
        \draw[knot] (0.325, 0.675) to (0.425, 0.575);
        \draw[knot] (0.575, 0.425) to (1,0);
        \draw[knot] (0.175, 0.175) to[out=45, in=135] (0.325, 0.175);
        \draw[knot] (0.175, 0.325) to[out=-45, in=-135] (0.325, 0.325);
        \draw[knot] (0.425, 0.425) to[out=45, in=135] (0.575, 0.425);
        \draw[knot] (0.425, 0.575) to[out=-45, in=-135] (0.575, 0.575);
        \draw[knot] (0.175, 0.675) to[out=45, in=-45] (0.175, 0.825);
        \draw[knot] (0.325, 0.675) to[out=135, in=-135] (0.325, 0.825);
}
    $}
    (111);
    \draw[knot, ->] (011) to
    node[pos=0.5, above, arrows=-] {$
\tikz[baseline={([yshift=-.5ex]current bounding box.center)}, scale=0.6]
{
    \draw[dotted] (.5,.5) circle(0.707);
        \draw[knot, red] (0.2, 0.25) -- (0.3, 0.25);
        \draw[knot] (0.5, 1.207) to[out=-90, in=45] (0.325, 0.825);
        \draw[knot] (0.175, 0.675) to[out=225, in=135] (0.175, 0.325);
        \draw[knot] (0.325, 0.175) to[out=-45, in=90] (0.5, -0.207);
        \draw[knot] (0,0) to (0.175, 0.175);
        \draw[knot] (0.325, 0.325) to (0.425, 0.425);
        \draw[knot] (0.575, 0.575) to (1,1);
        \draw[knot] (0,1) to (0.175, 0.825);
        \draw[knot] (0.325, 0.675) to (0.425, 0.575);
        \draw[knot] (0.575, 0.425) to (1,0);
        \draw[knot] (0.175, 0.175) to[out=45, in=-45] (0.175, 0.325);
        \draw[knot] (0.325, 0.175) to[out=135, in=-135] (0.325, 0.325);
        \draw[knot] (0.425, 0.425) to[out=45, in=-45] (0.425, 0.575); 
        \draw[knot] (0.575, 0.425) to[out=135, in=-135] (0.575, 0.575);
        \draw[knot] (0.175, 0.675) to[out=45, in=-45] (0.175, 0.825);
        \draw[knot] (0.325, 0.675) to[out=135, in=-135] (0.325, 0.825);
} \circ \varphi_{H_4}
    $}
    (111);
    \draw[overcross] (010) to (110);
        \draw[->] (010) to
        node[pos=0.28, above, arrows=-] {$
        \tikz[baseline={([yshift=-.5ex]current bounding box.center)}, scale=0.6]
{
    \draw[dotted] (.5,.5) circle(0.707);
        \draw[knot, red] (0.2, 0.25) -- (0.3, 0.25);
        \draw[knot] (0.5, 1.207) to[out=-90, in=45] (0.325, 0.825);
        \draw[knot] (0.175, 0.675) to[out=225, in=135] (0.175, 0.325);
        \draw[knot] (0.325, 0.175) to[out=-45, in=90] (0.5, -0.207);
        \draw[knot] (0,0) to (0.175, 0.175);
        \draw[knot] (0.325, 0.325) to (0.425, 0.425);
        \draw[knot] (0.575, 0.575) to (1,1);
        \draw[knot] (0,1) to (0.175, 0.825);
        \draw[knot] (0.325, 0.675) to (0.425, 0.575);
        \draw[knot] (0.575, 0.425) to (1,0);
        \draw[knot] (0.175, 0.175) to[out=45, in=-45] (0.175, 0.325);
        \draw[knot] (0.325, 0.175) to[out=135, in=-135] (0.325, 0.325);
        \draw[knot] (0.425, 0.425) to[out=45, in=-45] (0.425, 0.575); 
        \draw[knot] (0.575, 0.425) to[out=135, in=-135] (0.575, 0.575);
        \draw[knot] (0.175, 0.825) to[out=-45, in=-135] (0.325, 0.825);
        \draw[knot] (0.175, 0.675) to[out=45, in=135] (0.325, 0.675);
    } \circ \varphi_{H_2}
    $}
(110);
    \draw[overcross] (010) to (011);
        \draw[->] (010) to 
        node[pos=0.2, right, arrows=-] {$
        \tikz[baseline={([yshift=-.5ex]current bounding box.center)}, scale=0.6]
{
    \draw[dotted] (.5,.5) circle(0.707);
        \draw[knot, red] (0.25, 0.8) -- (0.25, 0.7);
        \draw[knot] (0.5, 1.207) to[out=-90, in=45] (0.325, 0.825);
        \draw[knot] (0.175, 0.675) to[out=225, in=135] (0.175, 0.325);
        \draw[knot] (0.325, 0.175) to[out=-45, in=90] (0.5, -0.207);
        \draw[knot] (0,0) to (0.175, 0.175);
        \draw[knot] (0.325, 0.325) to (0.425, 0.425);
        \draw[knot] (0.575, 0.575) to (1,1);
        \draw[knot] (0,1) to (0.175, 0.825);
        \draw[knot] (0.325, 0.675) to (0.425, 0.575);
        \draw[knot] (0.575, 0.425) to (1,0);
        \draw[knot] (0.175, 0.175) to[out=45, in=-45] (0.175, 0.325);
        \draw[knot] (0.325, 0.175) to[out=135, in=-135] (0.325, 0.325);
        \draw[knot] (0.425, 0.425) to[out=45, in=-45] (0.425, 0.575); 
        \draw[knot] (0.575, 0.425) to[out=135, in=-135] (0.575, 0.575);
        \draw[knot] (0.175, 0.825) to[out=-45, in=-135] (0.325, 0.825);
        \draw[knot] (0.175, 0.675) to[out=45, in=135] (0.325, 0.675);
    }
    $}
        (011);
    \node[draw,dashed,fit=(110)] {};
}
\]

Eyeing the boxed vertex, we have that

\[
\varphi_{\left(
    \tikz[baseline={([yshift=-.5ex]current bounding box.center)}, scale=.75]
{
    \draw[dotted] (.5,.5) circle(0.707);
        \draw[knot, red] (0.25, 0.2) -- (0.25, 0.3);
        \draw[knot, red] (0.25, 0.8) -- (0.25, 0.7);
        \draw[knot] (0.5, 1.207) to[out=-90, in=45] (0.325, 0.825);
        \draw[knot] (0.175, 0.675) to[out=225, in=135] (0.175, 0.325);
        \draw[knot] (0.325, 0.175) to[out=-45, in=90] (0.5, -0.207);
        \draw[knot] (0,0) to (0.175, 0.175);
        \draw[knot] (0.325, 0.325) to (0.425, 0.425);
        \draw[knot] (0.575, 0.575) to (1,1);
        \draw[knot] (0,1) to (0.175, 0.825);
        \draw[knot] (0.325, 0.675) to (0.425, 0.575);
        \draw[knot] (0.575, 0.425) to (1,0);
        \draw[knot] (0.175, 0.175) to[out=45, in=135] (0.325, 0.175);
        \draw[knot] (0.175, 0.325) to[out=-45, in=-135] (0.325, 0.325);
        \draw[knot] (0.425, 0.425) to[out=45, in=-45] (0.425, 0.575); 
        \draw[knot] (0.575, 0.425) to[out=135, in=-135] (0.575, 0.575);
        \draw[knot] (0.175, 0.825) to[out=-45, in=-135] (0.325, 0.825);
        \draw[knot] (0.175, 0.675) to[out=45, in=135] (0.325, 0.675);
},\, (1,1) \right)}
~\cong~
\varphi_{\left(
    \tikz[baseline={([yshift=-.5ex]current bounding box.center)}, scale=.75]
{
    \draw[dotted] (.5,.5) circle(0.707);
        \draw[knot, red] (0.25, 0.2) -- (0.25, 0.8);
        \draw[knot] (0.5, 1.207) to[out=-90, in=45] (0.325, 0.825);
        \draw[knot] (0.325, 0.175) to[out=-45, in=90] (0.5, -0.207);
        \draw[knot] (0,0) to (0.175, 0.175);
        \draw[knot] (0.575, 0.575) to (1,1);
        \draw[knot] (0,1) to (0.175, 0.825);
        \draw[knot] (0.575, 0.425) to (1,0);
        \draw[knot] (0.175, 0.175) to[out=45, in=135] (0.325, 0.175);
        \draw[knot] (0.575, 0.425) to[out=135, in=-135] (0.575, 0.575);
        \draw[knot] (0.175, 0.825) to[out=-45, in=-135] (0.325, 0.825);
}, \, (0,1)\right)}
\]
and, moreover, the delooping isomorphism provides that
\[
\varphi_{\left(
    \tikz[baseline={([yshift=-.5ex]current bounding box.center)}, scale=.75]
{
    \draw[dotted] (.5,.5) circle(0.707);
        \draw[knot, red] (0.25, 0.2) -- (0.25, 0.8);
        \draw[knot] (0.5, 1.207) to[out=-90, in=45] (0.325, 0.825);
        \draw[knot] (0.325, 0.175) to[out=-45, in=90] (0.5, -0.207);
        \draw[knot] (0,0) to (0.175, 0.175);
        \draw[knot] (0.575, 0.575) to (1,1);
        \draw[knot] (0,1) to (0.175, 0.825);
        \draw[knot] (0.575, 0.425) to (1,0);
        \draw[knot] (0.175, 0.175) to[out=45, in=135] (0.325, 0.175);
        \draw[knot] (0.575, 0.425) to[out=135, in=-135] (0.575, 0.575);
        \draw[knot] (0.175, 0.825) to[out=-45, in=-135] (0.325, 0.825);
}, \, (0,1)\right)}
    \tikz[baseline={([yshift=-.5ex]current bounding box.center)}, scale=1]
{
    \draw[dotted] (.5,.5) circle(0.707);
        \draw[knot] (0.5, 1.207) to[out=-90, in=45] (0.325, 0.825);
        \draw[knot] (0.175, 0.675) to[out=225, in=135] (0.175, 0.325);
        \draw[knot] (0.325, 0.175) to[out=-45, in=90] (0.5, -0.207);
        \draw[knot] (0,0) to (0.175, 0.175);
        \draw[knot] (0.325, 0.325) to (0.425, 0.425);
        \draw[knot] (0.575, 0.575) to (1,1);
        \draw[knot] (0,1) to (0.175, 0.825);
        \draw[knot] (0.325, 0.675) to (0.425, 0.575);
        \draw[knot] (0.575, 0.425) to (1,0);
        \draw[knot] (0.175, 0.175) to[out=45, in=135] (0.325, 0.175);
        \draw[knot] (0.175, 0.325) to[out=-45, in=-135] (0.325, 0.325);
        \draw[knot] (0.425, 0.425) to[out=45, in=-45] (0.425, 0.575); 
        \draw[knot] (0.575, 0.425) to[out=135, in=-135] (0.575, 0.575);
        \draw[knot] (0.175, 0.825) to[out=-45, in=-135] (0.325, 0.825);
        \draw[knot] (0.175, 0.675) to[out=45, in=135] (0.325, 0.675);
}
~\cong~
\varphi_{
    \tikz[baseline={([yshift=-.5ex]current bounding box.center)}, scale=.75]
{
    \draw[dotted] (.5,.5) circle(0.707);
        \draw[knot, red] (0.25, 0.2) -- (0.25, 0.8);
        \draw[knot] (0.5, 1.207) to[out=-90, in=45] (0.325, 0.825);
        \draw[knot] (0.325, 0.175) to[out=-45, in=90] (0.5, -0.207);
        \draw[knot] (0,0) to (0.175, 0.175);
        \draw[knot] (0.575, 0.575) to (1,1);
        \draw[knot] (0,1) to (0.175, 0.825);
        \draw[knot] (0.575, 0.425) to (1,0);
        \draw[knot] (0.175, 0.175) to[out=45, in=135] (0.325, 0.175);
        \draw[knot] (0.575, 0.425) to[out=135, in=-135] (0.575, 0.575);
        \draw[knot] (0.175, 0.825) to[out=-45, in=-135] (0.325, 0.825);
}}
\tikz[baseline={([yshift=-.5ex]current bounding box.center)}, scale=1]
{
    \draw[dotted] (.5,.5) circle(0.707);
        \draw[knot] (0.5, 1.207) to[out=-90, in=45] (0.325, 0.825);
        \draw[knot] (0.325, 0.175) to[out=-45, in=90] (0.5, -0.207);
        \draw[knot] (0,0) to (0.175, 0.175);
        \draw[knot] (0.575, 0.575) to (1,1);
        \draw[knot] (0,1) to (0.175, 0.825);
        \draw[knot] (0.575, 0.425) to (1,0);
        \draw[knot] (0.175, 0.175) to[out=45, in=135] (0.325, 0.175);
        \draw[knot] (0.575, 0.425) to[out=135, in=-135] (0.575, 0.575);
        \draw[knot] (0.175, 0.825) to[out=-45, in=-135] (0.325, 0.825);
}
~\oplus~
\varphi_{\left(
    \tikz[baseline={([yshift=-.5ex]current bounding box.center)}, scale=.75]
{
    \draw[dotted] (.5,.5) circle(0.707);
        \draw[knot, red] (0.25, 0.2) -- (0.25, 0.8);
        \draw[knot] (0.5, 1.207) to[out=-90, in=45] (0.325, 0.825);
        \draw[knot] (0.325, 0.175) to[out=-45, in=90] (0.5, -0.207);
        \draw[knot] (0,0) to (0.175, 0.175);
        \draw[knot] (0.575, 0.575) to (1,1);
        \draw[knot] (0,1) to (0.175, 0.825);
        \draw[knot] (0.575, 0.425) to (1,0);
        \draw[knot] (0.175, 0.175) to[out=45, in=135] (0.325, 0.175);
        \draw[knot] (0.575, 0.425) to[out=135, in=-135] (0.575, 0.575);
        \draw[knot] (0.175, 0.825) to[out=-45, in=-135] (0.325, 0.825);
}, \, (1,1)\right)}
\tikz[baseline={([yshift=-.5ex]current bounding box.center)}, scale=1]
{
    \draw[dotted] (.5,.5) circle(0.707);
        \draw[knot] (0.5, 1.207) to[out=-90, in=45] (0.325, 0.825);
        \draw[knot] (0.325, 0.175) to[out=-45, in=90] (0.5, -0.207);
        \draw[knot] (0,0) to (0.175, 0.175);
        \draw[knot] (0.575, 0.575) to (1,1);
        \draw[knot] (0,1) to (0.175, 0.825);
        \draw[knot] (0.575, 0.425) to (1,0);
        \draw[knot] (0.175, 0.175) to[out=45, in=135] (0.325, 0.175);
        \draw[knot] (0.575, 0.425) to[out=135, in=-135] (0.575, 0.575);
        \draw[knot] (0.175, 0.825) to[out=-45, in=-135] (0.325, 0.825);
}\,.
\]
Now, we apply Gaussian elimination, so that the northwest and southeast vertices of the forward-facing face cancel with the northeast vertex which we just delooped. Here is the resulting complex.

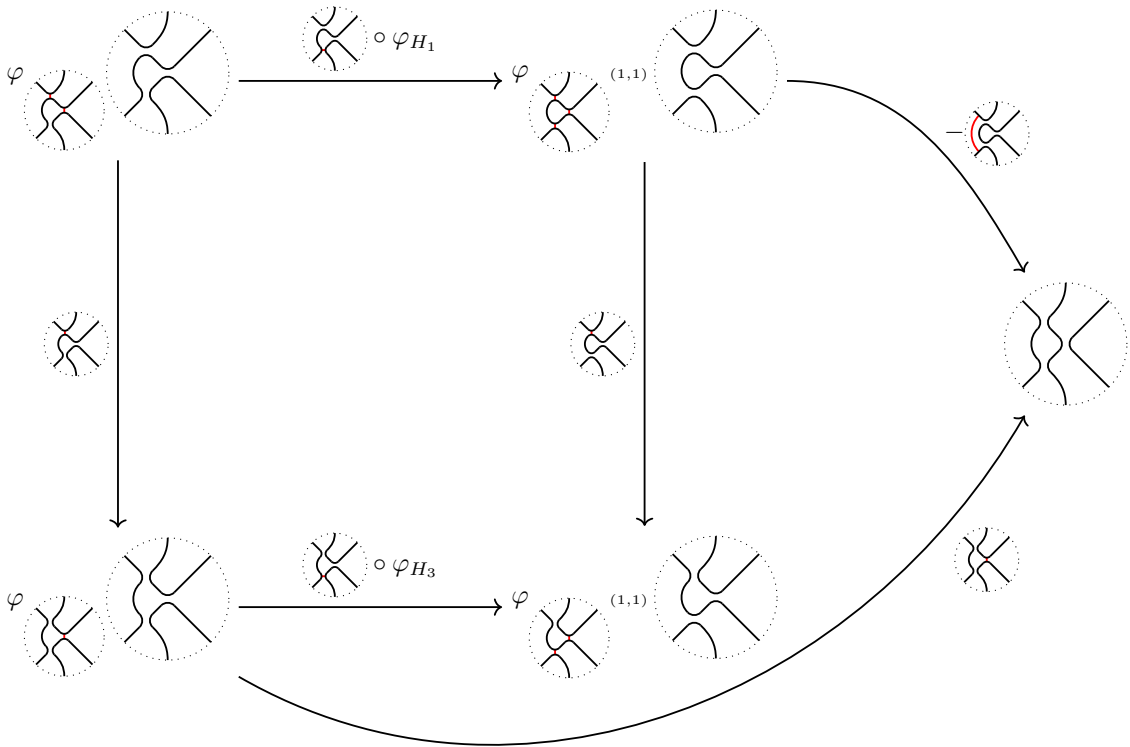
\[
\tikz[scale=1.4]
{
    \node(000) at (0,0) {$
    \varphi_{
    \tikz[baseline={([yshift=-.5ex]current bounding box.center)}, scale=.75]
{
    \draw[dotted] (.5,.5) circle(0.707);
        \draw[knot, red] (0.5, 0.45) -- (0.5, 0.55);
        \draw[knot, red] (0.25, 0.8) -- (0.25, 0.7);
        \draw[knot] (0.5, 1.207) to[out=-90, in=45] (0.325, 0.825);
        \draw[knot] (0.175, 0.675) to[out=225, in=135] (0.175, 0.325);
        \draw[knot] (0.325, 0.175) to[out=-45, in=90] (0.5, -0.207);
        \draw[knot] (0,0) to (0.175, 0.175);
        \draw[knot] (0.325, 0.325) to (0.425, 0.425);
        \draw[knot] (0.575, 0.575) to (1,1);
        \draw[knot] (0,1) to (0.175, 0.825);
        \draw[knot] (0.325, 0.675) to (0.425, 0.575);
        \draw[knot] (0.575, 0.425) to (1,0);
        \draw[knot] (0.175, 0.175) to[out=45, in=-45] (0.175, 0.325);
        \draw[knot] (0.325, 0.175) to[out=135, in=-135] (0.325, 0.325);
        \draw[knot] (0.425, 0.425) to[out=45, in=135] (0.575, 0.425);
        \draw[knot] (0.425, 0.575) to[out=-45, in=-135] (0.575, 0.575);
        \draw[knot] (0.175, 0.825) to[out=-45, in=-135] (0.325, 0.825);
        \draw[knot] (0.175, 0.675) to[out=45, in=135] (0.325, 0.675);
}
    }
    \tikz[baseline={([yshift=-.5ex]current bounding box.center)}, scale=1.15]
{
    \draw[dotted] (.5,.5) circle(0.707);
        \draw[knot] (0.5, 1.207) to[out=-90, in=45] (0.325, 0.825);
        \draw[knot] (0.175, 0.675) to[out=225, in=135] (0.175, 0.325);
        \draw[knot] (0.325, 0.175) to[out=-45, in=90] (0.5, -0.207);
        \draw[knot] (0,0) to (0.175, 0.175);
        \draw[knot] (0.325, 0.325) to (0.425, 0.425);
        \draw[knot] (0.575, 0.575) to (1,1);
        \draw[knot] (0,1) to (0.175, 0.825);
        \draw[knot] (0.325, 0.675) to (0.425, 0.575);
        \draw[knot] (0.575, 0.425) to (1,0);
        \draw[knot] (0.175, 0.175) to[out=45, in=-45] (0.175, 0.325);
        \draw[knot] (0.325, 0.175) to[out=135, in=-135] (0.325, 0.325);
        \draw[knot] (0.425, 0.425) to[out=45, in=135] (0.575, 0.425);
        \draw[knot] (0.425, 0.575) to[out=-45, in=-135] (0.575, 0.575);
        \draw[knot] (0.175, 0.825) to[out=-45, in=-135] (0.325, 0.825);
        \draw[knot] (0.175, 0.675) to[out=45, in=135] (0.325, 0.675);
}
    $};
    \node(100) at (5,0) {$
        \varphi_{
    \tikz[baseline={([yshift=-.5ex]current bounding box.center)}, scale=.75]
{
    \draw[dotted] (.5,.5) circle(0.707);
        \draw[knot, red] (0.25, 0.2) -- (0.25, 0.3);
        \draw[knot, red] (0.5, 0.45) -- (0.5, 0.55);
        \draw[knot, red] (0.25, 0.8) -- (0.25, 0.7);
        \draw[knot] (0.5, 1.207) to[out=-90, in=45] (0.325, 0.825);
        \draw[knot] (0.175, 0.675) to[out=225, in=135] (0.175, 0.325);
        \draw[knot] (0.325, 0.175) to[out=-45, in=90] (0.5, -0.207);
        \draw[knot] (0,0) to (0.175, 0.175);
        \draw[knot] (0.325, 0.325) to (0.425, 0.425);
        \draw[knot] (0.575, 0.575) to (1,1);
        \draw[knot] (0,1) to (0.175, 0.825);
        \draw[knot] (0.325, 0.675) to (0.425, 0.575);
        \draw[knot] (0.575, 0.425) to (1,0);
        \draw[knot] (0.175, 0.175) to[out=45, in=135] (0.325, 0.175);
        \draw[knot] (0.175, 0.325) to[out=-45, in=-135] (0.325, 0.325);
        \draw[knot] (0.425, 0.425) to[out=45, in=135] (0.575, 0.425);
        \draw[knot] (0.425, 0.575) to[out=-45, in=-135] (0.575, 0.575);
        \draw[knot] (0.175, 0.825) to[out=-45, in=-135] (0.325, 0.825);
        \draw[knot] (0.175, 0.675) to[out=45, in=135] (0.325, 0.675);
}^{(1,1)}
\,
    }
    \tikz[baseline={([yshift=-.5ex]current bounding box.center)}, scale=1.15]
{
    \draw[dotted] (.5,.5) circle(0.707);
        \draw[knot] (0.5, 1.207) to[out=-90, in=45] (0.325, 0.825);
        \draw[knot] (0.175, 0.675) to[out=225, in=135] (0.175, 0.325);
        \draw[knot] (0.325, 0.175) to[out=-45, in=90] (0.5, -0.207);
        \draw[knot] (0,0) to (0.175, 0.175);
        \draw[knot] (0.325, 0.325) to (0.425, 0.425);
        \draw[knot] (0.575, 0.575) to (1,1);
        \draw[knot] (0,1) to (0.175, 0.825);
        \draw[knot] (0.325, 0.675) to (0.425, 0.575);
        \draw[knot] (0.575, 0.425) to (1,0);
        \draw[knot] (0.175, 0.175) to[out=45, in=135] (0.325, 0.175);
        \draw[knot] (0.175, 0.325) to[out=-45, in=-135] (0.325, 0.325);
        \draw[knot] (0.425, 0.425) to[out=45, in=135] (0.575, 0.425);
        \draw[knot] (0.425, 0.575) to[out=-45, in=-135] (0.575, 0.575);
        \draw[knot] (0.175, 0.825) to[out=-45, in=-135] (0.325, 0.825);
        \draw[knot] (0.175, 0.675) to[out=45, in=135] (0.325, 0.675);
}
    $};
    \node(001) at (0,-5) {$
        \varphi_{
    \tikz[baseline={([yshift=-.5ex]current bounding box.center)}, scale=.75]
{
    \draw[dotted] (.5,.5) circle(0.707);
        \draw[knot, red] (0.5, 0.45) -- (0.5, 0.55);
        \draw[knot] (0.5, 1.207) to[out=-90, in=45] (0.325, 0.825);
        \draw[knot] (0.175, 0.675) to[out=225, in=135] (0.175, 0.325);
        \draw[knot] (0.325, 0.175) to[out=-45, in=90] (0.5, -0.207);
        \draw[knot] (0,0) to (0.175, 0.175);
        \draw[knot] (0.325, 0.325) to (0.425, 0.425);
        \draw[knot] (0.575, 0.575) to (1,1);
        \draw[knot] (0,1) to (0.175, 0.825);
        \draw[knot] (0.325, 0.675) to (0.425, 0.575);
        \draw[knot] (0.575, 0.425) to (1,0);
        \draw[knot] (0.175, 0.175) to[out=45, in=-45] (0.175, 0.325);
        \draw[knot] (0.325, 0.175) to[out=135, in=-135] (0.325, 0.325);
        \draw[knot] (0.425, 0.425) to[out=45, in=135] (0.575, 0.425);
        \draw[knot] (0.425, 0.575) to[out=-45, in=-135] (0.575, 0.575);
        \draw[knot] (0.175, 0.675) to[out=45, in=-45] (0.175, 0.825);
        \draw[knot] (0.325, 0.675) to[out=135, in=-135] (0.325, 0.825);
}
    }
    \tikz[baseline={([yshift=-.5ex]current bounding box.center)}, scale=1.15]
{
    \draw[dotted] (.5,.5) circle(0.707);
        \draw[knot] (0.5, 1.207) to[out=-90, in=45] (0.325, 0.825);
        \draw[knot] (0.175, 0.675) to[out=225, in=135] (0.175, 0.325);
        \draw[knot] (0.325, 0.175) to[out=-45, in=90] (0.5, -0.207);
        \draw[knot] (0,0) to (0.175, 0.175);
        \draw[knot] (0.325, 0.325) to (0.425, 0.425);
        \draw[knot] (0.575, 0.575) to (1,1);
        \draw[knot] (0,1) to (0.175, 0.825);
        \draw[knot] (0.325, 0.675) to (0.425, 0.575);
        \draw[knot] (0.575, 0.425) to (1,0);
        \draw[knot] (0.175, 0.175) to[out=45, in=-45] (0.175, 0.325);
        \draw[knot] (0.325, 0.175) to[out=135, in=-135] (0.325, 0.325);
        \draw[knot] (0.425, 0.425) to[out=45, in=135] (0.575, 0.425);
        \draw[knot] (0.425, 0.575) to[out=-45, in=-135] (0.575, 0.575);
        \draw[knot] (0.175, 0.675) to[out=45, in=-45] (0.175, 0.825);
        \draw[knot] (0.325, 0.675) to[out=135, in=-135] (0.325, 0.825);
}
    $};
    \node(101) at (5,-5) {$
        \varphi_{
    \tikz[baseline={([yshift=-.5ex]current bounding box.center)}, scale=.75]
{
    \draw[dotted] (.5,.5) circle(0.707);
        \draw[knot, red] (0.25, 0.2) -- (0.25, 0.3);
        \draw[knot, red] (0.5, 0.45) -- (0.5, 0.55);
        \draw[knot] (0.5, 1.207) to[out=-90, in=45] (0.325, 0.825);
        \draw[knot] (0.175, 0.675) to[out=225, in=135] (0.175, 0.325);
        \draw[knot] (0.325, 0.175) to[out=-45, in=90] (0.5, -0.207);
        \draw[knot] (0,0) to (0.175, 0.175);
        \draw[knot] (0.325, 0.325) to (0.425, 0.425);
        \draw[knot] (0.575, 0.575) to (1,1);
        \draw[knot] (0,1) to (0.175, 0.825);
        \draw[knot] (0.325, 0.675) to (0.425, 0.575);
        \draw[knot] (0.575, 0.425) to (1,0);
        \draw[knot] (0.175, 0.175) to[out=45, in=135] (0.325, 0.175);
        \draw[knot] (0.175, 0.325) to[out=-45, in=-135] (0.325, 0.325);
        \draw[knot] (0.425, 0.425) to[out=45, in=135] (0.575, 0.425);
        \draw[knot] (0.425, 0.575) to[out=-45, in=-135] (0.575, 0.575);
        \draw[knot] (0.175, 0.675) to[out=45, in=-45] (0.175, 0.825);
        \draw[knot] (0.325, 0.675) to[out=135, in=-135] (0.325, 0.825);
}^{(1,1)}
\,
    }
    \tikz[baseline={([yshift=-.5ex]current bounding box.center)}, scale=1.15]
{
    \draw[dotted] (.5,.5) circle(0.707);
        \draw[knot] (0.5, 1.207) to[out=-90, in=45] (0.325, 0.825);
        \draw[knot] (0.175, 0.675) to[out=225, in=135] (0.175, 0.325);
        \draw[knot] (0.325, 0.175) to[out=-45, in=90] (0.5, -0.207);
        \draw[knot] (0,0) to (0.175, 0.175);
        \draw[knot] (0.325, 0.325) to (0.425, 0.425);
        \draw[knot] (0.575, 0.575) to (1,1);
        \draw[knot] (0,1) to (0.175, 0.825);
        \draw[knot] (0.325, 0.675) to (0.425, 0.575);
        \draw[knot] (0.575, 0.425) to (1,0);
        \draw[knot] (0.175, 0.175) to[out=45, in=135] (0.325, 0.175);
        \draw[knot] (0.175, 0.325) to[out=-45, in=-135] (0.325, 0.325);
        \draw[knot] (0.425, 0.425) to[out=45, in=135] (0.575, 0.425);
        \draw[knot] (0.425, 0.575) to[out=-45, in=-135] (0.575, 0.575);
        \draw[knot] (0.175, 0.675) to[out=45, in=-45] (0.175, 0.825);
        \draw[knot] (0.325, 0.675) to[out=135, in=-135] (0.325, 0.825);
}
    $};
    \node(011) at (9, -2.5) {$
    \tikz[baseline={([yshift=-.5ex]current bounding box.center)}, scale=1.15]
{
    \draw[dotted] (.5,.5) circle(0.707);
        \draw[knot] (0.5, 1.207) to[out=-90, in=45] (0.325, 0.825);
        \draw[knot] (0.175, 0.675) to[out=225, in=135] (0.175, 0.325);
        \draw[knot] (0.325, 0.175) to[out=-45, in=90] (0.5, -0.207);
        \draw[knot] (0,0) to (0.175, 0.175);
        \draw[knot] (0.325, 0.325) to (0.425, 0.425);
        \draw[knot] (0.575, 0.575) to (1,1);
        \draw[knot] (0,1) to (0.175, 0.825);
        \draw[knot] (0.325, 0.675) to (0.425, 0.575);
        \draw[knot] (0.575, 0.425) to (1,0);
        \draw[knot] (0.175, 0.175) to[out=45, in=-45] (0.175, 0.325);
        \draw[knot] (0.325, 0.175) to[out=135, in=-135] (0.325, 0.325);
        \draw[knot] (0.425, 0.425) to[out=45, in=-45] (0.425, 0.575); 
        \draw[knot] (0.575, 0.425) to[out=135, in=-135] (0.575, 0.575);
        \draw[knot] (0.175, 0.675) to[out=45, in=-45] (0.175, 0.825);
        \draw[knot] (0.325, 0.675) to[out=135, in=-135] (0.325, 0.825);
}
    $};
    \draw[knot, ->] (000) to
    node[pos=0.5, above, arrows=-] {$
\tikz[baseline={([yshift=-.5ex]current bounding box.center)}, scale=0.6]
{
    \draw[dotted] (.5,.5) circle(0.707);
        \draw[knot, red] (0.2, 0.25) -- (0.3, 0.25);
        \draw[knot] (0.5, 1.207) to[out=-90, in=45] (0.325, 0.825);
        \draw[knot] (0.175, 0.675) to[out=225, in=135] (0.175, 0.325);
        \draw[knot] (0.325, 0.175) to[out=-45, in=90] (0.5, -0.207);
        \draw[knot] (0,0) to (0.175, 0.175);
        \draw[knot] (0.325, 0.325) to (0.425, 0.425);
        \draw[knot] (0.575, 0.575) to (1,1);
        \draw[knot] (0,1) to (0.175, 0.825);
        \draw[knot] (0.325, 0.675) to (0.425, 0.575);
        \draw[knot] (0.575, 0.425) to (1,0);
        \draw[knot] (0.175, 0.175) to[out=45, in=-45] (0.175, 0.325);
        \draw[knot] (0.325, 0.175) to[out=135, in=-135] (0.325, 0.325);
        \draw[knot] (0.425, 0.425) to[out=45, in=135] (0.575, 0.425);
        \draw[knot] (0.425, 0.575) to[out=-45, in=-135] (0.575, 0.575);
        \draw[knot] (0.175, 0.825) to[out=-45, in=-135] (0.325, 0.825);
        \draw[knot] (0.175, 0.675) to[out=45, in=135] (0.325, 0.675);
} \circ \varphi_{H_1}
    $}
    (100);
    \draw[knot, ->] (000) to
    node[pos=0.5, left, arrows=-] {$
\tikz[baseline={([yshift=-.5ex]current bounding box.center)}, scale=0.6]
{
    \draw[dotted] (.5,.5) circle(0.707);
        \draw[knot, red] (0.25, 0.8) -- (0.25, 0.7);
        \draw[knot] (0.5, 1.207) to[out=-90, in=45] (0.325, 0.825);
        \draw[knot] (0.175, 0.675) to[out=225, in=135] (0.175, 0.325);
        \draw[knot] (0.325, 0.175) to[out=-45, in=90] (0.5, -0.207);
        \draw[knot] (0,0) to (0.175, 0.175);
        \draw[knot] (0.325, 0.325) to (0.425, 0.425);
        \draw[knot] (0.575, 0.575) to (1,1);
        \draw[knot] (0,1) to (0.175, 0.825);
        \draw[knot] (0.325, 0.675) to (0.425, 0.575);
        \draw[knot] (0.575, 0.425) to (1,0);
        \draw[knot] (0.175, 0.175) to[out=45, in=-45] (0.175, 0.325);
        \draw[knot] (0.325, 0.175) to[out=135, in=-135] (0.325, 0.325);
        \draw[knot] (0.425, 0.425) to[out=45, in=135] (0.575, 0.425);
        \draw[knot] (0.425, 0.575) to[out=-45, in=-135] (0.575, 0.575);
        \draw[knot] (0.175, 0.825) to[out=-45, in=-135] (0.325, 0.825);
        \draw[knot] (0.175, 0.675) to[out=45, in=135] (0.325, 0.675);
}
    $}
    (001);
    \draw[knot, ->] (100) to
    node[pos=0.5, left, arrows=-] {$
\tikz[baseline={([yshift=-.5ex]current bounding box.center)}, scale=0.6]
{
    \draw[dotted] (.5,.5) circle(0.707);
        \draw[knot, red] (0.25, 0.8) -- (0.25, 0.7);
        \draw[knot] (0.5, 1.207) to[out=-90, in=45] (0.325, 0.825);
        \draw[knot] (0.175, 0.675) to[out=225, in=135] (0.175, 0.325);
        \draw[knot] (0.325, 0.175) to[out=-45, in=90] (0.5, -0.207);
        \draw[knot] (0,0) to (0.175, 0.175);
        \draw[knot] (0.325, 0.325) to (0.425, 0.425);
        \draw[knot] (0.575, 0.575) to (1,1);
        \draw[knot] (0,1) to (0.175, 0.825);
        \draw[knot] (0.325, 0.675) to (0.425, 0.575);
        \draw[knot] (0.575, 0.425) to (1,0);
        \draw[knot] (0.175, 0.175) to[out=45, in=135] (0.325, 0.175);
        \draw[knot] (0.175, 0.325) to[out=-45, in=-135] (0.325, 0.325);
        \draw[knot] (0.425, 0.425) to[out=45, in=135] (0.575, 0.425);
        \draw[knot] (0.425, 0.575) to[out=-45, in=-135] (0.575, 0.575);
        \draw[knot] (0.175, 0.825) to[out=-45, in=-135] (0.325, 0.825);
        \draw[knot] (0.175, 0.675) to[out=45, in=135] (0.325, 0.675);
}
    $}
    (101);
    \draw[knot, ->] (001) to
    node[pos=0.5, above, arrows=-] {$
\tikz[baseline={([yshift=-.5ex]current bounding box.center)}, scale=0.6]
{
    \draw[dotted] (.5,.5) circle(0.707);
        \draw[knot, red] (0.2, 0.25) -- (0.3, 0.25);
        \draw[knot] (0.5, 1.207) to[out=-90, in=45] (0.325, 0.825);
        \draw[knot] (0.175, 0.675) to[out=225, in=135] (0.175, 0.325);
        \draw[knot] (0.325, 0.175) to[out=-45, in=90] (0.5, -0.207);
        \draw[knot] (0,0) to (0.175, 0.175);
        \draw[knot] (0.325, 0.325) to (0.425, 0.425);
        \draw[knot] (0.575, 0.575) to (1,1);
        \draw[knot] (0,1) to (0.175, 0.825);
        \draw[knot] (0.325, 0.675) to (0.425, 0.575);
        \draw[knot] (0.575, 0.425) to (1,0);
        \draw[knot] (0.175, 0.175) to[out=45, in=-45] (0.175, 0.325);
        \draw[knot] (0.325, 0.175) to[out=135, in=-135] (0.325, 0.325);
        \draw[knot] (0.425, 0.425) to[out=45, in=135] (0.575, 0.425);
        \draw[knot] (0.425, 0.575) to[out=-45, in=-135] (0.575, 0.575);
        \draw[knot] (0.175, 0.675) to[out=45, in=-45] (0.175, 0.825);
        \draw[knot] (0.325, 0.675) to[out=135, in=-135] (0.325, 0.825);
} \circ \varphi_{H_3}
    $}
    (101);
    \draw[knot, ->] (001) to[out=-30, in=-120] (011);
    \draw[knot, ->] (100) to[out=0, in=120] (011);
    \node at (8.25, -4.55) {$\tikz[baseline={([yshift=-.5ex]current bounding box.center)}, scale=0.6]
{
    \draw[dotted] (.5,.5) circle(0.707);
        \draw[knot, red] (0.5, 0.45) -- (0.5, 0.55);
        \draw[knot] (0.5, 1.207) to[out=-90, in=45] (0.325, 0.825);
        \draw[knot] (0.175, 0.675) to[out=225, in=135] (0.175, 0.325);
        \draw[knot] (0.325, 0.175) to[out=-45, in=90] (0.5, -0.207);
        \draw[knot] (0,0) to (0.175, 0.175);
        \draw[knot] (0.325, 0.325) to (0.425, 0.425);
        \draw[knot] (0.575, 0.575) to (1,1);
        \draw[knot] (0,1) to (0.175, 0.825);
        \draw[knot] (0.325, 0.675) to (0.425, 0.575);
        \draw[knot] (0.575, 0.425) to (1,0);
        \draw[knot] (0.175, 0.175) to[out=45, in=-45] (0.175, 0.325);
        \draw[knot] (0.325, 0.175) to[out=135, in=-135] (0.325, 0.325);
        \draw[knot] (0.425, 0.425) to[out=45, in=135] (0.575, 0.425);
        \draw[knot] (0.425, 0.575) to[out=-45, in=-135] (0.575, 0.575);
        \draw[knot] (0.175, 0.675) to[out=45, in=-45] (0.175, 0.825);
        \draw[knot] (0.325, 0.675) to[out=135, in=-135] (0.325, 0.825);
}
    $};
    \node at (8.25,-0.5) {$-
    \tikz[baseline={([yshift=-.5ex]current bounding box.center)}, scale=0.6]
{
    \draw[dotted] (.5,.5) circle(0.707);
        \draw[knot, red] (0.1, 0.1) to[out=135, in=-135] (0.1, 0.9);
        \draw[knot] (0.5, 1.207) to[out=-90, in=45] (0.325, 0.825);
        \draw[knot] (0.175, 0.675) to[out=225, in=135] (0.175, 0.325);
        \draw[knot] (0.325, 0.175) to[out=-45, in=90] (0.5, -0.207);
        \draw[knot] (0,0) to (0.175, 0.175);
        \draw[knot] (0.325, 0.325) to (0.425, 0.425);
        \draw[knot] (0.575, 0.575) to (1,1);
        \draw[knot] (0,1) to (0.175, 0.825);
        \draw[knot] (0.325, 0.675) to (0.425, 0.575);
        \draw[knot] (0.575, 0.425) to (1,0);
        \draw[knot] (0.175, 0.175) to[out=45, in=135] (0.325, 0.175);
        \draw[knot] (0.175, 0.325) to[out=-45, in=-135] (0.325, 0.325);
        \draw[knot] (0.425, 0.425) to[out=45, in=135] (0.575, 0.425);
        \draw[knot] (0.425, 0.575) to[out=-45, in=-135] (0.575, 0.575);
        \draw[knot] (0.175, 0.825) to[out=-45, in=-135] (0.325, 0.825);
        \draw[knot] (0.175, 0.675) to[out=45, in=135] (0.325, 0.675);
}
    $}
}
\]

Now we do the same thing for
$\mathrm{Kh}\left(
\tikz[baseline={([yshift=-.5ex]current bounding box.center)}, scale=.85] {
        \draw[dotted] (.5,.5) circle(0.707);
        \draw[knot] (0.5, -0.207) to[out=90, in=-90] (0.915, 0.5);
        \draw[knot] (0.915, 0.5) to[out=90, in=-90] (0.5, 1.207);
        \draw[knot, ->] (0.915, 0.499) -- (0.915, 0.501);
        \draw[knot] (0,0) -- (0.45, 0.45);
        \draw[knot] (0.45, 0.45) -- (0.55, 0.55);
        \draw[knot, overcross] (0.55, 0.55) -- (1,1);
        \draw[knot, ->] (0.99, 0.99) -- (1,1);
        \draw[knot, <-] (0,1) -- (0.45, 0.55);
        \draw[knot, overcross] (0.45, 0.55) -- (0.55, 0.45);
        \draw[knot, overcross] (0.55, 0.45) -- (1,0);
        }\right)$. 
We will refrain from writing out the initial cube this time. Mirroring the previous argument---delooping and then applying Gaussian elimination to toss three of the four terms appearing in the forward-facing face---this complex is homotopy equivalent to the following.

\[
\tikz[scale=1.4]
{
    \node(000) at (0,0) {$
    \varphi_{
\tikz[baseline={([yshift=-.5ex]current bounding box.center)}, scale=0.75]
{
    \draw[dotted] (.5,.5) circle(0.707);
        \draw[knot, red] (0.75, 0.2) -- (0.75, 0.3);
        \draw[knot, red] (0.5, 0.45) -- (0.5, 0.55);
        \draw[knot] (0.5, 1.207) to[out=-90, in=135] (0.675, 0.825);
        \draw[knot] (0.825, 0.675) to[out=-45, in=45] (0.825, 0.325);
        \draw[knot] (0.675, 0.175) to[out=-135, in=90] (0.5, -0.207);
        \draw[knot] (0,0) to (0.425, 0.425);
        \draw[knot] (0.575, 0.575) -- (0.675, 0.675);
        \draw[knot] (0.825, 0.825) -- (1,1);
        \draw[knot] (0,1) -- (0.425, 0.575);
        \draw[knot] (0.575, 0.425) -- (0.675, 0.325);
        \draw[knot] (0.825, 0.175) -- (1,0);
        \draw[knot] (0.675, 0.175) to[out=45, in=135] (0.825, 0.175);
        \draw[knot] (0.675, 0.325) to[out=-45, in=-135] (0.825, 0.325);
        \draw[knot] (0.425, 0.425) to[out=45, in=135] (0.575, 0.425);
        \draw[knot] (0.425, 0.575) to[out=-45, in=-135] (0.575, 0.575);
        \draw[knot] (0.675, 0.675) to[out=45, in=-45] (0.675, 0.825);
        \draw[knot] (0.825, 0.675) to[out=135, in=-135] (0.825, 0.825);
}
}
\tikz[baseline={([yshift=-.5ex]current bounding box.center)}, scale=1.15]
{
    \draw[dotted] (.5,.5) circle(0.707);
        \draw[knot] (0.5, 1.207) to[out=-90, in=135] (0.675, 0.825);
        \draw[knot] (0.825, 0.675) to[out=-45, in=45] (0.825, 0.325);
        \draw[knot] (0.675, 0.175) to[out=-135, in=90] (0.5, -0.207);
        \draw[knot] (0,0) to (0.425, 0.425);
        \draw[knot] (0.575, 0.575) -- (0.675, 0.675);
        \draw[knot] (0.825, 0.825) -- (1,1);
        \draw[knot] (0,1) -- (0.425, 0.575);
        \draw[knot] (0.575, 0.425) -- (0.675, 0.325);
        \draw[knot] (0.825, 0.175) -- (1,0);
        \draw[knot] (0.675, 0.175) to[out=45, in=135] (0.825, 0.175);
        \draw[knot] (0.675, 0.325) to[out=-45, in=-135] (0.825, 0.325);
        \draw[knot] (0.425, 0.425) to[out=45, in=135] (0.575, 0.425);
        \draw[knot] (0.425, 0.575) to[out=-45, in=-135] (0.575, 0.575);
        \draw[knot] (0.675, 0.675) to[out=45, in=-45] (0.675, 0.825);
        \draw[knot] (0.825, 0.675) to[out=135, in=-135] (0.825, 0.825);
}
    $};
    \node(100) at (5,0) {$
    \varphi_{
\tikz[baseline={([yshift=-.5ex]current bounding box.center)}, scale=0.75]
{
    \draw[dotted] (.5,.5) circle(0.707);
        \draw[knot, red] (0.5, 0.45) -- (0.5, 0.55);
        \draw[knot] (0.5, 1.207) to[out=-90, in=135] (0.675, 0.825);
        \draw[knot] (0.825, 0.675) to[out=-45, in=45] (0.825, 0.325);
        \draw[knot] (0.675, 0.175) to[out=-135, in=90] (0.5, -0.207);
        \draw[knot] (0,0) to (0.425, 0.425);
        \draw[knot] (0.575, 0.575) -- (0.675, 0.675);
        \draw[knot] (0.825, 0.825) -- (1,1);
        \draw[knot] (0,1) -- (0.425, 0.575);
        \draw[knot] (0.575, 0.425) -- (0.675, 0.325);
        \draw[knot] (0.825, 0.175) -- (1,0);
        \draw[knot] (0.675, 0.175) to[out=45, in=-45] (0.675, 0.325);
        \draw[knot] (0.825, 0.175) to[out=135, in=-135] (0.825, 0.325);
        \draw[knot] (0.425, 0.425) to[out=45, in=135] (0.575, 0.425);
        \draw[knot] (0.425, 0.575) to[out=-45, in=-135] (0.575, 0.575);
        \draw[knot] (0.675, 0.675) to[out=45, in=-45] (0.675, 0.825);
        \draw[knot] (0.825, 0.675) to[out=135, in=-135] (0.825, 0.825);
}
}
\tikz[baseline={([yshift=-.5ex]current bounding box.center)}, scale=1.15]
{
    \draw[dotted] (.5,.5) circle(0.707);
        \draw[knot] (0.5, 1.207) to[out=-90, in=135] (0.675, 0.825);
        \draw[knot] (0.825, 0.675) to[out=-45, in=45] (0.825, 0.325);
        \draw[knot] (0.675, 0.175) to[out=-135, in=90] (0.5, -0.207);
        \draw[knot] (0,0) to (0.425, 0.425);
        \draw[knot] (0.575, 0.575) -- (0.675, 0.675);
        \draw[knot] (0.825, 0.825) -- (1,1);
        \draw[knot] (0,1) -- (0.425, 0.575);
        \draw[knot] (0.575, 0.425) -- (0.675, 0.325);
        \draw[knot] (0.825, 0.175) -- (1,0);
        \draw[knot] (0.675, 0.175) to[out=45, in=-45] (0.675, 0.325);
        \draw[knot] (0.825, 0.175) to[out=135, in=-135] (0.825, 0.325);
        \draw[knot] (0.425, 0.425) to[out=45, in=135] (0.575, 0.425);
        \draw[knot] (0.425, 0.575) to[out=-45, in=-135] (0.575, 0.575);
        \draw[knot] (0.675, 0.675) to[out=45, in=-45] (0.675, 0.825);
        \draw[knot] (0.825, 0.675) to[out=135, in=-135] (0.825, 0.825);
}
    $};
    \node(001) at (0,-5) {$\varphi_{
\tikz[baseline={([yshift=-.5ex]current bounding box.center)}, scale=0.75]
{
    \draw[dotted] (.5,.5) circle(0.707);
        \draw[knot, red] (0.75, 0.2) -- (0.75, 0.3);
        \draw[knot, red] (0.5, 0.45) -- (0.5, 0.55);
        \draw[knot, red] (0.75, 0.8) -- (0.75, 0.7);
        \draw[knot] (0.5, 1.207) to[out=-90, in=135] (0.675, 0.825);
        \draw[knot] (0.825, 0.675) to[out=-45, in=45] (0.825, 0.325);
        \draw[knot] (0.675, 0.175) to[out=-135, in=90] (0.5, -0.207);
        \draw[knot] (0,0) to (0.425, 0.425);
        \draw[knot] (0.575, 0.575) -- (0.675, 0.675);
        \draw[knot] (0.825, 0.825) -- (1,1);
        \draw[knot] (0,1) -- (0.425, 0.575);
        \draw[knot] (0.575, 0.425) -- (0.675, 0.325);
        \draw[knot] (0.825, 0.175) -- (1,0);
        \draw[knot] (0.675, 0.175) to[out=45, in=135] (0.825, 0.175);
        \draw[knot] (0.675, 0.325) to[out=-45, in=-135] (0.825, 0.325);
        \draw[knot] (0.425, 0.425) to[out=45, in=135] (0.575, 0.425);
        \draw[knot] (0.425, 0.575) to[out=-45, in=-135] (0.575, 0.575);
        \draw[knot] (0.675, 0.825) to[out=-45, in=-135] (0.825, 0.825);
        \draw[knot] (0.675, 0.675) to[out=45, in=135] (0.825, 0.675);
}^{(1,1)}
\,
}
\tikz[baseline={([yshift=-.5ex]current bounding box.center)}, scale=1.15]
{
    \draw[dotted] (.5,.5) circle(0.707);
        \draw[knot] (0.5, 1.207) to[out=-90, in=135] (0.675, 0.825);
        \draw[knot] (0.825, 0.675) to[out=-45, in=45] (0.825, 0.325);
        \draw[knot] (0.675, 0.175) to[out=-135, in=90] (0.5, -0.207);
        \draw[knot] (0,0) to (0.425, 0.425);
        \draw[knot] (0.575, 0.575) -- (0.675, 0.675);
        \draw[knot] (0.825, 0.825) -- (1,1);
        \draw[knot] (0,1) -- (0.425, 0.575);
        \draw[knot] (0.575, 0.425) -- (0.675, 0.325);
        \draw[knot] (0.825, 0.175) -- (1,0);
        \draw[knot] (0.675, 0.175) to[out=45, in=135] (0.825, 0.175);
        \draw[knot] (0.675, 0.325) to[out=-45, in=-135] (0.825, 0.325);
        \draw[knot] (0.425, 0.425) to[out=45, in=135] (0.575, 0.425);
        \draw[knot] (0.425, 0.575) to[out=-45, in=-135] (0.575, 0.575);
        \draw[knot] (0.675, 0.825) to[out=-45, in=-135] (0.825, 0.825);
        \draw[knot] (0.675, 0.675) to[out=45, in=135] (0.825, 0.675);
}
    $};
    \node(101) at (5,-5) {$
    \varphi_{
\tikz[baseline={([yshift=-.5ex]current bounding box.center)}, scale=0.75]
{
    \draw[dotted] (.5,.5) circle(0.707);
        \draw[knot, red] (0.5, 0.45) -- (0.5, 0.55);
        \draw[knot, red] (0.75, 0.8) -- (0.75, 0.7);
        \draw[knot] (0.5, 1.207) to[out=-90, in=135] (0.675, 0.825);
        \draw[knot] (0.825, 0.675) to[out=-45, in=45] (0.825, 0.325);
        \draw[knot] (0.675, 0.175) to[out=-135, in=90] (0.5, -0.207);
        \draw[knot] (0,0) to (0.425, 0.425);
        \draw[knot] (0.575, 0.575) -- (0.675, 0.675);
        \draw[knot] (0.825, 0.825) -- (1,1);
        \draw[knot] (0,1) -- (0.425, 0.575);
        \draw[knot] (0.575, 0.425) -- (0.675, 0.325);
        \draw[knot] (0.825, 0.175) -- (1,0);
        \draw[knot] (0.675, 0.175) to[out=45, in=-45] (0.675, 0.325);
        \draw[knot] (0.825, 0.175) to[out=135, in=-135] (0.825, 0.325);
        \draw[knot] (0.425, 0.425) to[out=45, in=135] (0.575, 0.425);
        \draw[knot] (0.425, 0.575) to[out=-45, in=-135] (0.575, 0.575);
        \draw[knot] (0.675, 0.825) to[out=-45, in=-135] (0.825, 0.825);
        \draw[knot] (0.675, 0.675) to[out=45, in=135] (0.825, 0.675);
}^{(1,1)}
}\,
\tikz[baseline={([yshift=-.5ex]current bounding box.center)}, scale=1.15]
{
    \draw[dotted] (.5,.5) circle(0.707);
        \draw[knot] (0.5, 1.207) to[out=-90, in=135] (0.675, 0.825);
        \draw[knot] (0.825, 0.675) to[out=-45, in=45] (0.825, 0.325);
        \draw[knot] (0.675, 0.175) to[out=-135, in=90] (0.5, -0.207);
        \draw[knot] (0,0) to (0.425, 0.425);
        \draw[knot] (0.575, 0.575) -- (0.675, 0.675);
        \draw[knot] (0.825, 0.825) -- (1,1);
        \draw[knot] (0,1) -- (0.425, 0.575);
        \draw[knot] (0.575, 0.425) -- (0.675, 0.325);
        \draw[knot] (0.825, 0.175) -- (1,0);
        \draw[knot] (0.675, 0.175) to[out=45, in=-45] (0.675, 0.325);
        \draw[knot] (0.825, 0.175) to[out=135, in=-135] (0.825, 0.325);
        \draw[knot] (0.425, 0.425) to[out=45, in=135] (0.575, 0.425);
        \draw[knot] (0.425, 0.575) to[out=-45, in=-135] (0.575, 0.575);
        \draw[knot] (0.675, 0.825) to[out=-45, in=-135] (0.825, 0.825);
        \draw[knot] (0.675, 0.675) to[out=45, in=135] (0.825, 0.675);
}
    $};
    \node(spe) at (9, -2.5) {$
\tikz[baseline={([yshift=-.5ex]current bounding box.center)}, scale=1.15]
{
    \draw[dotted] (.5,.5) circle(0.707);
        \draw[knot] (0.5, 1.207) to[out=-90, in=135] (0.675, 0.825);
        \draw[knot] (0.825, 0.675) to[out=-45, in=45] (0.825, 0.325);
        \draw[knot] (0.675, 0.175) to[out=-135, in=90] (0.5, -0.207);
        \draw[knot] (0,0) to (0.425, 0.425);
        \draw[knot] (0.575, 0.575) -- (0.675, 0.675);
        \draw[knot] (0.825, 0.825) -- (1,1);
        \draw[knot] (0,1) -- (0.425, 0.575);
        \draw[knot] (0.575, 0.425) -- (0.675, 0.325);
        \draw[knot] (0.825, 0.175) -- (1,0);
        \draw[knot] (0.675, 0.175) to[out=45, in=-45] (0.675, 0.325);
        \draw[knot] (0.825, 0.175) to[out=135, in=-135] (0.825, 0.325);
        \draw[knot] (0.425, 0.425) to[out=45, in=-45] (0.425, 0.575); 
        \draw[knot] (0.575, 0.425) to[out=135, in=-135] (0.575, 0.575);
        \draw[knot] (0.675, 0.675) to[out=45, in=-45] (0.675, 0.825);
        \draw[knot] (0.825, 0.675) to[out=135, in=-135] (0.825, 0.825);
}
    $};
    \draw[knot, ->] (000) to
    node[pos=0.5, above, arrows=-] {$
    \tikz[baseline={([yshift=-.5ex]current bounding box.center)}, scale=0.6]
{
    \draw[dotted] (.5,.5) circle(0.707);
        \draw[knot, red] (0.75, 0.2) -- (0.75, 0.3);
        \draw[knot] (0.5, 1.207) to[out=-90, in=135] (0.675, 0.825);
        \draw[knot] (0.825, 0.675) to[out=-45, in=45] (0.825, 0.325);
        \draw[knot] (0.675, 0.175) to[out=-135, in=90] (0.5, -0.207);
        \draw[knot] (0,0) to (0.425, 0.425);
        \draw[knot] (0.575, 0.575) -- (0.675, 0.675);
        \draw[knot] (0.825, 0.825) -- (1,1);
        \draw[knot] (0,1) -- (0.425, 0.575);
        \draw[knot] (0.575, 0.425) -- (0.675, 0.325);
        \draw[knot] (0.825, 0.175) -- (1,0);
        \draw[knot] (0.675, 0.175) to[out=45, in=135] (0.825, 0.175);
        \draw[knot] (0.675, 0.325) to[out=-45, in=-135] (0.825, 0.325);
        \draw[knot] (0.425, 0.425) to[out=45, in=135] (0.575, 0.425);
        \draw[knot] (0.425, 0.575) to[out=-45, in=-135] (0.575, 0.575);
        \draw[knot] (0.675, 0.675) to[out=45, in=-45] (0.675, 0.825);
        \draw[knot] (0.825, 0.675) to[out=135, in=-135] (0.825, 0.825);
}
    $}
    (100);
    \draw[knot, ->] (000) to
    node[pos=0.5, left, arrows=-] {$
    \tikz[baseline={([yshift=-.5ex]current bounding box.center)}, scale=0.6]
{
    \draw[dotted] (.5,.5) circle(0.707);
        \draw[knot, red] (0.7, 0.75) -- (0.8, 0.75);
        \draw[knot] (0.5, 1.207) to[out=-90, in=135] (0.675, 0.825);
        \draw[knot] (0.825, 0.675) to[out=-45, in=45] (0.825, 0.325);
        \draw[knot] (0.675, 0.175) to[out=-135, in=90] (0.5, -0.207);
        \draw[knot] (0,0) to (0.425, 0.425);
        \draw[knot] (0.575, 0.575) -- (0.675, 0.675);
        \draw[knot] (0.825, 0.825) -- (1,1);
        \draw[knot] (0,1) -- (0.425, 0.575);
        \draw[knot] (0.575, 0.425) -- (0.675, 0.325);
        \draw[knot] (0.825, 0.175) -- (1,0);
        \draw[knot] (0.675, 0.175) to[out=45, in=135] (0.825, 0.175);
        \draw[knot] (0.675, 0.325) to[out=-45, in=-135] (0.825, 0.325);
        \draw[knot] (0.425, 0.425) to[out=45, in=135] (0.575, 0.425);
        \draw[knot] (0.425, 0.575) to[out=-45, in=-135] (0.575, 0.575);
        \draw[knot] (0.675, 0.675) to[out=45, in=-45] (0.675, 0.825);
        \draw[knot] (0.825, 0.675) to[out=135, in=-135] (0.825, 0.825);
} \circ \varphi_{H_1'}
    $}
    (001);
    \draw[knot, ->] (100) to
    node[pos=0.5, left, arrows=-] {$
    \tikz[baseline={([yshift=-.5ex]current bounding box.center)}, scale=0.6]
{
    \draw[dotted] (.5,.5) circle(0.707);
        \draw[knot, red] (0.7, 0.75) -- (0.8, 0.75);
        \draw[knot] (0.5, 1.207) to[out=-90, in=135] (0.675, 0.825);
        \draw[knot] (0.825, 0.675) to[out=-45, in=45] (0.825, 0.325);
        \draw[knot] (0.675, 0.175) to[out=-135, in=90] (0.5, -0.207);
        \draw[knot] (0,0) to (0.425, 0.425);
        \draw[knot] (0.575, 0.575) -- (0.675, 0.675);
        \draw[knot] (0.825, 0.825) -- (1,1);
        \draw[knot] (0,1) -- (0.425, 0.575);
        \draw[knot] (0.575, 0.425) -- (0.675, 0.325);
        \draw[knot] (0.825, 0.175) -- (1,0);
        \draw[knot] (0.675, 0.175) to[out=45, in=-45] (0.675, 0.325);
        \draw[knot] (0.825, 0.175) to[out=135, in=-135] (0.825, 0.325);
        \draw[knot] (0.425, 0.425) to[out=45, in=135] (0.575, 0.425);
        \draw[knot] (0.425, 0.575) to[out=-45, in=-135] (0.575, 0.575);
        \draw[knot] (0.675, 0.675) to[out=45, in=-45] (0.675, 0.825);
        \draw[knot] (0.825, 0.675) to[out=135, in=-135] (0.825, 0.825);
} \circ \varphi_{H_2'}
    $}
    (101);
    \draw[knot, ->] (001) to
    node[pos=0.5, above, arrows=-] {$
    \tikz[baseline={([yshift=-.5ex]current bounding box.center)}, scale=0.6]
{
    \draw[dotted] (.5,.5) circle(0.707);
        \draw[knot, red] (0.75, 0.2) -- (0.75, 0.3);
        \draw[knot] (0.5, 1.207) to[out=-90, in=135] (0.675, 0.825);
        \draw[knot] (0.825, 0.675) to[out=-45, in=45] (0.825, 0.325);
        \draw[knot] (0.675, 0.175) to[out=-135, in=90] (0.5, -0.207);
        \draw[knot] (0,0) to (0.425, 0.425);
        \draw[knot] (0.575, 0.575) -- (0.675, 0.675);
        \draw[knot] (0.825, 0.825) -- (1,1);
        \draw[knot] (0,1) -- (0.425, 0.575);
        \draw[knot] (0.575, 0.425) -- (0.675, 0.325);
        \draw[knot] (0.825, 0.175) -- (1,0);
        \draw[knot] (0.675, 0.175) to[out=45, in=135] (0.825, 0.175);
        \draw[knot] (0.675, 0.325) to[out=-45, in=-135] (0.825, 0.325);
        \draw[knot] (0.425, 0.425) to[out=45, in=135] (0.575, 0.425);
        \draw[knot] (0.425, 0.575) to[out=-45, in=-135] (0.575, 0.575);
        \draw[knot] (0.675, 0.825) to[out=-45, in=-135] (0.825, 0.825);
        \draw[knot] (0.675, 0.675) to[out=45, in=135] (0.825, 0.675);
}
    $}
    (101);
    \draw[knot, ->] (001) to[out=-30, in=-120] (011);
    \draw[knot, ->] (100) to[out=0, in=120] (011);
    \node at (8.25, -4.55) {$-
    \tikz[baseline={([yshift=-.5ex]current bounding box.center)}, scale=0.6]
{
    \draw[dotted] (.5,.5) circle(0.707);
        \draw[knot, red] (0.9, 0.1) to[out=45, in=-45] (0.9, 0.9);
        \draw[knot] (0.5, 1.207) to[out=-90, in=135] (0.675, 0.825);
        \draw[knot] (0.825, 0.675) to[out=-45, in=45] (0.825, 0.325);
        \draw[knot] (0.675, 0.175) to[out=-135, in=90] (0.5, -0.207);
        \draw[knot] (0,0) to (0.425, 0.425);
        \draw[knot] (0.575, 0.575) -- (0.675, 0.675);
        \draw[knot] (0.825, 0.825) -- (1,1);
        \draw[knot] (0,1) -- (0.425, 0.575);
        \draw[knot] (0.575, 0.425) -- (0.675, 0.325);
        \draw[knot] (0.825, 0.175) -- (1,0);
        \draw[knot] (0.675, 0.175) to[out=45, in=135] (0.825, 0.175);
        \draw[knot] (0.675, 0.325) to[out=-45, in=-135] (0.825, 0.325);
        \draw[knot] (0.425, 0.425) to[out=45, in=135] (0.575, 0.425);
        \draw[knot] (0.425, 0.575) to[out=-45, in=-135] (0.575, 0.575);
        \draw[knot] (0.675, 0.825) to[out=-45, in=-135] (0.825, 0.825);
        \draw[knot] (0.675, 0.675) to[out=45, in=135] (0.825, 0.675);
}
    $};
    \node at (8.25,-0.5) {$
    \tikz[baseline={([yshift=-.5ex]current bounding box.center)}, scale=0.6]
{
    \draw[dotted] (.5,.5) circle(0.707);
        \draw[knot, red] (0.5, 0.45) -- (0.5, 0.55);
        \draw[knot] (0.5, 1.207) to[out=-90, in=135] (0.675, 0.825);
        \draw[knot] (0.825, 0.675) to[out=-45, in=45] (0.825, 0.325);
        \draw[knot] (0.675, 0.175) to[out=-135, in=90] (0.5, -0.207);
        \draw[knot] (0,0) to (0.425, 0.425);
        \draw[knot] (0.575, 0.575) -- (0.675, 0.675);
        \draw[knot] (0.825, 0.825) -- (1,1);
        \draw[knot] (0,1) -- (0.425, 0.575);
        \draw[knot] (0.575, 0.425) -- (0.675, 0.325);
        \draw[knot] (0.825, 0.175) -- (1,0);
        \draw[knot] (0.675, 0.175) to[out=45, in=-45] (0.675, 0.325);
        \draw[knot] (0.825, 0.175) to[out=135, in=-135] (0.825, 0.325);
        \draw[knot] (0.425, 0.425) to[out=45, in=135] (0.575, 0.425);
        \draw[knot] (0.425, 0.575) to[out=-45, in=-135] (0.575, 0.575);
        \draw[knot] (0.675, 0.675) to[out=45, in=-45] (0.675, 0.825);
        \draw[knot] (0.825, 0.675) to[out=135, in=-135] (0.825, 0.825);
}
    $}
}
\]
Finally, notice that
\[
\varphi_{\left(
    \tikz[baseline={([yshift=-.5ex]current bounding box.center)}, scale=.75]
{
    \draw[dotted] (.5,.5) circle(0.707);
        \draw[knot, red] (0.25, 0.2) -- (0.25, 0.3);
        \draw[knot, red] (0.5, 0.45) -- (0.5, 0.55);
        \draw[knot, red] (0.25, 0.8) -- (0.25, 0.7);
        \draw[knot] (0.5, 1.207) to[out=-90, in=45] (0.325, 0.825);
        \draw[knot] (0.175, 0.675) to[out=225, in=135] (0.175, 0.325);
        \draw[knot] (0.325, 0.175) to[out=-45, in=90] (0.5, -0.207);
        \draw[knot] (0,0) to (0.175, 0.175);
        \draw[knot] (0.325, 0.325) to (0.425, 0.425);
        \draw[knot] (0.575, 0.575) to (1,1);
        \draw[knot] (0,1) to (0.175, 0.825);
        \draw[knot] (0.325, 0.675) to (0.425, 0.575);
        \draw[knot] (0.575, 0.425) to (1,0);
        \draw[knot] (0.175, 0.175) to[out=45, in=135] (0.325, 0.175);
        \draw[knot] (0.175, 0.325) to[out=-45, in=-135] (0.325, 0.325);
        \draw[knot] (0.425, 0.425) to[out=45, in=135] (0.575, 0.425);
        \draw[knot] (0.425, 0.575) to[out=-45, in=-135] (0.575, 0.575);
        \draw[knot] (0.175, 0.825) to[out=-45, in=-135] (0.325, 0.825);
        \draw[knot] (0.175, 0.675) to[out=45, in=135] (0.325, 0.675);
    }
,\, (1,1)\right)}
\cong 
\varphi_{
\tikz[baseline={([yshift=-.5ex]current bounding box.center)}, scale=0.75]
{
    \draw[dotted] (.5,.5) circle(0.707);
        \draw[knot, red] (0.1, 0.1) to[out=135, in=-135] (0.1, 0.9);
        \draw[knot] (0.5, 1.207) to[out=-90, in=45] (0.325, 0.825);
        \draw[knot] (0.175, 0.675) to[out=225, in=135] (0.175, 0.325);
        \draw[knot] (0.325, 0.175) to[out=-45, in=90] (0.5, -0.207);
        \draw[knot] (0,0) to (0.175, 0.175);
        \draw[knot] (0.325, 0.325) to (0.425, 0.425);
        \draw[knot] (0.575, 0.575) to (1,1);
        \draw[knot] (0,1) to (0.175, 0.825);
        \draw[knot] (0.325, 0.675) to (0.425, 0.575);
        \draw[knot] (0.575, 0.425) to (1,0);
        \draw[knot] (0.175, 0.175) to[out=45, in=135] (0.325, 0.175);
        \draw[knot] (0.175, 0.325) to[out=-45, in=-135] (0.325, 0.325);
        \draw[knot] (0.425, 0.425) to[out=45, in=135] (0.575, 0.425);
        \draw[knot] (0.425, 0.575) to[out=-45, in=-135] (0.575, 0.575);
        \draw[knot] (0.175, 0.825) to[out=-45, in=-135] (0.325, 0.825);
        \draw[knot] (0.175, 0.675) to[out=45, in=135] (0.325, 0.675);
}
}
\]
and
\[
\varphi_{\left(
\tikz[baseline={([yshift=-.5ex]current bounding box.center)}, scale=0.75]
{
    \draw[dotted] (.5,.5) circle(0.707);
        \draw[knot, red] (0.75, 0.2) -- (0.75, 0.3);
        \draw[knot, red] (0.5, 0.45) -- (0.5, 0.55);
        \draw[knot, red] (0.75, 0.8) -- (0.75, 0.7);
        \draw[knot] (0.5, 1.207) to[out=-90, in=135] (0.675, 0.825);
        \draw[knot] (0.825, 0.675) to[out=-45, in=45] (0.825, 0.325);
        \draw[knot] (0.675, 0.175) to[out=-135, in=90] (0.5, -0.207);
        \draw[knot] (0,0) to (0.425, 0.425);
        \draw[knot] (0.575, 0.575) -- (0.675, 0.675);
        \draw[knot] (0.825, 0.825) -- (1,1);
        \draw[knot] (0,1) -- (0.425, 0.575);
        \draw[knot] (0.575, 0.425) -- (0.675, 0.325);
        \draw[knot] (0.825, 0.175) -- (1,0);
        \draw[knot] (0.675, 0.175) to[out=45, in=135] (0.825, 0.175);
        \draw[knot] (0.675, 0.325) to[out=-45, in=-135] (0.825, 0.325);
        \draw[knot] (0.425, 0.425) to[out=45, in=135] (0.575, 0.425);
        \draw[knot] (0.425, 0.575) to[out=-45, in=-135] (0.575, 0.575);
        \draw[knot] (0.675, 0.825) to[out=-45, in=-135] (0.825, 0.825);
        \draw[knot] (0.675, 0.675) to[out=45, in=135] (0.825, 0.675);
}, \, (1,1) \right)}
\cong
\varphi_{
\tikz[baseline={([yshift=-.5ex]current bounding box.center)}, scale=0.75]
{
    \draw[dotted] (.5,.5) circle(0.707);
        \draw[knot, red] (0.9, 0.1) to[out=45, in=-45] (0.9, 0.9);
        \draw[knot] (0.5, 1.207) to[out=-90, in=135] (0.675, 0.825);
        \draw[knot] (0.825, 0.675) to[out=-45, in=45] (0.825, 0.325);
        \draw[knot] (0.675, 0.175) to[out=-135, in=90] (0.5, -0.207);
        \draw[knot] (0,0) to (0.425, 0.425);
        \draw[knot] (0.575, 0.575) -- (0.675, 0.675);
        \draw[knot] (0.825, 0.825) -- (1,1);
        \draw[knot] (0,1) -- (0.425, 0.575);
        \draw[knot] (0.575, 0.425) -- (0.675, 0.325);
        \draw[knot] (0.825, 0.175) -- (1,0);
        \draw[knot] (0.675, 0.175) to[out=45, in=135] (0.825, 0.175);
        \draw[knot] (0.675, 0.325) to[out=-45, in=-135] (0.825, 0.325);
        \draw[knot] (0.425, 0.425) to[out=45, in=135] (0.575, 0.425);
        \draw[knot] (0.425, 0.575) to[out=-45, in=-135] (0.575, 0.575);
        \draw[knot] (0.675, 0.825) to[out=-45, in=-135] (0.825, 0.825);
        \draw[knot] (0.675, 0.675) to[out=45, in=135] (0.825, 0.675);
}
}
\]
as grading shift functors. From here, it is straight forward to verify that the complexes are homotopy equivalent, showing that
\[
\mathrm{Kh}
\left(
\tikz[baseline={([yshift=-.5ex]current bounding box.center)}, scale=.85]
{
        \draw[dotted] (.5,.5) circle(0.707);
        \draw[knot] (0.5, -0.207) to[out=90, in=-90] (0.085, 0.5);
        \draw[knot] (0.085, 0.5) to[out=90, in=-90] (0.5, 1.207);
        \draw[knot, ->] (0.085, 0.499) -- (0.085, 0.501);
        \draw[knot, overcross] (0,0) -- (0.45, 0.45);
        \draw[knot] (0.45, 0.45) -- (0.55, 0.55);
        \draw[knot, ->] (0.55, 0.55) -- (1,1);
        \draw[knot, overcross] (0,1) -- (0.45, 0.55);
        \draw[knot, <-] (0,1) -- (0.01, 0.99);
        \draw[knot, overcross] (0.45, 0.55) -- (0.55, 0.45);
        \draw[knot] (0.55, 0.45) -- (1,0);
}
\right)
\cong
\mathrm{Kh}
\left(
\tikz[baseline={([yshift=-.5ex]current bounding box.center)}, scale=.85]
{
        \draw[dotted] (.5,.5) circle(0.707);
        \draw[knot] (0.5, -0.207) to[out=90, in=-90] (0.915, 0.5);
        \draw[knot] (0.915, 0.5) to[out=90, in=-90] (0.5, 1.207);
        \draw[knot, ->] (0.915, 0.499) -- (0.915, 0.501);
        \draw[knot] (0,0) -- (0.45, 0.45);
        \draw[knot] (0.45, 0.45) -- (0.55, 0.55);
        \draw[knot, overcross] (0.55, 0.55) -- (1,1);
        \draw[knot, ->] (0.99, 0.99) -- (1,1);
        \draw[knot, <-] (0,1) -- (0.45, 0.55);
        \draw[knot, overcross] (0.45, 0.55) -- (0.55, 0.45);
        \draw[knot, overcross] (0.55, 0.45) -- (1,0);
}
\right).
\]

On the other hand, working the same program for 
$
\mathrm{Kh}
\left(
\tikz[baseline={([yshift=-.5ex]current bounding box.center)}, scale=.85]
{
        \draw[dotted] (.5,.5) circle(0.707);
        \draw[knot] (0.5, -0.207) to[out=90, in=-90] (0.085, 0.5);
        \draw[knot] (0.085, 0.5) to[out=90, in=-90] (0.5, 1.207);
        \draw[knot, <-] (0.085, 0.499) -- (0.085, 0.501);
        \draw[knot, overcross] (0,0) -- (0.45, 0.45);
        \draw[knot] (0.45, 0.45) -- (0.55, 0.55);
        \draw[knot, ->] (0.55, 0.55) -- (1,1);
        \draw[knot, overcross] (0,1) -- (0.45, 0.55);
        \draw[knot, overcross] (0.45, 0.55) -- (0.55, 0.45);
        \draw[knot] (0.55, 0.45) -- (1,0);
        \draw[knot, <-] (0,1) -- (0.01, 0.99);
}
\right)
$
and
$
\mathrm{Kh}
\left(
\tikz[baseline={([yshift=-.5ex]current bounding box.center)}, scale=.85]
{
        \draw[dotted] (.5,.5) circle(0.707);
        \draw[knot] (0.5, -0.207) to[out=90, in=-90] (0.915, 0.5);
        \draw[knot] (0.915, 0.5) to[out=90, in=-90] (0.5, 1.207);
        \draw[knot, <-] (0.915, 0.499) -- (0.915, 0.501);
        \draw[knot] (0,0) -- (0.45, 0.45);
        \draw[knot] (0.45, 0.45) -- (0.55, 0.55);
        \draw[knot, overcross] (0.55, 0.55) -- (1,1);
        \draw[knot, ->] (0.99, 0.99) -- (1,1);
        \draw[knot, <-] (0,1) -- (0.45, 0.55);
        \draw[knot, overcross] (0.45, 0.55) -- (0.55, 0.45);
        \draw[knot, overcross] (0.55, 0.45) -- (1,0);
}
\right)
$
we obtain the following complexes.
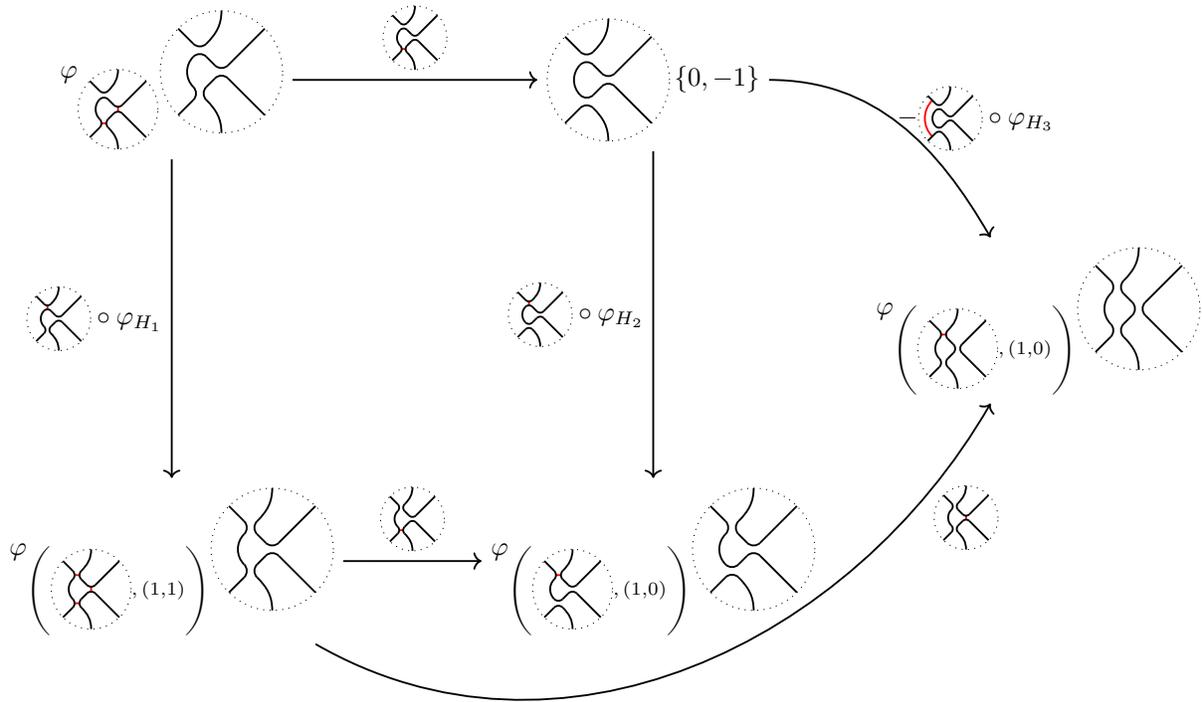
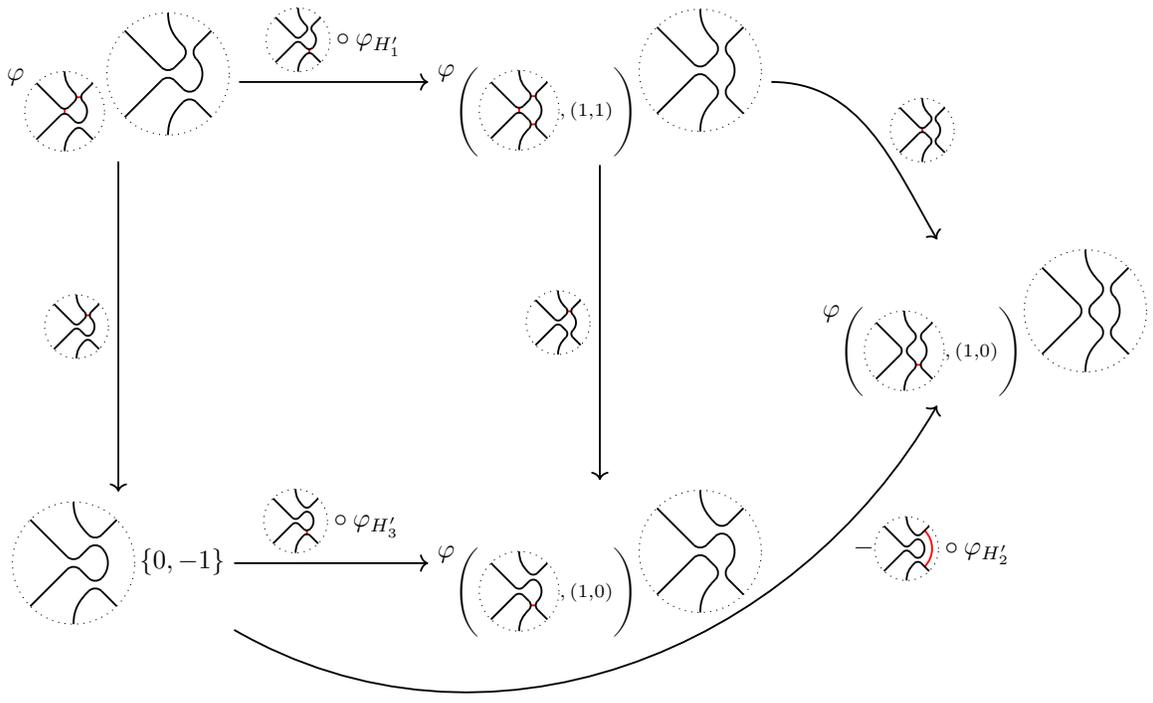
\[
\tikz[scale=1.28]
{
    \node(000) at (0,0) {$
    \varphi_{
    \tikz[baseline={([yshift=-.5ex]current bounding box.center)}, scale=.75]
{
    \draw[dotted] (.5,.5) circle(0.707);
        \draw[knot, red] (0.2, 0.25) -- (0.3, 0.25);
        \draw[knot, red] (0.5, 0.45) -- (0.5, 0.55);
        \draw[knot] (0.5, 1.207) to[out=-90, in=45] (0.325, 0.825);
        \draw[knot] (0.175, 0.675) to[out=225, in=135] (0.175, 0.325);
        \draw[knot] (0.325, 0.175) to[out=-45, in=90] (0.5, -0.207);
        \draw[knot] (0,0) to (0.175, 0.175);
        \draw[knot] (0.325, 0.325) to (0.425, 0.425);
        \draw[knot] (0.575, 0.575) to (1,1);
        \draw[knot] (0,1) to (0.175, 0.825);
        \draw[knot] (0.325, 0.675) to (0.425, 0.575);
        \draw[knot] (0.575, 0.425) to (1,0);
        \draw[knot] (0.175, 0.175) to[out=45, in=-45] (0.175, 0.325);
        \draw[knot] (0.325, 0.175) to[out=135, in=-135] (0.325, 0.325);
        \draw[knot] (0.425, 0.425) to[out=45, in=135] (0.575, 0.425);
        \draw[knot] (0.425, 0.575) to[out=-45, in=-135] (0.575, 0.575);
        \draw[knot] (0.175, 0.825) to[out=-45, in=-135] (0.325, 0.825);
        \draw[knot] (0.175, 0.675) to[out=45, in=135] (0.325, 0.675);
}
    }
    \tikz[baseline={([yshift=-.5ex]current bounding box.center)}, scale=1.15]
{
    \draw[dotted] (.5,.5) circle(0.707);
        \draw[knot] (0.5, 1.207) to[out=-90, in=45] (0.325, 0.825);
        \draw[knot] (0.175, 0.675) to[out=225, in=135] (0.175, 0.325);
        \draw[knot] (0.325, 0.175) to[out=-45, in=90] (0.5, -0.207);
        \draw[knot] (0,0) to (0.175, 0.175);
        \draw[knot] (0.325, 0.325) to (0.425, 0.425);
        \draw[knot] (0.575, 0.575) to (1,1);
        \draw[knot] (0,1) to (0.175, 0.825);
        \draw[knot] (0.325, 0.675) to (0.425, 0.575);
        \draw[knot] (0.575, 0.425) to (1,0);
        \draw[knot] (0.175, 0.175) to[out=45, in=-45] (0.175, 0.325);
        \draw[knot] (0.325, 0.175) to[out=135, in=-135] (0.325, 0.325);
        \draw[knot] (0.425, 0.425) to[out=45, in=135] (0.575, 0.425);
        \draw[knot] (0.425, 0.575) to[out=-45, in=-135] (0.575, 0.575);
        \draw[knot] (0.175, 0.825) to[out=-45, in=-135] (0.325, 0.825);
        \draw[knot] (0.175, 0.675) to[out=45, in=135] (0.325, 0.675);
}
    $};
    \node(100) at (5,0) {$
    \tikz[baseline={([yshift=-.5ex]current bounding box.center)}, scale=1.15]
{
    \draw[dotted] (.5,.5) circle(0.707);
        \draw[knot] (0.5, 1.207) to[out=-90, in=45] (0.325, 0.825);
        \draw[knot] (0.175, 0.675) to[out=225, in=135] (0.175, 0.325);
        \draw[knot] (0.325, 0.175) to[out=-45, in=90] (0.5, -0.207);
        \draw[knot] (0,0) to (0.175, 0.175);
        \draw[knot] (0.325, 0.325) to (0.425, 0.425);
        \draw[knot] (0.575, 0.575) to (1,1);
        \draw[knot] (0,1) to (0.175, 0.825);
        \draw[knot] (0.325, 0.675) to (0.425, 0.575);
        \draw[knot] (0.575, 0.425) to (1,0);
        \draw[knot] (0.175, 0.175) to[out=45, in=135] (0.325, 0.175);
        \draw[knot] (0.175, 0.325) to[out=-45, in=-135] (0.325, 0.325);
        \draw[knot] (0.425, 0.425) to[out=45, in=135] (0.575, 0.425);
        \draw[knot] (0.425, 0.575) to[out=-45, in=-135] (0.575, 0.575);
        \draw[knot] (0.175, 0.825) to[out=-45, in=-135] (0.325, 0.825);
        \draw[knot] (0.175, 0.675) to[out=45, in=135] (0.325, 0.675);
} \,\{0, -1\}
    $};
    \node(001) at (0,-5) {$
        \varphi_{\left(
    \tikz[baseline={([yshift=-.5ex]current bounding box.center)}, scale=.75]
{
    \draw[dotted] (.5,.5) circle(0.707);
        \draw[knot, red] (0.2, 0.25) -- (0.3, 0.25);
        \draw[knot, red] (0.5, 0.45) -- (0.5, 0.55);
        \draw[knot, red] (0.2, 0.75) -- (0.3, 0.75);
        \draw[knot] (0.5, 1.207) to[out=-90, in=45] (0.325, 0.825);
        \draw[knot] (0.175, 0.675) to[out=225, in=135] (0.175, 0.325);
        \draw[knot] (0.325, 0.175) to[out=-45, in=90] (0.5, -0.207);
        \draw[knot] (0,0) to (0.175, 0.175);
        \draw[knot] (0.325, 0.325) to (0.425, 0.425);
        \draw[knot] (0.575, 0.575) to (1,1);
        \draw[knot] (0,1) to (0.175, 0.825);
        \draw[knot] (0.325, 0.675) to (0.425, 0.575);
        \draw[knot] (0.575, 0.425) to (1,0);
        \draw[knot] (0.175, 0.175) to[out=45, in=-45] (0.175, 0.325);
        \draw[knot] (0.325, 0.175) to[out=135, in=-135] (0.325, 0.325);
        \draw[knot] (0.425, 0.425) to[out=45, in=135] (0.575, 0.425);
        \draw[knot] (0.425, 0.575) to[out=-45, in=-135] (0.575, 0.575);
        \draw[knot] (0.175, 0.675) to[out=45, in=-45] (0.175, 0.825);
        \draw[knot] (0.325, 0.675) to[out=135, in=-135] (0.325, 0.825);
},\, (1,1)\right)}
    \tikz[baseline={([yshift=-.5ex]current bounding box.center)}, scale=1.15]
{
    \draw[dotted] (.5,.5) circle(0.707);
        \draw[knot] (0.5, 1.207) to[out=-90, in=45] (0.325, 0.825);
        \draw[knot] (0.175, 0.675) to[out=225, in=135] (0.175, 0.325);
        \draw[knot] (0.325, 0.175) to[out=-45, in=90] (0.5, -0.207);
        \draw[knot] (0,0) to (0.175, 0.175);
        \draw[knot] (0.325, 0.325) to (0.425, 0.425);
        \draw[knot] (0.575, 0.575) to (1,1);
        \draw[knot] (0,1) to (0.175, 0.825);
        \draw[knot] (0.325, 0.675) to (0.425, 0.575);
        \draw[knot] (0.575, 0.425) to (1,0);
        \draw[knot] (0.175, 0.175) to[out=45, in=-45] (0.175, 0.325);
        \draw[knot] (0.325, 0.175) to[out=135, in=-135] (0.325, 0.325);
        \draw[knot] (0.425, 0.425) to[out=45, in=135] (0.575, 0.425);
        \draw[knot] (0.425, 0.575) to[out=-45, in=-135] (0.575, 0.575);
        \draw[knot] (0.175, 0.675) to[out=45, in=-45] (0.175, 0.825);
        \draw[knot] (0.325, 0.675) to[out=135, in=-135] (0.325, 0.825);
}
    $};
    \node(101) at (5,-5) {$
        \varphi_{\left(
    \tikz[baseline={([yshift=-.5ex]current bounding box.center)}, scale=.75]
{
    \draw[dotted] (.5,.5) circle(0.707);
        \draw[knot, red] (0.2, 0.75) -- (0.3, 0.75);
        \draw[knot] (0.5, 1.207) to[out=-90, in=45] (0.325, 0.825);
        \draw[knot] (0.175, 0.675) to[out=225, in=135] (0.175, 0.325);
        \draw[knot] (0.325, 0.175) to[out=-45, in=90] (0.5, -0.207);
        \draw[knot] (0,0) to (0.175, 0.175);
        \draw[knot] (0.325, 0.325) to (0.425, 0.425);
        \draw[knot] (0.575, 0.575) to (1,1);
        \draw[knot] (0,1) to (0.175, 0.825);
        \draw[knot] (0.325, 0.675) to (0.425, 0.575);
        \draw[knot] (0.575, 0.425) to (1,0);
        \draw[knot] (0.175, 0.175) to[out=45, in=135] (0.325, 0.175);
        \draw[knot] (0.175, 0.325) to[out=-45, in=-135] (0.325, 0.325);
        \draw[knot] (0.425, 0.425) to[out=45, in=135] (0.575, 0.425);
        \draw[knot] (0.425, 0.575) to[out=-45, in=-135] (0.575, 0.575);
        \draw[knot] (0.175, 0.675) to[out=45, in=-45] (0.175, 0.825);
        \draw[knot] (0.325, 0.675) to[out=135, in=-135] (0.325, 0.825);
},\, (1,0)\right)} 
    \tikz[baseline={([yshift=-.5ex]current bounding box.center)}, scale=1.15]
{
    \draw[dotted] (.5,.5) circle(0.707);
        \draw[knot] (0.5, 1.207) to[out=-90, in=45] (0.325, 0.825);
        \draw[knot] (0.175, 0.675) to[out=225, in=135] (0.175, 0.325);
        \draw[knot] (0.325, 0.175) to[out=-45, in=90] (0.5, -0.207);
        \draw[knot] (0,0) to (0.175, 0.175);
        \draw[knot] (0.325, 0.325) to (0.425, 0.425);
        \draw[knot] (0.575, 0.575) to (1,1);
        \draw[knot] (0,1) to (0.175, 0.825);
        \draw[knot] (0.325, 0.675) to (0.425, 0.575);
        \draw[knot] (0.575, 0.425) to (1,0);
        \draw[knot] (0.175, 0.175) to[out=45, in=135] (0.325, 0.175);
        \draw[knot] (0.175, 0.325) to[out=-45, in=-135] (0.325, 0.325);
        \draw[knot] (0.425, 0.425) to[out=45, in=135] (0.575, 0.425);
        \draw[knot] (0.425, 0.575) to[out=-45, in=-135] (0.575, 0.575);
        \draw[knot] (0.175, 0.675) to[out=45, in=-45] (0.175, 0.825);
        \draw[knot] (0.325, 0.675) to[out=135, in=-135] (0.325, 0.825);
}
    $};
    \node(011) at (9, -2.5) {$
    \varphi_{\left(
    \tikz[baseline={([yshift=-.5ex]current bounding box.center)}, scale=0.75]
{
    \draw[dotted] (.5,.5) circle(0.707);
        \draw[knot, red] (0.2, 0.75) -- (0.3, 0.75);
        \draw[knot] (0.5, 1.207) to[out=-90, in=45] (0.325, 0.825);
        \draw[knot] (0.175, 0.675) to[out=225, in=135] (0.175, 0.325);
        \draw[knot] (0.325, 0.175) to[out=-45, in=90] (0.5, -0.207);
        \draw[knot] (0,0) to (0.175, 0.175);
        \draw[knot] (0.325, 0.325) to (0.425, 0.425);
        \draw[knot] (0.575, 0.575) to (1,1);
        \draw[knot] (0,1) to (0.175, 0.825);
        \draw[knot] (0.325, 0.675) to (0.425, 0.575);
        \draw[knot] (0.575, 0.425) to (1,0);
        \draw[knot] (0.175, 0.175) to[out=45, in=-45] (0.175, 0.325);
        \draw[knot] (0.325, 0.175) to[out=135, in=-135] (0.325, 0.325);
        \draw[knot] (0.425, 0.425) to[out=45, in=-45] (0.425, 0.575); 
        \draw[knot] (0.575, 0.425) to[out=135, in=-135] (0.575, 0.575);
        \draw[knot] (0.175, 0.675) to[out=45, in=-45] (0.175, 0.825);
        \draw[knot] (0.325, 0.675) to[out=135, in=-135] (0.325, 0.825);
}
    ,\, (1,0)
    \right)}
    \tikz[baseline={([yshift=-.5ex]current bounding box.center)}, scale=1.15]
{
    \draw[dotted] (.5,.5) circle(0.707);
        \draw[knot] (0.5, 1.207) to[out=-90, in=45] (0.325, 0.825);
        \draw[knot] (0.175, 0.675) to[out=225, in=135] (0.175, 0.325);
        \draw[knot] (0.325, 0.175) to[out=-45, in=90] (0.5, -0.207);
        \draw[knot] (0,0) to (0.175, 0.175);
        \draw[knot] (0.325, 0.325) to (0.425, 0.425);
        \draw[knot] (0.575, 0.575) to (1,1);
        \draw[knot] (0,1) to (0.175, 0.825);
        \draw[knot] (0.325, 0.675) to (0.425, 0.575);
        \draw[knot] (0.575, 0.425) to (1,0);
        \draw[knot] (0.175, 0.175) to[out=45, in=-45] (0.175, 0.325);
        \draw[knot] (0.325, 0.175) to[out=135, in=-135] (0.325, 0.325);
        \draw[knot] (0.425, 0.425) to[out=45, in=-45] (0.425, 0.575); 
        \draw[knot] (0.575, 0.425) to[out=135, in=-135] (0.575, 0.575);
        \draw[knot] (0.175, 0.675) to[out=45, in=-45] (0.175, 0.825);
        \draw[knot] (0.325, 0.675) to[out=135, in=-135] (0.325, 0.825);
}
    $};
    \draw[knot, ->] (000) to
    node[pos=0.5, above, arrows=-] {$
\tikz[baseline={([yshift=-.5ex]current bounding box.center)}, scale=0.6]
{
    \draw[dotted] (.5,.5) circle(0.707);
        \draw[knot, red] (0.2, 0.25) -- (0.3, 0.25);
        \draw[knot] (0.5, 1.207) to[out=-90, in=45] (0.325, 0.825);
        \draw[knot] (0.175, 0.675) to[out=225, in=135] (0.175, 0.325);
        \draw[knot] (0.325, 0.175) to[out=-45, in=90] (0.5, -0.207);
        \draw[knot] (0,0) to (0.175, 0.175);
        \draw[knot] (0.325, 0.325) to (0.425, 0.425);
        \draw[knot] (0.575, 0.575) to (1,1);
        \draw[knot] (0,1) to (0.175, 0.825);
        \draw[knot] (0.325, 0.675) to (0.425, 0.575);
        \draw[knot] (0.575, 0.425) to (1,0);
        \draw[knot] (0.175, 0.175) to[out=45, in=-45] (0.175, 0.325);
        \draw[knot] (0.325, 0.175) to[out=135, in=-135] (0.325, 0.325);
        \draw[knot] (0.425, 0.425) to[out=45, in=135] (0.575, 0.425);
        \draw[knot] (0.425, 0.575) to[out=-45, in=-135] (0.575, 0.575);
        \draw[knot] (0.175, 0.825) to[out=-45, in=-135] (0.325, 0.825);
        \draw[knot] (0.175, 0.675) to[out=45, in=135] (0.325, 0.675);
}
    $}
    (100);
    \draw[knot, ->] (000) to
    node[pos=0.5, left, arrows=-] {$
\tikz[baseline={([yshift=-.5ex]current bounding box.center)}, scale=0.6]
{
    \draw[dotted] (.5,.5) circle(0.707);
        \draw[knot, red] (0.25, 0.8) -- (0.25, 0.7);
        \draw[knot] (0.5, 1.207) to[out=-90, in=45] (0.325, 0.825);
        \draw[knot] (0.175, 0.675) to[out=225, in=135] (0.175, 0.325);
        \draw[knot] (0.325, 0.175) to[out=-45, in=90] (0.5, -0.207);
        \draw[knot] (0,0) to (0.175, 0.175);
        \draw[knot] (0.325, 0.325) to (0.425, 0.425);
        \draw[knot] (0.575, 0.575) to (1,1);
        \draw[knot] (0,1) to (0.175, 0.825);
        \draw[knot] (0.325, 0.675) to (0.425, 0.575);
        \draw[knot] (0.575, 0.425) to (1,0);
        \draw[knot] (0.175, 0.175) to[out=45, in=-45] (0.175, 0.325);
        \draw[knot] (0.325, 0.175) to[out=135, in=-135] (0.325, 0.325);
        \draw[knot] (0.425, 0.425) to[out=45, in=135] (0.575, 0.425);
        \draw[knot] (0.425, 0.575) to[out=-45, in=-135] (0.575, 0.575);
        \draw[knot] (0.175, 0.825) to[out=-45, in=-135] (0.325, 0.825);
        \draw[knot] (0.175, 0.675) to[out=45, in=135] (0.325, 0.675);
}\circ \varphi_{H_1}
    $}
    (001);
    \draw[knot, ->] (100) to
    node[pos=0.5, left, arrows=-] {$
\tikz[baseline={([yshift=-.5ex]current bounding box.center)}, scale=0.6]
{
    \draw[dotted] (.5,.5) circle(0.707);
        \draw[knot, red] (0.25, 0.8) -- (0.25, 0.7);
        \draw[knot] (0.5, 1.207) to[out=-90, in=45] (0.325, 0.825);
        \draw[knot] (0.175, 0.675) to[out=225, in=135] (0.175, 0.325);
        \draw[knot] (0.325, 0.175) to[out=-45, in=90] (0.5, -0.207);
        \draw[knot] (0,0) to (0.175, 0.175);
        \draw[knot] (0.325, 0.325) to (0.425, 0.425);
        \draw[knot] (0.575, 0.575) to (1,1);
        \draw[knot] (0,1) to (0.175, 0.825);
        \draw[knot] (0.325, 0.675) to (0.425, 0.575);
        \draw[knot] (0.575, 0.425) to (1,0);
        \draw[knot] (0.175, 0.175) to[out=45, in=135] (0.325, 0.175);
        \draw[knot] (0.175, 0.325) to[out=-45, in=-135] (0.325, 0.325);
        \draw[knot] (0.425, 0.425) to[out=45, in=135] (0.575, 0.425);
        \draw[knot] (0.425, 0.575) to[out=-45, in=-135] (0.575, 0.575);
        \draw[knot] (0.175, 0.825) to[out=-45, in=-135] (0.325, 0.825);
        \draw[knot] (0.175, 0.675) to[out=45, in=135] (0.325, 0.675);
} \circ \varphi_{H_2}
    $}
    (101);
    \draw[knot, ->] (001) to
    node[pos=0.5, above, arrows=-] {$
\tikz[baseline={([yshift=-.5ex]current bounding box.center)}, scale=0.6]
{
    \draw[dotted] (.5,.5) circle(0.707);
        \draw[knot, red] (0.2, 0.25) -- (0.3, 0.25);
        \draw[knot] (0.5, 1.207) to[out=-90, in=45] (0.325, 0.825);
        \draw[knot] (0.175, 0.675) to[out=225, in=135] (0.175, 0.325);
        \draw[knot] (0.325, 0.175) to[out=-45, in=90] (0.5, -0.207);
        \draw[knot] (0,0) to (0.175, 0.175);
        \draw[knot] (0.325, 0.325) to (0.425, 0.425);
        \draw[knot] (0.575, 0.575) to (1,1);
        \draw[knot] (0,1) to (0.175, 0.825);
        \draw[knot] (0.325, 0.675) to (0.425, 0.575);
        \draw[knot] (0.575, 0.425) to (1,0);
        \draw[knot] (0.175, 0.175) to[out=45, in=-45] (0.175, 0.325);
        \draw[knot] (0.325, 0.175) to[out=135, in=-135] (0.325, 0.325);
        \draw[knot] (0.425, 0.425) to[out=45, in=135] (0.575, 0.425);
        \draw[knot] (0.425, 0.575) to[out=-45, in=-135] (0.575, 0.575);
        \draw[knot] (0.175, 0.675) to[out=45, in=-45] (0.175, 0.825);
        \draw[knot] (0.325, 0.675) to[out=135, in=-135] (0.325, 0.825);
}
    $}
    (101);
    \draw[knot, ->] (001) to[out=-30, in=-120] (011);
    \draw[knot, ->] (100) to[out=0, in=120] (011);
    \node at (8.25, -4.55) {$\tikz[baseline={([yshift=-.5ex]current bounding box.center)}, scale=0.6]
{
    \draw[dotted] (.5,.5) circle(0.707);
        \draw[knot, red] (0.5, 0.45) -- (0.5, 0.55);
        \draw[knot] (0.5, 1.207) to[out=-90, in=45] (0.325, 0.825);
        \draw[knot] (0.175, 0.675) to[out=225, in=135] (0.175, 0.325);
        \draw[knot] (0.325, 0.175) to[out=-45, in=90] (0.5, -0.207);
        \draw[knot] (0,0) to (0.175, 0.175);
        \draw[knot] (0.325, 0.325) to (0.425, 0.425);
        \draw[knot] (0.575, 0.575) to (1,1);
        \draw[knot] (0,1) to (0.175, 0.825);
        \draw[knot] (0.325, 0.675) to (0.425, 0.575);
        \draw[knot] (0.575, 0.425) to (1,0);
        \draw[knot] (0.175, 0.175) to[out=45, in=-45] (0.175, 0.325);
        \draw[knot] (0.325, 0.175) to[out=135, in=-135] (0.325, 0.325);
        \draw[knot] (0.425, 0.425) to[out=45, in=135] (0.575, 0.425);
        \draw[knot] (0.425, 0.575) to[out=-45, in=-135] (0.575, 0.575);
        \draw[knot] (0.175, 0.675) to[out=45, in=-45] (0.175, 0.825);
        \draw[knot] (0.325, 0.675) to[out=135, in=-135] (0.325, 0.825);
}
    $};
    \node at (8.35,-0.4) {$-
    \tikz[baseline={([yshift=-.5ex]current bounding box.center)}, scale=0.6]
{
    \draw[dotted] (.5,.5) circle(0.707);
        \draw[knot, red] (0.1, 0.1) to[out=135, in=-135] (0.1, 0.9);
        \draw[knot] (0.5, 1.207) to[out=-90, in=45] (0.325, 0.825);
        \draw[knot] (0.175, 0.675) to[out=225, in=135] (0.175, 0.325);
        \draw[knot] (0.325, 0.175) to[out=-45, in=90] (0.5, -0.207);
        \draw[knot] (0,0) to (0.175, 0.175);
        \draw[knot] (0.325, 0.325) to (0.425, 0.425);
        \draw[knot] (0.575, 0.575) to (1,1);
        \draw[knot] (0,1) to (0.175, 0.825);
        \draw[knot] (0.325, 0.675) to (0.425, 0.575);
        \draw[knot] (0.575, 0.425) to (1,0);
        \draw[knot] (0.175, 0.175) to[out=45, in=135] (0.325, 0.175);
        \draw[knot] (0.175, 0.325) to[out=-45, in=-135] (0.325, 0.325);
        \draw[knot] (0.425, 0.425) to[out=45, in=135] (0.575, 0.425);
        \draw[knot] (0.425, 0.575) to[out=-45, in=-135] (0.575, 0.575);
        \draw[knot] (0.175, 0.825) to[out=-45, in=-135] (0.325, 0.825);
        \draw[knot] (0.175, 0.675) to[out=45, in=135] (0.325, 0.675);
} \circ \varphi_{H_3}
    $}
}
\]
\[
\tikz[scale=1.28]
{
    \node(000) at (0,0) {$
    \varphi_{
\tikz[baseline={([yshift=-.5ex]current bounding box.center)}, scale=0.75]
{
    \draw[dotted] (.5,.5) circle(0.707);
        \draw[knot, red] (0.5, 0.45) -- (0.5, 0.55);
        \draw[knot, red] (0.7, 0.75) -- (0.8, 0.75);
        \draw[knot] (0.5, 1.207) to[out=-90, in=135] (0.675, 0.825);
        \draw[knot] (0.825, 0.675) to[out=-45, in=45] (0.825, 0.325);
        \draw[knot] (0.675, 0.175) to[out=-135, in=90] (0.5, -0.207);
        \draw[knot] (0,0) to (0.425, 0.425);
        \draw[knot] (0.575, 0.575) -- (0.675, 0.675);
        \draw[knot] (0.825, 0.825) -- (1,1);
        \draw[knot] (0,1) -- (0.425, 0.575);
        \draw[knot] (0.575, 0.425) -- (0.675, 0.325);
        \draw[knot] (0.825, 0.175) -- (1,0);
        \draw[knot] (0.675, 0.175) to[out=45, in=135] (0.825, 0.175);
        \draw[knot] (0.675, 0.325) to[out=-45, in=-135] (0.825, 0.325);
        \draw[knot] (0.425, 0.425) to[out=45, in=135] (0.575, 0.425);
        \draw[knot] (0.425, 0.575) to[out=-45, in=-135] (0.575, 0.575);
        \draw[knot] (0.675, 0.675) to[out=45, in=-45] (0.675, 0.825);
        \draw[knot] (0.825, 0.675) to[out=135, in=-135] (0.825, 0.825);
}
}
\tikz[baseline={([yshift=-.5ex]current bounding box.center)}, scale=1.15]
{
    \draw[dotted] (.5,.5) circle(0.707);
        \draw[knot] (0.5, 1.207) to[out=-90, in=135] (0.675, 0.825);
        \draw[knot] (0.825, 0.675) to[out=-45, in=45] (0.825, 0.325);
        \draw[knot] (0.675, 0.175) to[out=-135, in=90] (0.5, -0.207);
        \draw[knot] (0,0) to (0.425, 0.425);
        \draw[knot] (0.575, 0.575) -- (0.675, 0.675);
        \draw[knot] (0.825, 0.825) -- (1,1);
        \draw[knot] (0,1) -- (0.425, 0.575);
        \draw[knot] (0.575, 0.425) -- (0.675, 0.325);
        \draw[knot] (0.825, 0.175) -- (1,0);
        \draw[knot] (0.675, 0.175) to[out=45, in=135] (0.825, 0.175);
        \draw[knot] (0.675, 0.325) to[out=-45, in=-135] (0.825, 0.325);
        \draw[knot] (0.425, 0.425) to[out=45, in=135] (0.575, 0.425);
        \draw[knot] (0.425, 0.575) to[out=-45, in=-135] (0.575, 0.575);
        \draw[knot] (0.675, 0.675) to[out=45, in=-45] (0.675, 0.825);
        \draw[knot] (0.825, 0.675) to[out=135, in=-135] (0.825, 0.825);
}
    $};
    \node(100) at (5,0) {$
    \varphi_{\left(
\tikz[baseline={([yshift=-.5ex]current bounding box.center)}, scale=0.75]
{
    \draw[dotted] (.5,.5) circle(0.707);
        \draw[knot, red] (0.7, 0.25) -- (0.8, 0.25);
        \draw[knot, red] (0.5, 0.45) -- (0.5, 0.55);
        \draw[knot, red] (0.7, 0.75) -- (0.8, 0.75);
        \draw[knot] (0.5, 1.207) to[out=-90, in=135] (0.675, 0.825);
        \draw[knot] (0.825, 0.675) to[out=-45, in=45] (0.825, 0.325);
        \draw[knot] (0.675, 0.175) to[out=-135, in=90] (0.5, -0.207);
        \draw[knot] (0,0) to (0.425, 0.425);
        \draw[knot] (0.575, 0.575) -- (0.675, 0.675);
        \draw[knot] (0.825, 0.825) -- (1,1);
        \draw[knot] (0,1) -- (0.425, 0.575);
        \draw[knot] (0.575, 0.425) -- (0.675, 0.325);
        \draw[knot] (0.825, 0.175) -- (1,0);
        \draw[knot] (0.675, 0.175) to[out=45, in=-45] (0.675, 0.325);
        \draw[knot] (0.825, 0.175) to[out=135, in=-135] (0.825, 0.325);
        \draw[knot] (0.425, 0.425) to[out=45, in=135] (0.575, 0.425);
        \draw[knot] (0.425, 0.575) to[out=-45, in=-135] (0.575, 0.575);
        \draw[knot] (0.675, 0.675) to[out=45, in=-45] (0.675, 0.825);
        \draw[knot] (0.825, 0.675) to[out=135, in=-135] (0.825, 0.825);
}
, \, (1,1)
\right)
}
\tikz[baseline={([yshift=-.5ex]current bounding box.center)}, scale=1.15]
{
    \draw[dotted] (.5,.5) circle(0.707);
        \draw[knot] (0.5, 1.207) to[out=-90, in=135] (0.675, 0.825);
        \draw[knot] (0.825, 0.675) to[out=-45, in=45] (0.825, 0.325);
        \draw[knot] (0.675, 0.175) to[out=-135, in=90] (0.5, -0.207);
        \draw[knot] (0,0) to (0.425, 0.425);
        \draw[knot] (0.575, 0.575) -- (0.675, 0.675);
        \draw[knot] (0.825, 0.825) -- (1,1);
        \draw[knot] (0,1) -- (0.425, 0.575);
        \draw[knot] (0.575, 0.425) -- (0.675, 0.325);
        \draw[knot] (0.825, 0.175) -- (1,0);
        \draw[knot] (0.675, 0.175) to[out=45, in=-45] (0.675, 0.325);
        \draw[knot] (0.825, 0.175) to[out=135, in=-135] (0.825, 0.325);
        \draw[knot] (0.425, 0.425) to[out=45, in=135] (0.575, 0.425);
        \draw[knot] (0.425, 0.575) to[out=-45, in=-135] (0.575, 0.575);
        \draw[knot] (0.675, 0.675) to[out=45, in=-45] (0.675, 0.825);
        \draw[knot] (0.825, 0.675) to[out=135, in=-135] (0.825, 0.825);
}
    $};
    \node(001) at (0,-5) {$
\tikz[baseline={([yshift=-.5ex]current bounding box.center)}, scale=1.15]
{
    \draw[dotted] (.5,.5) circle(0.707);
        \draw[knot] (0.5, 1.207) to[out=-90, in=135] (0.675, 0.825);
        \draw[knot] (0.825, 0.675) to[out=-45, in=45] (0.825, 0.325);
        \draw[knot] (0.675, 0.175) to[out=-135, in=90] (0.5, -0.207);
        \draw[knot] (0,0) to (0.425, 0.425);
        \draw[knot] (0.575, 0.575) -- (0.675, 0.675);
        \draw[knot] (0.825, 0.825) -- (1,1);
        \draw[knot] (0,1) -- (0.425, 0.575);
        \draw[knot] (0.575, 0.425) -- (0.675, 0.325);
        \draw[knot] (0.825, 0.175) -- (1,0);
        \draw[knot] (0.675, 0.175) to[out=45, in=135] (0.825, 0.175);
        \draw[knot] (0.675, 0.325) to[out=-45, in=-135] (0.825, 0.325);
        \draw[knot] (0.425, 0.425) to[out=45, in=135] (0.575, 0.425);
        \draw[knot] (0.425, 0.575) to[out=-45, in=-135] (0.575, 0.575);
        \draw[knot] (0.675, 0.825) to[out=-45, in=-135] (0.825, 0.825);
        \draw[knot] (0.675, 0.675) to[out=45, in=135] (0.825, 0.675);
} \, \{0,-1\}
    $};
    \node(101) at (5,-5) {$
    \varphi_{\left(
\tikz[baseline={([yshift=-.5ex]current bounding box.center)}, scale=0.75]
{
    \draw[dotted] (.5,.5) circle(0.707);
        \draw[knot, red] (0.7, 0.25) -- (0.8, 0.25);
        \draw[knot] (0.5, 1.207) to[out=-90, in=135] (0.675, 0.825);
        \draw[knot] (0.825, 0.675) to[out=-45, in=45] (0.825, 0.325);
        \draw[knot] (0.675, 0.175) to[out=-135, in=90] (0.5, -0.207);
        \draw[knot] (0,0) to (0.425, 0.425);
        \draw[knot] (0.575, 0.575) -- (0.675, 0.675);
        \draw[knot] (0.825, 0.825) -- (1,1);
        \draw[knot] (0,1) -- (0.425, 0.575);
        \draw[knot] (0.575, 0.425) -- (0.675, 0.325);
        \draw[knot] (0.825, 0.175) -- (1,0);
        \draw[knot] (0.675, 0.175) to[out=45, in=-45] (0.675, 0.325);
        \draw[knot] (0.825, 0.175) to[out=135, in=-135] (0.825, 0.325);
        \draw[knot] (0.425, 0.425) to[out=45, in=135] (0.575, 0.425);
        \draw[knot] (0.425, 0.575) to[out=-45, in=-135] (0.575, 0.575);
        \draw[knot] (0.675, 0.825) to[out=-45, in=-135] (0.825, 0.825);
        \draw[knot] (0.675, 0.675) to[out=45, in=135] (0.825, 0.675);
}, \, (1,0)
\right)
}
\tikz[baseline={([yshift=-.5ex]current bounding box.center)}, scale=1.15]
{
    \draw[dotted] (.5,.5) circle(0.707);
        \draw[knot] (0.5, 1.207) to[out=-90, in=135] (0.675, 0.825);
        \draw[knot] (0.825, 0.675) to[out=-45, in=45] (0.825, 0.325);
        \draw[knot] (0.675, 0.175) to[out=-135, in=90] (0.5, -0.207);
        \draw[knot] (0,0) to (0.425, 0.425);
        \draw[knot] (0.575, 0.575) -- (0.675, 0.675);
        \draw[knot] (0.825, 0.825) -- (1,1);
        \draw[knot] (0,1) -- (0.425, 0.575);
        \draw[knot] (0.575, 0.425) -- (0.675, 0.325);
        \draw[knot] (0.825, 0.175) -- (1,0);
        \draw[knot] (0.675, 0.175) to[out=45, in=-45] (0.675, 0.325);
        \draw[knot] (0.825, 0.175) to[out=135, in=-135] (0.825, 0.325);
        \draw[knot] (0.425, 0.425) to[out=45, in=135] (0.575, 0.425);
        \draw[knot] (0.425, 0.575) to[out=-45, in=-135] (0.575, 0.575);
        \draw[knot] (0.675, 0.825) to[out=-45, in=-135] (0.825, 0.825);
        \draw[knot] (0.675, 0.675) to[out=45, in=135] (0.825, 0.675);
}
    $};
    \node(spe) at (9, -2.5) {$
    \varphi_{\left(
\tikz[baseline={([yshift=-.5ex]current bounding box.center)}, scale=0.75]
{
    \draw[dotted] (.5,.5) circle(0.707);
        \draw[knot, red] (0.7, 0.25) -- (0.8, 0.25);
        \draw[knot] (0.5, 1.207) to[out=-90, in=135] (0.675, 0.825);
        \draw[knot] (0.825, 0.675) to[out=-45, in=45] (0.825, 0.325);
        \draw[knot] (0.675, 0.175) to[out=-135, in=90] (0.5, -0.207);
        \draw[knot] (0,0) to (0.425, 0.425);
        \draw[knot] (0.575, 0.575) -- (0.675, 0.675);
        \draw[knot] (0.825, 0.825) -- (1,1);
        \draw[knot] (0,1) -- (0.425, 0.575);
        \draw[knot] (0.575, 0.425) -- (0.675, 0.325);
        \draw[knot] (0.825, 0.175) -- (1,0);
        \draw[knot] (0.675, 0.175) to[out=45, in=-45] (0.675, 0.325);
        \draw[knot] (0.825, 0.175) to[out=135, in=-135] (0.825, 0.325);
        \draw[knot] (0.425, 0.425) to[out=45, in=-45] (0.425, 0.575); 
        \draw[knot] (0.575, 0.425) to[out=135, in=-135] (0.575, 0.575);
        \draw[knot] (0.675, 0.675) to[out=45, in=-45] (0.675, 0.825);
        \draw[knot] (0.825, 0.675) to[out=135, in=-135] (0.825, 0.825);
}
,\, (1,0)
\right)
}
\tikz[baseline={([yshift=-.5ex]current bounding box.center)}, scale=1.15]
{
    \draw[dotted] (.5,.5) circle(0.707);
        \draw[knot] (0.5, 1.207) to[out=-90, in=135] (0.675, 0.825);
        \draw[knot] (0.825, 0.675) to[out=-45, in=45] (0.825, 0.325);
        \draw[knot] (0.675, 0.175) to[out=-135, in=90] (0.5, -0.207);
        \draw[knot] (0,0) to (0.425, 0.425);
        \draw[knot] (0.575, 0.575) -- (0.675, 0.675);
        \draw[knot] (0.825, 0.825) -- (1,1);
        \draw[knot] (0,1) -- (0.425, 0.575);
        \draw[knot] (0.575, 0.425) -- (0.675, 0.325);
        \draw[knot] (0.825, 0.175) -- (1,0);
        \draw[knot] (0.675, 0.175) to[out=45, in=-45] (0.675, 0.325);
        \draw[knot] (0.825, 0.175) to[out=135, in=-135] (0.825, 0.325);
        \draw[knot] (0.425, 0.425) to[out=45, in=-45] (0.425, 0.575); 
        \draw[knot] (0.575, 0.425) to[out=135, in=-135] (0.575, 0.575);
        \draw[knot] (0.675, 0.675) to[out=45, in=-45] (0.675, 0.825);
        \draw[knot] (0.825, 0.675) to[out=135, in=-135] (0.825, 0.825);
}
    $};
    \draw[knot, ->] (000) to
    node[pos=0.5, above, arrows=-] {$
    \tikz[baseline={([yshift=-.5ex]current bounding box.center)}, scale=0.6]
{
    \draw[dotted] (.5,.5) circle(0.707);
        \draw[knot, red] (0.75, 0.2) -- (0.75, 0.3);
        \draw[knot] (0.5, 1.207) to[out=-90, in=135] (0.675, 0.825);
        \draw[knot] (0.825, 0.675) to[out=-45, in=45] (0.825, 0.325);
        \draw[knot] (0.675, 0.175) to[out=-135, in=90] (0.5, -0.207);
        \draw[knot] (0,0) to (0.425, 0.425);
        \draw[knot] (0.575, 0.575) -- (0.675, 0.675);
        \draw[knot] (0.825, 0.825) -- (1,1);
        \draw[knot] (0,1) -- (0.425, 0.575);
        \draw[knot] (0.575, 0.425) -- (0.675, 0.325);
        \draw[knot] (0.825, 0.175) -- (1,0);
        \draw[knot] (0.675, 0.175) to[out=45, in=135] (0.825, 0.175);
        \draw[knot] (0.675, 0.325) to[out=-45, in=-135] (0.825, 0.325);
        \draw[knot] (0.425, 0.425) to[out=45, in=135] (0.575, 0.425);
        \draw[knot] (0.425, 0.575) to[out=-45, in=-135] (0.575, 0.575);
        \draw[knot] (0.675, 0.675) to[out=45, in=-45] (0.675, 0.825);
        \draw[knot] (0.825, 0.675) to[out=135, in=-135] (0.825, 0.825);
} \circ \varphi_{H_1'}
    $}
    (100);
    \draw[knot, ->] (000) to
    node[pos=0.5, left, arrows=-] {$
    \tikz[baseline={([yshift=-.5ex]current bounding box.center)}, scale=0.6]
{
    \draw[dotted] (.5,.5) circle(0.707);
        \draw[knot, red] (0.7, 0.75) -- (0.8, 0.75);
        \draw[knot] (0.5, 1.207) to[out=-90, in=135] (0.675, 0.825);
        \draw[knot] (0.825, 0.675) to[out=-45, in=45] (0.825, 0.325);
        \draw[knot] (0.675, 0.175) to[out=-135, in=90] (0.5, -0.207);
        \draw[knot] (0,0) to (0.425, 0.425);
        \draw[knot] (0.575, 0.575) -- (0.675, 0.675);
        \draw[knot] (0.825, 0.825) -- (1,1);
        \draw[knot] (0,1) -- (0.425, 0.575);
        \draw[knot] (0.575, 0.425) -- (0.675, 0.325);
        \draw[knot] (0.825, 0.175) -- (1,0);
        \draw[knot] (0.675, 0.175) to[out=45, in=135] (0.825, 0.175);
        \draw[knot] (0.675, 0.325) to[out=-45, in=-135] (0.825, 0.325);
        \draw[knot] (0.425, 0.425) to[out=45, in=135] (0.575, 0.425);
        \draw[knot] (0.425, 0.575) to[out=-45, in=-135] (0.575, 0.575);
        \draw[knot] (0.675, 0.675) to[out=45, in=-45] (0.675, 0.825);
        \draw[knot] (0.825, 0.675) to[out=135, in=-135] (0.825, 0.825);
}
    $}
    (001);
    \draw[knot, ->] (100) to
    node[pos=0.5, left, arrows=-] {$
    \tikz[baseline={([yshift=-.5ex]current bounding box.center)}, scale=0.6]
{
    \draw[dotted] (.5,.5) circle(0.707);
        \draw[knot, red] (0.7, 0.75) -- (0.8, 0.75);
        \draw[knot] (0.5, 1.207) to[out=-90, in=135] (0.675, 0.825);
        \draw[knot] (0.825, 0.675) to[out=-45, in=45] (0.825, 0.325);
        \draw[knot] (0.675, 0.175) to[out=-135, in=90] (0.5, -0.207);
        \draw[knot] (0,0) to (0.425, 0.425);
        \draw[knot] (0.575, 0.575) -- (0.675, 0.675);
        \draw[knot] (0.825, 0.825) -- (1,1);
        \draw[knot] (0,1) -- (0.425, 0.575);
        \draw[knot] (0.575, 0.425) -- (0.675, 0.325);
        \draw[knot] (0.825, 0.175) -- (1,0);
        \draw[knot] (0.675, 0.175) to[out=45, in=-45] (0.675, 0.325);
        \draw[knot] (0.825, 0.175) to[out=135, in=-135] (0.825, 0.325);
        \draw[knot] (0.425, 0.425) to[out=45, in=135] (0.575, 0.425);
        \draw[knot] (0.425, 0.575) to[out=-45, in=-135] (0.575, 0.575);
        \draw[knot] (0.675, 0.675) to[out=45, in=-45] (0.675, 0.825);
        \draw[knot] (0.825, 0.675) to[out=135, in=-135] (0.825, 0.825);
}
    $}
    (101);
    \draw[knot, ->] (001) to
    node[pos=0.5, above, arrows=-] {$
    \tikz[baseline={([yshift=-.5ex]current bounding box.center)}, scale=0.6]
{
    \draw[dotted] (.5,.5) circle(0.707);
        \draw[knot, red] (0.75, 0.2) -- (0.75, 0.3);
        \draw[knot] (0.5, 1.207) to[out=-90, in=135] (0.675, 0.825);
        \draw[knot] (0.825, 0.675) to[out=-45, in=45] (0.825, 0.325);
        \draw[knot] (0.675, 0.175) to[out=-135, in=90] (0.5, -0.207);
        \draw[knot] (0,0) to (0.425, 0.425);
        \draw[knot] (0.575, 0.575) -- (0.675, 0.675);
        \draw[knot] (0.825, 0.825) -- (1,1);
        \draw[knot] (0,1) -- (0.425, 0.575);
        \draw[knot] (0.575, 0.425) -- (0.675, 0.325);
        \draw[knot] (0.825, 0.175) -- (1,0);
        \draw[knot] (0.675, 0.175) to[out=45, in=135] (0.825, 0.175);
        \draw[knot] (0.675, 0.325) to[out=-45, in=-135] (0.825, 0.325);
        \draw[knot] (0.425, 0.425) to[out=45, in=135] (0.575, 0.425);
        \draw[knot] (0.425, 0.575) to[out=-45, in=-135] (0.575, 0.575);
        \draw[knot] (0.675, 0.825) to[out=-45, in=-135] (0.825, 0.825);
        \draw[knot] (0.675, 0.675) to[out=45, in=135] (0.825, 0.675);
} \circ \varphi_{H_3'}
    $}
    (101);
    \draw[knot, ->] (001) to[out=-30, in=-120] (011);
    \draw[knot, ->] (100) to[out=0, in=120] (011);
    \node at (8.45, -4.85) {$-
    \tikz[baseline={([yshift=-.5ex]current bounding box.center)}, scale=0.6]
{
    \draw[dotted] (.5,.5) circle(0.707);
        \draw[knot, red] (0.9, 0.1) to[out=45, in=-45] (0.9, 0.9);
        \draw[knot] (0.5, 1.207) to[out=-90, in=135] (0.675, 0.825);
        \draw[knot] (0.825, 0.675) to[out=-45, in=45] (0.825, 0.325);
        \draw[knot] (0.675, 0.175) to[out=-135, in=90] (0.5, -0.207);
        \draw[knot] (0,0) to (0.425, 0.425);
        \draw[knot] (0.575, 0.575) -- (0.675, 0.675);
        \draw[knot] (0.825, 0.825) -- (1,1);
        \draw[knot] (0,1) -- (0.425, 0.575);
        \draw[knot] (0.575, 0.425) -- (0.675, 0.325);
        \draw[knot] (0.825, 0.175) -- (1,0);
        \draw[knot] (0.675, 0.175) to[out=45, in=135] (0.825, 0.175);
        \draw[knot] (0.675, 0.325) to[out=-45, in=-135] (0.825, 0.325);
        \draw[knot] (0.425, 0.425) to[out=45, in=135] (0.575, 0.425);
        \draw[knot] (0.425, 0.575) to[out=-45, in=-135] (0.575, 0.575);
        \draw[knot] (0.675, 0.825) to[out=-45, in=-135] (0.825, 0.825);
        \draw[knot] (0.675, 0.675) to[out=45, in=135] (0.825, 0.675);
} \circ \varphi_{H_2'}
    $};
    \node at (8.35,-0.5) {$
    \tikz[baseline={([yshift=-.5ex]current bounding box.center)}, scale=0.6]
{
    \draw[dotted] (.5,.5) circle(0.707);
        \draw[knot, red] (0.5, 0.45) -- (0.5, 0.55);
        \draw[knot] (0.5, 1.207) to[out=-90, in=135] (0.675, 0.825);
        \draw[knot] (0.825, 0.675) to[out=-45, in=45] (0.825, 0.325);
        \draw[knot] (0.675, 0.175) to[out=-135, in=90] (0.5, -0.207);
        \draw[knot] (0,0) to (0.425, 0.425);
        \draw[knot] (0.575, 0.575) -- (0.675, 0.675);
        \draw[knot] (0.825, 0.825) -- (1,1);
        \draw[knot] (0,1) -- (0.425, 0.575);
        \draw[knot] (0.575, 0.425) -- (0.675, 0.325);
        \draw[knot] (0.825, 0.175) -- (1,0);
        \draw[knot] (0.675, 0.175) to[out=45, in=-45] (0.675, 0.325);
        \draw[knot] (0.825, 0.175) to[out=135, in=-135] (0.825, 0.325);
        \draw[knot] (0.425, 0.425) to[out=45, in=135] (0.575, 0.425);
        \draw[knot] (0.425, 0.575) to[out=-45, in=-135] (0.575, 0.575);
        \draw[knot] (0.675, 0.675) to[out=45, in=-45] (0.675, 0.825);
        \draw[knot] (0.825, 0.675) to[out=135, in=-135] (0.825, 0.825);
}
    $}
}
\]
Again, we know these complexes are \textit{not} homotopy equivalent by, for example, Remark \ref{rem:R3noticing}. Instead, we will show that the latter is taken to the former by a global grading shift of
\[
\varphi_{
    \tikz[baseline={([yshift=-.5ex]current bounding box.center)}, scale=0.75]
{
    \draw[dotted] (.5,.5) circle(0.707);
        \draw[knot, red] (0.2, 0.75) -- (0.3, 0.75);
        \draw[knot] (0.5, 1.207) to[out=-90, in=45] (0.325, 0.825);
        \draw[knot] (0.175, 0.675) to[out=225, in=135] (0.175, 0.325);
        \draw[knot] (0.325, 0.175) to[out=-45, in=90] (0.5, -0.207);
        \draw[knot] (0,0) to (0.175, 0.175);
        \draw[knot] (0.325, 0.325) to (0.425, 0.425);
        \draw[knot] (0.575, 0.575) to (1,1);
        \draw[knot] (0,1) to (0.175, 0.825);
        \draw[knot] (0.325, 0.675) to (0.425, 0.575);
        \draw[knot] (0.575, 0.425) to (1,0);
        \draw[knot] (0.175, 0.175) to[out=45, in=-45] (0.175, 0.325);
        \draw[knot] (0.325, 0.175) to[out=135, in=-135] (0.325, 0.325);
        \draw[knot] (0.425, 0.425) to[out=45, in=-45] (0.425, 0.575); 
        \draw[knot] (0.575, 0.425) to[out=135, in=-135] (0.575, 0.575);
        \draw[knot] (0.175, 0.675) to[out=45, in=-45] (0.175, 0.825);
        \draw[knot] (0.325, 0.675) to[out=135, in=-135] (0.325, 0.825);
}}
\circ
\varphi_{
\tikz[baseline={([yshift=-.5ex]current bounding box.center)}, scale=0.75]
{
    \draw[dotted] (.5,.5) circle(0.707);
        \draw[knot, red] (0.7, 0.25) -- (0.8, 0.25);
        \draw[knot] (0.5, 1.207) to[out=-90, in=135] (0.675, 0.825);
        \draw[knot] (0.825, 0.675) to[out=-45, in=45] (0.825, 0.325);
        \draw[knot] (0.675, 0.175) to[out=-135, in=90] (0.5, -0.207);
        \draw[knot] (0,0) to (0.425, 0.425);
        \draw[knot] (0.575, 0.575) -- (0.675, 0.675);
        \draw[knot] (0.825, 0.825) -- (1,1);
        \draw[knot] (0,1) -- (0.425, 0.575);
        \draw[knot] (0.575, 0.425) -- (0.675, 0.325);
        \draw[knot] (0.825, 0.175) -- (1,0);
        \draw[knot] (0.675, 0.175) to[out=45, in=-45] (0.675, 0.325);
        \draw[knot] (0.825, 0.175) to[out=135, in=-135] (0.825, 0.325);
        \draw[knot] (0.425, 0.425) to[out=45, in=-45] (0.425, 0.575); 
        \draw[knot] (0.575, 0.425) to[out=135, in=-135] (0.575, 0.575);
        \draw[knot] (0.675, 0.675) to[out=45, in=-45] (0.675, 0.825);
        \draw[knot] (0.825, 0.675) to[out=135, in=-135] (0.825, 0.825);
}
}^{-1}. 
\]
First, recall that $\varphi_{
\tikz[baseline={([yshift=-.5ex]current bounding box.center)}, scale=0.75]
{
    \draw[dotted] (.5,.5) circle(0.707);
        \draw[knot, red] (0.7, 0.25) -- (0.8, 0.25);
        \draw[knot] (0.5, 1.207) to[out=-90, in=135] (0.675, 0.825);
        \draw[knot] (0.825, 0.675) to[out=-45, in=45] (0.825, 0.325);
        \draw[knot] (0.675, 0.175) to[out=-135, in=90] (0.5, -0.207);
        \draw[knot] (0,0) to (0.425, 0.425);
        \draw[knot] (0.575, 0.575) -- (0.675, 0.675);
        \draw[knot] (0.825, 0.825) -- (1,1);
        \draw[knot] (0,1) -- (0.425, 0.575);
        \draw[knot] (0.575, 0.425) -- (0.675, 0.325);
        \draw[knot] (0.825, 0.175) -- (1,0);
        \draw[knot] (0.675, 0.175) to[out=45, in=-45] (0.675, 0.325);
        \draw[knot] (0.825, 0.175) to[out=135, in=-135] (0.825, 0.325);
        \draw[knot] (0.425, 0.425) to[out=45, in=-45] (0.425, 0.575); 
        \draw[knot] (0.575, 0.425) to[out=135, in=-135] (0.575, 0.575);
        \draw[knot] (0.675, 0.675) to[out=45, in=-45] (0.675, 0.825);
        \draw[knot] (0.825, 0.675) to[out=135, in=-135] (0.825, 0.825);
}
}^{-1}$
may be written, up to equivalence of grading shift functors, as $\varphi_{\left(
\tikz[baseline={([yshift=-.5ex]current bounding box.center)}, scale=0.75]
{
    \draw[dotted] (.5,.5) circle(0.707);
        \draw[knot, red] (0.75, 0.2) -- (0.75, 0.3);
        \draw[knot] (0.5, 1.207) to[out=-90, in=135] (0.675, 0.825);
        \draw[knot] (0.825, 0.675) to[out=-45, in=45] (0.825, 0.325);
        \draw[knot] (0.675, 0.175) to[out=-135, in=90] (0.5, -0.207);
        \draw[knot] (0,0) to (0.425, 0.425);
        \draw[knot] (0.575, 0.575) -- (0.675, 0.675);
        \draw[knot] (0.825, 0.825) -- (1,1);
        \draw[knot] (0,1) -- (0.425, 0.575);
        \draw[knot] (0.575, 0.425) -- (0.675, 0.325);
        \draw[knot] (0.825, 0.175) -- (1,0);
        \draw[knot] (0.675, 0.175) to[out=45, in=135] (0.825, 0.175);
        \draw[knot] (0.675, 0.325) to[out=-45, in=-135] (0.825, 0.325);
        \draw[knot] (0.425, 0.425) to[out=45, in=-45] (0.425, 0.575); 
        \draw[knot] (0.575, 0.425) to[out=135, in=-135] (0.575, 0.575);
        \draw[knot] (0.675, 0.675) to[out=45, in=-45] (0.675, 0.825);
        \draw[knot] (0.825, 0.675) to[out=135, in=-135] (0.825, 0.825);
},\, (1,1)\right)
}$. On the other hand, $\varphi_{
    \tikz[baseline={([yshift=-.5ex]current bounding box.center)}, scale=0.75]
{
    \draw[dotted] (.5,.5) circle(0.707);
        \draw[knot, red] (0.2, 0.75) -- (0.3, 0.75);
        \draw[knot] (0.5, 1.207) to[out=-90, in=45] (0.325, 0.825);
        \draw[knot] (0.175, 0.675) to[out=225, in=135] (0.175, 0.325);
        \draw[knot] (0.325, 0.175) to[out=-45, in=90] (0.5, -0.207);
        \draw[knot] (0,0) to (0.175, 0.175);
        \draw[knot] (0.325, 0.325) to (0.425, 0.425);
        \draw[knot] (0.575, 0.575) to (1,1);
        \draw[knot] (0,1) to (0.175, 0.825);
        \draw[knot] (0.325, 0.675) to (0.425, 0.575);
        \draw[knot] (0.575, 0.425) to (1,0);
        \draw[knot] (0.175, 0.175) to[out=45, in=-45] (0.175, 0.325);
        \draw[knot] (0.325, 0.175) to[out=135, in=-135] (0.325, 0.325);
        \draw[knot] (0.425, 0.425) to[out=45, in=-45] (0.425, 0.575); 
        \draw[knot] (0.575, 0.425) to[out=135, in=-135] (0.575, 0.575);
        \draw[knot] (0.175, 0.675) to[out=45, in=-45] (0.175, 0.825);
        \draw[knot] (0.325, 0.675) to[out=135, in=-135] (0.325, 0.825);
}}$ and $\varphi_{\left(
    \tikz[baseline={([yshift=-.5ex]current bounding box.center)}, scale=0.75]
{
    \draw[dotted] (.5,.5) circle(0.707);
        \draw[knot, red] (0.25, 0.7) -- (0.25, 0.8);
        \draw[knot] (0.5, 1.207) to[out=-90, in=45] (0.325, 0.825);
        \draw[knot] (0.175, 0.675) to[out=225, in=135] (0.175, 0.325);
        \draw[knot] (0.325, 0.175) to[out=-45, in=90] (0.5, -0.207);
        \draw[knot] (0,0) to (0.175, 0.175);
        \draw[knot] (0.325, 0.325) to (0.425, 0.425);
        \draw[knot] (0.575, 0.575) to (1,1);
        \draw[knot] (0,1) to (0.175, 0.825);
        \draw[knot] (0.325, 0.675) to (0.425, 0.575);
        \draw[knot] (0.575, 0.425) to (1,0);
        \draw[knot] (0.175, 0.175) to[out=45, in=-45] (0.175, 0.325);
        \draw[knot] (0.325, 0.175) to[out=135, in=-135] (0.325, 0.325);
        \draw[knot] (0.425, 0.425) to[out=45, in=-45] (0.425, 0.575); 
        \draw[knot] (0.575, 0.425) to[out=135, in=-135] (0.575, 0.575);
        \draw[knot] (0.175, 0.675) to[out=45, in=135] (0.325, 0.675);
        \draw[knot] (0.175, 0.825) to[out=-45, in=-135] (0.325, 0.825);
}, \, (-1,-1)\right)}^{-1}$ are isomorphic as grading shift functors.  
\begin{enumerate}[label=(\roman*)]
    \item Northwest vertex. As a warm-up, notice that $\varphi_{
\tikz[baseline={([yshift=-.5ex]current bounding box.center)}, scale=0.75]
{
    \draw[dotted] (.5,.5) circle(0.707);
        \draw[knot, red] (0.5, 0.45) -- (0.5, 0.55);
        \draw[knot, red] (0.7, 0.75) -- (0.8, 0.75);
        \draw[knot] (0.5, 1.207) to[out=-90, in=135] (0.675, 0.825);
        \draw[knot] (0.825, 0.675) to[out=-45, in=45] (0.825, 0.325);
        \draw[knot] (0.675, 0.175) to[out=-135, in=90] (0.5, -0.207);
        \draw[knot] (0,0) to (0.425, 0.425);
        \draw[knot] (0.575, 0.575) -- (0.675, 0.675);
        \draw[knot] (0.825, 0.825) -- (1,1);
        \draw[knot] (0,1) -- (0.425, 0.575);
        \draw[knot] (0.575, 0.425) -- (0.675, 0.325);
        \draw[knot] (0.825, 0.175) -- (1,0);
        \draw[knot] (0.675, 0.175) to[out=45, in=135] (0.825, 0.175);
        \draw[knot] (0.675, 0.325) to[out=-45, in=-135] (0.825, 0.325);
        \draw[knot] (0.425, 0.425) to[out=45, in=135] (0.575, 0.425);
        \draw[knot] (0.425, 0.575) to[out=-45, in=-135] (0.575, 0.575);
        \draw[knot] (0.675, 0.675) to[out=45, in=-45] (0.675, 0.825);
        \draw[knot] (0.825, 0.675) to[out=135, in=-135] (0.825, 0.825);
}}$ has two isomorphic representatives important to understanding the intermediate complex. They are hardly different, but making a choice here is one way to describe two representatives of the $\mathscr{G}$-grading shift obtained after the first global shift:
\[
\varphi_{\left(
\tikz[baseline={([yshift=-.5ex]current bounding box.center)}, scale=0.75]
{
    \draw[dotted] (.5,.5) circle(0.707);
        \draw[knot, red] (0.75, 0.2) -- (0.75, 0.3);
        \draw[knot] (0.5, 1.207) to[out=-90, in=135] (0.675, 0.825);
        \draw[knot] (0.825, 0.675) to[out=-45, in=45] (0.825, 0.325);
        \draw[knot] (0.675, 0.175) to[out=-135, in=90] (0.5, -0.207);
        \draw[knot] (0,0) to (0.425, 0.425);
        \draw[knot] (0.575, 0.575) -- (0.675, 0.675);
        \draw[knot] (0.825, 0.825) -- (1,1);
        \draw[knot] (0,1) -- (0.425, 0.575);
        \draw[knot] (0.575, 0.425) -- (0.675, 0.325);
        \draw[knot] (0.825, 0.175) -- (1,0);
        \draw[knot] (0.675, 0.175) to[out=45, in=135] (0.825, 0.175);
        \draw[knot] (0.675, 0.325) to[out=-45, in=-135] (0.825, 0.325);
        \draw[knot] (0.425, 0.425) to[out=45, in=-45] (0.425, 0.575); 
        \draw[knot] (0.575, 0.425) to[out=135, in=-135] (0.575, 0.575);
        \draw[knot] (0.675, 0.675) to[out=45, in=-45] (0.675, 0.825);
        \draw[knot] (0.825, 0.675) to[out=135, in=-135] (0.825, 0.825);
},\, (1,1)\right)
}
\circ
\begin{Bmatrix}
\varphi_{\left(
\tikz[baseline={([yshift=-.5ex]current bounding box.center)}, scale=0.75]
{
    \draw[dotted] (.5,.5) circle(0.707);
        \draw[knot, red] (0.7, 0.75) -- (0.8, 0.75);
        \draw[knot] (0.5, 1.207) to[out=-90, in=135] (0.675, 0.825);
        \draw[knot] (0.825, 0.675) to[out=-45, in=45] (0.825, 0.325);
        \draw[knot] (0.675, 0.175) to[out=-135, in=90] (0.5, -0.207);
        \draw[knot] (0,0) to (0.425, 0.425);
        \draw[knot] (0.575, 0.575) -- (0.675, 0.675);
        \draw[knot] (0.825, 0.825) -- (1,1);
        \draw[knot] (0,1) -- (0.425, 0.575);
        \draw[knot] (0.575, 0.425) -- (0.675, 0.325);
        \draw[knot] (0.825, 0.175) -- (1,0);
        \draw[knot] (0.675, 0.175) to[out=45, in=135] (0.825, 0.175);
        \draw[knot] (0.675, 0.325) to[out=-45, in=-135] (0.825, 0.325);
        \draw[knot] (0.425, 0.425) to[out=45, in=135] (0.575, 0.425);
        \draw[knot] (0.425, 0.575) to[out=-45, in=-135] (0.575, 0.575);
        \draw[knot] (0.675, 0.675) to[out=45, in=-45] (0.675, 0.825);
        \draw[knot] (0.825, 0.675) to[out=135, in=-135] (0.825, 0.825);
}, \, (0,-1)\right)}
\\
\varphi_{\left(
\tikz[baseline={([yshift=-.5ex]current bounding box.center)}, scale=0.75]
{
    \draw[dotted] (.5,.5) circle(0.707);
        \draw[knot, red] (0.5, 0.45) -- (0.5, 0.55);
        \draw[knot] (0.5, 1.207) to[out=-90, in=135] (0.675, 0.825);
        \draw[knot] (0.825, 0.675) to[out=-45, in=45] (0.825, 0.325);
        \draw[knot] (0.675, 0.175) to[out=-135, in=90] (0.5, -0.207);
        \draw[knot] (0,0) to (0.425, 0.425);
        \draw[knot] (0.575, 0.575) -- (0.675, 0.675);
        \draw[knot] (0.825, 0.825) -- (1,1);
        \draw[knot] (0,1) -- (0.425, 0.575);
        \draw[knot] (0.575, 0.425) -- (0.675, 0.325);
        \draw[knot] (0.825, 0.175) -- (1,0);
        \draw[knot] (0.675, 0.175) to[out=45, in=135] (0.825, 0.175);
        \draw[knot] (0.675, 0.325) to[out=-45, in=-135] (0.825, 0.325);
        \draw[knot] (0.425, 0.425) to[out=45, in=135] (0.575, 0.425);
        \draw[knot] (0.425, 0.575) to[out=-45, in=-135] (0.575, 0.575);
        \draw[knot] (0.675, 0.675) to[out=45, in=-45] (0.675, 0.825);
        \draw[knot] (0.825, 0.675) to[out=135, in=-135] (0.825, 0.825);
},\, (0,-1)\right)}
\end{Bmatrix}
\cong
\begin{Bmatrix}
\varphi_{\left(
\tikz[baseline={([yshift=-.5ex]current bounding box.center)}, scale=0.75]
{
    \draw[dotted] (.5,.5) circle(0.707);
        \draw[knot, red] (0.75, 0.2) -- (0.75, 0.3);
        \draw[knot, red] (0.7, 0.75) -- (0.8, 0.75);
        \draw[knot] (0.5, 1.207) to[out=-90, in=135] (0.675, 0.825);
        \draw[knot] (0.825, 0.675) to[out=-45, in=45] (0.825, 0.325);
        \draw[knot] (0.675, 0.175) to[out=-135, in=90] (0.5, -0.207);
        \draw[knot] (0,0) to (0.425, 0.425);
        \draw[knot] (0.575, 0.575) -- (0.675, 0.675);
        \draw[knot] (0.825, 0.825) -- (1,1);
        \draw[knot] (0,1) -- (0.425, 0.575);
        \draw[knot] (0.575, 0.425) -- (0.675, 0.325);
        \draw[knot] (0.825, 0.175) -- (1,0);
        \draw[knot] (0.675, 0.175) to[out=45, in=135] (0.825, 0.175);
        \draw[knot] (0.675, 0.325) to[out=-45, in=-135] (0.825, 0.325);
        \draw[knot] (0.425, 0.425) to[out=45, in=135] (0.575, 0.425);
        \draw[knot] (0.425, 0.575) to[out=-45, in=-135] (0.575, 0.575);
        \draw[knot] (0.675, 0.675) to[out=45, in=-45] (0.675, 0.825);
        \draw[knot] (0.825, 0.675) to[out=135, in=-135] (0.825, 0.825);
}, \, (1,0)\right)}
\\
\varphi_{\left(
\tikz[baseline={([yshift=-.5ex]current bounding box.center)}, scale=0.75]
{
    \draw[dotted] (.5,.5) circle(0.707);
        \draw[knot, red] (0.75, 0.2) -- (0.75, 0.3);
        \draw[knot, red] (0.5, 0.45) -- (0.5, 0.55);
        \draw[knot] (0.5, 1.207) to[out=-90, in=135] (0.675, 0.825);
        \draw[knot] (0.825, 0.675) to[out=-45, in=45] (0.825, 0.325);
        \draw[knot] (0.675, 0.175) to[out=-135, in=90] (0.5, -0.207);
        \draw[knot] (0,0) to (0.425, 0.425);
        \draw[knot] (0.575, 0.575) -- (0.675, 0.675);
        \draw[knot] (0.825, 0.825) -- (1,1);
        \draw[knot] (0,1) -- (0.425, 0.575);
        \draw[knot] (0.575, 0.425) -- (0.675, 0.325);
        \draw[knot] (0.825, 0.175) -- (1,0);
        \draw[knot] (0.675, 0.175) to[out=45, in=135] (0.825, 0.175);
        \draw[knot] (0.675, 0.325) to[out=-45, in=-135] (0.825, 0.325);
        \draw[knot] (0.425, 0.425) to[out=45, in=135] (0.575, 0.425);
        \draw[knot] (0.425, 0.575) to[out=-45, in=-135] (0.575, 0.575);
        \draw[knot] (0.675, 0.675) to[out=45, in=-45] (0.675, 0.825);
        \draw[knot] (0.825, 0.675) to[out=135, in=-135] (0.825, 0.825);
},\, (1,0)\right)}
\end{Bmatrix}.
\]
Of course, yet another representative of this grading shift, encapsulating both of these representatives, is $\varphi_{\left(
\tikz[baseline={([yshift=-.5ex]current bounding box.center)}, scale=0.75]
{
    \draw[dotted] (.5,.5) circle(0.707);
        \draw[knot, red] (0.75, 0.2) -- (0.75, 0.3);
        \draw[knot, red] (0.5, 0.45) -- (0.5, 0.55);
        \draw[knot, red] (0.7, 0.75) -- (0.8, 0.75);
        \draw[knot] (0.5, 1.207) to[out=-90, in=135] (0.675, 0.825);
        \draw[knot] (0.825, 0.675) to[out=-45, in=45] (0.825, 0.325);
        \draw[knot] (0.675, 0.175) to[out=-135, in=90] (0.5, -0.207);
        \draw[knot] (0,0) to (0.425, 0.425);
        \draw[knot] (0.575, 0.575) -- (0.675, 0.675);
        \draw[knot] (0.825, 0.825) -- (1,1);
        \draw[knot] (0,1) -- (0.425, 0.575);
        \draw[knot] (0.575, 0.425) -- (0.675, 0.325);
        \draw[knot] (0.825, 0.175) -- (1,0);
        \draw[knot] (0.675, 0.175) to[out=45, in=135] (0.825, 0.175);
        \draw[knot] (0.675, 0.325) to[out=-45, in=-135] (0.825, 0.325);
        \draw[knot] (0.425, 0.425) to[out=45, in=135] (0.575, 0.425);
        \draw[knot] (0.425, 0.575) to[out=-45, in=-135] (0.575, 0.575);
        \draw[knot] (0.675, 0.675) to[out=45, in=-45] (0.675, 0.825);
        \draw[knot] (0.825, 0.675) to[out=135, in=-135] (0.825, 0.825);
},\, (1,1)\right)}$. Anyway, applying the final global shift to the second representative above, we obtain the grading shift $\varphi_{\left(
\tikz[baseline={([yshift=-.5ex]current bounding box.center)}, scale=0.75]
{
    \draw[dotted] (.5,.5) circle(0.707);
        \draw[knot, red] (0.75, 0.2) -- (0.75, 0.3);
        \draw[knot] (0.5, 1.207) to[out=-90, in=135] (0.675, 0.825);
        \draw[knot] (0.825, 0.675) to[out=-45, in=45] (0.825, 0.325);
        \draw[knot] (0.675, 0.175) to[out=-135, in=90] (0.5, -0.207);
        \draw[knot] (0,0) to (0.425, 0.425);
        \draw[knot] (0.575, 0.575) -- (0.675, 0.675);
        \draw[knot] (0.825, 0.825) -- (1,1);
        \draw[knot] (0,1) -- (0.425, 0.575);
        \draw[knot] (0.575, 0.425) -- (0.675, 0.325);
        \draw[knot] (0.825, 0.175) -- (1,0);
        \draw[knot] (0.675, 0.175) to[out=45, in=135] (0.825, 0.175);
        \draw[knot] (0.675, 0.325) to[out=-45, in=-135] (0.825, 0.325);
        \draw[knot] (0.425, 0.425) to[out=45, in=135] (0.575, 0.425);
        \draw[knot] (0.425, 0.575) to[out=-45, in=-135] (0.575, 0.575);
        \draw[knot] (0.675, 0.675) to[out=45, in=-45] (0.675, 0.825);
        \draw[knot] (0.825, 0.675) to[out=135, in=-135] (0.825, 0.825);
},\, (0,-1)\right)} = \varphi_{\left(
    \tikz[baseline={([yshift=-.5ex]current bounding box.center)}, scale=.75]
{
    \draw[dotted] (.5,.5) circle(0.707);
        \draw[knot, red] (0.2, 0.25) -- (0.3, 0.25);
        \draw[knot] (0.5, 1.207) to[out=-90, in=45] (0.325, 0.825);
        \draw[knot] (0.175, 0.675) to[out=225, in=135] (0.175, 0.325);
        \draw[knot] (0.325, 0.175) to[out=-45, in=90] (0.5, -0.207);
        \draw[knot] (0,0) to (0.175, 0.175);
        \draw[knot] (0.325, 0.325) to (0.425, 0.425);
        \draw[knot] (0.575, 0.575) to (1,1);
        \draw[knot] (0,1) to (0.175, 0.825);
        \draw[knot] (0.325, 0.675) to (0.425, 0.575);
        \draw[knot] (0.575, 0.425) to (1,0);
        \draw[knot] (0.175, 0.175) to[out=45, in=-45] (0.175, 0.325);
        \draw[knot] (0.325, 0.175) to[out=135, in=-135] (0.325, 0.325);
        \draw[knot] (0.425, 0.425) to[out=45, in=135] (0.575, 0.425);
        \draw[knot] (0.425, 0.575) to[out=-45, in=-135] (0.575, 0.575);
        \draw[knot] (0.175, 0.825) to[out=-45, in=-135] (0.325, 0.825);
        \draw[knot] (0.175, 0.675) to[out=45, in=135] (0.325, 0.675);
}, \, (0, -1)\right)}$ which, similarly, is a representative of the grading shift $\varphi_{
    \tikz[baseline={([yshift=-.5ex]current bounding box.center)}, scale=.75]
{
    \draw[dotted] (.5,.5) circle(0.707);
        \draw[knot, red] (0.2, 0.25) -- (0.3, 0.25);
        \draw[knot, red] (0.5, 0.45) -- (0.5, 0.55);
        \draw[knot] (0.5, 1.207) to[out=-90, in=45] (0.325, 0.825);
        \draw[knot] (0.175, 0.675) to[out=225, in=135] (0.175, 0.325);
        \draw[knot] (0.325, 0.175) to[out=-45, in=90] (0.5, -0.207);
        \draw[knot] (0,0) to (0.175, 0.175);
        \draw[knot] (0.325, 0.325) to (0.425, 0.425);
        \draw[knot] (0.575, 0.575) to (1,1);
        \draw[knot] (0,1) to (0.175, 0.825);
        \draw[knot] (0.325, 0.675) to (0.425, 0.575);
        \draw[knot] (0.575, 0.425) to (1,0);
        \draw[knot] (0.175, 0.175) to[out=45, in=-45] (0.175, 0.325);
        \draw[knot] (0.325, 0.175) to[out=135, in=-135] (0.325, 0.325);
        \draw[knot] (0.425, 0.425) to[out=45, in=135] (0.575, 0.425);
        \draw[knot] (0.425, 0.575) to[out=-45, in=-135] (0.575, 0.575);
        \draw[knot] (0.175, 0.825) to[out=-45, in=-135] (0.325, 0.825);
        \draw[knot] (0.175, 0.675) to[out=45, in=135] (0.325, 0.675);
}}$.
    \item Southwest vertex. This is the trickiest since it is the vertex with one of its arrows altered by Gaussian elimination. On one hand, obviously if we apply $\varphi_{\left(
\tikz[baseline={([yshift=-.5ex]current bounding box.center)}, scale=0.75]
{
    \draw[dotted] (.5,.5) circle(0.707);
        \draw[knot, red] (0.75, 0.2) -- (0.75, 0.3);
        \draw[knot] (0.5, 1.207) to[out=-90, in=135] (0.675, 0.825);
        \draw[knot] (0.825, 0.675) to[out=-45, in=45] (0.825, 0.325);
        \draw[knot] (0.675, 0.175) to[out=-135, in=90] (0.5, -0.207);
        \draw[knot] (0,0) to (0.425, 0.425);
        \draw[knot] (0.575, 0.575) -- (0.675, 0.675);
        \draw[knot] (0.825, 0.825) -- (1,1);
        \draw[knot] (0,1) -- (0.425, 0.575);
        \draw[knot] (0.575, 0.425) -- (0.675, 0.325);
        \draw[knot] (0.825, 0.175) -- (1,0);
        \draw[knot] (0.675, 0.175) to[out=45, in=135] (0.825, 0.175);
        \draw[knot] (0.675, 0.325) to[out=-45, in=-135] (0.825, 0.325);
        \draw[knot] (0.425, 0.425) to[out=45, in=-45] (0.425, 0.575); 
        \draw[knot] (0.575, 0.425) to[out=135, in=-135] (0.575, 0.575);
        \draw[knot] (0.675, 0.675) to[out=45, in=-45] (0.675, 0.825);
        \draw[knot] (0.825, 0.675) to[out=135, in=-135] (0.825, 0.825);
},\, (1,1)\right)
}$ to $\tikz[baseline={([yshift=-.5ex]current bounding box.center)}, scale=0.85]
{
    \draw[dotted] (.5,.5) circle(0.707);
        \draw[knot] (0.5, 1.207) to[out=-90, in=135] (0.675, 0.825);
        \draw[knot] (0.825, 0.675) to[out=-45, in=45] (0.825, 0.325);
        \draw[knot] (0.675, 0.175) to[out=-135, in=90] (0.5, -0.207);
        \draw[knot] (0,0) to (0.425, 0.425);
        \draw[knot] (0.575, 0.575) -- (0.675, 0.675);
        \draw[knot] (0.825, 0.825) -- (1,1);
        \draw[knot] (0,1) -- (0.425, 0.575);
        \draw[knot] (0.575, 0.425) -- (0.675, 0.325);
        \draw[knot] (0.825, 0.175) -- (1,0);
        \draw[knot] (0.675, 0.175) to[out=45, in=135] (0.825, 0.175);
        \draw[knot] (0.675, 0.325) to[out=-45, in=-135] (0.825, 0.325);
        \draw[knot] (0.425, 0.425) to[out=45, in=135] (0.575, 0.425);
        \draw[knot] (0.425, 0.575) to[out=-45, in=-135] (0.575, 0.575);
        \draw[knot] (0.675, 0.825) to[out=-45, in=-135] (0.825, 0.825);
        \draw[knot] (0.675, 0.675) to[out=45, in=135] (0.825, 0.675);
} \, \{0,-1\}$ we are left with $\varphi_{\left(\tikz[baseline={([yshift=-.5ex]current bounding box.center)}, scale=0.75]
{
    \draw[dotted] (.5,.5) circle(0.707);
        \draw[knot, red] (0.75, 0.2) -- (0.75, 0.3);
        \draw[knot] (0.5, 1.207) to[out=-90, in=135] (0.675, 0.825);
        \draw[knot] (0.825, 0.675) to[out=-45, in=45] (0.825, 0.325);
        \draw[knot] (0.675, 0.175) to[out=-135, in=90] (0.5, -0.207);
        \draw[knot] (0,0) to (0.425, 0.425);
        \draw[knot] (0.575, 0.575) -- (0.675, 0.675);
        \draw[knot] (0.825, 0.825) -- (1,1);
        \draw[knot] (0,1) -- (0.425, 0.575);
        \draw[knot] (0.575, 0.425) -- (0.675, 0.325);
        \draw[knot] (0.825, 0.175) -- (1,0);
        \draw[knot] (0.675, 0.175) to[out=45, in=135] (0.825, 0.175);
        \draw[knot] (0.675, 0.325) to[out=-45, in=-135] (0.825, 0.325);
        \draw[knot] (0.425, 0.425) to[out=45, in=135] (0.575, 0.425);
        \draw[knot] (0.425, 0.575) to[out=-45, in=-135] (0.575, 0.575);
        \draw[knot] (0.675, 0.825) to[out=-45, in=-135] (0.825, 0.825);
        \draw[knot] (0.675, 0.675) to[out=45, in=135] (0.825, 0.675);
},\, (1,0)\right)}$. At first, this may not seem to square with the other arrow out of the vertex. To see that this questionable arrow is still a graded map, one may draw the original cube and trace it through the Gaussian elimination; we leave this as an exercise. Moving on, rewrite the grading shift as $\varphi_{\left(
    \tikz[baseline={([yshift=-.5ex]current bounding box.center)}, scale=.75]
{
    \draw[dotted] (.5,.5) circle(0.707);
        \draw[knot, red] (0.2, 0.25) -- (0.3, 0.25);
        \draw[knot] (0.5, 1.207) to[out=-90, in=45] (0.325, 0.825);
        \draw[knot] (0.175, 0.675) to[out=225, in=135] (0.175, 0.325);
        \draw[knot] (0.325, 0.175) to[out=-45, in=90] (0.5, -0.207);
        \draw[knot] (0,0) to (0.175, 0.175);
        \draw[knot] (0.325, 0.325) to (0.425, 0.425);
        \draw[knot] (0.575, 0.575) to (1,1);
        \draw[knot] (0,1) to (0.175, 0.825);
        \draw[knot] (0.325, 0.675) to (0.425, 0.575);
        \draw[knot] (0.575, 0.425) to (1,0);
        \draw[knot] (0.175, 0.175) to[out=45, in=-45] (0.175, 0.325);
        \draw[knot] (0.325, 0.175) to[out=135, in=-135] (0.325, 0.325);
        \draw[knot] (0.425, 0.425) to[out=45, in=135] (0.575, 0.425);
        \draw[knot] (0.425, 0.575) to[out=-45, in=-135] (0.575, 0.575);
        \draw[knot] (0.175, 0.675) to[out=45, in=-45] (0.175, 0.825);
        \draw[knot] (0.325, 0.675) to[out=135, in=-135] (0.325, 0.825);
},\, (1,0)\right)}$ and apply $\varphi_{
    \tikz[baseline={([yshift=-.5ex]current bounding box.center)}, scale=0.75]
{
    \draw[dotted] (.5,.5) circle(0.707);
        \draw[knot, red] (0.2, 0.75) -- (0.3, 0.75);
        \draw[knot] (0.5, 1.207) to[out=-90, in=45] (0.325, 0.825);
        \draw[knot] (0.175, 0.675) to[out=225, in=135] (0.175, 0.325);
        \draw[knot] (0.325, 0.175) to[out=-45, in=90] (0.5, -0.207);
        \draw[knot] (0,0) to (0.175, 0.175);
        \draw[knot] (0.325, 0.325) to (0.425, 0.425);
        \draw[knot] (0.575, 0.575) to (1,1);
        \draw[knot] (0,1) to (0.175, 0.825);
        \draw[knot] (0.325, 0.675) to (0.425, 0.575);
        \draw[knot] (0.575, 0.425) to (1,0);
        \draw[knot] (0.175, 0.175) to[out=45, in=-45] (0.175, 0.325);
        \draw[knot] (0.325, 0.175) to[out=135, in=-135] (0.325, 0.325);
        \draw[knot] (0.425, 0.425) to[out=45, in=-45] (0.425, 0.575); 
        \draw[knot] (0.575, 0.425) to[out=135, in=-135] (0.575, 0.575);
        \draw[knot] (0.175, 0.675) to[out=45, in=-45] (0.175, 0.825);
        \draw[knot] (0.325, 0.675) to[out=135, in=-135] (0.325, 0.825);
}}$ to obtain $\varphi_{\left(
    \tikz[baseline={([yshift=-.5ex]current bounding box.center)}, scale=.75]
{
    \draw[dotted] (.5,.5) circle(0.707);
        \draw[knot, red] (0.2, 0.25) -- (0.3, 0.25);
        \draw[knot, red] (0.2, 0.75) -- (0.3, 0.75);
        \draw[knot] (0.5, 1.207) to[out=-90, in=45] (0.325, 0.825);
        \draw[knot] (0.175, 0.675) to[out=225, in=135] (0.175, 0.325);
        \draw[knot] (0.325, 0.175) to[out=-45, in=90] (0.5, -0.207);
        \draw[knot] (0,0) to (0.175, 0.175);
        \draw[knot] (0.325, 0.325) to (0.425, 0.425);
        \draw[knot] (0.575, 0.575) to (1,1);
        \draw[knot] (0,1) to (0.175, 0.825);
        \draw[knot] (0.325, 0.675) to (0.425, 0.575);
        \draw[knot] (0.575, 0.425) to (1,0);
        \draw[knot] (0.175, 0.175) to[out=45, in=-45] (0.175, 0.325);
        \draw[knot] (0.325, 0.175) to[out=135, in=-135] (0.325, 0.325);
        \draw[knot] (0.425, 0.425) to[out=45, in=135] (0.575, 0.425);
        \draw[knot] (0.425, 0.575) to[out=-45, in=-135] (0.575, 0.575);
        \draw[knot] (0.175, 0.675) to[out=45, in=-45] (0.175, 0.825);
        \draw[knot] (0.325, 0.675) to[out=135, in=-135] (0.325, 0.825);
},\, (1,0)\right)}$. This is a representative of the grading shift $\varphi_{\left(
    \tikz[baseline={([yshift=-.5ex]current bounding box.center)}, scale=.75]
{
    \draw[dotted] (.5,.5) circle(0.707);
        \draw[knot, red] (0.2, 0.25) -- (0.3, 0.25);
        \draw[knot, red] (0.5, 0.45) -- (0.5, 0.55);
        \draw[knot, red] (0.2, 0.75) -- (0.3, 0.75);
        \draw[knot] (0.5, 1.207) to[out=-90, in=45] (0.325, 0.825);
        \draw[knot] (0.175, 0.675) to[out=225, in=135] (0.175, 0.325);
        \draw[knot] (0.325, 0.175) to[out=-45, in=90] (0.5, -0.207);
        \draw[knot] (0,0) to (0.175, 0.175);
        \draw[knot] (0.325, 0.325) to (0.425, 0.425);
        \draw[knot] (0.575, 0.575) to (1,1);
        \draw[knot] (0,1) to (0.175, 0.825);
        \draw[knot] (0.325, 0.675) to (0.425, 0.575);
        \draw[knot] (0.575, 0.425) to (1,0);
        \draw[knot] (0.175, 0.175) to[out=45, in=-45] (0.175, 0.325);
        \draw[knot] (0.325, 0.175) to[out=135, in=-135] (0.325, 0.325);
        \draw[knot] (0.425, 0.425) to[out=45, in=135] (0.575, 0.425);
        \draw[knot] (0.425, 0.575) to[out=-45, in=-135] (0.575, 0.575);
        \draw[knot] (0.175, 0.675) to[out=45, in=-45] (0.175, 0.825);
        \draw[knot] (0.325, 0.675) to[out=135, in=-135] (0.325, 0.825);
},\, (1,1)\right)}$ as we hoped.
    \item Northeast vertex. From $\varphi_{\left(
\tikz[baseline={([yshift=-.5ex]current bounding box.center)}, scale=0.75]
{
    \draw[dotted] (.5,.5) circle(0.707);
        \draw[knot, red] (0.7, 0.25) -- (0.8, 0.25);
        \draw[knot, red] (0.5, 0.45) -- (0.5, 0.55);
        \draw[knot, red] (0.7, 0.75) -- (0.8, 0.75);
        \draw[knot] (0.5, 1.207) to[out=-90, in=135] (0.675, 0.825);
        \draw[knot] (0.825, 0.675) to[out=-45, in=45] (0.825, 0.325);
        \draw[knot] (0.675, 0.175) to[out=-135, in=90] (0.5, -0.207);
        \draw[knot] (0,0) to (0.425, 0.425);
        \draw[knot] (0.575, 0.575) -- (0.675, 0.675);
        \draw[knot] (0.825, 0.825) -- (1,1);
        \draw[knot] (0,1) -- (0.425, 0.575);
        \draw[knot] (0.575, 0.425) -- (0.675, 0.325);
        \draw[knot] (0.825, 0.175) -- (1,0);
        \draw[knot] (0.675, 0.175) to[out=45, in=-45] (0.675, 0.325);
        \draw[knot] (0.825, 0.175) to[out=135, in=-135] (0.825, 0.325);
        \draw[knot] (0.425, 0.425) to[out=45, in=135] (0.575, 0.425);
        \draw[knot] (0.425, 0.575) to[out=-45, in=-135] (0.575, 0.575);
        \draw[knot] (0.675, 0.675) to[out=45, in=-45] (0.675, 0.825);
        \draw[knot] (0.825, 0.675) to[out=135, in=-135] (0.825, 0.825);
}
, \, (1,1)
\right)
}$, we will consider the representatives $\varphi_{\left(
\tikz[baseline={([yshift=-.5ex]current bounding box.center)}, scale=0.75]
{
    \draw[dotted] (.5,.5) circle(0.707);
        \draw[knot, red] (0.7, 0.25) -- (0.8, 0.25);
        \draw[knot, red] (0.7, 0.75) -- (0.8, 0.75);
        \draw[knot] (0.5, 1.207) to[out=-90, in=135] (0.675, 0.825);
        \draw[knot] (0.825, 0.675) to[out=-45, in=45] (0.825, 0.325);
        \draw[knot] (0.675, 0.175) to[out=-135, in=90] (0.5, -0.207);
        \draw[knot] (0,0) to (0.425, 0.425);
        \draw[knot] (0.575, 0.575) -- (0.675, 0.675);
        \draw[knot] (0.825, 0.825) -- (1,1);
        \draw[knot] (0,1) -- (0.425, 0.575);
        \draw[knot] (0.575, 0.425) -- (0.675, 0.325);
        \draw[knot] (0.825, 0.175) -- (1,0);
        \draw[knot] (0.675, 0.175) to[out=45, in=-45] (0.675, 0.325);
        \draw[knot] (0.825, 0.175) to[out=135, in=-135] (0.825, 0.325);
        \draw[knot] (0.425, 0.425) to[out=45, in=135] (0.575, 0.425);
        \draw[knot] (0.425, 0.575) to[out=-45, in=-135] (0.575, 0.575);
        \draw[knot] (0.675, 0.675) to[out=45, in=-45] (0.675, 0.825);
        \draw[knot] (0.825, 0.675) to[out=135, in=-135] (0.825, 0.825);
}
, \, (1,0)
\right)
}$ and $\varphi_{\left(
\tikz[baseline={([yshift=-.5ex]current bounding box.center)}, scale=0.75]
{
    \draw[dotted] (.5,.5) circle(0.707);
        \draw[knot, red] (0.7, 0.25) -- (0.8, 0.25);
        \draw[knot, red] (0.5, 0.45) -- (0.5, 0.55);
        \draw[knot] (0.5, 1.207) to[out=-90, in=135] (0.675, 0.825);
        \draw[knot] (0.825, 0.675) to[out=-45, in=45] (0.825, 0.325);
        \draw[knot] (0.675, 0.175) to[out=-135, in=90] (0.5, -0.207);
        \draw[knot] (0,0) to (0.425, 0.425);
        \draw[knot] (0.575, 0.575) -- (0.675, 0.675);
        \draw[knot] (0.825, 0.825) -- (1,1);
        \draw[knot] (0,1) -- (0.425, 0.575);
        \draw[knot] (0.575, 0.425) -- (0.675, 0.325);
        \draw[knot] (0.825, 0.175) -- (1,0);
        \draw[knot] (0.675, 0.175) to[out=45, in=-45] (0.675, 0.325);
        \draw[knot] (0.825, 0.175) to[out=135, in=-135] (0.825, 0.325);
        \draw[knot] (0.425, 0.425) to[out=45, in=135] (0.575, 0.425);
        \draw[knot] (0.425, 0.575) to[out=-45, in=-135] (0.575, 0.575);
        \draw[knot] (0.675, 0.675) to[out=45, in=-45] (0.675, 0.825);
        \draw[knot] (0.825, 0.675) to[out=135, in=-135] (0.825, 0.825);
}
, \, (1,0)
\right)
}$. Then,
\[
\varphi_{
\tikz[baseline={([yshift=-.5ex]current bounding box.center)}, scale=0.75]
{
    \draw[dotted] (.5,.5) circle(0.707);
        \draw[knot, red] (0.7, 0.25) -- (0.8, 0.25);
        \draw[knot] (0.5, 1.207) to[out=-90, in=135] (0.675, 0.825);
        \draw[knot] (0.825, 0.675) to[out=-45, in=45] (0.825, 0.325);
        \draw[knot] (0.675, 0.175) to[out=-135, in=90] (0.5, -0.207);
        \draw[knot] (0,0) to (0.425, 0.425);
        \draw[knot] (0.575, 0.575) -- (0.675, 0.675);
        \draw[knot] (0.825, 0.825) -- (1,1);
        \draw[knot] (0,1) -- (0.425, 0.575);
        \draw[knot] (0.575, 0.425) -- (0.675, 0.325);
        \draw[knot] (0.825, 0.175) -- (1,0);
        \draw[knot] (0.675, 0.175) to[out=45, in=-45] (0.675, 0.325);
        \draw[knot] (0.825, 0.175) to[out=135, in=-135] (0.825, 0.325);
        \draw[knot] (0.425, 0.425) to[out=45, in=-45] (0.425, 0.575); 
        \draw[knot] (0.575, 0.425) to[out=135, in=-135] (0.575, 0.575);
        \draw[knot] (0.675, 0.675) to[out=45, in=-45] (0.675, 0.825);
        \draw[knot] (0.825, 0.675) to[out=135, in=-135] (0.825, 0.825);
}
}^{-1}
\circ
\begin{Bmatrix}
    \varphi_{\left(
\tikz[baseline={([yshift=-.5ex]current bounding box.center)}, scale=0.75]
{
    \draw[dotted] (.5,.5) circle(0.707);
        \draw[knot, red] (0.7, 0.25) -- (0.8, 0.25);
        \draw[knot, red] (0.7, 0.75) -- (0.8, 0.75);
        \draw[knot] (0.5, 1.207) to[out=-90, in=135] (0.675, 0.825);
        \draw[knot] (0.825, 0.675) to[out=-45, in=45] (0.825, 0.325);
        \draw[knot] (0.675, 0.175) to[out=-135, in=90] (0.5, -0.207);
        \draw[knot] (0,0) to (0.425, 0.425);
        \draw[knot] (0.575, 0.575) -- (0.675, 0.675);
        \draw[knot] (0.825, 0.825) -- (1,1);
        \draw[knot] (0,1) -- (0.425, 0.575);
        \draw[knot] (0.575, 0.425) -- (0.675, 0.325);
        \draw[knot] (0.825, 0.175) -- (1,0);
        \draw[knot] (0.675, 0.175) to[out=45, in=-45] (0.675, 0.325);
        \draw[knot] (0.825, 0.175) to[out=135, in=-135] (0.825, 0.325);
        \draw[knot] (0.425, 0.425) to[out=45, in=135] (0.575, 0.425);
        \draw[knot] (0.425, 0.575) to[out=-45, in=-135] (0.575, 0.575);
        \draw[knot] (0.675, 0.675) to[out=45, in=-45] (0.675, 0.825);
        \draw[knot] (0.825, 0.675) to[out=135, in=-135] (0.825, 0.825);
}
, \, (1,0)
\right)
}
\\
\varphi_{\left(
\tikz[baseline={([yshift=-.5ex]current bounding box.center)}, scale=0.75]
{
    \draw[dotted] (.5,.5) circle(0.707);
        \draw[knot, red] (0.7, 0.25) -- (0.8, 0.25);
        \draw[knot, red] (0.5, 0.45) -- (0.5, 0.55);
        \draw[knot] (0.5, 1.207) to[out=-90, in=135] (0.675, 0.825);
        \draw[knot] (0.825, 0.675) to[out=-45, in=45] (0.825, 0.325);
        \draw[knot] (0.675, 0.175) to[out=-135, in=90] (0.5, -0.207);
        \draw[knot] (0,0) to (0.425, 0.425);
        \draw[knot] (0.575, 0.575) -- (0.675, 0.675);
        \draw[knot] (0.825, 0.825) -- (1,1);
        \draw[knot] (0,1) -- (0.425, 0.575);
        \draw[knot] (0.575, 0.425) -- (0.675, 0.325);
        \draw[knot] (0.825, 0.175) -- (1,0);
        \draw[knot] (0.675, 0.175) to[out=45, in=-45] (0.675, 0.325);
        \draw[knot] (0.825, 0.175) to[out=135, in=-135] (0.825, 0.325);
        \draw[knot] (0.425, 0.425) to[out=45, in=135] (0.575, 0.425);
        \draw[knot] (0.425, 0.575) to[out=-45, in=-135] (0.575, 0.575);
        \draw[knot] (0.675, 0.675) to[out=45, in=-45] (0.675, 0.825);
        \draw[knot] (0.825, 0.675) to[out=135, in=-135] (0.825, 0.825);
}
, \, (1,0)
\right)
}
\end{Bmatrix}
\cong
\begin{Bmatrix}
    \varphi_{\left(
\tikz[baseline={([yshift=-.5ex]current bounding box.center)}, scale=0.75]
{
    \draw[dotted] (.5,.5) circle(0.707);
        \draw[knot, red] (0.7, 0.75) -- (0.8, 0.75);
        \draw[knot] (0.5, 1.207) to[out=-90, in=135] (0.675, 0.825);
        \draw[knot] (0.825, 0.675) to[out=-45, in=45] (0.825, 0.325);
        \draw[knot] (0.675, 0.175) to[out=-135, in=90] (0.5, -0.207);
        \draw[knot] (0,0) to (0.425, 0.425);
        \draw[knot] (0.575, 0.575) -- (0.675, 0.675);
        \draw[knot] (0.825, 0.825) -- (1,1);
        \draw[knot] (0,1) -- (0.425, 0.575);
        \draw[knot] (0.575, 0.425) -- (0.675, 0.325);
        \draw[knot] (0.825, 0.175) -- (1,0);
        \draw[knot] (0.675, 0.175) to[out=45, in=-45] (0.675, 0.325);
        \draw[knot] (0.825, 0.175) to[out=135, in=-135] (0.825, 0.325);
        \draw[knot] (0.425, 0.425) to[out=45, in=135] (0.575, 0.425);
        \draw[knot] (0.425, 0.575) to[out=-45, in=-135] (0.575, 0.575);
        \draw[knot] (0.675, 0.675) to[out=45, in=-45] (0.675, 0.825);
        \draw[knot] (0.825, 0.675) to[out=135, in=-135] (0.825, 0.825);
}
, \, (1,0)
\right)
}
\\
\varphi_{\left(
\tikz[baseline={([yshift=-.5ex]current bounding box.center)}, scale=0.75]
{
    \draw[dotted] (.5,.5) circle(0.707);
        \draw[knot, red] (0.5, 0.45) -- (0.5, 0.55);
        \draw[knot] (0.5, 1.207) to[out=-90, in=135] (0.675, 0.825);
        \draw[knot] (0.825, 0.675) to[out=-45, in=45] (0.825, 0.325);
        \draw[knot] (0.675, 0.175) to[out=-135, in=90] (0.5, -0.207);
        \draw[knot] (0,0) to (0.425, 0.425);
        \draw[knot] (0.575, 0.575) -- (0.675, 0.675);
        \draw[knot] (0.825, 0.825) -- (1,1);
        \draw[knot] (0,1) -- (0.425, 0.575);
        \draw[knot] (0.575, 0.425) -- (0.675, 0.325);
        \draw[knot] (0.825, 0.175) -- (1,0);
        \draw[knot] (0.675, 0.175) to[out=45, in=-45] (0.675, 0.325);
        \draw[knot] (0.825, 0.175) to[out=135, in=-135] (0.825, 0.325);
        \draw[knot] (0.425, 0.425) to[out=45, in=135] (0.575, 0.425);
        \draw[knot] (0.425, 0.575) to[out=-45, in=-135] (0.575, 0.575);
        \draw[knot] (0.675, 0.675) to[out=45, in=-45] (0.675, 0.825);
        \draw[knot] (0.825, 0.675) to[out=135, in=-135] (0.825, 0.825);
}
, \, (1,0)
\right)
}
\end{Bmatrix}.
\]
The reader is invited to check that both representatives are used in the intermediary complex. Picking the latter and composing with $\varphi_{\left(
    \tikz[baseline={([yshift=-.5ex]current bounding box.center)}, scale=0.75]
{
    \draw[dotted] (.5,.5) circle(0.707);
        \draw[knot, red] (0.25, 0.7) -- (0.25, 0.8);
        \draw[knot] (0.5, 1.207) to[out=-90, in=45] (0.325, 0.825);
        \draw[knot] (0.175, 0.675) to[out=225, in=135] (0.175, 0.325);
        \draw[knot] (0.325, 0.175) to[out=-45, in=90] (0.5, -0.207);
        \draw[knot] (0,0) to (0.175, 0.175);
        \draw[knot] (0.325, 0.325) to (0.425, 0.425);
        \draw[knot] (0.575, 0.575) to (1,1);
        \draw[knot] (0,1) to (0.175, 0.825);
        \draw[knot] (0.325, 0.675) to (0.425, 0.575);
        \draw[knot] (0.575, 0.425) to (1,0);
        \draw[knot] (0.175, 0.175) to[out=45, in=-45] (0.175, 0.325);
        \draw[knot] (0.325, 0.175) to[out=135, in=-135] (0.325, 0.325);
        \draw[knot] (0.425, 0.425) to[out=45, in=-45] (0.425, 0.575); 
        \draw[knot] (0.575, 0.425) to[out=135, in=-135] (0.575, 0.575);
        \draw[knot] (0.175, 0.675) to[out=45, in=135] (0.325, 0.675);
        \draw[knot] (0.175, 0.825) to[out=-45, in=-135] (0.325, 0.825);
}, \, (-1,-1)\right)}^{-1}$, we obtain the grading shift $\{0, -1\}$.
    \item Southeast vertex. This is the most straightforward: applying the first global shift to $\varphi_{\left(
\tikz[baseline={([yshift=-.5ex]current bounding box.center)}, scale=0.75]
{
    \draw[dotted] (.5,.5) circle(0.707);
        \draw[knot, red] (0.7, 0.25) -- (0.8, 0.25);
        \draw[knot] (0.5, 1.207) to[out=-90, in=135] (0.675, 0.825);
        \draw[knot] (0.825, 0.675) to[out=-45, in=45] (0.825, 0.325);
        \draw[knot] (0.675, 0.175) to[out=-135, in=90] (0.5, -0.207);
        \draw[knot] (0,0) to (0.425, 0.425);
        \draw[knot] (0.575, 0.575) -- (0.675, 0.675);
        \draw[knot] (0.825, 0.825) -- (1,1);
        \draw[knot] (0,1) -- (0.425, 0.575);
        \draw[knot] (0.575, 0.425) -- (0.675, 0.325);
        \draw[knot] (0.825, 0.175) -- (1,0);
        \draw[knot] (0.675, 0.175) to[out=45, in=-45] (0.675, 0.325);
        \draw[knot] (0.825, 0.175) to[out=135, in=-135] (0.825, 0.325);
        \draw[knot] (0.425, 0.425) to[out=45, in=135] (0.575, 0.425);
        \draw[knot] (0.425, 0.575) to[out=-45, in=-135] (0.575, 0.575);
        \draw[knot] (0.675, 0.825) to[out=-45, in=-135] (0.825, 0.825);
        \draw[knot] (0.675, 0.675) to[out=45, in=135] (0.825, 0.675);
}, \, (1,0)
\right)
}$ yields a shift by $\{1,0\}$. Redrawing $\tikz[baseline={([yshift=-.5ex]current bounding box.center)}, scale=0.75]
{
    \draw[dotted] (.5,.5) circle(0.707);
        \draw[knot] (0.5, 1.207) to[out=-90, in=135] (0.675, 0.825);
        \draw[knot] (0.825, 0.675) to[out=-45, in=45] (0.825, 0.325);
        \draw[knot] (0.675, 0.175) to[out=-135, in=90] (0.5, -0.207);
        \draw[knot] (0,0) to (0.425, 0.425);
        \draw[knot] (0.575, 0.575) -- (0.675, 0.675);
        \draw[knot] (0.825, 0.825) -- (1,1);
        \draw[knot] (0,1) -- (0.425, 0.575);
        \draw[knot] (0.575, 0.425) -- (0.675, 0.325);
        \draw[knot] (0.825, 0.175) -- (1,0);
        \draw[knot] (0.675, 0.175) to[out=45, in=-45] (0.675, 0.325);
        \draw[knot] (0.825, 0.175) to[out=135, in=-135] (0.825, 0.325);
        \draw[knot] (0.425, 0.425) to[out=45, in=135] (0.575, 0.425);
        \draw[knot] (0.425, 0.575) to[out=-45, in=-135] (0.575, 0.575);
        \draw[knot] (0.675, 0.825) to[out=-45, in=-135] (0.825, 0.825);
        \draw[knot] (0.675, 0.675) to[out=45, in=135] (0.825, 0.675);
}$ as $\tikz[baseline={([yshift=-.5ex]current bounding box.center)}, scale=.75]
{
    \draw[dotted] (.5,.5) circle(0.707);
        \draw[knot] (0.5, 1.207) to[out=-90, in=45] (0.325, 0.825);
        \draw[knot] (0.175, 0.675) to[out=225, in=135] (0.175, 0.325);
        \draw[knot] (0.325, 0.175) to[out=-45, in=90] (0.5, -0.207);
        \draw[knot] (0,0) to (0.175, 0.175);
        \draw[knot] (0.325, 0.325) to (0.425, 0.425);
        \draw[knot] (0.575, 0.575) to (1,1);
        \draw[knot] (0,1) to (0.175, 0.825);
        \draw[knot] (0.325, 0.675) to (0.425, 0.575);
        \draw[knot] (0.575, 0.425) to (1,0);
        \draw[knot] (0.175, 0.175) to[out=45, in=135] (0.325, 0.175);
        \draw[knot] (0.175, 0.325) to[out=-45, in=-135] (0.325, 0.325);
        \draw[knot] (0.425, 0.425) to[out=45, in=135] (0.575, 0.425);
        \draw[knot] (0.425, 0.575) to[out=-45, in=-135] (0.575, 0.575);
        \draw[knot] (0.175, 0.675) to[out=45, in=-45] (0.175, 0.825);
        \draw[knot] (0.325, 0.675) to[out=135, in=-135] (0.325, 0.825);
}$, it is apparent that applying the second global shift provides $\varphi_{\left(
    \tikz[baseline={([yshift=-.5ex]current bounding box.center)}, scale=.75]
{
    \draw[dotted] (.5,.5) circle(0.707);
        \draw[knot, red] (0.2, 0.75) -- (0.3, 0.75);
        \draw[knot] (0.5, 1.207) to[out=-90, in=45] (0.325, 0.825);
        \draw[knot] (0.175, 0.675) to[out=225, in=135] (0.175, 0.325);
        \draw[knot] (0.325, 0.175) to[out=-45, in=90] (0.5, -0.207);
        \draw[knot] (0,0) to (0.175, 0.175);
        \draw[knot] (0.325, 0.325) to (0.425, 0.425);
        \draw[knot] (0.575, 0.575) to (1,1);
        \draw[knot] (0,1) to (0.175, 0.825);
        \draw[knot] (0.325, 0.675) to (0.425, 0.575);
        \draw[knot] (0.575, 0.425) to (1,0);
        \draw[knot] (0.175, 0.175) to[out=45, in=135] (0.325, 0.175);
        \draw[knot] (0.175, 0.325) to[out=-45, in=-135] (0.325, 0.325);
        \draw[knot] (0.425, 0.425) to[out=45, in=135] (0.575, 0.425);
        \draw[knot] (0.425, 0.575) to[out=-45, in=-135] (0.575, 0.575);
        \draw[knot] (0.175, 0.675) to[out=45, in=-45] (0.175, 0.825);
        \draw[knot] (0.325, 0.675) to[out=135, in=-135] (0.325, 0.825);
},\, (1,0)\right)}$, as desired.
\end{enumerate}
\end{proof}

\begin{remark}
In light of Lemma \ref{lem:r3}, the sequence of Remark \ref{rem:R3noticing} is recitfied: notice that
\begin{align*}
\mathrm{Kh}
\left(
\tikz[baseline={([yshift=-.5ex]current bounding box.center)}, scale=.85]
{
        \draw[dotted] (.5,.5) circle(0.707);
        \draw[knot, ->] (0,0) to[out=30, in=-90] (0.8, 0.5);
        \draw[knot] (0.8, 0.5) to[out=90, in=-30] (0,1);
        \draw[knot, overcross] (1,0) to[out=180, in=0] (0.15, 0.7);
        \draw[knot] (0.15, 0.7) to[out=180, in=90] (0, 0.5);
        \draw[knot, ->] (0.15, 0.3) to[out=180, in=-90] (0, 0.5);
        \draw[knot, overcross] (0.15, 0.3) to[out=0, in=180] (1,1);
}
\right)
&
\cong
\varphi_{
\tikz[baseline={([yshift=-.5ex]current bounding box.center)}, scale=.45]
{
        \draw[dotted] (.5,.5) circle(0.707); 
        \draw[knot, red] (0.8, 0.5) -- (0.5, 0.5);
    %
        \draw[knot] (0,0.5) circle(0.2);
        \draw[knot] (1,0) to[out=135, in=-135] (1,1);
        \draw[knot] (0.5, -0.207) -- (0.5, 1.207);
}
}
\circ
\varphi_{
\tikz[baseline={([yshift=-.5ex]current bounding box.center)}, scale=.45]
{
        \draw[dotted] (.5,.5) circle(0.707); 
        \draw[knot, red] (0.2, 0.5) -- (0.5, 0.5);
    %
        \draw[knot] (0,0.5) circle(0.2);
        \draw[knot] (1,0) to[out=135, in=-135] (1,1);
        \draw[knot] (0.5, -0.207) -- (0.5, 1.207);
}
}
^{-1}
\mathrm{Kh}
\left(
\tikz[baseline={([yshift=-.5ex]current bounding box.center)}, scale=.85]
{
        \draw[dotted] (.5,.5) circle(0.707);
        \draw[knot] (0, 0) to[out=30, in=-90] (0.55, 0.5);
        \draw[knot] (0.55, 0.5) to[out=90, in=-30] (0, 1);
        \draw[knot] (1, 0) to[out=150, in=0] (0.5, 0.8);
        \draw[knot, overcross] (0.5, 0.8) to[out=180, in=90] (0.25, 0.5);
        \draw[knot, overcross] (0.25, 0.5) to[out=-90, in=180] (0.5, 0.2);
        \draw[knot, overcross] (0.5, 0.2) to[out=0, in=210] (1, 1);
        \draw[knot, ->] (0.25, 0.499) -- (0.25, 0.501);
        \draw[knot, ->] (0.55, 0.499) -- (0.55, 0.501);
}
\right)
\\
&
\cong
\varphi_{
\left(
\tikz[baseline={([yshift=-.5ex]current bounding box.center)}, scale=.45]
{
        \draw[dotted] (.5,.5) circle(0.707); 
        \draw[knot, red] (0.8, 0.5) -- (0.5, 0.5);
    %
        \draw[knot] (0,0.5) circle(0.2);
        \draw[knot] (1,0) to[out=135, in=-135] (1,1);
        \draw[knot] (0.5, -0.207) -- (0.5, 1.207);
}
, \,
(1,0)
\right)
}
\mathrm{Kh}
\left(
\tikz[baseline={([yshift=-.5ex]current bounding box.center)}, scale=.85]
{
        \draw[dotted] (.5,.5) circle(0.707);
        \draw[knot] (0, 0) to[out=30, in=-90] (0.55, 0.5);
        \draw[knot] (0.55, 0.5) to[out=90, in=-30] (0, 1);
        \draw[knot] (1, 0) to[out=150, in=0] (0.5, 0.8);
        \draw[knot, overcross] (0.5, 0.8) to[out=180, in=90] (0.25, 0.5);
        \draw[knot, overcross] (0.25, 0.5) to[out=-90, in=180] (0.5, 0.2);
        \draw[knot, overcross] (0.5, 0.2) to[out=0, in=210] (1, 1);
        \draw[knot, ->] (0.25, 0.499) -- (0.25, 0.501);
        \draw[knot, ->] (0.55, 0.499) -- (0.55, 0.501);
}
\right).
\end{align*}
Composing with the grading shift by $\{-1,1\}$ and one last Reidemeister I move gives the desired grading shift by $\varphi_{\left(\tikz[baseline=8.25ex, scale=.45]
{
    \begin{scope}[rotate=90]
	\draw[dotted] (3,-2) circle(0.707);
	\draw[knot] (2.5,-1.5) .. controls (2.75,-1.75) and (3.25,-1.75) .. (3.5,-1.5);
	\draw[knot] (2.5,-2.5) .. controls  (2.75,-2.25) and (3.25,-2.25) .. (3.5,-2.5);
        \draw[red, knot] (3,-1.7) -- (3,-2.3);
    \end{scope}
},\, (0,1)\right)}$.
\end{remark}

\begin{theorem}
\label{thm:almosttangleinvt}
If $T$ and $S$ are isotopic diskular tangles, then there exists a grading shifting functor $\varphi_{\Delta^v}$ so that
\[
\varphi_{\Delta^v}\mathrm{Kh}(T) \cong \mathrm{Kh}(S).
\]
\end{theorem}

\begin{proof}
In general, if one decomposes a diskular tangle $T$ into $T_A(T_B)$, as pictured below, then Theorem \ref{thm:multigluing} tells us that $\mathcal{F}(T) \cong \mathcal{F}(T_B) \otimes_{H^n} \mathcal{F}(T_A)$. By Lemma \ref{lem:khcoherence}, there is a shifting functor $\varphi$ and $\ell\in\mathbb{Z}$ such that $\mathrm{Kh}(T) \cong \varphi(\mathcal{F}(T))[\ell]$. By the coherence isomorphisms $\beta$, we have that
\[
\varphi(\mathcal{F}(T_B) \otimes_{H^n} \mathcal{F}(T_A)) \cong \varphi_B\mathcal{F}(T_B) \otimes_{H^n} \varphi_A \mathcal{F}(T_A)
\]
for $\varphi_A$ and $\varphi_B$ restrictions of $\varphi$ to the regions $A$ and $B$. Moreover, as described in the proof of Lemma \ref{lem:khcoherence}, the $\mathscr{G}$-grading and homological-grading shifts here are determined by local crossing information, so it follows similarly that
\[
\varphi(\mathcal{F}(T))[\ell] \cong \varphi_B\mathcal{F}(T_B)[\ell_B] \otimes_{H^n} \varphi_A \mathcal{F}(T_A)[\ell_A]
\]
for those particular $\ell_A$, $\ell_B \in\mathbb{Z}$ satisfying $\ell_A + \ell_B = \ell$. Indeed, since each $\varphi_A$, $\varphi_B$, $\ell_A$, $\ell_B$ coming from $\varphi$ and $\ell$ are the same as the shifts coming from the proof of Lemma \ref{lem:khcoherence}, we have $ \varphi_B\mathcal{F}(T_B)[\ell_B] \cong \mathrm{Kh}(T_B)$ and $\varphi_A \mathcal{F}(T_A)[\ell_A] \cong \mathrm{Kh}(T_A)$. Summarizing, we have that
\[
\mathrm{Kh}(T) \cong \mathrm{Kh}(T_B) \otimes_{H^n} \mathrm{Kh}(T_A).
\]

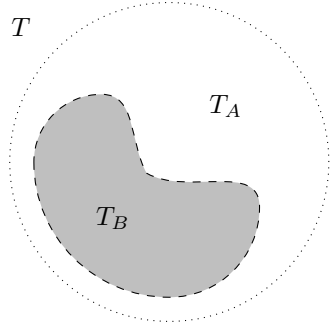
\[
\tikz[baseline={([yshift=-.5ex]current bounding box.center)}, scale=3]
{
    \draw[dotted] (.5,.5) circle(0.707);
    \draw[dashed, fill=gray!50!white] (-0.1,0.5) to[out=-90, in=180] (0.5, -0.1) to[out=0, in=270] (0.9, 0.3) to [out=90, in=-30] (0.4, 0.45) to[out=120, in=0] (0.2, 0.8) to[out=180, in=90] (-0.1,0.5);
    \node at (0.75, 0.75) {$T_A$};
    \node at (0.25, 0.25) {$T_B$};
    \node at (-0.157, 1.1) {$T$};
}
\]

If $T$ and $S$ are isotopic, then $S$ is obtained from $T$ by a finite sequence of Reidemeister moves. For each move in this sequence, apply the isomorphism above to the diskular region containing the Reidemeister move. Then, the theorem follows by applying this isomorphism, invoking one of Lemmas \ref{lem:r1}, \ref{lem:r2}, and \ref{lem:r3}, and repeating as needed.
\end{proof}

\subsubsection{Collapse to $q$-grading}
\label{sss:qcollapse}

To obtain a genuine tangle invariant, we will perform the same trick as is in Section 6.5 of \cite{naisse2020odd}. Define the \textit{degree collapsing map}
\begin{align*}
    \kappa:& \mathrm{Hom}_\mathscr{G} \to \mathbb{Z}\\
    & (D, (p_1, p_2)) \mapsto p_1 + p_2
\end{align*}
which forgets the planar arc diagram input of a $\mathscr{G}$-grading and sums the entries of the second coordinate. We will use $\kappa$ to notice that the $\mathscr{G}$-grading of any $\mathscr{G}$-graded object induces a coarser integral grading. First, by $\mathcal{F}_q(D)$, we mean the multimodule $\mathcal{F}(D)$ with an additional $\mathbb{Z}$-grading determined by its $\mathscr{G}$-grading: fix
\[
\deg_{\mathbb{Z} \times \mathscr{G}} (u) := \left(\kappa(\deg_\mathscr{G}(u)) + \sum_{i=1}^k m_i, \deg_\mathscr{G}(u)\right).
\]
This additional $\mathbb{Z}$-degree, determined by $\mathscr{G}$-degree, is called the \textit{quantum degree}, or \textit{$q$-degree}; we denoted it by $\deg_q(u)$. 

Notice that the composition maps $\mu$ preserve quantum degree. Furthermore, any cobordism $\Delta: Dd \to D'$ induces a map $\mathcal{F}(\Delta) : \mathcal{F}_q(D) \to \mathcal{F}_q(D')$ which is homogeneous of $q$-degree
\[
\deg_q(\mathcal{F}(\Delta)) = \#\text{births} + \#\text{deaths} - \#\text{saddles}.
\]
Sometimes we just write $\deg_q(\Delta)$ for $\deg_q(\mathcal{F}(\Delta))$.

Finally, we reinterpret a grading shift functor $\varphi_{\Delta^v}$ in the $\mathbb{Z} \times \mathscr{G}$-graded setting by
\[
\deg_q\left(\varphi_{\Delta^{(v_1,v_2)}}(m)\right) := \deg_q(m) + \deg_q(\Delta) + v_1 + v_2.
\]
This is to say that any cobordism $\Delta$ induces a $\mathbb{Z} \times \mathscr{G}$-graded map $\mathcal{F}(\Delta) : \varphi_{\Delta}\mathcal{F}_q(D) \to \mathcal{F}_q(D')$.

In conclusion, all results in the $\mathscr{G}$-graded setting extend to the $\mathbb{Z} \times \mathscr{G}$-graded setting with no change to the compatibility maps: all isomorphisms involved are graded with respect to quantum degree. In particular, each $d_{v,j}$ preserves $q$-degree, so we can define $\mathcal{F}_q(T)$ as a $\mathbb{Z} \times \mathscr{G}$-graded dg-multimodule, using $\mathcal{F}_q(T_v)$ in the place of $\mathcal{F}(T_v)$; define $\mathrm{Kh}_q(T)$ similarly.

Suppressing notation, we let $\mathrm{MultiMod}^q$ denote the category whose objects are the same as $\mathrm{MultiMod}^\mathscr{G}$ except we record the quantum degree (that is, objects are $\mathbb{Z}\times \mathscr{G}$-graded multimodules obtained from the regular $\mathscr{G}$-graded ones) but now, maps are only required to be homogeneous with respect to $\mathscr{G}$-degree, with the caveat that they must preserve quantum degree. By \textit{collapsing to $q$-degree}, we just mean that we are working in the category $\mathrm{MultiMod}^q$ rather than $\mathrm{MultiMod}^\mathscr{G}$. This is perhaps misleading, as the $\mathscr{G}$-degree is still present---what we mean to relay is that we have relaxed the requirement of $\mathscr{G}$-degree preservation to $\mathscr{G}$-degree homogeneity up to $q$-degree preservation.

We think of $\mathrm{Kh}_q(T)$ as an object of $\mathrm{Kom}(\mathrm{MultiMod}^q)$. In the final chapter, we are mostly interested in objects of $\mathrm{Kom}(H^n\mathrm{Mod}_R^\mathscr{G})$, which we say \textit{descend} to objects of $\mathrm{Kom}(H^n\mathrm{Mod}_R^q)$, and also to $\mathrm{Kom}(H^n\mathrm{Mod}_e^q)$ and $\mathrm{Kom}(H^n\mathrm{Mod}_o^q)$, specializing $X, Y, Z = 1$ and $X, Z = 1$, $Y=-1$ respectively. We call these objects of $\mathrm{Kom}(H^n \mathrm{Mod}^q)$ the \textit{image} of whatever object(s) of $\mathrm{Kom}(H^n \mathrm{Mod}^\mathscr{G})$ which descends to it.

Notice that a gluing property holds for $\mathcal{F}_q(T)$ and $\mathrm{Kh}_q(T)$, as before. Again, the benefit of working in $\mathrm{Kom}(\mathrm{MultiMod}^q)$ is that $\mathrm{Kh}_q$ becomes an honest tangle invariant.

\begin{theorem}
\label{thm:qinvt}
If $T$ and $S$ are isotopic diskular tangles, then 
\[
\mathrm{Kh}_q(T) \cong \mathrm{Kh}_q(S).
\]
\end{theorem}

\begin{proof}
This follows as long as the homotopy equivalences of Lemmas \ref{lem:r1}, \ref{lem:r2}, and \ref{lem:r3} are graded with respect to quantum degree. For Reidemeister I moves, this is trivial, as the homotopy equivalence was already graded with respect to $\mathscr{G}$-degree. For Reidemeister II moves, $\deg_q(\{-1, 1\}) = 0$ obviously, and 
\[
\deg_q\left(\varphi_{\left(\tikz[baseline=8.25ex, scale=.45]
{
    \begin{scope}[rotate=90]
	\draw[dotted] (3,-2) circle(0.707);
	\draw[knot] (2.5,-1.5) .. controls (2.75,-1.75) and (3.25,-1.75) .. (3.5,-1.5);
	\draw[knot] (2.5,-2.5) .. controls  (2.75,-2.25) and (3.25,-2.25) .. (3.5,-2.5);
        \draw[red, knot] (3,-1.7) -- (3,-2.3);
    \end{scope}
},\, (0,1)\right)}\right) = 1 + \deg_q\left(\tikz[baseline=8.25ex, scale=.45]
{
    \begin{scope}[rotate=90]
	\draw[dotted] (3,-2) circle(0.707);
	\draw[knot] (2.5,-1.5) .. controls (2.75,-1.75) and (3.25,-1.75) .. (3.5,-1.5);
	\draw[knot] (2.5,-2.5) .. controls  (2.75,-2.25) and (3.25,-2.25) .. (3.5,-2.5);
        \draw[red, knot] (3,-1.7) -- (3,-2.3);
    \end{scope}
} \right) = 1 + (-1) = 0.
\]
Similarly, it is clear that the $q$-degree of $\varphi_{
\tikz[baseline={([yshift=-.5ex]current bounding box.center)}, scale=.45]
{
        \draw[dotted] (.5,.5) circle(0.707); 
        \draw[knot, red] (0.2, 0.5) -- (0.5, 0.5);
        \draw[knot] (0,0) to[out=45, in=-45] (0,1);
        \draw[knot] (1,0) to[out=135, in=-135] (1,1);
        \draw[knot] (0.5, -0.207) -- (0.5, 1.207);
}
}
\circ
\varphi_{
\tikz[baseline={([yshift=-.5ex]current bounding box.center)}, scale=.45]
{
        \draw[dotted] (.5,.5) circle(0.707); 
        \draw[knot, red] (0.8, 0.5) -- (0.5, 0.5);
        \draw[knot] (0,0) to[out=45, in=-45] (0,1);
        \draw[knot] (1,0) to[out=135, in=-135] (1,1);
        \draw[knot] (0.5, -0.207) -- (0.5, 1.207);
}
}^{-1}$ is zero. Therefore, the grading shift appearing in Theorem \ref{thm:almosttangleinvt} has $\deg_q(\varphi_{W^v}) = 0$, and the result follows.
\end{proof}

\begin{remark}
If $T$ is a link, then the homology of $\mathrm{Kh}_q(T)$ is isomorphic to the unified Khovanov homology of $T$, as constructed in \cite{putyra20152categorychronologicalcobordismsodd}; see \cite{naisse2020odd} for a proof. In particular, setting $X = Z = 1$ and $Y = -1$ (before taking homology), we get a tangle invariant for odd Khovanov homology, as desired.
\end{remark}

\newpage

\section{Unified and odd projectors}
\label{Chapter:OddCKProj}

Finally, we apply Theorem \ref{thm:multigluing} (multigluing) to mimic the constructions of Cooper-Krushkal \cite{https://doi.org/10.48550/arxiv.1005.5117} and produce projectors living in $\mathrm{Kom}(H^n \mathrm{Mod}_R^\mathscr{G})$. Our work in this chapter follows an outline similar to \cite{stoffregen2024joneswenzlprojectorskhovanovhomotopy}, since we exploit the flexibility provided by diskular tangles, as Stoffregen-Willis do in the spectral setting. 

More explicitly, in \S \ref{ss:operationsviamultigluing}, we use multigluing to define stacking $\otimes$, juxtaposing $\sqcup$, and partial trace $\mathrm{Tr}$ operations, and the category $\mathrm{Chom}(n)^\mathscr{G}$ (which we conjecture is the same as $\mathrm{Kom}(H^n\mathrm{PMod}^\mathscr{G})$, as in \cite{Khovanov_2002}). We also take this opportunity to prove an adjunction, generalizing a theorem of Hogancamp \cite{https://doi.org/10.48550/arxiv.1405.2574}. The next section, \S \ref{ss:duality}, is mostly stand-alone: the main takeaway for this paper is Corollary \ref{cor:dualswap}, which we use in the proof of Lemma \ref{lem:pullbackHE}, itself used in the proofs of Proposition \ref{prop:unifiedprojectoruni} and Corollary \ref{cor:endothing}. In \S \ref{s:unifiedprops}, we define unified projectors as in \cite{https://doi.org/10.48550/arxiv.1405.2574}, though our proofs follow the methods outlined in \cite{stoffregen2024joneswenzlprojectorskhovanovhomotopy}, as their setting most resembles our own. We hope to illuminate preceding and successive work by computing the 2-stranded unified projector two different ways in \S \ref{ss:2strandedproj}. We also compute the homology of the closure of $P_2$ (cf. Section 4.3.1 of \cite{https://doi.org/10.48550/arxiv.1005.5117}), which we will use to show that our categorification of the colored Jones polynomial is distinct from that of \cite{https://doi.org/10.48550/arxiv.1005.5117}. Finally, we prove the existence of unified projectors (using the same procedure as \cite{stoffregen2024joneswenzlprojectorskhovanovhomotopy}) in \S \ref{ss:egunifiedprojectors}, and the existence of a unified colored link homology (which collapses to the categorification of the colored Jones polynomial of \cite{https://doi.org/10.48550/arxiv.1005.5117} on one hand, and to a new categorification on the other) in \S \ref{ss:unifiedcoloredhomology}.

We establish some notation. Proceeding, for $A, B \in \mathrm{Kom}(H^n\mathrm{Mod}^\mathscr{G})$, we will denote the HOM-complex of $A$ and $B$ by $\normalfont{\textsc{Hom}}_n(A, B)$. If $A$ and $B$ are (non-dg) $\mathscr{G}$-graded $H^n$ modules, we'll write $\mathrm{Hom}_n(A, B)$ as shorthand for $\mathrm{Hom}_{H^nMod^\mathscr{G}}(A,B)$.

\subsection{Operations defined via multigluing}
\label{ss:operationsviamultigluing}

As far as the existence of projectors is concerned, the main payoff of multigluing in the unified setting is that we can develop a notion for stacking and juxtaposing complexes of $\mathscr{G}$-graded modules. We can also use multigluing to define a partial trace for these complexes, allowing for an adjunction statement.

Given a diskular $n$-tangle $T$, we'll view it as a tangle in a rectangle as follows: traveling counter-clockwise from the basepoint along the boundary, place the first $n$ endpoints along the top of the rectangle and the last $n$ endpoints along the bottom. For this reason, flat diskular $n$-tangles are also called a Temperley-Lieb $n$-diagrams, i.e., each resolution of $T$ is a Temperley-Lieb diagram. 

\begin{definition}[Stacking]
\label{def:verttensor}
\textit{Vertical composition} is the operation 
\[
\otimes: \mathrm{Kom}(H^n\mathrm{Mod}_R^\mathscr{G}) \times \mathrm{Kom}(H^n\mathrm{Mod}_R^\mathscr{G}) \to \mathrm{Kom}(H^n\mathrm{Mod}_R^\mathscr{G})
\]
defined as follows: given complexes $A, B \in \mathrm{Kom}(H^n\mathrm{Mod}_R^\mathscr{G})$, $A\otimes B$ is the complex
\[
(A, B) \otimes_{\left(H^{n}, H^{n}\right)} \mathcal{F}(D_n^\otimes)
\]
where $D_n^\otimes$ is the $(n,n; n)$-planar arc diagram 
 
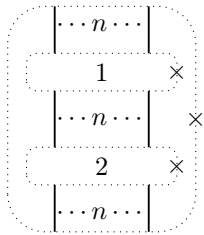
\[
\tikz{
    \draw[rounded corners=5mm, dotted] (0, 0) rectangle (2.5, 3) {};
    \draw[rounded corners, dotted] (.25, 1.875) rectangle (2.25, 2.375);
    \draw[rounded corners, dotted] (.25, .625) rectangle (2.25, 1.125);
    \node at (1.25, 2.125) {1};
    \node at (1.25, 0.875) {2};
    \draw[knot] (0.625, 0) -- (0.625, .625);
    \draw[knot] (0.625, 1.125) -- (0.625, 1.875);
    \draw[knot] (0.625, 2.375) -- (0.625, 3);
%
%
    \draw[knot] (1.875, 0) -- (1.875, .625);
    \draw[knot] (1.875, 1.125) -- (1.875, 1.875);
    \draw[knot] (1.875, 2.375) -- (1.875, 3);
    \node at (2.25, 2.125) {$\times$};
    \node at (2.25, 0.875) {$\times$};
    \node at (2.5, 1.5) {$\times$};
    \node at (1.25, 2.75) {$\cdots n \cdots$};
    \node at (1.25, 0.25) {$\cdots n \cdots$};
    \node at (1.25, 1.5) {$\cdots n \cdots$}; 
}
\]
with removed disks ordered as shown. In particular, if $T_1$ and $T_2$ are both diskular $n$-tangles, Theorem \ref{thm:multigluing} says that
\[
\mathcal{F}(T_1) \otimes \mathcal{F}(T_2) \cong
(\mathcal{F}(T_1), \mathcal{F}(T_2)) \otimes_{\left(H^{n}, H^{n}\right)} \mathcal{F}(D_n^\otimes) \cong \mathcal{F}(D_n^\otimes(T_1, T_2)).
\]
We say that this complex is the result of \textit{stacking} $\mathcal{F}(T_1)$ and $\mathcal{F}(T_2)$.
\end{definition}

\begin{definition}
\label{def:chom}
Consider the full subcategory $\mathrm{Chom}(n)^\mathscr{G}$ of $\mathrm{Kom}(H^n\mathrm{Mod}^\mathscr{G})$ consisting of (partially unbounded) $\mathscr{G}$-graded dg-modules whose entries are all direct sums of $\mathscr{G}$-graded modules associated to flat diskular $n$-tangles. 

In analogy with \cite{10.1215/S0012-7094-00-10131-7}, we expect that the subcategory $\mathrm{Chom}(n)^\mathscr{G}$ is just the category $\mathrm{Kom}(H^n\mathrm{PMod}^\mathscr{G})$ for $H^n\mathrm{PMod}^\mathscr{G}$ the category of projective $\mathscr{G}$-graded $H^n$-modules, although this seems worthy of further study. Additionally, we expect that vertical composition $\otimes$ for this subcategory is a monoidal product with monoidal identity
\[
\mathcal{I}_n := \tikz[baseline={([yshift=-.5ex]current bounding box.center)}, scale=1.5]{
    \draw[dotted, rounded corners] (0,0) rectangle (1.5,1);
    \draw[knot] (0.25, 0) -- (0.25, 1);
    \draw[knot] (1.25, 0) -- (1.25, 1);
    \node at (0.75, 0.5) {$\cdots n \cdots$};
    \node at (1.5, 0.5) {$\times$};
}
\]
(that is, $\mathcal{I}_n$ is the dg-module associated to the picture above), with monoidal structure provided by multigluing (Theorem \ref{thm:multigluing}). We let $\mathrm{Chom}(n)^q$ denote the image of $\mathrm{Chom}(n)^\mathscr{G}$ in $\mathrm{Kom}(H^n\mathrm{Mod}^q)$ after collapsing to $q$-degree, \S \ref{sss:qcollapse}. By definition, for $K_0^q$ the Grothendieck group which records only the $q$-degree of $\mathscr{G}$-graded objects, we have that 
\[
K_0^q(\mathrm{Chom}(n)^\mathscr{G}) \cong K_0^q(\mathrm{Chom}(n)^q) \cong TL_n.
\]
\end{definition}

Just as stacking can be realized as a multigluing operation, the horizontal juxtaposition can as well.

\begin{definition}[Juxtaposing]
\textit{Horizontal composition} is the operation 
\[
\sqcup: \mathrm{Kom}(H^{n_1}\mathrm{Mod}_R^\mathscr{G}) \times \mathrm{Kom}(H^{n_2}\mathrm{Mod}_R^\mathscr{G}) \to \mathrm{Kom}(H^{n_1+n_2}\mathrm{Mod}_R^\mathscr{G})
\]
defined as follows: for complexes $A \in \mathrm{Kom}(H^{n_1}\mathrm{Mod}_R^\mathscr{G})$ and $B \in \mathrm{Kom}(H^{n_2}\mathrm{Mod}_R^\mathscr{G})$, $A\sqcup B$ is the complex
\[
(\mathcal{F}(T_1), \mathcal{F}(T_2)) \otimes_{\left(H^{n_1}, H^{n_2}\right)} \mathcal{F}(D_{(n_1, n_2)}^\sqcup)
\]
where $D_{(n_1, n_2)}^\sqcup$ is the $(n_1, n_2; n_1 + n_2)$-planar arc diagram
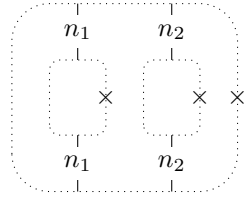
\[
\tikz{
    \draw[rounded corners=5mm, dotted] (0, 0) rectangle (3, 2.5) {};
    \draw[rounded corners, dotted] (0.5, 0.75) rectangle (1.25, 1.75);
    \draw[rounded corners, dotted] (1.75, 0.75) rectangle (2.5, 1.75);
    \draw (0.875, 0) -- (0.875, 0.75);
    \draw (0.875, 1.75) -- (0.875, 2.5);
    \draw (2.125, 0) -- (2.125, 0.75);
    \draw (2.125, 1.75) -- (2.125, 2.5);
    \node at (1.25, 1.25) {$\times$};
    \node at (2.5, 1.25) {$\times$};
    \node at (3, 1.25) {$\times$};
    \node[fill=white] at (0.875, 0.375) {$n_1$};
    \node[fill=white] at (0.875, 2.125) {$n_1$};
    \node[fill=white] at (2.125, 0.375) {$n_2$};
    \node[fill=white] at (2.125, 2.125) {$n_2$};
}
\]
If $T_i$ a diskular $n_i$-tangle, we'll write $\mathcal{F}(T_1) \sqcup \mathcal{F}(T_2)$ to denote the tensor product 
\[
(\mathcal{F}(T_1), \mathcal{F}(T_2)) \otimes_{\left(H^{n_1}, H^{n_2}\right)} \mathcal{F}(D_{(n_1, n_2)}^\sqcup) \cong \mathcal{F}(D_{(n_1, n_2)}^\sqcup(T_1, T_2)).
\]
We say that this complex is the result of \textit{juxtaposing} $\mathcal{F}(T_1)$ and $\mathcal{F}(T_2)$.
\end{definition}

\subsubsection{Adjunction}

First, consider the following operation on complexes in $\mathrm{Kom}(H^n\mathrm{Mod}_R^\mathscr{G})$.

\begin{definition}(Trace)
The \textit{trace} is an operation
\[
\mathrm{Tr}: \mathrm{Kom}(H^n \mathrm{Mod}_R^\mathscr{G}) \to \mathrm{Kom}(H^{n-1} \mathrm{Mod}_R^\mathscr{G})
\]
defined as follows: for $A \in \mathrm{Kom}(H^n \mathrm{Mod}_R^\mathscr{G})$, $\mathrm{Tr}(A)$ is the complex
\[
A \otimes_{H^n} \mathcal{F}(D_n^{\mathrm{Tr}})
\]
where $D_n^{\mathrm{Tr}}$ is the $(n; n-1)$-planar arc diagram
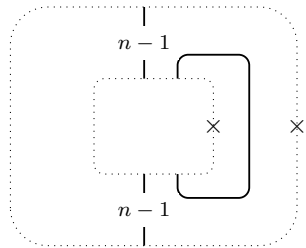
\[
\tikz[baseline={([yshift=-.5ex]current bounding box.center)}, scale=1.27]{
    \draw[rounded corners=5mm, dotted] (0, 0) rectangle (3, 2.5) {};
    \draw[rounded corners, dotted] (0.875, 0.75) rectangle (2.125, 1.75);
    \draw[knot] (1.4, 0) -- (1.4, 0.75);
    \draw[knot] (1.4, 1.75) -- (1.4, 2.5);
    \draw[knot, rounded corners] (1.75, 0.75) -- (1.75, 0.5) -- (2.5, 0.5) -- (2.5, 2) -- (1.75, 2) -- (1.75, 1.75);
    \node at (2.125, 1.25) {$\times$};
    \node at (3, 1.25) {$\times$};
    \node[fill=white] at (1.4, 0.375) {\footnotesize$n-1$};
    \node[fill=white] at (1.4, 2.125) {\footnotesize$n-1$};
}
\]
If $T$ is a diskular $n$-tangle, we'll write $\mathrm{Tr}(\mathcal{F}(T))$ to denote the complex
\[
\mathcal{F}(T) \otimes_{H^n} \mathcal{F}(D_n^{\mathrm{Tr}}) \cong \mathcal{F}(D_n^{\mathrm{Tr}}(T)).
\]
By the \textit{$k$th partial trace} of $A$, we mean the complex obtained from applying the partial trace $k$ times to obtain $\mathrm{Tr}^k(A) \in \mathrm{Kom}(H^{n-k}\mathrm{Mod}_R^\mathscr{G})$. The $n$th partial trace of $A$ is known simply as the \textit{trace} or \textit{closure} of $A$.
\end{definition}

In \cite{https://doi.org/10.48550/arxiv.1405.2574}, we saw that the operations $-\sqcup 1$ and $\mathrm{Tr}(-)$ were adjoint. Impressively, we can prove that a generalization of this adjunction exists in the $\mathscr{G}$-graded setting!

\begin{theorem}
\label{thm:adjunction1}
Suppose $A\in \mathrm{Kom}(H^{n-1}\mathrm{Mod}_R^\mathscr{G})$ and $B \in \mathrm{Kom}(H^n \mathrm{Mod}_R^\mathscr{G})$. Then we have the following isomorphisms of complexes.
\[
\normalfont{\textsc{Hom}}_n\left(\tikz[baseline={([yshift=-.5ex]current bounding box.center)}, scale=.3]{
            \draw[knot] (0,0) rectangle (3, 1.5);
            \draw[knot] (0.5, 0) -- (0.5, -1.25);
            \draw[knot] (1.5, 0) -- (1.5, -1.25);
            \draw[knot] (2.5, 0) -- (2.5, -1.25);
            \draw[knot] (0.5, 1.5) -- (0.5, 2.75);
            \draw[knot] (1.5, 1.5) -- (1.5, 2.75);
            \draw[knot] (2.5, 1.5) -- (2.5, 2.75);
            \draw[knot] (3.5, -1.25) -- (3.5, 2.75);
            \node at (1.5, 0.75) {$A$};
        }\,,~ \varphi_{\left(\tikz[baseline={([yshift=-.5ex]current bounding box.center)}, scale=.15]{
            \draw[fill=white, knot] (0,0) rectangle (4, 1.5);
            \draw[knot, red, rounded corners] (3.5, 2.125) -- (4.5, 2.125) -- (4.5, -0.625) -- (3.5, -0.625);
            \draw[knot] (0.5, 0) -- (0.5, -1.25);
            \draw[knot] (1.5, 0) -- (1.5, -1.25);
            \draw[knot] (2.5, 0) -- (2.5, -1.25);
            \draw[knot] (0.5, 1.5) -- (0.5, 2.75);
            \draw[knot] (1.5, 1.5) -- (1.5, 2.75);
            \draw[knot] (2.5, 1.5) -- (2.5, 2.75);
            \draw[knot] (3.5, 0) -- (3.5, -1.25);
            \draw[knot] (3.5, 1.5) -- (3.5, 2.75);
            \node at (2, 0.75) {\scriptsize \textit{B}};
        },\, (0,1)\right)} \tikz[baseline={([yshift=-.5ex]current bounding box.center)}, scale=.3]{
            \draw[fill=white, knot] (0,0) rectangle (4, 1.5);
            \draw[knot] (0.5, 0) -- (0.5, -1.25);
            \draw[knot] (1.5, 0) -- (1.5, -1.25);
            \draw[knot] (2.5, 0) -- (2.5, -1.25);
            \draw[knot] (0.5, 1.5) -- (0.5, 2.75);
            \draw[knot] (1.5, 1.5) -- (1.5, 2.75);
            \draw[knot] (2.5, 1.5) -- (2.5, 2.75);
            \draw[knot] (3.5, 0) -- (3.5, -1.25);
            \draw[knot] (3.5, 1.5) -- (3.5, 2.75);
            \node at (2, 0.75) {$B$};
        } \right)
        \cong
        \normalfont{\textsc{Hom}}_{n-1}\left( \tikz[baseline={([yshift=-.5ex]current bounding box.center)}, scale=.3]{
            \draw[knot] (0,0) rectangle (3, 1.5);
            \draw[knot] (0.5, 0) -- (0.5, -1.25);
            \draw[knot] (1.5, 0) -- (1.5, -1.25);
            \draw[knot] (2.5, 0) -- (2.5, -1.25);
            \draw[knot] (0.5, 1.5) -- (0.5, 2.75);
            \draw[knot] (1.5, 1.5) -- (1.5, 2.75);
            \draw[knot] (2.5, 1.5) -- (2.5, 2.75);
            \node at (1.5, 0.75) {$A$};
        }\,,~ \tikz[baseline={([yshift=-.5ex]current bounding box.center)}, scale=.3]{
            \draw[rounded corners, knot] (3.5, -1) rectangle (4.5, 2.5);
            \draw[fill=white, knot] (0,0) rectangle (4, 1.5);
            \draw[knot] (0.5, 0) -- (0.5, -1.25);
            \draw[knot] (1.5, 0) -- (1.5, -1.25);
            \draw[knot] (2.5, 0) -- (2.5, -1.25);
            \draw[knot] (0.5, 1.5) -- (0.5, 2.75);
            \draw[knot] (1.5, 1.5) -- (1.5, 2.75);
            \draw[knot] (2.5, 1.5) -- (2.5, 2.75);
            \node at (2, 0.75) {$B$};
        }\, \{-1,0\}\right)
\]
and 
\[
\normalfont{\textsc{Hom}}_n\left(
\varphi_{\left(\tikz[baseline={([yshift=-.5ex]current bounding box.center)}, scale=.15]{
            \draw[fill=white, knot] (0,0) rectangle (4, 1.5);
            \draw[knot, red, rounded corners] (3.5, 2.125) -- (4.5, 2.125) -- (4.5, -0.625) -- (3.5, -0.625);
            \draw[knot] (0.5, 0) -- (0.5, -1.25);
            \draw[knot] (1.5, 0) -- (1.5, -1.25);
            \draw[knot] (2.5, 0) -- (2.5, -1.25);
            \draw[knot] (0.5, 1.5) -- (0.5, 2.75);
            \draw[knot] (1.5, 1.5) -- (1.5, 2.75);
            \draw[knot] (2.5, 1.5) -- (2.5, 2.75);
            \draw[knot] (3.5, 0) -- (3.5, -1.25);
            \draw[knot] (3.5, 1.5) -- (3.5, 2.75);
            \node at (2, 0.75) {\scriptsize \textit{B}};
        },\, (0,1)\right)} \tikz[baseline={([yshift=-.5ex]current bounding box.center)}, scale=.3]{
            \draw[fill=white, knot] (0,0) rectangle (4, 1.5);
            \draw[knot] (0.5, 0) -- (0.5, -1.25);
            \draw[knot] (1.5, 0) -- (1.5, -1.25);
            \draw[knot] (2.5, 0) -- (2.5, -1.25);
            \draw[knot] (0.5, 1.5) -- (0.5, 2.75);
            \draw[knot] (1.5, 1.5) -- (1.5, 2.75);
            \draw[knot] (2.5, 1.5) -- (2.5, 2.75);
            \draw[knot] (3.5, 0) -- (3.5, -1.25);
            \draw[knot] (3.5, 1.5) -- (3.5, 2.75);
            \node at (2, 0.75) {$B$};
        } 
        \, ,~
\tikz[baseline={([yshift=-.5ex]current bounding box.center)}, scale=.3]{
            \draw[knot] (0,0) rectangle (3, 1.5);
            \draw[knot] (0.5, 0) -- (0.5, -1.25);
            \draw[knot] (1.5, 0) -- (1.5, -1.25);
            \draw[knot] (2.5, 0) -- (2.5, -1.25);
            \draw[knot] (0.5, 1.5) -- (0.5, 2.75);
            \draw[knot] (1.5, 1.5) -- (1.5, 2.75);
            \draw[knot] (2.5, 1.5) -- (2.5, 2.75);
            \draw[knot] (3.5, -1.25) -- (3.5, 2.75);
            \node at (1.5, 0.75) {$A$};
        }\right)
        \cong
        \normalfont{\textsc{Hom}}_{n-1}\left( \tikz[baseline={([yshift=-.5ex]current bounding box.center)}, scale=.3]{
            \draw[rounded corners, knot] (3.5, -1) rectangle (4.5, 2.5);
            \draw[fill=white, knot] (0,0) rectangle (4, 1.5);
            \draw[knot] (0.5, 0) -- (0.5, -1.25);
            \draw[knot] (1.5, 0) -- (1.5, -1.25);
            \draw[knot] (2.5, 0) -- (2.5, -1.25);
            \draw[knot] (0.5, 1.5) -- (0.5, 2.75);
            \draw[knot] (1.5, 1.5) -- (1.5, 2.75);
            \draw[knot] (2.5, 1.5) -- (2.5, 2.75);
            \node at (2, 0.75) {$B$};
        } 
        \, ,~
        \tikz[baseline={([yshift=-.5ex]current bounding box.center)}, scale=.3]{
            \draw[knot] (0,0) rectangle (3, 1.5);
            \draw[knot] (0.5, 0) -- (0.5, -1.25);
            \draw[knot] (1.5, 0) -- (1.5, -1.25);
            \draw[knot] (2.5, 0) -- (2.5, -1.25);
            \draw[knot] (0.5, 1.5) -- (0.5, 2.75);
            \draw[knot] (1.5, 1.5) -- (1.5, 2.75);
            \draw[knot] (2.5, 1.5) -- (2.5, 2.75);
            \node at (1.5, 0.75) {$A$};
        }\,
        \{0,-1\}\right).
\]
\end{theorem}

\begin{proof}
Unlike the analogues of this result for even Khovanov homology \cite{https://doi.org/10.48550/arxiv.1405.2574} and even Khovanov spectra \cite{stoffregen2024joneswenzlprojectorskhovanovhomotopy}, the fact that certain maps occur in disjoint disks does not mean that they commute, but rather that swapping the two changes the overall composition by an isomorphism induced by a locally vertical change of chronology. We will see that our $\mathscr{G}$-shifting 2-system accounts for this difference, so that the above result holds with little alterations to the aforementioned proofs.

We will prove the first isomorphism, leaving the second to the reader---notice that the grading shift by $\{-1,0\}$ in the former is replaced by a grading shift by $\{0,-1\}$ in the latter. Suppose that $f\in \normalfont{\textsc{Hom}}_n\left(A \sqcup 1, \varphi_{\left(\tikz[baseline={([yshift=-.5ex]current bounding box.center)}, scale=.15]{
            \draw[fill=white, knot] (0,0) rectangle (4, 1.5);
            \draw[knot, red, rounded corners] (3.5, 2.125) -- (4.5, 2.125) -- (4.5, -0.625) -- (3.5, -0.625);
            \draw[knot] (0.5, 0) -- (0.5, -1.25);
            \draw[knot] (1.5, 0) -- (1.5, -1.25);
            \draw[knot] (2.5, 0) -- (2.5, -1.25);
            \draw[knot] (0.5, 1.5) -- (0.5, 2.75);
            \draw[knot] (1.5, 1.5) -- (1.5, 2.75);
            \draw[knot] (2.5, 1.5) -- (2.5, 2.75);
            \draw[knot] (3.5, 0) -- (3.5, -1.25);
            \draw[knot] (3.5, 1.5) -- (3.5, 2.75);
            \node at (2, 0.75) {\scriptsize \textit{B}};
        },\, (0,1)\right)} B\right)$ has homogenous $\mathscr{I}$-degree $\Delta^v$, so it is realized as a $\mathscr{G}$-graded map $f: \varphi_{\Delta^v} A \sqcup 1 \to \varphi_{\left(\tikz[baseline={([yshift=-.5ex]current bounding box.center)}, scale=.15]{
            \draw[fill=white, knot] (0,0) rectangle (4, 1.5);
            \draw[knot, red, rounded corners] (3.5, 2.125) -- (4.5, 2.125) -- (4.5, -0.625) -- (3.5, -0.625);
            \draw[knot] (0.5, 0) -- (0.5, -1.25);
            \draw[knot] (1.5, 0) -- (1.5, -1.25);
            \draw[knot] (2.5, 0) -- (2.5, -1.25);
            \draw[knot] (0.5, 1.5) -- (0.5, 2.75);
            \draw[knot] (1.5, 1.5) -- (1.5, 2.75);
            \draw[knot] (2.5, 1.5) -- (2.5, 2.75);
            \draw[knot] (3.5, 0) -- (3.5, -1.25);
            \draw[knot] (3.5, 1.5) -- (3.5, 2.75);
            \node at (2, 0.75) {\scriptsize \textit{B}};
        },\, (0,1)\right)} B$ (we do not have to pay attention to the homological degree). Define $\phi(f) \in \normalfont{\textsc{Hom}}_{n-1}(A, \mathrm{Tr}(B)\{-1,0\})$ to be the composition
\[
\tikz{
     \node(A) at (0,0) 
        {$\varphi_{\Delta^v}
        \tikz[baseline={([yshift=-.5ex]current bounding box.center)}, scale=.35]{
            \draw[knot] (0,0) rectangle (3, 1.5);
            \draw[knot] (0.5, 0) -- (0.5, -1.25);
            \draw[knot] (1.5, 0) -- (1.5, -1.25);
            \draw[knot] (2.5, 0) -- (2.5, -1.25);
            \draw[knot] (0.5, 1.5) -- (0.5, 2.75);
            \draw[knot] (1.5, 1.5) -- (1.5, 2.75);
            \draw[knot] (2.5, 1.5) -- (2.5, 2.75);
            \node at (1.5, 0.75) {$A$};
            \draw[rounded corners, dashed, red] (3.5, -1) rectangle (4.5, 2.5);
        }$};
    \node(B) at (7,0) 
        {$\varphi_{\Delta^{v + (-1,0)}}
        \tikz[baseline={([yshift=-.5ex]current bounding box.center)}, scale=.35]{
            \draw[knot] (0,0) rectangle (3, 1.5);
            \draw[knot] (0.5, 0) -- (0.5, -1.25);
            \draw[knot] (1.5, 0) -- (1.5, -1.25);
            \draw[knot] (2.5, 0) -- (2.5, -1.25);
            \draw[knot] (0.5, 1.5) -- (0.5, 2.75);
            \draw[knot] (1.5, 1.5) -- (1.5, 2.75);
            \draw[knot] (2.5, 1.5) -- (2.5, 2.75);
            \node at (1.5, 0.75) {$A$};
            \draw[rounded corners, knot] (3.5, -1) rectangle (4.5, 2.5);
            \draw[rounded corners, dashed, red] (-.5, -.5) rectangle (4, 2);
        }$};
    \node(C) at (12,0) 
        {$\tikz[baseline={([yshift=-.5ex]current bounding box.center)}, scale=.35]{
            \draw[rounded corners, knot] (3.5, -1) rectangle (4.5, 2.5);
            \draw[fill=white, knot] (0,0) rectangle (4, 1.5);
            \draw[knot] (0.5, 0) -- (0.5, -1.25);
            \draw[knot] (1.5, 0) -- (1.5, -1.25);
            \draw[knot] (2.5, 0) -- (2.5, -1.25);
            \draw[knot] (0.5, 1.5) -- (0.5, 2.75);
            \draw[knot] (1.5, 1.5) -- (1.5, 2.75);
            \draw[knot] (2.5, 1.5) -- (2.5, 2.75);
            \node at (2, 0.75) {$B$};
        }\, \{-1,0\}\,$};
    \draw[->] (A) to[out=0,in=180] node[pos=.5,above,arrows=-]
        {$\lambda_{\phi(f)} \circ \mathcal{F}\left(\tikz[baseline={([yshift=-.5ex]current bounding box.center)}, scale=.4]{
	\draw (1,2) .. controls (1,1) and (2,1) .. (2,2);
	\draw (1,2) .. controls (1,1.75) and (2,1.75) .. (2,2);
	\draw (1,2) .. controls (1,2.25) and (2,2.25) .. (2,2);
}\right) \circ \varphi_{H_B}$} 
        (B);
    \draw[->] (B) to[out=0,in=180] node[pos=.5,above,arrows=-]
        {$\mathrm{Tr}(f)$} 
        (C);    
}
\]
where
\[
\varphi_{H_B}: \mathrm{Id} \Rightarrow \varphi_{\tikz[baseline={([yshift=-.5ex]current bounding box.center)}, scale=.4]{
	\draw (1,2) .. controls (1,1) and (2,1) .. (2,2);
	\draw (1,2) .. controls (1,1.75) and (2,1.75) .. (2,2);
	\draw (1,2) .. controls (1,2.25) and (2,2.25) .. (2,2);
}}^{-1} \circ \varphi_{\tikz[baseline={([yshift=-.5ex]current bounding box.center)}, scale=.4]{
	\draw (1,2) .. controls (1,1) and (2,1) .. (2,2);
	\draw (1,2) .. controls (1,1.75) and (2,1.75) .. (2,2);
	\draw (1,2) .. controls (1,2.25) and (2,2.25) .. (2,2);
}}^{\textcolor{white}{-1}} \cong \{-1,0\} \circ \{1,0\}.
\]
and $\lambda_{\phi(f)}$ is shorthand for the isomorphism which pushes the $\{-1,0\}$ shift after $\Delta^v$; that is, 
$\lambda_{\phi(f)} = \gamma_{(-1,0), \Delta^v} \circ \lambda(v, (-1,0))$. Schematically, 
\[
\lambda_{\phi(f)}: \tikz[baseline={([yshift=-.5ex]current bounding box.center)}, yscale=.5,xscale=.75]{
    \draw (0,-2.5) -- (0,2.5);
    \node[fill=white,draw,rounded corners,scale=1] at (0,-1.25) {$(-1,0)$};
    \node[fill=white,draw,rounded corners,scale=1] at (0,0) {$\Delta$};
    \node[fill=white,draw,rounded corners,scale=1] at (0,1.25) {$v$};
}
\Rightarrow
\tikz[baseline={([yshift=-.5ex]current bounding box.center)}, yscale=.5,xscale=.75]{
    \draw (0,-2.5) -- (0,2.5);
    \node[fill=white,draw,rounded corners,scale=1] at (0,-1.25) {$\Delta$};
    \node[fill=white,draw,rounded corners,scale=1] at (0,0) {$v$};
    \node[fill=white,draw,rounded corners,scale=1] at (0,1.25) {$(-1,0)$};
}.
\]
Lastly, by $\mathrm{Tr}(f)$ we just mean $f \otimes \mathbbm{1}_{D_n^{\mathrm{Tr}}}$. Notice that $\phi(f)$ has the desired form since
\[
\tikz[baseline={([yshift=-.5ex]current bounding box.center)}, scale=.3]{
            \draw[fill=white, knot] (0,0) rectangle (4, 1.5);
            \draw[knot, red, rounded corners] (3.5, 2.125) -- (4.5, 2.125) -- (4.5, -0.625) -- (3.5, -0.625);
            \draw[knot] (0.5, 0) -- (0.5, -1.25);
            \draw[knot] (1.5, 0) -- (1.5, -1.25);
            \draw[knot] (2.5, 0) -- (2.5, -1.25);
            \draw[knot] (0.5, 1.5) -- (0.5, 2.75);
            \draw[knot] (1.5, 1.5) -- (1.5, 2.75);
            \draw[knot] (2.5, 1.5) -- (2.5, 2.75);
            \draw[knot] (3.5, 0) -- (3.5, -1.25);
            \draw[knot] (3.5, 1.5) -- (3.5, 2.75);
            \node at (2, 0.75) {\textit{B}};
        } \bullet \mathbbm{1}_{D_{n}^{\mathrm{Tr}}} = \tikz[baseline={([yshift=-.5ex]current bounding box.center)}, scale=.3]{
            \draw[fill=white, knot] (0,0) rectangle (4, 1.5);
            \draw[knot, red, rounded corners] (3.5, 2.125) -- (4.5, 2.125) -- (4.5, -0.625) -- (3.5, -0.625);
            \draw[knot] (0.5, 0) -- (0.5, -1.25);
            \draw[knot] (1.5, 0) -- (1.5, -1.25);
            \draw[knot] (2.5, 0) -- (2.5, -1.25);
            \draw[knot] (0.5, 1.5) -- (0.5, 2.75);
            \draw[knot] (1.5, 1.5) -- (1.5, 2.75);
            \draw[knot] (2.5, 1.5) -- (2.5, 2.75);
            \draw[knot] (3.5, 0) -- (3.5, -1.25);
            \draw[knot] (3.5, 1.5) -- (3.5, 2.75);
            \node at (2, 0.75) {\textit{B}};
            \draw[knot] (3.5, 2.75) to[out=90, in=90] (5, 2.75) to[out=-90, in=90] (5, -1.25) to[out=-90, in=-90] (3.5, -1.25);
            }
\]
is a split, so the shifting functor associated to it is the $\mathbb{Z}\times\mathbb{Z}$-grading shift $\{0,-1\}$, thus canceling with the original $\mathbb{Z}\times\mathbb{Z}$-grading shift of $\{0,1\}$. Said another way, $\mathrm{Tr}(f) \in \normalfont{\textsc{Hom}}_n\left(A \sqcup \bigcirc, \mathrm{Tr}(B)\right)$.

Next, let $g\in \normalfont{\textsc{Hom}}_{n-1}(A, \mathrm{Tr}(B)\{-1,0\})$ and denote the $\mathscr{I}$-degree of $g$ by $\mathcal{E}^w$. We define $\psi(g)\in \normalfont{\textsc{Hom}}_n\left(A \sqcup 1, \varphi_{\left(\tikz[baseline={([yshift=-.5ex]current bounding box.center)}, scale=.15]{
            \draw[fill=white, knot] (0,0) rectangle (4, 1.5);
            \draw[knot, red, rounded corners] (3.5, 2.125) -- (4.5, 2.125) -- (4.5, -0.625) -- (3.5, -0.625);
            \draw[knot] (0.5, 0) -- (0.5, -1.25);
            \draw[knot] (1.5, 0) -- (1.5, -1.25);
            \draw[knot] (2.5, 0) -- (2.5, -1.25);
            \draw[knot] (0.5, 1.5) -- (0.5, 2.75);
            \draw[knot] (1.5, 1.5) -- (1.5, 2.75);
            \draw[knot] (2.5, 1.5) -- (2.5, 2.75);
            \draw[knot] (3.5, 0) -- (3.5, -1.25);
            \draw[knot] (3.5, 1.5) -- (3.5, 2.75);
            \node at (2, 0.75) {\scriptsize \textit{B}};
},\, (0,1)\right)} B\right)$ to be the composition
\[
\tikz{
     \node(A) at (0,0) 
        {$\varphi_{\mathcal{E}^w}
        \tikz[baseline={([yshift=-.5ex]current bounding box.center)}, scale=.35]{
            \draw[knot] (0,0) rectangle (3, 1.5);
            \draw[knot] (0.5, 0) -- (0.5, -1.25);
            \draw[knot] (1.5, 0) -- (1.5, -1.25);
            \draw[knot] (2.5, 0) -- (2.5, -1.25);
            \draw[knot] (0.5, 1.5) -- (0.5, 2.75);
            \draw[knot] (1.5, 1.5) -- (1.5, 2.75);
            \draw[knot] (2.5, 1.5) -- (2.5, 2.75);
            \draw[knot] (3.5, -1.25) -- (3.5, 2.75);
            \node at (1.5, 0.75) {$A$};
            \draw[rounded corners, dashed, red] (-0.25, -0.625) rectangle (3.25, 2.125);

        }$};
    \node(B) at (5,0) 
        {$\tikz[baseline={([yshift=-.5ex]current bounding box.center)}, scale=.35]{
            \draw[rounded corners, knot] (3.5, -1) rectangle (4.5, 2.5);
            \draw[fill=white, knot] (0,0) rectangle (4, 1.5);
            \draw[knot] (0.5, 0) -- (0.5, -1.25);
            \draw[knot] (1.5, 0) -- (1.5, -1.25);
            \draw[knot] (2.5, 0) -- (2.5, -1.25);
            \draw[knot] (0.5, 1.5) -- (0.5, 2.75);
            \draw[knot] (1.5, 1.5) -- (1.5, 2.75);
            \draw[knot] (2.5, 1.5) -- (2.5, 2.75);
            \draw[knot] (5.5, -1.25) -- (5.5, 2.75);
            \draw[rounded corners, dashed, red] (4.25, 0.1) rectangle (5.75, 1.4);
            \node at (2, 0.75) {$B$};
        }\,\{-1,0\}$};
    \node(C) at (12.5,0) 
        {$\{-1,0\} \circ 
        \varphi_{\tikz[baseline={([yshift=-.5ex]current bounding box.center)}, scale=.25]{
            \draw[knot, red] (4.5, 0.75) -- (5.5, 0.75);
            \draw[rounded corners, knot] (3.5, -1) rectangle (4.5, 2.5);
            \draw[fill=white, knot] (0,0) rectangle (4, 1.5);
            \draw[knot] (0.5, 0) -- (0.5, -1.25);
            \draw[knot] (1.5, 0) -- (1.5, -1.25);
            \draw[knot] (2.5, 0) -- (2.5, -1.25);
            \draw[knot] (0.5, 1.5) -- (0.5, 2.75);
            \draw[knot] (1.5, 1.5) -- (1.5, 2.75);
            \draw[knot] (2.5, 1.5) -- (2.5, 2.75);
            \draw[knot] (5.5, -1.25) -- (5.5, 2.75);
            \node at (2, 0.75) {\scriptsize \textit{B}};
        }}^{-1}~
        \tikz[baseline={([yshift=-.5ex]current bounding box.center)}, scale=.35]{
            \draw[fill=white, knot] (0,0) rectangle (4, 1.5);
            \draw[knot] (0.5, 0) -- (0.5, -1.25);
            \draw[knot] (1.5, 0) -- (1.5, -1.25);
            \draw[knot] (2.5, 0) -- (2.5, -1.25);
            \draw[knot] (0.5, 1.5) -- (0.5, 2.75);
            \draw[knot] (1.5, 1.5) -- (1.5, 2.75);
            \draw[knot] (2.5, 1.5) -- (2.5, 2.75);
            \draw[knot] (3.5, 1.5) to[out=90, in=180] (4, 2.5) to[out=0, in=90] (4.5, 1.5) to[out=-90, in=-90] (5.5, 1.5) to[out=90, in=-90] (5.5, 2.75);
            \draw[knot] (3.5, 0) to[out=-90, in=180] (4, -1) to[out=0, in=-90] (4.5, 0) to[out=90, in=90] (5.5,0) to[out=-90, in=90] (5.5, -1.25);
            \node at (2, 0.75) {$B$};
        }$};
    \draw[->] (A) to[out=0,in=180] node[pos=.5,above,arrows=-]
        {$g \sqcup 1$} 
        (B);
    \draw[->] (B) to[out=0,in=180] node[pos=.5,above,arrows=-]
        {$\tikz[baseline={([yshift=-.5ex]current bounding box.center)}, scale=.25]{
            \draw[knot, red] (4.5, 0.75) -- (5.5, 0.75);
            \draw[rounded corners, knot] (3.5, -1) rectangle (4.5, 2.5);
            \draw[fill=white, knot] (0,0) rectangle (4, 1.5);
            \draw[knot] (0.5, 0) -- (0.5, -1.25);
            \draw[knot] (1.5, 0) -- (1.5, -1.25);
            \draw[knot] (2.5, 0) -- (2.5, -1.25);
            \draw[knot] (0.5, 1.5) -- (0.5, 2.75);
            \draw[knot] (1.5, 1.5) -- (1.5, 2.75);
            \draw[knot] (2.5, 1.5) -- (2.5, 2.75);
            \draw[knot] (5.5, -1.25) -- (5.5, 2.75);
        } \circ \varphi_{H_S}$} 
        (C);    
}
\]
where
\[
\varphi_{H_{S}}: 
\mathrm{Id}
\Rightarrow
\varphi_{
\tikz[baseline={([yshift=-.5ex]current bounding box.center)}, scale=.25]{
            \draw[knot, red] (4.5, 0.75) -- (5.5, 0.75);
            \draw[rounded corners, knot] (3.5, -1) rectangle (4.5, 2.5);
            \draw[fill=white, knot] (0,0) rectangle (4, 1.5);
            \draw[knot] (0.5, 0) -- (0.5, -1.25);
            \draw[knot] (1.5, 0) -- (1.5, -1.25);
            \draw[knot] (2.5, 0) -- (2.5, -1.25);
            \draw[knot] (0.5, 1.5) -- (0.5, 2.75);
            \draw[knot] (1.5, 1.5) -- (1.5, 2.75);
            \draw[knot] (2.5, 1.5) -- (2.5, 2.75);
            \draw[knot] (5.5, -1.25) -- (5.5, 2.75);
            \node at (2, 0.75) {$B$};
        }}^{-1}
\circ
\varphi_{
\tikz[baseline={([yshift=-.5ex]current bounding box.center)}, scale=.25]{
            \draw[knot, red] (4.5, 0.75) -- (5.5, 0.75);
            \draw[rounded corners, knot] (3.5, -1) rectangle (4.5, 2.5);
            \draw[fill=white, knot] (0,0) rectangle (4, 1.5);
            \draw[knot] (0.5, 0) -- (0.5, -1.25);
            \draw[knot] (1.5, 0) -- (1.5, -1.25);
            \draw[knot] (2.5, 0) -- (2.5, -1.25);
            \draw[knot] (0.5, 1.5) -- (0.5, 2.75);
            \draw[knot] (1.5, 1.5) -- (1.5, 2.75);
            \draw[knot] (2.5, 1.5) -- (2.5, 2.75);
            \draw[knot] (5.5, -1.25) -- (5.5, 2.75);
            \node at (2, 0.75) {$B$};
        }}~.
\]
Then, $\psi(g)$ has the desired form since
\[
\varphi_{
\tikz[baseline={([yshift=-.5ex]current bounding box.center)}, scale=.25]{
            \draw[knot, red] (4.5, 0.75) -- (5.5, 0.75);
            \draw[rounded corners, knot] (3.5, -1) rectangle (4.5, 2.5);
            \draw[fill=white, knot] (0,0) rectangle (4, 1.5);
            \draw[knot] (0.5, 0) -- (0.5, -1.25);
            \draw[knot] (1.5, 0) -- (1.5, -1.25);
            \draw[knot] (2.5, 0) -- (2.5, -1.25);
            \draw[knot] (0.5, 1.5) -- (0.5, 2.75);
            \draw[knot] (1.5, 1.5) -- (1.5, 2.75);
            \draw[knot] (2.5, 1.5) -- (2.5, 2.75);
            \draw[knot] (5.5, -1.25) -- (5.5, 2.75);
            \node at (2, 0.75) {$B$};
        }}^{-1} = 
\varphi_{\left(
\tikz[baseline={([yshift=-.5ex]current bounding box.center)}, scale=.25]{
            \draw[fill=white, knot] (0,0) rectangle (4, 1.5);
            \draw[knot, red] (5, 1.25) -- (5, 0.25);
            \draw[knot] (0.5, 0) -- (0.5, -1.25);
            \draw[knot] (1.5, 0) -- (1.5, -1.25);
            \draw[knot] (2.5, 0) -- (2.5, -1.25);
            \draw[knot] (0.5, 1.5) -- (0.5, 2.75);
            \draw[knot] (1.5, 1.5) -- (1.5, 2.75);
            \draw[knot] (2.5, 1.5) -- (2.5, 2.75);
            \draw[knot] (3.5, 1.5) to[out=90, in=180] (4, 2.5) to[out=0, in=90] (4.5, 1.5) to[out=-90, in=-90] (5.5, 1.5) to[out=90, in=-90] (5.5, 2.75);
            \draw[knot] (3.5, 0) to[out=-90, in=180] (4, -1) to[out=0, in=-90] (4.5, 0) to[out=90, in=90] (5.5,0) to[out=-90, in=90] (5.5, -1.25);
            \node at (2, 0.75) {$B$};
        },\,(1,1)\right)}
\]
composed with $\{-1,0\}$ is $\varphi_{\left(\tikz[baseline={([yshift=-.5ex]current bounding box.center)}, scale=.15]{
            \draw[fill=white, knot] (0,0) rectangle (4, 1.5);
            \draw[knot, red, rounded corners] (3.5, 2.125) -- (4.5, 2.125) -- (4.5, -0.625) -- (3.5, -0.625);
            \draw[knot] (0.5, 0) -- (0.5, -1.25);
            \draw[knot] (1.5, 0) -- (1.5, -1.25);
            \draw[knot] (2.5, 0) -- (2.5, -1.25);
            \draw[knot] (0.5, 1.5) -- (0.5, 2.75);
            \draw[knot] (1.5, 1.5) -- (1.5, 2.75);
            \draw[knot] (2.5, 1.5) -- (2.5, 2.75);
            \draw[knot] (3.5, 0) -- (3.5, -1.25);
            \draw[knot] (3.5, 1.5) -- (3.5, 2.75);
            \node at (2, 0.75) {\scriptsize \textit{B}};
        },\, (0,1)\right)}$.

Now, we compute $\psi(\phi(f))$ as the composition
\[
\tikz{
     \node(A) at (0,0) 
        {$\varphi_{\Delta^v}
        \tikz[baseline={([yshift=-.5ex]current bounding box.center)}, scale=.3]{
            \draw[knot] (0,0) rectangle (3, 1.5);
            \draw[knot] (0.5, 0) -- (0.5, -1.25);
            \draw[knot] (1.5, 0) -- (1.5, -1.25);
            \draw[knot] (2.5, 0) -- (2.5, -1.25);
            \draw[knot] (0.5, 1.5) -- (0.5, 2.75);
            \draw[knot] (1.5, 1.5) -- (1.5, 2.75);
            \draw[knot] (2.5, 1.5) -- (2.5, 2.75);
            \draw[knot] (5.5, -1.25) -- (5.5, 2.75);
            \node at (1.5, 0.75) {$A$};
            \draw[rounded corners, dashed, red] (3.5, -1) rectangle (4.5, 2.5);
        }$};
    \node(B) at (7,0) 
        {$\varphi_{\Delta^{v + (-1,0)}}
        \tikz[baseline={([yshift=-.5ex]current bounding box.center)}, scale=.3]{
            \draw[knot] (0,0) rectangle (3, 1.5);
            \draw[knot] (0.5, 0) -- (0.5, -1.25);
            \draw[knot] (1.5, 0) -- (1.5, -1.25);
            \draw[knot] (2.5, 0) -- (2.5, -1.25);
            \draw[knot] (0.5, 1.5) -- (0.5, 2.75);
            \draw[knot] (1.5, 1.5) -- (1.5, 2.75);
            \draw[knot] (2.5, 1.5) -- (2.5, 2.75);
            \draw[knot] (5.5, -1.25) -- (5.5, 2.75);
            \node at (1.5, 0.75) {$A$};
            \draw[rounded corners, knot] (3.5, -1) rectangle (4.5, 2.5);
            \draw[rounded corners, dashed, red] (-.5, -.5) rectangle (4, 2);
        }$};
    \node(C) at (12.5,0) 
        {$\tikz[baseline={([yshift=-.5ex]current bounding box.center)}, scale=.3]{
            \draw[rounded corners, knot] (3.5, -1) rectangle (4.5, 2.5);
            \draw[fill=white, knot] (0,0) rectangle (4, 1.5);
            \draw[knot] (0.5, 0) -- (0.5, -1.25);
            \draw[knot] (1.5, 0) -- (1.5, -1.25);
            \draw[knot] (2.5, 0) -- (2.5, -1.25);
            \draw[knot] (0.5, 1.5) -- (0.5, 2.75);
            \draw[knot] (1.5, 1.5) -- (1.5, 2.75);
            \draw[knot] (2.5, 1.5) -- (2.5, 2.75);
            \draw[knot] (5.5, -1.25) -- (5.5, 2.75);
            \draw[rounded corners, dashed, red] (4.25, 0.1) rectangle (5.75, 1.4);
            \node at (2, 0.75) {$B$};
        }\,\{-1,0\}$};
    \node (D) at (12.5,-3)
        {$\varphi_{\left(\tikz[baseline={([yshift=-.5ex]current bounding box.center)}, scale=.15]{
            \draw[fill=white, knot] (0,0) rectangle (4, 1.5);
            \draw[knot, red, rounded corners] (3.5, 2.125) -- (4.5, 2.125) -- (4.5, -0.625) -- (3.5, -0.625);
            \draw[knot] (0.5, 0) -- (0.5, -1.25);
            \draw[knot] (1.5, 0) -- (1.5, -1.25);
            \draw[knot] (2.5, 0) -- (2.5, -1.25);
            \draw[knot] (0.5, 1.5) -- (0.5, 2.75);
            \draw[knot] (1.5, 1.5) -- (1.5, 2.75);
            \draw[knot] (2.5, 1.5) -- (2.5, 2.75);
            \draw[knot] (3.5, 0) -- (3.5, -1.25);
            \draw[knot] (3.5, 1.5) -- (3.5, 2.75);
            \node at (2, 0.75) {\scriptsize \textit{B}};
        },\, (0,1)\right)}
        \tikz[baseline={([yshift=-.5ex]current bounding box.center)}, scale=.3]{
            \draw[fill=white, knot] (0,0) rectangle (4, 1.5);
            \draw[knot] (0.5, 0) -- (0.5, -1.25);
            \draw[knot] (1.5, 0) -- (1.5, -1.25);
            \draw[knot] (2.5, 0) -- (2.5, -1.25);
            \draw[knot] (3.5, 0) -- (3.5, -1.25);
            \draw[knot] (0.5, 1.5) -- (0.5, 2.75);
            \draw[knot] (1.5, 1.5) -- (1.5, 2.75);
            \draw[knot] (2.5, 1.5) -- (2.5, 2.75);
            \draw[knot] (3.5, 1.5) -- (3.5, 2.75);
            \node at (2, 0.75) {$B$};
        }
        $};
    \draw[->] (A) to[out=0,in=180] node[pos=.5,above,arrows=-]
        {$\lambda_{\phi(f)} \circ \mathcal{F}\left(\tikz[baseline={([yshift=-.5ex]current bounding box.center)}, scale=.4]{
	\draw (1,2) .. controls (1,1) and (2,1) .. (2,2);
	\draw (1,2) .. controls (1,1.75) and (2,1.75) .. (2,2);
	\draw (1,2) .. controls (1,2.25) and (2,2.25) .. (2,2);
}\right) \circ \varphi_{H_B}$} 
        (B);
    \draw[->] (B) to[out=0,in=180] node[pos=.5,above,arrows=-]
        {$\mathrm{Tr}(f) \sqcup 1$} 
        (C);
    \draw[->] (C) to[out=-90,in=90] node[pos=.5,right,arrows=-]
        {$\tikz[baseline={([yshift=-.5ex]current bounding box.center)}, scale=.22]{
            \draw[knot, red] (4.5, 0.75) -- (5.5, 0.75);
            \draw[rounded corners, knot] (3.5, -1) rectangle (4.5, 2.5);
            \draw[fill=white, knot] (0,0) rectangle (4, 1.5);
            \draw[knot] (0.5, 0) -- (0.5, -1.25);
            \draw[knot] (1.5, 0) -- (1.5, -1.25);
            \draw[knot] (2.5, 0) -- (2.5, -1.25);
            \draw[knot] (0.5, 1.5) -- (0.5, 2.75);
            \draw[knot] (1.5, 1.5) -- (1.5, 2.75);
            \draw[knot] (2.5, 1.5) -- (2.5, 2.75);
            \draw[knot] (5.5, -1.25) -- (5.5, 2.75);
        } \circ \varphi_{H_S}$} (D)
}
\]
If we slide $f$ past the saddle, then the above complex is equivalent to the following one, where we have compensated for the slide by a change of chronology $\varphi_{H_1}$.
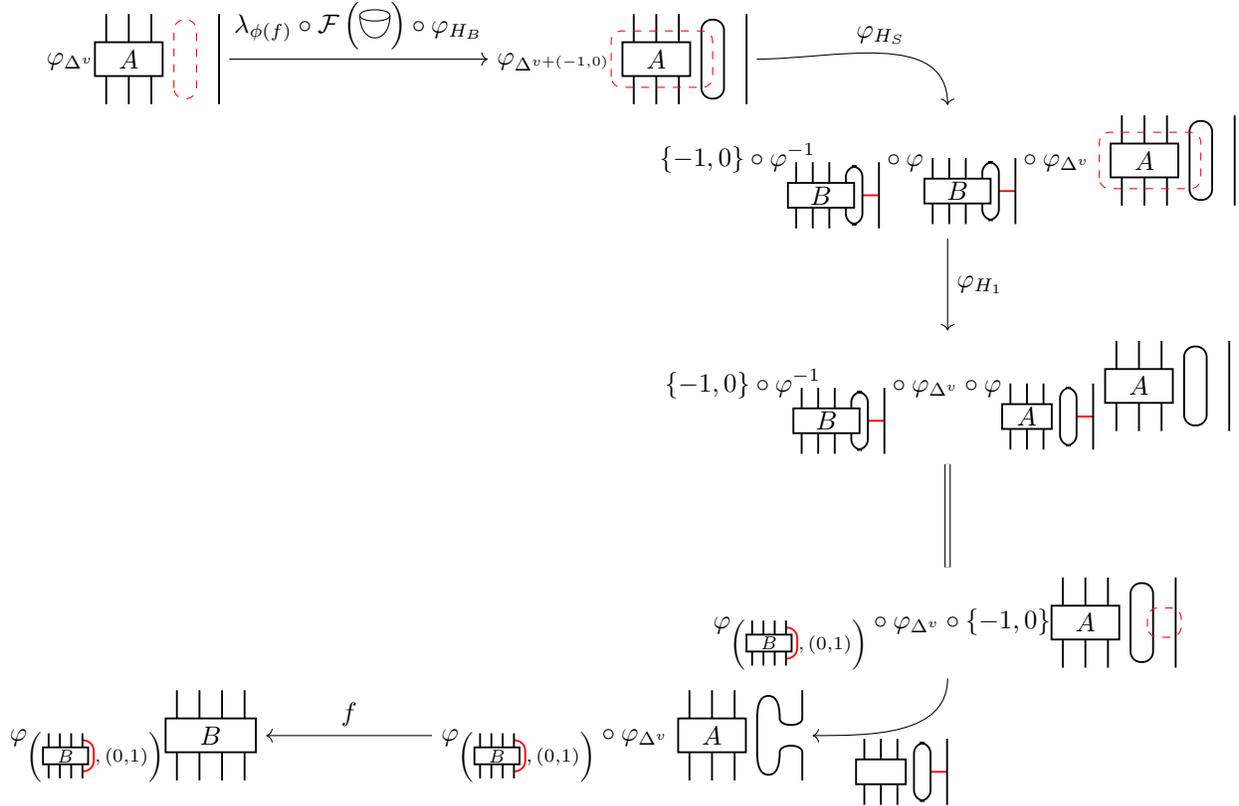
\[
\tikz{
     \node(A) at (0,0) 
        {$\varphi_{\Delta^v}
        \tikz[baseline={([yshift=-.5ex]current bounding box.center)}, scale=.3]{
            \draw[knot] (0,0) rectangle (3, 1.5);
            \draw[knot] (0.5, 0) -- (0.5, -1.25);
            \draw[knot] (1.5, 0) -- (1.5, -1.25);
            \draw[knot] (2.5, 0) -- (2.5, -1.25);
            \draw[knot] (0.5, 1.5) -- (0.5, 2.75);
            \draw[knot] (1.5, 1.5) -- (1.5, 2.75);
            \draw[knot] (2.5, 1.5) -- (2.5, 2.75);
            \draw[knot] (5.5, -1.25) -- (5.5, 2.75);
            \node at (1.5, 0.75) {$A$};
            \draw[rounded corners, dashed, red] (3.5, -1) rectangle (4.5, 2.5);
        }$};
    \node(B) at (6.5,0) 
        {$\varphi_{\Delta^{v + (-1,0)}}
        \tikz[baseline={([yshift=-.5ex]current bounding box.center)}, scale=.3]{
            \draw[knot] (0,0) rectangle (3, 1.5);
            \draw[knot] (0.5, 0) -- (0.5, -1.25);
            \draw[knot] (1.5, 0) -- (1.5, -1.25);
            \draw[knot] (2.5, 0) -- (2.5, -1.25);
            \draw[knot] (0.5, 1.5) -- (0.5, 2.75);
            \draw[knot] (1.5, 1.5) -- (1.5, 2.75);
            \draw[knot] (2.5, 1.5) -- (2.5, 2.75);
            \draw[knot] (5.5, -1.25) -- (5.5, 2.75);
            \node at (1.5, 0.75) {$A$};
            \draw[rounded corners, knot] (3.5, -1) rectangle (4.5, 2.5);
            \draw[rounded corners, dashed, red] (-.5, -.5) rectangle (4, 2);
        }$};
    \node(C) at (10.825,-1.5) 
        {$\{-1,0\} \circ \varphi_{
\tikz[baseline={([yshift=-.5ex]current bounding box.center)}, scale=.22]{
            \draw[knot, red] (4.5, 0.75) -- (5.5, 0.75);
            \draw[rounded corners, knot] (3.5, -1) rectangle (4.5, 2.5);
            \draw[fill=white, knot] (0,0) rectangle (4, 1.5);
            \draw[knot] (0.5, 0) -- (0.5, -1.25);
            \draw[knot] (1.5, 0) -- (1.5, -1.25);
            \draw[knot] (2.5, 0) -- (2.5, -1.25);
            \draw[knot] (0.5, 1.5) -- (0.5, 2.75);
            \draw[knot] (1.5, 1.5) -- (1.5, 2.75);
            \draw[knot] (2.5, 1.5) -- (2.5, 2.75);
            \draw[knot] (5.5, -1.25) -- (5.5, 2.75);
            \node at (2, 0.75) {$B$};
        }}^{-1}
\circ
\varphi_{
\tikz[baseline={([yshift=-.5ex]current bounding box.center)}, scale=.22]{
            \draw[knot, red] (4.5, 0.75) -- (5.5, 0.75);
            \draw[rounded corners, knot] (3.5, -1) rectangle (4.5, 2.5);
            \draw[fill=white, knot] (0,0) rectangle (4, 1.5);
            \draw[knot] (0.5, 0) -- (0.5, -1.25);
            \draw[knot] (1.5, 0) -- (1.5, -1.25);
            \draw[knot] (2.5, 0) -- (2.5, -1.25);
            \draw[knot] (0.5, 1.5) -- (0.5, 2.75);
            \draw[knot] (1.5, 1.5) -- (1.5, 2.75);
            \draw[knot] (2.5, 1.5) -- (2.5, 2.75);
            \draw[knot] (5.5, -1.25) -- (5.5, 2.75);
            \node at (2, 0.75) {$B$};
        }} \circ \varphi_{\Delta^{v}} ~
        \tikz[baseline={([yshift=-.5ex]current bounding box.center)}, scale=.3]{
            \draw[knot] (0,0) rectangle (3, 1.5);
            \draw[knot] (0.5, 0) -- (0.5, -1.25);
            \draw[knot] (1.5, 0) -- (1.5, -1.25);
            \draw[knot] (2.5, 0) -- (2.5, -1.25);
            \draw[knot] (0.5, 1.5) -- (0.5, 2.75);
            \draw[knot] (1.5, 1.5) -- (1.5, 2.75);
            \draw[knot] (2.5, 1.5) -- (2.5, 2.75);
            \draw[knot] (5.5, -1.25) -- (5.5, 2.75);
            \node at (1.5, 0.75) {$A$};
            \draw[rounded corners, knot] (3.5, -1) rectangle (4.5, 2.5);
            \draw[rounded corners, dashed, red] (-.5, -.5) rectangle (4, 2);
        }$};
    \node (D) at (10.825,-4.5)
        {$\{-1,0\} \circ \varphi_{
\tikz[baseline={([yshift=-.5ex]current bounding box.center)}, scale=.22]{
            \draw[knot, red] (4.5, 0.75) -- (5.5, 0.75);
            \draw[rounded corners, knot] (3.5, -1) rectangle (4.5, 2.5);
            \draw[fill=white, knot] (0,0) rectangle (4, 1.5);
            \draw[knot] (0.5, 0) -- (0.5, -1.25);
            \draw[knot] (1.5, 0) -- (1.5, -1.25);
            \draw[knot] (2.5, 0) -- (2.5, -1.25);
            \draw[knot] (0.5, 1.5) -- (0.5, 2.75);
            \draw[knot] (1.5, 1.5) -- (1.5, 2.75);
            \draw[knot] (2.5, 1.5) -- (2.5, 2.75);
            \draw[knot] (5.5, -1.25) -- (5.5, 2.75);
            \node at (2, 0.75) {$B$};
        }}^{-1}
        \circ
        \varphi_{\Delta^{v}}
        \circ
        \varphi_{
        \tikz[baseline={([yshift=-.5ex]current bounding box.center)}, scale=.22]{
            \draw[knot, red] (4.5, 0.75) -- (5.5, 0.75);
            \draw[knot] (0,0) rectangle (3, 1.5);
            \draw[knot] (0.5, 0) -- (0.5, -1.25);
            \draw[knot] (1.5, 0) -- (1.5, -1.25);
            \draw[knot] (2.5, 0) -- (2.5, -1.25);
            \draw[knot] (0.5, 1.5) -- (0.5, 2.75);
            \draw[knot] (1.5, 1.5) -- (1.5, 2.75);
            \draw[knot] (2.5, 1.5) -- (2.5, 2.75);
            \draw[knot] (5.5, -1.25) -- (5.5, 2.75);
            \node at (1.5, 0.75) {$A$};
            \draw[rounded corners, knot] (3.5, -1) rectangle (4.5, 2.5);
        }}~
        \tikz[baseline={([yshift=-.5ex]current bounding box.center)}, scale=.3]{
            \draw[knot] (0,0) rectangle (3, 1.5);
            \draw[knot] (0.5, 0) -- (0.5, -1.25);
            \draw[knot] (1.5, 0) -- (1.5, -1.25);
            \draw[knot] (2.5, 0) -- (2.5, -1.25);
            \draw[knot] (0.5, 1.5) -- (0.5, 2.75);
            \draw[knot] (1.5, 1.5) -- (1.5, 2.75);
            \draw[knot] (2.5, 1.5) -- (2.5, 2.75);
            \draw[knot] (5.5, -1.25) -- (5.5, 2.75);
            \node at (1.5, 0.75) {$A$};
            \draw[rounded corners, knot] (3.5, -1) rectangle (4.5, 2.5);
        }$};
    \node (E) at (10.825, -7.5) {$\varphi_{\left(\tikz[baseline={([yshift=-.5ex]current bounding box.center)}, scale=.15]{
            \draw[fill=white, knot] (0,0) rectangle (4, 1.5);
            \draw[knot, red, rounded corners] (3.5, 2.125) -- (4.5, 2.125) -- (4.5, -0.625) -- (3.5, -0.625);
            \draw[knot] (0.5, 0) -- (0.5, -1.25);
            \draw[knot] (1.5, 0) -- (1.5, -1.25);
            \draw[knot] (2.5, 0) -- (2.5, -1.25);
            \draw[knot] (0.5, 1.5) -- (0.5, 2.75);
            \draw[knot] (1.5, 1.5) -- (1.5, 2.75);
            \draw[knot] (2.5, 1.5) -- (2.5, 2.75);
            \draw[knot] (3.5, 0) -- (3.5, -1.25);
            \draw[knot] (3.5, 1.5) -- (3.5, 2.75);
            \node at (2, 0.75) {\scriptsize \textit{B}};
        },\, (0,1)\right)}
        \circ
        \varphi_{\Delta^{v}}
        \circ
        \{-1,0\}        %
        \tikz[baseline={([yshift=-.5ex]current bounding box.center)}, scale=.3]{
            \draw[knot] (0,0) rectangle (3, 1.5);
            \draw[knot] (0.5, 0) -- (0.5, -1.25);
            \draw[knot] (1.5, 0) -- (1.5, -1.25);
            \draw[knot] (2.5, 0) -- (2.5, -1.25);
            \draw[knot] (0.5, 1.5) -- (0.5, 2.75);
            \draw[knot] (1.5, 1.5) -- (1.5, 2.75);
            \draw[knot] (2.5, 1.5) -- (2.5, 2.75);
            \draw[knot] (5.5, -1.25) -- (5.5, 2.75);
            \node at (1.5, 0.75) {$A$};
            \draw[rounded corners, knot] (3.5, -1) rectangle (4.5, 2.5);
            \draw[rounded corners, dashed, red] (4.25, 0.1) rectangle (5.75, 1.4);
        }$};
    \node (F) at (6.5, -9) {$\varphi_{\left(\tikz[baseline={([yshift=-.5ex]current bounding box.center)}, scale=.15]{
            \draw[fill=white, knot] (0,0) rectangle (4, 1.5);
            \draw[knot, red, rounded corners] (3.5, 2.125) -- (4.5, 2.125) -- (4.5, -0.625) -- (3.5, -0.625);
            \draw[knot] (0.5, 0) -- (0.5, -1.25);
            \draw[knot] (1.5, 0) -- (1.5, -1.25);
            \draw[knot] (2.5, 0) -- (2.5, -1.25);
            \draw[knot] (0.5, 1.5) -- (0.5, 2.75);
            \draw[knot] (1.5, 1.5) -- (1.5, 2.75);
            \draw[knot] (2.5, 1.5) -- (2.5, 2.75);
            \draw[knot] (3.5, 0) -- (3.5, -1.25);
            \draw[knot] (3.5, 1.5) -- (3.5, 2.75);
            \node at (2, 0.75) {\scriptsize \textit{B}};
        },\, (0,1)\right)}
        \circ
        \varphi_{\Delta^{v}}~
        \tikz[baseline={([yshift=-.5ex]current bounding box.center)}, scale=.3]{
            \draw[knot] (0,0) rectangle (3, 1.5);
            \draw[knot] (0.5, 0) -- (0.5, -1.25);
            \draw[knot] (1.5, 0) -- (1.5, -1.25);
            \draw[knot] (2.5, 0) -- (2.5, -1.25);
            \draw[knot] (0.5, 1.5) -- (0.5, 2.75);
            \draw[knot] (1.5, 1.5) -- (1.5, 2.75);
            \draw[knot] (2.5, 1.5) -- (2.5, 2.75);
            \node at (1.5, 0.75) {$A$};
            \draw[knot] (5.5, -1.25) to[out=90, in =-90] (5.5,0) to[out=90, in=90] (4.5, 0) to[out=-90, in=0] (4, -1) to[out=180, in=-90] (3.5, 0) to[out=90, in=-90] (3.5, 1.5) to[out=90, in=180] (4, 2.5) to[out=0, in=90] (4.5, 1.5) to[out=-90, in=-90] (5.5, 1.5) to[out=90, in=-90] (5.5, 2.75)
        }$};
    \node (G) at (0, -9) {$\varphi_{\left(\tikz[baseline={([yshift=-.5ex]current bounding box.center)}, scale=.15]{
            \draw[fill=white, knot] (0,0) rectangle (4, 1.5);
            \draw[knot, red, rounded corners] (3.5, 2.125) -- (4.5, 2.125) -- (4.5, -0.625) -- (3.5, -0.625);
            \draw[knot] (0.5, 0) -- (0.5, -1.25);
            \draw[knot] (1.5, 0) -- (1.5, -1.25);
            \draw[knot] (2.5, 0) -- (2.5, -1.25);
            \draw[knot] (0.5, 1.5) -- (0.5, 2.75);
            \draw[knot] (1.5, 1.5) -- (1.5, 2.75);
            \draw[knot] (2.5, 1.5) -- (2.5, 2.75);
            \draw[knot] (3.5, 0) -- (3.5, -1.25);
            \draw[knot] (3.5, 1.5) -- (3.5, 2.75);
            \node at (2, 0.75) {\scriptsize \textit{B}};
        },\, (0,1)\right)}
        \tikz[baseline={([yshift=-.5ex]current bounding box.center)}, scale=.3]{
            \draw[fill=white, knot] (0,0) rectangle (4, 1.5);
            \draw[knot] (0.5, 0) -- (0.5, -1.25);
            \draw[knot] (1.5, 0) -- (1.5, -1.25);
            \draw[knot] (2.5, 0) -- (2.5, -1.25);
            \draw[knot] (3.5, 0) -- (3.5, -1.25);
            \draw[knot] (0.5, 1.5) -- (0.5, 2.75);
            \draw[knot] (1.5, 1.5) -- (1.5, 2.75);
            \draw[knot] (2.5, 1.5) -- (2.5, 2.75);
            \draw[knot] (3.5, 1.5) -- (3.5, 2.75);
            \node at (2, 0.75) {$B$};
        }
        $};
    \draw[->] (A) to[out=0,in=180] node[pos=.5,above,arrows=-]
        {$\lambda_{\phi(f)} \circ \mathcal{F}\left(\tikz[baseline={([yshift=-.5ex]current bounding box.center)}, scale=.4]{
	\draw (1,2) .. controls (1,1) and (2,1) .. (2,2);
	\draw (1,2) .. controls (1,1.75) and (2,1.75) .. (2,2);
	\draw (1,2) .. controls (1,2.25) and (2,2.25) .. (2,2);
}\right) \circ \varphi_{H_B}$} 
        (B);
    \draw[->] (B) to[out=0,in=90] node[pos=.5,above,arrows=-]
        {$\varphi_{H_S}$} 
        (C);
    \draw[->] (C) to[out=-90,in=90] node[pos=.5,right,arrows=-]
        {$\varphi_{H_1}$} (D);
    \draw[double equal sign distance] (D) to[out=-90,in=90] (E);
    \draw[->] (E) to[out=-90, in=0] node[pos=.5,below,arrows=-]
        {$\tikz[baseline={([yshift=-.5ex]current bounding box.center)}, scale=.22]{
            \draw[knot, red] (4.5, 0.75) -- (5.5, 0.75);
            \draw[knot] (0,0) rectangle (3, 1.5);
            \draw[knot] (0.5, 0) -- (0.5, -1.25);
            \draw[knot] (1.5, 0) -- (1.5, -1.25);
            \draw[knot] (2.5, 0) -- (2.5, -1.25);
            \draw[knot] (0.5, 1.5) -- (0.5, 2.75);
            \draw[knot] (1.5, 1.5) -- (1.5, 2.75);
            \draw[knot] (2.5, 1.5) -- (2.5, 2.75);
            \draw[knot] (5.5, -1.25) -- (5.5, 2.75);
            \draw[rounded corners, knot] (3.5, -1) rectangle (4.5, 2.5);
        }$} (F);
    \draw[->] (F) to[out=180, in=0] node[pos=.5,above,arrows=-] {$f$} (G);
}
\]
The key observation is that $\lambda_{\phi(f)}$---which corresponds to sliding a shift by $\{-1,0\}$ through $\Delta^v$---and $\varphi_{H_1}$---which corresponds to a change of chronology which pushes a saddle through $\Delta^v$, at which point it is realized as a merge (and the grading shift associated to merges is $\{-1,0\}$)---are inverse to one another. After this, the birth and merge cancel with one another, and we conclude that $\psi(\phi(f)) = f$.

We play a similar game for $\phi(\psi(g))$: it is computed as
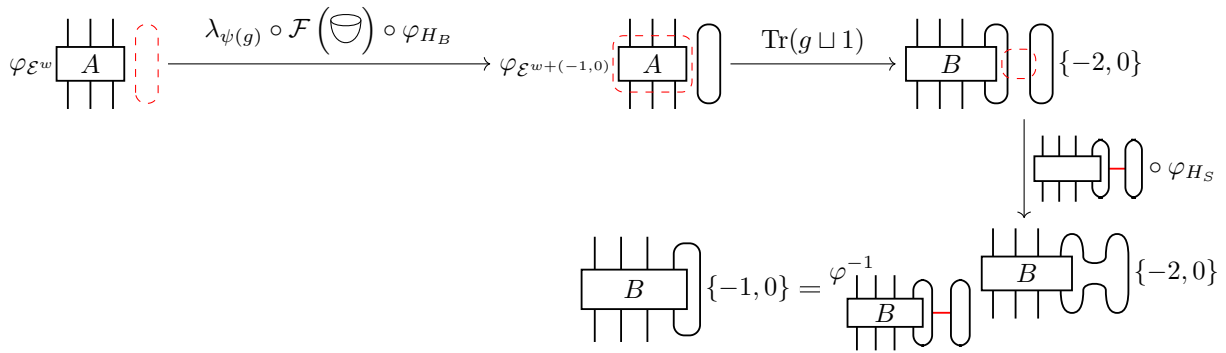
\[
\tikz{
     \node(A) at (0,0) 
        {$\varphi_{\mathcal{E}^w}
        \tikz[baseline={([yshift=-.5ex]current bounding box.center)}, scale=.3]{
            \draw[knot] (0,0) rectangle (3, 1.5);
            \draw[knot] (0.5, 0) -- (0.5, -1.25);
            \draw[knot] (1.5, 0) -- (1.5, -1.25);
            \draw[knot] (2.5, 0) -- (2.5, -1.25);
            \draw[knot] (0.5, 1.5) -- (0.5, 2.75);
            \draw[knot] (1.5, 1.5) -- (1.5, 2.75);
            \draw[knot] (2.5, 1.5) -- (2.5, 2.75);
            \node at (1.5, 0.75) {$A$};
            \draw[rounded corners, dashed, red] (3.5, -1) rectangle (4.5, 2.5);
        }$};
    \node(B) at (7,0) 
        {$\varphi_{\mathcal{E}^{w + (-1,0)}}
        \tikz[baseline={([yshift=-.5ex]current bounding box.center)}, scale=.3]{
            \draw[knot] (0,0) rectangle (3, 1.5);
            \draw[knot] (0.5, 0) -- (0.5, -1.25);
            \draw[knot] (1.5, 0) -- (1.5, -1.25);
            \draw[knot] (2.5, 0) -- (2.5, -1.25);
            \draw[knot] (0.5, 1.5) -- (0.5, 2.75);
            \draw[knot] (1.5, 1.5) -- (1.5, 2.75);
            \draw[knot] (2.5, 1.5) -- (2.5, 2.75);
            \node at (1.5, 0.75) {$A$};
            \draw[rounded corners, knot] (3.5, -1) rectangle (4.5, 2.5);
            \draw[rounded corners, dashed, red] (-.25, -.5) rectangle (3.25, 2);
        }$};
    \node(C) at (12.5,0) 
        {$\tikz[baseline={([yshift=-.5ex]current bounding box.center)}, scale=.3]{
            \draw[rounded corners, knot] (3.5, -1) rectangle (4.5, 2.5);
            \draw[fill=white, knot] (0,0) rectangle (4, 1.5);
            \draw[knot] (0.5, 0) -- (0.5, -1.25);
            \draw[knot] (1.5, 0) -- (1.5, -1.25);
            \draw[knot] (2.5, 0) -- (2.5, -1.25);
            \draw[knot] (0.5, 1.5) -- (0.5, 2.75);
            \draw[knot] (1.5, 1.5) -- (1.5, 2.75);
            \draw[knot] (2.5, 1.5) -- (2.5, 2.75);
            \draw[rounded corners, knot] (5.5, -1) rectangle (6.5, 2.5);
            \draw[rounded corners, dashed, red] (4.25, 0.1) rectangle (5.75, 1.4);
            \node at (2, 0.75) {$B$};
        }\,\{-2,0\}$};
    \node (D) at (12.5,-3)
        {$\varphi_{
\tikz[baseline={([yshift=-.5ex]current bounding box.center)}, scale=.25]{
            \draw[knot, red] (4.5, 0.75) -- (5.5, 0.75);
            \draw[rounded corners, knot] (3.5, -1) rectangle (4.5, 2.5);
            \draw[fill=white, knot] (0,0) rectangle (4, 1.5);
            \draw[knot] (0.5, 0) -- (0.5, -1.25);
            \draw[knot] (1.5, 0) -- (1.5, -1.25);
            \draw[knot] (2.5, 0) -- (2.5, -1.25);
            \draw[knot] (0.5, 1.5) -- (0.5, 2.75);
            \draw[knot] (1.5, 1.5) -- (1.5, 2.75);
            \draw[knot] (2.5, 1.5) -- (2.5, 2.75);
            \draw[rounded corners, knot] (5.5, -1) rectangle (6.5, 2.5);
            \node at (2, 0.75) {$B$};
        }
        }^{-1}~
        \tikz[baseline={([yshift=-.5ex]current bounding box.center)}, scale=.3]{
            \draw[fill=white, knot] (0,0) rectangle (4, 1.5);
            \draw[knot] (0.5, 0) -- (0.5, -1.25);
            \draw[knot] (1.5, 0) -- (1.5, -1.25);
            \draw[knot] (2.5, 0) -- (2.5, -1.25);
            \draw[knot] (0.5, 1.5) -- (0.5, 2.75);
            \draw[knot] (1.5, 1.5) -- (1.5, 2.75);
            \draw[knot] (2.5, 1.5) -- (2.5, 2.75);
            \draw[knot] (3.5, 1.5) to[out=90, in=180] (4, 2.5) to[out=0, in=90] (4.5, 1.5) to[out=-90, in=-90] (5.5, 1.5) to[out=90, in=180] (6, 2.5) to[out=0, in=90] (6.5, 1.5) to[out=-90, in=90] (6.5, 0) to[out=-90, in=0] (6, -1) to[out=180, in=-90] (5.5, 0) to[out=90, in=90] (4.5, 0) to[out=-90, in=0] (4, -1) to[out=180, in=-90] (3.5, 0);
            \node at (2, 0.75) {$B$};
        }\, \{-2, 0\}$};
    \node(E) at (8,-3) {$\tikz[baseline={([yshift=-.5ex]current bounding box.center)}, scale=.35]{
            \draw[rounded corners, knot] (3.5, -1) rectangle (4.5, 2.5);
            \draw[fill=white, knot] (0,0) rectangle (4, 1.5);
            \draw[knot] (0.5, 0) -- (0.5, -1.25);
            \draw[knot] (1.5, 0) -- (1.5, -1.25);
            \draw[knot] (2.5, 0) -- (2.5, -1.25);
            \draw[knot] (0.5, 1.5) -- (0.5, 2.75);
            \draw[knot] (1.5, 1.5) -- (1.5, 2.75);
            \draw[knot] (2.5, 1.5) -- (2.5, 2.75);
            \node at (2, 0.75) {$B$};
        }\,\{-1,0\}$};
    \draw[->] (A) to[out=0,in=180] node[pos=.5,above,arrows=-]
        {$\lambda_{\psi(g)} \circ \mathcal{F}\left(\tikz[baseline={([yshift=-.5ex]current bounding box.center)}, scale=.4]{
	\draw (1,2) .. controls (1,1) and (2,1) .. (2,2);
	\draw (1,2) .. controls (1,1.75) and (2,1.75) .. (2,2);
	\draw (1,2) .. controls (1,2.25) and (2,2.25) .. (2,2);
}\right) \circ \varphi_{H_B}$} 
        (B);
    \draw[->] (B) to[out=0,in=180] node[pos=.5,above,arrows=-]
        {$\mathrm{Tr}(g \sqcup 1)$} 
        (C);
    \draw[->] (C) to[out=-90,in=90] node[pos=.5,right,arrows=-]
        {$\tikz[baseline={([yshift=-.5ex]current bounding box.center)}, scale=.22]{
            \draw[knot, red] (4.5, 0.75) -- (5.5, 0.75);
            \draw[rounded corners, knot] (3.5, -1) rectangle (4.5, 2.5);
            \draw[fill=white, knot] (0,0) rectangle (4, 1.5);
            \draw[knot] (0.5, 0) -- (0.5, -1.25);
            \draw[knot] (1.5, 0) -- (1.5, -1.25);
            \draw[knot] (2.5, 0) -- (2.5, -1.25);
            \draw[knot] (0.5, 1.5) -- (0.5, 2.75);
            \draw[knot] (1.5, 1.5) -- (1.5, 2.75);
            \draw[knot] (2.5, 1.5) -- (2.5, 2.75);
            \draw[rounded corners, knot] (5.5, -1) rectangle (6.5, 2.5);
        } \circ \varphi_{H_S}$} (D);
    \draw[double equal sign distance] (D) to[out=180, in=0] (E);
}
\]
where the last equality follows because $\varphi_{
\tikz[baseline={([yshift=-.5ex]current bounding box.center)}, scale=.25]{
            \draw[knot, red] (4.5, 0.75) -- (5.5, 0.75);
            \draw[rounded corners, knot] (3.5, -1) rectangle (4.5, 2.5);
            \draw[fill=white, knot] (0,0) rectangle (4, 1.5);
            \draw[knot] (0.5, 0) -- (0.5, -1.25);
            \draw[knot] (1.5, 0) -- (1.5, -1.25);
            \draw[knot] (2.5, 0) -- (2.5, -1.25);
            \draw[knot] (0.5, 1.5) -- (0.5, 2.75);
            \draw[knot] (1.5, 1.5) -- (1.5, 2.75);
            \draw[knot] (2.5, 1.5) -- (2.5, 2.75);
            \draw[rounded corners, knot] (5.5, -1) rectangle (6.5, 2.5);
            \node at (2, 0.75) {$B$};
        }
        }^{-1} \cong \{1,0\}$. 
Now, slide $g$ before the birth; as before, to do so, we have to compensate by $\varphi_{H_2}: \varphi_{\mathcal{E}^w} \circ \{1,0\} \Rightarrow \{1,0\} \circ \varphi_{\mathcal{E}^w}$. Here is the resulting composition.
\[
\tikz{
     \node(A) at (0,0) 
        {$\varphi_{\mathcal{E}^w}~
        \tikz[baseline={([yshift=-.5ex]current bounding box.center)}, scale=.3]{
            \draw[knot] (0,0) rectangle (3, 1.5);
            \draw[knot] (0.5, 0) -- (0.5, -1.25);
            \draw[knot] (1.5, 0) -- (1.5, -1.25);
            \draw[knot] (2.5, 0) -- (2.5, -1.25);
            \draw[knot] (0.5, 1.5) -- (0.5, 2.75);
            \draw[knot] (1.5, 1.5) -- (1.5, 2.75);
            \draw[knot] (2.5, 1.5) -- (2.5, 2.75);
            \node at (1.5, 0.75) {$A$};
        }$};
    \node(B) at (6.5,0) 
        {$\varphi_{\mathcal{E}^w} \circ \{-1,0\} \circ \{1,0\}~
        \tikz[baseline={([yshift=-.5ex]current bounding box.center)}, scale=.3]{
            \draw[knot] (0,0) rectangle (3, 1.5);
            \draw[knot] (0.5, 0) -- (0.5, -1.25);
            \draw[knot] (1.5, 0) -- (1.5, -1.25);
            \draw[knot] (2.5, 0) -- (2.5, -1.25);
            \draw[knot] (0.5, 1.5) -- (0.5, 2.75);
            \draw[knot] (1.5, 1.5) -- (1.5, 2.75);
            \draw[knot] (2.5, 1.5) -- (2.5, 2.75);
            \node at (1.5, 0.75) {$A$};
        }$};
    \node(C) at (11.5,-1.5) 
        {$\{-1,0\} \circ \varphi_{\mathcal{E}^w} \circ \{1,0\}~
        \tikz[baseline={([yshift=-.5ex]current bounding box.center)}, scale=.3]{
            \draw[knot] (0,0) rectangle (3, 1.5);
            \draw[knot] (0.5, 0) -- (0.5, -1.25);
            \draw[knot] (1.5, 0) -- (1.5, -1.25);
            \draw[knot] (2.5, 0) -- (2.5, -1.25);
            \draw[knot] (0.5, 1.5) -- (0.5, 2.75);
            \draw[knot] (1.5, 1.5) -- (1.5, 2.75);
            \draw[knot] (2.5, 1.5) -- (2.5, 2.75);
            \node at (1.5, 0.75) {$A$};
        }$};
    \node (D) at (11.5,-4.5)
        {$\{-1,0\} \circ \{1,0\} \circ \varphi_{\mathcal{E}^w}~
        \tikz[baseline={([yshift=-.5ex]current bounding box.center)}, scale=.3]{
            \draw[knot] (0,0) rectangle (3, 1.5);
            \draw[knot] (0.5, 0) -- (0.5, -1.25);
            \draw[knot] (1.5, 0) -- (1.5, -1.25);
            \draw[knot] (2.5, 0) -- (2.5, -1.25);
            \draw[knot] (0.5, 1.5) -- (0.5, 2.75);
            \draw[knot] (1.5, 1.5) -- (1.5, 2.75);
            \draw[knot] (2.5, 1.5) -- (2.5, 2.75);
            \node at (1.5, 0.75) {$A$};
            \draw[rounded corners, dashed, red] (-.25, -.5) rectangle (3.25, 2);
        }$};
    \node (E) at (11.5, -7.5) {$\{-2,0\} \circ \{1,0\} ~\tikz[baseline={([yshift=-.5ex]current bounding box.center)}, scale=.35]{
            \draw[rounded corners, knot] (3.5, -1) rectangle (4.5, 2.5);
            \draw[fill=white, knot] (0,0) rectangle (4, 1.5);
            \draw[knot] (0.5, 0) -- (0.5, -1.25);
            \draw[knot] (1.5, 0) -- (1.5, -1.25);
            \draw[knot] (2.5, 0) -- (2.5, -1.25);
            \draw[knot] (0.5, 1.5) -- (0.5, 2.75);
            \draw[knot] (1.5, 1.5) -- (1.5, 2.75);
            \draw[knot] (2.5, 1.5) -- (2.5, 2.75);
            \node at (2, 0.75) {$B$};
            \draw[rounded corners, red, dashed] (5.5, -1) rectangle (6.5, 2.5);
        }$};
    \node (F) at (6.5, -9) {$\tikz[baseline={([yshift=-.5ex]current bounding box.center)}, scale=.3]{
            \draw[rounded corners, knot] (3.5, -1) rectangle (4.5, 2.5);
            \draw[fill=white, knot] (0,0) rectangle (4, 1.5);
            \draw[knot] (0.5, 0) -- (0.5, -1.25);
            \draw[knot] (1.5, 0) -- (1.5, -1.25);
            \draw[knot] (2.5, 0) -- (2.5, -1.25);
            \draw[knot] (0.5, 1.5) -- (0.5, 2.75);
            \draw[knot] (1.5, 1.5) -- (1.5, 2.75);
            \draw[knot] (2.5, 1.5) -- (2.5, 2.75);
            \draw[rounded corners, knot] (5.5, -1) rectangle (6.5, 2.5);
            \draw[rounded corners, dashed, red] (4.25, 0.1) rectangle (5.75, 1.4);
            \node at (2, 0.75) {$B$};
        }\,\{-2,0\}$};
    \node (G) at (0, -9) {$\tikz[baseline={([yshift=-.5ex]current bounding box.center)}, scale=.35]{
            \draw[rounded corners, knot] (3.5, -1) rectangle (4.5, 2.5);
            \draw[fill=white, knot] (0,0) rectangle (4, 1.5);
            \draw[knot] (0.5, 0) -- (0.5, -1.25);
            \draw[knot] (1.5, 0) -- (1.5, -1.25);
            \draw[knot] (2.5, 0) -- (2.5, -1.25);
            \draw[knot] (0.5, 1.5) -- (0.5, 2.75);
            \draw[knot] (1.5, 1.5) -- (1.5, 2.75);
            \draw[knot] (2.5, 1.5) -- (2.5, 2.75);
            \node at (2, 0.75) {$B$};
        }\,\{-1,0\}$};
    \draw[->] (A) to[out=0,in=180] node[pos=.5,above,arrows=-]
        {$\varphi_{H_B}$} 
        (B);
    \draw[->] (B) to[out=0,in=90] node[pos=.5,above,arrows=-]
        {$\lambda_{\phi(f)}$} 
        (C);
    \draw[->] (C) to[out=-90,in=90] node[pos=.5,right,arrows=-]
        {$\varphi_{H_2}$} (D);
    \draw[->] (D) to[out=-90,in=90] node[pos=.5,right,arrows=-]
        {$g$} (E);
    \draw[->] (E) to[out=-90, in=0] node[pos=.5,below,arrows=-]
        {$\mathcal{F}\left(\tikz[baseline={([yshift=-.5ex]current bounding box.center)}, scale=.4]{
	\draw (1,2) .. controls (1,1) and (2,1) .. (2,2);
	\draw (1,2) .. controls (1,1.75) and (2,1.75) .. (2,2);
	\draw (1,2) .. controls (1,2.25) and (2,2.25) .. (2,2);
}\right)$} (F);
    \draw[->] (F) to[out=180, in=0] node[pos=.5,above,arrows=-] {$\tikz[baseline={([yshift=-.5ex]current bounding box.center)}, scale=.22]{
            \draw[knot, red] (4.5, 0.75) -- (5.5, 0.75);
            \draw[rounded corners, knot] (3.5, -1) rectangle (4.5, 2.5);
            \draw[fill=white, knot] (0,0) rectangle (4, 1.5);
            \draw[knot] (0.5, 0) -- (0.5, -1.25);
            \draw[knot] (1.5, 0) -- (1.5, -1.25);
            \draw[knot] (2.5, 0) -- (2.5, -1.25);
            \draw[knot] (0.5, 1.5) -- (0.5, 2.75);
            \draw[knot] (1.5, 1.5) -- (1.5, 2.75);
            \draw[knot] (2.5, 1.5) -- (2.5, 2.75);
            \draw[rounded corners, knot] (5.5, -1) rectangle (6.5, 2.5);
        } \circ \varphi_{H_S}$} (G);
}
\]
Now, $\lambda_{\phi(f)}$ and $\varphi_{H_2}$ are inverse to one another, since $\lambda((a,b),(-1,0)) = X^{-a}Z^b$ and $\lambda((a,b),(1,0)) = X^a Z^{-b}$. Again, the birth cancels with the merge, and we have that $\phi(\psi(g)) = g$, concluding the proof.
\end{proof}

\begin{remark}
Since
\[
\deg_q\left(\varphi_{\left(\tikz[baseline={([yshift=-.5ex]current bounding box.center)}, scale=.15]{
            \draw[fill=white, knot] (0,0) rectangle (4, 1.5);
            \draw[knot, red, rounded corners] (3.5, 2.125) -- (4.5, 2.125) -- (4.5, -0.625) -- (3.5, -0.625);
            \draw[knot] (0.5, 0) -- (0.5, -1.25);
            \draw[knot] (1.5, 0) -- (1.5, -1.25);
            \draw[knot] (2.5, 0) -- (2.5, -1.25);
            \draw[knot] (0.5, 1.5) -- (0.5, 2.75);
            \draw[knot] (1.5, 1.5) -- (1.5, 2.75);
            \draw[knot] (2.5, 1.5) -- (2.5, 2.75);
            \draw[knot] (3.5, 0) -- (3.5, -1.25);
            \draw[knot] (3.5, 1.5) -- (3.5, 2.75);
            \node at (2, 0.75) {\scriptsize \textit{B}};
        },\, (0,1)\right)}\right)=0
\]
this result descends to Theorem 2.31 of \cite{https://doi.org/10.48550/arxiv.1405.2574} if we collapse the $\mathscr{G}$-grading to the $q$-grading.
\end{remark}

\subsection{Duals and mirrors}
\label{ss:duality}

Suppose $R$, $S$, and $T$ are $\mathscr{C}$-graded algebras. Per usual, we expect that if $M$ is a $\mathscr{C}$-graded $(R; S)$-multimodule, and $N$ is a $\mathscr{C}$-graded $(R; T)$-multimodule, then $\mathrm{Hom}_R(M, N)$ is an $(S; T)$-multimodule by
\[
\rho_L^{\mathrm{Hom}}(s, f) (m) := f(\rho_R^M(m,s)) \qquad \text{and} \qquad \rho_R^{\mathrm{Hom}}(f,t)(m) := \rho_R^N(f(m), t)
\]
for each $f\in \mathrm{Hom}_R(M, N)$, $m\in M$, $s\in S$ and $t\in T$. \textit{However}, $\mathrm{Hom}_R(M, N)$ does not satisfy the axioms of a $\mathscr{C}$-graded multimodule: by definition, $\mathrm{Hom}_R(M, N)$ is graded by $\widetilde{\mathscr{I}} \times \mathbb{Z} = (\mathscr{I} \sqcup \{\mathrm{Id}\}) \times \mathbb{Z}$, and the reader is invited to verify that
\begin{itemize}
    \item $\rho_L^{\mathrm{Hom}} (s_1 \cdot s_2, f) (m) = \alpha\left( \abs{m}, \abs{s_1}, \abs{s_2} \right)^{-1} \rho_L^{\mathrm{Hom}}(s_1, \rho_L^{\mathrm{Hom}}(s_2, f))(m)$, 
    \item $\rho_R^{\mathrm{Hom}}(\rho_R^{\mathrm{Hom}}(f, t_1), t_2)(m) = \alpha \left(\abs{f(m)}, \abs{t_1}, \abs{t_2} \right) \rho_R^{\mathrm{Hom}}(f, t_1 \cdot t_2)(m)$, and
    \item $\rho_R^{\mathrm{Hom}} (\rho_L^\mathrm{Hom}(s, f), t) (m) = \rho_L^\mathrm{Hom}(s, \rho_R^\mathrm{Hom}(f,t)) (m)$.
\end{itemize}

Despite this ambiguity, we are able to give a type of duality statement which turns out to be a generalization of Theorem 4.12 in \cite{hogancamp2012morphismscategorifiedspinnetworks}. This implies a unified analogue to Lemma 4.14 of \cite{stoffregen2024joneswenzlprojectorskhovanovhomotopy}, which is all we will need to prove the uniqueness of unified Cooper-Krushkal projectors.

We \textit{dualize} a flat diskular $(m;n)$-tangle $T$ by the following operation, flipping radially,
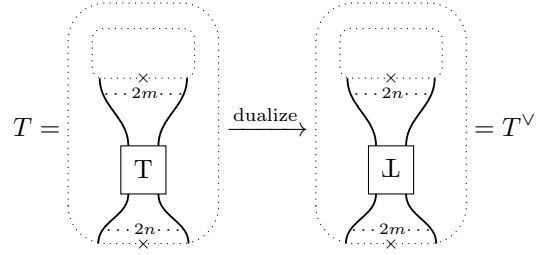
\[
T = 
\tikz[baseline={([yshift=-.5ex]current bounding box.center)}, scale=.8]{
    \draw[rounded corners=5mm, dotted] (0, 0) rectangle (2.5, 4) {};
    \draw[rounded corners, dotted] (.4, 2.75) rectangle (2.1, 3.575);
    \draw (0.875, 0.825) rectangle (1.625, 1.625);
    \node at (1.25, 1.25) {T};
    \draw[knot] (1, 0.825) to[out=-90, in=90] (0.5, 0);
    \draw[knot] (1.5, 0.825) to[out=-90, in=90] (2, 0);
    \draw[knot] (1, 1.625) to[out=90, in=-90] (0.525, 2.75);
    \draw[knot] (1.5, 1.625) to[out=90, in=-90] (1.975, 2.75);
    \node at (1.25, 0.25) {\tiny $\dotsm 2n \dotsi$};
    \node at (1.25, 2.5) {\tiny $\dotsm 2m \dotsi$};
    \node at (1.25, 0) {\tiny$\times$};
    \node at (1.25, 2.75) {\tiny $\times$};
}
\xrightarrow{\text{dualize}}
\tikz[baseline={([yshift=-.5ex]current bounding box.center)}, scale=.8]{
    \draw[rounded corners=5mm, dotted] (0, 0) rectangle (2.5, 4) {};
    \draw[rounded corners, dotted] (.4, 2.75) rectangle (2.1, 3.575);
    \draw (0.875, 0.825) rectangle (1.625, 1.625);
    \node at (1.25, 1.25) {\rotatebox{180}{T}};
    \draw[knot] (1, 0.825) to[out=-90, in=90] (0.5, 0);
    \draw[knot] (1.5, 0.825) to[out=-90, in=90] (2, 0);
    \draw[knot] (1, 1.625) to[out=90, in=-90] (0.525, 2.75);
    \draw[knot] (1.5, 1.625) to[out=90, in=-90] (1.975, 2.75);
    \node at (1.25, 0.25) {\tiny $\dotsm 2m \dotsi$};
    \node at (1.25, 2.5) {\tiny $\dotsm 2n \dotsi$};
    \node at (1.25, 0) {\tiny$\times$};
    \node at (1.25, 2.75) {\tiny $\times$};
}
= T^\vee
\]
to obtain a diskular $(n;m)$-tangle. Notice that if $T$ is a flat diskular $n$-tangle, then $T^\vee$ is a flat diskular $(n; 0)$-tangle; this is the case we are most interested in.
On cobordisms of $T$ embedded in $[0,1]^3$, $(-)^\vee$ acts by the transformation $(x,y,z) \mapsto (x, 1-y, 1-z)$. 

Now we describe how $(-)^\vee$ establishes a contravariant functor $\mathrm{Chom}(n)^\mathscr{G} \to \mathrm{Chom}(n)^\mathscr{G}$. On objects (which are chain complexes of summands of $\mathscr{G}$-graded $H^n$-modules associated to flat diskular $n$-tangles with a differential of matrices of cobordisms), $(-)^\vee$ applies $(-)^\vee$ on each entry, reverses homological degree (i.e., $(A^\vee)^k := (A^{-k})^\vee$), applies $(-)^\vee$ on each cobordism and takes the transpose of each matrix of cobordisms, and reverses $\mathscr{G}$-degree. By this last point, we mean that each cobordism shift $W$ is dualized (note that if $W: a\to b$, then $W^\vee: b^\vee \to a^\vee$) and $\mathbb{Z}\times\mathbb{Z}$-degree is reversed: $\{v_1, v_2\}^\vee = \{-v_2, -v_1\}$.

In particular, if $d_A$ is the differential for $A\in \mathrm{Chom}(n)^\mathscr{G}$, then (abusing notation), fix
\[
d_{A^\vee} = -(d_A)^\vee \circ \varphi_H
\]
where $\varphi_H$ means that we are applying the change of chronology
\[
\varphi_H: \mathrm{Id} \Rightarrow \varphi_{(d_A)^\vee}^{-1} \circ \varphi_{(d_A)^\vee}
\]
on each entry of each matrix comprising $d_{A^\vee}$. For example, the dual of the complex
\[
 \varphi_{\tikz[baseline=-6.5ex, scale=.45]
{
	\draw[dotted] (3,-2) circle(0.707);
        \draw[knot, red,thick] (3,-1.7) -- (3,-2.3);
	\draw[knot] (2.5,-1.5) .. controls (2.75,-1.75) and (3.25,-1.75) .. (3.5,-1.5);
	\draw[knot] (2.5,-2.5) .. controls  (2.75,-2.25) and (3.25,-2.25) .. (3.5,-2.5);
}}
\mathcal{F}\left(
\tikz[baseline={([yshift=-.5ex]current bounding box.center)}, scale=.75]
{
	\draw[dotted] (3,-2) circle(0.707);
	\draw[knot] (2.5,-1.5) .. controls (2.75,-1.75) and (3.25,-1.75) .. (3.5,-1.5);
	\draw[knot] (2.5,-2.5) .. controls  (2.75,-2.25) and (3.25,-2.25) .. (3.5,-2.5);
}
  \right) \{-1,0\}
\xrightarrow{\mathcal{F}\left( \tikz[baseline=-6.5ex, scale=.45]
{
	\draw[dotted] (3,-2) circle(0.707);
        \draw[knot, red,thick] (3,-1.7) -- (3,-2.3);
	\draw[knot] (2.5,-1.5) .. controls (2.75,-1.75) and (3.25,-1.75) .. (3.5,-1.5);
	\draw[knot] (2.5,-2.5) .. controls  (2.75,-2.25) and (3.25,-2.25) .. (3.5,-2.5);
} \right)}
\mathcal{F}\left(
\tikz[baseline={([yshift=-.5ex]current bounding box.center)}, scale=.75]
{
	\draw[dotted] (-2,-2) circle(0.707);
	\draw[knot] (-1.5,-1.5) .. controls (-1.75,-1.75) and (-1.75,-2.25) .. (-1.5,-2.5);
	\draw[knot] (-2.5,-1.5) .. controls (-2.25,-1.75) and (-2.25,-2.25) ..  (-2.5,-2.5);
}
\right) \{-1,0\}
\]
is the complex
\[
\mathcal{F}\left(
\tikz[baseline={([yshift=-.5ex]current bounding box.center)}, scale=.75]
{
	\draw[dotted] (-2,-2) circle(0.707);
	\draw[knot] (-1.5,-1.5) .. controls (-1.75,-1.75) and (-1.75,-2.25) .. (-1.5,-2.5);
	\draw[knot] (-2.5,-1.5) .. controls (-2.25,-1.75) and (-2.25,-2.25) ..  (-2.5,-2.5);
} 
  \right) \{0,1\}
\xrightarrow{
\mathcal{F}\left( \tikz[baseline=8.25ex, scale=.45]
{
    \begin{scope}[rotate=90]
	\draw[dotted] (3,-2) circle(0.707);
        \draw[knot, red,thick] (3,-1.7) -- (3,-2.3);
	\draw[knot] (2.5,-1.5) .. controls (2.75,-1.75) and (3.25,-1.75) .. (3.5,-1.5);
	\draw[knot] (2.5,-2.5) .. controls  (2.75,-2.25) and (3.25,-2.25) .. (3.5,-2.5);
    \end{scope}
} \right)
\circ \varphi_H}
\varphi_{\tikz[baseline=8.25ex, scale=.45]
{
    \begin{scope}[rotate=90]
	\draw[dotted] (3,-2) circle(0.707);
        \draw[knot, red,thick] (3,-1.7) -- (3,-2.3);
	\draw[knot] (2.5,-1.5) .. controls (2.75,-1.75) and (3.25,-1.75) .. (3.5,-1.5);
	\draw[knot] (2.5,-2.5) .. controls  (2.75,-2.25) and (3.25,-2.25) .. (3.5,-2.5);
    \end{scope}
}} ^{-1}
\mathcal{F}\left(
\tikz[baseline={([yshift=-.5ex]current bounding box.center)}, scale=.75]
{
	\draw[dotted] (3,-2) circle(0.707);
	\draw[knot] (2.5,-1.5) .. controls (2.75,-1.75) and (3.25,-1.75) .. (3.5,-1.5);
	\draw[knot] (2.5,-2.5) .. controls  (2.75,-2.25) and (3.25,-2.25) .. (3.5,-2.5);
}
\right) \{0,1\}
\]
for $\varphi_H: \mathrm{Id} \Rightarrow \varphi_{\tikz[baseline=8.25ex, scale=.45]
{
    \begin{scope}[rotate=90]
	\draw[dotted] (3,-2) circle(0.707);
        \draw[knot, red,thick] (3,-1.7) -- (3,-2.3);
	\draw[knot] (2.5,-1.5) .. controls (2.75,-1.75) and (3.25,-1.75) .. (3.5,-1.5);
	\draw[knot] (2.5,-2.5) .. controls  (2.75,-2.25) and (3.25,-2.25) .. (3.5,-2.5);
    \end{scope}
}} ^{-1} \circ \varphi_{\tikz[baseline=8.25ex, scale=.45]
{
    \begin{scope}[rotate=90]
	\draw[dotted] (3,-2) circle(0.707);
        \draw[knot, red,thick] (3,-1.7) -- (3,-2.3);
	\draw[knot] (2.5,-1.5) .. controls (2.75,-1.75) and (3.25,-1.75) .. (3.5,-1.5);
	\draw[knot] (2.5,-2.5) .. controls  (2.75,-2.25) and (3.25,-2.25) .. (3.5,-2.5);
    \end{scope}
}}$. In particular, this is to say that 
\[
\mathrm{Kh}\left( 
\tikz[baseline={([yshift=-.5ex]current bounding box.center)}, scale=.75]
{
	\draw[dotted] (.5,.5) circle(0.707);
	\draw[knot, ->](0,0) -- (1,1);
	\fill[fill=white] (.5,.5) circle (.15);
	\draw[knot, ->](1,0) -- (0,1);
         \node at (1.207, 0.5) {$\times$};
}\right) ^\vee
=
\mathrm{Kh}\left(
\tikz[baseline={([yshift=-.5ex]current bounding box.center)}, scale=.75]
{
	\draw[dotted] (.5,.5) circle(0.707);
	\draw[knot, <-](1,0) -- (0,1);
	\fill[fill=white] (.5,.5) circle (.15);
	\draw[knot, <-](0,0) -- (1,1);
        \node at (1.207, 0.5) {$\times$};
}
 \right) 
\]
as one might hope.

Finally, on morphisms, to $f\in \mathrm{Hom}^k_{\mathrm{Chom}(n)}(\varphi_{W, (v_1, v_2)}A, B)$ (where $k$ is the homological degree and $(W, (v_1, v_2))$ is the $\widetilde{\mathscr{I}}$-degree) we define $f^\vee\in \mathrm{Hom}^k_{\mathrm{Chom}(n)}(\varphi_{W^\vee, (v_2, v_1)}B^\vee, A^\vee)$ to be
\[
(f^\vee)_i = (-1)^{ik}(f_{-i-k})^\vee
\]
following the commutativity of the square 
\[
\begin{tikzcd}[column sep = huge, row sep = large]
(B^\vee)_i \arrow[d, equals] \arrow[r, "(f^\vee)_i"] & (A^\vee)_{i+k} \arrow[d, equals] \\
(B_{-i})^\vee \arrow[r, "(-1)^{ik} (f_{-i-k})^\vee"] & (A_{-i-k})^\vee
\end{tikzcd}
\]
As consequence of reversing $\mathscr{G}$-degree, the $\widetilde{\mathscr{I}}$-degree of compositions of morphisms is also reversed; this is to say that (for, say, maps of homological degree zero) $(g \circ_\mathscr{G} f)^\vee = \varphi_{H_{f,g}} g^\vee \circ_\mathscr{G} f^\vee$ where $\varphi_{H_{f,g}}$ denotes the change of chronology prioritizing the degree shift of $g$ before that of $f$. Then, we have the following standard lemma.

\begin{lemma} For $A, B$ and $f, g$ as above,
\begin{enumerate}
    \item $(-)^\vee$ induces a degree-zero chain map
\[
\mathrm{Hom}_{\mathrm{Chom}(n)}(A, B) \to \mathrm{Hom}_{\mathrm{Chom}(n)}(B^\vee, A^\vee);
\]
    \item $(g\circ_\mathscr{G} f)^\vee = \varphi_{H_{f,g}} (-1)^{\abs{f}_h\abs{g}_h} f^\vee \circ_\mathscr{G} g^\vee.$
\end{enumerate}
\end{lemma}

The purpose of the rest of this section is to prove that
\[
\normalfont{\textsc{Hom}}_n(A \otimes \delta,B) \cong \normalfont{\textsc{Hom}}_n(A, B \otimes \delta^\vee) 
\]
for any $A, B\in \mathrm{Chom}(n)$ and $\delta$ any flat diskular $n$-tangle. In order to describe our logical process for proving this statement, we will introduce yet another tensor product which will not reappear anywhere else in this paper.

\begin{definition}
Suppose $A, B\in \mathrm{Chom}(n)^\mathscr{G}$. Recall that we may represent, for example, $A$ and $A^\vee$ as
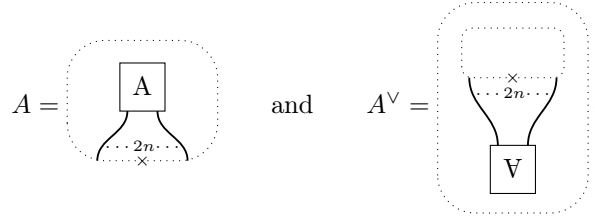
\[
A = 
\tikz[baseline={([yshift=-.5ex]current bounding box.center)}, scale=.8]{
    \draw[rounded corners=5mm, dotted] (0, 0) rectangle (2.5, 2) {};
    \draw (0.875, 0.825) rectangle (1.625, 1.625);
    \node at (1.25, 1.25) {A};
    \draw[knot] (1, 0.825) to[out=-90, in=90] (0.5, 0);
    \draw[knot] (1.5, 0.825) to[out=-90, in=90] (2, 0);
    \node at (1.25, 0.25) {\tiny $\dotsm 2n \dotsi$};
    \node at (1.25, 0) {\tiny$\times$};
}
\qquad \text{and} \qquad
A^\vee = 
\tikz[baseline={([yshift=-.5ex]current bounding box.center)}, scale=.8]{
    \draw[rounded corners=5mm, dotted] (0, 0.5) rectangle (2.5, 4);
    \draw[rounded corners, dotted] (.4, 2.75) rectangle (2.1, 3.575);
    \draw (0.875, 0.825) rectangle (1.625, 1.625);
    \node at (1.25, 1.25) {\rotatebox{180}{A}};

    \draw[knot] (1, 1.625) to[out=90, in=-90] (0.525, 2.75);
    \draw[knot] (1.5, 1.625) to[out=90, in=-90] (1.975, 2.75);
    \node at (1.25, 2.5) {\tiny $\dotsm 2n \dotsi$};
    \node at (1.25, 2.75) {\tiny $\times$};
}\,.
\]
We define two natural operations. By $A \mid B^\vee$, we mean the tensor $A \otimes_{H^n} B^\vee$; on the other hand, by $A^\vee \mid B$, we mean the tensor $A^\vee \otimes_{H^0} B$. Diagramatically,
\[
A \mid B^\vee \cong 
\tikz[baseline={([yshift=-.5ex]current bounding box.center)}, scale=.8]{
    \draw (0.875, 0.825) rectangle (1.625, 1.625);
    \node at (1.25, 1.25) {A};
    \draw[knot] (1, 0.825) to[out=-90, in=90] (0.5, 0);
    \draw[knot] (1.5, 0.825) to[out=-90, in=90] (2, 0);
    \node at (1.25, 0) {\tiny $\dotsm 2n \dotsi$};
    \draw[knot] (0.5, 0) to[out=-90, in=90] (1, -0.825);
    \draw[knot] (2, 0) to[out=-90, in=90] (1.5, -0.825);
    \draw (0.875, -0.825) rectangle (1.625, -1.625);
    \node at (1.25, -1.25) {\scalebox{-1}[1]{B}};
    \draw[rounded corners, dotted] (0,-2) rectangle (2.5, 2);
}
\qquad \text{and} \qquad
A^\vee \mid B \cong 
\tikz[baseline={([yshift=-.5ex]current bounding box.center)}, scale=.8]{
    \draw (0.875, 0.825) rectangle (1.625, 1.625);
    \node at (1.25, 1.25) {B};
    \draw[knot] (1, 0.825) to[out=-90, in=90] (0.5, 0);
    \draw[knot] (1.5, 0.825) to[out=-90, in=90] (2, 0);
    \node at (1.25, 0.25) {\tiny $\dotsm 2n \dotsi$};
    \draw[knot] (1, 2.625) to[out=90, in=-90] (0.5, 3.45);
    \draw[knot] (1.5, 2.625) to[out=90, in=-90] (2, 3.45);
    \draw (0.875, 1.825) rectangle (1.625, 2.625);
    \node at (1.25, 2.25) {\rotatebox{180}{A}};
    \node at (1.25, 3.2) {\tiny $\dotsm 2n \dotsi$};
    \draw[rounded corners, dotted] (0,0) rectangle (2.5, 4.25);
    \draw[rounded corners, dotted] (0.25, 3.45) rectangle (2.25, 4);
    \node at (1.25, 0) {$\times$};
    \node at (1.25, 3.45) {$\times$}; 
}
\]
by Theorem \ref{thm:multigluing}.
\end{definition}

\begin{theorem}[cf. Theorem 4.12, \cite{hogancamp2012morphismscategorifiedspinnetworks}]
\label{thm:dualclosing}
Suppose $A, B\in \mathrm{Chom}(n)^\mathscr{G}$. Then there is an isomorphism of complexes
\[
\normalfont{\textsc{Hom}}_n(A, B) \cong \normalfont{\textsc{Hom}}_0(\emptyset, B\mid A^\vee \{-n,0\}).
\]
\end{theorem}

This Theorem implies our goal for the section.

\begin{corollary}
\label{cor:dualswap}
Suppose $\delta$ is a flat diskular $n$-tangle. Then 
\[
\normalfont{\textsc{Hom}}_n (A \otimes \mathcal{F}(\delta), B) \cong \normalfont{\textsc{Hom}}_n (A, B\otimes \mathcal{F}(\delta^\vee)).
\]
\end{corollary}

\begin{proof}
Writing $\delta$ for $\mathcal{F}(\delta)$, we have 
\begin{align*}
    \normalfont{\textsc{Hom}}_n(A \otimes \delta, B) & \cong \normalfont{\textsc{Hom}}_0(\emptyset, B \mid (A \otimes \delta)^\vee \{-n, 0\}) \\
    &\cong \normalfont{\textsc{Hom}}_0 (\emptyset, B \mid (\delta^\vee \otimes A^\vee)\{-n, 0\}) \\
    &\cong \normalfont{\textsc{Hom}}_0 (\emptyset, (B \otimes \delta^\vee) \mid A^\vee\{-n, 0\})\\
    & \cong \normalfont{\textsc{Hom}}_n(A, B \otimes \delta^\vee).
\end{align*}
The first and last isomorphisms are provided by Theorem \ref{thm:dualclosing}, while the second follows from the definition of $(-)^\vee$ and the third is an application of Theorem \ref{thm:multigluing}.
\end{proof}

We prove Theorem \ref{thm:dualclosing} in two steps. First, we prove an analogue of Theorem \ref{thm:dualclosing} for crossingless matchings. Then, we argue that this implies the general statement.

\begin{definition}
Suppose $a$ is a crossingless matching on $2n$ points; i.e., a planar diskular $n$-tangle. In this definition, we will assume $a$ is \textit{indecomposable}; that is, $a$ is void of circle components.
\begin{enumerate}
    \item Define $\eta_a$ as the map
    \[
        \eta_a: \emptyset \xrightarrow{\varphi_{H_{\eta_a}}} \{-n, 0\} \circ \{n,0\} \, \emptyset \longrightarrow \{-n, 0\} \, a \mid a^\vee
    \]
    consisting of $n$-many births (since $a \mid a^\vee$ is exactly $n$-many circles).
    \item Let $s_a$ denote the map
    \[
        s_a: \varphi_{\Sigma_a} a^\vee \mid a \to 1_n
    \]
    defined by the minimal chronological cobordism $\Sigma_a$ given by contracting symmetric arcs, right-to-left, with framing pointed upwards.
\end{enumerate}
\end{definition}

The following lemma is apparent.

\begin{lemma}
\label{lem:CPmincob}
Fix indecomposable crossingless parings on $2n$-points $a, b$. 
\begin{enumerate}
    \item $(\mathbbm{1}_a \mid s_a) \circ (\eta_a \mid \mathbbm{1}_a) = \mathbbm{1}_a$ and $(s_a \mid \mathbbm{1}_{a^\vee}) \circ (\mathbbm{1}_{a^\vee} \circ \eta_a) = \mathbbm{1}_{a^\vee}$.
    \item Suppose $b \mid a^\vee$ consists of $\ell$-many circles, $1 \le \ell \le n$ (note that $\ell = n$ if and only if $b=a$). Then
    \[
        \mathbbm{1}_b \mid s_a : b \mid a^\vee \mid a \to b
    \]
    consists of $\ell$-many merges followed by a minimal cobordism $W: a\to b$.
\end{enumerate}
\end{lemma}

Note that $W$ consists of $(n-\ell)$-many saddles. We'll write $\abs{b \mid a^\vee}$ to denote the number of loops in $b \mid a^\vee$ (above, $\abs{b \mid a^\vee} = \ell$). We'll denote crossingless matchings, pictorially, as
\[
a = \tikz[baseline={([yshift=-.5ex]current bounding box.center)}, scale=1]{
    \draw[knot] (0,0) -- (0, 0.25) to[out=90, in=90] (1, 0.25) -- (1,0);
    \draw[fill=white] (0.5, 0.55) circle (7.5pt);
    \node at (0.5, 0.55) {$a$};
}
\qquad\text{and}\qquad
a^\vee = \tikz[baseline={([yshift=-.5ex]current bounding box.center)}, scale=1]{
    \draw[knot] (0,1) -- (0, 0.75) to[out=-90, in=-90] (1, 0.75) -- (1,1);
    \draw[fill=white] (0.5, 0.45) circle (7.5pt);
    \node at (0.5, 0.45) {$a^\vee$};
}\,.
\]
For example, part 2 of Lemma \ref{lem:CPmincob} describes a cobordism
\[
\tikz[baseline={([yshift=-.5ex]current bounding box.center)}, scale=1]{
    \draw[knot] (0,0) -- (0, 0.25) to[out=90, in=90] (1, 0.25) -- (1,0);
    \draw[fill=white] (0.5, 0.55) circle (7.5pt);
    \node at (0.5, 0.55) {$a$};
    \begin{scope}[yshift=.75cm]
        \draw[knot] (0,1) -- (0, 0.75) to[out=-90, in=-90] (1, 0.75) -- (1,1);
        \draw[fill=white] (0.5, 0.45) circle (7.5pt);
        \node at (0.5, 0.45) {$a^\vee$};
    \end{scope}
    \begin{scope}[yshift=1.75cm]
        \draw [knot] (0,0) -- (0, 0.25) to[out=90, in=90] (1, 0.25) -- (1,0);
        \draw[fill=white] (0.5, 0.55) circle (7.5pt);
        \node at (0.5, 0.55) {$b$};
    \end{scope}
    \draw[red, rounded corners, dashed] (-0.1, 0.125) rectangle (1.1, 1.75);
}~
\xrightarrow{\mathbbm{1}_b \mid s_a}~
\tikz[baseline={([yshift=-.5ex]current bounding box.center)}, scale=1]{
    \draw[knot] (0,0) -- (0, 1.75);
    \draw[knot] (1,0) -- (1, 1.75);
    \begin{scope}[yshift=1.75cm]
        \draw [knot] (0,0) -- (0, 0.25) to[out=90, in=90] (1, 0.25) -- (1,0);
        \draw[fill=white] (0.5, 0.55) circle (7.5pt);
        \node at (0.5, 0.55) {$b$};
    \end{scope}
}
~=~ 
\tikz[baseline={([yshift=-.5ex]current bounding box.center)}, scale=1]{
        \draw [knot] (0,0) -- (0, 0.25) to[out=90, in=90] (1, 0.25) -- (1,0);
        \draw[fill=white] (0.5, 0.55) circle (7.5pt);
        \node at (0.5, 0.55) {$b$};
}\,.
\]
While these pictures are a departure from the planar arc diagrams we are accustomed to, they are a little more natural for the proof of the following proposition.

\begin{proposition}[cf. Proposition 4.8, \cite{hogancamp2012morphismscategorifiedspinnetworks}]
\label{prop:CPadjunction2}
Suppose $a$ and $b$ are crossingless matchings on $2n$ points (not necessarily indecomposable) and fix a minimal cobordism $W: \widehat{a} \to \widehat{b}$, where $\widehat{a}$, $\widehat{b}$ are $a$ and $b$ with circle components removed. Then
\[
\mathrm{Hom}_n\left(\varphi_{\left(W, (n - \abs{\widehat{b} \mid \widehat{a}}, 0)\right)}\, a, b \right)
\cong
\mathrm{Hom}_0 \left( \emptyset, b \mid a^\vee \{-n, 0\}\right).
\]
In pictures, 
\[
\mathrm{Hom}_n\left(\varphi_{\left(W, (n - \abs{\widehat{b} \mid \widehat{a}}, 0)\right)}\,
\tikz[baseline={([yshift=-.5ex]current bounding box.center)}, scale=1]{
        \draw [knot] (0,0) -- (0, 0.25) to[out=90, in=90] (1, 0.25) -- (1,0);
        \draw[fill=white] (0.5, 0.55) circle (7.5pt);
        \node at (0.5, 0.55) {$a$};
}
\,, \,
\tikz[baseline={([yshift=-.5ex]current bounding box.center)}, scale=1]{
        \draw [knot] (0,0) -- (0, 0.25) to[out=90, in=90] (1, 0.25) -- (1,0);
        \draw[fill=white] (0.5, 0.55) circle (7.5pt);
        \node at (0.5, 0.55) {$b$};
}
\right)
\cong
\mathrm{Hom}_0 \left( \emptyset, \,
\tikz[baseline={([yshift=-.5ex]current bounding box.center)}, scale=1]{
        \draw[knot] (0,1) -- (0, 0.75) to[out=-90, in=-90] (1, 0.75) -- (1,1);
        \draw[fill=white] (0.5, 0.45) circle (7.5pt);
        \node at (0.5, 0.45) {$a^\vee$};
    \begin{scope}[yshift=1cm]
        \draw [knot] (0,0) -- (0, 0.25) to[out=90, in=90] (1, 0.25) -- (1,0);
        \draw[fill=white] (0.5, 0.55) circle (7.5pt);
        \node at (0.5, 0.55) {$b$};
    \end{scope}
}
\{-n, 0\}\right).
\]
\end{proposition}

\begin{proof}
First, we can assume without loss of generality that $a$ and $b$ are both indecomposable---the general result follows immediately by delooping.

Proceeding, we will frequently denote $\varphi_{\left(W, (n - \abs{\widehat{b} \mid \widehat{a}}, 0)\right)}$ by $\varphi_{W^N}$. Notice that the $\widetilde{\mathscr{I}}$-degree of any $f\in \mathrm{Hom}_n\left(\varphi_{W^N}a, b\right)$ can be chosen to be described purely by a $\mathbb{Z} \times \mathbb{Z}$-degree, since $W:a\to b$ is minimal. Recall that this is also the case for any $g\in \mathrm{Hom}_0\left(\emptyset, b \mid a^\vee \, \{-n,0\}\right)$ since any grading shift associated to a cobordism between closed diagrams is canonically isomorphic to a pure $\mathbb{Z}\times\mathbb{Z}$-shift. Thus, we will denote the homogeneous degree of $f$ and $g$ by $v_f$ and $w_g\in\mathbb{Z} \times \mathbb{Z}$ respectively.

The rest of this proof proceeds like the proof of Theorem \ref{thm:adjunction1}. To any $f\in \mathrm{Hom}_n\left(\varphi_{W^n} a, b\right)$, define $\phi(f) \in \mathrm{Hom}_0 (\emptyset, b \mid a^\vee\, \{-n, 0\})$ as the composition
\[
\tikz{
     \node(A) at (0,0) 
        {$v_f \, \emptyset$};
    \node(B) at (4,0) 
        {$v_f \circ \{-n, 0\} \, 
        \tikz[baseline={([yshift=-.5ex]current bounding box.center)}, scale=1]{
        \draw[knot] (0,1) -- (0, 0.75) to[out=-90, in=-90] (1, 0.75) -- (1,1);
        \draw[fill=white] (0.5, 0.45) circle (7.5pt);
        \node at (0.5, 0.45) {$a^\vee$};
    \begin{scope}[yshift=1cm]
        \draw [knot] (0,0) -- (0, 0.25) to[out=90, in=90] (1, 0.25) -- (1,0);
        \draw[fill=white] (0.5, 0.55) circle (7.5pt);
        \node at (0.5, 0.55) {$a$};
    \end{scope}
    }
        $};
    \node(C) at (10,0) 
        {$\{-n, 0\} \circ v_f \, 
        \tikz[baseline={([yshift=-.5ex]current bounding box.center)}, scale=1]{
        \draw[knot] (0,1) -- (0, 0.75) to[out=-90, in=-90] (1, 0.75) -- (1,1);
        \draw[fill=white] (0.5, 0.45) circle (7.5pt);
        \node at (0.5, 0.45) {$a^\vee$};
    \begin{scope}[yshift=1cm]
        \draw [knot] (0,0) -- (0, 0.25) to[out=90, in=90] (1, 0.25) -- (1,0);
        \draw[fill=white] (0.5, 0.55) circle (7.5pt);
        \node at (0.5, 0.55) {$a$};
        \draw[red, rounded corners, dashed] (-0.1, 0) rectangle (1.1, 1);
    \end{scope}
    }
        $};
    \node(D) at (14.5,0) 
        {$\{-n, 0\}\, 
        \tikz[baseline={([yshift=-.5ex]current bounding box.center)}, scale=1]{
        \draw[knot] (0,1) -- (0, 0.75) to[out=-90, in=-90] (1, 0.75) -- (1,1);
        \draw[fill=white] (0.5, 0.45) circle (7.5pt);
        \node at (0.5, 0.45) {$a^\vee$};
    \begin{scope}[yshift=1cm]
        \draw [knot] (0,0) -- (0, 0.25) to[out=90, in=90] (1, 0.25) -- (1,0);
        \draw[fill=white] (0.5, 0.55) circle (7.5pt);
        \node at (0.5, 0.55) {$b$};
    \end{scope}
    }
        $};
    \draw[->] (A) to[out=0,in=180] node[pos=.5,above,arrows=-]
        {$\eta_a$} 
        (B);
    \draw[->] (B) to[out=0,in=180] node[pos=.5,above,arrows=-]
        {$\lambda(v_f, (-n,0))$} 
        (C);  
    \draw[->] (C) to[out=0,in=180] node[pos=.5,above,arrows=-]
        {$f \mid \mathbbm{1}_{a^\vee}$} 
        (D);  
}
\]
To clear up any confusion, notice that the minimal cobordism $W: a\to b$, which has $\left(n - \abs{b \mid a^\vee}\right)$-many saddles, extends to a cobordism $W \bullet \mathbbm{1}_{a^\vee}: a \mid a^\vee \to b \mid a^\vee$ in which all saddles are realized as merges. Thus $\varphi_{(W\bullet \mathbbm{1}_{a^\vee})^N} \cong \mathrm{Id}$.

Next, to $g\in \mathrm{Hom}_0 (\emptyset, b\mid a^\vee\, \{-n, 0\})$, define $\psi(g) \in \mathrm{Hom}_n(\varphi_{W^N} a, b)$ by
\[
\tikz{
     \node(A) at (0,0) 
        {$w_g \circ \varphi_{W^N} \, 
        \tikz[baseline={([yshift=-.5ex]current bounding box.center)}, scale=1]{
            \draw[knot] (0,0) -- (0, 0.25) to[out=90, in=90] (1, 0.25) -- (1,0);
            \draw[fill=white] (0.5, 0.55) circle (7.5pt);
            \node at (0.5, 0.55) {$a$};
        }
    $};
    \node(B) at (3.7,0) 
        {$\varphi_{W^N} \circ w_g \, 
        \tikz[baseline={([yshift=-.5ex]current bounding box.center)}, scale=1]{
            \draw[knot] (0,0) -- (0, 0.25) to[out=90, in=90] (1, 0.25) -- (1,0);
            \draw[fill=white] (0.5, 0.55) circle (7.5pt);
            \node at (0.5, 0.55) {$a$};
            \draw[red, rounded corners, dashed] (-0.1, 0.9) rectangle (1.1, 1.9);
        }
    $};
    \node(C) at (8.15,0) 
        {$\varphi_{W^N} \circ \{-n,0\} \,
        \tikz[baseline={([yshift=-.5ex]current bounding box.center)}, scale=1]{
            \draw[knot] (0,0) -- (0, 0.25) to[out=90, in=90] (1, 0.25) -- (1,0);
            \draw[fill=white] (0.5, 0.55) circle (7.5pt);
            \node at (0.5, 0.55) {$a$};
            \begin{scope}[yshift=.75cm]
                \draw[knot] (0,1) -- (0, 0.75) to[out=-90, in=-90] (1, 0.75) -- (1,1);
                \draw[fill=white] (0.5, 0.45) circle (7.5pt);
                \node at (0.5, 0.45) {$a^\vee$};
            \end{scope}
            \begin{scope}[yshift=1.75cm]
                \draw [knot] (0,0) -- (0, 0.25) to[out=90, in=90] (1, 0.25) -- (1,0);
                \draw[fill=white] (0.5, 0.55) circle (7.5pt);
                \node at (0.5, 0.55) {$b$};
            \end{scope}
        }
    $};
    \node(D) at (13.25,0) 
        {$\varphi_{\left(W, (-\abs{b\mid a^\vee}, 0)\right)} \,
        \tikz[baseline={([yshift=-.5ex]current bounding box.center)}, scale=1]{
            \draw[knot] (0,0) -- (0, 0.25) to[out=90, in=90] (1, 0.25) -- (1,0);
            \draw[fill=white] (0.5, 0.55) circle (7.5pt);
            \node at (0.5, 0.55) {$a$};
            \begin{scope}[yshift=.75cm]
                \draw[knot] (0,1) -- (0, 0.75) to[out=-90, in=-90] (1, 0.75) -- (1,1);
                \draw[fill=white] (0.5, 0.45) circle (7.5pt);
                \node at (0.5, 0.45) {$a^\vee$};
            \end{scope}
            \begin{scope}[yshift=1.75cm]
                \draw [knot] (0,0) -- (0, 0.25) to[out=90, in=90] (1, 0.25) -- (1,0);
                \draw[fill=white] (0.5, 0.55) circle (7.5pt);
                \node at (0.5, 0.55) {$b$};
            \end{scope}
            \draw[red, rounded corners, dashed] (-0.1, 0.125) rectangle (1.1, 1.75);
        }
    $};
    \node(E) at (13.25,-3) 
        {$\tikz[baseline={([yshift=-.5ex]current bounding box.center)}, scale=1]{
            \draw [knot] (0,0) -- (0, 0.25) to[out=90, in=90] (1, 0.25) -- (1,0);
            \draw[fill=white] (0.5, 0.55) circle (7.5pt);
            \node at (0.5, 0.55) {$b$};
        }
    $};
    \draw[->] (A) to[out=0,in=180] node[pos=.5,above,arrows=-]
        {$\lambda_{\psi(g)}$} 
        (B);
    \draw[->] (B) to[out=0,in=180] node[pos=.5,above,arrows=-]
        {$g \mid \mathbbm{1}_{a}$} 
        (C);  
    \draw[->] (C) to[out=0,in=180] node[pos=.5,above,arrows=-]
        {$\gamma_{W, \{-n,0\}}$} 
        (D);
    \draw[->] (D) to[out=-90,in=90] node[pos=.5,right,arrows=-]
        {$\mathbbm{1}_b \mid s_a$} 
        (E);
}
\]
where we set
\[
    \lambda_{\psi(g)} := \gamma_{W^N, w_g} \circ \lambda(N, w_g).
\]
Note that the last map is $\mathscr{G}$-graded by part 2 of Lemma \ref{lem:CPmincob}.

We compute $\psi(\phi(f))$ as the composition
\[
\tikz{
     \node(A) at (0,0) 
        {$v_f \circ \varphi_{W^N} \,
        \tikz[baseline={([yshift=-.5ex]current bounding box.center)}, scale=1]{
            \draw[knot] (0,0) -- (0, 0.25) to[out=90, in=90] (1, 0.25) -- (1,0);
            \draw[fill=white] (0.5, 0.55) circle (7.5pt);
            \node at (0.5, 0.55) {$a$};
        }
        $};
    \node(B) at (4.5,0) 
        {$ \varphi_{W^N} \circ v_f \,
        \tikz[baseline={([yshift=-.5ex]current bounding box.center)}, scale=1]{
            \draw[knot] (0,0) -- (0, 0.25) to[out=90, in=90] (1, 0.25) -- (1,0);
            \draw[fill=white] (0.5, 0.55) circle (7.5pt);
            \node at (0.5, 0.55) {$a$};
            \draw[red, rounded corners, dashed] (-0.1, 0.9) rectangle (1.1, 1.9);
        }
    $};    
    \node(C) at (9,0) 
        {$\varphi_{W^N} \circ v_f \circ \{-n,0\} \,
        \tikz[baseline={([yshift=-.5ex]current bounding box.center)}, scale=1]{
            \draw[knot] (0,0) -- (0, 0.25) to[out=90, in=90] (1, 0.25) -- (1,0);
            \draw[fill=white] (0.5, 0.55) circle (7.5pt);
            \node at (0.5, 0.55) {$a$};
            \begin{scope}[yshift=.75cm]
                \draw[knot] (0,1) -- (0, 0.75) to[out=-90, in=-90] (1, 0.75) -- (1,1);
                \draw[fill=white] (0.5, 0.45) circle (7.5pt);
                \node at (0.5, 0.45) {$a^\vee$};
            \end{scope}
            \begin{scope}[yshift=1.75cm]
                \draw [knot] (0,0) -- (0, 0.25) to[out=90, in=90] (1, 0.25) -- (1,0);
                \draw[fill=white] (0.5, 0.55) circle (7.5pt);
                \node at (0.5, 0.55) {$a$};
            \end{scope}
        }
    $};
    \node(D) at (12.5,-3) 
        {$ \varphi_{W^N} \circ \{-n,0\} \circ v_f \,
        \tikz[baseline={([yshift=-.5ex]current bounding box.center)}, scale=1]{
            \draw[knot] (0,0) -- (0, 0.25) to[out=90, in=90] (1, 0.25) -- (1,0);
            \draw[fill=white] (0.5, 0.55) circle (7.5pt);
            \node at (0.5, 0.55) {$a$};
            \begin{scope}[yshift=.75cm]
                \draw[knot] (0,1) -- (0, 0.75) to[out=-90, in=-90] (1, 0.75) -- (1,1);
                \draw[fill=white] (0.5, 0.45) circle (7.5pt);
                \node at (0.5, 0.45) {$a^\vee$};
            \end{scope}
            \begin{scope}[yshift=1.75cm]
                \draw [knot] (0,0) -- (0, 0.25) to[out=90, in=90] (1, 0.25) -- (1,0);
                \draw[fill=white] (0.5, 0.55) circle (7.5pt);
                \node at (0.5, 0.55) {$a$};
                \draw[red, rounded corners, dashed] (-0.1, 0) rectangle (1.1, 1);
            \end{scope}
        }
    $};
    \node(E) at (9,-6) 
        {$\varphi_{W^N} \circ \{-n,0\} \, 
        \tikz[baseline={([yshift=-.5ex]current bounding box.center)}, scale=1]{
            \draw[knot] (0,0) -- (0, 0.25) to[out=90, in=90] (1, 0.25) -- (1,0);
            \draw[fill=white] (0.5, 0.55) circle (7.5pt);
            \node at (0.5, 0.55) {$a$};
            \begin{scope}[yshift=.75cm]
                \draw[knot] (0,1) -- (0, 0.75) to[out=-90, in=-90] (1, 0.75) -- (1,1);
                \draw[fill=white] (0.5, 0.45) circle (7.5pt);
                \node at (0.5, 0.45) {$a^\vee$};
            \end{scope}
            \begin{scope}[yshift=1.75cm]
                \draw [knot] (0,0) -- (0, 0.25) to[out=90, in=90] (1, 0.25) -- (1,0);
                \draw[fill=white] (0.5, 0.55) circle (7.5pt);
                \node at (0.5, 0.55) {$b$};
            \end{scope}
        }
    $};
    \node(F) at (3.8,-6) 
        {$\varphi_{\left(W, (-\abs{b \mid a^\vee}, 0)\right)}
        \tikz[baseline={([yshift=-.5ex]current bounding box.center)}, scale=1]{
            \draw[knot] (0,0) -- (0, 0.25) to[out=90, in=90] (1, 0.25) -- (1,0);
            \draw[fill=white] (0.5, 0.55) circle (7.5pt);
            \node at (0.5, 0.55) {$a$};
            \begin{scope}[yshift=.75cm]
                \draw[knot] (0,1) -- (0, 0.75) to[out=-90, in=-90] (1, 0.75) -- (1,1);
                \draw[fill=white] (0.5, 0.45) circle (7.5pt);
                \node at (0.5, 0.45) {$a^\vee$};
            \end{scope}
            \begin{scope}[yshift=1.75cm]
                \draw [knot] (0,0) -- (0, 0.25) to[out=90, in=90] (1, 0.25) -- (1,0);
                \draw[fill=white] (0.5, 0.55) circle (7.5pt);
                \node at (0.5, 0.55) {$b$};
            \end{scope}
            \draw[red, rounded corners, dashed] (-0.1, 0.125) rectangle (1.1, 1.75);
        }
    $};
    \node(G) at (-0.2,-6) 
        {$
        \tikz[baseline={([yshift=-.5ex]current bounding box.center)}, scale=1]{
            \draw[knot] (0,0) -- (0, 0.25) to[out=90, in=90] (1, 0.25) -- (1,0);
            \draw[fill=white] (0.5, 0.55) circle (7.5pt);
            \node at (0.5, 0.55) {$b$};
        }
        $};
        \draw[->] (A) to[out=0,in=180] node[pos=.5,above,arrows=-]
        {$\lambda_{\psi(\phi(f))}$} 
        (B);
    \draw[->] (B) to[out=0,in=180] node[pos=.5,above,arrows=-]
        {$\eta_a$} 
        (C);  
    \draw[->] (C) to[out=0, in=90] node[pos=.5,right,arrows=-]
        {$\lambda(v_f, (-n,0))$} 
        (D);
    \draw[->] (D) to[out=-90, in=0] node[pos=.5,right,arrows=-]
        {$f \mid \mathbbm{1}_{a\mid a^\vee}$} 
        (E);
    \draw[->] (E) to[out=180,in=0] node[pos=.55,above,arrows=-]
        {$\gamma_{W, \{-n,0\}}$} 
        (F);
    \draw[->] (F) to[out=-180,in=0] node[pos=.5,above,arrows=-]
        {$\mathbbm{1}_b \mid s_a$} 
        (G);
}
\]
Sliding $f$ past the saddle introduces a change of chronology $\varphi_{H_1}$ which cancels with $\lambda(v_f, (-n,0))$ and $\lambda_{\psi(\phi(f))}$. A discerning eye notices that this change of chronology also kills the $\gamma_{W, \{-n,0\}}$ term, since the roles of $\varphi_{W^N}$ and $\{-n,0\}$ are interchanged during this change of chronology. After this slide, the composition looks like
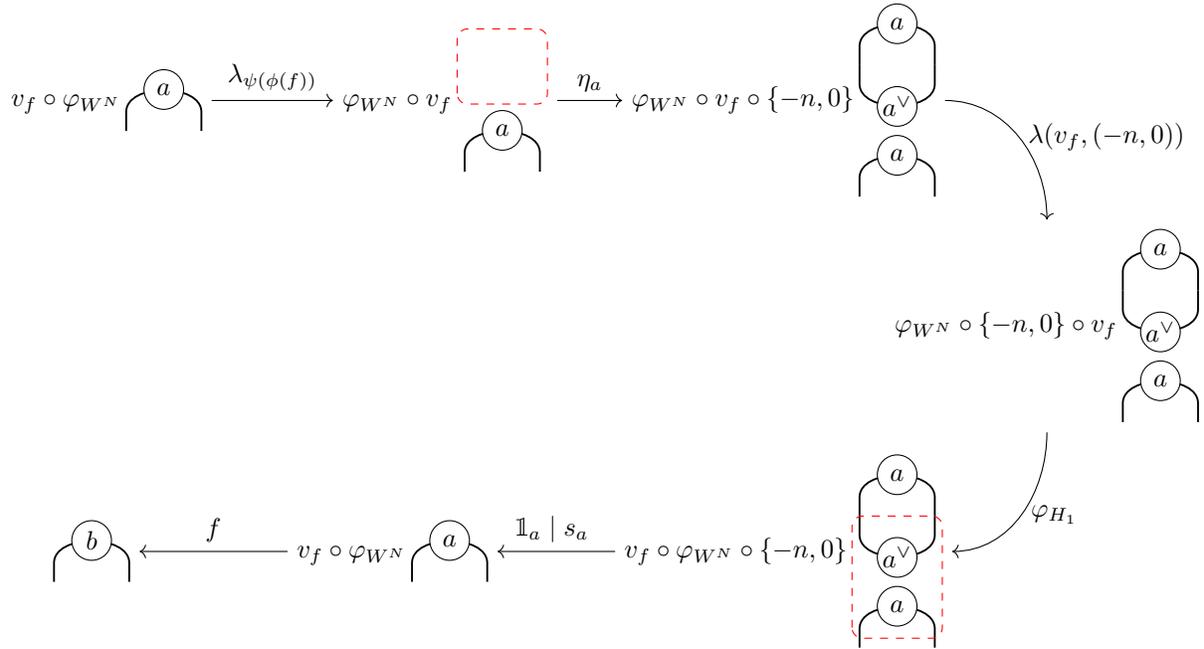
\[
\tikz{
     \node(A) at (0,0) 
        {$v_f \circ \varphi_{W^N} \,
        \tikz[baseline={([yshift=-.5ex]current bounding box.center)}, scale=1]{
            \draw[knot] (0,0) -- (0, 0.25) to[out=90, in=90] (1, 0.25) -- (1,0);
            \draw[fill=white] (0.5, 0.55) circle (7.5pt);
            \node at (0.5, 0.55) {$a$};
        }
        $};
    \node(B) at (4.5,0) 
        {$ \varphi_{W^N} \circ v_f \,
        \tikz[baseline={([yshift=-.5ex]current bounding box.center)}, scale=1]{
            \draw[knot] (0,0) -- (0, 0.25) to[out=90, in=90] (1, 0.25) -- (1,0);
            \draw[fill=white] (0.5, 0.55) circle (7.5pt);
            \node at (0.5, 0.55) {$a$};
            \draw[red, rounded corners, dashed] (-0.1, 0.9) rectangle (1.1, 1.9);
        }
    $};    
    \node(C) at (9,0) 
        {$\varphi_{W^N} \circ v_f \circ \{-n,0\} \,
        \tikz[baseline={([yshift=-.5ex]current bounding box.center)}, scale=1]{
            \draw[knot] (0,0) -- (0, 0.25) to[out=90, in=90] (1, 0.25) -- (1,0);
            \draw[fill=white] (0.5, 0.55) circle (7.5pt);
            \node at (0.5, 0.55) {$a$};
            \begin{scope}[yshift=.75cm]
                \draw[knot] (0,1) -- (0, 0.75) to[out=-90, in=-90] (1, 0.75) -- (1,1);
                \draw[fill=white] (0.5, 0.45) circle (7.5pt);
                \node at (0.5, 0.45) {$a^\vee$};
            \end{scope}
            \begin{scope}[yshift=1.75cm]
                \draw [knot] (0,0) -- (0, 0.25) to[out=90, in=90] (1, 0.25) -- (1,0);
                \draw[fill=white] (0.5, 0.55) circle (7.5pt);
                \node at (0.5, 0.55) {$a$};
            \end{scope}
        }
    $};
    \node(D) at (12.5,-3) 
        {$ \varphi_{W^N} \circ \{-n,0\} \circ v_f \,
        \tikz[baseline={([yshift=-.5ex]current bounding box.center)}, scale=1]{
            \draw[knot] (0,0) -- (0, 0.25) to[out=90, in=90] (1, 0.25) -- (1,0);
            \draw[fill=white] (0.5, 0.55) circle (7.5pt);
            \node at (0.5, 0.55) {$a$};
            \begin{scope}[yshift=.75cm]
                \draw[knot] (0,1) -- (0, 0.75) to[out=-90, in=-90] (1, 0.75) -- (1,1);
                \draw[fill=white] (0.5, 0.45) circle (7.5pt);
                \node at (0.5, 0.45) {$a^\vee$};
            \end{scope}
            \begin{scope}[yshift=1.75cm]
                \draw [knot] (0,0) -- (0, 0.25) to[out=90, in=90] (1, 0.25) -- (1,0);
                \draw[fill=white] (0.5, 0.55) circle (7.5pt);
                \node at (0.5, 0.55) {$a$};
            \end{scope}
        }
    $};
    \node(E) at (9,-6) 
        {$v_f \circ \varphi_{W^N} \circ \{-n,0\}\,
        \tikz[baseline={([yshift=-.5ex]current bounding box.center)}, scale=1]{
            \draw[knot] (0,0) -- (0, 0.25) to[out=90, in=90] (1, 0.25) -- (1,0);
            \draw[fill=white] (0.5, 0.55) circle (7.5pt);
            \node at (0.5, 0.55) {$a$};
            \begin{scope}[yshift=.75cm]
                \draw[knot] (0,1) -- (0, 0.75) to[out=-90, in=-90] (1, 0.75) -- (1,1);
                \draw[fill=white] (0.5, 0.45) circle (7.5pt);
                \node at (0.5, 0.45) {$a^\vee$};
            \end{scope}
            \begin{scope}[yshift=1.75cm]
                \draw [knot] (0,0) -- (0, 0.25) to[out=90, in=90] (1, 0.25) -- (1,0);
                \draw[fill=white] (0.5, 0.55) circle (7.5pt);
                \node at (0.5, 0.55) {$a$};
            \end{scope}
            \draw[red, rounded corners, dashed] (-0.1, 0.125) rectangle (1.1, 1.75);
        }
    $};
    \node(F) at (3.8,-6) 
        {$v_f \circ \varphi_{W^N} \,
        \tikz[baseline={([yshift=-.5ex]current bounding box.center)}, scale=1]{
            \draw[knot] (0,0) -- (0, 0.25) to[out=90, in=90] (1, 0.25) -- (1,0);
            \draw[fill=white] (0.5, 0.55) circle (7.5pt);
            \node at (0.5, 0.55) {$a$};
        }
    $};
    \node(G) at (-0.2,-6) 
        {$
        \tikz[baseline={([yshift=-.5ex]current bounding box.center)}, scale=1]{
            \draw[knot] (0,0) -- (0, 0.25) to[out=90, in=90] (1, 0.25) -- (1,0);
            \draw[fill=white] (0.5, 0.55) circle (7.5pt);
            \node at (0.5, 0.55) {$b$};
        }
        $};
        \draw[->] (A) to[out=0,in=180] node[pos=.5,above,arrows=-]
        {$\lambda_{\psi(\phi(f))}$} 
        (B);
    \draw[->] (B) to[out=0,in=180] node[pos=.5,above,arrows=-]
        {$\eta_a$} 
        (C);  
    \draw[->] (C) to[out=0, in=90] node[pos=.5,right,arrows=-]
        {$\lambda(v_f, (-n,0))$} 
        (D);
    \draw[->] (D) to[out=-90, in=0] node[pos=.5,right,arrows=-]
        {$\varphi_{H_1}$} 
        (E);
    \draw[->] (E) to[out=180,in=0] node[pos=.55,above,arrows=-]
        {$\mathbbm{1}_a\mid s_a$} 
        (F);
    \draw[->] (F) to[out=-180,in=0] node[pos=.5,above,arrows=-]
        {$f$} 
        (G);
}
\]
Notice that $\mathbbm{1}_a \mid s_a$ consists of $n$ merges, so the penultimate arrow makes sense. Then, 1 of Lemma \ref{lem:CPmincob} gives us that $\psi(\phi(f)) = f$.

On the other hand, $\phi(\psi(g))$ is rather easy to compute; the reader is invited to verify that this composition simplifies to
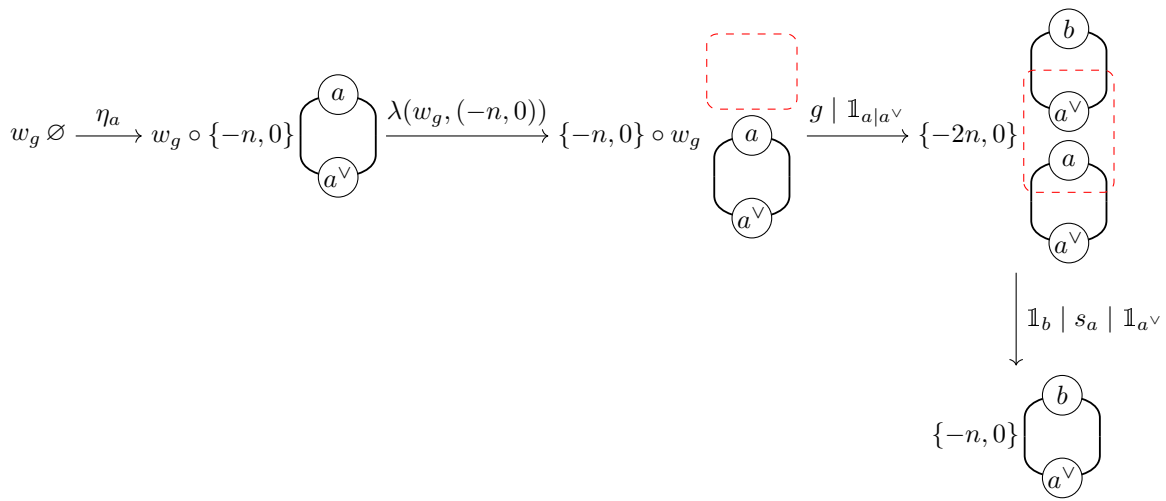
\[
\tikz{
     \node(A) at (0,0) 
        {$w_g  \, 
        \emptyset
        $};
    \node(B) at (3,0) 
        {$w_g \circ \{-n,0\} \,
        \tikz[baseline={([yshift=-.5ex]current bounding box.center)}, scale=1]{
            \draw[knot] (0,1) -- (0, 0.75) to[out=-90, in=-90] (1, 0.75) -- (1,1);
            \draw[fill=white] (0.5, 0.45) circle (7.5pt);
            \node at (0.5, 0.45) {$a^\vee$};
        \begin{scope}[yshift=1cm]
            \draw [knot] (0,0) -- (0, 0.25) to[out=90, in=90] (1, 0.25) -- (1,0);
            \draw[fill=white] (0.5, 0.55) circle (7.5pt);
            \node at (0.5, 0.55) {$a$};
        \end{scope}
        }
        $};
    \node(C) at (8.5,0) 
        {$\{-n,0\} \circ w_g \,
        \tikz[baseline={([yshift=-.5ex]current bounding box.center)}, scale=1]{
            \draw[knot] (0,1) -- (0, 0.75) to[out=-90, in=-90] (1, 0.75) -- (1,1);
            \draw[fill=white] (0.5, 0.45) circle (7.5pt);
            \node at (0.5, 0.45) {$a^\vee$};
        \begin{scope}[yshift=1cm]
            \draw [knot] (0,0) -- (0, 0.25) to[out=90, in=90] (1, 0.25) -- (1,0);
            \draw[fill=white] (0.5, 0.55) circle (7.5pt);
            \node at (0.5, 0.55) {$a$};
            \draw[red, rounded corners, dashed] (-0.1, 0.9) rectangle (1.1, 1.9);
        \end{scope}
        }
        $};
    \node(D) at (13,0) 
        {$\{-2n,0\} \,
        \tikz[baseline={([yshift=-.5ex]current bounding box.center)}, scale=1]{
            \draw[knot] (0,1) -- (0, 0.75) to[out=-90, in=-90] (1, 0.75) -- (1,1);
            \draw[fill=white] (0.5, 0.45) circle (7.5pt);
            \node at (0.5, 0.45) {$a^\vee$};
        \begin{scope}[yshift=1cm]
            \draw [knot] (0,0) -- (0, 0.25) to[out=90, in=90] (1, 0.25) -- (1,0);
            \draw[fill=white] (0.5, 0.55) circle (7.5pt);
            \node at (0.5, 0.55) {$a$};
            \begin{scope}[yshift=.75cm]
                \draw[knot] (0,1) -- (0, 0.75) to[out=-90, in=-90] (1, 0.75) -- (1,1);
                \draw[fill=white] (0.5, 0.45) circle (7.5pt);
                \node at (0.5, 0.45) {$a^\vee$};
            \end{scope}
            \begin{scope}[yshift=1.75cm]
                \draw [knot] (0,0) -- (0, 0.25) to[out=90, in=90] (1, 0.25) -- (1,0);
                \draw[fill=white] (0.5, 0.55) circle (7.5pt);
                \node at (0.5, 0.55) {$b$};
            \end{scope}
            \draw[red, rounded corners, dashed] (-0.1, 0.125) rectangle (1.1, 1.75);
        \end{scope}
        }
        $};
    \node(E) at (13,-4) 
        {$\{-n, 0\}\, 
        \tikz[baseline={([yshift=-.5ex]current bounding box.center)}, scale=1]{
        \draw[knot] (0,1) -- (0, 0.75) to[out=-90, in=-90] (1, 0.75) -- (1,1);
        \draw[fill=white] (0.5, 0.45) circle (7.5pt);
        \node at (0.5, 0.45) {$a^\vee$};
    \begin{scope}[yshift=1cm]
        \draw [knot] (0,0) -- (0, 0.25) to[out=90, in=90] (1, 0.25) -- (1,0);
        \draw[fill=white] (0.5, 0.55) circle (7.5pt);
        \node at (0.5, 0.55) {$b$};
    \end{scope}
    }
        $};
        \draw[->] (A) to[out=0,in=180] node[pos=.5,above,arrows=-]
        {$\eta_a$} 
        (B);
    \draw[->] (B) to[out=0,in=180] node[pos=.5,above,arrows=-]
        {$\lambda(w_g, (-n,0))$} 
        (C);  
    \draw[->] (C) to[out=0,in=180] node[pos=.5,above,arrows=-]
        {$g \mid \mathbbm{1}_{a\mid a^\vee}$} 
        (D);
    \draw[->] (D) to[out=-90,in=90] node[pos=.5,right,arrows=-]
        {$\mathbbm{1}_b \mid s_a \mid \mathbbm{1}_{a^\vee}$} 
        (E);
}
\]
Then, pushing $g$ before the birth introduces a change of chronology $\varphi_{H_2}$ equal to $\lambda((-n,0), w_g)$. This is inverse to $\lambda(w_g, (-n,0))$, so that the new composition is
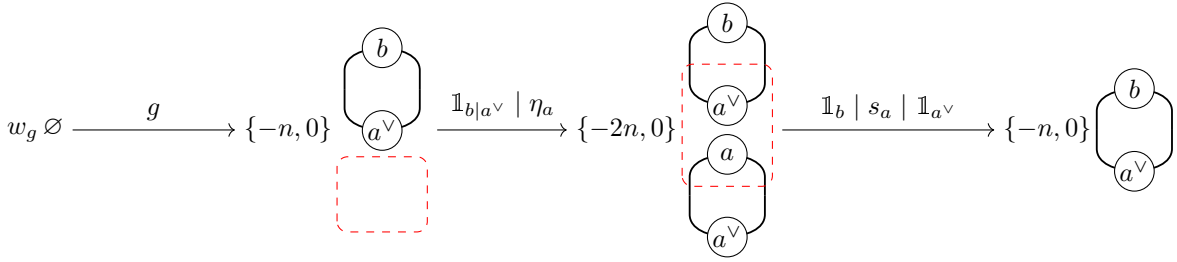
\[
\tikz{
     \node(A) at (0,0) 
        {$w_g \, \emptyset$};
    \node(B) at (4,0) 
        {$\{-n, 0\} \, 
        \tikz[baseline={([yshift=-.5ex]current bounding box.center)}, scale=1]{
        \draw[knot] (0,1) -- (0, 0.75) to[out=-90, in=-90] (1, 0.75) -- (1,1);
        \draw[fill=white] (0.5, 0.45) circle (7.5pt);
        \node at (0.5, 0.45) {$a^\vee$};
    \begin{scope}[yshift=1cm]
        \draw [knot] (0,0) -- (0, 0.25) to[out=90, in=90] (1, 0.25) -- (1,0);
        \draw[fill=white] (0.5, 0.55) circle (7.5pt);
        \node at (0.5, 0.55) {$b$};
    \end{scope}
    \begin{scope}[yshift=-1.8cm]
    \draw[red, rounded corners, dashed] (-0.1, 0.9) rectangle (1.1, 1.9);
    \end{scope}
    }
        $};
    \node(C) at (8.5,0) 
        {$\{-2n,0\} \,
        \tikz[baseline={([yshift=-.5ex]current bounding box.center)}, scale=1]{
            \draw[knot] (0,1) -- (0, 0.75) to[out=-90, in=-90] (1, 0.75) -- (1,1);
            \draw[fill=white] (0.5, 0.45) circle (7.5pt);
            \node at (0.5, 0.45) {$a^\vee$};
        \begin{scope}[yshift=1cm]
            \draw [knot] (0,0) -- (0, 0.25) to[out=90, in=90] (1, 0.25) -- (1,0);
            \draw[fill=white] (0.5, 0.55) circle (7.5pt);
            \node at (0.5, 0.55) {$a$};
            \begin{scope}[yshift=.75cm]
                \draw[knot] (0,1) -- (0, 0.75) to[out=-90, in=-90] (1, 0.75) -- (1,1);
                \draw[fill=white] (0.5, 0.45) circle (7.5pt);
                \node at (0.5, 0.45) {$a^\vee$};
            \end{scope}
            \begin{scope}[yshift=1.75cm]
                \draw [knot] (0,0) -- (0, 0.25) to[out=90, in=90] (1, 0.25) -- (1,0);
                \draw[fill=white] (0.5, 0.55) circle (7.5pt);
                \node at (0.5, 0.55) {$b$};
            \end{scope}
            \draw[red, rounded corners, dashed] (-0.1, 0.125) rectangle (1.1, 1.75);
        \end{scope}
        }
        $};
    \node(D) at (14,0) 
        {$\{-n, 0\}\, 
        \tikz[baseline={([yshift=-.5ex]current bounding box.center)}, scale=1]{
        \draw[knot] (0,1) -- (0, 0.75) to[out=-90, in=-90] (1, 0.75) -- (1,1);
        \draw[fill=white] (0.5, 0.45) circle (7.5pt);
        \node at (0.5, 0.45) {$a^\vee$};
    \begin{scope}[yshift=1cm]
        \draw [knot] (0,0) -- (0, 0.25) to[out=90, in=90] (1, 0.25) -- (1,0);
        \draw[fill=white] (0.5, 0.55) circle (7.5pt);
        \node at (0.5, 0.55) {$b$};
    \end{scope}
    }
        $};
    \draw[->] (A) to[out=0,in=180] node[pos=.5,above,arrows=-]
        {$g$} 
        (B);
    \draw[->] (B) to[out=0,in=180] node[pos=.5,above,arrows=-]
        {$\mathbbm{1}_{b \mid a^\vee} \mid \eta_a$} 
        (C);  
    \draw[->] (C) to[out=0,in=180] node[pos=.5,above,arrows=-]
        {$\mathbbm{1}_b \mid s_a \mid \mathbbm{1}_{a^\vee}$} 
        (D);  
}
\]
which simplifies to $g$ by Lemma \ref{lem:CPmincob}. This concludes the proof.
\end{proof}

\begin{remark}
Since a minimal cobordism $a\to b$ consists of $(n - \abs{b\mid a^\vee})$-many saddles, 
\[
\deg_q\left(
\varphi_{\left(W, (n - \abs{\widehat{b} \mid \widehat{a}}, 0)\right)}
\right) = 0
\]
and we obtain a generalization of Proposition 4.8 in \cite{hogancamp2012morphismscategorifiedspinnetworks}.
\end{remark}

\begin{proof}[Proof of Theorem \ref{thm:dualclosing}]
Recall that $\normalfont{\textsc{Hom}}_n (A, B)$, for $A$ and $B$ $\mathscr{G}$-graded dg-$H^n$-modules, is the chain complex of bihomogeneous (that is, homogeneous in homological degree and purely homogeneous in $\widetilde{\mathscr{I}}$-degree) maps $f$ of arbitrary $(\mathbb{Z} \times \widetilde{\mathscr{I}})$-degree. So, we can view $\normalfont{\textsc{Hom}}_n$-complexes as bigraded abelian groups
\[
\normalfont{\textsc{Hom}}_n(A, B)^k_{(W, v)}
\cong
\prod_{\ell \in \mathbb{Z}} \mathrm{Hom}_n \left(\varphi_{(W, v)} A^\ell, B^{\ell + k}\right).
\]
However, notice that for each $\ell$, $k$, and $(W, v)$, Proposition \ref{ShiftingTubes} says that $\varphi_{(W, v)} \cong \varphi_{(W_\ell^k, v')}$ for $W_\ell^k$ a minimal cobordism $A^\ell \to B^{\ell + k}$ and $v' = v + \tau_W(-1,-1)$. This means that $\mathrm{Hom}_n(\varphi_{(W, v)} A^\ell, B^{\ell + k})$ is canonically isomorphic to $\mathrm{Hom}_n\left(\varphi_{(W_\ell^k, v')}A^\ell, B^{\ell + k}\right)$. Set $v_\ell^k := (n - \abs{A^\ell \mid (B^{\ell + k})^\vee}, 0)$; we conclude that
\[
\mathrm{Hom}_n\left(\varphi_{(W, v)} A^\ell, B^{\ell + k}\right)
\cong
\mathrm{Hom}_n \left(\varphi_{(W_\ell^k, v_\ell^k)} A^\ell \, \{v' - v_\ell^k\}, B^{\ell + k}\right).
\]
Thus, in the $\mathscr{G}$-graded case, we can absorb the first coordinate of the $\widetilde{\mathscr{I}}$-grading into the homological degree and view $\normalfont{\textsc{Hom}}_n(A, B)$ as bigraded by $\mathbb{Z} \times \mathbb{Z}^2$. Then
\begin{align*}
    \normalfont{\textsc{Hom}}_n\left(A, B\right)^k & \cong \prod_{\ell \in \mathbb{Z}} \mathrm{Hom}_n\left( \varphi_{(W_\ell^k, v_\ell^k)} A^\ell, B^{\ell + k}\right)
    \\ &\cong \prod_{\ell \in \mathbb{Z}} \mathrm{Hom}_0\left(\emptyset, B^{\ell + k} \mid (A^\ell)^\vee \, \{-n, 0\}\right)
    \\ &\cong \normalfont{\textsc{Hom}}_0 \left(\emptyset, B \mid A^\vee \, \{-n, 0\}\right)
\end{align*}
where the second isomorphism follows from Proposition \ref{prop:CPadjunction2}.

This proves the isomorphism on the level of bigraded abelian groups. The rest of the statement follows from the argument provided in the proof of Theorem 4.12 in \cite{hogancamp2012morphismscategorifiedspinnetworks}. We will not review the proof here, but for the argument to apply we must show that
\[
(g \mid \mathbbm{1}_{a^\vee}) \circ_\mathscr{G} \phi(f) = \phi(g\circ_\mathscr{G} f) = (\mathbbm{1}_c \mid f^\vee) \circ_\mathscr{G} \phi(g)
\]
where $f\in \mathrm{Hom}_n(\varphi_{(W_1, N_1)}a, b)$, $g\in \mathrm{Hom}_n(\varphi_{(W_2, N_2)}b, c)$, and  $\phi: \mathrm{Hom}_n(\varphi_{W^N}a, c) \to \mathrm{Hom}_0(\emptyset, c \mid a^\vee \, \{-n, 0\})$ is the isomorphism from the proof of Proposition \ref{prop:CPadjunction2}. Here, $W_1:a\to b$ and $W_2: b \to c$ are minimal cobordisms, and $N_1 = (n - \abs{b \mid a^\vee}, 0)$ and $N_2 = (n - \abs{c \mid b^\vee}, 0)$, thus $g\circ_\mathscr{G} f \in \mathrm{Hom}_n\left(\varphi_{(W_2\circ W_1, N_1 + N_2)}a, c\right)$. The equality on the left-hand side is immediate. We will content ourselves by proving the right-hand side.

To start, we claim that 
\[
(f \mid \mathbbm{1}_{a^\vee}) \circ \eta_a = (\mathbbm{1}_b \mid f^\vee) \circ \eta_a.
\]
Notice that the claim holds trivially when $f$ is a dot. When $f$ is a saddle, both $f$ and $f^\vee$ are necessarily merge, and their $\widetilde{\mathscr{I}}$-degree is $\mathrm{Id}$. Thus, in this case, isotopy invariance implies the equality. Indeed, for any $f \in \mathrm{Hom}_n(\varphi_{(W_1, N_1)} a, b)$, the $\widetilde{\mathscr{I}}$-degree of $f \mid \mathbbm{1}_{a^\vee}$ is supported entirely in the $\mathbb{Z}\times \mathbb{Z}$-coordinate; the same is true for $\mathbbm{1}_{b^\vee} \mid f^\vee$. We denote this degree by $v_f$ and, in this case, we have that $v_f = v_{f^\vee}$. To conclude the proof of the claim, we have to show the equality holds for compositions $g\circ_\mathscr{G} f$, for $f$ and $g$ as above. First, notice that 
\[
(g \circ_\mathscr{G} f)  \mid \mathbbm{1}_{a^\vee}  = (g\mid \mathbbm{1}_{a^\vee}) \circ_\mathscr{G} (f \mid \mathbbm{1}_{a^\vee})
\]
by Proposition \ref{prop:vertical&horizontalcomp} (here, $\Xi = 1$ since $\mathbbm{1}_{a^\vee}$ is two of the four inputted maps). On the other hand,
\[
(g\mid \mathbbm{1}_{a^\vee}) \circ_\mathscr{G} (f \mid \mathbbm{1}_{a^\vee}) = (g\mid \mathbbm{1}_{a^\vee}) \circ (f \mid \mathbbm{1}_{a^\vee})
\]
since each map in the composite has trivial $\widetilde{\mathscr{I}}$-degree. So, we have
\begin{align*}
    (g\mid \mathbbm{1}_{a^\vee}) \circ (f \mid \mathbbm{1}_{a^\vee}) \circ \eta_a &= 
    (g\mid \mathbbm{1}_{a^\vee}) \circ (\mathbbm{1}_{b} \mid f^\vee) \circ \eta_b \\
    &= (\mathbbm{1}_c \mid f^\vee) \circ (g \mid \mathbbm{1}_{b^\vee}) \circ \lambda(w_g, v_f) \circ \eta_b \\
    &= (\mathbbm{1}_c \mid f^\vee) \circ (\mathbbm{1}_c \mid g^\vee) \circ \lambda(w_g, v_f) \circ \eta_c.
\end{align*}
The first and last equalities are by assumption. The second equality follows from applying a change of chronology. Notice that $\lambda(w_g, v_f)$ is, in this setting, equal to the value $\varphi_{H_{f,g}}$. Then, again applying Proposition \ref{prop:vertical&horizontalcomp}, we conclude that
\begin{align*}
((g \circ_\mathscr{G} f) \mid \mathbbm{1}_{a^\vee}) \circ \eta_a &= (\mathbbm{1}_{c} \mid (f^\vee \circ_\mathscr{G} g^\vee)) \circ \varphi_{H_{f,g}} \circ \eta_c \\ & = (\mathbbm{1}_{c} \mid (g\circ_\mathscr{G} f)^\vee) \circ \eta_c.
\end{align*}
We leave it to the reader to verify that one application of this claim implies that 
\[
\phi(g\circ_\mathscr{G} f) = (\mathbbm{1}_c \mid f^\vee) \circ_\mathscr{G} \phi(g)
\]
concluding our proof.
\end{proof}

\subsection{Definition and properties of unified projectors}
\label{s:unifiedprops}

Recall that the \textit{through-degree} of a Temperley-Lieb diagram $\delta$, denoted $\tau(\delta)$, is the number of strands with endpoints on opposite ends of the disk. We say that $A \in \mathrm{Chom}(n)^\mathscr{G}$ has \textit{through-degree} less than $k$ if $A$ is homotopy equivalent to a colimit of $\mathscr{G}$-graded dg-modules $\mathcal{F}(\delta)$ for Temperley-Lieb diagramas $\delta$ with $\tau(\delta) < k$. In this case, we also write $\tau(A) < k$. Since the tensor product commutes with colimits, we have that $\tau(A \otimes B) \le \min\{\tau(A), \tau(B)\}$.

\begin{definition}
We say that $A \in \mathrm{Chom}(n)^\mathscr{G}$ \textit{kills turnbacks from above} if, for each $B\in \mathrm{Chom}(n)^\mathscr{G}$ with $\tau(B) < n$, we have $B \otimes A \simeq *$. Similarly, $A \in \mathrm{Chom}(n)$ \textit{kills turnbacks from below} if, for each $B$ with $\tau(B) < n$, $A \otimes B \simeq *$.
\end{definition}

Since all Temperley-Lieb diagrams with through-degree less than $k$ can be built by stacking various generators $e_i$ of the Temperley-Lieb algebra, we have the following (stated without proof).

\begin{proposition}
Let $e_i$ denote a standard generator of the Temperley-Lieb algebra. Then any object $A$ of $\mathrm{Chom}(n)^\mathscr{G}$ kills turnbacks from above (resp. below) if and only if $\mathcal{F}(e_i) \otimes A \simeq *$ (resp. $A \otimes \mathcal{F}(e_i) \simeq *$.
\end{proposition}

\begin{definition}
\label{def:unifiedprojector}
A \textit{unified Cooper-Krushkal projector} (or simply \textit{unified projector}) is a pair $(P_n, \iota)$ consisting of an object $P_n \in \mathrm{Chom}(n)^\mathscr{G}$ and a morphism $\iota: \mathcal{I}_n \to P_n$, called the \textit{unit} of the projector, so that
\begin{enumerate}
    \item[(CK1)] $\mathrm{Cone}(\iota)$ has through-degree less than $n$, and
    \item[(CK2)] the $\mathscr{G}$-graded dg-module $P_n$ kills turnbacks (from above and below).
\end{enumerate}
\end{definition}

\begin{lemma}
\label{lem:pullbackHE}
If $(P_n, \iota)$ is a unified projector, there is a homotopy equivalence
\[
\normalfont{\textsc{Hom}}_n(P_n, P_n) \to  \normalfont{\textsc{Hom}}_n(\mathcal{I}_n, P_n)
\]
induced by $\iota$.
\end{lemma}

\begin{proof}
Specifically, we will show that the pullback $\iota^*: \normalfont{\textsc{Hom}}_n(P_n, P_n) \to  \normalfont{\textsc{Hom}}_n(\mathcal{I}_n, P_n)$ is a homotopy equivalence. It suffices to show that $\mathrm{Cone}(\iota^*)$ is contractible. We compute
\begin{align*}
    \mathrm{Cone}(\iota^*) & \simeq \normalfont{\textsc{Hom}}_n (\mathrm{Cone}(\iota), P_n) & \\
    & \simeq \normalfont{\textsc{Hom}}_n(\mathrm{colim}(\mathcal{F}(\delta)_i), P_n) & \text{(CK1), $\tau(\delta) < n$ for all $i$} \\
    & \simeq \varprojlim (\normalfont{\textsc{Hom}}_n(\mathcal{F}(\delta)_i, P_n)) & \\
    & \simeq \varprojlim (\normalfont{\textsc{Hom}}_n(\mathcal{I}_n, P_n \otimes \mathcal{F}(\delta^\vee)_i)) & \text{Corollary \ref{cor:dualswap}} \\
    & \simeq \varprojlim (\normalfont{\textsc{Hom}}_n(\mathcal{I}_n, *)) & \text{(CK2)} \\
    & \simeq *
\end{align*}
as desired.
\end{proof}

\begin{proposition}[Properties of unified projectors]
\label{prop:unifiedprojectoruni}
Suppose $(P_n, \iota)$ and $(P_n', \iota')$ are two unified projectors of $\mathrm{Chom}(n)^\mathscr{G}$.
\begin{enumerate}
    \item (Uniqueness) $P_n \simeq P_n' \otimes P_n \simeq P_n'$, and there is a homotopy equivalence $h: P_n \to P_n'$ satisfying $h\circ \iota \simeq \iota'$.
    \item (Idempotence) $(P_n\otimes P_n, \iota \otimes \iota)$ is a projector; thus, by uniqueness, $P_n \otimes P_n \simeq P_n$.
    \item (Generalized absorbtion) More generally, for $\ell \le n$
    \[
        P_n \otimes \left(P_\ell \sqcup \mathcal{I}_{n-\ell}\right) \simeq P_n \simeq \left(P_\ell \sqcup \mathcal{I}_{n-\ell}\right) \otimes P_n.
    \]
\end{enumerate}
\end{proposition}

\begin{proof}
Consider the following $\mathscr{G}$-graded commutative diagram.
\[
\begin{tikzcd}
    &  P_n  \arrow[rr] & & \mathcal{I}_n \otimes P_n \arrow[dr, "\iota' \otimes \mathrm{id}_{P_n}"] & \\
    \mathcal{I}_n \arrow[ur, "\iota"] \arrow[rr] \arrow[dr, "\iota'"'] & & \mathcal{I}_n \otimes \mathcal{I}_n \arrow[ur, "\mathrm{id}_{\mathcal{I}_n} \otimes \iota"] \arrow[rr, "\iota' \otimes \iota"] \arrow[dr, "\iota' \otimes \mathrm{id}_{\mathcal{I}_n}"'] & & P_n' \otimes P_n \\
    & P_n' \arrow[rr] & & P_n' \otimes \mathcal{I}_n \arrow[ur, "\mathrm{id}_{P_n'} \otimes \iota"'] &
\end{tikzcd}
\]
The unmarked arrows are isomorphisms coming from multigluing (or, if one likes, the probable the monoidal structure of $\mathrm{Chom}(n)^\mathscr{G}$). Since this diagram is $\mathscr{G}$-graded commutative, it commutes up to homotopy, which is all we need going forward.

For the proof of uniqueness, notice that $\iota' \otimes \mathrm{id}_{P_n}$ is a homotopy equivalence, as $\mathrm{Cone}(\iota' \otimes \mathrm{id}_{P_n}) \simeq \mathrm{Cone}(\iota') \otimes P_n \simeq *$, using (CK1) and (CK2). By the same reasoning, $\mathrm{id}_{P_n'}\otimes \iota$ is a homotopy equivalence, thus 
\[
P_n \simeq P_n'\otimes P_n \simeq P_n'.
\]
Then, since both these maps are homotopy equivalences, choosing a homotopy inverse for, say, $\mathrm{id}_{P_n'}\otimes \iota$ induces a (class of) homotopy equivalence(s) $h: P_n \to P_n'$ satisfying $h \circ \iota \simeq \iota'$. To see that $h$ is unique up to homotopy, suppose $h_1, h_2$ are two homotopy equivalences satisfying, for $i=1,2$,  $h_i \circ \iota \simeq \iota'$, and that $\overline{h}_2$ is a homotopy inverse for $h_2$. Then $(\iota - \overline{h}_2 \circ h_1 \circ \iota) = (\mathrm{id}_{P_n} - \overline{h}_2 \circ h_1) \circ \iota \in \normalfont{\textsc{Hom}}_n(\mathcal{I}_n, P_n)$ is nullhomotopic, so Lemma \ref{lem:pullbackHE} implies that $\mathrm{id}_{P_n} - \overline{h}_2 \circ h_1$ is as well; thus $h_1\simeq h_2$.

For idempotence, replace $P_n'$ in the diagram with $P_n$ everywhere. Then we have that $P_n \otimes P_n \simeq P_n$. More generally, that $P_n \otimes P_n$ kills turnbacks is clear by the monoidal structure of $\mathrm{Chom}(n)$. Then, since $\iota \otimes \mathrm{id}_{P_n}$ is a homotopy equivalence, the homotopy commutativity of the diagram implies that $\mathrm{Cone}(\iota \otimes \iota) \simeq \mathrm{Cone}(\mathrm{id}_{\mathcal{I}_n} \otimes \iota) \simeq *$.

More generally, for $\ell < n$, $P_\ell$ comes equipped with unit $\iota_\ell: \mathcal{I}_\ell \to P_\ell$. Then, it is clear that 
\[
\mathrm{id}_{P_n} \otimes (\iota_\ell \sqcup \mathrm{id}_{\mathcal{I}_{n-\ell}}): P_n \otimes (P_\ell \sqcup \mathcal{I}_{n-\ell}) \longrightarrow P_n\otimes \mathcal{I}_n \simeq P_n
\]
is a homotopy equivalence (its cone is contractible by (CK2)). The other homotopy equivalence is analogous.
\end{proof}

\begin{remark}
We can define projectors for the category $\mathrm{Chom}(n)^q$ similarly. Notice that projectors of $\mathrm{Chom}(n)^\mathscr{G}$ descend to projectors of $\mathrm{Chom}(n)^q$; in addition, given any $(W,v)\in\mathscr{I}$ with $\deg_q(\varphi_{W^v}) = 0$, $\varphi_{W^v}P_n$ defines a projector of $\mathrm{Chom}(n)^q$.
\end{remark}

In future work, we hope to find particular elements $U_n \in \normalfont{\textsc{Hom}}_n(P_n, P_n)$ coming from an action on $P_n$, as in \cite{https://doi.org/10.48550/arxiv.1405.2574}. Fundamental to this study is the homotopy equivalence between the endomorphism complex of $P_n$ and (a shift of) the closure of $P_n$. We point out that Theorem \ref{thm:adjunction1} and Lemma \ref{lem:pullbackHE} imply a generalization of this result in the unified setting; we state it for the $q$-graded category. For $A, B$ $\mathscr{G}$-graded dg $H^n$-modules let $\normalfont{\textsc{Hom}}_n^q(A,B)$ denote the HOM-complex $\normalfont{\textsc{Hom}}_n(A,B)$ obtained by collapsing $\mathscr{G}$-grading.

\begin{corollary}
\label{cor:endothing}
If $P_n$ is a unified projector, it descends to a projector in $\mathrm{Chom}(n)^q$. We have that
\[
\normalfont{\textsc{Hom}}_n^q(P_n, P_n) \cong q^{-n} \mathrm{Tr}^n (P_n).
\]
\end{corollary}

\begin{proof}
Apply Lemma \ref{lem:pullbackHE} and then apply Theorem \ref{thm:adjunction1} $n$-times.
\end{proof}

\subsection{Explicit computations for the 2-stranded projector}
\label{ss:2strandedproj}

Finally, our previous work allows us to mimic \cite{https://doi.org/10.48550/arxiv.1005.5117} in the $\mathscr{G}$-graded (that is, unified) setting. Consider the complex we will call $P_2$, which has the form
\[
\begin{tikzcd}[sep=scriptsize]
\cdots \arrow[r, "C_{-4}"] & 
\varphi_{\left(\tikz[baseline=-4.33ex, scale=.3]
{
	\draw[dotted] (3,-2) circle(0.707);
	\draw[knot] (2.5,-1.5) .. controls (2.75,-1.75) and (3.25,-1.75) .. (3.5,-1.5);
	\draw[knot] (2.5,-2.5) .. controls  (2.75,-2.25) and (3.25,-2.25) .. (3.5,-2.5);
        \draw[red,thick] (3,-1.7) -- (3,-2.3);
}\,,~ (-2,-2)\right)} \tikz[baseline=-6.5ex, scale=.45]
{
	\draw[dotted] (3,-2) circle(0.707);
	\draw[knot] (2.5,-1.5) .. controls (2.75,-1.75) and (3.25,-1.75) .. (3.5,-1.5);
	\draw[knot] (2.5,-2.5) .. controls  (2.75,-2.25) and (3.25,-2.25) .. (3.5,-2.5);
} \arrow[r, "C_{-3}"] &
\varphi_{\left(\tikz[baseline=-4.33ex, scale=.3]
{
	\draw[dotted] (3,-2) circle(0.707);
	\draw[knot] (2.5,-1.5) .. controls (2.75,-1.75) and (3.25,-1.75) .. (3.5,-1.5);
	\draw[knot] (2.5,-2.5) .. controls  (2.75,-2.25) and (3.25,-2.25) .. (3.5,-2.5);
        \draw[red,thick] (3,-1.7) -- (3,-2.3);
}\,,~ (-1,-1)\right)} \tikz[baseline=-6.5ex, scale=.45]
{
	\draw[dotted] (3,-2) circle(0.707);
	\draw[knot] (2.5,-1.5) .. controls (2.75,-1.75) and (3.25,-1.75) .. (3.5,-1.5);
	\draw[knot] (2.5,-2.5) .. controls  (2.75,-2.25) and (3.25,-2.25) .. (3.5,-2.5);
} \arrow[r, "C_{-2}"] &
\varphi_{\tikz[baseline=-4.33ex, scale=.3]
{
	\draw[dotted] (3,-2) circle(0.707);
	\draw[knot] (2.5,-1.5) .. controls (2.75,-1.75) and (3.25,-1.75) .. (3.5,-1.5);
	\draw[knot] (2.5,-2.5) .. controls  (2.75,-2.25) and (3.25,-2.25) .. (3.5,-2.5);
        \draw[red,thick] (3,-1.7) -- (3,-2.3);
}} ~ \tikz[baseline=-6.5ex, scale=.45]
{
	\draw[dotted] (3,-2) circle(0.707);
	\draw[knot] (2.5,-1.5) .. controls (2.75,-1.75) and (3.25,-1.75) .. (3.5,-1.5);
	\draw[knot] (2.5,-2.5) .. controls  (2.75,-2.25) and (3.25,-2.25) .. (3.5,-2.5);
} \arrow[r, "C_{-1}"] &
\tikz[baseline=8.25ex, scale=.45]
{
    \begin{scope}[rotate=90]
	\draw[dotted] (3,-2) circle(0.707);
	\draw[knot] (2.5,-1.5) .. controls (2.75,-1.75) and (3.25,-1.75) .. (3.5,-1.5);
	\draw[knot] (2.5,-2.5) .. controls  (2.75,-2.25) and (3.25,-2.25) .. (3.5,-2.5);
    \end{scope}
}
\end{tikzcd}
\]
where
\[
C_i = \begin{cases} \tikz[baseline=-6.5ex, scale=.45]
{
	\draw[dotted] (3,-2) circle(0.707);
	\draw[knot] (2.5,-1.5) .. controls (2.75,-1.75) and (3.25,-1.75) .. (3.5,-1.5);
	\draw[knot] (2.5,-2.5) .. controls  (2.75,-2.25) and (3.25,-2.25) .. (3.5,-2.5);
        \draw[red,thick] (3,-1.7) -- (3,-2.3);
} & i=-1 \\[10pt] 
\tikz[baseline=-6.5ex, scale=.45]
{
	\draw[dotted] (3,-2) circle(0.707);
	\draw[knot] (2.5,-1.5) .. controls (2.75,-1.75) and (3.25,-1.75) .. (3.5,-1.5);
	\draw[knot] (2.5,-2.5) .. controls  (2.75,-2.25) and (3.25,-2.25) .. (3.5,-2.5);
        \node at (3,-2.35) {$\bullet$};
}~-~\tikz[baseline=-6.5ex, scale=.45]
{
	\draw[dotted] (3,-2) circle(0.707);
	\draw[knot] (2.5,-1.5) .. controls (2.75,-1.75) and (3.25,-1.75) .. (3.5,-1.5);
	\draw[knot] (2.5,-2.5) .. controls  (2.75,-2.25) and (3.25,-2.25) .. (3.5,-2.5);
         \node at (3,-1.7) {$\bullet$};
} & i=-2k\\[10pt]
\tikz[baseline=-6.5ex, scale=.45]
{
	\draw[dotted] (3,-2) circle(0.707);
	\draw[knot] (2.5,-1.5) .. controls (2.75,-1.75) and (3.25,-1.75) .. (3.5,-1.5);
	\draw[knot] (2.5,-2.5) .. controls  (2.75,-2.25) and (3.25,-2.25) .. (3.5,-2.5);
         \node at (3,-2.35) {$\bullet$};
}~+ XY ~\tikz[baseline=-6.5ex, scale=.45]
{
	\draw[dotted] (3,-2) circle(0.707);
	\draw[knot] (2.5,-1.5) .. controls (2.75,-1.75) and (3.25,-1.75) .. (3.5,-1.5);
	\draw[knot] (2.5,-2.5) .. controls  (2.75,-2.25) and (3.25,-2.25) .. (3.5,-2.5);
         \node at (3,-1.7) {$\bullet$};
} & i=-2k -1\end{cases}
\]
for all $i<0$. Notice that taking $X, Y, Z \mapsto 1$ recovers a 2-strand projector of \cite{https://doi.org/10.48550/arxiv.1005.5117}; taking $X, Z \mapsto 1$ and $Y \mapsto -1$ recovers the one of \cite{Sch_tz_2022}.

\begin{proposition}
$P_2 \in \mathrm{Chom}(2)^\mathscr{G}$.
\end{proposition}

\begin{proof}
For the first case, notice that $C_{-1} \circ C_{-2} = 0$ just as in the even case: passing a dot below a saddle and then back up the opposing side introduces two changes of chronology whose evaluations are inverse to one another, since $\lambda(v,u) = \lambda(u,v)^{-1}$.

The other two cases are slightly different since dots may not move past each other freely, but rather by multiplication by $XY$: 
\begin{align*}
& \left(\tikz[baseline=-6.5ex, scale=.45]
{
	\draw[dotted] (3,-2) circle(0.707);
	\draw[knot] (2.5,-1.5) .. controls (2.75,-1.75) and (3.25,-1.75) .. (3.5,-1.5);
	\draw[knot] (2.5,-2.5) .. controls  (2.75,-2.25) and (3.25,-2.25) .. (3.5,-2.5);
         \node at (3,-2.35) {$\bullet$};
}~+ XY ~\tikz[baseline=-6.5ex, scale=.45]
{
	\draw[dotted] (3,-2) circle(0.707);
	\draw[knot] (2.5,-1.5) .. controls (2.75,-1.75) and (3.25,-1.75) .. (3.5,-1.5);
	\draw[knot] (2.5,-2.5) .. controls  (2.75,-2.25) and (3.25,-2.25) .. (3.5,-2.5);
         \node at (3,-1.7) {$\bullet$};
}\right)
\circ
\left(\tikz[baseline=-6.5ex, scale=.45]
{
	\draw[dotted] (3,-2) circle(0.707);
	\draw[knot] (2.5,-1.5) .. controls (2.75,-1.75) and (3.25,-1.75) .. (3.5,-1.5);
	\draw[knot] (2.5,-2.5) .. controls  (2.75,-2.25) and (3.25,-2.25) .. (3.5,-2.5);
        \node at (3,-2.35) {$\bullet$};
}~-~\tikz[baseline=-6.5ex, scale=.45]
{
	\draw[dotted] (3,-2) circle(0.707);
	\draw[knot] (2.5,-1.5) .. controls (2.75,-1.75) and (3.25,-1.75) .. (3.5,-1.5);
	\draw[knot] (2.5,-2.5) .. controls  (2.75,-2.25) and (3.25,-2.25) .. (3.5,-2.5);
         \node at (3,-1.7) {$\bullet$};
}\right) \\
&=
\tikz[baseline=6ex, scale = .475]{
    \draw (-2,0) arc (180:360:2 and 0.6);
    \draw[dashed] (2,0) arc (0:180:2 and 0.6);
    \draw (-2,4) arc (180:360:2 and 0.6);
    \draw (2,4) arc (0:180:2 and 0.6);
    \draw (-2,0) -- (-2,4);
    \draw (2,0) -- (2,4);
    \draw[rounded corners = 3mm] (-0.395777629218,4.588134683673) -- (-.4,4) -- (-0.80159834416,4-0.54969956205);
    \draw[rounded corners = 3mm] (0.398899441995,4.587944836839) -- (.4,4) -- (0.802005173956,4-0.549646152616);
    \draw (-0.395777629218,4.588134683673) -- (-0.395777629218,0.588134683673);
    \draw (-0.80159834416,4-0.54969956205) -- (-0.80159834416,-0.54969956205);
    \draw (0.398899441995,4.587944836839) -- (0.398899441995,0.587944836839); 
    \draw (0.802005173956,4-0.549646152616) -- (0.802005173956,-0.549646152616);
    \draw[rounded corners = 3mm] (-0.395777629218,0.588134683673) -- (-.4,0) -- (-0.80159834416,-0.54969956205);
    \draw[rounded corners = 3mm] (0.398899441995,0.587944836839) -- (.4,0) -- (0.802005173956,-0.549646152616);
    \node at (.6,1.3) {$\bullet$};
    \node at (.6, 2.6) {$\bullet$};
    \node[left] at (-2,0) {top};
    \node[right] at (2,0) {bot};
}
- ~
\tikz[baseline=6ex, scale = .475]{
    \draw (-2,0) arc (180:360:2 and 0.6);
    \draw[dashed] (2,0) arc (0:180:2 and 0.6);
    \draw (-2,4) arc (180:360:2 and 0.6);
    \draw (2,4) arc (0:180:2 and 0.6);
    \draw (-2,0) -- (-2,4);
    \draw (2,0) -- (2,4);
    \draw[rounded corners = 3mm] (-0.395777629218,4.588134683673) -- (-.4,4) -- (-0.80159834416,4-0.54969956205);
    \draw[rounded corners = 3mm] (0.398899441995,4.587944836839) -- (.4,4) -- (0.802005173956,4-0.549646152616);
    \draw (-0.395777629218,4.588134683673) -- (-0.395777629218,0.588134683673);
    \draw (-0.80159834416,4-0.54969956205) -- (-0.80159834416,-0.54969956205);
    \draw (0.398899441995,4.587944836839) -- (0.398899441995,0.587944836839); 
    \draw (0.802005173956,4-0.549646152616) -- (0.802005173956,-0.549646152616);
    \draw[rounded corners = 3mm] (-0.395777629218,0.588134683673) -- (-.4,0) -- (-0.80159834416,-0.54969956205);
    \draw[rounded corners = 3mm] (0.398899441995,0.587944836839) -- (.4,0) -- (0.802005173956,-0.549646152616);
    \node at (-.6,1.3) {$\bullet$};
    \node at (.6, 2.6) {$\bullet$};
}
~+ XY~ 
\tikz[baseline=6ex, scale = .475]{
    \draw (-2,0) arc (180:360:2 and 0.6);
    \draw[dashed] (2,0) arc (0:180:2 and 0.6);
    \draw (-2,4) arc (180:360:2 and 0.6);
    \draw (2,4) arc (0:180:2 and 0.6);
    \draw (-2,0) -- (-2,4);
    \draw (2,0) -- (2,4);
    \draw[rounded corners = 3mm] (-0.395777629218,4.588134683673) -- (-.4,4) -- (-0.80159834416,4-0.54969956205);
    \draw[rounded corners = 3mm] (0.398899441995,4.587944836839) -- (.4,4) -- (0.802005173956,4-0.549646152616);
    \draw (-0.395777629218,4.588134683673) -- (-0.395777629218,0.588134683673);
    \draw (-0.80159834416,4-0.54969956205) -- (-0.80159834416,-0.54969956205);
    \draw (0.398899441995,4.587944836839) -- (0.398899441995,0.587944836839); 
    \draw (0.802005173956,4-0.549646152616) -- (0.802005173956,-0.549646152616);
    \draw[rounded corners = 3mm] (-0.395777629218,0.588134683673) -- (-.4,0) -- (-0.80159834416,-0.54969956205);
    \draw[rounded corners = 3mm] (0.398899441995,0.587944836839) -- (.4,0) -- (0.802005173956,-0.549646152616);
    \node at (.6,1.3) {$\bullet$};
    \node at (-.6, 2.6) {$\bullet$};
}
~- XY~
\tikz[baseline=6ex, scale = .475]{
    \draw (-2,0) arc (180:360:2 and 0.6);
    \draw[dashed] (2,0) arc (0:180:2 and 0.6);
    \draw (-2,4) arc (180:360:2 and 0.6);
    \draw (2,4) arc (0:180:2 and 0.6);
    \draw (-2,0) -- (-2,4);
    \draw (2,0) -- (2,4);
    \draw[rounded corners = 3mm] (-0.395777629218,4.588134683673) -- (-.4,4) -- (-0.80159834416,4-0.54969956205);
    \draw[rounded corners = 3mm] (0.398899441995,4.587944836839) -- (.4,4) -- (0.802005173956,4-0.549646152616);
    \draw (-0.395777629218,4.588134683673) -- (-0.395777629218,0.588134683673);
    \draw (-0.80159834416,4-0.54969956205) -- (-0.80159834416,-0.54969956205);
    \draw (0.398899441995,4.587944836839) -- (0.398899441995,0.587944836839); 
    \draw (0.802005173956,4-0.549646152616) -- (0.802005173956,-0.549646152616);
    \draw[rounded corners = 3mm] (-0.395777629218,0.588134683673) -- (-.4,0) -- (-0.80159834416,-0.54969956205);
    \draw[rounded corners = 3mm] (0.398899441995,0.587944836839) -- (.4,0) -- (0.802005173956,-0.549646152616);
    \node at (-.6,1.3) {$\bullet$};
    \node at (-.6, 2.6) {$\bullet$};
} \\
& =
- XY ~ 
\tikz[baseline=6ex, scale = .475]{
    \draw (-2,0) arc (180:360:2 and 0.6);
    \draw[dashed] (2,0) arc (0:180:2 and 0.6);
    \draw (-2,4) arc (180:360:2 and 0.6);
    \draw (2,4) arc (0:180:2 and 0.6);
    \draw (-2,0) -- (-2,4);
    \draw (2,0) -- (2,4);
    \draw[rounded corners = 3mm] (-0.395777629218,4.588134683673) -- (-.4,4) -- (-0.80159834416,4-0.54969956205);
    \draw[rounded corners = 3mm] (0.398899441995,4.587944836839) -- (.4,4) -- (0.802005173956,4-0.549646152616);
    \draw (-0.395777629218,4.588134683673) -- (-0.395777629218,0.588134683673);
    \draw (-0.80159834416,4-0.54969956205) -- (-0.80159834416,-0.54969956205);
    \draw (0.398899441995,4.587944836839) -- (0.398899441995,0.587944836839); 
    \draw (0.802005173956,4-0.549646152616) -- (0.802005173956,-0.549646152616);
    \draw[rounded corners = 3mm] (-0.395777629218,0.588134683673) -- (-.4,0) -- (-0.80159834416,-0.54969956205);
    \draw[rounded corners = 3mm] (0.398899441995,0.587944836839) -- (.4,0) -- (0.802005173956,-0.549646152616);
    \node at (.6,1.3) {$\bullet$};
    \node at (-.6, 2.6) {$\bullet$};
}
~+ XY~
\tikz[baseline=6ex, scale = .475]{
    \draw (-2,0) arc (180:360:2 and 0.6);
    \draw[dashed] (2,0) arc (0:180:2 and 0.6);
    \draw (-2,4) arc (180:360:2 and 0.6);
    \draw (2,4) arc (0:180:2 and 0.6);
    \draw (-2,0) -- (-2,4);
    \draw (2,0) -- (2,4);
    \draw[rounded corners = 3mm] (-0.395777629218,4.588134683673) -- (-.4,4) -- (-0.80159834416,4-0.54969956205);
    \draw[rounded corners = 3mm] (0.398899441995,4.587944836839) -- (.4,4) -- (0.802005173956,4-0.549646152616);
    \draw (-0.395777629218,4.588134683673) -- (-0.395777629218,0.588134683673);
    \draw (-0.80159834416,4-0.54969956205) -- (-0.80159834416,-0.54969956205);
    \draw (0.398899441995,4.587944836839) -- (0.398899441995,0.587944836839); 
    \draw (0.802005173956,4-0.549646152616) -- (0.802005173956,-0.549646152616);
    \draw[rounded corners = 3mm] (-0.395777629218,0.588134683673) -- (-.4,0) -- (-0.80159834416,-0.54969956205);
    \draw[rounded corners = 3mm] (0.398899441995,0.587944836839) -- (.4,0) -- (0.802005173956,-0.549646152616);
    \node at (.6,1.3) {$\bullet$};
    \node at (-.6, 2.6) {$\bullet$};
}
\end{align*}
as desired, recalling that two dots are evaluated as zero. The other composition is the same.
\end{proof}

\begin{proposition}
The chain complex $P_2\in \mathrm{Chom}(2)^\mathscr{G}$ is a unified Cooper-Krushkal projector.
\end{proposition}

\begin{proof}

(CK1) is satisfied clearly. We must check (CK2), that $P_2$ is killed by turnbacks. We will show that $e_1 \otimes P_2 \simeq 0$; the other direction is totally similar.

We have
\[e_1 \otimes \varphi_{\left(\tikz[baseline={([yshift=-.5ex]current bounding box.center)}, scale=.35]{
    \draw[red, thick] (.5, -.25) -- (.5, .25);
    \draw[knot] (0,.5) to[out=-70,in=250] (1,.5);
    \draw[knot] (0,-.5) to[out=70,in=110] (1,-.5);
    }, (-n, -n) \right)}
\tikz[baseline={([yshift=-.5ex]current bounding box.center)}, scale=.7]{
    \draw[knot] (0,.5) to[out=-70,in=250] (1,.5);
    \draw[knot] (0,-.5) to[out=70,in=110] (1,-.5);
    }
\cong  
\varphi_{\left(\tikz[baseline={([yshift=-.5ex]current bounding box.center)}, scale=.35]{
        \draw[red, thick] (0.5, -0.35) -- (0.5, -0.75);
        \draw[knot] (0,1) to[out=-70,in=250] (1,1);
        \draw[knot] (0,-1) to[out=70,in=110] (1,-1);
        \draw[knot] (.5,0) circle (.35);
        }, (-n, -n) \right)} 
\tikz[baseline={([yshift=-.5ex]current bounding box.center)}, scale=.7]{
        \draw[knot] (0,1) to[out=-70,in=250] (1,1);
        \draw[knot] (0,-1) to[out=70,in=110] (1,-1);
        \draw[knot] (.5,0) circle (.35);
}.\]
    Thus, the previously ambiguous saddles appearing in the shifting functors of $P_2$ are seen to be a merge upon tensoring with $e_1$. Merges have the effect of shifting $\mathbb{Z}\times \mathbb{Z}$ degree by $(-1, 0)$, so we conclude that
\[\varphi_{\left(\tikz[baseline={([yshift=-.5ex]current bounding box.center)}, scale=.35]{
        \draw[red, thick] (0.5, -0.35) -- (0.5, -0.75);
        \draw[knot] (0,1) to[out=-70,in=250] (1,1);
        \draw[knot] (0,-1) to[out=70,in=110] (1,-1);
        \draw[knot] (.5,0) circle (.35);
        }, (-n, -n) \right)} \cong \{-(n+1), -n\}.
\]
Consequently, the chain complex $e_1\otimes P_2$ has the form
\[
\tikz[xscale=3.5]{
    \node(dots) at (.8,0) {$\cdots$};
    \node(res1) at (1,0) {$\tikz[anchor=base, baseline, scale=.65]{
        \draw[knot] (0,1) to[out=-70,in=250] (1,1);
        \draw[knot] (0,-1) to[out=70,in=110] (1,-1);
        \draw[knot] (.5,0) circle (.35);
        }$};
    \draw(res1.south east) node[anchor=south west] {$\{-3,-2\}$};
    \node(res2) at (2,0) {$\tikz[anchor=base, baseline, scale=.65]{
        \draw[knot] (0,1) to[out=-70,in=250] (1,1);
        \draw[knot] (0,-1) to[out=70,in=110] (1,-1);
        \draw[knot] (.5,0) circle (.35);
        }$};
    \draw(res2.south east) node[anchor=south west] {$\{-2,-1\}$};
    \node(res3) at (3,0) {$\tikz[anchor=base, baseline, scale=.65]{
        \draw[knot] (0,1) to[out=-70,in=250] (1,1);
        \draw[knot] (0,-1) to[out=70,in=110] (1,-1);
        \draw[knot] (.5,0) circle (.35);
        }$};
    \draw(res3.south east) node[anchor=south west] {$\{-1,0\}$};
    \node(res4) at (4,0) {$\tikz[anchor=base, baseline, scale=.65]{
        \draw[knot] (0,2) to[out=-70,in=250] (1,2);
         \draw [knot, rounded corners = 2mm] (0,0) -- (.45, .5) -- (0,1) -- (.5,1.45) -- (1,1) -- (.55,.5) -- (1,0);       
        }$};
    \draw[->] (res1) to[out=0,in=180] node[pos=.5,above,arrows=-]
        {$\tikz[anchor=base, baseline, scale=.35]{
            \draw[knot] (0,1) to[out=-70,in=250] (1,1);
            \draw[knot] (0,-1) to[out=70,in=110] node[pos=.5,inner sep=0](bot){} (1,-1);
            \draw[knot] (.5,0) circle (.35);
            \node at (.5, -.96) {$\bullet$};
        }
        + XY
        \tikz[anchor=base, baseline, scale=.35]{
            \draw[knot] (0,1) to[out=-70,in=250] (1,1);
            \draw[knot] (0,-1) to[out=70,in=110] node[pos=.5,inner sep=0](bot){} (1,-1);
            \draw[knot] (.5,0) circle (.35);
            \node at (.5,-.6) {$\bullet$};
        }$} (res2);
    \draw[->] (res2) to[out=0,in=180] node[pos=.5,above,arrows=-]      
        {$\tikz[anchor=base, baseline, scale=.35]{
            \draw[knot] (0,1) to[out=-70,in=250] (1,1);
            \draw[knot] (0,-1) to[out=70,in=110] node[pos=.5,inner sep=0](bot){} (1,-1);
            \draw[knot] (.5,0) circle (.35);
            \node at (.5, -.96) {$\bullet$};
        }
        -
        \tikz[anchor=base, baseline, scale=.35]{
            \draw[knot] (0,1) to[out=-70,in=250] (1,1);
            \draw[knot] (0,-1) to[out=70,in=110] node[pos=.5,inner sep=0](bot){} (1,-1);
            \draw[knot] (.5,0) circle (.35);
            \node at (.5,-.6) {$\bullet$};
        }$} (res3);
    \draw[->] (res3) to[out=0,in=180] node[pos=.5,above,arrows=-]
        {$\tikz[anchor=base, baseline, scale=.35]{
            \draw[red, thick] (0.5, -0.35) -- (0.5, -0.75);
            \draw[knot] (0,1) to[out=-70,in=250] (1,1);
            \draw[knot] (0,-1) to[out=70,in=110] (1,-1);
            \draw[knot] (.5,0) circle (.35);
        }$} (res4);
    }.
\]
Delooping yields the complex
\[
\tikz[scale=3.5]{
    \node(dots) at (.8,0) {$\cdots$};
    \node(res1) at (1,.3) {$\tikz[anchor=base, baseline, scale=.65]{
        \draw[knot] (0,1) to[out=-70,in=250] (1,1);
        \draw[knot] (0,-1) to[out=70,in=110] (1,-1);
        }$};
    \draw(res1.south east) node[anchor=south west] {$\{-3,-3\}$};
    \node(oplus1) at (1,0) {$\oplus$};
    \node(res1') at (1,-.3) {$\tikz[anchor=base, baseline, scale=.65]{
        \draw[knot] (0,1) to[out=-70,in=250] (1,1);
        \draw[knot] (0,-1) to[out=70,in=110] (1,-1);
        }$};
    \draw(res1'.south east) node[anchor=south west] {$\{-2,-2\}$};
    \node(res2) at (2,.3) {$\tikz[anchor=base, baseline, scale=.65]{
        \draw[knot] (0,1) to[out=-70,in=250] (1,1);
        \draw[knot] (0,-1) to[out=70,in=110] (1,-1);
        }$};
    \draw(res2.south east) node[anchor=south west] {$\{-2,-2\}$};
    \node(oplus2) at (2,0) {$\oplus$};
    \node(res2') at (2,-.3) {$\tikz[anchor=base, baseline, scale=.65]{
        \draw[knot] (0,1) to[out=-70,in=250] (1,1);
        \draw[knot] (0,-1) to[out=70,in=110] (1,-1);
        }$};
    \draw(res2'.south east) node[anchor=south west] {$\{-1,-1\}$};
    \node(res3) at (3,.3) {$\tikz[anchor=base, baseline, scale=.65]{
        \draw[knot] (0,1) to[out=-70,in=250] (1,1);
        \draw[knot] (0,-1) to[out=70,in=110] (1,-1);
        }$};
    \draw(res3.south east) node[anchor=south west] {$\{-1,-1\}$};
    \node(oplus3) at (3,0) {$\oplus$};
    \node(res3') at (3,-.3) {$\tikz[anchor=base, baseline, scale=.65]{
        \draw[knot] (0,1) to[out=-70,in=250] (1,1);
        \draw[knot] (0,-1) to[out=70,in=110] (1,-1);
        }$};
    \node(res4) at (4,0) {$\tikz[anchor=base, baseline, scale=.65]{
        \draw[knot] (0,2) to[out=-70,in=250] (1,2);
         \draw[knot, rounded corners = 2mm] (0,0) -- (.45, .5) -- (0,1) -- (.5,1.45) -- (1,1) -- (.55,.5) -- (1,0);       
        }$};
    \draw[->] (res1) to[out=0,in=180] node[pos=.5,above,arrows=-]
        {$YZ^{-1} \tikz[baseline={([yshift=-.5ex]current bounding box.center)}, scale=.35]{
        \draw[knot] (0,.5) to[out=-70,in=250] (1,.5);
        \draw[knot] (0,-.5) to[out=70,in=110] (1,-.5);
        \node at (.5, -.8) {$\bullet$};
        }$} (res2);
    \draw[->] (res1') to[out=0,in=180] node[pos=.5,above,arrows=-]
        {$XZ^{-1} \tikz[baseline={([yshift=-.5ex]current bounding box.center)}, scale=.35]{
        \draw[knot] (0,.5) to[out=-70,in=250] (1,.5);
        \draw[knot] (0,-.5) to[out=70,in=110] (1,-.5);
        \node at (.5, -.8) {$\bullet$};
        }$} (res2');
    \draw[->] (res1') to node[pos=.35,above,arrows=-]
        {$XY$} (res2);
    \draw[->] (res2) to[out=0,in=180] node[pos=.5,above,arrows=-]
        {$YZ^{-1} \tikz[baseline={([yshift=-.5ex]current bounding box.center)}, scale=.35]{
        \draw[knot] (0,.5) to[out=-70,in=250] (1,.5);
        \draw[knot] (0,-.5) to[out=70,in=110] (1,-.5);
        \node at (.5, -.8) {$\bullet$};
        }$} (res3);
    \draw[->] (res2') to[out=0,in=180] node[pos=.5,above,arrows=-]
        {$XZ^{-1} \tikz[baseline={([yshift=-.5ex]current bounding box.center)}, scale=.35]{
        \draw[knot] (0,.5) to[out=-70,in=250] (1,.5);
        \draw[knot] (0,-.5) to[out=70,in=110] (1,-.5);
        \node at (.5, -.8) {$\bullet$};
        }$} (res3');
    \draw[->] (res2') to node[pos=.35,above,arrows=-]
        {$-1$} (res3);
    \draw[->] (res3) to[out=0,in=135] node[pos=.4,above,arrows=-]
        {$XZ^{-1} \tikz[baseline={([yshift=-.5ex]current bounding box.center)}, scale=.35]{
        \draw[knot] (0,.5) to[out=-70,in=250] (1,.5);
        \draw[knot] (0,-.5) to[out=70,in=110] (1,-.5);
        \node at (.5, -.8) {$\bullet$};
        }$} (res4);
    \draw[->] (res3') to node[pos=.4,above,arrows=-]
        {$1$} (res4);
}
\]
where each of the maps down and to the right are zero and are therefore not pictured. Simplifying the maps after delooping is not difficult---one need only take caution when applying the S1 relation. Noting that each of the nonzero, diagonal maps are invertible, simultaneous Gaussian elimination (Proposition \ref{SimGE}) implies that this complex is homotopy equivalent to the zero complex.
\end{proof}

\subsubsection{Homology of the trace}
\label{sss:homologyoftrace}

As in the even case, the unified projector satisfies a categorification of the closure property $\langle \mathrm{Tr}(p_n)\rangle = [n+1]$. In the $n=2$ case, notice that 
\[\varphi_{\left(\tikz[baseline={([yshift=-.5ex]current bounding box.center)}, scale=.35]{
    \draw[red, thick] (.5, -.25) -- (.5, .25);
    \draw[knot] (0,.5) to[out=-70,in=250] (1,.5);
    \draw[knot] (0,-.5) to[out=70,in=110] (1,-.5);
    \draw[knot, rounded corners=.5mm] (0,.5) -- (0,1) -- (1.5, 1) -- (1.5,-1) -- (0,-1) -- (0,-.5);
    \draw[knot, rounded corners=.5mm] (1,.5) -- (1,.75) -- (1.25,.75) -- (1.25,-.75) -- (1,-.75) -- (1,-.5);
    }, (-n, -n) \right)}
\tikz[baseline={([yshift=-.5ex]current bounding box.center)}, scale=.6]{
    \draw[knot] (0,.5) to[out=-70,in=250] (1,.5);
    \draw[knot] (0,-.5) to[out=70,in=110] (1,-.5);
    \draw[knot, rounded corners=.5mm] (0,.5) -- (0,1) -- (1.5, 1) -- (1.5,-1) -- (0,-1) -- (0,-.5);
    \draw[knot, rounded corners=.5mm] (1,.5) -- (1,.75) -- (1.25,.75) -- (1.25,-.75) -- (1,-.75) -- (1,-.5);
    } ~=~ \tikz[baseline={([yshift=-.5ex]current bounding box.center)}, scale=.6]{
    \draw[knot] (0,.5) to[out=-70,in=250] (1,.5);
    \draw[knot] (0,-.5) to[out=70,in=110] (1,-.5);
    \draw[knot, rounded corners=.5mm] (0,.5) -- (0,1) -- (1.5, 1) -- (1.5,-1) -- (0,-1) -- (0,-.5);
    \draw[knot, rounded corners=.5mm] (1,.5) -- (1,.75) -- (1.25,.75) -- (1.25,-.75) -- (1,-.75) -- (1,-.5);
    } ~\{-n, -(n+1)\}\]
because the typically ambiguous saddle is a split after taking closure. Then, we see that the complex $\mathrm{Tr}^2(P_2)$ has the form
\begin{equation}
\label{eq:2coloredunknotcomplex}
\tikz[baseline={([yshift=-.5ex]current bounding box.center)},xscale=3.5]{
    \node(dots) at (.8,0) {$\cdots$};
    \node(res1) at (1,0) {$\tikz[baseline={([yshift=-.5ex]current bounding box.center)}, scale=.6]{
    \draw[knot] (0,.5) to[out=-70,in=250] (1,.5);
    \draw[knot] (0,-.5) to[out=70,in=110] (1,-.5);
    \draw[knot, rounded corners=.5mm] (0,.5) -- (0,1) -- (1.5, 1) -- (1.5,-1) -- (0,-1) -- (0,-.5);
    \draw[knot, rounded corners=.5mm] (1,.5) -- (1,.75) -- (1.25,.75) -- (1.25,-.75) -- (1,-.75) -- (1,-.5);
    }$};
    \draw(res1.south east) node[anchor=south west] {$\{-2,-3\}$};
    \node(res2) at (2,0) {$\tikz[baseline={([yshift=-.5ex]current bounding box.center)}, scale=.6]{
    \draw[knot] (0,.5) to[out=-70,in=250] (1,.5);
    \draw[knot] (0,-.5) to[out=70,in=110] (1,-.5);
    \draw[knot, rounded corners=.5mm] (0,.5) -- (0,1) -- (1.5, 1) -- (1.5,-1) -- (0,-1) -- (0,-.5);
    \draw[knot, rounded corners=.5mm] (1,.5) -- (1,.75) -- (1.25,.75) -- (1.25,-.75) -- (1,-.75) -- (1,-.5);
    }$};
    \draw(res2.south east) node[anchor=south west] {$\{-1,-2\}$};
    \node(res3) at (3,0) {$\tikz[baseline={([yshift=-.5ex]current bounding box.center)}, scale=.6]{
    \draw[knot] (0,.5) to[out=-70,in=250] (1,.5);
    \draw[knot] (0,-.5) to[out=70,in=110] (1,-.5);
    \draw[knot, rounded corners=.5mm] (0,.5) -- (0,1) -- (1.5, 1) -- (1.5,-1) -- (0,-1) -- (0,-.5);
    \draw[knot, rounded corners=.5mm] (1,.5) -- (1,.75) -- (1.25,.75) -- (1.25,-.75) -- (1,-.75) -- (1,-.5);
    }$};
    \draw(res3.south east) node[anchor=south west] {$\{0,-1\}$};
    \node(res4) at (4,0) {$\tikz[baseline={([yshift=-.5ex]current bounding box.center)}, scale=.6]{
    \draw[knot] (0,.5) to[out=-70,in=70] (0,-.5);
    \draw[knot] (1,.5) to[out=250,in=110] (1,-.5);
    \draw[knot, rounded corners=.5mm] (0,.5) -- (0,1) -- (1.5, 1) -- (1.5,-1) -- (0,-1) -- (0,-.5);
    \draw[knot, rounded corners=.5mm] (1,.5) -- (1,.75) -- (1.25,.75) -- (1.25,-.75) -- (1,-.75) -- (1,-.5);
    }$};
    \draw[->] (res1) to[out=0,in=180] node[pos=.5,above,arrows=-]
        {$(1+XY)~\tikz[baseline={([yshift=-.5ex]current bounding box.center)}, scale=.35]{
    \draw[knot] (0,.5) to[out=-70,in=250] (1,.5);
    \draw[knot] (0,-.5) to[out=70,in=110] (1,-.5);
    \draw[knot, rounded corners=.5mm] (0,.5) -- (0,1) -- (1.5, 1) -- (1.5,-1) -- (0,-1) -- (0,-.5);
    \draw[knot, rounded corners=.5mm] (1,.5) -- (1,.75) -- (1.25,.75) -- (1.25,-.75) -- (1,-.75) -- (1,-.5);
    \node at (1.25,-.6125) {$\bullet$};
    }$} (res2);
    \draw[->] (res2) to[out=0,in=180] node[pos=.5,above,arrows=-]      
        {$0$} (res3);
    \draw[->] (res3) to[out=0,in=180] node[pos=.5,above,arrows=-]
        {$\tikz[baseline={([yshift=-.5ex]current bounding box.center)}, scale=.35]{
    \draw[red, thick] (.5, -.25) -- (.5, .25);
    \draw[knot] (0,.5) to[out=-70,in=250] (1,.5);
    \draw[knot] (0,-.5) to[out=70,in=110] (1,-.5);
    \draw[knot, rounded corners=.5mm] (0,.5) -- (0,1) -- (1.5, 1) -- (1.5,-1) -- (0,-1) -- (0,-.5);
    \draw[knot, rounded corners=.5mm] (1,.5) -- (1,.75) -- (1.25,.75) -- (1.25,-.75) -- (1,-.75) -- (1,-.5);
    }$} (res4);
    }.
\end{equation}
Then we compute
\begin{equation}
\label{eq:2coloredunknot}
H_n(\mathrm{Tr}^2(P_2)) = 
    \begin{cases} 
        R \{2,0\} \oplus R \{1,-1\} & n=0 \\
        0 & n = -1 \\
        R \{-2k+2, -2k\} \oplus \frac{R}{(1+XY)R} \{-2k+1, -2k-1\} & n = -2k \\
        (1+XY) R \{-2k, -2k-2\} & n = -2k-1
    \end{cases}
\end{equation}
whenever $k>0$. Note that we recover the solution in the even case (see Section 4.3.1 of \cite{https://doi.org/10.48550/arxiv.1005.5117}) when $X, Y, Z \mapsto 1$. In the odd case, we see that it is important to specialize coefficients before taking homology, since the dotted map is killed by setting $Y=-1$. In either case , the Euler characteristic reproduces $[3] = q^{2} + 1 + q^{-2}$, despite infinite homology.

\subsubsection{Unified Khovanov homology of the infinite 2-twist}

While we have succeeded in constructing a representative for the second projector by guessing based on the result in the even case, we will prove the existence of unified projectors in the following section based on the suspicion that it ought to correspond to the Khovanov complex of an infinite twist (\cite{rozansky2010infinitetorusbraidyields, Willis_2018, stoffregen2024joneswenzlprojectorskhovanovhomotopy}).

We'll illustrate this fact in the $n=2$ case, using multigluing to compute the Khovanov complex for 2-strand torus braids, yielding a unified Cooper-Krushkal projector. Perhaps it is interesting that the projector obtained in this way has a slightly different appearance compared to $P_2$ in the previous sections, although the homotopy equivalence is obvious.

To a single (negative) crossing we associate the complex
\[
\tikz[xscale=3.5]{
    \node(res0) at (1,0) 
        {$\varphi_{\tikz[baseline={([yshift=-.5ex]current bounding box.center)}, scale=.3]{
        \draw[red,thick] (0.5, -0.225) -- (0.5, 0.225); 
        \draw[knot] (0,.5) to[out=-70,in=250] (1,.5);
        \draw[knot] (0,-.5) to[out=70,in=110] (1,-.5);
        }}\tikz[baseline={([yshift=-.5ex]current bounding box.center)}, scale=.6]{
        \draw[knot] (0,.5) to[out=-70,in=250] (1,.5);
        \draw[knot] (0,-.5) to[out=70,in=110] (1,-.5);
        }$};
    \node(res1) at (2,0) 
        {$\tikz[baseline={([yshift=-.5ex]current bounding box.center)}, scale=.6]{
        \draw[knot] (0,.5) to[out=-70,in=70] (0,-.5);
        \draw[knot] (1,.5) to[out=250,in=110] (1,-.5);
        }$};
    \draw[->] (res0) to[out=0,in=180] node[pos=.5,above,arrows=-]
        {$\tikz[baseline={([yshift=-.5ex]current bounding box.center)}, scale=.35]{
        \draw[red,thick] (0.5, -0.225) -- (0.5, 0.225);
        \draw[knot] (0,.5) to[out=-70,in=250] (1,.5);
        \draw[knot] (0,-.5) to[out=70,in=110] (1,-.5);
        }$} (res1);
}
\]
Thus, to the torus 2-braid with two negative crossings we assoiciate the complex
\[
\left(
\tikz[baseline={([yshift=-.5ex]current bounding box.center)}, xscale=3.5]{
    \node(res0) at (1,0) 
        {$\varphi_{\tikz[baseline={([yshift=-.5ex]current bounding box.center)}, scale=.3]{
        \draw[red,thick] (0.5, -0.225) -- (0.5, 0.225); 
        \draw[knot] (0,.5) to[out=-70,in=250] (1,.5);
        \draw[knot] (0,-.5) to[out=70,in=110] (1,-.5);
        }}\tikz[baseline={([yshift=-.5ex]current bounding box.center)}, scale=.6]{
        \draw[knot] (0,.5) to[out=-70,in=250] (1,.5);
        \draw[knot] (0,-.5) to[out=70,in=110] (1,-.5);
        }$};
    \node(res1) at (2,0) 
        {$\tikz[baseline={([yshift=-.5ex]current bounding box.center)}, scale=.6]{
        \draw[knot] (0,.5) to[out=-70,in=70] (0,-.5);
        \draw[knot] (1,.5) to[out=250,in=110] (1,-.5);
        }$};
    \draw[->] (res0) to[out=0,in=180] node[pos=.5,above,arrows=-]
        {$\tikz[baseline={([yshift=-.5ex]current bounding box.center)}, scale=.35]{
        \draw[red,thick] (0.5, -0.225) -- (0.5, 0.225);
        \draw[knot] (0,.5) to[out=-70,in=250] (1,.5);
        \draw[knot] (0,-.5) to[out=70,in=110] (1,-.5);
        }$} (res1);
}
\right)
\otimes
\left( 
\tikz[baseline={([yshift=-.5ex]current bounding box.center)}, xscale=3.5]{
    \node(res0) at (1,0) 
        {$\varphi_{\tikz[baseline={([yshift=-.5ex]current bounding box.center)}, scale=.3]{
        \draw[red,thick] (0.5, -0.225) -- (0.5, 0.225); 
        \draw[knot] (0,.5) to[out=-70,in=250] (1,.5);
        \draw[knot] (0,-.5) to[out=70,in=110] (1,-.5);
        }}\tikz[baseline={([yshift=-.5ex]current bounding box.center)}, scale=.6]{
        \draw[knot] (0,.5) to[out=-70,in=250] (1,.5);
        \draw[knot] (0,-.5) to[out=70,in=110] (1,-.5);
        }$};
    \node(res1) at (2,0) 
        {$\tikz[baseline={([yshift=-.5ex]current bounding box.center)}, scale=.6]{
        \draw[knot] (0,.5) to[out=-70,in=70] (0,-.5);
        \draw[knot] (1,.5) to[out=250,in=110] (1,-.5);
        }$};
    \draw[->] (res0) to[out=0,in=180] node[pos=.5,above,arrows=-]
        {$\tikz[baseline={([yshift=-.5ex]current bounding box.center)}, scale=.35]{
        \draw[red,thick] (0.5, -0.225) -- (0.5, 0.225);
        \draw[knot] (0,.5) to[out=-70,in=250] (1,.5);
        \draw[knot] (0,-.5) to[out=70,in=110] (1,-.5);
        }$} (res1);
}
\right),
\]
which is isomorphic to
\[
\tikz[xscale=3.5]{
    \node(res00) at (0,1.25) {$\varphi_{\tikz[baseline={([yshift=-.5ex]current bounding box.center)}, scale=.3]{
        \draw[red, thick] (.5,-.35)--(.5,-.7);
        \draw[red, thick] (.5, .35) -- (.5, .7);
        \draw[knot] (0,1) to[out=-70,in=250] (1,1);
        \draw[knot] (0,-1) to[out=70,in=110] (1,-1);
        \draw[knot] (.5,0) circle (.35);
        }}\tikz[baseline={([yshift=-.5ex]current bounding box.center)}, scale=.65]{
        \draw[knot] (0,1) to[out=-70,in=250] (1,1);
        \draw[knot] (0,-1) to[out=70,in=110] (1,-1);
        \draw[knot] (.5,0) circle (.35);
        }$};
    \node(res10) at (1.25,1.25) {$\varphi_{\tikz[baseline={([yshift=-.5ex]current bounding box.center)}, scale=.3]{
        \draw[red, thick] (.5, .3) -- (.5, .6);
        \draw[knot] (0,0) to[out=70,in=-250] (1,0);
        \draw[knot, rounded corners = .75mm] (0,2) -- (.45, 1.5) -- (0,1) -- (.5,0.55) -- (1,1) -- (.55,1.5) -- (1,2);
        }}\tikz[baseline={([yshift=-.5ex]current bounding box.center)}, scale=.65]{
        \draw[knot] (0,0) to[out=70,in=-250] (1,0);
        \draw[knot, rounded corners = 2mm] (0,2) -- (.45, 1.5) -- (0,1) -- (.5,0.55) -- (1,1) -- (.55,1.5) -- (1,2);       
        }$};
    \node(res01) at (1.25,-1.25) {$\varphi_{\tikz[baseline={([yshift=-.5ex]current bounding box.center)}, scale=.3]{
        \draw[red, thick] (.5, 1.7) -- (.5,1.4);
        \draw[knot] (0,2) to[out=-70,in=250] (1,2);
        \draw[knot, rounded corners = .75mm] (0,0) -- (.45, .5) -- (0,1) -- (.5,1.45) -- (1,1) -- (.55,.5) -- (1,0);
        }}\tikz[baseline={([yshift=-.5ex]current bounding box.center)}, scale=.65]{
        \draw[knot] (0,2) to[out=-70,in=250] (1,2);
        \draw[knot, rounded corners = 2mm] (0,0) -- (.45, .5) -- (0,1) -- (.5,1.45) -- (1,1) -- (.55,.5) -- (1,0);       
        }$};
    \node(res11) at (2.5,-1.25) {$\tikz[baseline={([yshift=-.5ex]current bounding box.center)}, scale=.65]{
        \draw[knot, rounded corners = 2mm] (0,0) -- (.45, .5) -- (0,1) -- (.45, 1.5) -- (0,2);
        \draw[knot, rounded corners = 2mm] (1,0) -- (.55,.5) -- (1,1) -- (.55,1.5) -- (1,2);
        }$};
    \node at (1.25, 0) {$\oplus$};
    \draw[->] (res00) to[out=0,in=180] node[pos=.5,above,arrows=-] {$\tikz[baseline={([yshift=-.5ex]current bounding box.center)}, scale=.3]{
        \draw[red, thick] (.5, .35) -- (.5, .7);
        \draw[knot] (0,1) to[out=-70,in=250] (1,1);
        \draw[knot] (0,-1) to[out=70,in=110] node[pos=.5,inner sep=0](bot){} (1,-1);
        \draw[knot]  (.5,0) circle (.35);
        }$} (res10);
    \draw[->] (res00) to node[pos=.5,above,arrows=-] {$\tikz[baseline={([yshift=-.5ex]current bounding box.center)}, scale=.3]{
        \draw[red, thick] (.5,-.35)--(0.5, -0.7);
        \draw[knot] (0,1) to[out=-70,in=250] (1,1);
        \draw[knot] (0,-1) to[out=70,in=110] (1,-1);
        \draw[knot] (.5,0) circle (.35);
        }$} (res01);
    \draw[->] (res10) to node[pos=.5,above,arrows=-] {$\tikz[baseline={([yshift=-.5ex]current bounding box.center)}, scale=.3]{
        \draw[red, thick] (.5, .3) -- (.5, .6);
        \draw[knot] (0,0) to[out=70,in=-250] (1,0);
        \draw[knot, rounded corners = .75mm] (0,2) -- (.45, 1.5) -- (0,1) -- (.5,0.55) -- (1,1) -- (.55,1.5) -- (1,2);
        }$} (res11);
     \draw[->] (res01) to[out=0,in=180] node[pos=.5,above,arrows=-] {$\tikz[baseline={([yshift=-.5ex]current bounding box.center)}, scale=.3]{
        \draw[red, thick] (.5, 1.7) -- (.5,1.4);
        \draw[knot] (0,2) to[out=-70,in=250] (1,2);
        \draw[knot, rounded corners = .75mm] (0,0) -- (.45, .5) -- (0,1) -- (.5,1.45) -- (1,1) -- (.55,.5) -- (1,0);
        }$} (res11);
}
\]
Focusing on the leftmost vertex, we've shown that $\varphi_{\tikz[baseline={([yshift=-.5ex]current bounding box.center)}, scale=.35]{
        \draw[red, thick] (.5,-.35)--(0.5, -0.7);
        \draw[red, thick] (.5, .35) -- (.5, .7);
        \draw[knot] (0,1) to[out=-70,in=250] (1,1);
        \draw[knot] (0,-1) to[out=70,in=110] (1,-1);
        \draw[knot] (.5,0) circle (.35);
        }}
\cong
\varphi_{\left(\tikz[baseline={([yshift=-.5ex]current bounding box.center)}, scale=.35]{
    \draw[red, thick] (.5, -.25) -- (.5, .25);
    \draw[knot] (0,.5) to[out=-70,in=250] (1,.5);
    \draw[knot] (0,-.5) to[out=70,in=110] (1,-.5);
    }, (-1, 0) \right)}.$ 
Moreover, delooping tells us that we have the following isomorphism for any $n$:
\[
\tikz[xscale=3.5]{
    \node(loop) at (0,0) {$\varphi_{\left(\tikz[baseline={([yshift=-.5ex]current bounding box.center)}, scale=.35]{
    \draw[red, thick] (.5, -.25) -- (.5, .25);
    \draw[knot] (0,.5) to[out=-70,in=250] (1,.5);
    \draw[knot] (0,-.5) to[out=70,in=110] (1,-.5);
    }, (-n, 1-n) \right)} \tikz[baseline={([yshift=-.5ex]current bounding box.center)}, scale=.65]{
        \draw[knot] (0,1) to[out=-70,in=250] (1,1);
        \draw[knot] (0,-1) to[out=70,in=110] (1,-1);
        \draw[knot] (.5,0) circle (.35);
        } $};
    \node(Ldeloop) at (-1,-3) {$\varphi_{\left(\tikz[baseline={([yshift=-.5ex]current bounding box.center)}, scale=.35]{
    \draw[red, thick] (.5, -.25) -- (.5, .25);
    \draw[knot] (0,.5) to[out=-70,in=250] (1,.5);
    \draw[knot] (0,-.5) to[out=70,in=110] (1,-.5);
    }, (-n, -n) \right)} \tikz[baseline={([yshift=-.5ex]current bounding box.center)}, scale=.65]{
        \draw[knot] (0,1) to[out=-70,in=250] (1,1);
        \draw[knot] (0,-1) to[out=70,in=110] (1,-1);
        }$};
    \node(Rdeloop) at (1,-3) {$\varphi_{\left(\tikz[baseline={([yshift=-.5ex]current bounding box.center)}, scale=.35]{
    \draw[red, thick] (.5, -.25) -- (.5, .25);
    \draw[knot] (0,.5) to[out=-70,in=250] (1,.5);
    \draw[knot] (0,-.5) to[out=70,in=110] (1,-.5);
    }, (1-n, 1-n) \right)} \tikz[baseline={([yshift=-.5ex]current bounding box.center)}, scale=.65]{
        \draw[knot] (0,1) to[out=-70,in=250] (1,1);
        \draw[knot] (0,-1) to[out=70,in=110] (1,-1);
        }$};
    \node at (0,-3) {$\oplus$};
    \draw[->] (loop) to[out=225,in=75] node[pos=.7,above,arrows=-] {$\tikz[baseline=-4ex, scale=.2]{\draw[domain=0:180] plot ({2*cos(\x)}, {2*sin(\x)});
        \draw (-2,0) arc (180:360:2 and 0.6);
        \draw[dashed] (2,0) arc (0:180:2 and 0.6);
        \draw[->] (0,2.8) [partial ellipse=0:270:5ex and 2ex];}$} (Ldeloop);
    \draw[->] (Ldeloop) to[out=45,in=255] node[pos=.5,below=1.5mm,arrows=-] {$\tikz[baseline=-4ex, scale=.2]{
        \draw [domain=180:360] plot ({2*cos(\x)}, {2*sin(\x)});
        \draw (-2,0) arc (180:360:2 and 0.6);
        \draw (2,0) arc (0:180:2 and 0.6);
        \node at (0,-.35) {$\bullet$};}$} (loop);
    \draw[->] (loop) to[out=315,in=105] node[pos=.7,above,arrows=-] {$\tikz[baseline=-4ex, scale=.2]{\draw[domain=0:180] plot ({2*cos(\x)}, {2*sin(\x)});
        \draw (-2,0) arc (180:360:2 and 0.6);
        \draw[dashed] (2,0) arc (0:180:2 and 0.6);
        \node at (0,.35) {$\bullet$};
        \draw[->] (0,2.8) [partial ellipse=0:270:5ex and 2ex];}$} (Rdeloop);
    \draw[->] (Rdeloop) to[out=135,in=285] node[pos=.5,below=1.5mm,arrows=-]{$\tikz[baseline=-4ex, scale=.2]{
        \draw [domain=180:360] plot ({2*cos(\x)}, {2*sin(\x)});
        \draw (-2,0) arc (180:360:2 and 0.6);
        \draw (2,0) arc (0:180:2 and 0.6);}$} (loop);
}
\]
Thus the original complex is isomorphic to the following complex.
\[
\tikz[xscale=3.5]{
    \node(res00) at (0,1.25) {$\varphi_{\left(\tikz[baseline={([yshift=-.5ex]current bounding box.center)}, scale=.35]{
    \draw[red, thick] (.5, -.25) -- (.5, .25);
    \draw[knot] (0,.5) to[out=-70,in=250] (1,.5);
    \draw[knot] (0,-.5) to[out=70,in=110] (1,-.5);
    }, (-1, -1) \right)} \tikz[baseline={([yshift=-.5ex]current bounding box.center)}, scale=.5]{
    \draw[knot] (0,.5) to[out=-70,in=250] (1,.5);
    \draw[knot] (0,-.5) to[out=70,in=110] (1,-.5);
    }$};
    \node(resx) at (.2, -1.25) {$\varphi_{\tikz[baseline={([yshift=-.5ex]current bounding box.center)}, scale=.35]{
    \draw[red, thick] (.5, -.25) -- (.5, .25);
    \draw[knot] (0,.5) to[out=-70,in=250] (1,.5);
    \draw[knot] (0,-.5) to[out=70,in=110] (1,-.5);
    }} \tikz[baseline={([yshift=-.5ex]current bounding box.center)}, scale=.5]{
    \draw[knot] (0,.5) to[out=-70,in=250] (1,.5);
    \draw[knot] (0,-.5) to[out=70,in=110] (1,-.5);
    }$};
    \node(res10) at (1.25,1.25) {$\varphi_{\tikz[baseline={([yshift=-.5ex]current bounding box.center)}, scale=.35]{
    \draw[red, thick] (.5, -.25) -- (.5, .25);
    \draw[knot] (0,.5) to[out=-70,in=250] (1,.5);
    \draw[knot] (0,-.5) to[out=70,in=110] (1,-.5);
    }} \tikz[baseline={([yshift=-.5ex]current bounding box.center)}, scale=.5]{
    \draw[knot] (0,.5) to[out=-70,in=250] (1,.5);
    \draw[knot] (0,-.5) to[out=70,in=110] (1,-.5);
    }$};
    \node(res01) at (1.25,-1.25) {$\varphi_{\tikz[baseline={([yshift=-.5ex]current bounding box.center)}, scale=.35]{
    \draw[red, thick] (.5, -.25) -- (.5, .25);
    \draw[knot] (0,.5) to[out=-70,in=250] (1,.5);
    \draw[knot] (0,-.5) to[out=70,in=110] (1,-.5);
    }} \tikz[baseline={([yshift=-.5ex]current bounding box.center)}, scale=.5]{
    \draw[knot] (0,.5) to[out=-70,in=250] (1,.5);
    \draw[knot] (0,-.5) to[out=70,in=110] (1,-.5);
    }$};
    \node(res11) at (2.5,-1.25) {$\tikz[baseline={([yshift=-.5ex]current bounding box.center)}, scale=.6]{
        \draw[knot] (0,.5) to[out=-70,in=70] (0,-.5);
        \draw[knot] (1,.5) to[out=250,in=110] (1,-.5);
        }$};
    \node at (1.25, 0) {$\oplus$};
    \draw[->] (res00) to[out=0,in=180] node[pos=.5,above,arrows=-] {$\tikz[baseline={([yshift=-.5ex]current bounding box.center)}, scale=.35]{
    \node at (-1,-.75) {$XZ$};
    \draw[knot] (0,.5) to[out=-70,in=250] (1,.5);
    \draw[knot] (0,-.5) to[out=70,in=110] (1,-.5);
    \node at (.5, -0.4) {$\bullet$};
    }$} (res10);
    \draw[->] (res00) to node[pos=.4,above,arrows=-] {$\tikz[baseline={([yshift=-.5ex]current bounding box.center)}, scale=.35]{
    \node at (-1,-.75) {$XZ$};
    \draw[knot] (0,.5) to[out=-70,in=250] (1,.5);
    \draw[knot] (0,-.5) to[out=70,in=110] (1,-.5);
    \node at (.5, -.85) {$\bullet$};
    }$} (res01);
    \draw[->] (res10) to node[pos=.5,above,arrows=-] {$\tikz[baseline={([yshift=-.5ex]current bounding box.center)}, scale=.35]{
    \draw[red, thick] (.5, -.25) -- (.5, .25);
    \draw[knot] (0,.5) to[out=-70,in=250] (1,.5);
    \draw[knot] (0,-.5) to[out=70,in=110] (1,-.5);
    }$} (res11);
     \draw[->] (res01) to[out=0,in=180] node[pos=.5,above,arrows=-] {$\tikz[baseline={([yshift=-.5ex]current bounding box.center)}, scale=.35]{
    \draw[red, knot] (.5, -.25) -- (.5, .25);
    \draw[knot] (0,.5) to[out=-70,in=250] (1,.5);
    \draw[knot] (0,-.5) to[out=70,in=110] (1,-.5);
    }$} (res11);
    \draw[->] (resx) to node[pos=.5,above,arrows=-] {$1$} (res01);
    \draw (.23, -.75) -- (.55, -.145);
    \draw[->] (.725, .185) -- (1.05, .8);
    \node at (.3, -.4) {$1$};
    \node at (0, 0) {$\oplus$};
}
\]
The $XZ = \lambda((-1,0), (-1,-1))$ factor comes from sliding a dot past a merge:
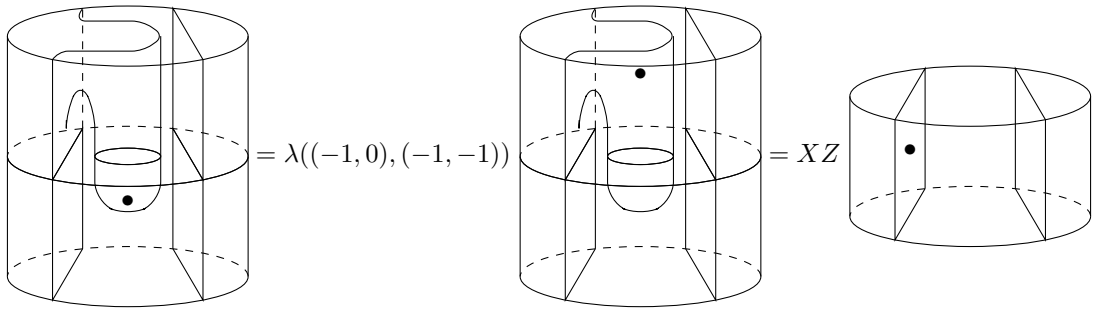
\[
\tikz[baseline={([yshift=-.5ex]current bounding box.center)}, scale = .4]{
        \draw (-4,0) arc (180:360:4 and 1);
        \draw[dashed] (4,0) arc (0:180:4 and 1);
        \draw (-4,4) arc (180:360:4 and 1);
        \draw (-4,4) -- (-4,0);
        \draw (4,4) -- (4,0);
        \draw (-1.5,0.927) -- (-2.5, -0.781) -- (-2.5, -0.781+4) -- (-1.5,4.927) -- (-1.5,0.927);
        \draw (1.5,0.927) -- (2.5, -0.781) -- (2.5, -0.781+4) -- (1.5,4.927) -- (1.5,0.927);
        \draw (-1.095,4) arc (180:360:1.095 and 0.274);
        \draw (1.095,4) arc (0:180:1.095 and 0.274);
        \draw[rounded corners=7.5pt] (-1.095,4) -- (-1.095,2.5) -- (0,2) -- (1.095,2.5) -- (1.095,4);
        \node at (0,2.5) {$\bullet$};
\begin{scope}[shift={(0,4)}]
        \draw (-4,0) arc (180:360:4 and 1);
        \draw[dashed] (4,0) arc (0:180:4 and 1);
        \draw (-4,4) arc (180:360:4 and 1);
        \draw (4,4) arc (0:180:4 and 1);
        \draw (-4,4) -- (-4,0);
        \draw (4,4) -- (4,0);
        \draw (-1.5,0.927) -- (-2.5, -0.781) -- (-2.5, -0.781+4);
        \draw[dashed] (-1.5,0.927) -- (-1.5,4.927);
        \draw[rounded corners] (-1.5,4.927) -- (-1.75, 4.5) -- (0.5,4.5) -- (1.2175,4) -- (0.5, 3.5) -- (-2.335, 3.5) -- (-2.5, -0.781+4);
        \draw (1.5,0.927) -- (2.5, -0.781) -- (2.5, -0.781+4) -- (1.5,4.927) -- (1.5,0.927);
        \draw (-1.095,0) arc (180:360:1.095 and 0.274);
        \draw (1.095,0) arc (0:180:1.095 and 0.274);
        \draw (1.095,0) -- (1.095,4);
        \draw[rounded corners=6.5pt] (-1.095,0) -- (-1.095,1.5) -- (-1.5695, 2.5) -- (-2.043, 1.5) -- (-2.043, 1.2);
\end{scope}
}
= \lambda((-1,0), (-1,-1))~\tikz[baseline={([yshift=-.5ex]current bounding box.center)}, scale = .4]{
        \draw (-4,0) arc (180:360:4 and 1);
        \draw[dashed] (4,0) arc (0:180:4 and 1);
        \draw (-4,4) arc (180:360:4 and 1);
        \draw (-4,4) -- (-4,0);
        \draw (4,4) -- (4,0);
        \draw (-1.5,0.927) -- (-2.5, -0.781) -- (-2.5, -0.781+4) -- (-1.5,4.927) -- (-1.5,0.927);
        \draw (1.5,0.927) -- (2.5, -0.781) -- (2.5, -0.781+4) -- (1.5,4.927) -- (1.5,0.927);
        \draw (-1.095,4) arc (180:360:1.095 and 0.274);
        \draw (1.095,4) arc (0:180:1.095 and 0.274);
        \draw[rounded corners=7.5pt] (-1.095,4) -- (-1.095,2.5) -- (0,2) -- (1.095,2.5) -- (1.095,4);
        \node at (0,6.75) {$\bullet$};
\begin{scope}[shift={(0,4)}]
        \draw (-4,0) arc (180:360:4 and 1);
        \draw[dashed] (4,0) arc (0:180:4 and 1);
        \draw (-4,4) arc (180:360:4 and 1);
        \draw (4,4) arc (0:180:4 and 1);
        \draw (-4,4) -- (-4,0);
        \draw (4,4) -- (4,0);
        \draw (-1.5,0.927) -- (-2.5, -0.781) -- (-2.5, -0.781+4);
        \draw[dashed] (-1.5,0.927) -- (-1.5,4.927);
        \draw[rounded corners] (-1.5,4.927) -- (-1.75, 4.5) -- (0.5,4.5) -- (1.2175,4) -- (0.5, 3.5) -- (-2.335, 3.5) -- (-2.5, -0.781+4);
        \draw (1.5,0.927) -- (2.5, -0.781) -- (2.5, -0.781+4) -- (1.5,4.927) -- (1.5,0.927);
        \draw (-1.095,0) arc (180:360:1.095 and 0.274);
        \draw (1.095,0) arc (0:180:1.095 and 0.274);
        \draw (1.095,0) -- (1.095,4);
        \draw[rounded corners=6.5pt] (-1.095,0) -- (-1.095,1.5) -- (-1.5695, 2.5) -- (-2.043, 1.5) -- (-2.043, 1.2);
\end{scope}
}
=
XZ~\tikz[baseline={([yshift=-.5ex]current bounding box.center)}, scale = .4]{
        \draw (-4,0) arc (180:360:4 and 1);
        \draw[dashed] (4,0) arc (0:180:4 and 1);
        \draw (-4,4) arc (180:360:4 and 1);
        \draw (4,4) arc (0:180:4 and 1);
        \draw (-4,4) -- (-4,0);
        \draw (4,4) -- (4,0);
        \draw (-1.5,0.927) -- (-2.5, -0.781) -- (-2.5, -0.781+4) -- (-1.5,4.927) -- (-1.5,0.927);
        \draw (1.5,0.927) -- (2.5, -0.781) -- (2.5, -0.781+4) -- (1.5,4.927) -- (1.5,0.927);
        \node at (-2,2.2) {$\bullet$};
        }
\]
Then, applying Gaussian elimination, we obtain the following complex.
\[
\tikz[xscale=3.5]{
    \node(res1) at (0,0) {$\varphi_{\left(\tikz[baseline={([yshift=-.5ex]current bounding box.center)}, scale=.35]{
        \draw[red, thick] (.5, -.25) -- (.5, .25);
        \draw[knot] (0,.5) to[out=-70,in=250] (1,.5);
        \draw[knot] (0,-.5) to[out=70,in=110] (1,-.5);
    }, (-1, -1) \right)} \tikz[baseline={([yshift=-.5ex]current bounding box.center)}, scale=.5]{
        \draw[knot] (0,.5) to[out=-70,in=250] (1,.5);
        \draw[knot] (0,-.5) to[out=70,in=110] (1,-.5);
        }$};
    \node(res2) at (1.5,0) {$\varphi_{\tikz[baseline={([yshift=-.5ex]current bounding box.center)}, scale=.35]{
        \draw[red, thick] (.5, -.25) -- (.5, .25);
        \draw[knot] (0,.5) to[out=-70,in=250] (1,.5);
        \draw[knot] (0,-.5) to[out=70,in=110] (1,-.5);
    }} \tikz[baseline={([yshift=-.5ex]current bounding box.center)}, scale=.5]{
        \draw[knot] (0,.5) to[out=-70,in=250] (1,.5);
        \draw[knot] (0,-.5) to[out=70,in=110] (1,-.5);
        }$};
    \node(res3) at (2.5,0) {$\tikz[baseline={([yshift=-.5ex]current bounding box.center)}, scale=.6]{
        \draw[knot] (0,.5) to[out=-70,in=70] (0,-.5);
        \draw[knot] (1,.5) to[out=250,in=110] (1,-.5);
        }$};
    \draw[->] (res1) to node[pos=.5,above,arrows=-] {$\tikz[baseline={([yshift=-.5ex]current bounding box.center)}, scale=.35]{
    \node at (-1,-.75) {$XZ$};
    \draw[knot] (0,.5) to[out=-70,in=250] (1,.5);
    \draw[knot] (0,-.5) to[out=70,in=110] (1,-.5);
    \node at (.5, -0.4) {$\bullet$};
    } - \tikz[baseline={([yshift=-.5ex]current bounding box.center)}, scale=.35]{
    \node at (-1,-.75) {$XZ$};
    \draw[knot] (0,.5) to[out=-70,in=250] (1,.5);
    \draw[knot] (0,-.5) to[out=70,in=110] (1,-.5);
    \node at (.5, -.85) {$\bullet$};
    }$} (res2);
    \draw[->] (res2) to node[pos=.5,above,arrows=-] {$\tikz[baseline={([yshift=-.5ex]current bounding box.center)}, scale=.35]{
    \draw[red, thick] (.5, -.25) -- (.5, .25);
    \draw[knot] (0,.5) to[out=-70,in=250] (1,.5);
    \draw[knot] (0,-.5) to[out=70,in=110] (1,-.5);
    }$} (res3);
}
\]
To stack with another crossing means to tensor this complex with the original single crossing complex. After delooping, this complex has the following form (arrows which are not pictured are zero; dotted arrows are ones which die during Gaussian elimination).
\[
\tikz[xscale=3.5]{
    \node(resA1) at (-.25,1) {$\varphi_{\left(\tikz[baseline={([yshift=-.5ex]current bounding box.center)}, scale=.35]{
        \draw[red, thick] (.5, -.25) -- (.5, .25);
        \draw[knot] (0,.5) to[out=-70,in=250] (1,.5);
        \draw[knot] (0,-.5) to[out=70,in=110] (1,-.5);
        }, (-2, -2) \right)}   \tikz[baseline={([yshift=-.5ex]current bounding box.center)}, scale=.5]{
        \draw[knot] (0,.5) to[out=-70,in=250] (1,.5);
        \draw[knot] (0,-.5) to[out=70,in=110] (1,-.5);
        }$};
    \node(resA2) at (-.25,-1) {$\varphi_{\left(\tikz[baseline={([yshift=-.5ex]current bounding box.center)}, scale=.35]{
        \draw[red, thick] (.5, -.25) -- (.5, .25);
        \draw[knot] (0,.5) to[out=-70,in=250] (1,.5);
        \draw[knot] (0,-.5) to[out=70,in=110] (1,-.5);
        }, (-1, -1) \right)} \tikz[baseline={([yshift=-.5ex]current bounding box.center)}, scale=.5]{
        \draw[knot] (0,.5) to[out=-70,in=250] (1,.5);
        \draw[knot] (0,-.5) to[out=70,in=110] (1,-.5);
        }$};
    \node(resB1) at (1.25, 1) {$\varphi_{\left(\tikz[baseline={([yshift=-.5ex]current bounding box.center)}, scale=.35]{
        \draw[red, thick] (.5, -.25) -- (.5, .25);
        \draw[knot] (0,.5) to[out=-70,in=250] (1,.5);
        \draw[knot] (0,-.5) to[out=70,in=110] (1,-.5);
        }, (-1, -1) \right)} \tikz[baseline={([yshift=-.5ex]current bounding box.center)}, scale=.5]{
        \draw[knot] (0,.5) to[out=-70,in=250] (1,.5);
        \draw[knot] (0,-.5) to[out=70,in=110] (1,-.5);
        }$};
    \node(resB2) at (1.25, -1) {$\varphi_{\tikz[baseline={([yshift=-.5ex]current bounding box.center)}, scale=.35]{
        \draw[red, thick] (.5, -.25) -- (.5, .25);
        \draw[knot] (0,.5) to[out=-70,in=250] (1,.5);
        \draw[knot] (0,-.5) to[out=70,in=110] (1,-.5);
        }} \tikz[baseline={([yshift=-.5ex]current bounding box.center)}, scale=.5]{
        \draw[knot] (0,.5) to[out=-70,in=250] (1,.5);
        \draw[knot] (0,-.5) to[out=70,in=110] (1,-.5);
        }$};
    \node(resB3) at (1.25, -3) {$\varphi_{\left(\tikz[baseline={([yshift=-.5ex]current bounding box.center)}, scale=.35]{
        \draw[red, thick] (.5, -.25) -- (.5, .25);
        \draw[knot] (0,.5) to[out=-70,in=250] (1,.5);
        \draw[knot] (0,-.5) to[out=70,in=110] (1,-.5);
        }, (-1, -1) \right)} \tikz[baseline={([yshift=-.5ex]current bounding box.center)}, scale=.5]{
        \draw[knot] (0,.5) to[out=-70,in=250] (1,.5);
        \draw[knot] (0,-.5) to[out=70,in=110] (1,-.5);
        }$};
    \node(resC1) at (2.75, -1)  {$\varphi_{\tikz[baseline={([yshift=-.5ex]current bounding box.center)}, scale=.35]{
        \draw[red, thick] (.5, -.25) -- (.5, .25);
        \draw[knot] (0,.5) to[out=-70,in=250] (1,.5);
        \draw[knot] (0,-.5) to[out=70,in=110] (1,-.5);
        }}   \tikz[baseline={([yshift=-.5ex]current bounding box.center)}, scale=.5]{
        \draw[knot] (0,.5) to[out=-70,in=250] (1,.5);
        \draw[knot] (0,-.5) to[out=70,in=110] (1,-.5);
        }$};
    \node(resC2) at (2.75, -3)  {$\varphi_{\tikz[baseline={([yshift=-.5ex]current bounding box.center)}, scale=.35]{
        \draw[red, thick] (.5, -.25) -- (.5, .25);
        \draw[knot] (0,.5) to[out=-70,in=250] (1,.5);
        \draw[knot] (0,-.5) to[out=70,in=110] (1,-.5);
        }} \tikz[baseline={([yshift=-.5ex]current bounding box.center)}, scale=.5]{
        \draw[knot] (0,.5) to[out=-70,in=250] (1,.5);
        \draw[knot] (0,-.5) to[out=70,in=110] (1,-.5);
        }$};
    \node(resD) at (3.85, -3) {$\tikz[baseline={([yshift=-.5ex]current bounding box.center)}, scale=.6]{
        \draw[knot] (0,.5) to[out=-70,in=70] (0,-.5);
        \draw[knot] (1,.5) to[out=250,in=110] (1,-.5);
        }$};
    \draw[->] (resA1) to node[pos=.5,above,arrows=-] {$\tikz[baseline={([yshift=-.5ex]current bounding box.center)}, scale=.35]{
        \node at (-1.45,-.75) {$XYZ^2$};
        \draw[knot] (0,.5) to[out=-70,in=250] (1,.5);
        \draw[knot] (0,-.5) to[out=70,in=110] (1,-.5);
        \node at (.5, -0.4) {$\bullet$};
        }$} (resB1);
    \draw[->] (resA2) to node[pos=.65,below,arrows=-] {$-XZ$} (resB1);
    \draw[->, dotted] (resA2) to (resB2);
    \draw[->] (resA2) to node[pos=.5,above, arrows=-] {$1$} (resB3);
    \draw[->, white, line width=2mm] (resA1) to (resB3);
    \draw[->] (resA1) to node[pos=.25,above=1.1mm,arrows=-] {$\tikz[baseline={([yshift=-.5ex]current bounding box.center)}, scale=.35]{
        \node at (-1.35,-.75) {$XZ$};
        \draw[knot] (0,.5) to[out=-70,in=250] (1,.5);
        \draw[knot] (0,-.5) to[out=70,in=110] (1,-.5);
        \node at (.5, -.85) {$\bullet$};
        }$} (resB3);
    \draw[->] (resB1) to node[pos=.5,above,arrows=-] {$\tikz[baseline={([yshift=-.5ex]current bounding box.center)}, scale=.35]{
    \node at (-1,-.75) {$XZ$};
    \draw[knot] (0,.5) to[out=-70,in=250] (1,.5);
    \draw[knot] (0,-.5) to[out=70,in=110] (1,-.5);
    \node at (.5, -0.4) {$\bullet$};
    }$} (resC1);
    \draw[->] (resB2) to node[pos=.2,above,arrows=-] {$1$} (resC1);
    \draw[->, dotted] (resB3) to (resC1);
    \draw[->, dotted] (resB3) to (resC2);
    \draw[->, white, line width=2mm] (resB1) to (resC2);
    \draw[->] (resB1) to node[pos=.5,above=1.1mm,arrows=-] {$\tikz[baseline={([yshift=-.5ex]current bounding box.center)}, scale=.35]{
    \node at (-1,-.75) {$XZ$};
    \draw[knot] (0,.5) to[out=-70,in=250] (1,.5);
    \draw[knot] (0,-.5) to[out=70,in=110] (1,-.5);
    \node at (.5, -.85) {$\bullet$};
    }$} (resC2);
    \draw[->, white, line width=.8mm] (resB2) to (resC2);
    \draw[->] (resB2) to node[pos=.4,above,arrows=-] {$1$} (resC2);
    \draw[->] (resC1) to node[pos=.5,above,arrows=-] {$\tikz[baseline={([yshift=-.5ex]current bounding box.center)}, scale=.35]{
    \draw[red, thick] (.5, -.25) -- (.5, .25);
    \draw[knot] (0,.5) to[out=-70,in=250] (1,.5);
    \draw[knot] (0,-.5) to[out=70,in=110] (1,-.5);
    }$} (resD);
    \draw[->] (resC2) to node[pos=.5,above,arrows=-] {$\tikz[baseline={([yshift=-.5ex]current bounding box.center)}, scale=.35]{
    \draw[red, thick] (.5, -.25) -- (.5, .25);
    \draw[knot] (0,.5) to[out=-70,in=250] (1,.5);
    \draw[knot] (0,-.5) to[out=70,in=110] (1,-.5);
    }$} (resD);
}
\]
The arrows in the left-most column are obtained by applying sphere relations. Note that we can apply S1 and dot-sliding relations to obtain the following equivalences. 
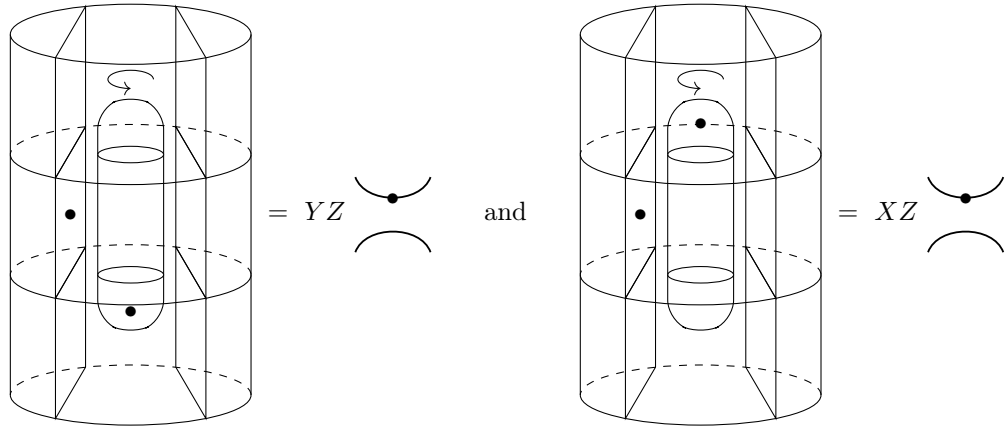
\[
\tikz[baseline={([yshift=-.5ex]current bounding box.center)}, scale = .4]{
        \draw (-4,0) arc (180:360:4 and 1);
        \draw[dashed] (4,0) arc (0:180:4 and 1);
        \draw (-4,4) -- (-4,0);
        \draw (4,4) -- (4,0);
        \draw (-1.5,0.927) -- (-2.5, -0.781) -- (-2.5, -0.781+4) -- (-1.5,4.927) -- (-1.5,0.927);
        \draw (1.5,0.927) -- (2.5, -0.781) -- (2.5, -0.781+4) -- (1.5,4.927) -- (1.5,0.927);
        \draw (-1.095,4) arc (180:360:1.095 and 0.274);
        \draw (1.095,4) arc (0:180:1.095 and 0.274);
        \draw[rounded corners=7.5pt] (-1.095,4) -- (-1.095,2.5) -- (0,2) -- (1.095,2.5) -- (1.095,4);
        \node at (0,2.75) {$\bullet$};
    \begin{scope}[shift={(0,4)}]
        \draw (-4,0) arc (180:360:4 and 1);
        \draw[dashed] (4,0) arc (0:180:4 and 1);
        \draw (-4,4) -- (-4,0);
        \draw (4,4) -- (4,0);
        \draw (-1.5,0.927) -- (-2.5, -0.781) -- (-2.5, -0.781+4) -- (-1.5,4.927) -- (-1.5,0.927);
        \draw (1.5,0.927) -- (2.5, -0.781) -- (2.5, -0.781+4) -- (1.5,4.927) -- (1.5,0.927);
        \draw (-1.095,4) arc (180:360:1.095 and 0.274);
        \draw (1.095,4) arc (0:180:1.095 and 0.274);
        \draw(-1.095,4) -- (-1.095,0);
        \draw (1.095,4) -- (1.095, 0);
        \node at (-2,2) {$\bullet$};
    \end{scope}
    \begin{scope}[shift={(0,8)}]
        \draw (-4,0) arc (180:360:4 and 1);
        \draw[dashed] (4,0) arc (0:180:4 and 1);
        \draw (-4,4) arc (180:360:4 and 1);
        \draw (4,4) arc (0:180:4 and 1);
        \draw (-4,4) -- (-4,0);
        \draw (4,4) -- (4,0);
        \draw (-1.5,0.927) -- (-2.5, -0.781) -- (-2.5, -0.781+4) -- (-1.5,4.927) -- (-1.5,0.927);
        \draw (1.5,0.927) -- (2.5, -0.781) -- (2.5, -0.781+4) -- (1.5,4.927) -- (1.5,0.927);
        \draw[rounded corners=7.5pt] (-1.095,0) -- (-1.095,1.5) -- (0,2) -- (1.095,1.5) -- (1.095,0);
        \draw[->] (0,2.5) [partial ellipse=0:270:5ex and 2ex];
    \end{scope}
}
~=~YZ~ \tikz[baseline={([yshift=-.5ex]current bounding box.center)}, scale=1]{
        \draw[knot] (0,.5) to[out=-70,in=250] (1,.5);
        \draw[knot] (0,-.5) to[out=70,in=110] (1,-.5);
        \node at (.5, 0.21) {$\bullet$};
        }
\qquad \text{and} \qquad
\tikz[baseline={([yshift=-.5ex]current bounding box.center)}, scale = .4]{
        \draw (-4,0) arc (180:360:4 and 1);
        \draw[dashed] (4,0) arc (0:180:4 and 1);
        \draw (-4,4) -- (-4,0);
        \draw (4,4) -- (4,0);
        \draw (-1.5,0.927) -- (-2.5, -0.781) -- (-2.5, -0.781+4) -- (-1.5,4.927) -- (-1.5,0.927);
        \draw (1.5,0.927) -- (2.5, -0.781) -- (2.5, -0.781+4) -- (1.5,4.927) -- (1.5,0.927);
        \draw (-1.095,4) arc (180:360:1.095 and 0.274);
        \draw (1.095,4) arc (0:180:1.095 and 0.274);
        \draw[rounded corners=7.5pt] (-1.095,4) -- (-1.095,2.5) -- (0,2) -- (1.095,2.5) -- (1.095,4);
    \begin{scope}[shift={(0,4)}]
        \draw (-4,0) arc (180:360:4 and 1);
        \draw[dashed] (4,0) arc (0:180:4 and 1);
        \draw (-4,4) -- (-4,0);
        \draw (4,4) -- (4,0);
        \draw (-1.5,0.927) -- (-2.5, -0.781) -- (-2.5, -0.781+4) -- (-1.5,4.927) -- (-1.5,0.927);
        \draw (1.5,0.927) -- (2.5, -0.781) -- (2.5, -0.781+4) -- (1.5,4.927) -- (1.5,0.927);
        \draw (-1.095,4) arc (180:360:1.095 and 0.274);
        \draw (1.095,4) arc (0:180:1.095 and 0.274);
        \draw(-1.095,4) -- (-1.095,0);
        \draw (1.095,4) -- (1.095, 0);
        \node at (-2,2) {$\bullet$};
    \end{scope}
    \begin{scope}[shift={(0,8)}]
        \draw (-4,0) arc (180:360:4 and 1);
        \draw[dashed] (4,0) arc (0:180:4 and 1);
        \draw (-4,4) arc (180:360:4 and 1);
        \draw (4,4) arc (0:180:4 and 1);
        \draw (-4,4) -- (-4,0);
        \draw (4,4) -- (4,0);
        \draw (-1.5,0.927) -- (-2.5, -0.781) -- (-2.5, -0.781+4) -- (-1.5,4.927) -- (-1.5,0.927);
        \draw (1.5,0.927) -- (2.5, -0.781) -- (2.5, -0.781+4) -- (1.5,4.927) -- (1.5,0.927);
        \draw[rounded corners=7.5pt] (-1.095,0) -- (-1.095,1.5) -- (0,2) -- (1.095,1.5) -- (1.095,0);
        \node at (0,1) {$\bullet$};
        \draw[->] (0,2.5) [partial ellipse=0:270:5ex and 2ex];
    \end{scope}
}
~=~XZ~ \tikz[baseline={([yshift=-.5ex]current bounding box.center)}, scale=1]{
        \draw[knot] (0,.5) to[out=-70,in=250] (1,.5);
        \draw[knot] (0,-.5) to[out=70,in=110] (1,-.5);
        \node at (.5, 0.21) {$\bullet$};
        }
\]
Applying Gaussian elimination, the above complex is homotopy equivalent to the following.
\[
\tikz[xscale=3.5]{
    \node(res0) at (-1.7, 0) {$\varphi_{\left(\tikz[baseline={([yshift=-.5ex]current bounding box.center)}, scale=.35]{
        \draw[red, thick] (.5, -.25) -- (.5, .25);
        \draw[knot] (0,.5) to[out=-70,in=250] (1,.5);
        \draw[knot] (0,-.5) to[out=70,in=110] (1,-.5);
    }, (-2, -2) \right)} \tikz[baseline={([yshift=-.5ex]current bounding box.center)}, scale=.5]{
        \draw[knot] (0,.5) to[out=-70,in=250] (1,.5);
        \draw[knot] (0,-.5) to[out=70,in=110] (1,-.5);
        }$};
    \node(res1) at (0,0) {$\varphi_{\left(\tikz[baseline={([yshift=-.5ex]current bounding box.center)}, scale=.35]{
        \draw[red, thick] (.5, -.25) -- (.5, .25);
        \draw[knot] (0,.5) to[out=-70,in=250] (1,.5);
        \draw[knot] (0,-.5) to[out=70,in=110] (1,-.5);
    }, (-1, -1) \right)} \tikz[baseline={([yshift=-.5ex]current bounding box.center)}, scale=.5]{
        \draw[knot] (0,.5) to[out=-70,in=250] (1,.5);
        \draw[knot] (0,-.5) to[out=70,in=110] (1,-.5);
        }$};
    \node(res2) at (1.5,0) {$\varphi_{\tikz[baseline={([yshift=-.5ex]current bounding box.center)}, scale=.35]{
        \draw[red, thick] (.5, -.25) -- (.5, .25);
        \draw[knot] (0,.5) to[out=-70,in=250] (1,.5);
        \draw[knot] (0,-.5) to[out=70,in=110] (1,-.5);
    }} \tikz[baseline={([yshift=-.5ex]current bounding box.center)}, scale=.5]{
        \draw[knot] (0,.5) to[out=-70,in=250] (1,.5);
        \draw[knot] (0,-.5) to[out=70,in=110] (1,-.5);
        }$};
    \node(res3) at (2.5,0) {$\tikz[baseline={([yshift=-.5ex]current bounding box.center)}, scale=.6]{
        \draw[knot] (0,.5) to[out=-70,in=70] (0,-.5);
        \draw[knot] (1,.5) to[out=250,in=110] (1,-.5);
        }$};
    \draw[->] (res0) to node[pos=.5,above,arrows=-] {$\tikz[baseline={([yshift=-.5ex]current bounding box.center)}, scale=.35]{
    \node at (-1.5,-.75) {$XYZ^2$};
    \draw[knot] (0,.5) to[out=-70,in=250] (1,.5);
    \draw[knot] (0,-.5) to[out=70,in=110] (1,-.5);
    \node at (.5, -0.4) {$\bullet$};
    } + \tikz[baseline={([yshift=-.5ex]current bounding box.center)}, scale=.35]{
    \node at (-.8,-.75) {$Z^2$};
    \draw[knot] (0,.5) to[out=-70,in=250] (1,.5);
    \draw[knot] (0,-.5) to[out=70,in=110] (1,-.5);
    \node at (.5, -.85) {$\bullet$};
    }$} (res1);
    \draw[->] (res1) to node[pos=.5,above,arrows=-] {$\tikz[baseline={([yshift=-.5ex]current bounding box.center)}, scale=.35]{
    \node at (-1,-.75) {$XZ$};
    \draw[knot] (0,.5) to[out=-70,in=250] (1,.5);
    \draw[knot] (0,-.5) to[out=70,in=110] (1,-.5);
    \node at (.5, -0.4) {$\bullet$};
    } - \tikz[baseline={([yshift=-.5ex]current bounding box.center)}, scale=.35]{
    \node at (-1,-.75) {$XZ$};
    \draw[knot] (0,.5) to[out=-70,in=250] (1,.5);
    \draw[knot] (0,-.5) to[out=70,in=110] (1,-.5);
    \node at (.5, -.85) {$\bullet$};
    }$} (res2);
    \draw[->] (res2) to node[pos=.5,above,arrows=-] {$\tikz[baseline={([yshift=-.5ex]current bounding box.center)}, scale=.35]{
    \draw[red, thick] (.5, -.25) -- (.5, .25);
    \draw[knot] (0,.5) to[out=-70,in=250] (1,.5);
    \draw[knot] (0,-.5) to[out=70,in=110] (1,-.5);
    }$} (res3);
}
\]
At this point a pattern emerges which controls the complex for any two stranded braid (although this might be easier to see computing the next case; we leave it to the reader). The complex has the form
\[
\tikz[xscale=3.5]{
    \node at (-.5,0) {$\cdots$};
    \node(res0) at (0, 0) {$\varphi_{\left(\tikz[baseline={([yshift=-.5ex]current bounding box.center)}, scale=.35]{
        \draw[red, thick] (.5, -.25) -- (.5, .25);
        \draw[knot] (0,.5) to[out=-70,in=250] (1,.5);
        \draw[knot] (0,-.5) to[out=70,in=110] (1,-.5);
    }, (-3, -3) \right)} \tikz[baseline={([yshift=-.5ex]current bounding box.center)}, scale=.5]{
        \draw[knot] (0,.5) to[out=-70,in=250] (1,.5);
        \draw[knot] (0,-.5) to[out=70,in=110] (1,-.5);
        }$};
    \node(res1) at (1, 0) {$\varphi_{\left(\tikz[baseline={([yshift=-.5ex]current bounding box.center)}, scale=.35]{
        \draw[red, thick] (.5, -.25) -- (.5, .25);
        \draw[knot] (0,.5) to[out=-70,in=250] (1,.5);
        \draw[knot] (0,-.5) to[out=70,in=110] (1,-.5);
    }, (-2, -2) \right)} \tikz[baseline={([yshift=-.5ex]current bounding box.center)}, scale=.5]{
        \draw[knot] (0,.5) to[out=-70,in=250] (1,.5);
        \draw[knot] (0,-.5) to[out=70,in=110] (1,-.5);
        }$};
    \node(res2) at (2,0) {$\varphi_{\left(\tikz[baseline={([yshift=-.5ex]current bounding box.center)}, scale=.35]{
        \draw[red, thick] (.5, -.25) -- (.5, .25);
        \draw[knot] (0,.5) to[out=-70,in=250] (1,.5);
        \draw[knot] (0,-.5) to[out=70,in=110] (1,-.5);
    }, (-1, -1) \right)} \tikz[baseline={([yshift=-.5ex]current bounding box.center)}, scale=.5]{
        \draw[knot] (0,.5) to[out=-70,in=250] (1,.5);
        \draw[knot] (0,-.5) to[out=70,in=110] (1,-.5);
        }$};
    \node(res3) at (3,0) {$\varphi_{\tikz[baseline={([yshift=-.5ex]current bounding box.center)}, scale=.35]{
        \draw[red, thick] (.5, -.25) -- (.5, .25);
        \draw[knot] (0,.5) to[out=-70,in=250] (1,.5);
        \draw[knot] (0,-.5) to[out=70,in=110] (1,-.5);
    }} \tikz[baseline={([yshift=-.5ex]current bounding box.center)}, scale=.5]{
        \draw[knot] (0,.5) to[out=-70,in=250] (1,.5);
        \draw[knot] (0,-.5) to[out=70,in=110] (1,-.5);
        }$};
    \node(res4) at (4,0) {$\tikz[baseline={([yshift=-.5ex]current bounding box.center)}, scale=.6]{
        \draw[knot] (0,.5) to[out=-70,in=70] (0,-.5);
        \draw[knot] (1,.5) to[out=250,in=110] (1,-.5);
        }$};
    \draw[->] (res0) to node[pos=.5,above,arrows=-] {$C_{-4}$} (res1);
    \draw[->] (res1) to node[pos=.5,above,arrows=-] {$C_{-3}$} (res2);
    \draw[->] (res2) to node[pos=.5,above,arrows=-] {$C_{-2}$} (res3);
    \draw[->] (res3) to node[pos=.5,above,arrows=-] {$C_{-1}$} (res4);
}
\]
where
\[
C_i = \begin{cases}
    \tikz[baseline={([yshift=-.5ex]current bounding box.center)}, scale=.5]{
    \draw[red, thick] (.5, -.25) -- (.5, .25);
    \draw[knot] (0,.5) to[out=-70,in=250] (1,.5);
    \draw[knot] (0,-.5) to[out=70,in=110] (1,-.5);
    } & i=-1 \\
    XZ^{2k-1}\left(\tikz[baseline={([yshift=-.5ex]current bounding box.center)}, scale=.5]{
    \draw[knot] (0,.5) to[out=-70,in=250] (1,.5);
    \draw[knot] (0,-.5) to[out=70,in=110] (1,-.5);
    \node at (.5, .21) {$\bullet$};
    } - \tikz[baseline={([yshift=-.5ex]current bounding box.center)}, scale=.5]{
    \draw[knot] (0,.5) to[out=-70,in=250] (1,.5);
    \draw[knot] (0,-.5) to[out=70,in=110] (1,-.5);
    \node at (.5, -.26) {$\bullet$};
    }\right) & i=-2k \\
    Z^{2k}\left(XY
    \tikz[baseline={([yshift=-.5ex]current bounding box.center)}, scale=.5]{
    \draw[knot] (0,.5) to[out=-70,in=250] (1,.5);
    \draw[knot] (0,-.5) to[out=70,in=110] (1,-.5);
    \node at (.5, .21) {$\bullet$};
    } + \tikz[baseline={([yshift=-.5ex]current bounding box.center)}, scale=.5]{
    \draw[knot] (0,.5) to[out=-70,in=250] (1,.5);
    \draw[knot] (0,-.5) to[out=70,in=110] (1,-.5);
    \node at (.5, -.26) {$\bullet$};
    }
    \right) & i=-2k-1
\end{cases}
\]
for all $i<0$. As promised, this complex is homotopy equivalent to the $P_2$ we guessed earlier on.

\subsection{Existence of unified projectors}
\label{ss:egunifiedprojectors}

In \cite{rozansky2010infinitetorusbraidyields}, Rozansky showed that the Khovanov complex associated to an infinte twist on $n$ strands is a Cooper-Krushkal projector. In \cite{Willis_2018}, Willis generalized this argument to the spectral setting. His argument was further generalized in \cite{stoffregen2024joneswenzlprojectorskhovanovhomotopy} for the setting of spectral multimodules.

We will adapt the arguments of \cite{stoffregen2024joneswenzlprojectorskhovanovhomotopy} to prove that unified Cooper-Krushkal projectors exist. As in the work of Stoffregen-Willis, the \textit{left-handed fractional twist complex}, denoted $\mathcal{T}_n$, is the complex associated to the diskular $n$-tangle shown below.
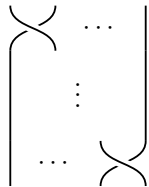
\[
\tikz[baseline={([yshift=-.5ex]current bounding box.center)}, scale=0.6]
{
        \draw[knot] (3,1) -- (3,0);
        \node at (2, 0.5) {$\cdots$};
        \draw[knot] (1,1) to[out=-90, in=90] (0,0);
        \draw[knot, overcross] (0,1) to[out=-90, in=90] (1,0);
        \draw[knot] (0,0) -- (0,-2);
        \draw[knot] (3,0) -- (3, -2);
        \node at (1.5, -0.8) {$\vdots$};
    \begin{scope}[yshift=-3cm]
        \draw[knot] (3,1) to[out=-90, in=90] (2,0);
        \draw[knot, overcross] (2,1) to[out=-90, in=90] (3,0);
        \node at (1, 0.5) {$\cdots$};
        \draw[knot] (0,1) -- (0,0);
    \end{scope}
}
\]
Superscripts will indicate stacking:
\[
\mathcal{T}_n^m = \underbrace{\mathcal{T}_n \otimes \cdots \otimes \mathcal{T}_n}_{m\text{-times}}
\]
with $\mathcal{T}_n^0 = \mathcal{I}_n$. Notice that $\mathcal{T}_n^n$ can be viewed as a pure braid; we call this the \textit{left-handed full twist complex}. Finally, for any $n\in \mathbb{N}$, the \textit{left-handed infinite twist complex}, denoted $\mathcal{T}_n^\infty$, is defined as the colimit of the sequence
\[
\mathcal{T}_n^\infty = \mathrm{colim} \left(\mathcal{T}_n^0 \to \mathcal{T}_n^1 \to \cdots  \mathcal{T}_n^m \to \cdots \right)
\]
where each arrow comes from compositions of maps arising from the cofibration sequence
\[
\mathcal{F}\left(
\tikz[baseline={([yshift=-.5ex]current bounding box.center)}, scale=.55]
{
	\draw[dotted] (-2,-2) circle(0.707);
	\draw[knot] (-1.5,-1.5) .. controls (-1.75,-1.75) and (-1.75,-2.25) .. (-1.5,-2.5);
	\draw[knot] (-2.5,-1.5) .. controls (-2.25,-1.75) and (-2.25,-2.25) ..  (-2.5,-2.5);
}  
 \right)\longrightarrow \mathcal{F}\left(
\tikz[baseline={([yshift=-.5ex]current bounding box.center)}, scale=.55]
{
	\draw[dotted] (.5,.5) circle(0.707);
	\draw[knot](0,0) -- (1,1);
	\fill[fill=white] (.5,.5) circle (.15);
	\draw[knot](1,0) -- (0,1);
}
\right) \longrightarrow 
\varphi_{\tikz[baseline=-6.5ex, scale=.35]
{
	\draw[dotted] (3,-2) circle(0.707);
	\draw[knot] (2.5,-1.5) .. controls (2.75,-1.75) and (3.25,-1.75) .. (3.5,-1.5);
	\draw[knot] (2.5,-2.5) .. controls  (2.75,-2.25) and (3.25,-2.25) .. (3.5,-2.5);
        \draw[red,thick] (3,-1.7) -- (3,-2.3);
}}
\mathcal{F}\left(
\tikz[baseline={([yshift=-.5ex]current bounding box.center)}, scale=.55]
{
	\draw[dotted] (3,-2) circle(0.707);
	\draw[knot] (2.5,-1.5) .. controls (2.75,-1.75) and (3.25,-1.75) .. (3.5,-1.5);
	\draw[knot] (2.5,-2.5) .. controls  (2.75,-2.25) and (3.25,-2.25) .. (3.5,-2.5);
}
  \right) 
[1]
\]
of Proposition \ref{Prop:cofibration}. By the same proposition,
\[
\tikz[baseline={([yshift=-.5ex]current bounding box.center)}, scale=1]{
    \node at (0.85, 0.6) {$\mathcal{T}_n^{nk + r}$};
    \draw[knot] (0.2, 1.2) -- (0.2, 1.5);
    \draw[knot] (1.5, 1.2) -- (1.5, 1.5);
    \node at (0.85, 1.35) {$\cdots$};
    \draw[knot] (0.2, 0) -- (0.2,-0.75);
    \draw[knot] (1.2, 0) to[out=-90, in=90] (0.5, -0.75);
    \draw[knot, overcross] (0.5, 0) to[out=-90, in=90] (1.2, -0.75);
    \draw[knot] (1.5, 0) -- (1.5,-0.75);
    \draw (0,0) rectangle (1.7, 1.2);
}\,
=\, \mathrm{Cone}\left(
\varphi_{
\tikz[baseline={([yshift=-.5ex]current bounding box.center)}, scale=.6]{
    \draw[knot,red] (0.85, -0.2) -- (0.85, -0.55);
    \draw[knot] (0.2, 1.2) -- (0.2, 1.5);
    \draw[knot] (1.5, 1.2) -- (1.5, 1.5);
    \node at (0.85, 1.35) {$\cdots$};
    \draw[knot] (0.2, 0) -- (0.2,-0.75);
    \draw[knot] (0.5, 0) to[out=-90, in=-90] (1.2, 0);
    \draw[knot] (0.5, -0.75) to[out=90, in=90]  (1.2, -0.75);
    \draw[knot] (1.5, 0) -- (1.5,-0.75);
    \draw (0,0) rectangle (1.7, 1.2);
}
}
~
\tikz[baseline={([yshift=-.5ex]current bounding box.center)}, scale=1]{
    \node at (0.85, 0.6) {$\mathcal{T}_n^{nk + r}$};
    \draw[knot] (0.2, 1.2) -- (0.2, 1.5);
    \draw[knot] (1.5, 1.2) -- (1.5, 1.5);
    \node at (0.85, 1.35) {$\cdots$};
    \draw[knot] (0.2, 0) -- (0.2,-0.75);
    \draw[knot] (0.5, 0) -- (0.5, -0.1) to[out=-90, in=-90] (1.2, -0.1) -- (1.2, 0);
    \draw[knot] (0.5, -0.75) -- (0.5, -0.65) to[out=90, in=90] (1.2, -0.65) -- (1.2, -0.75);
    \draw[knot] (1.5, 0) -- (1.5,-0.75);
    \draw (0,0) rectangle (1.7, 1.2);
}\,
\xrightarrow{\qquad \qquad}\,
\tikz[baseline={([yshift=-.5ex]current bounding box.center)}, scale=1]{
    \node at (0.85, 0.6) {$\mathcal{T}_n^{nk + r}$};
    \draw[knot] (0.2, 1.2) -- (0.2, 1.5);
    \draw[knot] (1.5, 1.2) -- (1.5, 1.5);
    \node at (0.85, 1.35) {$\cdots$};
    \draw[knot] (0.2, 0) -- (0.2, -0.75);
    \draw[knot] (1.5, 0) -- (1.5, -0.75);
    \node at (0.85, -0.325) {$\cdots$};
    \draw (0,0) rectangle (1.7, 1.2);
}
\right)
\]
We start our argument by computing a simplification of the term $\mathcal{T}_n^{nk+r} \otimes e_i$, for $0\le r < n$ and $1\le i \le n-1$. Note that
\[
\mathcal{F}\left(
\tikz[baseline={([yshift=-.5ex]current bounding box.center)}, scale=.65]
{
        \draw[dotted] (.5,.5) circle(0.707);
        \draw[knot] (1, 0.5) to[out=-90, in=-60] (0,1);
        \draw[knot, overcross] (0,0) to[out=60, in =90] (1, 0.5);

}
\right)
\cong 
\mathcal{F}\left(
\tikz[baseline={([yshift=-.5ex]current bounding box.center)}, scale=.65]
{
        \draw[dotted] (.5,.5) circle(0.707);
        \draw[knot] (0,0) to[out=60, in=-90] (.7, .5) to[out=90, in=-60] (0,1);
}
\right) \{-1,-1\}
\qquad \text{and} \qquad
\mathcal{F}\left(
\tikz[baseline={([yshift=-.5ex]current bounding box.center)}, scale=.65]
{
        \draw[dotted] (.5,.5) circle(0.707);
        \draw[knot] (0,0) to[out=30, in=-90] (0.8, 0.5);
        \draw[knot] (0.8, 0.5) to[out=90, in=-30] (0,1);
        \draw[knot, overcross] (1,0) to[out=150, in=-90] (0.3, 0.5);
        \draw[knot, overcross] (0.3, 0.5) to[out=90, in=210] (1,1);
        \draw[knot] (0.3, 0.499) -- (0.3, 0.501);
}
\right)
\cong
\varphi_{
\tikz[baseline={([yshift=-.5ex]current bounding box.center)}, scale=.45]
{
        \draw[dotted] (.5,.5) circle(0.707);
        \draw[knot, red] (0.325, 0.5) -- (0.675, 0.5);
        \draw[knot] (0, 0) to[out=30, in=-90] (0.35, 0.5);
        \draw[knot] (0.35, 0.5) to[out=90, in=-30] (0, 1);
        \draw[knot] (1, 0) to[out=150, in=-90] (0.65, 0.5);
        \draw[knot] (0.65, 0.5) to[out=90, in=210] (1, 1);
}
}
~
\mathcal{F}\left(
\tikz[baseline={([yshift=-.5ex]current bounding box.center)}, scale=.65]
{
        \draw[dotted] (.5,.5) circle(0.707);
        \draw[knot] (0, 0) to[out=30, in=-90] (0.35, 0.5);
        \draw[knot] (0.35, 0.5) to[out=90, in=-30] (0, 1);
        \draw[knot] (1, 0) to[out=150, in=-90] (0.65, 0.5);
        \draw[knot] (0.65, 0.5) to[out=90, in=210] (1, 1);
}
\right)
\]
by delooping and Gaussian elimination.

\begin{figure}[ht]
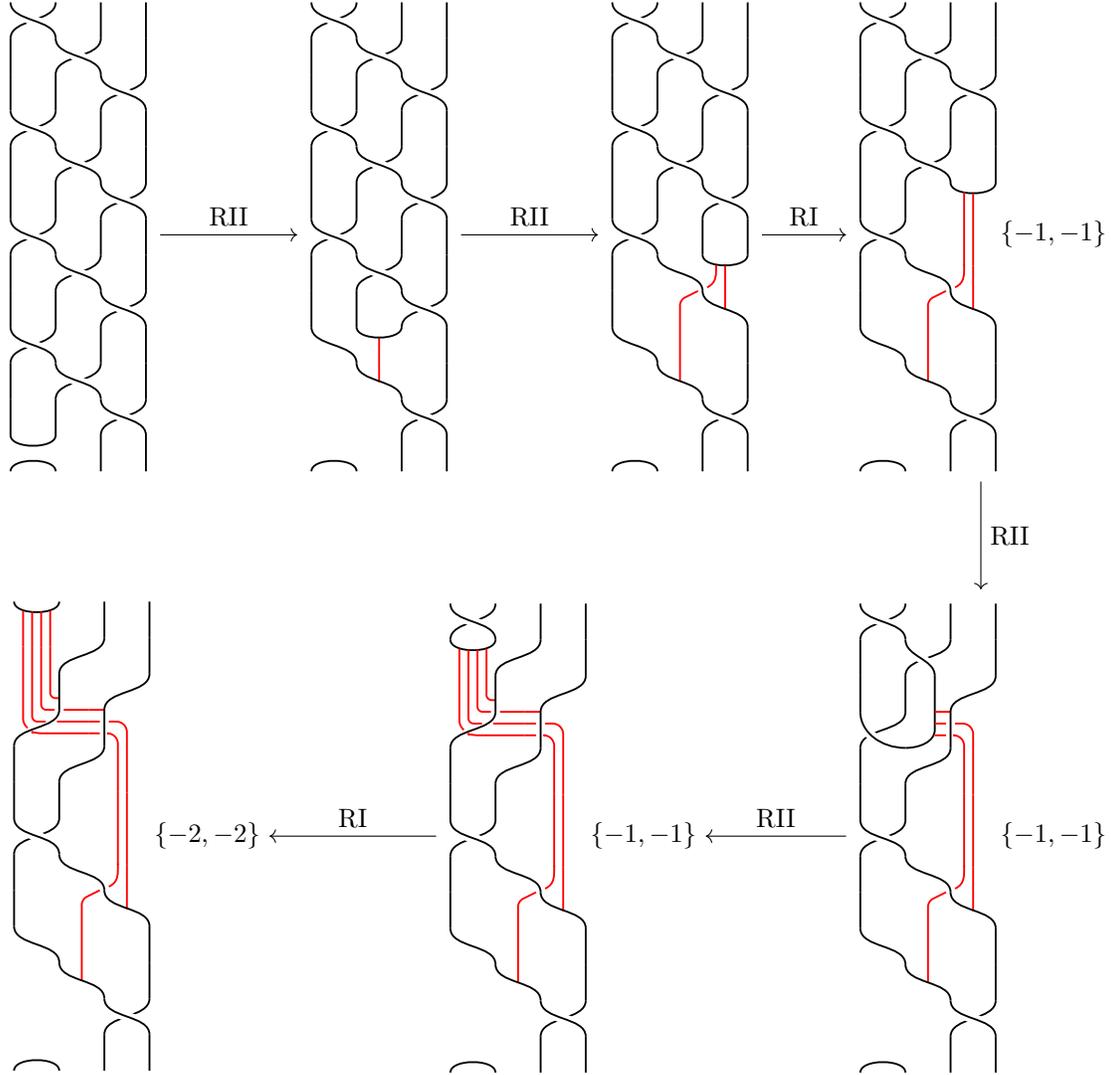

\[
\tikz{
    \node(A) at (0,0) {$\tikz[baseline={([yshift=-.5ex]current bounding box.center)}, scale=0.6, yscale=0.8]
{
        \draw[knot] (3,1) -- (3,0);
        \draw[knot] (2,1) -- (2,0);
        \draw[knot] (1,1) to[out=-90, in=90] (0,0);
        \draw[knot, overcross] (0,1) to[out=-90, in=90] (1,0);
    \begin{scope}[yshift=-1cm]
        \draw[knot] (3,1) -- (3,0);
        \draw[knot] (2,1) to[out=-90, in=90] (1,0);
        \draw[knot, overcross] (1,1) to[out=-90, in=90] (2,0);
        \draw[knot] (0,1) -- (0,0);
    \end{scope}
    \begin{scope}[yshift=-2cm]
        \draw[knot] (3,1) to[out=-90, in=90] (2,0);
        \draw[knot, overcross] (2,1) to[out=-90, in=90] (3,0);
        \draw[knot] (1,1) -- (1,0);
        \draw[knot] (0,1) -- (0,0);
    \end{scope}
    \begin{scope}[yshift=-3cm]
        \draw[knot] (3,1) -- (3,0);
        \draw[knot] (2,1) -- (2,0);
        \draw[knot] (1,1) to[out=-90, in=90] (0,0);
        \draw[knot, overcross] (0,1) to[out=-90, in=90] (1,0);
    \begin{scope}[yshift=-1cm]
        \draw[knot] (3,1) -- (3,0);
        \draw[knot] (2,1) to[out=-90, in=90] (1,0);
        \draw[knot, overcross] (1,1) to[out=-90, in=90] (2,0);
        \draw[knot] (0,1) -- (0,0);
    \end{scope}
    \begin{scope}[yshift=-2cm]
        \draw[knot] (3,1) to[out=-90, in=90] (2,0);
        \draw[knot, overcross] (2,1) to[out=-90, in=90] (3,0);
        \draw[knot] (1,1) -- (1,0);
        \draw[knot] (0,1) -- (0,0);
    \end{scope}
    \end{scope}
    \begin{scope}[yshift=-6cm]
        \draw[knot] (3,1) -- (3,0);
        \draw[knot] (2,1) -- (2,0);
        \draw[knot] (1,1) to[out=-90, in=90] (0,0);
        \draw[knot, overcross] (0,1) to[out=-90, in=90] (1,0);
    \begin{scope}[yshift=-1cm]
        \draw[knot] (3,1) -- (3,0);
        \draw[knot] (2,1) to[out=-90, in=90] (1,0);
        \draw[knot, overcross] (1,1) to[out=-90, in=90] (2,0);
        \draw[knot] (0,1) -- (0,0);
    \end{scope}
    \begin{scope}[yshift=-2cm]
        \draw[knot] (3,1) to[out=-90, in=90] (2,0);
        \draw[knot, overcross] (2,1) to[out=-90, in=90] (3,0);
        \draw[knot] (1,1) -- (1,0);
        \draw[knot] (0,1) -- (0,0);
    \end{scope}
    \end{scope}
    \begin{scope}[yshift=-9cm]
        \draw[knot] (3,1) -- (3,0);
        \draw[knot] (2,1) -- (2,0);
        \draw[knot] (1,1) to[out=-90, in=90] (0,0);
        \draw[knot, overcross] (0,1) to[out=-90, in=90] (1,0);
    \begin{scope}[yshift=-1cm]
        \draw[knot] (3,1) -- (3,0);
        \draw[knot] (2,1) to[out=-90, in=90] (1,0);
        \draw[knot, overcross] (1,1) to[out=-90, in=90] (2,0);
        \draw[knot] (0,1) -- (0,0);
    \end{scope}
    \begin{scope}[yshift=-2cm]
        \draw[knot] (3,1) to[out=-90, in=90] (2,0);
        \draw[knot, overcross] (2,1) to[out=-90, in=90] (3,0);
        \draw[knot] (1,1) -- (1,0);
        \draw[knot] (0,1) -- (0,0);
        \draw[knot] (0,0) to[out=-90, in=-90] (1,0);
        \draw[knot] (0, -1) to[out=90, in=90] (1, -1);
        \draw[knot] (2,0) -- (2,-1);
        \draw[knot] (3,0) -- (3,-1);
    \end{scope}
    \end{scope}
}$};
    \node(B) at (4,0) {$\tikz[baseline={([yshift=-.5ex]current bounding box.center)}, scale=0.6, yscale=0.8]
{
        \draw[knot] (3,1) -- (3,0);
        \draw[knot] (2,1) -- (2,0);
        \draw[knot] (1,1) to[out=-90, in=90] (0,0);
        \draw[knot, overcross] (0,1) to[out=-90, in=90] (1,0);
    \begin{scope}[yshift=-1cm]
        \draw[knot] (3,1) -- (3,0);
        \draw[knot] (2,1) to[out=-90, in=90] (1,0);
        \draw[knot, overcross] (1,1) to[out=-90, in=90] (2,0);
        \draw[knot] (0,1) -- (0,0);
    \end{scope}
    \begin{scope}[yshift=-2cm]
        \draw[knot] (3,1) to[out=-90, in=90] (2,0);
        \draw[knot, overcross] (2,1) to[out=-90, in=90] (3,0);
        \draw[knot] (1,1) -- (1,0);
        \draw[knot] (0,1) -- (0,0);
    \end{scope}
    \begin{scope}[yshift=-3cm]
        \draw[knot] (3,1) -- (3,0);
        \draw[knot] (2,1) -- (2,0);
        \draw[knot] (1,1) to[out=-90, in=90] (0,0);
        \draw[knot, overcross] (0,1) to[out=-90, in=90] (1,0);
    \begin{scope}[yshift=-1cm]
        \draw[knot] (3,1) -- (3,0);
        \draw[knot] (2,1) to[out=-90, in=90] (1,0);
        \draw[knot, overcross] (1,1) to[out=-90, in=90] (2,0);
        \draw[knot] (0,1) -- (0,0);
    \end{scope}
    \begin{scope}[yshift=-2cm]
        \draw[knot] (3,1) to[out=-90, in=90] (2,0);
        \draw[knot, overcross] (2,1) to[out=-90, in=90] (3,0);
        \draw[knot] (1,1) -- (1,0);
        \draw[knot] (0,1) -- (0,0);
    \end{scope}
    \end{scope}
    \begin{scope}[yshift=-6cm]
        \draw[knot] (3,1) -- (3,0);
        \draw[knot] (2,1) -- (2,0);
        \draw[knot] (1,1) to[out=-90, in=90] (0,0);
        \draw[knot, overcross] (0,1) to[out=-90, in=90] (1,0);
    \begin{scope}[yshift=-1cm]
        \draw[knot] (3,1) -- (3,0);
        \draw[knot] (2,1) to[out=-90, in=90] (1,0);
        \draw[knot, overcross] (1,1) to[out=-90, in=90] (2,0);
        \draw[knot] (0,1) -- (0,0);
    \end{scope}
    \begin{scope}[yshift=-2cm]
        \draw[knot] (3,1) to[out=-90, in=90] (2,0);
        \draw[knot, overcross] (2,1) to[out=-90, in=90] (3,0);
        \draw[knot] (1,1) -- (1,0);
        \draw[knot] (0,1) -- (0,0);
    \end{scope}
    \end{scope}
    \begin{scope}[yshift=-9cm]
        \draw[knot, red] (1.5, 0.7) -- (1.5, -0.525);
        \draw[knot] (3,1) -- (3,0);
        \draw[knot, overcross] (0,1) to[out=-90, in=90] (1,0);
        \draw[knot] (1,1) to[out=-90, in=-90] (2,1);
    \begin{scope}[yshift=-1cm]
        \draw[knot] (3,1) -- (3,0);
        \draw[knot] (1,1) to[out=-90, in=90] (2,0);
    \end{scope}
    \begin{scope}[yshift=-2cm]
        \draw[knot] (3,1) to[out=-90, in=90] (2,0);
        \draw[knot, overcross] (2,1) to[out=-90, in=90] (3,0);
        \draw[knot] (0, -1) to[out=90, in=90] (1, -1);
        \draw[knot] (2,0) -- (2,-1);
        \draw[knot] (3,0) -- (3,-1);
    \end{scope}
    \end{scope}
}$};
    \node(C) at (8,0) {$\tikz[baseline={([yshift=-.5ex]current bounding box.center)}, scale=0.6, yscale=0.8]
{
        \draw[knot] (3,1) -- (3,0);
        \draw[knot] (2,1) -- (2,0);
        \draw[knot] (1,1) to[out=-90, in=90] (0,0);
        \draw[knot, overcross] (0,1) to[out=-90, in=90] (1,0);
    \begin{scope}[yshift=-1cm]
        \draw[knot] (3,1) -- (3,0);
        \draw[knot] (2,1) to[out=-90, in=90] (1,0);
        \draw[knot, overcross] (1,1) to[out=-90, in=90] (2,0);
        \draw[knot] (0,1) -- (0,0);
    \end{scope}
    \begin{scope}[yshift=-2cm]
        \draw[knot] (3,1) to[out=-90, in=90] (2,0);
        \draw[knot, overcross] (2,1) to[out=-90, in=90] (3,0);
        \draw[knot] (1,1) -- (1,0);
        \draw[knot] (0,1) -- (0,0);
    \end{scope}
    \begin{scope}[yshift=-3cm]
        \draw[knot] (3,1) -- (3,0);
        \draw[knot] (2,1) -- (2,0);
        \draw[knot] (1,1) to[out=-90, in=90] (0,0);
        \draw[knot, overcross] (0,1) to[out=-90, in=90] (1,0);
    \begin{scope}[yshift=-1cm]
        \draw[knot] (3,1) -- (3,0);
        \draw[knot] (2,1) to[out=-90, in=90] (1,0);
        \draw[knot, overcross] (1,1) to[out=-90, in=90] (2,0);
        \draw[knot] (0,1) -- (0,0);
    \end{scope}
    \begin{scope}[yshift=-2cm]
        \draw[knot] (3,1) to[out=-90, in=90] (2,0);
        \draw[knot, overcross] (2,1) to[out=-90, in=90] (3,0);
        \draw[knot] (1,1) -- (1,0);
        \draw[knot] (0,1) -- (0,0);
    \end{scope}
    \end{scope}
    \begin{scope}[yshift=-6cm]
        \draw[knot] (3,1) -- (3,0);
        \draw[knot] (2,1) -- (2,0);
        \draw[knot] (1,1) to[out=-90, in=90] (0,0);
        \draw[knot, overcross] (0,1) to[out=-90, in=90] (1,0);
    \begin{scope}[yshift=-1cm]
        \draw[knot, overcross] (1,1) to[out=-90, in=90] (2,0);
        \draw[knot] (0,1) -- (0,0);
        \draw[knot, red] ;
        \draw[knot] (2,1) to[out=-90, in=-90] (3,1);   \draw[knot, red] (2.5, 0.7) -- (2.5, -0.525);
        \draw[knot, red] (2.3, 0.72) to[out=-90, in=30] (2.1,0.1);
        \draw[knot, red, rounded corners=0.2em] (1.9, 0) -- (1.5, -0.25) -- (1.5, -2.5);
    \end{scope}
    \begin{scope}[yshift=-2cm]
        \draw[knot] (2,1) to[out=-90, in=90] (3,0);
        \draw[knot] (0,1) -- (0,0);
    \end{scope}
    \end{scope}
    \begin{scope}[yshift=-9cm]
        \draw[knot] (3,1) -- (3,0);
        \draw[knot, overcross] (0,1) to[out=-90, in=90] (1,0);
    \begin{scope}[yshift=-1cm]
        \draw[knot] (3,1) -- (3,0);
        \draw[knot] (1,1) to[out=-90, in=90] (2,0);
    \end{scope}
    \begin{scope}[yshift=-2cm]
        \draw[knot] (3,1) to[out=-90, in=90] (2,0);
        \draw[knot, overcross] (2,1) to[out=-90, in=90] (3,0);
        \draw[knot] (0, -1) to[out=90, in=90] (1, -1);
        \draw[knot] (2,0) -- (2,-1);
        \draw[knot] (3,0) -- (3,-1);
    \end{scope}
    \end{scope}
}$};
    \node(D) at (12,0) {$\tikz[baseline={([yshift=-.5ex]current bounding box.center)}, scale=0.6, yscale=0.8]
{
        \draw[knot] (3,1) -- (3,0);
        \draw[knot] (2,1) -- (2,0);
        \draw[knot] (1,1) to[out=-90, in=90] (0,0);
        \draw[knot, overcross] (0,1) to[out=-90, in=90] (1,0);
    \begin{scope}[yshift=-1cm]
        \draw[knot] (3,1) -- (3,0);
        \draw[knot] (2,1) to[out=-90, in=90] (1,0);
        \draw[knot, overcross] (1,1) to[out=-90, in=90] (2,0);
        \draw[knot] (0,1) -- (0,0);
    \end{scope}
    \begin{scope}[yshift=-2cm]
        \draw[knot] (3,1) to[out=-90, in=90] (2,0);
        \draw[knot, overcross] (2,1) to[out=-90, in=90] (3,0);
        \draw[knot] (1,1) -- (1,0);
        \draw[knot] (0,1) -- (0,0);
    \end{scope}
    \begin{scope}[yshift=-3cm]
        \draw[knot] (3,1) -- (3,0);
        \draw[knot] (2,1) -- (2,0);
        \draw[knot] (1,1) to[out=-90, in=90] (0,0);
        \draw[knot, overcross] (0,1) to[out=-90, in=90] (1,0);
    \begin{scope}[yshift=-1cm]
        \draw[knot] (3,1) -- (3,0);
        \draw[knot] (2,1) to[out=-90, in=90] (1,0);
        \draw[knot, overcross] (1,1) to[out=-90, in=90] (2,0);
        \draw[knot] (0,1) -- (0,0);
    \end{scope}
    \begin{scope}[yshift=-2cm]
        \draw[knot] (1,1) -- (1,0);
        \draw[knot] (0,1) -- (0,0);
        \draw[knot] (2,1) to[out=-90, in=-90] (3,1);
        \draw[knot, red] (2.5, 0.7) -- (2.5, -1.3);
        \draw[knot, red] (2.3, 0.74) -- (2.3, -1.3);
    \end{scope}
    \end{scope}
    \begin{scope}[yshift=-6cm]
        \draw[knot] (1,1) to[out=-90, in=90] (0,0);
        \draw[knot, overcross] (0,1) to[out=-90, in=90] (1,0);
    \begin{scope}[yshift=-1cm]
        \draw[knot, overcross] (1,1) to[out=-90, in=90] (2,0);
        \draw[knot] (0,1) -- (0,0);
        \draw[knot, red] (2.5, 0.7) -- (2.5, -0.525);
        \draw[knot, red] (2.3, 0.74) to[out=-90, in=30] (2.1,0.1);
        \draw[knot, red, rounded corners=0.2em] (1.9, 0) -- (1.5, -0.25) -- (1.5, -2.5);
    \end{scope}
    \begin{scope}[yshift=-2cm]
        \draw[knot] (2,1) to[out=-90, in=90] (3,0);
        \draw[knot] (0,1) -- (0,0);
    \end{scope}
    \end{scope}
    \begin{scope}[yshift=-9cm]
        \draw[knot] (3,1) -- (3,0);
        \draw[knot, overcross] (0,1) to[out=-90, in=90] (1,0);
    \begin{scope}[yshift=-1cm]
        \draw[knot] (3,1) -- (3,0);
        \draw[knot] (1,1) to[out=-90, in=90] (2,0);
    \end{scope}
    \begin{scope}[yshift=-2cm]
        \draw[knot] (3,1) to[out=-90, in=90] (2,0);
        \draw[knot, overcross] (2,1) to[out=-90, in=90] (3,0);
        \draw[knot] (0, -1) to[out=90, in=90] (1, -1);
        \draw[knot] (2,0) -- (2,-1);
        \draw[knot] (3,0) -- (3,-1);
    \end{scope}
    \end{scope}
}
\{-1,-1\}$};
    \node (E) at (12,-8) {$\tikz[baseline={([yshift=-.5ex]current bounding box.center)}, scale=0.6, yscale=0.8]
{
        \draw[knot] (3,1) -- (3,0);
        \draw[knot] (2,1) -- (2,0);
        \draw[knot] (1,1) to[out=-90, in=90] (0,0);
        \draw[knot, overcross] (0,1) to[out=-90, in=90] (1,0);
    \begin{scope}[yshift=-1cm]
        \draw[knot] (3,1) -- (3,0);
        \draw[knot] (2,1) to[out=-90, in=90] (1,0);
        \draw[knot, overcross] (1,1) to[out=-90, in=90] (1.65,0);
        \draw[knot] (1.65, 0) -- (1.65, -1.5) to[out=-90, in=0] (1, -2);
        \draw[knot, red] (1.65, -1) -- (2, -1);
        \draw[knot, red] (1.65, -1.325) -- (1.9, -1.325);
            \draw[knot, red ,rounded corners] (2.1, -1.325) -- (2.5, -1.325) -- (2.5, -3.3);
        \draw[knot, red] (1.65, -1.65) -- (1.9, -1.65);
            \draw[knot, red, rounded corners] (2.1, -1.65) -- (2.3, -1.65) -- (2.3, -3.26);
        \draw[knot] (0,1) -- (0,0);
    \end{scope}
    \begin{scope}[yshift=-2cm]
        \draw[knot] (3,1) to[out=-90, in=90] (2,0);
        \draw[knot] (1,1) -- (1,0);
        \draw[knot] (0,1) -- (0,0);
    \end{scope}
    \begin{scope}[yshift=-3cm]
        \draw[knot] (2,1) -- (2,0);
        \draw[knot] (1,1) to[out=-90, in=90] (0,0);
        \draw[knot, overcross] (0,1) to[out=-90, in=180] (1,0);
    \begin{scope}[yshift=-1cm]
        \draw[knot] (2,1) to[out=-90, in=90] (1,0);
        \draw[knot] (0,1) -- (0,0);
    \end{scope}
    \begin{scope}[yshift=-2cm]
        \draw[knot] (1,1) -- (1,0);
        \draw[knot] (0,1) -- (0,0);
        \draw[knot, red] (2.5, 0.7) -- (2.5, -1.3);
        \draw[knot, red] (2.3, 0.74) -- (2.3, -1.3);
    \end{scope}
    \end{scope}
    \begin{scope}[yshift=-6cm]
        \draw[knot] (1,1) to[out=-90, in=90] (0,0);
        \draw[knot, overcross] (0,1) to[out=-90, in=90] (1,0);
    \begin{scope}[yshift=-1cm]
        \draw[knot, overcross] (1,1) to[out=-90, in=90] (2,0);
        \draw[knot] (0,1) -- (0,0);
        \draw[knot, red] (2.5, 0.7) -- (2.5, -0.525);
        \draw[knot, red] (2.3, 0.74) to[out=-90, in=30] (2.1,0.1);
        \draw[knot, red, rounded corners=0.2em] (1.9, 0) -- (1.5, -0.25) -- (1.5, -2.5);
    \end{scope}
    \begin{scope}[yshift=-2cm]
        \draw[knot] (2,1) to[out=-90, in=90] (3,0);
        \draw[knot] (0,1) -- (0,0);
    \end{scope}
    \end{scope}
    \begin{scope}[yshift=-9cm]
        \draw[knot] (3,1) -- (3,0);
        \draw[knot, overcross] (0,1) to[out=-90, in=90] (1,0);
    \begin{scope}[yshift=-1cm]
        \draw[knot] (3,1) -- (3,0);
        \draw[knot] (1,1) to[out=-90, in=90] (2,0);
    \end{scope}
    \begin{scope}[yshift=-2cm]
        \draw[knot] (3,1) to[out=-90, in=90] (2,0);
        \draw[knot, overcross] (2,1) to[out=-90, in=90] (3,0);
        \draw[knot] (0, -1) to[out=90, in=90] (1, -1);
        \draw[knot] (2,0) -- (2,-1);
        \draw[knot] (3,0) -- (3,-1);
    \end{scope}
    \end{scope}
}
\{-1,-1\}$};
    \node(F) at (6.55, -8) {$\tikz[baseline={([yshift=-.5ex]current bounding box.center)}, scale=0.6, yscale=0.8]
{
        \draw[knot] (3,1) -- (3,0);
        \draw[knot] (2,1) -- (2,0);
        \draw[knot] (1,1) to[out=-90, in=90] (0,0);
        \draw[knot, overcross] (0,1) to[out=-90, in=90] (1,0);
        \draw[knot] (0,0) to[out=-90, in=-90] (1,0);
    \begin{scope}[yshift=-1cm]
        \draw[knot] (3,1) -- (3,0);
        \draw[knot] (2,1) to[out=-90, in=90] (1,0);
    \draw[knot, red, rounded corners=0.2em] (0.6, 0.725) -- (0.6, -1) -- (0.9,-1);
    \draw[knot, red] (1.1, -1) -- (1.65, -1);
        \draw[knot, red] (1.65, -1) -- (2, -1);
    \draw[knot, red] (0.95, -1.325) -- (1.65, -1.325);
    \draw[knot, red, rounded corners=0.2em] (0.7, -1.325) -- (0.4, -1.325) -- (0.4, 0.725);
        \draw[knot, red] (1.65, -1.325) -- (1.9, -1.325);
            \draw[knot, red ,rounded corners] (2.1, -1.325) -- (2.5, -1.325) -- (2.5, -3.3);
    \draw[knot, red, rounded corners=0.2em] (1.65, -1.65) -- (0.5, -1.65) -- (0.425, -1.625);
    \draw[knot, red, rounded corners=0.2em] (0.2, 0.75) -- (0.2, -1.4) -- (0.3, -1.5);
        \draw[knot, red] (1.65, -1.65) -- (1.9, -1.65);
            \draw[knot, red, rounded corners] (2.1, -1.65) -- (2.3, -1.65) -- (2.3, -3.26);
        \draw[knot, red, rounded corners=0.2em] (0.8, 0.75) -- (0.8, -0.675) -- (1, -0.675);
    \end{scope}
    \begin{scope}[yshift=-2cm]
        \draw[knot] (3,1) to[out=-90, in=90] (2,0);
        \draw[knot] (1,1) -- (1,0);
    \end{scope}
    \begin{scope}[yshift=-3cm]
        \draw[knot] (2,1) -- (2,0);
        \draw[knot] (1,1) to[out=-90, in=90] (0,0);
    \begin{scope}[yshift=-1cm]
        \draw[knot] (2,1) to[out=-90, in=90] (1,0);
        \draw[knot] (0,1) -- (0,0);
    \end{scope}
    \begin{scope}[yshift=-2cm]
        \draw[knot] (1,1) -- (1,0);
        \draw[knot] (0,1) -- (0,0);
        \draw[knot, red] (2.5, 0.7) -- (2.5, -1.3);
        \draw[knot, red] (2.3, 0.74) -- (2.3, -1.3);
    \end{scope}
    \end{scope}
    \begin{scope}[yshift=-6cm]
        \draw[knot] (1,1) to[out=-90, in=90] (0,0);
        \draw[knot, overcross] (0,1) to[out=-90, in=90] (1,0);
    \begin{scope}[yshift=-1cm]
        \draw[knot, overcross] (1,1) to[out=-90, in=90] (2,0);
        \draw[knot] (0,1) -- (0,0);
        \draw[knot, red] (2.5, 0.7) -- (2.5, -0.525);
        \draw[knot, red] (2.3, 0.74) to[out=-90, in=30] (2.1,0.1);
        \draw[knot, red, rounded corners=0.2em] (1.9, 0) -- (1.5, -0.25) -- (1.5, -2.5);
    \end{scope}
    \begin{scope}[yshift=-2cm]
        \draw[knot] (2,1) to[out=-90, in=90] (3,0);
        \draw[knot] (0,1) -- (0,0);
    \end{scope}
    \end{scope}
    \begin{scope}[yshift=-9cm]
        \draw[knot] (3,1) -- (3,0);
        \draw[knot, overcross] (0,1) to[out=-90, in=90] (1,0);
    \begin{scope}[yshift=-1cm]
        \draw[knot] (3,1) -- (3,0);
        \draw[knot] (1,1) to[out=-90, in=90] (2,0);
    \end{scope}
    \begin{scope}[yshift=-2cm]
        \draw[knot] (3,1) to[out=-90, in=90] (2,0);
        \draw[knot, overcross] (2,1) to[out=-90, in=90] (3,0);
        \draw[knot] (0, -1) to[out=90, in=90] (1, -1);
        \draw[knot] (2,0) -- (2,-1);
        \draw[knot] (3,0) -- (3,-1);
    \end{scope}
    \end{scope}
}
\{-1,-1\}$};
    \node(G) at (0.75, -8) {$\tikz[baseline={([yshift=-.5ex]current bounding box.center)}, scale=0.6, yscale=0.8]
{
        \draw[knot] (3,1) -- (3,0);
        \draw[knot] (2,1) -- (2,0);
        \draw[knot] (0,1) to[out=-90, in=-90] (1,1);
        \draw[knot, red] (0.2, 0.75) -- (0.2, -0.25);
        \draw[knot, red] (0.4, 0.72) -- (0.4, -0.28);
        \draw[knot, red] (0.6, 0.72) -- (0.6, -0.28);
        \draw[knot, red] (0.8, 0.75) -- (0.8, -0.25);
    \begin{scope}[yshift=-1cm]
        \draw[knot] (3,1) -- (3,0);
        \draw[knot] (2,1) to[out=-90, in=90] (1,0);
    \draw[knot, red, rounded corners=0.2em] (0.6, 0.725) -- (0.6, -1) -- (0.9,-1);
    \draw[knot, red] (1.1, -1) -- (1.65, -1);
        \draw[knot, red] (1.65, -1) -- (2, -1);
    \draw[knot, red] (0.95, -1.325) -- (1.65, -1.325);
    \draw[knot, red, rounded corners=0.2em] (0.7, -1.325) -- (0.4, -1.325) -- (0.4, 0.725);
        \draw[knot, red] (1.65, -1.325) -- (1.9, -1.325);
            \draw[knot, red ,rounded corners] (2.1, -1.325) -- (2.5, -1.325) -- (2.5, -3.3);
    \draw[knot, red, rounded corners=0.2em] (1.65, -1.65) -- (0.5, -1.65) -- (0.425, -1.625);
    \draw[knot, red, rounded corners=0.2em] (0.2, 0.75) -- (0.2, -1.4) -- (0.3, -1.5);
        \draw[knot, red] (1.65, -1.65) -- (1.9, -1.65);
            \draw[knot, red, rounded corners] (2.1, -1.65) -- (2.3, -1.65) -- (2.3, -3.26);
        \draw[knot, red, rounded corners=0.2em] (0.8, 0.75) -- (0.8, -0.675) -- (1, -0.675);
    \end{scope}
    \begin{scope}[yshift=-2cm]
        \draw[knot] (3,1) to[out=-90, in=90] (2,0);
        \draw[knot] (1,1) -- (1,0);
    \end{scope}
    \begin{scope}[yshift=-3cm]
        \draw[knot] (2,1) -- (2,0);
        \draw[knot] (1,1) to[out=-90, in=90] (0,0);
    \begin{scope}[yshift=-1cm]
        \draw[knot] (2,1) to[out=-90, in=90] (1,0);
        \draw[knot] (0,1) -- (0,0);
    \end{scope}
    \begin{scope}[yshift=-2cm]
        \draw[knot] (1,1) -- (1,0);
        \draw[knot] (0,1) -- (0,0);
        \draw[knot, red] (2.5, 0.7) -- (2.5, -1.3);
        \draw[knot, red] (2.3, 0.74) -- (2.3, -1.3);
    \end{scope}
    \end{scope}
    \begin{scope}[yshift=-6cm]
        \draw[knot] (1,1) to[out=-90, in=90] (0,0);
        \draw[knot, overcross] (0,1) to[out=-90, in=90] (1,0);
    \begin{scope}[yshift=-1cm]
        \draw[knot, overcross] (1,1) to[out=-90, in=90] (2,0);
        \draw[knot] (0,1) -- (0,0);
        \draw[knot, red] (2.5, 0.7) -- (2.5, -0.525);
        \draw[knot, red] (2.3, 0.74) to[out=-90, in=30] (2.1,0.1);
        \draw[knot, red, rounded corners=0.2em] (1.9, 0) -- (1.5, -0.25) -- (1.5, -2.5);
    \end{scope}
    \begin{scope}[yshift=-2cm]
        \draw[knot] (2,1) to[out=-90, in=90] (3,0);
        \draw[knot] (0,1) -- (0,0);
    \end{scope}
    \end{scope}
    \begin{scope}[yshift=-9cm]
        \draw[knot] (3,1) -- (3,0);
        \draw[knot, overcross] (0,1) to[out=-90, in=90] (1,0);
    \begin{scope}[yshift=-1cm]
        \draw[knot] (3,1) -- (3,0);
        \draw[knot] (1,1) to[out=-90, in=90] (2,0);
    \end{scope}
    \begin{scope}[yshift=-2cm]
        \draw[knot] (3,1) to[out=-90, in=90] (2,0);
        \draw[knot, overcross] (2,1) to[out=-90, in=90] (3,0);
        \draw[knot] (0, -1) to[out=90, in=90] (1, -1);
        \draw[knot] (2,0) -- (2,-1);
        \draw[knot] (3,0) -- (3,-1);
    \end{scope}
    \end{scope}
}
\{-2,-2\}$};
    \draw[->] (A) -- node[pos=.5,above,arrows=-]
        {RII}  (B);
    \draw[->] (B) -- node[pos=.5,above,arrows=-]
        {RII}  (C);
    \draw[->] (C) -- node[pos=.5,above,arrows=-]
        {RI}  (D);
    \draw[->] (D) -- node[pos=.5,right,arrows=-]
        {RII}  (E);
    \draw[->] (E) -- node[pos=.5,above,arrows=-]
        {RII}  (F);
    \draw[->] (F) -- node[pos=.5,above,arrows=-]
        {RI}  (G);
}
\]
\caption{Computing the grading shift on $\mathcal{T}_4^4 \otimes e_1$.}
\label{fig:fulltwistGS}
\end{figure}

We'll write $e_i$ as $e_i^{\mathrm{top}} \otimes e_i^{\mathrm{bot}}$, although this tensor product is not exactly the same as the one in Definition \ref{def:verttensor}; we do not belabor the point. Assume that $r=0$. Then $\mathcal{T}_n^{nk}$ is $k$-full twists, and we have that
\[
\mathcal{T}_n^{nk} \otimes e_i = 
e_{i'}^{\mathrm{top}} \otimes \varphi_{W_n^{nk}} \mathcal{T}_{n-2}^{(n-2)k} \{-2k, -2k\} \otimes e_i^{\mathrm{bot}}
\]
where $W_{n}^{nk}$ is a cobordism consisting of $2k(n-2)$ saddles (for the $2k(n-2)$-many Reidemeister II moves performed) and $i' = i+r \mod n$. There are also $2k$ Reidemeister I moves, accounting for the $\mathbb{Z}\times \mathbb{Z}$-shift. To aid in comprehending $\varphi_{W_n^{nk}}$, consider Figure \ref{fig:fulltwistGS}. We remark that the tensor on the left is vertical stacking as in definition \ref{def:verttensor}, and the one on the right is as in the writing of $e_i^{\mathrm{top}} \otimes e_i^{\mathrm{bot}}$. Notice that $e_0^{\mathrm{top}}$ is allowed; by this we mean the following picture.
\[
e_0^{\mathrm{top}} :=
\tikz[baseline={([yshift=-.5ex]current bounding box.center)}, scale=0.6, yscale=0.8]{
    \draw[knot] (0.7,1) -- (0.7,-0.75);
    \draw[knot] (2.3,1) -- (2.3,-0.75);
    \node at (1.5, 0.75) {$\cdots$};
    \draw[knot, overcross] (0,1) to[out=-90, in=-90] (3,1);
}
\]

Now, for $1\le r < n$, there are three cases.
\begin{enumerate}
    \item If $i< n-r$, the extra isotopy contains no Reidemeister I moves, but it does consist of $r$-many Reidemeister II moves.
    \item If $i=n-r$, the isotopy contains $(r-1)$ more Reidemeister II moves and exactly 1 more Reidemeister I move. Note that $i'=0$ in this case.
    \item If $i>n-r$, the isotopy contains the addition of a sequence of $(r-2)$ many Reidemeister II moves, then 1 Reidemeister I move, followed by $(n-2)$ more Reidemeister II moves, and another lone Reidemeister I move; that is, $(n+r-4)$ Reidemeister II moves and 2 Reidemeister I moves.
\end{enumerate}
So, we have proven the following.

\begin{lemma}
\label{lem:toplemma}
For any $0\le r < n$ and $0 < i < n$, 
\[
    \mathcal{T}_n^{nk+r} \otimes e_i \simeq e_{i'}^{\mathrm{top}} \otimes \varphi_{W_n^{nk+r}} \mathcal{T}_{n-2}^{(n-2)k + r_i} \{-(2k+k_i), -(2k+k_i)\} \otimes e_i^{\mathrm{bot}}
\]
where $i' = i+r \mod n$ and
\begin{enumerate}
    \item if $i<n-r$, $W_n^{nk+r}$ consists of $2k(n-2) + r$ saddles, $r_i = r$, and $k_i = 0$;
    \item if $i = n-r$, $W_n^{nk+r}$ consists of $2k(n-2) + (r - 1)$ saddles, $r_i = r-1$, and $k_i = 1$;
    \item if $i>n-r$, $W_n^{nk+r}$ consists of $2k(n-2) + (n+r-4)$ saddles, $r_i = r-2$, and $k_i=2$.
\end{enumerate}
In each of these cases, $W_n^{nk+r}$ is a cobordism in the style of Figure \ref{fig:fulltwistGS}.
\end{lemma}

We'll denote by $s_i$ the number of additional saddles depending on $r$. That is, $W_{n}^{nk+r}$ consists of $2k(n-2) + s_i$ saddles, where 
\begin{enumerate}
    \item $s_i=r$ if $i < n-r$;
    \item $s_i = r-1$ if $i = n-r$;
    \item $s_i = n+r-4$ if $i > n-r$.
\end{enumerate}
Note that our $s_i$ is not the same as the one appearing in \cite{stoffregen2024joneswenzlprojectorskhovanovhomotopy}.

We will use this Lemma to prove the existence of projectors. First, we would like to draw some connections between our work and computations found in Section 5 of \cite{stoffregen2024joneswenzlprojectorskhovanovhomotopy}. Consider the complex $\mathcal{C}_{m+1}$ defined as the cone
\[
\mathcal{C}_{m+1} := \mathrm{Cone}(\mathcal{T}_n^m \to \mathcal{T}_n^{m+1}).
\]
Then, $\mathcal{C}_{m+1}$ looks like (that is, is homotopy equivalent to) a cube of resolutions for $\mathcal{T}_n^1$ with $\mathcal{T}_n^m$ stacked on top, modulo the identity term, which is taken to be zero. Any entry of the cube of resolutions for $\mathcal{T}_n$ (apart from the identity entry, which we have avoided) is isomorphic to $\mathcal{F}(e_i) \otimes \mathcal{F}(\delta)$ for some flat diskular $n$-tangle $\delta$ and $1\le i \le n-1$. Dropping the $\mathcal{F}$ notation, this is to say that $\mathcal{C}_{m+1}$ is homotopy equivalent to a colimit in which all nontrivial terms are of the form
\[
\varphi_{\alpha_n^1} \mathcal{T}_n^{m} \otimes e_i \otimes \delta
\]
where $\varphi_{\alpha_n^1}$ denotes the grading shift coming from the cube of resolutions for $\mathcal{T}_n^1$. Writing $m = nk +r$, Lemma \ref{lem:toplemma} says that this term is equivalent to 
\begin{equation}
\label{eq:inftwistterms}
\varphi_{\alpha_n^1} \left(e_{i'}^{\mathrm{top}} \otimes \varphi_{W_n^{nk+r}} \mathcal{T}_{n-2}^{(n-2)k + r_i} \{-(2k+k_i),-(2k+k_i)\} \otimes e_i^{\mathrm{bot}}\right) \otimes \delta.
\end{equation}

As in \cite{stoffregen2024joneswenzlprojectorskhovanovhomotopy}, we want to provide a bound on grading shifts. On one hand, given a $\mathscr{G}$-graded dg $H^n$-module $A$, by a \textit{global upper $q$-bound on $\mathscr{G}$-grading shifts}, we mean some $B \in \mathbb{Z}$ so that, for each entry $A_i$ of $A$ with grading shift $\varphi_{W^v}$, $\deg_q(\varphi_{W^v})\le B$. For example, we can compute an upper bound of a complex with $\mathscr{G}$-grading finding the minimum number of saddles appearing in each grading shift and maximizing the $\mathbb{Z}\times \mathbb{Z}$-degree. We define global lower bounds similarly. This definition extends to a stricter notion on objects of $\mathrm{Kom}(H^n\mathrm{Mod}^\mathscr{G})$ by taking the minimum (resp. maximum) among all global upper (resp. global lower) bounds for each complex $A'$ homotopy equivalent to $A$.

Referring again to Proposition \ref{Prop:cofibration}, to any diskular tangle $T$, $\mathcal{F}(T)$ has an entry with trivial $\mathscr{G}$-grading; this is to say that a global upper bound on $\mathcal{F}(T)$ is 0. Similarly, a global lower bound is given by $-c(T)$, for $c(T)$ the number of crossings in the diagram for $T$. 

Note that $\varphi_{\alpha_n^1}$ always consists of at least one saddle, by construction. Then, we can compute the $q$-grading shift on (\ref{eq:inftwistterms}) on a case-by-case basis via Lemma \ref{lem:toplemma} and conclude that $\mathcal{C}_{m+1}$ is homotopy equivalent to a complex with global upper bound on $\mathscr{G}$-grading
\[
b_\epsilon \le B_{m+1} := -2nk-r-1.
\]
Observe that this bound is similar to the one provided in \cite{stoffregen2024joneswenzlprojectorskhovanovhomotopy}.

\begin{remark}
As in \cite{stoffregen2024joneswenzlprojectorskhovanovhomotopy}, we can present a model in which $\mathcal{T}_n^\infty$ is an iterated mapping cone. Start by setting $\mathcal{A}^1 = \mathcal{T}_n^1$ and, inductively, assume $\mathcal{A}^2, \ldots, \mathcal{A}^m$ have been constructed, each satisfying $\mathcal{A}^\ell \simeq \mathcal{T}_n^\ell$. We construct $\mathcal{A}^{m+1}$ as follows. From the definition of $\mathcal{C}_{m+1}$, there is an exact triangle
\[
\begin{tikzcd}
    \mathcal{T}_n^m \arrow[r] & \mathcal{T}_n^{m+1} \arrow[r] & \mathcal{C}_{m+1},
\end{tikzcd}
\]
thus there is a map $\psi_m$ so that $\mathcal{T}_n^{m+1} \simeq \mathrm{Cone}(\mathcal{C}_{m+1} \xrightarrow{\psi_m} \mathcal{T}_n^m)$.

Now, using Lemma \ref{lem:toplemma}, we have argued that $\mathcal{C}_{m+1}$ is homotopy equivalent to a complex we'll call $\mathcal{C}_{m+1}'$ with glabal upper bound $B_{m+1}$. Let $\psi_m'$ denote the map defined by the commutative square
\[
\begin{tikzcd}
    \mathcal{C}_{m+1} \arrow[r, "\psi_m"] & \mathcal{T}_n^m \\
    \mathcal{C}_{m+1}' \arrow[u, "\sim" {anchor=south, rotate=90}] \arrow[r, "\psi_m'"] & \mathcal{A}^m \arrow[u, "\sim"'{anchor=north, rotate=90}] 
\end{tikzcd}
\]
where each vertical arrow is a homotopy equivalence. Then, set
\[
\mathcal{A}^{m+1} := \mathrm{Cone}(\mathcal{C}_{m+1}' \xrightarrow{\psi_m'} \mathcal{A}^m).
\]
Unfurling definitions and homotopy equivalences, it follows that $\mathcal{A}^{m+1} \simeq \mathcal{T}_n^{m+1}$.

In particular, $\mathcal{A}^{m+1}$ is obtained from $\mathcal{A}^m$ by including finitely many new entries with $\mathscr{G}$-grading shifts bounded from above by $B_{m+1}$. As $m \to \infty$, $B_{m+1} \to -\infty$, and we obtain a model for $\mathcal{T}_n^\infty \simeq \mathcal{A}^\infty$ as an iterated mapping cone.
\end{remark}

On the other hand, define a \textit{global upper $\mathbb{Z} \times \mathbb{Z}$-bound on $\mathscr{G}$-grading shifts} to be some $(B_1, B_2) \in \mathbb{Z} \times \mathbb{Z}$ so that, for each $A_i$ of $A$ with grading shift $\varphi_{W^{(v_1, v_2)}}$, we can find a simplification of $\varphi_{W^{(v_1, v_2)}}$, written $\varphi_{\check{W}^{(v_1', v_2')}}$, in which $v_1' \le B_1$ and $v_2' \le B_2$. By a simplification, we that $\varphi_{W^{(v_1, v_2)}} \cong \varphi_{\check{W}^{(v_1', v_2')}}$ for $\check{W}$ a minimal cobordism void of births, deaths, and unambiguous saddles.

Notice that, since $\varphi_{W_n^{nk+r}}$ consists only of saddles, we have that $(-2k, -2k)$ provides a global upper $\mathbb{Z}\times\mathbb{Z}$-bound on $\mathscr{G}$-grading shifts for a complex homotopy equivalent to $\mathcal{T}_n^{nk+r} \otimes e_i$.

\begin{theorem}
\label{thm:existenceofprojectors}
For each $n$, $\mathcal{T}_n^\infty$ is a unified projector.
\end{theorem}

\begin{proof}
Recall that $\mathcal{T}_n^\infty$ is defined as the colimit
\[
\mathcal{T}_n^\infty = \mathrm{colim} \left(\mathcal{T}_n^0 \to \mathcal{T}_n^1 \to \cdots \mathcal{T}_n^m \to \cdots \right)
\]
which we'll write $\mathrm{colim}(\mathcal{T}_n^{nk+r})$. Axiom (CK1) is apparent by definition, so we will content ourselves with a proof of (CK2). First, notice that
\[
\mathrm{colim}(\mathcal{T}_n^{nk+r}) \otimes e_i \simeq \mathrm{colim}(\mathcal{T}_n^{nk+r} \otimes e_i)
\]
so if the homology of the colimit on the right-hand side is trivial, we can conclude that the colimit itself is contractible, thus $\mathcal{T}_n^\infty \simeq *$. Recall that any homology class of the colimit arises as a homology class in a piece of the colimit. However, by Lemma \ref{lem:toplemma}, this colimit is built from complexes with a global upper $\mathbb{Z}\times\mathbb{Z}$-bound of $(-2k, -2k)$. As $m\to \infty$, $k\to \infty$, and the global upper bound goes to $(-\infty, -\infty)$, so any nontrivial homology class must die in the colimit.
\end{proof}

\subsection{A unified colored link homology}
\label{ss:unifiedcoloredhomology}

With very little work, the existence of unified projectors (together with multigluing) implies the existence of a  unified colored link homology specializing to an even one (\cite{https://doi.org/10.48550/arxiv.1005.5117}, see also \cite{khovanov2003categorificationscoloredjonespolynomial, beliakova2006categorificationcoloredjonespolynomial} by way of \cite{beliakova2023unificationcoloredannularsl2}), but also specializing to a new odd version. Recall the following definition, adapted from Definition 5.1 of \cite{https://doi.org/10.48550/arxiv.1005.5117}.
\begin{definition}
For any $n\in \mathbb{N}$ and $\mathbf{m} = (m_1,\ldots, m_n) \in \mathbb{N}^n$, we denote by $\mathrm{Chom}_{\mathbf{m}}(n)^\mathscr{G}$ the category where
\begin{itemize}
    \item $\mathrm{ob}(\mathrm{Chom}_{\mathbf{m}}(n)^\mathscr{G}) = \mathrm{ob}(\mathrm{Chom}(n)^\mathscr{G})$ and
    \item $\mathrm{Hom}_{\mathrm{Chom}_{\mathbf{m}}(n)^\mathscr{G}} (A, B) = \mathrm{Hom}_{\mathrm{Chom}(Mn)^\mathscr{G}} (\Pi^{\mathbf{m}}(A), \Pi^{\mathbf{m}}(B))$
\end{itemize}
where $M = \sum_i m_i$ and  $\Pi^{\mathbf{m}}$ replaces the $i$th strand in each diagram with its $m_i$th parallel composed with a copy of the $m_i$th projector. We define $\mathrm{Chom}_{\mathbf{m}}(n)^q$ by taking objects and morphisms of $\mathrm{Chom}_{\mathbf{m}}(n)^\mathscr{G}$ and collapsing degree, as usual.
\end{definition}

We will represent projectors by small boxes, e.g., $P_n = 
\tikz[baseline={([yshift=-.5ex]current bounding box.center)}, scale=0.6, yscale=0.8]{
    \draw (0,0.5) -- (0,2.5);
    \draw[fill=white] (-0.6,1) rectangle (0.6,2);
    \node at (0,1.5) {$n$};
}$~. We will define the operation $\Pi^m$ on links, via operations on diskular tangles, as follows. As an example, if $K$ is a knot, let $\mathring{K}$ denote the diskular 1-tangle 
$\tikz[baseline={([yshift=-.5ex]current bounding box.center)}, scale=0.6, yscale=0.95]{
    \draw[dotted, rounded corners] (0,0) rectangle (1.5, 1.5);
    \draw (0.75, 0) -- (0.75, 1.5);
    \draw[fill=white] (0.375, 0.375) rectangle (1.125,1.125);
    \node at (0.75, 0.75) {$K$};
    \node at (1.5, 0.75) {$\times$};
}$ 
and suppose $\mathring{K}^m$ denotes its $m$th parallel. Then
\[
\Pi^m(K) = \mathrm{Tr}^m(\mathrm{Kh}_q(\mathring{K}^m) \otimes P_m)
\]
More generally, if $L$ is an $n$-component link, we use multigluing. Let $\mathbf{m} = (m_1,\ldots, m_n) \in \mathbb{N}^n$, and denote by $T_L^\mathbf{m}$ the result of taking $m_i$ parallel copies of the $i$th component of $L$ and then removing a small diskular region from each of the original components (again, see Figure \ref{fig:coloredcomplex}). Then, set
\[
\Pi^\mathbf{m}(L) := (P_{m_1}, \ldots, P_{m_n}) \otimes_{(H^{m_1}, \ldots, H^{m_n})} \mathrm{Kh}_q(T_L^{\mathbf{m}})
\]
where each of the $P_{m_i}$ is viewed as an object of $\mathrm{Chom}(m_i)_R^q$.

\begin{lemma}
\label{lem:projsliding}
We have the following isomorphisms in $\mathrm{Kom}(H^{m+n}\mathrm{Mod})^\mathscr{G}$:
\[
\tikz[baseline={([yshift=-.5ex]current bounding box.center)}, scale=1]
{
    \begin{scope}[rotate=-45]
        \draw[dotted] (.5,.5) circle(0.707);
        \draw[knot] (0.5, -0.207) -- (0.5, 1.207);
        \draw[knot, overcross] (-0.207, 0.5) -- (1.207, 0.5);
        \draw[fill=white] (0.3, 0.7535) rectangle (0.7, 0.9535);
        \node[transform shape] at (0.5, 0.8535) {\small $m$};
    \end{scope}
}
~\cong~
\tikz[baseline={([yshift=-.5ex]current bounding box.center)}, scale=1]
{
    \begin{scope}[rotate=-45]
        \draw[dotted] (.5,.5) circle(0.707);
        \draw[knot] (0.5, -0.207) -- (0.5, 1.207);
        \draw[knot, overcross] (-0.207, 0.5) -- (1.207, 0.5);
        \draw[fill=white] (0.3, 0.7535) rectangle (0.7, 0.9535);
        \node[transform shape] at (0.5, 0.8535) {\small $m$};
        \draw[fill=white] (0.3, 0.0465) rectangle (0.7, 0.2465);
        \node[transform shape] at (0.5, 0.1465) {$m$};
    \end{scope}
}
~\cong~
\tikz[baseline={([yshift=-.5ex]current bounding box.center)}, scale=1]
{
    \begin{scope}[rotate=-45]
        \draw[dotted] (.5,.5) circle(0.707);
        \draw[knot] (0.5, -0.207) -- (0.5, 1.207);
        \draw[knot, overcross] (-0.207, 0.5) -- (1.207, 0.5);
        \draw[fill=white] (0.3, 0.0465) rectangle (0.7, 0.2465);
        \node[transform shape] at (0.5, 0.1465) {$m$};
    \end{scope}
}
\]
and
\[
\tikz[baseline={([yshift=-.5ex]current bounding box.center)}, scale=1]
{
    \begin{scope}[rotate=-45]
        \draw[dotted] (.5,.5) circle(0.707);
        \draw[knot] (0.5, -0.207) -- (0.5, 1.207);
        \draw[knot, overcross] (-0.207, 0.5) -- (1.207, 0.5);
        \draw[fill=white] (0.0465, 0.3) rectangle (0.2465, 0.7);
        \node[transform shape, rotate=90] at (0.1465, 0.5) {$n$};
    \end{scope}
}
~\cong~
\tikz[baseline={([yshift=-.5ex]current bounding box.center)}, scale=1]
{
    \begin{scope}[rotate=-45]
        \draw[dotted] (.5,.5) circle(0.707);
        \draw[knot] (0.5, -0.207) -- (0.5, 1.207);
        \draw[knot, overcross] (-0.207, 0.5) -- (1.207, 0.5);
        \draw[fill=white] (0.753, 0.3) rectangle (0.9535, 0.7);
        \node[transform shape, rotate=90] at (0.8535, 0.5) {\small $n$};
        \draw[fill=white] (0.0465, 0.3) rectangle (0.2465, 0.7);
        \node[transform shape, rotate=90] at (0.1465, 0.5) {$n$};
    \end{scope}
}
~\cong~
\tikz[baseline={([yshift=-.5ex]current bounding box.center)}, scale=1]
{
    \begin{scope}[rotate=-45]
        \draw[dotted] (.5,.5) circle(0.707);
        \draw[knot] (0.5, -0.207) -- (0.5, 1.207);
        \draw[knot, overcross] (-0.207, 0.5) -- (1.207, 0.5);
        \draw[fill=white] (0.753, 0.3) rectangle (0.9535, 0.7);
        \node[transform shape, rotate=90] at (0.8535, 0.5) {\small $n$};
    \end{scope}
}
\]
That is, free (parallel) strands can be moved over or under projectors in $\mathrm{Kom}(H^{2n}\mathrm{Mod})^\mathscr{G}$.
\end{lemma}

\begin{proof}
We'll explain the first homotopy equivalence; the others are proven with the same procedure. The trick is to start with the middle complex: using (CK1), $\tikz[baseline={([yshift=-.5ex]current bounding box.center)}, scale=.8]
{
    \begin{scope}[rotate=-45]
        \draw[dotted] (.5,.5) circle(0.707);
        \draw[knot] (0.5, -0.207) -- (0.5, 1.207);
        \draw[knot, overcross] (-0.207, 0.5) -- (1.207, 0.5);
        \draw[fill=white] (0.3, 0.7535) rectangle (0.7, 0.9535);
        \node[transform shape] at (0.5, 0.8535) {\small $m$};
        \draw[fill=white] (0.3, 0.0465) rectangle (0.7, 0.2465);
        \node[transform shape] at (0.5, 0.1465) {$m$};
    \end{scope}
}$ is homotopy equivalent to the complex of complexes
\[
\tikz[baseline={([yshift=-.5ex]current bounding box.center)}, scale=1]
{
    \begin{scope}[rotate=-45]
        \draw[dotted] (.5,.5) circle(0.707);
        \draw[knot] (0.5, -0.207) -- (0.5, 1.207);
        \draw[knot, overcross] (-0.207, 0.5) -- (1.207, 0.5);
        \draw[fill=white] (0.3, 0.7535) rectangle (0.7, 0.9535);
        \node[transform shape] at (0.5, 0.8535) {\small $m$};
    \end{scope}
}
\longrightarrow
\tikz[baseline={([yshift=-.5ex]current bounding box.center)}, scale=1]
{
    \begin{scope}[rotate=-45]
        \draw[dotted] (.5,.5) circle(0.707);
        \draw[knot] (0.5, -0.207) -- (0.5, 1.207);
        \draw[knot, overcross] (-0.207, 0.5) -- (1.207, 0.5);
        \draw[fill=white] (0.3, 0.7535) rectangle (0.7, 0.9535);
        \node[transform shape] at (0.5, 0.8535) {\small $m$};
        \draw[fill=white] (0.3, 0.0465) rectangle (0.7, 0.2465);
        \node[transform shape] at (0.5, 0.1465) {$c$};
    \end{scope}
}
\]
where $c = \mathrm{Cone}(\iota)$. Again by (CK1), $c$ has through degree $<m$, so it contains some turnback. Pushing the turnback through the parallel overstrands induces nontrivial $\mathscr{G}$-grading shifts (see Lemma \ref{lem:r2}), but after it passes through all $n$ overstrands, (CK2) tells us that that the entire complex on the right is contractible, and we're done.
\end{proof}

Using this Lemma, together with multiguling (Theorem \ref{thm:multigluing}) and idempotence (Proposition \ref{prop:unifiedprojectoruni}), $\Pi^\mathbf{m}$ can be described up to homotopy as sending 
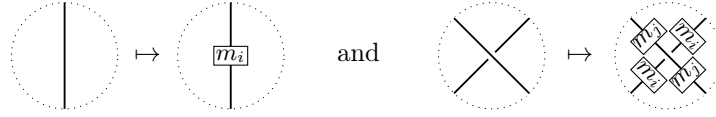
\[
\tikz[baseline={([yshift=-.5ex]current bounding box.center)}, scale=1]
{
        \draw[dotted] (.5,.5) circle(0.707);
        \draw[knot] (0.5, -0.207) -- (0.5, 1.207);
}
~\mapsto~
\tikz[baseline={([yshift=-.5ex]current bounding box.center)}, scale=1]
{
        \draw[dotted] (.5,.5) circle(0.707);
        \draw[knot] (0.5, -0.207) -- (0.5, 1.207);
        \draw[fill=white] (0.275, 0.375) rectangle (0.725, 0.625);
        \node[transform shape] at (0.5, 0.5) {\small $m_i$};
}
\qquad
\text{and}
\qquad
\tikz[baseline={([yshift=-.5ex]current bounding box.center)}, scale=1]
{
    \begin{scope}[rotate=-45]
        \draw[dotted] (.5,.5) circle(0.707);
        \draw[knot] (0.5, -0.207) -- (0.5, 1.207);
        \draw[knot, overcross] (-0.207, 0.5) -- (1.207, 0.5);
    \end{scope}
}
~\mapsto~
\tikz[baseline={([yshift=-.5ex]current bounding box.center)}, scale=1]
{
    \begin{scope}[rotate=-45]
        \draw[dotted] (.5,.5) circle(0.707);
        \draw[knot] (0.5, -0.207) -- (0.5, 1.207);
        \draw[knot, overcross] (-0.207, 0.5) -- (1.207, 0.5);
        \draw[fill=white] (0.3-0.02, 0.7535-0.02 + 0.01) rectangle (0.7+0.02, 0.9535+0.02+0.01);
        \node[transform shape] at (0.5, 0.8535+0.01) {\small $m_i$};
        \draw[fill=white] (0.753-0.02+0.01, 0.3-0.02) rectangle (0.9535+0.02+0.01, 0.7+0.02);
        \node[transform shape, rotate=90] at (0.8535+0.01, 0.5) {\small $m_j$};
        \draw[fill=white] (0.3-0.02, 0.0465-0.02-0.01) rectangle (0.7+0.02, 0.2465+0.02-0.01);
        \node[transform shape] at (0.5, 0.1465-0.01) {\small $m_i$};
        \draw[fill=white] (0.0465-0.02-0.01, 0.3-0.02) rectangle (0.2465+0.02-0.01, 0.7+0.02);
        \node[transform shape, rotate=90] at (0.1465-0.01, 0.5) {\small $m_j$};
    \end{scope}
}
\]
on the $i$th strand and each crossing of the $i$th strand under the $j$th.

\begin{theorem}
The category $\mathrm{Kom}_{\mathbf{m}}(n)^q$ contains invariants of framed tangles.
\end{theorem}

\begin{proof}
Applying $\Pi^\mathbf{m}$ to the following typical diskular 2-tangle and applying idempotence $(P_n \otimes P_n \simeq P_n$) and Lemma \ref{lem:projsliding}, we obtain
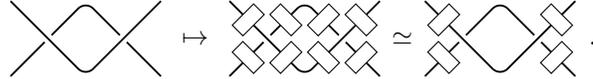
\[
\tikz[baseline={([yshift=-.5ex]current bounding box.center)}, scale=1]
{
    \begin{scope}[rotate=-45]
        \draw[knot, rounded corners] (0.5, -0.207) -- (0.5, 1.207) -- (1.914, 1.207);
        \draw[knot, overcross, rounded corners] (-0.207, 0.5) -- (1.207, 0.5) -- (1.207, 1.914);
    \end{scope}
}
~\mapsto~
\tikz[baseline={([yshift=-.5ex]current bounding box.center)}, scale=1]
{
    \begin{scope}[rotate=-45]
        \draw[knot, rounded corners] (0.5, -0.207) -- (0.5, 1.207) -- (1.914, 1.207);
        \draw[knot, overcross, rounded corners] (-0.207, 0.5) -- (1.207, 0.5) -- (1.207, 1.914);
        \draw[fill=white] (0.3, 0.7535) rectangle (0.7, 0.9535);
        \draw[fill=white] (0.753, 0.3) rectangle (0.9535, 0.7);
        \draw[fill=white] (0.3, 0.0465) rectangle (0.7, 0.2465);
        \draw[fill=white] (0.0465, 0.3) rectangle (0.2465, 0.7);
        \draw[fill=white] (1.007, 1.4605) rectangle (0.7+0.707, 0.9535+0.707);
        \draw[fill=white] (0.753+0.707, 0.3+0.707) rectangle (0.9535+0.707, 0.7+0.707);
        \draw[fill=white] (0.3+0.707, 0.0465+0.707) rectangle (0.7+0.707, 0.2465+0.707);
        \draw[fill=white] (0.0465+0.707, 0.3+0.707) rectangle (0.2465+0.707, 0.7+0.707);
    \end{scope}
}
\simeq
\tikz[baseline={([yshift=-.5ex]current bounding box.center)}, scale=1]
{
    \begin{scope}[rotate=-45]
        \draw[knot, rounded corners] (0.5, -0.207) -- (0.5, 1.207) -- (1.914, 1.207);
        \draw[knot, overcross, rounded corners] (-0.207, 0.5) -- (1.207, 0.5) -- (1.207, 1.914);
    %
        \draw[fill=white] (0.3, 0.0465) rectangle (0.7, 0.2465);
        \draw[fill=white] (0.0465, 0.3) rectangle (0.2465, 0.7);
        \draw[fill=white] (1.007, 1.4605) rectangle (0.7+0.707, 0.9535+0.707);
        \draw[fill=white] (0.753+0.707, 0.3+0.707) rectangle (0.9535+0.707, 0.7+0.707);
    \end{scope}
}~.
\]
Taking $\mathrm{Kh}$ (after picking any orientation), we know that
\[
\mathrm{Kh}\left(
\tikz[baseline={([yshift=-.5ex]current bounding box.center)}, scale=1]
{
    \begin{scope}[rotate=-45]
        \draw[knot, rounded corners] (0.5, -0.207) -- (0.5, 1.207) -- (1.914, 1.207);
        \draw[knot, overcross, rounded corners] (-0.207, 0.5) -- (1.207, 0.5) -- (1.207, 1.914);
    %
        \draw[fill=white] (0.3, 0.0465) rectangle (0.7, 0.2465);
        \draw[fill=white] (0.0465, 0.3) rectangle (0.2465, 0.7);
        \draw[fill=white] (1.007, 1.4605) rectangle (0.7+0.707, 0.9535+0.707);
        \draw[fill=white] (0.753+0.707, 0.3+0.707) rectangle (0.9535+0.707, 0.7+0.707);
    \end{scope}
}
\right)
\simeq
\varphi_{W^v}
\mathrm{Kh}\left(
\tikz[baseline={([yshift=-.5ex]current bounding box.center)}, scale=1]
{
    \begin{scope}[rotate=-45]
        \draw[knot, rounded corners=1em] (0.5, -0.207) -- (0.5, 0.4) -- (1.307,1.207) -- (1.914, 1.207);
        \draw[knot, rounded corners=1em] (-0.207, 0.5) -- (0.4, 0.5) -- (1.207,1.307) -- (1.207, 1.914);
    %
        \draw[fill=white] (0.3, 0.0465) rectangle (0.7, 0.2465);
        \draw[fill=white] (0.0465, 0.3) rectangle (0.2465, 0.7);
        \draw[fill=white] (1.007, 1.4605) rectangle (0.7+0.707, 0.9535+0.707);
        \draw[fill=white] (0.753+0.707, 0.3+0.707) rectangle (0.9535+0.707, 0.7+0.707);
    \end{scope}
}
\right)
\simeq
\varphi_{W^v}
\mathrm{Kh}\left(
\tikz[baseline={([yshift=-.5ex]current bounding box.center)}, scale=1]
{
    \begin{scope}[rotate=-45]
        \draw[knot, rounded corners=1em] (0.5, -0.207) -- (0.5, 0.4) -- (1.307,1.207) -- (1.914, 1.207);
        \draw[knot, rounded corners=1em] (-0.207, 0.5) -- (0.4, 0.5) -- (1.207,1.307) -- (1.207, 1.914);
    %
        \draw[fill=white] (0.0465, 0.3) rectangle (0.2465, 0.7);
    %
        \draw[fill=white] (0.753+0.707, 0.3+0.707) rectangle (0.9535+0.707, 0.7+0.707);
    \end{scope}
}\right)
\]
where $\varphi_{W^v}$ is the grading shift obtained by $m_im_j$ Reidemeister II moves (appeal to Lemma \ref{lem:r2} for an exact value, if desired). We have that $\deg_q(\varphi_{W^v}) = 0$ by Theorem \ref{thm:qinvt} which concludes the argument for the first framed tangle move. The argument for Reidemeister III moves is similar and left to the reader.
\end{proof}

If $L$ is a link, we denote by $\mathcal{H}(L;\mathbf{m})$ the homology of $\Pi^\mathbf{m}(L)$. Moreover, denote by $\Pi_e^\mathbf{m}(L)$ and $\Pi_o^\mathbf{m}(L)$ the complexes obtained from $\Pi^\mathbf{m}(L)$ by taking $X, Y, Z \mapsto 1$ and $X, Z \mapsto 1$ and $Y \mapsto -1$ respectively. These complexes are also invariants of the framed link $(L; \mathbf{m})$; denote their respective homology by $\mathcal{H}_e(L; \mathbf{m})$ and $\mathcal{H}_o(L; \mathbf{m})$. We write $\chi_q$ to denote the graded Euler characteristic which records only the $q$-grading associated to a particular $\mathscr{G}$-grading or $\mathscr{G}$-grading shift. By definition,
\[
\chi_q(\mathcal{H}_e(L; \mathbf{m})) = J(L;\mathbf{m})(q) = \chi_q(\mathcal{H}_o(L; \mathbf{m}))
\]
where $J(L;\mathbf{m})(q)$ denotes the colored Jones polynomial with indeterminate $q$. While $\mathcal{H}_e(L; \mathbf{m})$ is the colored link homology of \cite{https://doi.org/10.48550/arxiv.1005.5117}, $\mathcal{H}_o(L; \mathbf{m})$ provides a new categorification of the colored Jones polynomial of $L$. To verify that the two homologies are distinct, recall that the computation in \S \ref{sss:homologyoftrace} implies that $\mathcal{H}_e(U; 2) \not\cong \mathcal{H}_o(U; 2)$ for $U$ the unknot. 

\newpage

\bibliographystyle{amsalpha}
\bibliography{ref.bib}
\end{document}